\documentclass[11pt,oneside,english,final]{amsart}
\usepackage[T1]{fontenc}
\usepackage[notcite,notref]{showkeys}  
\usepackage{amstext,amsthm,amssymb,cases}
\usepackage[latin1]{inputenc}
\usepackage{geometry,xcolor}
\usepackage{mathrsfs,stmaryrd}
\usepackage{graphicx}
\usepackage{savesym}
\savesymbol{checkmark}
\usepackage{pifont}
\usepackage{multirow,textcase}

\usepackage{setspace}
\usepackage{xr}
\usepackage{enumitem}
\usepackage{bbm,bm}

\usepackage[colorlinks,bookmarks,linkcolor=black,citecolor=black]{hyperref} 
\usepackage[english,algoruled,lined,noresetcount,norelsize]{algorithm2e}

\usepackage{tikz}
\usetikzlibrary{shapes}
\usetikzlibrary{arrows.meta}

\newtheorem{theorem}{Theorem}
\newtheorem{lemma}[theorem]{Lemma}
\newcounter{LemmaDCounter}
\newcounter{LemmaECounter}
\newcounter{claimcounter}
\newcounter{claimcounterInLemmaD}
\newcounter{claimcounterInLemmaE}
\numberwithin{claimcounter}{theorem}
\numberwithin{claimcounterInLemmaD}{LemmaDCounter}
\numberwithin{claimcounterInLemmaE}{LemmaECounter}
\newtheorem{claim}[claimcounter]{Claim}
\newtheorem{claimInLemmaD}[claimcounterInLemmaD]{Claim}
\newtheorem{claimInLemmaE}[claimcounterInLemmaE]{Claim}
\newtheorem{proposition}[theorem]{Proposition}
\newtheorem{corollary}[theorem]{Corollary}
\newtheorem{fact}[theorem]{Fact}
\newtheorem{conjecture}[theorem]{Conjecture}

\newtheorem{setting}[theorem]{Setting}

\newtheorem*{InductiveStatementAst}{Inductive statement~$(\ast)_k$}

\newtheorem{definition}[theorem]{Definition}

\theoremstyle{remark}

\newcommand{\oldqed}{}
\newcommand{\qedClaim}{\hfill\scalebox{.6}{$\Box$}}
\newenvironment{claimproof}[1][Proof]{
	\renewcommand{\oldqed}{\qedsymbol}
	\renewcommand{\qedsymbol}{\qedClaim}
	\begin{proof}[#1]
	}{
	\end{proof}
	\renewcommand{\qedsymbol}{\oldqed}
} 

\def\sublem#1{\textrm{L}{\ref{#1}}}

\setlist{itemsep=2pt,parsep=1pt,topsep=3pt,partopsep=0pt}  
\setenumerate{leftmargin=*,labelindent=\parindent} 
\setdescription{style=sameline,leftmargin=.7cm,itemsep=2pt}

\def\itm#1{\rm ({#1})} 
\def\itmit#1{\itm{\it #1\,}} 
\def\rom{\itmit{\roman{*}}} 
\def\abc{\itmit{\alph{*}}}

\def\itmrom#1{\mbox{\itm{{\it #1.\roman{*}}\hspace{.1em}}}}
\def\itmarab#1{\mbox{\itm{{\it #1\,}\arabic{*}\hspace{.05em}}}}

\newcommand{\RandomEmbedding}{\emph{RandomEmbedding}}
\newcommand{\ModRandomEmbedding}{\emph{Modified RandomEmbedding}}
\newcommand{\PackingProcess}{\emph{PackingProcess}}

\newcommand{\RandomPathEmbedding}{\emph{RandomPathEmbedding}}
\newcommand{\PathPacking}{\emph{PathPacking}}
\newcommand{\AlgMap}{\hookrightarrow}
\newcommand{\deltnonspanning}{\delta}
\newcommand{\newdelta}{\tilde{\delta}} 
\newcommand{\newD}{\tilde{D}} 
\newcommand{\sigmKJ}{\lambda}
\newcommand{\sigmJnula}{\sigma_0}
\newcommand{\sigmJjedna}{\sigma_1}
\newcommand{\gamJ}{\gamma_\mathbf{F}}
\newcommand{\gamNew}{\gamma_\mathbf{E}}
\newcommand{\gamAnchor}{\gamma_\mathbf{D}}
\newcommand{\gamcore}{\gamma_\mathbf{A}}
\newcommand{\iniquasi}{\xi}
\newcommand{\DParity}{D_{\mathrm{odd}}}
\newcommand{\imA}{\im^\mathbf{A}}
\newcommand{\imB}{\im^\mathbf{B}}
\newcommand{\imC}{\im^\mathbf{C}}
\newcommand{\imCMinus}{\im^{\mathbf{C}^-}}
\newcommand{\imD}{\im^\mathbf{D}}
\newcommand{\imDMinus}{\im^{\mathbf{D}^-}}
\newcommand{\imE}{\im^\mathbf{E}}
\newcommand{\imF}{\im^\mathbf{F}}

\newcommand{\By}[2]{\overset{\mbox{\tiny{#1}}}{#2}}
\newcommand{\ByRef}[2]{   \By{\eqref{#1}}{#2} }
\newcommand{\eqBy}[1]{    \By{#1}{=} }

\newcommand{\leBy}[1]{    \By{#1}{\le} }
\newcommand{\geBy}[1]{    \By{#1}{\ge} }
\newcommand{\eqByRef}[1]{ \ByRef{#1}{=} }

\newcommand{\geByRef}[1]{ \ByRef{#1}{\ge} }

\DeclareMathOperator{\im}{im}

\def\JUSTIFY#1{\mbox{\fbox{\tiny#1}}\quad}

\newcommand{\eps}{\varepsilon}
\renewcommand{\rho}{\varrho}
\renewcommand{\subset}{\subseteq}
\newcommand{\dbigcup}{\dot\bigcup}
\newcommand{\dcup}{\dot\cup}

\newcommand{\NN}{\mathbb{N}}
\newcommand{\ZZ}{\mathbb{Z}}
\newcommand{\RR}{\mathbb{R}}

\newcommand{\cE}{\mathcal{E}}
\newcommand{\cF}{\mathcal{F}}
\newcommand{\cJ}{\mathcal{J}}
\newcommand{\cK}{\mathcal{K}}

\newcommand{\cD}{\mathcal{D}}
\newcommand{\cA}{\mathcal{A}}

\newcommand{\cG}{\mathcal{G}}

\newcommand{\cP}{\mathcal{P}}
\newcommand{\cQ}{\mathcal{Q}}

\newcommand{\cX}{\mathcal{X}}

\newcommand{\Vmin}{V_{\boxminus}}
\newcommand{\Vplus}{V_{\boxplus}}

\newcommand{\hist}{\mathscr{H}}
\newcommand{\histens}{\mathscr{L}}

\newcommand{\ONE}{\mathbbm{1}} 
\newcommand{\Exp}{\mathbf{E}} 
\newcommand{\Prob}{\mathbf{P}} 
\newcommand{\Var}{\mathbf{Var}}

\newcommand{\NBH}{\mathsf{N}}   
\newcommand{\NBHout}{\mathsf{N}^{\mathrm{out}}}   
\newcommand{\NBHin}{\mathsf{N}^{\mathrm{in}}}  
\newcommand{\degout}{\deg^{\mathrm{out}}}
\newcommand{\degin}{\deg^{\mathrm{in}}}
\newcommand{\LNBH}{\mathsf{N}^{-}} 
\newcommand{\CANDSET}{C}
\newcommand{\LEFTDEG}{\deg^{-}}
\newcommand{\OurPackingClass}{\mathsf{OurPackingClass}}
\newcommand{\SpecLeaves}{\mathsf{OddVert}}
\newcommand{\OddOut}{\mathrm{OddOut}}
\newcommand{\SpecPaths}{\mathsf{SpecPaths}}
\newcommand{\BasicPaths}{\mathsf{BasicPaths}}
\newcommand{\SpecShortPaths}{\mathsf{SpecShortPaths}}
\newcommand{\leftpath}{\mathsf{leftpath}}
\newcommand{\rightpath}{\mathsf{rightpath}}
\newcommand{\middlepath}{\mathsf{middlepath}}

\newcommand{\PathTerm}{\mathrm{PathTerm}}

\newcommand{\Paths}{\mathsf{Paths}}

\newcommand{\Parity}{\mathrm{Parity}}

\newcommand{\HStageA}{H_\mathbf{A}}
\newcommand{\HStageB}{H_\mathbf{B}}
\newcommand{\HStageC}{H_\mathbf{C}}
\newcommand{\HStageCMinus}{H_{\mathbf{C}^-}}
\newcommand{\HStageD}{H_\mathbf{D}}
\newcommand{\HStageE}{H_\mathbf{E}}
\newcommand{\HStageF}{H_\mathbf{F}}

\newcommand{\dStageA}{d_\mathbf{A}}
\newcommand{\dStageB}{d_\mathbf{B}}
\newcommand{\dStageC}{d_\mathbf{C}}
\newcommand{\dStageCMinus}{d_{\mathbf{C}^-}}
\newcommand{\dStageD}{d_\mathbf{D}}
\newcommand{\dStageE}{d_\mathbf{E}}
\newcommand{\dStageF}{d_\mathbf{F}}
\newcommand{\dStageEbar}{\bar{d}_{\mathbf{E}}}
\newcommand{\dStageFbar}{\bar{d}_{\mathbf{F}}}

\newcommand{\LC}{L_\mathbf{C}}
\newcommand{\LD}{L_\mathbf{D}}
\newcommand{\LE}{L_\mathbf{E}}
\newcommand{\LF}{L_\mathbf{F}}
\newcommand{\alphaCareful}{\delta}

\newcommand{\indexofvertex}{\bm{\iota}}

\DeclareRobustCommand{\SkipTocEntry}[5]{}  

\makeindex

\title{The tree packing conjecture for trees of almost linear maximum degree}

\author[P. Allen]{Peter Allen}
\address{(PA) London School of Economics, Department of Mathematics, Houghton Street, London WC2A 2AE, UK}
\thanks{(PA) Partially supported by EPSRC, grant number EP/P032125/1.}
\email{p.d.allen@lse.ac.uk}
\author[J. B\"ottcher]{Julia B\"ottcher}
\address{(JB) London School of Economics, Department of Mathematics, Houghton Street, London WC2A 2AE, UK}
\thanks{(JB) Partially supported by EPSRC, grant number EP/R00532X/1.}
\email{j.boettcher@lse.ac.uk}
\author[D. Clemens]{Dennis Clemens}
\address{(DC) Technische Universit\"at Hamburg, Institut f\"ur Mathematik, Am Schwarzenberg-Campus 3, 21073 Hamburg, Germany }
\email{dennis.clemens@tuhh.de}
\author[J. Hladk\'y]{Jan Hladk\'y}
\address{(JH, corresponding author) Institute of Computer Science of the Czech Academy of Sciences, Pod
	Vod\'arenskou v\v{e}\v{z}\'\i, 2, 182~07 Prague, Czechia. With institutional
	support RVO:67985807. \emph{This work was in part done while affiliated with:} Institute of Mathematics of the Czech Academy of Sciences, \v{Z}itn\'a 25, 115~67 Prague, Czechia. With institutional
	support RVO:67985840}
\thanks{(JH) Supported by the Czech Science Foundation, grant number 18-01472Y}
\email{hladky@cs.cas.cz}
\author[D. Piguet]{Diana Piguet}
\address{(DP) Institute of Computer Science of the Czech Academy of Sciences, Pod
	Vod\'arenskou v\v{e}\v{z}\'\i, 2, 182~07 Prague, Czechia. With institutional
	support RVO:67985807}
\thanks{(DP) Supported by the Czech Science Foundation, grant number 19-08740S}
\email{piguet@cs.cas.cz}
\author[A. Taraz]{Anusch Taraz}
\address{(AT) Technische Universit\"at Hamburg, Institut f\"ur Mathematik, Am Schwarzenberg-Campus 3, 21073 Hamburg, Germany}
\email{taraz@tuhh.de}

\date{\today}

\begin{document}
	
\begin{abstract}
  We prove that there is $c>0$ such that for all sufficiently large $n$,
  if $T_1,\dots,T_n$ are any trees such that $T_i$ has $i$ vertices and
  maximum degree at most $cn/\log n$, then $\{T_1,\dots,T_n\}$ packs into
  $K_n$. Our main result actually allows us to replace the host graph $K_n$
  by an arbitrary quasirandom graph, and to generalise from trees to
  graphs of bounded degeneracy that are rich in bare paths, contain some
  odd-degree vertices, and only satisfy much less stringent restrictions
  on their number of vertices.
\end{abstract}
	
\maketitle



\newpage
\tableofcontents
\newpage

\section{Introduction}
Let $G_1,\dots,G_t$ be a collection of graphs, and let $H$ be a graph. We say the family $\{G_1,\dots,G_t\}$ \emph{packs into $H$} if there are edge-disjoint copies of $G_1$,\dots, $G_t$ in $H$. The packing is called \emph{perfect} or \emph{exact} if $\sum_{i\in[t]}e(G_i)=e(H)$, so that every edge of $H$ is used exactly once.
The study of (perfect) packings is one of the oldest topics in graph theory. Indeed, the problem of the existence of designs --- one of the most fascinating questions of mathematics whose origins go back to the 19th century --- can be phrased as a perfect (hyper)graph packing problem. This problem was solved only recently, first by Keevash~\cite{Kee1}, and independently by Glock, K\"uhn, Lo and Osthus~\cite{GKLO:Designs}.

In this paper, we concentrate on packings with larger graphs. The two most influential conjectures in this area concern the packing of trees, one of which is the following.

\begin{conjecture}[Ringel's conjecture]
	For each $n\in\NN$ and for each tree $T$ of order $n+1$, we have that
	$2n+1$ copies of $T$ pack into the complete graph $K_{2n+1}$.
\end{conjecture}

Ringel's conjecture~\cite{Ringel} was stated in 1968, and for a long time was
only known to hold for very specific families of trees, such as stars, paths and
similar trees. The first general result on this conjecture was proved
in~\cite{BHPT}.  While the result in~\cite{BHPT} has several further restrictions,
the one we want to highlight here is that it is \emph{approximate}, by which we mean
that the total number of edges of the embedded trees must be at most
$\binom{n}{2}-\Omega(n^2)$. For obtaining a corresponding \emph{exact} result it remains to
remove the $\Omega(n^2)$ term, which is hard. Indeed, this gap in difficulty
between an approximate and an exact result is quite common in the area of
packing and is best illustrated by the increase in difficulty needed to get from
R\"odl's proof~\cite{RodlNibble} that approximate designs exist to the existence
of designs~\cite{GKLO:Designs,Kee1}.  Joos, Kim, K\"uhn and Osthus~\cite{JKKO}
proved that Ringel's conjecture holds exactly for large bounded degree
trees. Finally, Montgomery, Pokrovskiy and Sudakov~\cite{MPS:Ringel}, and later
Keevash and Staden~\cite{KS:Ringel}, proved that Ringel's conjecture holds for
all sufficiently large trees.

The result of Joos, Kim, K\"uhn and Osthus is much more general than just Ringel's conjecture, and in particular their result allows the packing of any collection of bounded degree trees $\{T_1,\dots,T_{2n+1}\}$ each on $n+1$ vertices into $K_{2n+1}$. The results of~\cite{KS:Ringel,MPS:Ringel} do not allow for such an extension.

The second influential conjecture in the area is the following.

\begin{conjecture}[Tree packing conjecture]
	For each $n\in\NN$ and for each family of trees $\left(T_s\right)_{s\in[n]}$,
	$v(T_s)=s$, we have that $\left(T_s\right)_{s\in[n]}$ packs
	into the complete graph $K_{n}$.
\end{conjecture}

Gy\'arf\'as~\cite{GyaLeh} formulated this conjecture in 1978, and again for a long time it was known only for specific path-like and star-like families. In 1983, Bollob\'as~\cite{MR710892} observed that the $\lfloor n/\sqrt{2}\rfloor$ smallest trees can be packed. Packing a few largest trees is already much more difficult. In particular, as a relaxation of the tree packing conjecture, Bollob\'as~\cite{MR1373679} asked in 1995 if for every $k\in \NN$ and $n>n_0(k)$, the $k$ largest trees $T_{n-k+1},\ldots,T_n$ can be packed into $K_n$. Until recently, this was known only for $k=5$ due to \.{Z}ak~\cite{ZAK2017252}, while
Balogh and Palmer~\cite{Balogh2013} showed that one can pack the largest $n^{1/4}/10$ trees into a somewhat larger $K_{n+1}$. After the first version of this manuscript was published, a breakthrough on Bollob\'as's question was achieved by Janzer and Montgomery in~\cite{janzer2024packinglargesttreestree}. Namely, they proved that the $cn$ largest trees can be packed, for some absolute constant $c>0$.

Again, the first general approximate result on the tree packing conjecture is~\cite{BHPT}, and the theorem of Joos, Kim, K\"uhn and Osthus~\cite{JKKO} mentioned above proves also that the tree packing conjecture holds for any family of bounded degree trees when $n$ is sufficiently large.

Moving away from trees, Messuti, R\"odl and Schacht~\cite{MRS} and Ferber, Lee and Mousset~\cite{FLM} considered approximate packings of bounded-degree graphs which are non-expanding in a suitable sense. Subsequently, Kim, K\"uhn, Osthus and Tyomkyn~\cite{KKOT} proved a packing version of the blow-up lemma, which in particular allows for approximate packings of general bounded-degree graphs. In another direction, moving away from bounded-degree graphs, Ferber and Samotij~\cite{FerSam} proved an approximate version of the tree packing conjecture for trees of maximum degree $O(n/\log n)$. We should remark that both these results apply in more generality than packing into $K_n$. Indeed, \cite{KKOT} in fact allows for packing into a Szemer\'edi partition, while~\cite{FerSam} works also in sparse random graphs.

 For packings into complete graphs, the following result from~\cite{DegPack} generalises both~\cite{FerSam,KKOT}.\footnote{Version~2 on the arXiv of~\cite{DegPack} appeared in a journal,~\cite{DegPack_Advances}. The current arXiv version~4 of~\cite{DegPack} corrects an error in the journal version (described in~\cite[Section 1.1]{DegPack}). Numbering of theorems and lemmas in these versions is identical.}
 A graph is said to be \emph{$D$-degenerate}\index{degenerate} if there is an ordering of its vertices such that each vertex has at most $D$ neighbours preceding it in the order. Many interesting families of graphs are degenerate---for example, trees are $1$-degenerate and planar graphs are $5$-degenerate.

\begin{theorem}[{\cite[Theorem~\ref{PAPERpackingdegenerate.thm:MAINunbounded}]{DegPack}}]
  \label{thm:ABHPDegenerateNonTechnical}
  For every $D\in\NN$ and $\eta>0$,
  there exist $n_{0}\in\NN$ and $c>0$ so that for each $n>n_{0}$
  the following holds. Suppose that $\left(G_s\right)_{s\in\mathcal{F}}$
  is a family of $D$-degenerate graphs of orders at most $n$ and maximum
  degrees at most $\frac{cn}{\log n}$, whose total number of edges is
  at most $\binom{n}{2}-\eta n^2$.
  Then $\left(G_s\right)_{s\in\mathcal{F}}$ packs into $K_{n}$.
\end{theorem}

Building on this, in~\cite{ABCT:PackingManyLeaves} a perfect packing result for
degenerate graphs was obtained with the additional condition that many of the
graphs are non-spanning and contain linearly many leaves (we state a generalisation of this result in Theorem~\ref{thm:ABCT}).
Here, a \emph{leaf}\index{leaf} in a graph is a vertex of degree~1.
As observed in~\cite{ABCT:PackingManyLeaves}, this result implies that the tree packing conjecture holds for almost all families of trees. However, observe that a tree necessarily either contains many leaves or many short bare paths. Here, a
subset $U$ of vertices of a graph $G$ induces a \emph{bare path}\index{bare
  path} if $G[U]$ is a path and for each vertex $u\in U$ we have
$\deg_{G}(u)=2$.
Hence to prove the tree packing conjecture, at least for trees with maximum degree at most $cn/\log n$, it only remains to consider trees with many short bare paths.
This paper resolves this case.

We remark that many results in this area allow for packing into more general graphs than~$K_n$. One can also pack into quasirandom graphs, which we now define. Given a graph $H$ with $n$ vertices, its \emph{density}\index{density} is the number ${e(H)}/{\binom{v(H)}{2}}$. Given a bipartite graph $K$ with parts of sizes $a$ and $b$, its \emph{bipartite density}\index{bipartite density} is the number ${e(K)}/ab$. Suppose that $v\in V(H)$ and $S\subset V(H)$. The \emph{neighbourhood}\index{neighbourhood} of $v$ is denoted by $\NBH_H(v)$\index{$\NBH_H(v)$, $\NBH_{H}(S)$, $\NBH_{H}(v;A)$, $\NBH_{H}(S;A)$} and the \emph{common neighbourhood} of~$S$ by $\NBH_{H}(S)=\bigcap_{v\in S}\NBH_{H}(v)$. Here, the convention is that $\NBH_{H}(\emptyset)=V(H)$. For $A\subset V(H)$, we write $\NBH_H(v;A)=\NBH_H(v)\cap A$ and $\NBH_H(S;A)=\NBH_H(S)\cap A$. We write $\deg_H(v):=|\NBH_H(v)|$\index{$\deg_H(v)$, $\deg_{H}(S)$, $\deg_{H}(v;A)$} and $\deg_H(S):=|\NBH_H(S)|$. For $A\subset V(H)$ we further write $\deg_H(v;A):=|\NBH_H(v;A)|$ for the number of neighbours of $v$ inside~$A$.
  \begin{definition}[quasirandom]\label{def:quasirandomness}
Suppose that $L\in\NN$ and $\gamma>0$.
Suppose that $H$ is a graph with $n$ vertices and with density~$p$. We say that $H$ is \emph{$(\gamma,L)$-quasirandom}\index{quasirandom} if for every set $S\subset V(H)$ of at most $L$ vertices we have
	$|\NBH_{H}(S)|=(1\pm\gamma)p^{|S|}n$.
\end{definition}

Theorem~\ref{thm:ABHPDegenerateNonTechnical} allows more generally packings in
quasirandom graphs. The same is true of the result of Joos, Kim, K\"uhn and
Osthus~\cite{JKKO}, and Keevash and Staden~\cite{KS:Ringel} formulate and prove
a version of Ringel's conjecture for quasirandom graphs. The following main result
from~\cite{ABCT:PackingManyLeaves} also handles quasirandom graphs.

\begin{theorem}[{\cite[Theorem~\ref{PAPERmanyleaves.thm:main}]{ABCT:PackingManyLeaves}}]
  \label{thm:ABCT}For every $D\in\NN$ and $d,\alpha>0$, there exist
	$n_{0},L\in\NN$ and $c,\iniquasi>0$ so that for each $n>n_{0}$ the
  following holds. Suppose that~$H$ is a $(\iniquasi,L)$-quasirandom graph of
  order~$n$ with at least $dn^2$ edges and that
  $\left(G_s\right)_{s\in\mathcal{F}}$ is a family
  of $D$-degenerate graphs of orders at most $n$, maximum degrees at most
  $\frac{cn}{\log n}$, and total number of edges at most $e(H)$. Suppose further that there exists an index set
  $\mathcal{B}\subset\mathcal{F}$ such that
  \begin{itemize}
  \item for each $s\in\mathcal{B}$ we have $v(G_s)\le(1-\alpha)n$, and
  \item the total number of leaves in the family $\left(G_s\right)_{s\in\mathcal{B}}$
    is at least $\alpha n^{2}$.
  \end{itemize}
  Then $\left(G_s\right)_{s\in\mathcal{F}}$ packs into~$H$.
\end{theorem}

\addtocontents{toc}{\SkipTocEntry}
\subsection{Our results}
Our main result, Theorem~\ref{thm:maintechnical}, states that we
can perfectly pack graphs from families we introduce in Definition~\ref{def:family} into a quasirandom graph. This definition is technical and tailored to the maximal possible generality allowed by our methods. However, combined with Theorem~\ref{thm:ABCT}, it allows us to prove that the tree packing conjecture holds for trees with maximum degree $O(n/\log n)$.

\begin{theorem}
  \label{thm:TreePackingForLargeDeg}
  There exist $c>0$ and $n_{0}\in\NN$
  such that for each $n>n_{0}$ any family of trees $\left(T_s\right)_{s\in[n]}$ with
  $v(T_s)=s$ and $\Delta(T_s)\le\frac{cn}{\log n}$ packs into $K_{n}$.
\end{theorem}

Similarly, we obtain an analogue of Ringel's conjecture for trees with degrees bounded by $O(n/\log n)$, where different trees are allowed.

\begin{theorem}
  \label{thm:RingelOurs}
  There exist $c>0$ and $n_{0}\in\NN$
  such that for each $n>n_{0}$ any family of trees $\left(T_s\right)_{s\in[2n+1]}$ with
  $v(T_s)=n+1$ and $\Delta
  (T_s)\le\frac{cn}{\log n}$ packs into $K_{2n+1}$.
\end{theorem}

In fact, Theorem~\ref{thm:TreePackingForLargeDeg} and
Theorem~\ref{thm:RingelOurs} both follow easily from the following more
general packing result for trees. Both this result and the two deductions are
given in Section~\ref{sec:reductionofMainTheorems}.

\begin{theorem}\label{thm:PackingTrees}
  For each $\delta,d>0$ there exist $c,\iniquasi>0$ and $n_{0},L\in\NN$ such that for
  each $n>n_{0}$ and any $(\iniquasi,L)$-quasirandom graph~$H$ with~$n$ vertices and at least $dn^2$ edges the
  following holds.
  Any family of trees $\left(T_s\right)_{s\in[N]}$ with $N\ge\left(\frac12+\delta\right)n$ satisfying
  \begin{enumerate}[label=\abc]
  \item $\sum_{s\in[N]} e(T_s)\le e(H)$,
  \item $\Delta(T_s)\le\frac{cn}{\log n}$ for all $s\in[N]$,
  \item $\delta n\le v(T_s)\le (1-\delta)n$ for all $1\le s\le(\frac12+\delta)n$ and $v(T_s)\le n$ for all $(\frac12+\delta)n<s\le N$, 
  \end{enumerate}
  packs into~$H$.
\end{theorem}

We now introduce the class of graphs we can pack in our main result.  The
\emph{length}\index{length of path} of a path is the number of its edges.

\begin{definition}[$\OurPackingClass$]
  \label{def:family}\mbox{}\\
  Given $n,m,D_0\in\NN$ and $\delta,c>0$,
  let $\OurPackingClass(n,m;\delta,c,D_0)$\index{$\OurPackingClass(n,m;\delta,c,D_0)$} be
  the set of all families $\left(G_s\right)_{s\in\cG}$ of
  graphs for which there are disjoint index sets \index{$\cK$}\index{$\cJ$}$\cK,\cJ\subset\cG$ with $\delta n\le |\cJ|\le n$, and an odd number
  \index{$\DParity$}$\DParity\le D_0$ such that
  \begin{enumerate}[label=\abc]
  \item \label{enu:familyEachG} for each $s\in\cG$, the graph~$G_s$ has $v\left(G_s\right)\le n$ vertices, maximum degree
    $\Delta\left(G_s\right)\le\frac{c n}{\log n}$, and is~$D_0$-degenerate,
  \item \label{enu:familyEdgesAll}$\sum_{s\in\cG}e\left(G_s\right)\le m$,
  \item \label{enu:familyOrderJK}for each $s\in\cJ\cup\cK$, we have $v\left(G_s\right)\le\left(1-\delta\right)n$,
  \item \label{enu:familyPathsJ} for each $s\in\cJ$, the graph $G_s$ contains
    a family $\BasicPaths_s$\index{$\BasicPaths_s$} of $\alphaCareful n$ vertex-disjoint bare paths of length $11$,
  \item\label{enu:defSpecLeaves} for each $s\in\cK$ there is a non-empty
    independent set \index{$\SpecLeaves_s$}$\SpecLeaves_s$ of vertices of $G_s$
    whose degree is $\DParity$ with $|\SpecLeaves_s|\le\tfrac{cn}{\log n}$, and such that
    $\sum_{s\in\cK}|\SpecLeaves_s|\ge(1+\delta)n$.
  \end{enumerate}
\end{definition}
To summarise, Definition~\ref{def:family} allows families of graphs that 
\begin{itemize}
\item by~\ref{enu:familyEachG} may be spanning with respect to the host graph $H$ into which the family is to be embedded, if $n=v(H)$,
  have degeneracy bounded by a constant and maximum degrees bounded by $O(n/\log n)$,
\item by~\ref{enu:familyEdgesAll}
 may in total have the same number of edges as the host graph $H$, if $m=e(H)$, 
\item by~\ref{enu:familyOrderJK} and~\ref{enu:familyPathsJ} contain a reasonably sized collection of non-spanning graphs, each of which
  contains linearly many constant-length bare paths, and
\item by~\ref{enu:familyOrderJK} and~\ref{enu:defSpecLeaves} contain a collection of non-spanning graphs, each of which contains at
  least one vertex of odd and not too large degree (which we will use to correct
  parities), not too many of which are in any one graph and which total slightly
  more than $n$.
\end{itemize}

Our main result states that any family of graphs from this class can be
perfectly packed into any dense and sufficiently quasirandom graph.

\begin{theorem}[main result]
	\label{thm:maintechnical}For every $D_0\in\NN$ and $\deltnonspanning,d>0$,
	there exist $n_{0},L\in\NN$ and $c,\iniquasi>0$ so that for each
	$n>n_{0}$ the following holds. Suppose that $H$ is a $(\iniquasi,L)$-quasirandom
	graph of order $n$ with at least $dn^{2}$ edges. Then each family of graphs from
	$\OurPackingClass(n,e(H);\deltnonspanning,c,D_0)$ packs into~$H$.
\end{theorem}

\addtocontents{toc}{\SkipTocEntry}
\subsection{Optimality}\label{ssec:opt}

Let us now briefly discuss the optimality of Theorem~\ref{thm:maintechnical}.
Firstly, we cannot relax the maximum degree condition in
Definition~\ref{def:family}\ref{enu:familyEachG}. Indeed, Koml\'os, S\'ark\"ozy
and Szemer\'edi show in~\cite{Komlos2001} that if we connect a disjoint union of
$\frac{\log n}{C}$ many stars of orders $\frac{Cn}{\log n}$ by a path going
through the centres of these stars, then we get a tree which asymptotically
almost surely does not appear in the random graph $\mathbb{G}(n,p)$ (here,
$p\in(0,1)$ is arbitrary, and $C$ is sufficiently large, depending on $p$). Of
course, $\mathbb{G}(n,p)$ asymptotically almost surely satisfies our
quasirandomness condition. This example shows that not only packing but even
finding a single copy of such a tree is impossible. This also shows the optimality of the
maximum degree condition in Theorem~\ref{thm:PackingTrees}.
Our maximum degree condition is apparently not optimal for the now-resolved Ringel's conjecture, and it is conjecturally not optimal for the tree packing conjecture either. We do not see an easy way to improve our result in this specific setting. The limitation that prohibits us from also packing trees with vertices of degree $\omega\left(\frac{n}{\log n}\right)$ lies in the randomised packing algorithm underlying our proof. Specifically, the concentration properties on which our algorithm relies are not strong enough in this regime, and some vertices of the host graph would be assigned substantially more such high-degree vertices than others. As a result, these heavily loaded host vertices would create a bottleneck for the remainder of any packing. This is a well-known limitation of similar methods, described for example in~\cite[Section~6]{FerSam}.

Secondly, some of the graphs to be packed must be non-spanning, even if we only pack trees. So,
condition~\ref{enu:familyOrderJK} in Definition~\ref{def:family} cannot be omitted. Indeed, in~Section~9.1 of~\cite{BHPT}, a family of bounded-degree trees of orders $n$ and total number of edges ${n \choose 2}$ is given that does not pack into $K_{n}$. Moving away from bounded degree graphs, suppose $\cG$ is the following family of graphs. For any sufficiently small $c>0$, we put $\tfrac13n$
stars with $\tfrac{cn}{\log n}$ leaves into $\cG$. We let the remaining graphs in $\cG$ be long paths on at most $n$ vertices, such that in total we have $\binom{n}{2}$ edges.
Suppose now that there is a packing of $\cG$ into $K_n$. Let $H$ be the subgraph of edges used by the stars. Consider the at least $\tfrac23n$ vertices to which no star centre is embedded. These vertices form an independent set in $H$ and are in total adjacent to at most $\tfrac{cn^2}{3\log n}$ edges, so in particular one must have degree at most
$\tfrac{cn}{2\log n}$ in $H$. Thus in order to have a perfect packing, $\cG$
must contain at least $\tfrac{n-1}{2}-\tfrac{cn}{4\log n}$ paths. If the
average length of these paths is $\ell$, then we have
$\binom{n}{2}=\tfrac{cn^2}{3\log n}+\ell\big(\tfrac{n-1}{2}-\tfrac{cn}{4\log
  n}\big)$ from which we conclude $\ell\le n-\tfrac{c n}{6\log n}$. In particular, $\Omega(n/\log n)$ of the paths must be $\Omega(n/\log n)$-far from spanning. This shows that our requirement on the number of non-spanning graphs
with many bare paths is sharp up to a log factor even for packing into $K_n$
(and the same argument works in any other regular graph). For packing into a
quasirandom $H$, our conditions permit one vertex of~$H$ to have $\Omega(n)$ more
neighbours than the average degree. Letting $\cG$ be a family of long
paths, a similar argument tells us that the average length of the long paths
is $n-\Omega(n)$, so that in particular $\Omega(n)$ of the graphs in $\cG$ must
be $\Omega(n)$-far from spanning, so that in this more general setting our
conditions are sharp up to the value of $\deltnonspanning$.

\addtocontents{toc}{\SkipTocEntry}
\subsection{Organisation}\label{ssec:org}

The remainder of this paper is structured as follows. In
Section~\ref{sec:reductionofMainTheorems} we deduce
Theorem~\ref{thm:PackingTrees} from our main result. In
Section~\ref{sec:proofsketch} we provide a sketch of the proof of our main
theorem, Theorem~\ref{thm:maintechnical}. In Section~\ref{sec:designs} we
provide a decomposition result that follows from the deep results on the
existence of designs by Keevash~\cite{Kee2, KeeColouredDirected} and that we shall use in the final stage of our packing. In Section~\ref{sec:tools} we collect some probabilistic tools and facts about certain quasirandom properties that we use in our proof. In Section~\ref{sec:theproof} we provide the main lemmas covering
the different stages (Stage~A--Stage~G) of our proof and show how they imply
Theorem~\ref{thm:maintechnical}. In
Section~\ref{sec:pathpack} we provide an auxiliary non-perfect packing result
for anchored paths that we need in Stages~D and~E of our proof.  In Sections
\ref{sec:StageA}--\ref{sec:StageG} we give the proofs of the lemmas for Stages~A--G. We finish the paper with some concluding remarks in
Section~\ref{sec:concl}.

Throughout the paper, we write $\dcup$ for disjoint union of sets.

\section{Deducing the tree packing results from the main theorem}
\label{sec:reductionofMainTheorems}

In this section, we first deduce Theorem~\ref{thm:PackingTrees} from our main result, Theorem~\ref{thm:maintechnical}, and Theorem~\ref{thm:ABCT}, and then deduce Theorems~\ref{thm:TreePackingForLargeDeg} and~\ref{thm:RingelOurs} from Theorem~\ref{thm:PackingTrees}.

\begin{proof}[Proof of Theorem~\ref{thm:PackingTrees}]
  Given~$\delta$ and~$d$, let $\deltnonspanning'=\delta/1000$. We choose $c,\iniquasi>0$ sufficiently small and $n_0,L$ sufficiently
  large for both Theorem~\ref{thm:maintechnical} with input $D_0=1$,
  $\deltnonspanning'$ and $d$, and Theorem~\ref{thm:ABCT} with input
  $D=1$, $d$ and $\alpha=(\deltnonspanning')^2$.

  Given~$H$ and trees $T_1,\dots,T_N$ satisfying the conditions of
  Theorem~\ref{thm:PackingTrees}, we distinguish two cases.  First, suppose that
  among the indices $\{1,\dots,(\frac12+\delta)n\}$ there is a subset~$\mathcal{B}$ of
  size $\deltnonspanning' n$, such that each tree $T_i$ with $i\in\mathcal B$
  contains at least $\deltnonspanning' n$ leaves. Then we apply
  Theorem~\ref{thm:ABCT}, with constants as above and this $\mathcal{B}$, and it
  returns a perfect packing of $(T_i)_{i\in[N]}$ into~$H$.

  Second, suppose no such $\mathcal{B}$ exists.
  Then there must be a set $\cJ$ of indices $i\le(\frac12+\delta)n$ such that $T_i$ has less than
  $\deltnonspanning' n$ leaves, with $|\cJ|=\deltnonspanning' n$.
  Given any $i\in\cJ$ such
  that $T_i$ has less than $\deltnonspanning' n$ leaves, observe that the sum of
  the degrees of $T_i$ is $2v(T_i)-2$. Since at most $\deltnonspanning' n$
  vertices have degree~$1$ and all other vertices have degree at least~$2$, we
  see that at most $\deltnonspanning' n$ vertices have degree exceeding~$2$
  in~$T_i$. Now let~$\tilde E$ be the set of edges in~$T_i$ that contain at
  least one vertex that is not of degree~$2$ in~$T_i$. We have
  \[|\tilde E|\le \sum_{\substack{x\in V(T_i),\\ \deg(x)=1}}1+\sum_{\substack{x\in V(T_i),\\ \deg(x)>2}}\deg(x)
  =2v(T_i)-2-\sum_{\substack{x\in V(T_i),\\ \deg(x)=2}}2
  \le 2v(T_i)-2-2\big(v(T_i)-2\deltnonspanning' n\big)
  \le 4\deltnonspanning' n.
  \]
  Removing all the edges in~$\tilde E$ from~$T_i$, we
  obtain a graph $F$ with at most $4\deltnonspanning' n+1\le 5\deltnonspanning'
  n$ components. The total number of vertices in components with less than $50$
  vertices is at most $250\deltnonspanning' n$, so the remaining at least
  $\delta n-250\deltnonspanning' n\ge 100\deltnonspanning' n$ vertices all lie in
  components with at least $50$ vertices. Note that each such component is the vertex set of a bare path in $T_i$ and that there are at most
  $5\deltnonspanning' n$ such components. From the components with at
  least $50$ vertices, we choose pairwise vertex-disjoint paths of
  length $11$ until we have a set $\BasicPaths_i$ of $\deltnonspanning' n$
  vertex-disjoint bare paths in $T_i$ of length $11$.
  Note that we can simply choose such paths by 
  splitting each component of $F$ into vertex-disjoint subpaths of length 11 and at most one subpath with at most 11 vertices,
  and then greedily choosing
  such subpaths of length 11
  until $\deltnonspanning' n$ paths are chosen.
  Also note that we can indeed reach the desired number $\deltnonspanning' n$ of bare paths, since in each of the at most $5\deltnonspanning' n$
  components we could choose subpaths until we leave a leftover of at most 11 vertices, 
 leading to a total of
 at most $55\deltnonspanning' n$ leftover vertices.

  We let $\cK$ be an arbitrary subset of $[(\frac12+\delta)n]\setminus \cJ$ of
  size $(\frac12+\frac12\delta)n$. For each $T_i$ with $i\in\cK$, let
  $\SpecLeaves_i$ be a set of two leaves in $T_i$. Set $\DParity=1$. Observe
  that $\sum_{i\in\cK}|\SpecLeaves_i|=(1+\delta)n\ge(1+\deltnonspanning')n$ and
  hence $(T_i)_{i\in[N]}$ is in
  $\OurPackingClass(n,e(H);\deltnonspanning',c,D_0)$.  Therefore, we can
  apply Theorem~\ref{thm:maintechnical} which returns a perfect packing of
  $(T_i)_{i\in[N]}$ into~$H$.
\end{proof}

We can now deduce our two tree packing results. Note that a complete graph $K_m$ has density~$1$ and satisfies $|\NBH_{K_m}(S)|=m-|S|$ for every $S\subseteq V(K_m)$, so that for each $\iniquasi>0$ and $L\in\NN$ it is $(\iniquasi,L)$-quasirandom once $m\ge L/\iniquasi$; moreover $e(K_m)=\binom m2\ge\tfrac14m^2$ for $m\ge2$. We may therefore apply Theorem~\ref{thm:PackingTrees} with $H=K_m$ and $d=\tfrac14$.

\begin{proof}[Proof of Theorem~\ref{thm:RingelOurs}]
  We apply Theorem~\ref{thm:PackingTrees} with $\delta=\tfrac14$ and $d=\tfrac14$, and with $m:=2n+1$ in the role of its~$n$; let $c$ be the constant it returns. The host graph is $H=K_m$. Since $\sum_{s\in[m]}e(T_s)=(2n+1)n=\binom{2n+1}{2}=e(H)$, condition~(a) holds with equality. Every tree of the family has $n+1$ vertices, and $\tfrac14m\le n+1\le\tfrac34m$, so condition~(c) holds for every index; also $N=m\ge\big(\tfrac12+\delta\big)m$. Finally, as $x\mapsto\tfrac{x}{\log x}$ is increasing for $x>e$, we have $\tfrac{cn}{\log n}\le\tfrac{cm}{\log m}$, so condition~(b) holds as well. Theorem~\ref{thm:PackingTrees} therefore packs $(T_s)_{s\in[m]}$ into $K_{2n+1}$.
\end{proof}

\begin{proof}[Proof of Theorem~\ref{thm:TreePackingForLargeDeg}]
  We apply Theorem~\ref{thm:PackingTrees} with $\delta=\tfrac18$ and $d=\tfrac14$; let $c$ be the constant it returns, and let $H=K_n$. Conditions~(a) and~(b) are immediate, since $\sum_{s\in[n]}e(T_s)=\sum_{s\in[n]}(s-1)=\binom n2=e(H)$, and $N=n\ge\big(\tfrac12+\delta\big)n$.

  Condition~(c), on the other hand, requires the first $\big(\tfrac12+\delta\big)n$ trees of the family to have order between $\delta n$ and $(1-\delta)n$, which is not the case for the given indexing. We therefore reorder the family: there are $(1-2\delta)n\ge\big(\tfrac12+\delta\big)n$ indices $s$ with $\delta n\le s\le(1-\delta)n$, and we let the corresponding trees occupy the first $\big(\tfrac12+\delta\big)n$ positions, in an arbitrary order, followed by the remaining trees in an arbitrary order. Since $v(T_s)\le n$ for every~$s$, condition~(c) now holds. As whether a family packs into~$H$ does not depend on the order in which the family is listed, Theorem~\ref{thm:PackingTrees} gives the required packing of $(T_s)_{s\in[n]}$ into $K_n$.
\end{proof}

\section{Outline of the proof of our main theorem}
\label{sec:proofsketch}

In this section, we give a rough sketch of our proof of Theorem~\ref{thm:maintechnical}. We will give a much more detailed discussion in Section~\ref{sec:theproof}, together with precise statements of the lemmas we need along the way. We note that Table~\ref{tab:stages} therein may already help to give an overview of the embedding stages
that we sketch in the following.

We need to embed the graphs $(G_s)_{s\in\cG}$ into $H$; we can without loss of generality suppose this will be a perfect packing, i.e.~\ref{enu:familyEdgesAll} of Definition~\ref{def:family} holds with equality,
as we could otherwise add
single edges as graphs to the 
given family and then pack everything.
 We first embed all the graphs $G_s$ with $s\in\cG\setminus(\cJ\cup\cK)$, and most vertices of all the remaining graphs, using a randomised packing algorithm from~\cite{DegPack}; this is the \emph{bulk embedding}, Stage~A.

What remains is the following: for the graphs $G_s$ with $s\in\cK$ we still need to embed (some of) the vertices $\SpecLeaves_s$. For the graphs $G_s$ with $s\in\cJ$ we still need to embed some bare paths contained in $\BasicPaths_s$, whose ends are already embedded (we say the paths are \emph{anchored}).

We next complete the embedding of the graphs $G_s$ with $s\in\cK$. Since, after this step, only bare paths from $(\BasicPaths_s)_{s\in\cJ}$ --- that is, vertices of degree~2 --- will be left to pack, there are some obvious parity restrictions ahead. So, we use the odd degrees of $\SpecLeaves_s$ (cf.~Definition~\ref{def:family}\ref{enu:defSpecLeaves}) to fulfil these. This is the \emph{parity correction}, Stage~B.

Now we need to embed the anchored bare paths of $G_s$ with $s\in\cJ$. At each
vertex $v$ of $H$, we have to use one edge per path anchored at $v$, and in
addition each time we choose to use $v$ as an interior vertex we use two edges;
this is why we needed to set the parity in the previous step. We split these
graphs up into three parts, $\cJ=\cJ_0\dcup\cJ_1\dcup\cJ_2$, and perform the
embedding of these bare paths over several stages, preparing for an application
of the results of Keevash~\cite{Kee2,KeeColouredDirected} on the existence of
designs in our final stage.  We now first explain in detail the setup we will use to
apply the design results in order to pack (parts of) our bare paths. We then
outline how we can apply the design results in this setup. Finally, we explain
which preparations are performed for achieving this setup.

Before the final stage we will be left with a subgraph of~$H$;
let us call this subgraph~$H^*$ here. We will ensure that~$H^*$ has an even
number of vertices, which come in \emph{terminal pairs}\index{terminal pair}
$\{\boxminus_i,\boxplus_i\}$. Thus we have a partition $V(H^*)=\Vmin\dcup\Vplus$
into equal parts, which we call \emph{sides}\index{side}.  We will also ensure that all that
remains is to embed some paths with three edges from some graphs $G_s$ with
$s\in\cJ_0$ such that each of these paths is anchored to a terminal pair: that
is, if $xyzw$ is such a $3$-edge path in $G_s$, then there exists some $i$ such
that $x$ is embedded to the vertex $\boxminus_i$ of $H^*$ and $w$ to
$\boxplus_i$. When packing these paths in the final stage, we shall insist
that~$y$ gets embedded to the same side as $x$, i.e.\ $\Vmin$, and $z$ to the
same side as $w$, i.e.\ $\Vplus$.

Let us now explain how we apply the design results.  We represent the embedding
of the path $xyzw$ into~$H^*$ as the embedding of a $4$-vertex configuration
called a \emph{diamond}\index{diamond} in an auxiliary coloured and partially directed multigraph
called a \emph{chest}\index{chest}, which we now describe. (For an illustration see also
Figure~\ref{fig:ChestsAndDiamonds}.) The vertices of the chest come in two
parts, the set $V=\{1,2,\dots,|\Vmin|\}$ and the set $U=\cJ_0$. We put coloured
edges into $V$ as follows (see Figure~\ref{fig:outline} for an illustration). For each edge $\boxminus_i\boxminus_j\in E(H^*)$, we
put a blue edge $ij$. For each edge $\boxplus_i\boxplus_j$, we put a red edge
$ij$. Finally, for each edge $\boxminus_i\boxplus_j$ we put a green arc directed
from $i$ to $j$. Thus a $3$-vertex cycle $ijk$ in the chest, in which $ij$ is
blue, $jk$ is green and directed from $j$ to $k$, and $ki$ is red, represents a
$3$-edge path $\boxminus_i\boxminus_j\boxplus_k\boxplus_i$ in~$H^*$. In addition, we put edges
from $s\in U$ to $V$ as follows. For each $i\in V$ such that
$G_s$ has a path anchored on $\{\boxminus_i,\boxplus_i\}$ we put a grey edge
$is$. If $j\in V$ is such that $\boxminus_j$ has not been used in the embedding
of $G_s$ before the final stage, we put a black edge $js$, and if $\boxplus_k$
has not been used in the embedding of $G_s$ before the final stage, we put a
purple edge $ks$. A \emph{diamond} is then the $4$-vertex configuration obtained
by adding a vertex $s\in U$ to the $3$-vertex cycle $ijk$ described above,
with $is$ grey, $js$ black and $ks$ purple. A copy of this configuration in the
chest represents a $3$-edge path $\boxminus_i\boxminus_j\boxplus_k\boxplus_i$ in
$H^*$, with the additional information that we can use this copy to embed a
$3$-edge path of $G_s$ anchored at $\{\boxminus_i,\boxplus_i\}$ and that
$\boxminus_j$ and $\boxplus_k$ have not been used in the embedding of $G_s$. In
other words, a diamond in the chest represents a valid way to embed one path of
one graph $G_s$ in~$H^*$. See Figure~\ref{fig:outline} for an illustration of
this translation from paths to diamonds.

\begin{figure}
  
\definecolor{xred}{rgb}{1,0.2,0.2}
\definecolor{xblue}{rgb}{0.3,0.5,1}
\definecolor{xgreen}{rgb}{0.4,0.8,0.2}
\definecolor{xpurple}{rgb}{0.7,0.2,1}
\definecolor{xgray}{rgb}{0.7,0.7,0.7}
\definecolor{xblack}{rgb}{0,0,0}
  
\tikzset{  vertex/.style={
    circle,
    inner sep=2pt,
    fill=black
  },
 boxedge/.style={
    shorten < = -1mm,
    shorten > = -1mm,
  },  
  >=Stealth
}

\begin{tikzpicture}[auto,thick]

\begin{scope}[scale=.8]
  \draw (0,0) node (m5) {$\boxminus_5$};
  \draw (0,1) node (m4) {$\boxminus_4$};
  \draw (0,2) node (m3) {$\boxminus_3$};
  \draw (0,3) node (m2) {$\boxminus_2$};
  \draw (0,4) node (m1) {$\boxminus_1$};

  \draw [dotted] (1,-0.5) to (1,4.5);
  
  \draw (2,0) node (p5) {$\boxplus_5$};
  \draw (2,1) node (p4) {$\boxplus_4$};
  \draw (2,2) node (p3) {$\boxplus_3$};
  \draw (2,3) node (p2) {$\boxplus_2$};
  \draw (2,4) node (p1) {$\boxplus_1$};

  \draw [boxedge] (m4) to [bend right=30] (m3);
  \draw [boxedge] (m3) to (p1);
  \draw [boxedge] (p1) to [bend left=30] (p4);

  \draw [boxedge,dashed] (m5) to [bend left=30] (m1);
  \draw [boxedge,dashed] (m1) to (p3);
  \draw [boxedge,dashed] (p3) to [bend right=30] (p5);

\end{scope}
  
\begin{scope}[xshift=4cm,scale=.8]
  \draw [label distance=5mm] (0,0) node[vertex,label=left:$5$] (5) {};
  \draw [label distance=5mm] (0,1) node[vertex,label=left:$4$] (4) {};
  \draw [label distance=5mm] (0,2) node[vertex,label=left:$3$] (3) {};
  \draw [label distance=5mm] (0,3) node[vertex,label=left:$2$] (2) {};
  \draw [label distance=5mm] (0,4) node[vertex,label=left:$1$] (1) {};

  \draw (0.5,-.5) node {$V$};
  
  \draw [dotted] (1.5,0) to (1.5,4);

  \draw (2.5,-.5) node {$U$};
  
  \draw (2,2) node[vertex,label=right:$j$] (j) {};
  
  \draw [color=xblue] (4) to [bend left=20] (3);
  \draw [->,color=xgreen] (3) to [bend left=20] (1);
  \draw [color=xred] (4) to [bend left=40] (1);

  \draw [color=xgray] (4) to (j);
  \draw [color=xblack] (3) [bend left=15] to (j);
  \draw [color=xpurple] (1) [bend left=15] to (j);

  \draw [color=xblue,dashed] (5) to [bend left=45] (1);
  \draw [->,color=xgreen,dashed] (1) to [bend left=20] (3);
  \draw [color=xred,dashed] (3) to [bend left=20] (5);

  \draw [color=xgray,dashed] (5) to  (j);
  \draw [color=xblack,dashed] (1) to (j);
  \draw [color=xpurple,dashed] (3) to (j);

\end{scope}

\begin{scope}[xshift=8cm, scale=.8]
  \draw (0,0) node (m5) {$\boxminus_5$};
  \draw (0,1) node (m4) {$\boxminus_4$};
  \draw (0,2) node (m3) {$\boxminus_3$};
  \draw (0,3) node (m2) {$\boxminus_2$};
  \draw (0,4) node (m1) {$\boxminus_1$};

  \draw [dotted] (1,-0.5) to (1,4.5);
  
  \draw (2,0) node (p5) {$\boxplus_5$};
  \draw (2,1) node (p4) {$\boxplus_4$};
  \draw (2,2) node (p3) {$\boxplus_3$};
  \draw (2,3) node (p2) {$\boxplus_2$};
  \draw (2,4) node (p1) {$\boxplus_1$};

  \draw [boxedge] (m5) to [bend right=30] (m4);
  \draw [boxedge] (m4) to (p1);
  \draw [boxedge] (p1) to [bend left=30] (p5);

  \draw [boxedge,dashed] (m5) to [bend left=30] (m1);
  \draw [boxedge,dashed] (m1) to (p3);
  \draw [boxedge,dashed] (p3) to [bend right=30] (p5);
\end{scope}

\begin{scope}[xshift=12cm,scale=.8]
  \draw [label distance=5mm] (0,0) node[vertex,label=left:$5$] (5) {};
  \draw [label distance=5mm] (0,1) node[vertex,label=left:$4$] (4) {};
  \draw [label distance=5mm] (0,2) node[vertex,label=left:$3$] (3) {};
  \draw [label distance=5mm] (0,3) node[vertex,label=left:$2$] (2) {};
  \draw [label distance=5mm] (0,4) node[vertex,label=left:$1$] (1) {};

  \draw (0.5,-.5) node {$V$};
  
  \draw [dotted] (1.5,0) to (1.5,4);

  \draw (2.5,-.5) node {$U$};
  
  \draw (2,2.5) node[vertex,label=right:$j$] (j) {};
  \draw (2,1.5) node[vertex,label=right:$j'$] (j') {};
  
  \draw [color=xblue] (5) to [bend left=20] (4);
  \draw [->,color=xgreen] (4) to [bend left=20] (1);
  \draw [color=xred] (1) to [bend right=35] (5);

  \draw [color=xgray] (5) to (j);
  \draw [color=xblack] (4) to (j);
  \draw [color=xpurple] (1) to (j);

  \draw [color=xblue,dashed] (5) to [bend left=45] (1);
  \draw [->,color=xgreen,dashed] (1) to [bend left=20] (3);
  \draw [color=xred,dashed] (3) to [bend left=20] (5);

  \draw [color=xgray,dashed] (5) to  (j');
  \draw [color=xblack,dashed] (1) to (j');
  \draw [color=xpurple,dashed] (3) to (j');
\end{scope}

\end{tikzpicture}

\begin{tikzpicture}[auto,thick]

\begin{scope}
\draw (0,-1) node {Corresponding edge types in the chest:
  \textcolor{xgreen}{$\vec{E}_1$},
  \textcolor{xblue}{$E_2$},
  \textcolor{xred}{$E_3$},
  \textcolor{xblack}{$E_4$},
  \textcolor{xpurple}{$E_5$},
  \textcolor{xgray}{$E_6$}.
};
\end{scope}

\end{tikzpicture}

\caption{The two left-hand pictures show an example of two paths from the same
  graph~$G_j$ embedded in~$H^*$ to the paths
  $\boxminus_4 \boxminus_3 \boxplus_1 \boxplus_4$ and
  $\boxminus_5 \boxminus_1 \boxplus_3 \boxplus_5$ and the corresponding diamonds
  in the chest.  The two right-hand pictures show an example of two paths from
  different graphs~$G_j$ and~$G_{j'}$ embedded in~$H^*$ to the paths
  $\boxminus_5 \boxminus_4 \boxplus_1 \boxplus_5$ and
  $\boxminus_5 \boxminus_1 \boxplus_3 \boxplus_5$, respectively, and the
  corresponding diamonds in the chest. The solid and dashed lines in these pictures are only used to distinguish the two paths/diamonds.}
\label{fig:outline}
\end{figure}

We should at this point note that the chest can and will have multiple edges
between its vertices, but we will ensure that $\boxminus_i\boxplus_i$ is not an
edge of $H^*$, meaning that the chest has no loops. Suppose now that we have a
collection of edge-disjoint diamonds in a chest, that is, we do not use any one
coloured edge or arc in two different diamonds. We remark that it may still be
the case that one red edge $ij$ is used in one diamond while a parallel blue
edge $ij$ is used in another. This condition of edge-disjointness of the
diamonds translated back to~$H^*$ means the following: the embeddings of
$3$-edge paths from various graphs $G_s$ with $s\in\cJ_0$ encoded by the
diamonds do not use any one edge of $H^*$ twice, and any two paths from any one
$G_s$ get embedded to disjoint sets of vertices. In other words, a collection of
edge-disjoint diamonds represents a way to extend the packing we have before the
final stage to a larger packing. In particular, we will ensure that the numbers of
red, blue, green and grey edges are all the same. A collection of edge-disjoint
diamonds which uses all of them---which we call a \emph{diamond
core-decomposition}---then represents an extension of the packing before the
final stage to the desired perfect packing. Note that our diamond
core-decompositions will not use all the black or purple edges, since the graphs
$G_s$ with $s\in\cJ$ are not spanning. A diamond core-decomposition is precisely
a generalised design, whose existence is proved by
Keevash~\cite{Kee2,KeeColouredDirected} under certain conditions on the chest
(see Section~\ref{sec:designs} for more details). This final stage of our
packing will be Stage~G, \emph{designs completion}.

The preparation we need to do before this is then simply to ensure that we end up with the above setting, and that the Keevash conditions for the existence of the required generalised design are met. In Stage~C, \emph{partite reduction}, we split the vertices of $H$ into two equal sides and pair them up into terminal pairs randomly. If $v(H)$ is odd, we have one leftover vertex $\boxdot$. We embed a few paths from the graphs $G_s$ with $s\in\cJ$, including all of those anchored at $\boxdot$, in order to use up all the edges of $H$ leaving $\boxdot$ and any edges of $H$ between terminal pairs.

In Stage~D, \emph{connecting to terminal pairs}, we embed some vertices from the graphs $G_s$ with $s\in\cJ_0$. Prior to this stage, the bare paths in these graphs all have $11$ edges, and they can be anchored anywhere in $H$ (except $\boxdot$), not necessarily to terminal pairs. In this stage, we embed the first four and last four vertices of each path so that what remains is to embed the middle $3$ edges between some terminal pair $\{\boxminus_i,\boxplus_i\}$ which we choose randomly.

As mentioned above, for our designs completion Stage~G, we need the same numbers of red, blue and green edges in the chest. That means we need the same number of edges of $H$ within $\Vmin$, within $\Vplus$ and crossing between $\Vmin$ and $\Vplus$. After Stage~D, the leftover graph $H$ has a fairly uniform density, and in particular there are about twice as many edges crossing from $\Vmin$ to $\Vplus$ as within $\Vmin$. In Stage~E, \emph{density correction}, we now make sure that we obtain precisely the correct number of edges within $\Vmin$ and $\Vplus$ and crossing. We will embed further paths in the following Stage~F, but the number of edges from each set we use is fixed and taken into account. We do this density correction by completing the embedding of the graphs $G_s$ with $s\in\cJ_2$.

In Stage~F, \emph{degree correction}, we ensure that some conditions of the following form hold. For any given $\boxminus_i\in\Vmin$, the number of edges of $H$ going from $\boxminus_i$ to $\Vmin$ is equal to the number of edges going from $\boxminus_i$ to $\Vplus$ plus the number of $s\in\cJ_0$ for which $G_s$ has a path anchored at $\boxminus_i$. Observe that if we embed a diamond using $i$, there is a blue edge at $i$ and either a green arc leaving $i$ or a grey edge at $i$, so that a diamond core-decomposition can only exist if the above equality holds. Following the language of Keevash~\cite{Kee2,KeeColouredDirected} we call these \emph{divisibility conditions}. All these conditions are very close to holding already after Stage~E, so only tiny corrections are necessary. We perform the degree correction by completing the embedding of the graphs $G_s$ with $s\in\cJ_1$.

We have $\cJ=\cJ_0\cup\cJ_1\cup\cJ_2$, so after Stage~F, what remains is precisely the embedding of $3$-edge paths from graphs $G_s$ with $s\in\cJ_0$ described above as Stage~G, \emph{designs completion}. To apply Keevash's designs result in this stage we need one further condition:
the chest has to have certain quasirandomness properties (which
we define precisely in Section~\ref{sec:designs}). In Stages~A,~D and~E (which is where we embed a
significant number of edges) we use randomised algorithms to perform our
embedding, and with high probability these algorithms give the required
quasirandom output. In Stages~B,~C and~F the number of edges we embed is tiny
compared to $e(H)$. We simply ensure that we do not embed to any one vertex too
often in these stages, and this is enough to ensure that quasirandomness cannot
be seriously affected. We should note that while we eventually need
quasirandomness to apply the Keevash machinery, we will also use it often in our
analysis of the preparatory stages; it is for instance the quasirandomness that
guarantees our divisibility conditions are all close to correct after
Stage~E.

\smallskip

Let us close this proof outline by describing the bulk embedding, Stage~A, in
some more detail.  This part of our proof is an extension of the proof of
Theorem~\ref{thm:ABHPDegenerateNonTechnical} given in~\cite{DegPack}, which is
built on the following idea. Enumerate the graphs as
$G_1,\ldots,G_{|\mathcal{G}|}$, with the graphs $G_s$ for $s\in\cJ\cup\cK$
last. We embed the graphs one after another, edge-disjointly, in this order,
except for the vertices in $\SpecLeaves_s$ and $\BasicPaths_s$.

For each $G_s$, we pack the first $(1-\deltnonspanning)n$ vertices in the $D_0$-degeneracy order as follows. We embed the first vertex of $G_s$ uniformly to
$V(H)$. Thereafter, when we need to embed vertex $i$, we look at the already
embedded neighbours $x_1,\dots,x_r$ of $i$. Suppose these are embedded to
$v_1,\dots,v_r$ in $H$. Then we need to embed $i$ to a vertex which is adjacent
to all of $v_1,\dots,v_r$ in $H$ and which we did not previously use in
embedding $G_s$. We pick such a vertex~$v$ uniformly at random, embed $i$ there,
and delete all edges $vv_j$ with $j\in[r]$ from~$H$. 

For the graphs $G_s$ with $s\in\cG\setminus(\cJ\cup\cK)$, we still need to embed the remaining $\deltnonspanning n$ vertices. This requires a separate argument which is not particularly relevant to this discussion.

To show that this random process succeeds, we need to show that after each $G_s$ is embedded, the remaining graph $H$ is still quasirandom. In turn, during the
embedding of $G_s$, we need to argue that the first $i$ vertices are embedded to
roughly an $i/n$ fraction of the common neighbourhood of any $D_0$ or fewer
vertices of $H$. This is called the~\emph{diet condition}, and together with the
quasirandomness of $H$ in particular it tells us that there will always be a
large set of vertices to which we can embed vertex $i+1$ of $G_s$.

The above analysis is detailed in~\cite{DegPack}. However, to make Stages~B--G
work, we need quite a few additional properties of this packing, which we prove
in this paper (in Section~\ref{sec:StageA}). In addition, the randomised
algorithm which we use in Stages~D and~E is related to the procedure described
above, but requires an entirely new analysis (which is given in
Section~\ref{sec:pathpack}).

We remark that another approach extending~\cite{DegPack} was used
in~\cite{ABCT:PackingManyLeaves} to prove Theorem~\ref{thm:ABCT}. However, the
structure that is put aside and packed later there consists of the leaves in
$\left(G_s\right)_{s\in\mathcal{B}}$. The subsequent perfect packing of these
leaves is much easier than that of our systems of paths.

\section{Designs}
\label{sec:designs}

The main purpose of this section is to formulate a general decomposition result
for directed coloured partite multigraphs
(Theorem~\ref{thm:DesignsIntermediate}), which is a special case of the deep
results on the existence of designs by Keevash~\cite{Kee2, KeeColouredDirected},
and to apply it to obtain a specific decomposition result of certain partially
directed coloured partite multigraphs
(Proposition~\ref{prop:Designs-oursetting}) that we shall use to complete the
packing in our proof of Theorem~\ref{thm:maintechnical}.
Here, a \emph{partially directed multigraph}\index{partially
  directed multigraph} consists of a vertex set, a set (or a collection of sets)
of directed edges and a set (or a collection of sets) of undirected edges, with
multi-edges allowed but loops not permitted.

We first introduce some basic notions. For a directed graph~$G$ and $v\in V(G)$
we define $\NBHout_G(v)$\index{$\NBHout_G(v)$} and
$\NBHin_G(v)$\index{$\NBHin_G(v)$} as the out-neighbourhood and the
in-neighbourhood of $v$, respectively. As before, $\NBHout_{G}(S)=\bigcap_{v\in
  S}\NBHout_{G}(v)$ and $\NBHin_{G}(S)=\bigcap_{v\in
  S}\NBHin_{G}(v)$\index{$\NBHout_{G}(S)$}\index{$\NBHin_{G}(S)$}. We write
$\degout_G(v)$\index{$\degout_G(v)$} and $\degin_G(v)$\index{$\degin_G(v)$}
for the corresponding degrees, that is, for the number of edges of~$G$ directed
away from~$v$ and the number of edges of~$G$ directed towards~$v$,
respectively. We shall be considering partially directed multigraphs~$\mathcal M$
with a collection of different edge sets $E_1,E_2,\dots,E_k$, where $E_i$ either
consists only of directed edges or only of undirected edges. We also write
$\NBHout_{E_i}(v)$, $\NBHin_{E_i}(v)$, $\degout_{E_i}(v)$, $\degin_{E_i}(v)$ in
the former case, and $\NBH_{E_i}(v)$, $\deg_{E_i}(v)$ in the latter case for the
neighbourhoods and degrees in the sub(di)graph of~$\mathcal M$ with edge
set~$E_i$.

Let us emphasise that degrees are always counted \emph{with multiplicity},
whereas neighbourhoods are sets; thus $\degout_G(v)\ge|\NBHout_G(v)|$ in
general, with equality if and only if no two parallel edges leave~$v$. This is
worth recording because the multi-digraphs of
Definition~\ref{def:general-digraph} below do admit parallel edges of the same
colour, while the divisibility conditions of Definition~\ref{def:divisibility}
are phrased in terms of degrees. In each (multi-)digraph to which we apply these
notions, however, degrees are taken with respect to a single edge set~$E_i$,
respectively a single colour class, and none of these contains parallel edges;
so for us the two readings in fact agree.

\begin{figure}
	\includegraphics[scale=0.8]{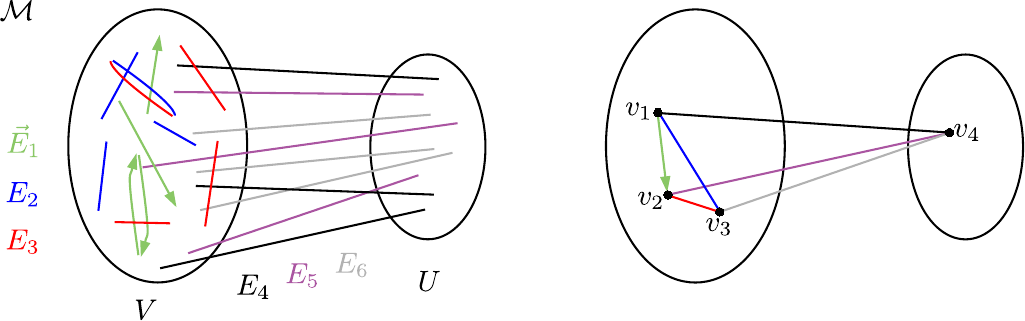}
	\caption{A chest (on the left) and a diamond in it (on the right).}
	\label{fig:ChestsAndDiamonds}
\end{figure}

We now define the setup in which we want to apply the decomposition results
mentioned above. The partially directed multigraph which we want to decompose,
and which we shall call \emph{chest}\index{chest}, has the following form.
See the left-hand side of Figure~\ref{fig:ChestsAndDiamonds} for an illustration
of this setup.

\begin{definition}[chest]\label{def:chest}
  We have disjoint sets~$V$ and~$U$ of vertices and the following collection of
  edge sets.  The edge set $\vec{E}_1$ contains directed edges in $V$, the edge
  sets $E_2$ and $E_3$ contain edges in $V$, and $E_4, E_5, E_6$ contain edges
  between $V$ and $U$. Parallel edges or loops are not allowed within any of
  these edge sets, though antiparallel edges are allowed in the set of
  directed edges $\vec{E}_1$, and different edge sets may have edges in common.
  If $V$, $U$, $\vec{E}_1$, $E_2,\dots,E_6$ satisfy these properties, we say
  that $\mathcal{M}= (V\dcup U; \vec{E}_1, E_2, E_3, E_4, E_5, E_6)$ is a
  \emph{chest}\index{chest}.
\end{definition}

In order to be able to apply the decomposition results, we need
our chest to have certain quasirandomness properties.

\begin{definition}[quasirandom chest]\label{def:quasichest}
  Define
  $d_1:=\frac{|\vec{E}_1|}{|V|^2}$, $d_i:=\frac{|E_i|}{\binom{|V|}{2}}$ for
  $i=2,3$, and $d_i:=\frac{|E_i|}{|V||U|}$ for $i=4,5,6$. We say that a chest
  $\mathcal{M}$ is \emph{$(\gamma,L)$-quasirandom}\index{quasirandom chest} if for
  every choice of $S_1,S'_1,S_2, S_3,S'_4, S'_5,S'_6\subseteq V$ and
  $S_4,S_5,S_6\subseteq U$ of mutually disjoint sets of total size at most $L$ we
  have that
  \begin{align*}
    \left|V\cap\NBHout_{\vec{E}_1}(S_1)\cap \NBHin_{\vec{E}_1}(S'_1)\cap \bigcap_{i=2}^6 \NBH_{E_i}(S_i)\right|&=(1\pm \gamma)\cdot d_1^{|S_1|+|S'_1|}\cdot \prod_{i=2}^6 d_i^{|S_i|}\cdot|V|\quad\mbox{and}\\
    \left|U\cap \bigcap_{i=4}^6 \NBH_{E_i}(S'_i)\right|&=(1\pm \gamma)\cdot \prod_{i=4}^6 d_i^{|S'_i|}\cdot|U|\;.
  \end{align*}
\end{definition}
  
We next define the partially directed graph into which we want to decompose a quasirandom chest.

\begin{definition}[diamond]\label{def:diamond}
  A \emph{diamond}\index{diamond} in $\mathcal{M}$ is a partially directed graph
  $(\{v_i\}_{i=1}^4,\{e_i\}_{i=1}^6)$ with $e_1\in \vec{E}_1$ and $e_i\in E_i$ for
  $i=2,\ldots,6$, such that $e_1=(v_1,v_2)$, $e_2=\{v_3,v_1\}$,
  $e_3=\{v_2,v_3\}$, $e_4=\{v_1,v_4\}$, $e_5=\{v_2,v_4\}$, and
  $e_6=\{v_3,v_4\}$.
\end{definition}

See the right-hand side of
Figure~\ref{fig:ChestsAndDiamonds} for an illustration of a diamond.
A collection $\cD$ of diamonds in a chest $\mathcal{M}$ is a
\emph{diamond core-decomposition}\index{diamond core-decomposition} of
$\mathcal{M}$ if each edge of $\vec{E}_1\dcup E_2\dcup E_3\dcup E_6$ is used by
exactly one diamond of $\cD$, and each edge of $E_4\dcup E_5$ is used by at most
one diamond of $\cD$.

It turns out that, in order to complete the packing of our trees in the proof of
Theorem~\ref{thm:maintechnical}, we need a diamond core-decomposition of a
quasirandom chest (this part of our proof is encapsulated in
Lemma~\ref{lem:StageG}). The following proposition provides conditions under
which we can obtain such a diamond core-decomposition.

\begin{proposition}[diamond core-decomposition]\label{prop:Designs-oursetting}
For every $d, \sigma >0$ there exist $L, n_0\in\NN$ and $\gamma>0$ with the following property. Let $\mathcal{M}= (V\dcup U; \vec{E}_1, E_2, E_3, E_4, E_5 , E_6)$ be a $(\gamma,L)$-quasirandom chest (see Definitions~\ref{def:chest} and~\ref{def:quasichest}) with $|U\cup V|=n>n_0$ and $\min\{|U|,|V|\}\ge \sigma n$, such that $d_1=\frac{|\vec{E}_1|}{|V|^2}>d$,  and $d_i=\frac{|E_i|}{|V||U|}>d$ for $i=4,5,6$. Suppose that $|\vec{E}_1|=|E_2|=|E_3|=|E_6|$. Suppose further that for any vertex $v\in V$ we have
	\begin{enumerate}[label=\rom]
		\item\label{itm:designs:i}  $\deg_{E_2}(v)=\degout_{\vec{E}_1}(v)+\deg_{E_6}(v)$,
		\item\label{itm:designs:ii} $\deg_{E_3}(v)=\degin_{\vec{E}_1}(v)+\deg_{E_6}(v)$,
		\item\label{itm:designs:iii} $\deg_{E_4}(v)\ge \degout_{\vec{E}_1}(v)+d^5\sigma^2|U|$,
		\item\label{itm:designs:iv} $\deg_{E_5}(v)\ge \degin_{\vec{E}_1}(v)+d^5\sigma^2|U|$,
	\end{enumerate}
	and for any vertex $u\in U$ we have
	\begin{enumerate}[label=\rom,resume]
		\item\label{itm:designs:v} $\deg_{E_4}(u)\ge \deg_{E_6}(u)+d^5\sigma^2|V|$, and
		\item\label{itm:designs:vi} $\deg_{E_5}(u)\ge \deg_{E_6}(u)+d^5\sigma^2|V|$.
	\end{enumerate}
Then $\mathcal{M}$ has a diamond core-decomposition.
\end{proposition}

The remainder of this subsection is dedicated to the deduction of
Proposition~\ref{prop:Designs-oursetting}
from~\cite[Theorem~19]{KeeColouredDirected}. This theorem is very general and
relies on heavy terminology specific to~\cite{KeeColouredDirected}. So, as an
intermediate step, we present in Theorem~\ref{thm:DesignsIntermediate} a
decomposition result which is fairly general yet relatively easy to state. In
order to state this result we need a number of definitions.

Theorem~\ref{thm:DesignsIntermediate} allows us to decompose multi-digraphs into
a family of digraphs in a coloured and partite setting.  For an integer~$D$, a
\emph{$[D]$-edge-coloured digraph} is a digraph whose edges are assigned colours
from the set $[D]$. We remark that in this theorem the roles of~$H$ and~$G$
are interchanged compared to the setting in the remainder of the paper: we
decompose a digraph~$G$ into graphs~$H$ from a family~$\mathcal H$. The reason
for this change is that this will make it easier for us to explain how
Theorem~\ref{thm:DesignsIntermediate} follows
from~\cite[Theorem~19]{KeeColouredDirected}.

\begin{definition}[decomposition]\label{def:decomposition}
 Let $\mathcal H$ be a family of $[D]$-edge-coloured digraphs on $[q]$, and let
 $G$ be a $[D]$-edge-coloured digraph on $[n]$. We say that $G$ has an
 \index{decomposition}\emph{$\mathcal H$-decomposition} if the edges of $G$ can
 be partitioned into copies of digraphs from $\mathcal H$ that preserve the
 colouring.
\end{definition}

We can only decompose into certain types of edge-coloured digraphs,
which we call simple canonical digraphs. 

\begin{definition}[simple canonical digraphs]\label{def:canonical}
  Let $D,q\in \NN$, and let $\mathcal P=(P_1, \ldots, P_t)$ be a partition of
  $[q]$.  Let $\mathcal H$ be a family of $[D]$-edge-coloured digraphs
  on~$[q]$. Further, for each colour $d\in[D]$ assume we are given a pair
  $(i,j)\in [t]^2$, which we call \emph{colour location of $d$}. For a colour
  $d\in[D]$ with colour location $(i,j)$ such that $i=j$ we assume further
  that~$d$ is specified as being an \emph{oriented colour}\index{oriented
    colour} or an \emph{unoriented colour}\index{unoriented colour}.
  
  In this case, we say that the family~$\mathcal H$ is \index{simple canonical
    family of digraphs}\emph{simple $\mathcal P$-canonical} if the following
  hold for every $H\in \mathcal H$ and every colour $d\in [D]$ and its colour
  location $(i,j)$.
  \begin{itemize}[leftmargin=*]
  \item There is no loop in~$H$, there are no parallel edges in~$H$, and there
    are no antiparallel edges of different colours in~$H$.
  \item All the edges of $H$ of colour $d$ start in $P_i$ and end in
    $P_j$.
  \item If $i=j$ and $d$ is an unoriented colour, then all the edges of $H$ of
    colour $d$ come in pairs that form directed $2$-cycles.
  \item If $i=j$ and $d$ is an
    oriented colour, then $H$ contains no directed $2$-cycles in colour $d$.
  \end{itemize}
\end{definition}

As an illustration, let us explain how this definition is used in our application.
Recall that in our setup we decompose into diamonds on vertex set $[4]$
(replacing~$v_i$ in Definition~\ref{def:diamond} with~$i$) containing directed
edges and undirected edges and with vertices $1,2,3$ mapped to~$V$ and
vertex~$4$ mapped to~$U$, as well as leftover edges from~$E_4$
and~$E_5$. Translating this to the setting of simple canonical digraphs, we
would choose $\mathcal{P}=(P_1,P_2)$ with $P_1=\{1,2,3\}$ and $P_2=\{4\}$ and
$\mathcal H=\{H_1,H_2,H_3\}$ containing the following digraphs.  The
digraph~$H_1$ is a directed version of the diamond: we give the~$6$ edges of the
diamond~$6$ different colours, and then replace each undirected edge~$\{i,4\}$ of
colour~$d$ with an edge directed towards~$4$ of colour~$d$ (this choice of
direction is arbitrary), and each other undirected edge (within~$P_1$) of
colour~$d$ with two antiparallel edges of colour~$d$.  The digraphs~$H_2$ and~$H_3$ can be
used for the leftover edges from~$E_4$ and~$E_5$: $H_2$ only contains the edge
$(1,4)$ in the same colour as in~$H_1$ and $H_3$ only contains the edge $(2,4)$
in the same colour as in~$H_1$.

We next define the types of digraphs that can be decomposed with the help of
Theorem~\ref{thm:DesignsIntermediate}.

\begin{definition}[general multi-digraph]\label{def:general-digraph}
  Let $\mathcal P'=(P'_1, \ldots, P'_t)$ be a partition of $[n]$. We say that a
  $[D]$-edge-coloured multi-digraph $G$ on~$[n]$ with partition $\mathcal P'$ is
  \index{general digraph}\emph{general} if $G$ has no loop, but multi-edges,
  parallel and antiparallel, of the same or different colours are allowed, as long as
  for any colour $d\in[D]$ between any two parts $P'_i$ and $P'_j$ with $i\neq
  j$, either all edges of colour~$d$ are directed towards~$P'_i$ or all edges of
  colour~$d$ are directed towards~$P'_j$.  
\end{definition}

Further, the digraphs we want to decompose need to satisfy a number of
conditions, which we define next, and which allow for a partite setting.  We
start with the divisibility conditions. For a $[D]$-edge-coloured
(multi-)digraph~$H$, a vertex~$v$ of~$H$ and a colour $d\in[D]$, we write~$H_d$
for the sub-(multi-)digraph of~$H$ containing exactly those edges of colour~$d$.
We also write $\degout_{H,d}(v)$ for $\degout_{H_d}(v)$ and $\degin_{H,d}(v)$
for $\degin_{H_d}(v)$, and analogously for $\NBHout_{H,d}(v)$,
$\NBHin_{H,d}(v)$, $\NBHout_{H,d}(S)$, and $\NBHin_{H,d}(S)$, where $S\subset
V(H)$.

\begin{definition}[divisibility]\label{def:divisibility}
  Let $\mathcal P=(P_1, \ldots, P_t)$ be a partition of $[q]$ and $\mathcal
  P'=(P'_1, \ldots, P'_t)$ be a partition of $[n]$. Let $G$ be a
  $[D]$-edge-coloured general multi-digraph on~$[n]$. Let $\mathcal H$ be a
  family of $[D]$-edge-coloured simple $\mathcal P$-canonical digraphs on
  $[q]$. We say that $(G,\mathcal P')$ is \index{divisible}\emph{$(\mathcal H,\mathcal
  P)$-divisible}  if the following hold.
  \begin{description}[font=\textnormal]
  \item[$0$-divisibility]\index{0-divisible} For $d\in [D]$, and $H\in \mathcal
    H$, let $c_{d,H}$ denote the number of edges in $H$ coloured~$d$. There
    are integers $(m_{H})_{H\in \mathcal H}$ such that for each $d\in [D]$ the
    number of edges of $G$ in colour $d$ equals $\sum_{H\in \mathcal H}m_H\cdot
    c_{d,H}$.
  \item[$1$-divisibility]\index{1-divisible} For any $i\in [t]$ and any vertex
    $v\in P'_i$, there exist integers $(m_{H,x})_{H\in \mathcal H, x\in P_i}$
    such that for each $d\in [D]$ we have
    \[
    \quad
    \degout_{G,d}(v)=\sum_{H\in \mathcal H, x\in P_i}m_{H, x}\cdot \degout_{H,d}(x)
    \quad\text{and}\quad
    \degin_{G,d}(v)=\sum_{H\in \mathcal H, x\in P_i}m_{H, x}\cdot \degin_{H,d}(x) \;.
    \] 
  \item[$2$-divisibility]\index{2-divisible} For each $d\in [D]$ with colour
    location $(i,j)$ (with respect to $\mathcal{H}$), all edges in $G$ of colour
    $d$ start in $P'_i$ and end in $P'_j$. Further, if $i=j$ and colour $d$ is
    unoriented then all the edges of $G$ of colour $d$ come in pairs that form
    directed $2$-cycles.
  \end{description}
\end{definition}

Observe that $2$-divisibility mandates that in~$G$ edges of colour~$d$ only run
between a unique pair $(P'_i,P'_j)$ of parts (with possibly $i=j$).
Since~$\mathcal H$ is canonical also in any $H\in\mathcal H$ we can have edges
of this colour~$d$ only between the corresponding pair $(P_i,P_j)$.

The next condition requires that copies of the digraphs~$H$ into which we seek
to decompose the digraph~$G$ are distributed regularly on edges of~$G$.

\begin{definition}[regularity]\label{def:regularity(designs)}
  Let $\mathcal H$ be a family of $[D]$-edge-coloured digraphs on $[q]$ and let
  $G$ be a general $[D]$-edge-coloured multi-digraph on $[n]$. Let $c,\omega>0$ be reals.  We
  say that $G$ is \index{regular}\emph{$(\mathcal H,c,\omega)$-regular} if we
  can assign a weight $w_{H'}\in [\omega\cdot n^{2-q},\omega^{-1}\cdot n^{2-q}]$ to each
  coloured copy $H'\subseteq G$ of any $H\in \mathcal H$, such that for every
  edge $e\in E(G)$ we have
  \begin{equation}\label{eq:reg1}
   \sum_{H\in \mathcal H}\quad \sum_{H'\subseteq G: e\in E(H'), H'\sim H}w_{H'}=(1\pm c)\,,  \end{equation}
  where we denote the fact that~$H'$ is a coloured copy of~$H$ by $H'\sim H$.
\end{definition}

In this definition, a \emph{coloured copy} of~$H$ in~$G$ is an injection from
the vertex set $[q]$ of~$H$ to the vertex set $[n]$ of~$G$ which maps every edge
of~$H$ to an edge of~$G$ of the same colour and orientation; in particular,
vertices which are isolated in~$H$ are also embedded. Further, whenever
partitions $\mathcal P$ of $[q]$ and $\mathcal P'$ of $[n]$ are given (as in
Theorem~\ref{thm:DesignsIntermediate} below), we additionally require that
copies map each part $P_i$ into $P'_i$, and we then say that $(G,\mathcal P')$
is $(\mathcal H,c,\omega)$-regular. The normalisation of the weights is
$n^{2-q}$ because we fix an edge in~$G$ and sum over copies of~$H$ in~$G$ which
contain this fixed edge.

The following final condition considers any vertex~$x$ in any $H\in\mathcal H$
and requires that whenever we choose for each vertex~$y$ in~$H$ other than~$x$ a
set~$A_y$ of at most~$h$ vertices in the part of~$G$ where~$y$ should be
embedded then the common neighbourhood of these sets~$A_y$ in the part where~$x$
should be embedded is of linear size, where for this common neighbourhood we
take edges directed as mandated by corresponding edges in~$H$.

\begin{definition}[vertex-extendability]\label{def:vertex-extensibility}
  Let $\mathcal P=(P_1, \ldots, P_t)$ be a partition of $[q]$ and let $\mathcal
  P'=(P'_1, \ldots, P'_t)$ be a partition of $[n]$.  Let $\mathcal H$ be a
  simple $\mathcal P$-canonical family of $[D]$-edge-coloured digraphs on
  $[q]$. Let $G$ be a general $[D]$-edge-coloured multi-digraph on $[n]$.

  Let $H\in\mathcal H$,
  let $x\in[q]$ be any vertex of $H$, and let~$P_i$ be the part of~$\mathcal P$
  containing~$x$.
  Let~$y$ be any other vertex of~$H$, that is, $y\in[q]\setminus\{x\}$.
  If $(x,y)$ or $(y,x)$ is an edge of~$H$, then let~$d$ be its
  colour. (Recall that since~$\mathcal H$ is canonical, if $(x,y)$ and $(y,x)$
  are both edges, then they have the same unoriented colour.)
  For any set $A\subset[n]$ of vertices in~$G$, we define
  \[
  \NBH^{(y,x,H)}_G(A)=\begin{cases}
   \NBHout_{G,d}(A)\cap P'_i & \text{if only } (y,x)\in E(H)\,, \\
   \NBHin_{G,d}(A)\cap P'_i & \text{if only } (x,y)\in E(H)\,, \\
   \NBHout_{G,d}(A)\cap\NBHin_{G,d}(A)\cap P'_i \quad & \text{if both } (y,x),(x,y)\in E(H)\,, \\
  P'_i & \text{otherwise\,.}
  \end{cases}
  \]

  Let $h$ be an integer and $\omega>0$.
  We  say that $(G,\mathcal P')$ is \index{vertex extendable}\emph{$(\mathcal
  H,\mathcal P,\omega,h)$-vertex-extendable}
  if the following holds for every $H\in \mathcal H$ and every $x\in [q]$.
  For every choice of pairwise disjoint sets $\{A_y\}_{y\in[q]\setminus \{x\}}$ of size $|A_y|\le h$
  with $A_y\subseteq P'_j$ whenever $y\in P_j$, the set $A_x:= \bigcap_{y\in [q]\setminus
    \{x\}}\NBH^{(y,x,H)}_G(A_y)$ has size $|A_x|\ge \omega n$.
\end{definition}

In our setting, regularity and vertex-extendability both follow from the quasirandomness properties of the chest.
We are now ready to state the decomposition result that implies
Proposition~\ref{prop:Designs-oursetting}.

\begin{theorem}[decomposition result]\label{thm:DesignsIntermediate}
  Given $q,D\in \NN$ and $\sigma>0$, there exist numbers $\omega_0>0$ and
  $n_0\in \NN$ such that with $q'=\max\{q,8+\lfloor\log_2(1/\sigma)\rfloor\}$,
  $h=2^{50 q'^3}$ and $\delta=2^{-10^3q'^5}$ the following holds for each
  $n>n_0$ and $\omega\in (n^{-\delta},\omega_0)$.

  Let $\mathcal P=(P_1,\ldots, P_t)$ be a partition of $[q]$ and $\mathcal
  P'=(P'_1,\ldots,P'_t)$ be a partition of $[n]$ such that $|P'_i|\ge \sigma n$
  for each $i\in[t]$. Let $\mathcal H$ be a simple $\mathcal P$-canonical family
  of $[D]$-edge-coloured digraphs on $[q]$, and let $G$ be a general
  $[D]$-edge-coloured multi-digraph on $[n]$.

  Suppose that $(G,\mathcal P')$ is $(\mathcal H,\mathcal P)$-divisible, $(\mathcal
  H,\omega^{h^{20}},2h^{q'-q}\omega)$-regular and $( \mathcal
  H,\mathcal P,\sqrt[q'h]{\omega},h)$-vertex-exten\-dable. Then $G$ has an $\mathcal
  H$-decomposition.
\end{theorem}

Theorem~\ref{thm:DesignsIntermediate} is a special case
of~\cite[Theorem~19]{KeeColouredDirected}. In Appendix~\ref{appendix:designs},
we explain the connections between the notion used in~\cite{KeeColouredDirected}
and the one used here. We now use Theorem~\ref{thm:DesignsIntermediate} to prove
Proposition~\ref{prop:Designs-oursetting}.

\begin{proof}[Proof of Proposition~\ref{prop:Designs-oursetting}]
  We first choose the constants.
  Let $d,\sigma>0$ be given.  Set $q:= 4$ and $D:=6$.  Let
  $\omega_0$ and $n'_0$ be the output parameters of
  Theorem~\ref{thm:DesignsIntermediate} for input $q$, $D$ and $\sigma$ and set
  \[q':=\max\{q,8+\lfloor\log_2(1/\sigma)\rfloor\}\,,\quad
  h:=2^{50 q'^3}\quad \text{and}\quad
  \delta:=2^{-10^3q'^5}\,,
  \]
  as specified in Theorem~\ref{thm:DesignsIntermediate}.
  Further, set
  \[L:=3h\,,\quad
  \omega:=\min\Big\{\Big(\frac{\sigma d}{2}\Big)^{5Lhq'},\frac{\omega_0}{2}\Big\}\,,\quad
  n_0:=\max\{n'_0,\omega^{-1/\delta}\}\,,\quad\text{and}\quad
  \gamma:=\frac{\omega^{h^{20}}}{6}>0\,.
  \]
  Let $n>n_0$ and observe that $n^{-\delta}<\omega<\omega_0$ as required by Theorem~\ref{thm:DesignsIntermediate}.

  Now let $\mathcal M$ be a $(\gamma,L)$-quasirandom chest with vertex set
  $[n]=U\dcup V$ satisfying the assumptions of
  Proposition~\ref{prop:Designs-oursetting}.  As
  Theorem~\ref{thm:DesignsIntermediate} only decomposes digraphs and~$\mathcal
  M$ has undirected edges, we need to transform $\mathcal M$ as follows. Replace
  every unoriented edge inside $V$ by an oriented $2$-cycle of the same colour
  and replace every unoriented edge between $V$ and $U$ by an oriented edge from
  $V$ to $U$ of the same colour. Abusing notation, we call this transformed
  digraph $\mathcal M$ as well. Observe that~$\mathcal M$ is a general
  multi-digraph. Naturally, we shall be using the partition $\mathcal
  P'=(P'_1,P'_2)$ with $P'_1=V$ and $P'_2=U$ for~$\mathcal M$, hence $t:=2$ in
  our application of Theorem~\ref{thm:DesignsIntermediate}.
  By assumption we have $|P'_i|\ge \sigma n$ for $i=1,2$.

  We next define the family~$\mathcal H$ of three simple canonical digraphs on vertex set $[q]=[4]$ we are
  decomposing~$\mathcal M$ into. (This was motivated already after Definition~\ref{def:canonical}.)
  Indeed, let $P_1:=\{1,2,3\}$
  and $P_2:=\{4\}$. Let $\mathcal H:=\{H_1, H_2, H_3\}$, where
  \begin{itemize}[leftmargin=*]
  \item $H_1$ has an edge coloured by $1$ directed from $1$ to $2$, a $2$-cycle
    of colour $2$ between $3$ and $1$, a $2$-cycle of colour $3$ between $2$ and
    $3$, an edge in colour $4$ directed from $1$ to $4$, an edge in colour $5$
    directed from $2$ to $4$, and an edge in colour $6$ directed from $3$ to
    $4$;
  \item the edge set of $H_2$ consists of a single edge coloured by $4$ directed from $1$ to $4$;
  \item the edge set of $H_3$ consists of a single edge coloured by $5$ directed from $2$ to $4$.
  \end{itemize}
  In other words, $H_1$ is a directed analogue of the diamond, with directions
  and colours consistent with directions and colours in~$\mathcal M$.  The
  digraphs $H_2$ and $H_3$ consist of a single edge and their purpose is
  that instead of using each edge $\vec{E}_1\dcup E_2\dcup E_3\dcup
  E_6$ exactly once and each edge of $E_4\dcup E_5$ at most once in the
  definition of the diamond core-decomposition, we now want to use each edge
  exactly once in a copy of a graph from $\mathcal H$ as in
  Definition~\ref{def:decomposition}. 

  By construction, an $\mathcal H$-decomposition of~$\mathcal M$ gives the
  desired diamond core-decomposition. Theorem~\ref{thm:DesignsIntermediate}
  provides us with such an $\mathcal H$-decomposition. Hence it remains to check
  that the conditions of this theorem are satisfied:
  we shall show that $(\mathcal M,\mathcal P')$ is $(\mathcal H,\mathcal
  P)$-divisible, $(\mathcal H,
  6\gamma, 2h^{q'-q}\omega)$-regular and $(\mathcal H,\mathcal P,\sqrt[q'h]{\omega}, h)$-vertex-extendable.

  \smallskip
  
  \emph{Divisibility:} The $0$-divisibility condition is implied by the fact
  that by assumption $|\vec{E}_1|=|E_2|=|E_3|=|E_6|\le |E_4|,|E_5|$, where the
  inequality follows from~\ref{itm:designs:iii}--\ref{itm:designs:vi} of
  Proposition~\ref{prop:Designs-oursetting}.  The $1$-divisibility condition
  follows from~\ref{itm:designs:i}--\ref{itm:designs:vi} of
  Proposition~\ref{prop:Designs-oursetting}, the fact that each colour in the
  (original) diamond appears only once, and the fact that~$H_2$ and~$H_3$ are
  single edges. The $2$-divisibility comes from the definition of the
  chest~$\mathcal M$ and the corresponding digraphs $H_1,H_2, H_3$.

  \smallskip

  \emph{Regularity:} For every copy $H'_1$ of~$H_1$ in $\mathcal M$ and
  every edge $H_2'$ of colour $4$ and every edge $H_3'$ of colour $5$ in $\mathcal M$, we
  choose the weights
 \begin{equation}
\label{eq:weightsE}
w_{H'_1}:=\frac{1}{d_2d_3d_4d_5d_6|V||U|}\,,\quad
  w_{H_2'}:=\frac{d_4-d_6}{d_4|V|^2}\,,\quad w_{H_3'}:=\frac{d_5-d_6}{d_5|V|^2}\,.
   \end{equation}
  Note first that $|E_2|=|E_3|=|\vec{E}_1|$ gives $d_2=d_3=\tfrac{|V|}{|V|-1}\cdot 2d_1$; in particular $d_2,d_3>d_1>d$.
  We shall show that all the weights satisfy $w_{H'_1}, w_{H'_2}, w_{H'_3}\in
  [(d^5\sigma^2)\cdot n^{-2},(d^5\sigma^2)^{-1}\cdot n^{-2}]$, and so they are in the allowed range (see Definition~\ref{def:regularity(designs)}) $[2h^{q'-q}\omega n^{-2},(2h^{q'-q}\omega)^{-1}n^{-2}]$ required by Theorem~\ref{thm:DesignsIntermediate}, since $2h^{q'-q}\omega\le d^5\sigma^2$ by the choice of~$\omega$. By plugging directly into~\eqref{eq:weightsE} we get $w_{H'_1}\le \frac{1}{d^5(\sigma n)^2}$, $w_{H'_1}\ge \frac{1}{n^2}$, and $w_{H_2'}\le \frac{1}{(\sigma n)^2}$. For a lower bound on $w_{H_2'}$, we write $$
  w_{H_2'}=\frac{d_4-d_6}{d_4|V|^2}\ge \frac{|E_4|-|E_6|}{|U||V|^3}=\frac{\sum_{u\in U}\left(\deg_{E_4}(u)-\deg_{E_6}(u)\right)}{|U||V|^3}\geBy{\ref{itm:designs:v}} \frac{d^5\sigma^2}{|V|^2}>\frac{d^5\sigma^2}{n^2}
  \;.
  $$
  The bounds for $w_{H_3'}$ follow similarly, using~\ref{itm:designs:vi}. 
  
  Let us now turn to~\eqref{eq:reg1}. By the $(\gamma,L)$-quasirandomness of
  the chest $\mathcal M$ we have that every edge $e$ in~$\mathcal M$ with colour $i\in
  \{1,2,3\}$ is contained in $(1\pm \gamma)^5d_2d_3|V|\cdot d_4d_5d_6 |U|$
  copies of $H_1$ and $0$ copies of $H_2$ and $H_3$. (For $i=1$ this is
  immediate. For $i\in\{2,3\}$, recall that the undirected edges of $E_i$ were
  replaced by oriented $2$-cycles, so that $e$ can play the role of either of the
  two edges of colour~$i$ of $H_1$; this gives $2d_1d_3|V|$, respectively
  $2d_1d_2|V|$, choices for the remaining vertex in~$V$, and this equals
  $(1\pm\gamma)d_2d_3|V|$ because $2d_1=\tfrac{|V|-1}{|V|}d_2=\tfrac{|V|-1}{|V|}d_3$.) Hence,
  \begin{align*}
    \sum_{H\in \mathcal H}\quad\sum_{H'\subseteq\mathcal M; e\in E(H'), H'\sim H}w_{H'}
    =&(1\pm \gamma)^5d_2d_3|V|\cdot d_4d_5d_6 |U|\cdot \frac{1}{d_2d_3d_4d_5d_6|V||U|}
    =(1\pm 6\gamma)\;.
  \end{align*} 
  Every edge~$e$ in~$\mathcal M$ with colour $6$ is contained in $(1\pm
  \gamma)^5d_1d_2d_3d_4d_5|V|^2$ copies of $H_1$ and $0$ copies of $H_2$ and
  $H_3$. Recall that $|\vec{E}_1|=d_1|V|^2$, $|E_6|=d_6|U||V|$, and
  $|\vec{E}_1|=|E_6|$ by assumption.
  Hence,
  \begin{multline*}
    \sum_{H\in \mathcal H}\quad\sum_{H'\subseteq\mathcal M; e\in E(H'), H'\sim H}w_{H'}
    =(1\pm \gamma)^5d_1d_2d_3d_4d_5|V|^2\cdot \frac{1}{d_2d_3d_4d_5d_6|V||U|}\\
    =(1\pm 6\gamma) \frac{d_1|V|}{d_6|U|}
    =(1\pm6\gamma)\frac{d_1|V|^2}{|E_6|}
    =(1\pm 6\gamma)\frac{|\vec{E}_1|}{|E_6|}
    =(1\pm 6\gamma)\;.
  \end{multline*} 
  Next, every edge~$e$ in~$\mathcal M$ of colour $4$ is contained in $(1\pm
  \gamma)^5d_1d_2d_3d_5d_6|V|^2$ copies of $H_1$, in $(|V|-1)(|V|-2)$ copies of
  $H_2$ corresponding to the placement of the two isolated vertices $\{2,3\}$
  in $\mathcal M$, and in~$0$ copies of~$H_3$. Hence, using $|E_4|=d_4|U||V|$
  and $d_1|V|^2=|\vec{E}_1|=|E_6|=d_6|U||V|$, we get
  \begin{align*}
    \sum_{H\in \mathcal H}\quad&\sum_{H'\subseteq\mathcal M; e\in E(H'), H'\sim H}w_{H'}
    =(1\pm \gamma)^5\frac{d_1d_2d_3d_5d_6|V|^2}{d_2d_3d_4d_5d_6|V||U|}+(1\pm \gamma)|V|^2\left(\frac{d_4-d_6}{d_4|V|^2}\right)\\
    &=(1\pm 6\gamma)\bigg(\frac{d_1|V|}{d_4|U|}+\frac{(d_4-d_6)|U|}{d_4|U|}\bigg)
    =(1\pm6\gamma)\frac{d_1|V|^2+d_4|U||V|-d_6|U||V|}{|E_4|} \\
  &=(1\pm 6\gamma)\frac{|E_4|}{|E_4|}
    =(1\pm 6\gamma)\;.
  \end{align*} 
  The calculation for an edge of colour $5$ is analogous to the previous case.
  
  \smallskip
  
  \emph{Extendability:} We provide the details for checking vertex-extendability
  only for $H=H_1$ and $x=1\in V(H_1)$. The other cases are analogous.  We
  choose any collection of pairwise disjoint sets $\{A_i\}_{i\in[4]\setminus
    \{1\}}$ of size $|A_i|\le h$ with $A_i\subset V$ if $i=2,3$ and
  $A_i\subseteq U$ if $i=4$.  Using the notation from the definition of
  vertex-extendability, we have $\NBH_{\mathcal M}^{(2,1,H_1)}(A_2)=\NBHin_{\mathcal
    M,1}(A_2)=\NBHin_{\vec{E}_1}(A_2)$. Similarly $\NBH_{\mathcal
    M}^{(3,1,H_1)}(A_3) =\NBHin_{\mathcal M,2}(A_3)\cap\NBHout_{\mathcal
    M,2}(A_3) =\NBH_{E_2}(A_3)$, where we use that (undirected) edges in~$E_2$
  were replaced by oriented $2$-cycles.  Analogously, $\NBH_{\mathcal
    M}^{(4,1,H_1)}(A_4)=\NBH_{E_4}(A_4)$. Accordingly, we want to lower-bound
  the size of
  \[
    A_1:= \bigcap_{y\in [4]\setminus
    \{1\}}\NBH^{(y,1,H_1)}_{\mathcal M}(A_y)=
    \NBHin_{\vec{E}_1}(A_2)\cap \NBH_{E_2}(A_3)\cap \NBH_{E_4}(A_4)\,.
  \]
  By the definition of $(\gamma,L)$-quasirandomness for the chest~$\mathcal M$, setting
  $S'_1:=A_2$, $S_2:=A_3$, $S_4:=A_4$ and
  $S_1=S_3=S_4'=S_5'=S_6'=S_5=S_6=\emptyset$ (note that these sets have total
  size at most $3h=L$, as required in Definition~\ref{def:quasichest}), and using that
  $\omega\le(\frac{\sigma d}{2})^{5Lhq'}$,
  we get that
  \begin{align*}
    |A_1|\ge (1-\gamma)d_1^{L}d_2^Ld_4^L|V|>\frac
    12d^{3L}\sigma n\ge \sqrt[q'h]{\omega} n\;,
  \end{align*}
  as required for $(\mathcal H,\mathcal P,\sqrt[q'h]{\omega}, h)$-vertex-extendability.
\end{proof}

\section{Probabilistic tools and quasirandomness}
\label{sec:tools}

In this section we collect a number of other tools that we shall need for the
proof of Theorem~\ref{thm:maintechnical}. We start with some standard
probabilistic results concerning properties of the hypergeometric distribution,
McDiarmid's inequality, and Freedman's inequality. We then turn to the
discussion of various quasirandomness properties.

\subsection{Probabilistic tools}\label{ssec:probabilistictools}
Let us quickly recall the hypergeometric distribution. 
\begin{definition}\label{def:hypergeometric}
	Suppose that $X$ is a set of order $n$, and $Y\subset X$ has $\ell$ elements. Let $Z$ be a uniformly random subset of $X$ of size $h$. Then $|Y\cap Z|$ has \emph{hypergeometric distribution}\index{hypergeometric distribution} with parameters $(n,\ell,h)$.
\end{definition}

It is a well-known fact that the hypergeometric distribution is at least as concentrated as the binomial distribution with the same mean. In particular, we have the following bounds which apply to the binomial, and hence also hypergeometric, distributions.

\begin{fact}[Corollary 2.4, Theorems 2.8 and 2.10, and (2.9) in~\cite{JLR}]\label{fact:hypergeometricBasicProperties}
	Given integers $(n,\ell,h)$, the hypergeometric distribution with parameters $(n,\ell,h)$ has mean $\tfrac{\ell h}{n}$. 
    Given $p\in [0,1]$, the binomial distribution with parameters $(n,p)$ has mean $np$.	
	Suppose that $N$ is a random variable whose distribution is either hypergeometric with parameters $(n,\ell,h)$, or binomial with parameters $(n,p)$. Then, for each $\vartheta\in(0,\frac32)$, we have
	\[\Prob\left[|N-\Exp[N]|>\vartheta\Exp[N]\right]<2\exp\left(-\frac{\vartheta^2}{3}\cdot\Exp[N]\right)\;.\]
	We also have
	\[\Prob\left[N-\Exp[N]>s\right]<\exp(-s) \]
	whenever $s\ge 6\Exp[N]$.
\end{fact}
Putting these two bounds together, we get the following observation, which we will use rather often. We have $N=\Exp[N]\pm n^{0.9}$ with probability at least $1-\exp(-n^{0.65})$ for all sufficiently large $n$. To see that this is true, observe that if $\Exp[N]\le n^{0.89}$ then we only need to consider the possibility $N\ge\Exp[N]+n^{0.9}$; since $n^{0.9}\ge6\Exp[N]$ the second bound applies and gives us the desired result. If on the other hand $n^{0.89}\le \Exp[N]\le n$, then we use the first bound with $\vartheta=n^{-0.1}$.

Suppose that $\Omega=\prod_{i=1}^k \Omega_i$ is a product probability space. A measurable function $f:\Omega\rightarrow \mathbb R$ is said to be \emph{$C$-Lipschitz}\index{Lipschitz function} if for each $\omega_1\in\Omega_1,\ldots,\omega_k\in\Omega_k$, for each $i\in[k]$ and each $\omega'_i\in \Omega_i$ we have $|f(\omega_1,\ldots,\omega_i,\ldots,\omega_k)-f(\omega_1,\ldots,\omega'_i,\ldots,\omega_k)|\le C$.
McDiarmid's inequality states that Lipschitz functions are concentrated around their expectation.
\begin{lemma}[McDiarmid's inequality,  \cite{McDiarmid:BoundedDiff}]\label{lem:McDiarmidInequ}
	Let $f:\prod_{i=1}^k \Omega_i\rightarrow \RR$ be a $C$-Lipschitz function. Then for each $t>0$ we have
	$$\Prob\big[|f-\Exp[f]|>t\big]\le 2\exp\left(-\frac{2t^2}{C^2 k}\right)\;.
	$$
\end{lemma}

Let
  $\Omega$ be a finite probability space. A \emph{filtration}\index{filtration} $\mathcal{F}_0$,
  $\mathcal{F}_1$,\dots, $\mathcal{F}_n$ is a sequence of partitions of~$\Omega$ such that $\mathcal{F}_i$ refines $\mathcal{F}_{i-1}$ for all $i\in[n]$. Note that in this finite setting, a function $f:\Omega\rightarrow\RR$ is \emph{$\mathcal{F}_i$-measurable} if $f$ is constant on each part of $\mathcal{F}_i$. Further, for any random variable
  $Y\colon\Omega\rightarrow\RR$ the \emph{conditional expectation}\index{conditional expectation} $\Exp(Y|\mathcal{F}_i)\colon\Omega\rightarrow\RR$ and the \emph{conditional variance} $\Var(Y|\mathcal{F}_i)\colon\Omega\rightarrow\RR$ of $Y$ with respect to $\mathcal{F}_i$ are defined by
  \begin{equation*}
    \begin{split}
      \Exp(Y|\mathcal{F}_i)(x)&=\Exp(Y|X),\\
      \Var(Y|\mathcal{F}_i)(x)&=\Var(Y|X),
    \end{split}
    \qquad
    \text{where $X\in\mathcal{F}_i$ is such that $x\in X$\,.}
  \end{equation*}

Suppose that we have an algorithm which proceeds in $m$ rounds using a new
source of randomness $\Omega_i$ in each round $i$. Then the probability space
underlying the run of the algorithm is $\prod_{i=1}^{m}\Omega_i$. By
\emph{history up to time $t$} we mean a set of the form
$\{\omega_1\}\times\cdots\times\{\omega_t\}\times
\Omega_{t+1}\times\cdots\times \Omega_m$, where $\omega_i\in\Omega_i$. We shall use the
symbol $\hist_{t}$ to denote any particular history of such a form. By a
\emph{history ensemble up to time $t$}\index{history ensemble} we mean any union of histories up to time
$t$; we shall use the symbol $\histens$ to denote any one such. Observe that
there are natural filtrations associated to such a probability space: given
times $t_1<t_2<\dots$ we let $\mathcal{F}_{t_i}$ denote the partition of
$\Omega$ into the histories up to time $t_i$.

The following inequality, a corollary of Freedman's inequality~\cite{Freedman} derived in~\cite{ABCT:PackingManyLeaves}, will be our main concentration tool for analysing random processes. The event $\cE$ will generally be an assertion that various good properties are maintained up to some stage in a random process, and in particular it is important that $\cE$ need not be measurable with respect to any elements of the given filtration.

\begin{corollary}[{\cite[Corollary~\ref{PAPERmanyleaves.cor:freedm}]{ABCT:PackingManyLeaves}}]\label{cor:freedm}
  Let $\Omega$ be a finite probability space, and $(\mathcal{F}_0,
  \mathcal{F}_{1},\dots,\mathcal{F}_{n})$ be a filtration. Suppose that we have
  $R>0$, and for each $1\le i\le n$ we have an $\mathcal{F}_i$-measurable
  non-negative random variable $Y_{i}$, non-negative real numbers $\tilde{\mu},\tilde{\nu}$ and
  an event $\cE$.
  \begin{enumerate}[label=\abc]
  \item\label{cor:freedm:uppertail} Suppose that either $\cE$ does not occur or
    we have $\sum_{i=1}^{n}\Exp\big[Y_{i}\big|\mathcal{F}_{i-1}\big]\le\tilde{\mu}$,
    and $0\le Y_i\le R$ for each $1\le i\le n$. Then
    \begin{equation*}
      \Prob\left[\mathcal{\cE}\text{ and
        }\sum_{i=1}^{n}Y_{i} >2\tilde{\mu}\right] \le\exp\Big(-\frac{\tilde{\mu}}{4R}\Big)\,.
    \end{equation*}
  \item\label{cor:freedm:tails} Suppose that either $\cE$ does not occur or we
    have $\sum_{i=1}^{n}\Exp\big[Y_{i}\big|\mathcal{F}_{i-1}\big]=\tilde{\mu}\pm\tilde{\nu}$,
    and $0\le Y_i\le R$ for each $1\le i\le n$. Then for each $\tilde{\rho}>0$ we have
    \begin{equation*}
      \Prob\left[\cE\text{ and
        }\sum_{i=1}^{n}Y_{i} \neq\tilde{\mu}\pm(\tilde{\nu}+\tilde{\rho})\right] \le2\exp\Big(-\frac{\tilde{\rho}^2}{2R(\tilde{\mu}+\tilde{\nu}+\tilde{\rho})}\Big)\,.
    \end{equation*}
    In particular, if $\tilde{\nu}=\tilde{\rho}=\tilde{\mu}\tilde{\eta}>0$ and $\tilde{\eta}\le\frac12$, then
    \begin{equation*}
      \Prob\left[\cE\text{ and
        }\sum_{i=1}^{n}Y_{i} \neq\tilde{\mu}(1\pm 2\tilde{\eta})\right] \le2\exp\Big(-\frac{\tilde{\mu}\tilde{\eta}^2}{4R}\Big)\,.
    \end{equation*}
  \end{enumerate}
\end{corollary}

\subsection{Quasirandomness}

Here we collect the definitions of various auxiliary quasirandomness properties
that we use in our proof. The last of these, index-quasirandomness, is the
strongest. However, a large fraction of this paper will be taken up by the
analysis of two randomised algorithms for which we need a less general notion of
quasirandomness; in order to keep notation in these parts manageable, we also
give weaker definitions at the level that these two pieces of analysis need.

We start with a generalisation of
$(\gamma,L)$-quasirandomness,
which was already introduced in~\cite{DegPack}. We will need this notion to state various lemmas from~\cite{DegPack} and to perform our own extra analysis of the randomised algorithm from~\cite{DegPack}.

\begin{definition}[diet condition\index{diet condition}]
	\label{def:dietcondition} 
  Let~$H$ be a graph with~$n$ vertices and $p\binom{n}{2}$ edges, and let
  $X\subseteq V(H)$ be any vertex set. We say that the pair $(H,X)$
  satisfies the \emph{$(\gamma,L)$-diet condition}\index{diet condition} if for every set
  $S\subset V(H)$ of at most~$L$ vertices we have 
  \[|\NBH_{H}(S)\setminus
  X|=(1\pm\gamma)p^{|S|}(n-|X|)\,.\]
\end{definition}
Observe that if $H$ is a $(\gamma,L)$-quasirandom graph, and $X$ is a randomly chosen subset of vertices whose size is not too close to $n$, then it is very likely that $(H,X)$ has the $((1+o(1))\gamma,L)$-diet condition. The randomised algorithm of~\cite{DegPack}, which we briefly described at the end of Section~\ref{sec:proofsketch}, has the property that $(H,X)$ is very likely to satisfy the diet condition, where $X$ is the image of the currently embedded vertices of $G_i$ at any given time in the embedding of $G_i$. This is our way of formalising the idea that the image of $G_i$ looks like a random set of vertices.

An important tool in this paper is a randomised algorithm for packing path-forests whose leaves are all embedded, which we state and prove in Section~\ref{sec:pathpack}. As discussed in the proof sketch of Section~\ref{sec:proofsketch}, when we need this algorithm we will be working with a graph whose vertices are split into two parts (the `sides' mentioned in the proof sketch) and densities within and between the sides are not necessarily equal. The following definition is a partite version of quasirandomness and the diet condition.

\begin{definition}[block-quasirandom, block-diet\index{block-quasirandom}\index{block-diet}]\label{def:block-quasi}
Suppose that $L\in \NN$ and $\gamma>0$ are given.
  Let~$H$ be a graph, and $U\subset V(H)$ be a subset. We say that $(H,U)$ is \emph{$(\gamma,L)$-block-diet} on $\Vmin\dcup
  \Vplus=V(H)$ if for every pair of sets $S_{\boxminus}\subset \Vmin$,
  $S_{\boxplus}\subset \Vplus$ with $|S_{\boxminus}|+|S_{\boxplus}|\le L$ we
  have
  \begin{align*}
    \Big|\big(\NBH_{H[\Vmin,\Vplus]}(S_{\boxminus}) \cap \NBH_{H[\Vplus]}(S_{\boxplus})\big)\setminus U\Big|
    &= (1\pm\gamma)|\Vplus\setminus U| \cdot d(H[\Vplus])^{|S_{\boxplus}|} \cdot d(H[\Vmin,\Vplus])^{|S_{\boxminus}|}\,,\\
    \Big|\big(\NBH_{H[\Vmin]}(S_{\boxminus}) \cap \NBH_{H[\Vplus,\Vmin]}(S_{\boxplus})\big)\setminus U\Big|
    &= (1\pm\gamma)|\Vmin\setminus U| \cdot d(H[\Vmin])^{|S_{\boxminus}|} \cdot d(H[\Vmin,\Vplus])^{|S_{\boxplus}|}\,,
  \end{align*}
  where $d(H[\Vplus])$ is the density of $H[\Vplus]$, $d(H[\Vmin])$ is the
  density of $H[\Vmin]$, and $d(H[\Vmin,\Vplus])$ is the bipartite density of
  $H[\Vmin,\Vplus]$. If $(H,\emptyset)$ is $(\gamma,L)$-block-diet, we say $H$ is \emph{$(\gamma,L)$-block-quasirandom}. 
\end{definition}

We remark that we shall only consider whether some $(H,U)$ is block-diet when we already know
that~$H$ is block-quasirandom.

Our final quasirandomness condition is a good deal more complicated. Recall from the proof sketch in Section~\ref{sec:proofsketch} that eventually in Stage~G we will draw an auxiliary \emph{chest} and argue that a generalised design in this chest corresponds to completing our perfect packing. We need the chest to be regular and extendable (as defined in Section~\ref{sec:designs}) and to obtain this we ask for a partite quasirandomness condition rather similar to that of block-quasirandomness. In particular, we should be able to control the size of common neighbourhoods (in specified colours) of several vertices of the chest. The following index-quasirandomness controls the sizes of these sets when they are within the part $V$ whose vertices $i\in[\tfrac{n}{2}]$ correspond to terminal pairs $\{\boxminus_i,\boxplus_i\}$ of $H$. 

Recall that a vertex $\ell\in U$ in the chest corresponds to a graph $G_\ell$ (and we still need to embed a path-forest in $G_\ell$). There are edges of three different colours leaving $\ell$ in the chest, which tell us which vertices of $\Vmin$ and of $\Vplus$ have not been used to embed vertices of $G_\ell$, and which terminal pairs need connecting by a path. We do not need to know the graph structure of $G_\ell$ in order to know where these edges go; we simply need to know the set $U_\ell$ of vertices of $H$ used in the embedding, and the set $A_\ell$ which contains the indices $i\in[\tfrac{n}{2}]$ of terminal pairs $\{\boxminus_i,\boxplus_i\}$ to which paths of $G_\ell$ are anchored.

In order to stick to this definition throughout the stages of our packing, we need to enhance it a little. From Stage~C onwards, we have a collection of path-forests in various different graphs, indexed by $\kappa$ in the following definition, to pack. Each of these graphs has a used set $U_\ell$. However, only some of them, indexed by $\kappa'$, will be the graphs whose packing we complete in Stage~G, and only these graphs have a set $A_\ell$ of terminal pairs.

In the following definition we have a size parameter $L$ and an error parameter $\gamma$, which play much the same role as the similar parameters in the above quasirandomness conditions; note that here, in contrast to those conditions, the size parameter comes first in the notation. In addition we have two density parameters $d_1$ and $d_2$. The density $d_1$ is (approximately) the density within each of $\Vmin$ and $\Vplus$; we will always have very nearly (but not necessarily exactly) the same number of edges within each of these sets. The density $d_2$ is the density between $\Vmin$ and $\Vplus$. We have sets $S_1$ and $S_2$ of vertices of $H$, whose common neighbourhoods in respectively $\Vmin$ and $\Vplus$ we want to consider. In addition, we have sets $T_1,T_2\subset\kappa$ and $T_3\subset\kappa'$. We want to know which vertices of $\Vmin$ are not used in the embedding of graphs $G_\ell$ with $\ell\in T_1$, and similarly vertices of $\Vplus$ for $\ell\in T_2$. In addition, we want to know which $i\in[\tfrac{n}{2}]$ are
indices of terminal pairs for each $G_\ell$ with $\ell\in T_3$.

A coloured common neighbourhood in the chest is thus the same thing as a set of vertices $i\in[\tfrac{n}{2}]$ which satisfy all of the following. $\boxminus_i$ is in the common neighbourhood of $S_1$, and not used by any graph $G_\ell$ with $\ell\in T_1$. Similarly, $\boxplus_i$ is in the common neighbourhood of $S_2$ and not used by any graph $G_\ell$ with $\ell\in T_2$. Furthermore, $\{\boxminus_i,\boxplus_i\}$ is a terminal pair of $G_\ell$ for each $\ell\in T_3$. The following definition controls the sizes of all such sets.

\begin{definition}[index-quasirandom\index{index-quasirandom}]\label{def:index-quasirandom}
  Let~$\tilde H$ be a graph on~$n$ vertices for~$n$ even, and
  $\Vmin=\{\boxminus_i:i\in[\frac{n}2]\}$,
  $\Vplus=\{\boxplus_i:i\in[\frac{n}2]\}$ be disjoint vertex sets such that
  $V(\tilde H)=\Vmin\dcup\Vplus$. Further, let $(U_\ell)_{\ell\in\kappa}$ be a
  collection of vertex sets $U_\ell\subset V(\tilde H)$, and
  $(A_\ell)_{\ell\in\kappa'}$ be a collection of index sets with
  $A_\ell\subset[\frac{n}{2}]$ and $\kappa'\subset\kappa$. 
  For a vertex set $S\subset V(\tilde H)$ and an index set $T\subset \kappa$ we define
  \[\mathbb{N}_{\tilde H}(S,T):=\NBH_{\tilde H}(S)\setminus\bigcup_{\ell\in T}U_\ell\,.\]
  One should think of $\mathbb{N}_{\tilde{H}}(S,T)$ as being the vertices of $\tilde H$ adjacent to all members of $S$ and not used in the embedding of any $G_\ell$ with $\ell\in T$. 
  
  We say that the triple $(\tilde H,\Vmin,\Vplus)$ is
  \emph{$(L,\gamma,d_1,d_2)$-index-quasirandom} with respect to
  $(U_\ell)_{\ell\in\kappa}$ and $(A_\ell)_{\ell\in\kappa'}$,
  if for all $S_1,S_2\subset V(\tilde H)$, for all $T_3\subset\kappa'$ and all $T_1,T_2\subset\kappa$ with $T_1,T_2,T_3$ pairwise disjoint
  such that \[|S_1|,|S_2|,|T_1|,|T_2|,|T_3|\le L\] the set
  \index{$\mathbb{U}_{\tilde H}(S_1,S_2,T_1,T_2,T_3)$}
  \begin{equation}\label{eq:U5param}
  \mathbb{U}_{\tilde H}(S_1,S_2,T_1,T_2,T_3):=
    \left\{i\in[\tfrac{n}{2}]:
      \boxminus_i\in\mathbb{N}_{\tilde H}(S_1,T_1),
      \boxplus_i\in\mathbb{N}_{\tilde H}(S_2,T_2),
      i\in\bigcap_{\ell\in T_3}A_\ell
    \right\}
 \end{equation}
  satisfies
  \[|\mathbb{U}_{\tilde H}(S_1,S_2,T_1,T_2,T_3)|=
    (1\pm\gamma)
    d_1^{|S_1\cap\Vmin|+|S_2\cap\Vplus|}
    d_2^{|S_1\cap\Vplus|+|S_2\cap\Vmin|}
    \cdot\frac{n}{2}
    \cdot\prod_{\ell\in T_1\cup T_2}\Big(1-\frac{|U_\ell|}{n}\Big)
    \cdot\prod_{\ell\in T_3}\frac{|A_\ell|}{n/2}\,.
  \]
\end{definition}

As mentioned, through the stages of our packing we will maintain index-quasirandomness. However, to apply our path packing theorem, we need the weaker block-quasirandomness and block-diet. The following lemma checks that we indeed get it.

\begin{lemma}\label{lem:indexquasi-implies-blockquasi}
	Suppose that the triple $(H, \Vmin, \Vplus)$ is $(L, \gamma, d_1,d_2)$-index-quasirandom with respect to sets $(U_{\ell})_{\ell\in \kappa}$ and $(A_\ell)_{\ell\in \kappa'}$. Suppose that $\gamma$ is sufficiently small compared to $L^{-1}$ and that $n=|V(H)|$ is sufficiently large given $L,\gamma,d_1,d_2$. Then we have
	\[d(H[\Vmin]), d(H[\Vplus])=(1\pm 2\gamma)d_1\quad\text{ and }\quad d(H[\Vmin, \Vplus])=(1\pm \gamma)d_2\,.\]
	Moreover, $H$ is $\big((L+3)\gamma,L\big)$-block-quasirandom. In addition, for each $\ell\in\kappa$, we have
	\[|\Vmin \setminus U_\ell|, |\Vplus\setminus U_\ell|=(1\pm \gamma)\left(\tfrac{n}{2}-\tfrac{|U_\ell|}{2}\right)\,,\]
	and the pair $(H, U_\ell)$ is $\big((L+3)\gamma,L\big)$-block-diet.
	\end{lemma}

\begin{proof}
	    Observe that $\big|\NBH_H(v)\cap\Vmin\big|=\big|\mathbb{U}_{H}(\{v\},\emptyset,\emptyset,\emptyset,\emptyset)\big|$. So we have
	    \begin{align*}
		d(H[\Vmin])\cdot \binom{n/2}{2}=\frac{1}{2}\sum_{v\in \Vmin}|\NBH_{H}(v)\cap\Vmin|=\frac{1}{2}\cdot \frac{n}{2}(1\pm \gamma)d_1\frac{n}{2}= (1\pm (\gamma+\tfrac{\gamma}{4L}))d_1\cdot \binom{n/2}{2}
		\end{align*}
		for large enough $n$,
leading to $d(H[\Vmin])=(1\pm (\gamma+\tfrac{\gamma}{4L}))d_1$.		
		Analogously, we obtain $d(H[\Vplus])=(1\pm (\gamma+\tfrac{\gamma}{4L}))d_1$, and 
	    \begin{align*}
		d(H[\Vmin,\Vplus])\cdot \left(\frac{n}{2}\right)^2 = \sum_{v\in \Vmin}|\NBH_{H}(v)\cap\Vplus|=\frac{n}{2}(1\pm \gamma)d_2\frac{n}{2} = (1\pm \gamma)d_2\cdot \left(\frac{n}{2}\right)^2
		\end{align*}
giving $d(H[\Vmin,\Vplus])=(1\pm \gamma)d_2$. This proves the first part of the lemma.

		For the second part, let $U=U_\ell$. We have
		\begin{align*}
			|\Vmin\setminus U|=|\mathbb{U}_{H}(\emptyset, \emptyset,\{\ell\},\emptyset,\emptyset)|=(1\pm \gamma)(1-\frac{|U|}{n})\cdot \frac{n}{2}=(1\pm \gamma)\left(\tfrac{n}{2}-\tfrac{|U|}{2}\right),
		\end{align*}
and $|\Vplus\setminus U|=(1\pm \gamma)\left(\tfrac{n}{2}-\tfrac{|U|}{2}\right)$ holds analogously.		
		
		Let $S_\boxplus$ and $S_\boxminus$ be as in the definition of block-diet. Setting $S_2:= S_{\boxplus}\cup S_{\boxminus}$ we have
		\begin{multline*}
			|(\NBH_{H[\Vmin,\Vplus]}(S_{\boxminus})\cap \NBH_{H[\Vplus]}(S_{\boxplus}))\setminus U|=|\mathbb{U}_{H}(\emptyset,S_2,\emptyset, \{\ell\}, \emptyset )|
			=(1\pm \gamma)d_1^{|S_{\boxplus}|}d_2^{|S_{\boxminus}|}\cdot (1-\frac{|U|}{n}) \frac{n}{2}\\
                  \begin{aligned}
			&=(1\pm \gamma)\left(\frac{d(H[\Vplus])}{(1\pm (\gamma+\tfrac{\gamma}{4L}))}\right)^{|S_{\boxplus}|}\left(\frac{d(H[\Vmin
				,\Vplus])}{(1\pm \gamma)}\right)^{|S_{\boxminus}|}\cdot \frac{|\Vplus\setminus U|}{1\pm \gamma}  .
				\end{aligned}
		\end{multline*}
Now, note that for any $\varepsilon \in (0,1)$ and $0<x<\frac{\varepsilon}{2}$ we have 
$\frac{1}{1\pm x} = 1 \pm (1+\varepsilon)x$.
Moreover, for any positive integer $K$, if $x>0$ is sufficiently small compared to $\varepsilon K^{-1}$,
then, by applying the binomial theorem,
we obtain $(1\pm x)^K = 1 \pm (1 + \varepsilon)K x$. Therefore, 
choosing $\varepsilon$ and then $\gamma$ sufficiently small,
and using $|S_{\boxplus}|+|S_{\boxminus}|\le L$,
we can continue our estimate as follows:
\begin{align*}
			& |(\NBH_{H[\Vmin,\Vplus]}(S_{\boxminus})\cap \NBH_{H[\Vplus]}(S_{\boxplus}))\setminus U| \\
			& = (1\pm (1+\varepsilon)(\gamma+\frac{\gamma}{4L}))^{L+2} d(H[\Vplus])^{|S_{\boxplus}|}\cdot d(H[\Vmin,\Vplus])^{|S_{\boxminus}|}|\Vplus\setminus U|\ \\
			& = (1\pm (1+\varepsilon)^2(\gamma+\frac{\gamma}{4L})(L+2)) \cdot d(H[\Vplus])^{|S_{\boxplus}|}\cdot d(H[\Vmin,\Vplus])^{|S_{\boxminus}|}|\Vplus\setminus U|\\
            & = (1\pm (L+3)\gamma) d(H[\Vplus])^{|S_{\boxplus}|}\cdot d(H[\Vmin,\Vplus])^{|S_{\boxminus}|}|\Vplus\setminus U|.
\end{align*}		
		The estimate for $|(\NBH_{H[\Vmin]}(S_{\boxminus})\cap \NBH_{H[\Vmin,\Vplus]}(S_{\boxplus}))\setminus U|$ is done analogously.

		Finally, we show that $H$ is $\big((L+3)\gamma,L\big)$-block-quasirandom, i.e.\ that $(H,\emptyset)$ is $\big((L+3)\gamma,L\big)$-block-diet; this is the case $U=\emptyset$ of the above, taking $T_1=T_2=\emptyset$ in place of $T_1$ or $T_2=\{\ell\}$, so that $|\Vplus\setminus\emptyset|=|\Vplus|=\tfrac n2$ exactly, with no error to absorb from this quantity. Let $S_\boxminus,S_\boxplus$ be as in the definition of block-diet, and set $S_2:=S_\boxplus\cup S_\boxminus$. Directly from the index-quasirandomness of $(H,\Vmin,\Vplus)$ (with $T_1=T_2=T_3=\emptyset$), we have
		\begin{align*}
			|\NBH_{H[\Vmin,\Vplus]}(S_{\boxminus})\cap \NBH_{H[\Vplus]}(S_{\boxplus})|=|\mathbb{U}_{H}(\emptyset,S_2,\emptyset, \emptyset, \emptyset )|
			=(1\pm \gamma)d_1^{|S_{\boxplus}|}d_2^{|S_{\boxminus}|}\cdot \frac{n}{2}\,.
		\end{align*}
		Substituting $d_1=\big(1\pm(1+\varepsilon)(\gamma+\tfrac{\gamma}{4L})\big)d(H[\Vplus])$ and $d_2=(1\pm(1+\varepsilon)\gamma)d(H[\Vmin,\Vplus])$ as above, and combining the leading $(1\pm\gamma)$ factor with the (at most $L$) factors coming from the exponents $|S_\boxplus|,|S_\boxminus|$, we obtain, exactly as in the estimate above but with one fewer factor since there is no $|\Vplus\setminus U|$ term to account for,
		\[|\NBH_{H[\Vmin,\Vplus]}(S_{\boxminus})\cap \NBH_{H[\Vplus]}(S_{\boxplus})|=\big(1\pm(1+\varepsilon)^2(\gamma+\tfrac{\gamma}{4L})(L+1)\big)d(H[\Vplus])^{|S_{\boxplus}|}\cdot d(H[\Vmin,\Vplus])^{|S_{\boxminus}|}\cdot \frac{n}{2}\,,\]
		which for $\varepsilon,\gamma$ chosen sufficiently small is $\big(1\pm(L+2)\gamma\big)d(H[\Vplus])^{|S_{\boxplus}|}\cdot d(H[\Vmin,\Vplus])^{|S_{\boxminus}|}\cdot \tfrac{n}{2}$, and since $|\Vplus|=n/2$ exactly and $(L+2)\gamma\le(L+3)\gamma$, this is what is required. The estimate for $|\NBH_{H[\Vmin]}(S_{\boxminus})\cap \NBH_{H[\Vmin,\Vplus]}(S_{\boxplus})|$ is done analogously.
	\end{proof}

\section{The proof of the main theorem}\label{sec:theproof}

In this section we give our proof of Theorem~\ref{thm:maintechnical}. We first
set up constants and preprocess the graphs $(G_s)_{s\in\cG}$. In the key
Section~\ref{ssec:packingstages} we introduce the seven packing stages (Stages~A--G) in seven
corresponding lemmas (Lemmas~\ref{lem:StageA}--\ref{lem:StageG}), which readily yield Theorem~\ref{thm:maintechnical}.

Before we start, we introduce some notation.
Suppose that $P$ is a path of length $\ell$. For any such path, we will
tacitly (and arbitrarily) fix a left-right orientation. Now, given
$h\in\{0,1,\ldots,\ell\}$, let $\leftpath_h(P)$\index{$\leftpath_h(P)$} be the
leftmost subpath of $P$ of length $h$ (that is, with $h$ edges and so $h+1$
vertices). If $h=0$, then $\leftpath_h(P)$ is the leftmost vertex of $P$ (which we
will typically view as a vertex rather than a $1$-vertex graph). We define
$\rightpath_h(P)$\index{$\rightpath_h(P)$} analogously.
If $h<\ell$ and $h$ has the same parity as $\ell$, we define
$\middlepath_h(P)$\index{$\middlepath_h(P)$} as
$P-(\leftpath_{(\ell-h-2)/2}(P)\cup\rightpath_{(\ell-h-2)/2}(P))$. Hence the
length of $\middlepath_h(P)$ is $h$.

\subsection{Constants}
Recall that we are given $D_0$, $d$ and $\deltnonspanning$. We may and do assume in addition that $\deltnonspanning\le\tfrac{1}{100}$: as we explain in the paragraph following~\eqref{eq:CONSTANTS}, the parameter $\deltnonspanning$ enters Definition~\ref{def:family} only through requirements which become weaker as $\deltnonspanning$ decreases, so we may replace $\deltnonspanning$ by $\min\big\{\deltnonspanning,\tfrac{1}{100}\big\}$ without loss of generality. Throughout Sections~\ref{sec:theproof}--\ref{sec:StageG}, we shall use constants satisfying the following hierarchy, where the numbers $D_0$, $\LF$, $\LE$, $\LD$, $\LC$ and $D$ are integers and the remaining constants are positive real numbers (we will describe how to calculate these constants at the end of this subsection):
\begin{multline}\label{eq:CONSTANTS}
	\index{$D_0$}
	\index{$\sigmKJ$}
	\index{$\LF,\LE,\LD,\LC$}
	\index{$\sigmJnula,\sigmJjedna$}
	\index{$\gamJ,\gamNew,\gamAnchor,\gamcore$}
	\frac{1}{D_0},d,\deltnonspanning\gg\sigmKJ\gg\sigmJnula\gg \frac{1}{\LF}\ge \frac{1}{\LE}\ge \frac{1}{\LD} \\
	\ge \frac{1}{\LC}\ge\frac{1}{D}\gg\gamJ\gg\sigmJjedna\gg\gamNew\gg\gamAnchor\gg\gamcore\gg c,\iniquasi\gg \frac{1}{n_0}\;.
\end{multline}
The relations $\ge$ between the size parameters $\LF$, $\LE$, $\LD$, $\LC$ and $D$ are genuinely only inequalities: no separation between these parameters is required anywhere in the proof, and the explicit choice we make at the end of this subsection satisfies them with small constant factors. The relations $\gg$ should be understood as follows. Given $D_0^{-1},d,\deltnonspanning$, there is a function $f(D_0^{-1},d,\deltnonspanning)$ which is monotone decreasing as any of the parameters decreases. We may choose any $0<\sigmKJ<f(D_0^{-1},d,\deltnonspanning)$. In fact we shall be marginally more careful and choose $\sigmKJ<f\big(D_0^{-1},d,\tfrac12\deltnonspanning\big)$, and likewise make every subsequent choice as though $\deltnonspanning$ were $\tfrac12\deltnonspanning$; this costs us nothing and gives us the freedom, used below, to decrease $\deltnonspanning$ slightly. Given this choice, we may then choose any $\sigmJnula$ sufficiently small (in the same sense), and so on. We will finally be supplied with a parameter $n\ge n_0$ by Theorem~\ref{thm:maintechnical}, and we will want to assume that $\sigmKJ n$, $\sigmJnula n$ and $\sigmJjedna n$ are all integers. Observe that for any given $x>0$ there is a real number $y$ with $x-n_0^{-1}<y\le x$ such that $yn$ is an integer. So what we will in fact do is the following: we do not at this stage specify $\sigmKJ$ precisely, but rather an interval of the form $\big[\tfrac12x,x\big]$ where $x$ is sufficiently small given $D_0^{-1},d,\deltnonspanning$. Again we do not choose $\sigmJnula$ precisely, but choose an interval $\big[\tfrac12y,y\big]$, where $y$ is sufficiently small for $\tfrac12x$. We then insist that $\LF^{-1}$ is sufficiently small for $\tfrac12y$; the remaining size parameters $\LE$, $\LD$, $\LC$ and $D$ are then fixed as described below, and we continue in the same way from $\gamJ$ onwards. Once we proceed to $\sigmJjedna$, we again choose an interval. With this construction we in particular have fixed values for all the constants listed above except $\sigmKJ$, $\sigmJnula$ and $\sigmJjedna$, and a guarantee that whatever values we choose in the specified intervals for these three constants, our calculations will work. 

The constants promised by Theorem~\ref{thm:maintechnical} are $n_0$, $L:=4D+3$,
$c$, $\iniquasi$. Suppose now that the remaining parameter~$n\ge n_0$ of
Theorem~\ref{thm:maintechnical} is given. We now fix $\sigmKJ$, $\sigmJnula$ and $\sigmJjedna$ in their specified intervals such that $\sigmKJ n$, $\sigmJnula n$ and $\sigmJjedna n$ are all integers.  Similarly, if necessary, we decrease $\deltnonspanning$ by less than $1/n$ such that $\deltnonspanning n$ is an integer. It is important that we decrease rather than increase $\deltnonspanning$ here. Indeed, $\deltnonspanning$ enters Definition~\ref{def:family} only through the requirements $|\cJ|\ge\deltnonspanning n$, \ref{enu:familyOrderJK}, \ref{enu:familyPathsJ} and~\ref{enu:defSpecLeaves}, each of which becomes weaker as $\deltnonspanning$ decreases; consequently the family we are given still lies in $\OurPackingClass(n,e(H);\deltnonspanning,c,D_0)$ for the new value of $\deltnonspanning$, whereas increasing $\deltnonspanning$ could destroy~\ref{enu:familyOrderJK}. Nor does this change disturb the hierarchy~\eqref{eq:CONSTANTS}: all the constants there were chosen sufficiently small given $\tfrac12\deltnonspanning$, and since $n\ge n_0\ge2\deltnonspanning^{-1}$ the new value of $\deltnonspanning$ is still at least $\tfrac12\deltnonspanning$.

We now summarise the meanings of these various constants. In this summary we shall refer to some sets which are only defined in Section~\ref{ssec:definingSubgraphs} below in order to keep these definitions in one place. More precisely, 
we refer to sets $\cJ_0,\cJ_1,\cJ_2\subset\cJ$ and use a set of paths $\SpecPaths_s$ which will be constructed from the paths in $\BasicPaths_s$; recall that the latter have length 11.

\begin{description}
 \item[$D_0,\DParity,d$] The graphs which we want to pack are all $D_0$-degenerate, and the graph $H$ into which we pack has at least $dn^2$ edges. All vertices in $\SpecLeaves$ have odd degree $\DParity\le D_0$.
 \item[$\deltnonspanning$] The graphs $G_s$ with $s\in\cJ\cup\cK$ have at most $(1-\deltnonspanning)n$ vertices.
 \item[$\sigmKJ$] The size of $\cJ$ (i.e.\ the number of graphs with many bare paths) is $\sigmKJ n$. For each $s\in\cJ_0\cup\cJ_2$ we have $|\SpecPaths_s|=\sigmKJ n$.
 \item[$\sigmJnula$] The size of $\cJ_0$ is $\sigmJnula n$.
 \item[$\sigmJjedna$] The size of $\cJ_1$ is $\sigmJjedna n$, and if $s\in\cJ_1$ then $|\SpecPaths_s|=\sigmJjedna n$.
 \item[$D,\LC,\LD,\LE,\LF$] In quasirandomness conditions in our various stages, we bound from above the number of vertices for which we control the size of a common neighbourhood and the number of guest graphs whose common image we look at. This bound is $D$ after Stage~A (for consistency with~\cite{DegPack}), and then $\LC$ after Stage~C, then $\LD$ after Stage~D, and so on.
 \item[$\gamcore,\gamAnchor,\gamNew,\gamJ$] The accuracy of our quasirandomness conditions is measured by $\gamcore$ after Stage~A (and by $2\gamcore$ after Stage~B and $100\gamcore$ after Stage~C), and is then $\gamAnchor$ after Stage~D, and so on.
 \item[$c$] We impose the bound $\Delta(G_s)\le\frac{cn}{\log n}$ for all guest graphs $G_s$.
 \item[$\iniquasi$] We require $H$ to be $(\iniquasi,4D+3)$-quasirandom.
 \item[$n_0$] We require $n\ge n_0$.
\end{description}

Finally, we explain how to calculate these constants in terms of the lemmas for Stages A--G given later in this section. We are given $D_0,d,\deltnonspanning$. We have various inequalities through the proof which (given these constants) require $\sigmKJ$ to be sufficiently small, and we choose the minimum of these requirements (for the upper bound of the interval in which $\sigmKJ$ is to be chosen). Similarly, $\sigmJnula$ needs to be sufficiently small for various further inequalities (which are monotone in $\sigmKJ$; we choose $\sigmJnula$ small enough for the lower bound on $\sigmKJ$).

Now Lemma~\ref{lem:StageG} (Stage~G), once these constants are fixed, gives us a lower bound for $\LF$; we take $\LF$ to be the maximum of this lower bound and $20$. We can now set $\LE=\LF$, $\LD=2\LE$, $\LC=6\LD$ and $D=\max(D_0,3\LC)$. Lemma~\ref{lem:StageG} in addition requires $\gamJ$ sufficiently small. Next, we look at Lemma~\ref{lem:StageFNew} (Stage~F). In the proof of this lemma, $\gamJ$ is assumed to be sufficiently small as well; we fix $\gamJ$ to be the minimum of these two requirements. Then Lemma~\ref{lem:StageFNew} requires that $\sigmJjedna$ is sufficiently small; in addition $\sigmJjedna$ is required to be small enough for various inequalities in the earlier stages (which do not involve constants coming later in the hierarchy). We fix an upper bound on $\sigmJjedna$ sufficiently small for these requirements, and as explained a lower bound of half the upper bound. Finally, Lemma~\ref{lem:StageFNew} gives us an upper bound on $\gamNew$. We now make a succession of similar choices: $\gamNew$ is chosen to satisfy this upper bound and in addition sufficiently small for Lemma~\ref{lem:StageENew}, and that lemma also gives us an upper bound on $\gamAnchor$. We repeat similar processes for Lemma~\ref{lem:StageDNew} (Stage~D) and Lemma~\ref{lem:StageCNew} (Stage~C) to obtain $\gamAnchor$, and an upper bound on $\gamcore$. Finally, both Lemmas~\ref{lem:StageA} and~\ref{lem:StageB} require $\gamcore$ to be sufficiently small, and we choose the minimum of these. At last, Lemma~\ref{lem:StageA} gives values for $c$ and $\iniquasi$.

Observe that this description explicitly chooses constants in the order of the hierarchy~\eqref{eq:CONSTANTS}. In each place, we specify a choice which works, but any smaller choice (and then continuing the constant choice process as above) will also work, as is implied by the hierarchy notation. From this point onwards we shall use these constants without providing explicit quantifiers (for example in the statements of the lemmas for the different packing stages below).

\subsection{Correcting inequalities and sizes}\label{ssec:IneqEqual}
It is convenient to assume in our proof that several of the inequalities in the definition of $\OurPackingClass$ (Definition~\ref{def:family}) hold with equality, and to reduce the size of some sets from what is guaranteed by Definition~\ref{def:family}.

We want $\sum_{s\in\cG}e(G_s)=e(H)$, so that our aim will always be a perfect packing. By Definition~\ref{def:family}\ref{enu:familyEdgesAll}, we have $\sum_{s\in\cG}e(G_s)\le e(H)$. So, we add $e(H)-\sum_{s\in\cG}e(G_s)$ copies of the graph $K_2$ to our family and correspondingly increase the size of $\cG$. Observe that a perfect packing of the new family contains a packing of the old family.

We want
\begin{equation}\label{eq:sizeJ}
|\cJ|=\sigmKJ n\;.
\end{equation}
Since we are given to start with $|\cJ|\ge\deltnonspanning n>\sigmKJ n$ by Definition~\ref{def:family}, we arbitrarily remove indices from $\cJ$ until we obtain this.

We also want
\begin{equation}\label{eq:sizeBasicPaths}
|\BasicPaths_s|=\sigmKJ n\quad\mbox{for each $s\in\cJ$}\;.
\end{equation}
By Definition~\ref{def:family}\ref{enu:familyPathsJ} we start with $|\BasicPaths_s|\ge\deltnonspanning n>\sigmKJ n$, so we arbitrarily remove paths from $\BasicPaths_s$ until we obtain this.

We want $(1+\deltnonspanning)n\le\big|\bigcup_{s\in\cK}\SpecLeaves_s\big|\le (1+2\deltnonspanning)n$. The lower bound is given by Definition~\ref{def:family}\ref{enu:defSpecLeaves}. If the upper bound is exceeded, we arbitrarily remove indices from $\cK$ one at a time until it is no longer exceeded. Since removing one index decreases (by Definition~\ref{def:family}\ref{enu:defSpecLeaves}) the size of $\bigcup_{s\in\cK}\SpecLeaves_s$ by at most $\tfrac{cn}{\log n}<\deltnonspanning n$, the result is the desired pair of inequalities.

We want that each $G_s$ with $s\in\cJ\cup\cK$ has exactly $(1-\deltnonspanning)n$ vertices. So, we add $(1-\deltnonspanning)n-v(G_s)$ isolated vertices to each $G_s$ with $s\in\cJ\cup\cK$. Observe that a packing of the new family of graphs immediately gives a packing of the old family by ignoring the added isolated vertices. Similarly, we obtain $v(G_s)=n$ for each $s\in\cG\setminus(\cJ\cup\cK)$.

Abusing notation slightly, we will continue to use the same letters $G_s,\cJ$ and so on for the modified family of graphs.

\subsection{Subgraphs of \texorpdfstring{$(G_s)_{s\in\cG}$}{(Gs)} and index sets used in the packing stages}\label{ssec:definingSubgraphs}
Our packing stages, detailed below, will refer to the following subgraphs of our guest graphs $(G_s)_{s\in\cG}$.
Fix two disjoint families $\cJ_0,\cJ_1\subset \cJ$\index{$\cJ_{0}$}\index{$\cJ_{1}$} with 
\begin{equation}\label{eq:sizeJ0J1}
|\cJ_0|=\sigmJnula n\quad \mbox{and}\quad |\cJ_1|=\sigmJjedna n\;.
\end{equation}
Set \index{$\cJ_2$}$\cJ_2:=\cJ\setminus (\cJ_0\cup \cJ_1)$, meaning that 
\begin{equation}\label{eq:sizeJ2}
|\cJ_2|=(\sigmKJ-\sigmJnula-\sigmJjedna)n\;.
\end{equation}
The graphs $G_s$ indexed by these three families will play quite different roles in the packing. 
Recall that each path in $\BasicPaths_s$, for $s\in \cJ$,
has length 11. For $s\in \cJ_1$, let $\BasicPaths'_s\subset \BasicPaths_s$ be an arbitrary family of $\sigmJjedna n$ paths.
For $s\in\cJ$, we define\index{$\SpecPaths_s$}
\[\SpecPaths_s=\begin{cases}
\BasicPaths_s & \text{ if } s\in \cJ_0\,, \\
\BasicPaths'_s & \text{ if } s\in \cJ_1\,, \\
\big\{\middlepath_7(P):P\in\BasicPaths_s\big\} \quad & \text{ if } s\in  \cJ_2\,.
\end{cases}
\]
For $s\in \cK$, let\index{$G_s^\spadesuit$}
\[G_s^\spadesuit:=G_s-\SpecLeaves_s\,.\]
Next, for each $s\in \cJ$, we are going to define a graph $G_s^\parallel\subset G_s$\index{$G_s^\parallel$} by trimming off parts of the bare paths, i.e., for $s\in\cJ$ let
\[G_s^\parallel:=G_s- \bigcup_{P\in\SpecPaths_s \atop P\text{ has length }\ell} V(\middlepath_{\ell-2}(P)).\]
Summarising, we have
\begin{align}
\label{eq:familyOrderAllEquality}
v(G_s)&=n \quad \mbox{for each $s\in \cG\setminus(\cJ\cup\cK)$}\;,\\
\label{eq:familyOrderKEquality}
v(G_s^\spadesuit)&=(1-\deltnonspanning)n-|\SpecLeaves_s| \quad \mbox{for each $s\in \cK$}\;,\\
\label{eq:familyOrderJ1Equality}
v(G_s^\parallel)&=(1-\deltnonspanning-10\sigmJjedna)n \quad \mbox{for each $s\in \cJ_1$}\;,\\
\label{eq:familyOrderJ0Equality}
v(G_s^\parallel)&=(1-\deltnonspanning-10\sigmKJ)n \quad \mbox{for each $s\in \cJ_0$}\;,\\
\label{eq:familyOrderJ2Equality}
v(G_s^\parallel)&=(1-\deltnonspanning-6\sigmKJ)n \quad \mbox{for each $s\in \cJ_2$}\;.
\end{align}

\subsection{The packing stages}\label{ssec:packingstages}
\begin{figure}
	\includegraphics[width=\textwidth]{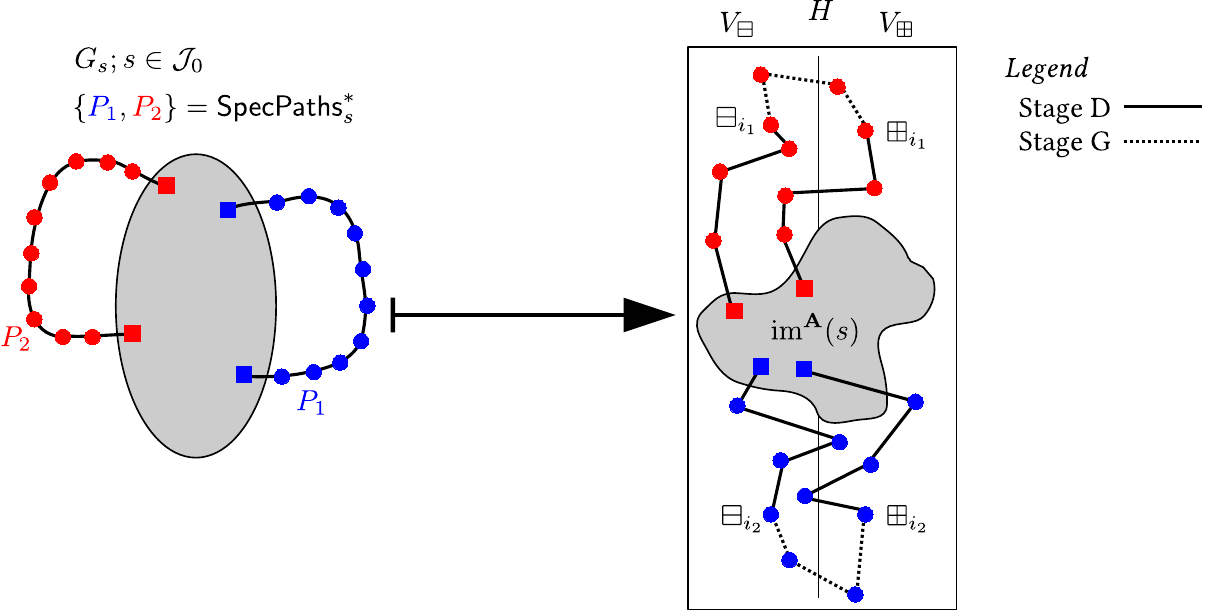}
	\caption{Embedding a graph $G_s$, $s\in \cJ_0$.}
	\label{fig:embedJ0}
\end{figure}
Our packing of $(G_s)_{s\in \cG}$ into~$H$ will be provided in seven stages,
called Stages~A--G. These stages are captured by Lemma~\ref{lem:StageA}
(Stage~A), Lemma~\ref{lem:StageB} (Stage~B), Lemma~\ref{lem:StageCNew}
(Stage~C), Lemma~\ref{lem:StageDNew} (Stage~D), Lemma~\ref{lem:StageENew}
(Stage~E), Lemma~\ref{lem:StageFNew} (Stage~F) and Lemma~\ref{lem:StageG}
(Stage~G). The proofs of these key lemmas are given in
Sections~\ref{sec:StageA}--\ref{sec:StageG}.

Before turning to these lemmas, let us briefly summarise which parts of which graphs are embedded in which stage.
The graphs $(G_s)_{s\in\cG\setminus (\cK\cup\cJ)}$ will be packed entirely in Stage~A. In Stage~A, we shall also pack $(G_s^\spadesuit)_{s\in\cK}$ and $(G_s^\parallel)_{s\in\cJ}$. We shall add the odd-degree vertices $(\SpecLeaves_s)_{s\in\cK}$ to the packing in Stage~B, hence finalising the packing of $(G_s)_{s\in \cK}$. Stage~C is devoted to splitting $V(H)$ into two vertex sets $\Vmin$ and $\Vplus$ of size $\lfloor\frac{n}{2}\rfloor$, and a single vertex $\boxdot$ if $n$ is odd. Further, for each $s\in\cJ$ we embed a subset of the paths $\SpecPaths_s$; the unembedded paths will be called $\SpecPaths^*_s$\index{$\SpecPaths^*_s$}. In doing so, we use all the edges of the form $\boxminus_i\boxplus_i$. If $n$ is odd, we further do away with having to deal with $\boxdot$ in the future by using all the edges incident with it.
From Stage~D to Stage~G we pack the remaining paths from $(G_s)_{s\in\cJ}$. These packings differ for $\cJ_0$ (see Figure~\ref{fig:embedJ0}), $\cJ_1$ and $\cJ_2$. In each of these stages, we shall have some family, say $\cP$, of paths we want to process during that stage. The paths in $\cP$ are subgraphs of $\bigcup_{s\in \cJ}\SpecPaths^*_s$. We emphasise that being a subgraph can mean either that for paths of $\bigcup_{s\in \cJ}\SpecPaths^*_s$ we restrict to some shorter paths, or that we consider subfamilies of $\{\SpecPaths^*_s\}_{s\in \cJ}$ and take the corresponding paths of the full length.
The paths of $\cP$ will be \index{anchor}\emph{anchored}, by which we mean that prior to that stage, for each $P\in\cP$ we have that the embedding of $\leftpath_0(P)$ and $\rightpath_0(P)$ is defined, while the rest of $P$ is unembedded. When $\{\leftpath_0(P),\rightpath_0(P)\}$ is embedded to $\{u,v\}\subset V(H)$, we say that \emph{$P$ is anchored at $u$ and $v$}. In that given stage we shall then embed $\cP$.
The progress of the embedding during the above stages is also outlined in Table~\ref{tab:stages}.

\begin{table}
\tiny{
\bgroup
\renewcommand{\arraystretch}{1.3}
\begin{tabular}{|c||c|c|c|c|c|}
\hline
 \!\!\textbf{Stage}\!\!& $\cG\setminus (\cK\cup \cJ)$ & $\cK$ & $\cJ_0$ & $\cJ_1$ & $\cJ_2$\tabularnewline
 \hline 
 \hline 
\multirow{2}{*}{\textbf{A}} & Packing $G_s$ & Packing $G^\spadesuit_s$ & Packing $G^\parallel_s$ &  Packing $G^\parallel_s$& Packing $G^\parallel_s$ \tabularnewline
 &\ding{51} Packing complete & \ding{56} $\SpecLeaves_s$ missing & \ding{56} 10-paths missing & \ding{56} 10-paths missing & \ding{56} 6-paths missing
\tabularnewline
 \hline 
\multirow{2}{*}{\textbf{B}} & & Packing $\SpecLeaves_s$ & & & \tabularnewline
& & \ding{51} Packing complete &  &  &
\tabularnewline
 \hline 
 \multirow{2}{*}{\textbf{C}}& & & Packing $\le 1$ path &\multicolumn{2}{c|}{Packing $\le n^{0.6}$ paths from $\SpecPaths_s$}\tabularnewline
 & & & from $\SpecPaths_s$ & \multicolumn{2}{c|}{}
\tabularnewline
\hline 
 \multirow{4}{*}{\textbf{D}} & & &Packing $\leftpath_4(P)$ &  &  \tabularnewline
& & & and $\rightpath_4(P)$ & &
\tabularnewline
& & & for $P\in\SpecPaths^*_s$ & &
\tabularnewline
& &  & \ding{56} 2-paths missing &  &
\tabularnewline
\hline 
 \multirow{2}{*}{\textbf{E}} & & &  &  &  Packing $\SpecPaths^*_s$
\tabularnewline
& &  & &  & \ding{51} Packing complete
\tabularnewline
\hline 
 \multirow{2}{*}{\textbf{F}} & &  &  &  Packing $\SpecPaths^*_s$ &
\tabularnewline
& &  & &  \ding{51} Packing complete &
\tabularnewline
\hline 
 \multirow{3}{*}{\textbf{G}} & &  & Packing $\middlepath_3(P)$ &  &
\tabularnewline
 &  & &  for $P\in\SpecPaths^*_s$ & &
\tabularnewline
 &  & &  \ding{51} Packing complete & &
\tabularnewline
\hline 
\end{tabular}
\vspace{1em}
\egroup
}
\caption{Packing of different parts of the graphs $\{G_s\}_{s\in\cG}$ in
  Stage~A-Stage~G depending on the type $s$. The description is slightly
  imprecise. For the purpose of the table, a \emph{missing $\ell$-path} means
  that $\ell$ vertices on a bare path need to be embedded; since these are
  surrounded by 2 anchors, this means $\ell+1$ edges to be packed.}
\label{tab:stages}
\end{table}

We remark that each of the seven stages starts with a quasirandomness assumption, and indeed none of the seven packing stages would be possible without such an assumption. Hence, the outcome of each stage needs to provide the quasirandomness conditions for the next step. These quasirandomness features are more complicated than the one formulated in Definition~\ref{def:quasirandomness} which talks merely about the structure of the host graph. In particular, we need to guarantee quasirandomness properties for
the images (as vertex sets) of the embedded graphs and for the anchors of the paths still to be embedded.

\smallskip

Let us now turn to the stages.
Given embeddings $(\phi_s)_{s\in\cF}$ of graphs $(F_s)_{s\in\cF}$ into a
graph~$H$, the \emph{leftover graph} is the graph on vertex set $V(H)$
containing precisely all edges of~$H$ not used by any $\phi_s$. The \emph{image
in~$H$} of the embedding $\phi_s$ is the set of vertices in $V(H)$ used by
$\phi_s$.

In Stage~A we shall pack the family $(G_s)_{s\in\cG\setminus (\cK\cup\cJ)}\cup (G_s^\spadesuit)_{s\in\cK}\cup (G_s^\parallel)_{s\in\cJ}$ into the graph~$H$ (see~\ref{A:packing}) with leftover graph~$\HStageA$.
As mentioned above, in addition to providing the packing of the family $(G_s^\spadesuit)_{s\in\cK}\cup(G_s^\parallel)_{s\in\cJ}$, we need to obtain certain quasirandomness conditions for~$\HStageA$ that will be used in later stages. 
Here, \ref{A:megaquasirandomness1} and~\ref{A:typicality} guarantee that the individual images of the graphs $(G_s^\spadesuit)_{s\in\cK}\cup(G_s^\parallel)_{s\in\cJ}$ are spread uniformly, \ref{A:parleavesSpread} and~\ref{A:SpecLeavesCover} guarantee that neighbours of the reserved odd-degree vertices are embedded uniformly,
while~\ref{A:vertex-covering}--\ref{A:PathsNotSquashed} assert that
anchors of paths that remain unembedded are spread uniformly.

Before stating the lemma, we also note that we will abuse notation slightly for the last three properties~\ref{A:PathsWithOneVertex}--\ref{A:PathsNotSquashed}: whenever a statement considers a path $P\in\bigcup_{s'\in \cJ}\SpecPaths_{s'}$, we let $s$ denote the unique index in $\cJ$ such that $P\in \SpecPaths_s$.

\begin{lemma}[Stage~A, bulk embedding]\label{lem:StageA}\index{bulk embedding!statement of Lemma~\ref{lem:StageA}}\index{Stage A!statement of Lemma~\ref{lem:StageA}}
  Given graphs $(G_s)_{s\in\cG}$ and their subgraphs as described
  above, and a graph $H$ which is $(\iniquasi,4D+3)$-quasirandom, there exist
  maps $(\phi^\mathbf{A}_s)_{s\in\cG}$ with the following properties.
  \begin{enumerate}[label=\itmrom{A},leftmargin=*]
  \item\label{A:packing} $(\phi^\mathbf{A}_s)_{s\in\cG}$ packs
    $(G_s)_{s\in\cG\setminus (\cK\cup\cJ)}\cup
    (G_s^\spadesuit)_{s\in\cK}\cup
    (G_s^\parallel)_{s\in\cJ}$ into~$H$ with leftover graph
    $\HStageA$\index{$\HStageA$} of density $\dStageA$\index{$\dStageA$} and with images $(\imA(s))_{s\in\cG}$ in~$H$.
    
  \item\label{A:megaquasirandomness1}
    For each $S\subset V(H)$, $|S|\le 2D$ and each $T\subset \cJ\cup \cK$, $|T|\le D$ we have
    \[
    \Big|\NBH_{\HStageA}(S) \setminus \bigcup_{s\in T}\imA(s)\Big|
    =
    (1\pm \gamcore)\dStageA^{|S|} n \prod_{s\in T}\Big(1-\frac{|\imA(s)|}{n}\Big)
    \;.
    \]

  \item \label{A:typicality} 
  For every pair of disjoint sets $S_1, S_2\subseteq V(H)$ with $|S_1|,|S_2|\le D$ we have with
  $\cX:=\{s\in\cJ:\imA(s)\cap S_1=\emptyset,  S_2\subseteq \imA(s)\}$ that
  \begin{align*}
    |\cJ_0\cap\cX|
    &=(1\pm \gamcore)(\deltnonspanning+10\sigmKJ)^{|S_1|}(1-\deltnonspanning-10\sigmKJ)^{|S_2|}|\cJ_0|\;,\\
    |\cJ_1\cap\cX|
    &=(1\pm \gamcore)(\deltnonspanning+10\sigmJjedna)^{|S_1|}(1-\deltnonspanning-10\sigmJjedna)^{|S_2|}|\cJ_1|\;,\\
    |\cJ_2\cap\cX|
    &=(1\pm \gamcore)(\deltnonspanning+6\sigmKJ)^{|S_1|}(1-\deltnonspanning-6\sigmKJ)^{|S_2|}|\cJ_2|\;.
  \end{align*}

\item\label{A:parleavesSpread} 
For every $v\in V(H)$ we have
$\sum_{s\in\cK}\deg_{G_s}\big((\phi_s^\mathbf{A})^{-1}(v);\SpecLeaves_s\big)\le \tfrac{20c n}{\log n}$,
where if $v\notin\imA(s)$ we count these degrees as zero.

\item\label{A:SpecLeavesCover} 
For every $v\in V(H)$ we have
\[\sum_{s\in \cK; v\notin \imA(s)}\ \sum_{x\in \SpecLeaves_s}\ \prod_{y\in \NBH_{G_s}(x)}\ONE_{v\phi^\mathbf{A}_s(y)\in E(\HStageA)}
\ge 4^{-200 D^{2}\dStageA^{-D}\deltnonspanning^{-1}}n\;.\]

  \item \label{A:vertex-covering}
    For each $v\in V(H)$ and for every $s\in \cJ$ such that $v\notin \imA(s)$  we have that
    \begin{multline*}
      \hspace{5em}\Big|\Big\{x\in 
      \bigcup_{P\in \SpecPaths_s}\{\leftpath_0(P), \rightpath_0(P)\}:
      v\in \NBH_{\HStageA}
      \big(\phi^\mathbf{A}_s(x)\big)\Big\}\Big| \\
      = (1\pm \gamcore)\dStageA\cdot 2|\SpecPaths_s|\;.\qquad
    \end{multline*}
    
\item\label{A:PathsWithOneVertex}
For every $v\in V(H)$ the following hold. The number of $P\in\bigcup_{s\in \cJ_1}\SpecPaths_s$ for which $v\in \{\phi^\mathbf{A}_s(\leftpath_0(P)),\phi^\mathbf{A}_s(\rightpath_0(P))\}$ is equal to $2(1\pm \frac{\gamcore}2) \sigmJjedna  |\cJ_1|$. 
For every $\cJ^*\in \{\cJ_0,\cJ_2\}$ we have that the number of $P\in\bigcup_{s\in \cJ^*}\SpecPaths_s$ for which $v\in \{\phi^\mathbf{A}_s(\leftpath_0(P)),\phi^\mathbf{A}_s(\rightpath_0(P))\}$ is equal to $2(1\pm \frac{\gamcore}2) \sigmKJ |\cJ^*|$.

\item\label{A:PathsOneVtxIm}
For every pair of distinct vertices $u,v\in V(H)$ and for every $\cJ^*\in \{\cJ_0,\cJ_1,\cJ_2\}$ the following hold. The number of $P\in\bigcup_{s\in \cJ^*}\SpecPaths_s$ for which $u\notin\imA(s)$ and $v\in \{\phi^\mathbf{A}_s(\leftpath_0(P)),\phi^\mathbf{A}_s(\rightpath_0(P))\}$ is equal to 
\begin{enumerate}[label=\abc]
\item \label{enu:PathsOneVtxIm0}$2(1\pm \frac{\gamcore}2) \sigmKJ(\deltnonspanning+10\sigmKJ) |\cJ_0|$, if $\cJ^*=\cJ_0$,
\item \label{enu:PathsOneVtxIm1}$2(1\pm \frac{\gamcore}2) \sigmJjedna(\deltnonspanning+10\sigmJjedna)|\cJ_1|$, if $\cJ^*=\cJ_1$,
\item \label{enu:PathsOneVtxIm2}$2(1\pm \frac{\gamcore}2) \sigmKJ(\deltnonspanning+6\sigmKJ) |\cJ_2|$, if $\cJ^*=\cJ_2$.
\end{enumerate}
If $u,u',v$ are distinct vertices of $V(H)$, the number of $P\in\bigcup_{s\in\cJ_0}\SpecPaths_s$ for which $u,u'\notin\imA(s)$ and $v\in\{\phi^\mathbf{A}_s(\leftpath_0(P)),\phi^\mathbf{A}_s(\rightpath_0(P))\}$ is equal to
\begin{enumerate}[label=\abc,start=4]
\item \label{enu:PathsTwoVtxIm0}$2(1\pm \frac{\gamcore}2) \sigmKJ(\deltnonspanning+10\sigmKJ)^2 |\cJ_0|$.
\end{enumerate}

\item\label{A:PathsNotSquashed}
For every two distinct vertices $v_1,v_2\in V(H)$ we have that the number of $P\in\bigcup_{s\in \cJ}\SpecPaths_s$ for which $\{v_1,v_2\}= \{\phi^\mathbf{A}_s(\leftpath_0(P)),\phi^\mathbf{A}_s(\rightpath_0(P))\}$ is at most $n^{0.3}$.
\end{enumerate}
\end{lemma}
The proof of this lemma, which we give in Section~\ref{sec:StageA}, relies on the analysis of the randomised algorithm \PackingProcess, which was introduced in~\cite{DegPack} and was also one of the key components in~\cite{ABCT:PackingManyLeaves}. We use the main technical result of~\cite{DegPack} as a black box to pack the graphs $(G_s)_{s\in\cG\setminus(\cK\cup\cJ)}$, after which the remaining edges still form a quasirandom graph. Unfortunately, we cannot use the results of~\cite{DegPack} as a black box to pack the remaining graphs $(G_s^\spadesuit)_{s\in\cK}\cup (G_s^\parallel)_{s\in\cJ}$, all of which have at most $(1-\deltnonspanning)n$ vertices. The reason is that, in addition to providing the packing of the family $(G_s^\spadesuit)_{s\in\cK}\cup (G_s^\parallel)_{s\in\cJ}$ we also need to ensure various quasirandomness properties listed in~\ref{A:megaquasirandomness1}--\ref{A:PathsNotSquashed}.
Note that the density $\dStageA$ of the leftover graph $\HStageA$ of this stage satisfies
\begin{align}\label{eq:dStageA}
\dStageA&=\binom{n}{2}^{-1}\cdot \Big(e(H)-\sum_{s\in\cG\setminus(\cK\cup\cJ)}e(G_s)-\sum_{s\in\cK}e(G_s^\spadesuit)-\sum_{s\in\cJ}e(G_s^\parallel)\Big)\\
\nonumber
\JUSTIFY{Sec~\ref{ssec:IneqEqual},\ref{ssec:definingSubgraphs}}&=\binom{n}{2}^{-1}\cdot \Big(\DParity \cdot\big|\bigcup_{s\in\cK}\SpecLeaves_s\big|+11\sigmKJ n|\cJ_0|+11\sigmJjedna n|\cJ_1|+7\sigmKJ n|\cJ_2|\Big)\\
\label{eq:dStageAapprox}
\JUSTIFY{Sec~\ref{ssec:definingSubgraphs},\eqref{eq:CONSTANTS}}&=(14\pm 0.1) \sigmKJ^2\;.
\end{align}

We now turn to Stage~B. By suitably embedding the reserved vertices $\bigcup_{s\in\cK}\SpecLeaves_s$ of odd degree, we shall ensure that the vertex degrees in the leftover graph of this stage have parities suitable for the remaining stages of our packing process. Let us provide more details on this now. Suppose that $v\in V(H)$ is an arbitrary vertex. After Stage~A, there are $\deg_{\HStageA}(v)$ edges incident to~$v$ which are available for the remaining stages. We know that once the packing is completed, all these edges will be used (by our assumptions on equalities from Section~\ref{ssec:IneqEqual}). These host graph edges will be used to accommodate guest graph edges of the following types.

\begin{enumerate}[label={[et\arabic*]}]
\item\label{et:1} Edges of the type $xy\in E(G_s)$, where $y\in\SpecLeaves_s$, $x\in \NBH_{G_s}(y)$, $\phi^\mathbf{A}_s(x)=v$ and $s\in\cK$. The number of such edges is equal to 
	\[
	\index{$\OddOut(v)$}
	\OddOut(v):=\sum_{s\in \cK}\deg_{G_s}\big((\phi_s^{\mathbf{A}})^{-1}(v);\SpecLeaves_s\big)\;.\]

\item\label{et:2} Edges of the type $xy\in E(G_s)$, where $y\in\SpecLeaves_s$, $x\in \NBH_{G_s}(y)$, $s\in\cK$, and where $y$ will be mapped to $v$.

\item\label{et:3} Edges $\leftpath_1(P)$ and $\rightpath_1(P)$ of the paths $P\in\SpecPaths_s$ (for $s\in\cJ$), for those paths $P$ for which $\phi^\mathbf{A}_s(\leftpath_0(P))=v$ or $\phi^\mathbf{A}_s(\rightpath_0(P))=v$.
The number of such edges is equal to 
\begin{align}
\label{eq:defPathTerm}
\index{$\PathTerm(v)$}
\PathTerm(v):=\sum_{s\in \cJ}
	\left|\{v\}\cap\{\phi^\mathbf{A}_s(\rightpath_0(P)),\phi^\mathbf{A}_s(\leftpath_0(P)):P\in\SpecPaths_s\}\right|
	\;.
\end{align}
	
\item\label{et:4} Edges $xy$ of the paths $P\in\SpecPaths_s$ (for $s\in\cJ$) with $x$ or~$y$ mapped to~$v$, for which $P$ is not anchored at~$v$, i.e., $(\phi^\mathbf{A}_s)^{-1}(v)\cap \{\leftpath_0(P), \rightpath_0(P)\}=\emptyset$.
\end{enumerate}

Each path counted in~\ref{et:4} will actually use~2 edges at $v$. In particular, the total number of edges counted in~\ref{et:4} will be even. Since edges counted in~\ref{et:1} and~\ref{et:3} are already determined after Stage~A, we need to adjust the number of edges counted in~\ref{et:2} so that it has the same parity as 
\[\index{$\Parity(v)$}
\Parity(v):=\deg_{\HStageA}(v)-\OddOut(v)-\PathTerm(v)\mod 2\;.\] In other words, we need to embed $\bigcup_{s\in\cK}\SpecLeaves_s$ in such a way that for each vertex $v$, the number of vertices from $\bigcup_{s\in\cK}\SpecLeaves_s$ embedded to $v$ has the same parity as $\Parity(v)$ (see~\ref{B:parity}).

The number of edges we are newly embedding in Stage~B is at most
$\DParity\big|\bigcup_{s\in\cK}\SpecLeaves_s\big|$ by our assumptions, which is
linear in $n$. So the density~$\dStageB$ of the leftover graph~$\HStageB$ after
the packing of $(G_s)_{s\in\cG\setminus
  (\cK\cup\cJ)}\cup(G_s^\parallel)_{s\in\cJ}\cup(G_s)_{s\in\cK}$ is very close
to $\dStageA$.  Because we are embedding so few edges in Stage~B, we do not need
to use a randomised algorithm for this stage but can rely on matching-type
arguments. We shall show (see~\ref{B:degrees}) that we can perform this packing
stage so that each host graph vertex loses only a tiny fraction of incident
edges. It follows that the leftover graph $\HStageB$ automatically inherits the
quasirandomness properties concerning the uniform spread of images and anchors
from $\HStageA$. Accordingly, the following lemma does not list most of these
quasirandomness properties again.
We prove this lemma in Section~\ref{sec:NEWStageB}.

\begin{lemma}[Stage~B, parity correction]\label{lem:StageB}
\index{parity correction!statement of Lemma~\ref{lem:StageB}}\index{Stage B!statement of Lemma~\ref{lem:StageB}}
  Given partial embeddings
  $(\phi^\mathbf{A}_s)_{s\in\cG}$ as described in Lemma~\ref{lem:StageA}, 
  there exist maps $(\phi^\mathbf{B}_s)_{s\in\cG}$ extending the $(\phi^\mathbf{A}_s)_{s\in\cG}$ with the following properties.
\begin{enumerate}[label=\itmrom{B}]
\item\label{B:one}
  For each $s\in \cK$, the map $\phi^\mathbf{B}_s$ is an embedding of the entire graph $G_s$ into $H$.
  For each $s\in \cG\setminus \cK$ we have $\phi^\mathbf{B}_s=\phi^\mathbf{A}_s$.
\item\label{B:two} $(\phi^\mathbf{B}_s)_{s\in\cG}$ is a packing
  of $(G_s)_{s\in\cG\setminus \cJ}\cup
  (G_s^\parallel)_{s\in\cJ}$ into the
  graph $H$ with leftover graph $\HStageB$\index{$\HStageB$} of density
  $\dStageB$\index{$\dStageB$} and images $\big(\imB(s)\big)_{s\in\cG}$\index{$\imB(s)$}.
\item\label{B:parity}
For each $v\in V(H)$ we have
\[
\deg_{\HStageB}(v)\equiv\PathTerm(v)
\mod 2\;.\]

\item\label{B:degrees}
For each $v\in V(H)$ we have
$\deg_{\HStageB}(v)\ge \deg_{\HStageA}(v)-8Dcn$.

\item\label{B:megaquasirandomness1}
For each $S\subset V(H)$, $|S|\le 2D$ and each $T\subset \cJ$, $|T|\le D$ we have
\[
\Big|\NBH_{\HStageB}(S) \setminus \bigcup_{s\in T}\imB(s)\Big|
=
(1\pm 2\gamcore)\dStageB^{|S|} n \prod_{s\in T}\left(1-\frac{|\imB(s)|}{n}\right)
\;.
\]
\end{enumerate}
\end{lemma}

In Stage~C, we randomly split $V(H)$ into two vertex sets $\Vmin$ and $\Vplus$ of size $\lfloor\frac{n}{2}\rfloor$, and a special vertex $\boxdot$ if $n$ is odd. We randomly label the vertices of $\Vmin$ as $\{\boxminus_i:i\in [\lfloor\frac{n}2\rfloor]\}$ and of $\Vplus$ as $\{\boxplus_i:i\in [\lfloor\frac{n}2\rfloor]\}$. We then embed a few of the paths in $\bigcup_{s\in\cJ}\SpecPaths_s$ so that the following three things happen, of which the last two apply only when $n$ is odd. Firstly, all the edges of $H$ which are of the form $\boxminus_i\boxplus_i$ are used. Secondly, all paths which are anchored at $\boxdot$ are embedded. Thirdly, all other edges leaving $\boxdot$ are used.
After this, if $n$ is odd then all edges incident to the vertex $\boxdot$ are used up. So, regardless of the parity of $n$, in the remaining Stages~D--G we effectively work with a host graph on an even number of vertices.

The key assumption of Theorem~\ref{thm:maintechnical}, which we are proving, is the quasirandomness of our host graph, formalised in Definition~\ref{def:quasirandomness}. In~\ref{A:megaquasirandomness1}--\ref{A:PathsNotSquashed} we saw more complicated quasirandomness conditions which dealt for example with mutual positions of the images of different graphs $G_s$. In the following stages we will encapsulate our running quasirandomness properties under a common \emph{quasirandom setup}, which we define now. That is, after each of the next four stages, we will argue that we have maintained the quasirandom setup (with quasirandomness parameters becoming a bit worse in each stage) and have gained certain additional good properties, which permit the final step to complete a perfect packing. In what follows, the reader should think of $U_s$ as the vertices in $V(H)$ used in the current partial embedding of $G_s$; of $F_s$ as the path-forest in $G_s$ which remains to be embedded, with $A_s$ the set of vertices to which the paths of $F_s$ are anchored and $\phi_s$ the anchoring; and of $I_s$ as the indices of terminal pairs
which have not been used in the partial embedding of~$G_s$ (after Stage~C) or which currently serve as anchors (after Stage~D).

\begin{definition}[Quasirandom setup]\label{def:quasisetup}\index{$(\gamma,L,d_1,d_2)$-quasirandom setup}
  Let $H$ be a graph on an even number of vertices partitioned into $\Vmin=\{\boxminus_i:i\in[v(H)/2]\}$ and $\Vplus=\{\boxplus_i:i\in[v(H)/2]\}$. Let $\cJ=\cJ_0\dcup\cJ_1\dcup\cJ_2$ be a set of indices. For each $s\in\cJ$, let $F_s$ be a path-forest whose components all have at least~4 vertices. Let $U_s$ be a subset of $V(H)$, and let $A_s\subseteq U_s$. Let $\phi_s$ be a bijection from the leaves of $F_s$ to~$A_s$. For each $s\in\cJ$ and $a\in\{\boxplus,\boxminus\}$, let $A_s^a$ denote the set of vertices $u\in A_s$ such that the path in~$F_s$ with endvertex $\phi_s^{-1}(u)$ has its other endvertex embedded by $\phi_s$ to $V_a$.
 Finally, for each $s\in\cJ_0$ let the set $I_s$ be a subset of $[v(H)/2]$. This is a \emph{$(\gamma,L,d_1,d_2)$-quasirandom setup} if the following conditions are satisfied.
 \begin{enumerate}[label=\itmarab{Quasi}]
  \item\label{quasi:indexquasi} The triple $(H,\Vmin,\Vplus)$ is $(L,\gamma,d_1,d_2)$-index-quasirandom with respect to the sets $(U_s)_{s\in\cJ}$ and $(I_s)_{s\in\cJ_0}$.
  \item\label{quasi:anchorsets} For each $a,b\in\{\boxminus,\boxplus\}$ and each $s\in\cJ_1\cup\cJ_2$ we have $|A^b_s\cap V_a|=(1\pm\gamma)\tfrac14|A_s|$. For each $s\in\cJ_0$ we have $|A_s\cap V_a|=(1\pm\gamma)\tfrac12|A_s|$.
  \item\label{quasi:PathsWithOneVertex}	For every $v\in V(H)$,  $a\in\{\boxminus,\boxplus\}$ and $i\in\{1,2\}$ we have that the number of $s\in\cJ_i$ such that $F_s$ has a path with leaves $x$ and $y$ such that $\phi_s(x)=v$ and $\phi_s(y)\in V_a$ is equal to $(1\pm \gamma)\tfrac{1}{2n}\sum_{s\in\cJ_i}|A_s|$.
  \item\label{quasi:termVtxNbs} For each $s\in\cJ_0$, each $u\in \Vmin\setminus U_s$ and each $v\in\Vplus\setminus U_s$ we have
   \begin{align*}
    \big|\{\boxminus_i\in \NBH_H(u):i\in I_s\}\big|,\big|\{\boxplus_i\in \NBH_H(v):i\in I_s\}\big|&=(1\pm\gamma)d_1|I_s|\quad     \text{, and}\\
 \big|\{\boxminus_i\in \NBH_H(v):i\in I_s\}\big|,\big|\{\boxplus_i\in \NBH_H(u):i\in I_s\}\big|&=(1\pm\gamma)d_2|I_s|\,.
   \end{align*}
  \item\label{quasi:distributionOfAnchors} For each $s\in\cJ_0\cup\cJ_1\cup\cJ_2$, each $a\in\{\boxplus,\boxminus\}$, each $u\in \Vmin\setminus U_s$ and each $v\in\Vplus\setminus U_s$ we have
   \begin{align}
    \label{eq:quasi:distAnchorsa1}\big|\NBH_H(u)\cap A^a_s\cap \Vmin\big|&=(1\pm\gamma)d_1|A^a_s\cap\Vmin|\,,\\
    \label{eq:quasi:distAnchorsa2}\big|\NBH_H(v)\cap A^a_s \cap \Vmin\big|&=(1\pm\gamma)d_2|A^a_s\cap\Vmin|\,,\\
    \label{eq:quasi:distAnchorsa3}\big|\NBH_H(v)\cap A^a_s\cap\Vplus\big|&=(1\pm\gamma)d_1|A^a_s\cap\Vplus|\,\text{, and}\\
    \label{eq:quasi:distAnchorsa4}\quad \big|\NBH_H(u)\cap A^a_s\cap \Vplus\big|&=(1\pm\gamma)d_2|A^a_s\cap\Vplus|\,.
   \end{align}
  \item\label{quasi:imagecaps} Given disjoint sets $S_1,S_2\subset V(H)$ with $|S_1|,|S_2|\le L$, for each $i=1,2$ we have
   \[\big|\{s\in\cJ_i\,:\,S_1\cap U_s=\emptyset, S_2\subset U_s\}\big|=(1\pm\gamma)\sum_{s\in\cJ_i}\frac{(n-|U_s|)^{|S_1|}|U_s|^{|S_2|}}{n^{|S_1|+|S_2|}}\,.\]
   Given additionally $T\subseteq [v(H)/2]$ with $|T|\le L$ such that $\{\boxminus_i,\boxplus_i:i\in T\}$ is disjoint from $S_1\cup S_2$, suppose that there is no $i$ such that $\{\boxminus_i,\boxplus_i\}\subset  S_1\cup S_2$. Then we have
   \[\big|\{s\in\cJ_0\,:\,S_1\cap U_s=\emptyset, S_2\subset U_s,T\subset I_s\}\big|=(1\pm\gamma)\sum_{s\in\cJ_0}\frac{(n-|U_s|)^{|S_1|}|U_s|^{|S_2|}(2|I_s|)^{|T|}}{n^{|S_1|+|S_2|+|T|}}\,.\]
  \item\label{quasi:NumAnchors} For each $i \in \{0,1,2\}$ and each $v\in V(H)$, the number of $s\in\cJ_i$ such that $v\in A_s$ is $(1\pm\gamma)\sum_{s\in\cJ_i}\tfrac{|A_s|}{n}$.
  \item\label{quasi:NumAnchorsNoIm} For each $i\in \{1,2\}$ and each $u,v\in V(H)$ with $u\neq v$, the number of $s\in\cJ_i$ such that $u\in A_s$ and $v\notin U_s$ is $(1\pm\gamma)\sum_{s\in\cJ_i}\frac{|A_s|(n-|U_s|)}{n^2}$. If in addition there is no $j$ such that $\{u,v\}=\{\boxminus_j,\boxplus_j\}$ then the number of $s\in\cJ_0$ such that $u\in A_s$ and $v\notin U_s$ is $(1\pm\gamma)\sum_{s\in\cJ_0}\tfrac{|A_s|(n-|U_s|)}{n^2}$.
 \end{enumerate}
\end{definition}

We remark that when we use this definition, in all the summations it is the value of the sum, and not merely each individual summand, that changes little if the summands are replaced by a single typical value; for example, even though $|A_s|$ may vary somewhat as $s$ ranges over any one $\cJ_i$, replacing $|A_s|$ throughout the sum $\sum_{s\in\cJ_i}|A_s|$ by a single representative value changes the sum by only a small relative amount.
We further remark that for formal reasons when we claim a quasirandom setup in the following lemmas, we will not always use our index sets $\cJ_0,\cJ_1,\cJ_2$ from Section~\ref{ssec:definingSubgraphs} as
index sets $\cJ_0,\cJ_1,\cJ_2$ in the quasirandom setup. More precisely, when we have already completely embedded each~$G_s$ with~$s\in\cJ_i$ we will use $\emptyset$ instead of $\cJ_i$ in the corresponding quasirandom setup. We shall always make this explicit.

The following lemma now encapsulates what we do in Stage~C (which was informally described before Definition~\ref{def:quasisetup}). We prove this lemma in Section~\ref{sec:StageC}.

\begin{lemma}[Stage~C, partite reduction]\label{lem:StageCNew}
	\index{partite reduction!statement of Lemma~\ref{lem:StageCNew}}\index{Stage C!statement of Lemma~\ref{lem:StageCNew}}
  Given partial embeddings $(\phi^\mathbf{B}_s)_{s\in\cG}$ as described in Lemma~\ref{lem:StageB} extending partial embeddings
  $(\phi^\mathbf{A}_s)_{s\in\cG}$ as described in Lemma~\ref{lem:StageA}, we can construct the following.
 
  If $n$ is even we find a labelling $V(H)=\{\boxminus_i,\boxplus_i\}_{i\in [n/2]}$, and if $n$ is odd we find a labelling $V(H)=\{\boxdot\}\cup\{\boxminus_i,\boxplus_i\}_{i\in [(n-1)/2]}$\index{$\boxdot$, $\boxminus_i$, $\boxplus_i$}. Set $\Vmin:=\{\boxminus_i:i\in [\lfloor\frac{n}2\rfloor]\}$ and $\Vplus:=\{\boxplus_i:i\in [\lfloor\frac{n}2\rfloor]\}$.\index{$\Vmin$}\index{$\Vplus$}
  Further, we find maps $(\phi^\mathbf{C}_s)_{s\in\cG}$ extending the
  $(\phi^\mathbf{B}_s)_{s\in\cG}$ with the following properties.
  \begin{enumerate}[label=\itmrom{C}]
  \item\label{C:one} $(\phi^\mathbf{C}_s)_{s\in\cG}$ is a packing of subgraphs
    of $\left(G_s\right)_{s\in\cG}$ into~$H$ with leftover graph
    $\HStageC$\index{$\HStageC$} on $\Vmin\dcup\Vplus$ and with
    density~$\dStageC$.\index{$\dStageC$}
  \item\label{C:two} For each $s\in\cJ$, the map $\phi^\mathbf{C}_s$ embeds
    $G_s^\parallel$ as well as some entire paths from $\SpecPaths_s$ (and no
    other vertices).  More precisely, if $s\in\cJ_0$, exactly one path from
    $\SpecPaths_s$ is embedded when $n$ is odd, and no paths when $n$
    is even.  If $s\in\cJ_1\cup\cJ_2$, at most $n^{0.6}$ paths from
    $\SpecPaths_s$ are embedded.
    We denote the remaining unembedded paths in $G_s$ by
    $\SpecPaths^*_s$\index{$\SpecPaths^*_s$} for each $s\in\cJ$. We write
    $\imC(s)=\im(\phi^\mathbf{C}_s)$.\index{$\imC(s)$}
  \item\label{C:three} Every edge of $H$ of the form $\boxminus_i\boxplus_i$
    with $i\in[\lfloor n/2\rfloor]$ is in the image of some
    $\phi^\mathbf{C}_s$. If $n$ is odd, every edge of $H$ incident with
    $\boxdot$ is in the image of some $\phi^\mathbf{C}_s$.
  \item\label{C:four} With the given $\cJ_0,\cJ_1,\cJ_2$, for $s\in\cJ_0\cup\cJ_1\cup\cJ_2$ let $F^\mathbf{C}_s=\SpecPaths^*_s$, let
    $U^\mathbf{C}_s=\imC(s)$, let $\phi_s$ be the restriction of
    $\phi^\mathbf{C}_s$ to the leaves of $F^\mathbf{C}_s$, and let
    $A^\mathbf{C}_s$ be the image of $\phi_s$.  For each $s\in\cJ_0$ let
    $I^\mathbf{C}_s=\{i\in[\lfloor\frac{n}2\rfloor]\,:\,\boxminus_i,\boxplus_i\notin
    U^\mathbf{C}_s\}$.
    Then for the graph $\HStageC$ we have a
    $(100\LC\gamcore,\LC,\dStageC,\dStageC)$-quasirandom setup, and    
    for each $s\in\cJ$ we have $A^\mathbf{C}_s\subset
    V(\HStageC)$, that is, all unembedded paths have their anchors in
    $\Vmin\dcup\Vplus$. 
\item\label{C:five}
For each $s\in\cJ_0$ we have \[|I^\mathbf{C}_s|=(1\pm
10\LC\gamcore)\tfrac{n}{2}\big(1-\tfrac{|U^\mathbf{C}_s|}{n}\big)^2\,.\]
For
each $u\in V(\HStageC)$ and $j\in[\lfloor\frac{n}2\rfloor]$ such that
$u\notin\{\boxminus_j,\boxplus_j\}$ we have
\[
\big|\big\{s\in\cJ_0:u\in A^\mathbf{C}_s \text{ and } j\in I^\mathbf{C}_s\big\}\big|
=(1\pm100\LC\gamcore)\sum_{s\in\cJ_0}\tfrac{2|A^\mathbf{C}_s||I^\mathbf{C}_s|}{n^2}\,.
\]
  \end{enumerate}
\end{lemma}

In Stage~D we shall extend our packing by outer parts of paths of $\bigcup_{s\in \cJ_0}\SpecPaths^*_s$ (recall Figure~\ref{fig:embedJ0}). More precisely, for each $s\in \cJ_0$ we shall extend the packing to  \index{$G_s^{\boxminus\boxplus}$}\[G_s^{\boxminus\boxplus}:=G_s-\bigcup_{P\in\SpecPaths^*_s}V(\middlepath_1(P))\;.\] 
That is, for each $P\in\SpecPaths^*_s$, the middle pair of vertices of $P$ is removed from $G_s$. This corresponds to a removal of 3~edges: the edge between the middle pair, as well as 2~edges surrounding it.

Observe that $\SpecPaths_s^*$ has the same size for every $s\in \cJ_0$. Consequently, we define
\begin{equation}\label{eq:sizeSpecPathStar}\index{$\sigmKJ^*$}
\sigmKJ^*=\frac{|\SpecPaths_s^*|}{n}\;.
\end{equation}
Recalling that $\SpecPaths_s=\BasicPaths_s$ (since $s\in \cJ_0$; see Section~\ref{ssec:definingSubgraphs}) has size $\sigmKJ n$ by~\eqref{eq:sizeBasicPaths}, we see that 
\begin{equation}\label{eq:sigmKJsigmKJ}
\mbox{$\sigmKJ^*=\sigmKJ$, if $n$ is even, and $\sigmKJ^*=\sigmKJ-\frac1n$, if $n$ is odd.}
\end{equation}
A crucial property guaranteed in the lemma for Stage~D (Lemma~\ref{lem:StageDNew}) is that the unembedded edges of $\middlepath_3(P)$ are always anchored at a \emph{terminal pair}\index{terminal pair}, that is, a pair of vertices of the form $\{\boxminus_i,\boxplus_i\}$.

We shall call the graph consisting of edges that remain after this step~$\HStageD$.
Recall that, after Stage~D, only the following paths remain to be embedded: the paths in $\SpecPaths_s^*$ for $s\in\cJ_1\cup\cJ_2$, 
and the paths 
$\middlepath_3(P)$ with $P\in\SpecPaths^*_s$ for $s\in\cJ_0$. We have $|\cJ_0|=\sigmJnula n$, $|\cJ_1|=\sigmJjedna n$ by~\eqref{eq:sizeJ0J1} and $|\cJ_2|=(\sigmKJ-\sigmJnula-\sigmJjedna)n$ by~\eqref{eq:sizeJ2}. Further, for $s\in\cJ_0$ there are $\sigmKJ n\pm 1$ paths in $\SpecPaths^*_s$, for each of which only~$3$ edges remain to be embedded after Stage~D; for $s\in\cJ_1$ there are $\sigmJjedna n \pm n^{0.6}$ paths in $\SpecPaths_s^*$ each with~$11$ edges, and for $s\in\cJ_2$ there are $\sigmKJ n\pm n^{0.6}$ paths in $\SpecPaths_s^*$ each with~$7$ edges.
Therefore,
\[
e(\HStageD)=3\sigmJnula n(\sigmKJ n\pm 1)+11\sigmJjedna n(\sigmJjedna n \pm n^{0.6})+7(\sigmKJ-\sigmJnula-\sigmJjedna)n(\sigmKJ n\pm n^{0.6})\,.
\]
Hence the density of $\HStageD$ will satisfy
\begin{align}\label{eq:dStageD}\index{$\dStageD$}
\dStageD:= \frac{e(\HStageD)}{\binom{|\Vmin\cup\Vplus|}{2}}
=(1\pm \iniquasi)\left(14(\sigmKJ-\sigmJnula-\sigmJjedna)\sigmKJ+6\sigmKJ\sigmJnula+22\sigmJjedna^2\right)\;.
\end{align}
The following lemma is proved in Section~\ref{sec:StageD}. For~\ref{D:two}
recall that each path in $\SpecPaths^*_s$ with $s\in\cJ_0$ has~$12$ vertices.

\begin{lemma}[Stage~D, connecting to terminal pairs]\label{lem:StageDNew}
\index{connecting to terminal pairs!statement of Lemma~\ref{lem:StageDNew}}\index{Stage D!statement of Lemma~\ref{lem:StageDNew}}
\addtocounter{LemmaDCounter}{\value{theorem}}
  Given partial embeddings $(\phi^\mathbf{C}_s)_{s\in\cG}$ as described in
  Lemma~\ref{lem:StageCNew}, there exist maps $(\phi^\mathbf{D}_s)_{s\in\cG}$
  extending the $(\phi^\mathbf{C}_s)_{s\in\cG}$ with the following
  properties.
  \begin{enumerate}[label=\itmrom{D}]
  \item\label{D:one} $(\phi^\mathbf{D}_s)_{s\in\cG}$ is a packing of subgraphs
    of $\left(G_s\right)_{s\in\cG}$ into~$H$ with leftover graph
    $\HStageD$\index{$\HStageD$} on $\Vmin\dcup\Vplus$ and with density
    $\dStageD$\index{$\dStageD$}.
  \item\label{D:two} For each $s\in\cJ\setminus\cJ_0$ we have
    $\phi^\mathbf{D}_s=\phi^\mathbf{C}_s$.
    For each $s\in\cJ_0$ and each $P\in\SpecPaths^*_s$, the map
    $\phi^\mathbf{D}_s$ embeds all vertices of $P$ except for the two vertices in
    $\middlepath_1(P)$. Furthermore, the fifth and
    eighth vertices
    of $P$ are embedded to $\{\boxminus_i,\boxplus_i\}$ for some
    $i\in[\lfloor\frac{n}2\rfloor]$ (not necessarily in this order).
  \item\label{D:three} With the given $\cJ_0,\cJ_1,\cJ_2$, for each $s\in\cJ_1\cup\cJ_2$ let $F_s=\SpecPaths^*_s$,
    let $U_s=\im(\phi^\mathbf{D}_s)$, let~$\phi_s$ be the restriction of
    $\phi^\mathbf{D}_s$ to the leaves of $F_s$, and
    let~$A_s$ be the image of~$\phi_s$. For each $s\in\cJ_0$ let~$F_s$ be
    \index{$\SpecShortPaths_s$}~\[\SpecShortPaths_s:=\{\middlepath_3(P):P\in\SpecPaths^*_s\}\,,\]
    let $U_s=\im(\phi^\mathbf{D}_s)$, let $\phi_s$ be the restriction of
    $\phi^\mathbf{D}_s$ to the leaves of $F_s$, let $A_s$ be the image of
    $\phi_s$, and let
    $I_s=\{i\in[\lfloor\frac{n}2\rfloor]\,:\,\boxminus_i,\boxplus_i\in
    A_s\}$. Then we have a $(\gamAnchor,\LD,\dStageD,\dStageD)$-quasirandom
    setup.
  \end{enumerate}  
\end{lemma}

In Stage~E we shall finish the packing of $(G_s)_{s\in \cJ_2}$. To this end we need to pack $\SpecPaths^*_s$ for each $s\in \cJ_2$. The main desired feature we want to achieve in this stage is a precise given count of edges of $\HStageE$ within $\Vmin$, within $\Vplus$, and from $\Vmin$ to $\Vplus$. To motivate where these count requirements come from, we need to outline what we do in Stage~F and Stage~G.
\begin{enumerate}[label={[fg\arabic*]}]
	\item In Stage~F, we pack $\{\SpecPaths^*_s\}_{s\in\cJ_1}$. Recall that each path $P$ from this set has~11 edges. If one endvertex of $P$ was mapped (in Stage~A) to $\Vmin$ and the other to $\Vplus$, then the packing in Stage~F will be such that~3 edges of $P$ will be mapped inside $\Vmin$, 3 edges will be mapped inside $\Vplus$ and~5 edges will be mapped into $H[\Vmin,\Vplus]$. If both endvertices of $P$ were mapped (in Stage~A) to $\Vmin$, then the packing in Stage~F will be such that~2 edges of $P$ will be mapped inside $\Vmin$, 3 edges will be mapped inside $\Vplus$ and~6 edges will be mapped into $H[\Vmin,\Vplus]$. The situation when both endvertices of $P$ were mapped (in Stage~A) to $\Vplus$ is symmetric. 
	\item In Stage~G, we pack $\{\SpecShortPaths_s\}_{s\in\cJ_0}$. Recall that each path $P$ from this set has~3 edges. The packing will be such that~1 edge of $P$ will be mapped inside $\Vmin$, 1 edge will be mapped inside $\Vplus$ and~1 edge will be mapped into $H[\Vmin,\Vplus]$. (Recall Figure~\ref{fig:outline}.)
\end{enumerate}
So, denote by \index{$j_{\boxminus \boxplus}$}\index{$j_{\boxminus \boxminus}$}\index{$j_{\boxplus \boxplus}$} $j_{\boxminus \boxplus}$ (resp.\ $j_{\boxminus\boxminus}$ and $j_{\boxplus\boxplus}$) the number of paths from $\bigcup_{s\in\cJ_1}\SpecPaths_s^*$ with one anchor in $\Vmin$ (resp.\ in $\Vmin$, and in $\Vplus$) and the other one in $\Vplus$ (resp.\ in $\Vmin$, and in $\Vplus$). The above description explains why we require the counts~\eqref{eq:E:edges}.

Define \index{$\dStageE^*$}\index{$\dStageEbar$}
\begin{equation}\label{eq:dStageE-def}
\dStageE^*:= 8\Big(\sigmKJ^*\sigmJnula+\frac{11}{4}\sigmJjedna^2\Big) \quad\text{and}\quad \dStageEbar:= 4\Big(\sigmKJ^*\sigmJnula+\frac{11}{2}\sigmJjedna^2\Big)\,.
\end{equation}
The following fact is immediate.
\begin{fact}\label{fact:min:d_E}
	We have $\min\{\dStageE^*, \dStageEbar\}\ge 4\sigmJnula\sigmKJ^*$ and $\dStageE^*-\dStageEbar= 4\sigmKJ^*\sigmJnula$.
\end{fact}

The numbers $\dStageEbar$ and $\dStageE^*$ (defined in~\eqref{eq:dStageE-def}) will be close to the edge densities between the vertex sets $\Vmin$ and $\Vplus$, and within either of the sets $\Vmin$ and $\Vplus$, respectively, after Stage~E finishes. This is expressed by the quasirandomness setup below. The proof of the following lemma can be found in Section~\ref{sec:StageE2}.

\begin{lemma}[Stage~E, density correction]\label{lem:StageENew}
\index{density correction!statement of Lemma~\ref{lem:StageENew}}\index{Stage E!statement of Lemma~\ref{lem:StageENew}}
\addtocounter{LemmaECounter}{\value{theorem}}	
  Given partial embeddings $(\phi^\mathbf{D}_s)_{s\in\cG}$ as described in
  Lemma~\ref{lem:StageDNew}, there exist maps $(\phi^\mathbf{E}_s)_{s\in\cG}$
  extending the $(\phi^\mathbf{D}_s)_{s\in\cG}$ with the following properties.
  \begin{enumerate}[label=\itmrom{E}]
    \item\label{E:one} 
      $(\phi^\mathbf{E}_s)_{s\in\cG}$ is a packing of subgraphs
      of $\left(G_s\right)_{s\in\cG}$ into~$H$ with leftover graph
      $\HStageE$\index{$\HStageE$} on $\Vmin\dcup\Vplus$. 
    \item\label{E:two} For each $s\in\cJ\setminus\cJ_2$ we have
      $\phi^\mathbf{E}_s=\phi^\mathbf{D}_s$. For each $s\in\cJ_2$, the map
      $\phi^\mathbf{E}_s$ embeds all vertices of $G_s$.
    \item\label{E:three}  
      We have
      \begin{align}
        \begin{split} \label{eq:E:edges}
          e_{\HStageE}(\Vmin)&=\sigmJnula\sigmKJ^* n^2+2 j_{\boxminus\boxminus}+3(j_{\boxminus\boxplus}+ j_{\boxplus\boxplus})\;,\\
          e_{\HStageE}(\Vplus)&=\sigmJnula\sigmKJ^* n^2+2 j_{\boxplus\boxplus}+ 3(j_{\boxminus\boxplus}+ j_{\boxminus\boxminus}) \;,\\
          e_{\HStageE}(\Vmin,\Vplus)&=\sigmJnula\sigmKJ^* n^2+6( j_{\boxminus\boxminus}+ j_{\boxplus\boxplus})+5 j_{\boxminus\boxplus}\;.
        \end{split}
      \end{align}
    \item\label{E:four}
      With the given $\cJ_0,\cJ_1$ and $\emptyset$ instead of $\cJ_2$, for each
      $s\in\cJ_0$ let $F_s$ be $\SpecShortPaths_s$, and for each $s\in\cJ_1$ let
      $F_s=\SpecPaths^*_s$. For each $i=0,1$ and $s\in\cJ_i$ let
      $U_s=\im(\phi^\mathbf{E}_s)$, let $\phi_s$ be the restriction of
      $\phi^\mathbf{E}_s$ to the leaves of $F_s$, and let $A_s$ be the image of
      $\phi_s$. For each $s\in\cJ_0$ let
      $I_s=\{i\in[\lfloor\frac{n}2\rfloor]\,:\,\boxminus_i,\boxplus_i\in
      A_s\}$. Then we have a $(\gamNew,\LE,\dStageE^*,\dStageEbar)$-quasirandom setup.
  \end{enumerate}
\end{lemma}

At this moment, it only remains to embed $\bigcup_{s\in\cJ_{1}}\SpecPaths^*_s$ and $\bigcup_{s\in\cJ_{0}}\SpecShortPaths_s$. 
In Stage~F we embed the former family in such a way that it will allow embedding the latter family in Stage~G. In order to see how the embedding should look in Stage~F, we thus need to look at the situation we will be in before and during Stage~G. For any $v\in \Vmin\dcup\Vplus$, let us write $t(v)$\index{$t(v)$} for the number of paths from $\bigcup_{s\in\cJ_0}\SpecShortPaths_s$ anchored at $v$. Observe that since for each $s\in\cJ_0$ the set $\SpecShortPaths_s$ is a set of disjoint paths whose endvertices were embedded in Stage~D, we have that $t(v)$ is equal to the number of $s\in\cJ_0$ such that $v\in A_s$, which by~\ref{D:three} and~\ref{quasi:NumAnchors} is
\begin{equation}\label{eq:t(v)}
t(v)=(1\pm \gamAnchor)2\sigmKJ^*\sigmJnula n\;.
\end{equation}
Since each path of $\SpecShortPaths_s$ is anchored at a terminal pair by~\ref{D:two}, for each $i\in[\lfloor\frac{n}2\rfloor]$ we have $t(\boxminus_i)=t(\boxplus_i)$.

In Stage~G, we will embed all paths as follows: given a path $abcd\in\SpecShortPaths_s$, with $\phi^\mathbf{F}_s(a)=\boxminus_i$ and $\phi^\mathbf{F}_s(d)=\boxplus_i$, we will embed $b$ into $\Vmin$ and $c$ into $\Vplus$. Consider some $\boxminus_i$. In Stage~G, we will use $t(\boxminus_i)$ edges from $\boxminus_i$ into $\Vmin$ to deal with the paths anchored at $\boxminus_i$, and (since we complete a perfect packing in Stage~G) all the remaining edges from $\boxminus_i$ to $\Vmin$ are used to accommodate inner vertices of paths in $\bigcup_{s\in\cJ_0}\SpecShortPaths_s$, and these paths also use all the edges from $\boxminus_i$ to $\Vplus$. A similar statement holds for $\boxplus_i$. It follows that in $\HStageF$ we need, for each $i\in[\lfloor\frac{n}2\rfloor]$,
\begin{equation}\label{eq:stageF:degtarget}
 \deg_{\HStageF}(\boxminus_i;\Vplus)=\deg_{\HStageF}(\boxminus_i;\Vmin) - t(\boxminus_i)\,\text{ and }\,\deg_{\HStageF}(\boxplus_i;\Vmin)=\deg_{\HStageF}(\boxplus_i;\Vplus) - t(\boxplus_i)\,.
\end{equation}

We will see that this is essentially the only requirement for Stage~F. That is, the fact that we have after Stage~E a quasirandom setup, together with the fact that we will use only very few edges at any given vertex in Stage~F, ensures that we still have a quasirandom setup after Stage~F, which is one of the key properties needed for Stage~G. The following lemma is proved in Section~\ref{sec:StageF}.

\begin{lemma}[Stage~F, degree correction]\label{lem:StageFNew}
\index{degree correction!statement of Lemma~\ref{lem:StageFNew}}\index{Stage F!statement of Lemma~\ref{lem:StageFNew}}
  Given partial embeddings $(\phi^\mathbf{E}_s)_{s\in\cG}$ as described in Lemma~\ref{lem:StageENew},
  there exist maps $(\phi^\mathbf{F}_s)_{s\in\cG}$
  extending the $(\phi^\mathbf{E}_s)_{s\in\cG}$ with the following properties.
  \begin{enumerate}[label=\itmrom{F}]
    \item\label{F:one} 
      $(\phi^\mathbf{F}_s)_{s\in\cG}$ is a packing of subgraphs
      of $\left(G_s\right)_{s\in\cG}$ into~$H$ with leftover graph
      $\HStageF$\index{$\HStageF$} on $\Vmin\dcup\Vplus$. 
    \item\label{F:two} For each $s\in\cJ_0$ we have
      $\phi^\mathbf{F}_s=\phi^\mathbf{E}_s$ and for each $s\notin\cJ_0$ the map
      $\phi^\mathbf{F}_s$ embeds all vertices of $G_s$.
    \item\label{F:three}
      Equation~\eqref{eq:stageF:degtarget} holds for each $i\in[\lfloor n/2\rfloor]$.
    \item\label{F:four} With the given $\cJ_0$, $\emptyset$ instead of $\cJ_1$,
      and $\emptyset$ instead of $\cJ_2$, for each $s\in\cJ_0$ let
      $F_s=\SpecShortPaths_s$, let $U_s=\im(\phi^\mathbf{F}_s)$, let $\phi_s$ be
      the restriction of $\phi^\mathbf{F}_s$ to the leaves of $F_s$, let $A_s$
      be the image of $\phi_s$, and let $I_s=\{i\in[\lfloor
        n/2\rfloor]\,:\,\boxminus_i,\boxplus_i\in A_s\}$.  Then we have a
      $(\gamJ,\LF,\dStageF^*,\dStageFbar)$-quasirandom setup, for suitable
      numbers $\dStageF^*,\dStageFbar>0$\index{$\dStageF^*$}\index{$\dStageFbar$}.
  \end{enumerate}
\end{lemma}

It remains to pack $\{\SpecShortPaths_s\}_{s\in \cJ_0}$. This is done in our last stage, encapsulated by the following lemma, which we prove in Section~\ref{sec:StageG}.

\begin{lemma}[Stage~G, designs completion]\label{lem:StageG}
\index{design completion!statement of Lemma~\ref{lem:StageG}}\index{Stage G!statement of Lemma~\ref{lem:StageG}}	
  Given partial embeddings $(\phi^\mathbf{F}_s)_{s\in\cG}$ as described in
  Lemma~\ref{lem:StageFNew}, there exist maps $(\phi^\mathbf{G}_s)_{s\in\cG}$
  extending the $(\phi^\mathbf{F}_s)_{s\in\cG}$ such that
  $(\phi^\mathbf{G}_s)_{s\in\cG}$ is a packing of $\left(G_s\right)_{s\in\cG}$
  into~$H$.
\end{lemma}

The key tool of the proof of Lemma~\ref{lem:StageG} is Keevash's machinery of generalised designs, which we introduced in a tailored setting in Proposition~\ref{prop:Designs-oursetting}. Again, recall that the key features for this machinery to work are the $(\gamJ,\LF,\dStageF^*,\dStageFbar)$-quasirandom setup, \eqref{eq:stageF:degtarget} and that no terminal pair forms an edge (which we obtained in Stage~C).

Once all these stages are complete, the proof of Theorem~\ref{thm:maintechnical} is finished.

\subsubsection{Evolution of leftover densities}\label{sssection:EvolutionOfDensities}
To guide the reader, let us briefly summarise how the density of the leftover host graph evolves throughout the packing stages and relates to the hierarchy of constants in~\eqref{eq:CONSTANTS}. The initial graph $H$ has density at least $d$. In Stage~A, we pack the bulk of the graphs, and the density of the leftover graph $\HStageA$ drops significantly to $\dStageA = (14\pm 0.1)\sigmKJ^2$ (as calculated in~\eqref{eq:dStageAapprox}). Because Stages~B and~C embed only $O(n^{1.6})$ edges in total, the densities $\dStageB$ and $\dStageC$ remain $\dStageA \pm o(1) \approx 14\sigmKJ^2$. In Stage~D, we embed the outer parts of the paths from $\cJ_0$, reducing the density to $\dStageD$, which is slightly smaller but still $\Theta(\sigmKJ^2)$ by~\eqref{eq:dStageD}. In Stage~E, we pack the paths from $\cJ_2$, which consume a significant number of edges, and the leftover densities drop to $\dStageE^*, \dStageEbar = \Theta(\sigmJnula\sigmKJ^*)$ (defined in~\eqref{eq:dStageE-def}; see also Fact~\ref{fact:min:d_E}). Finally, Stage~F only makes fine degree adjustments, so $\dStageF^*, \dStageFbar$ remain $\Theta(\sigmJnula\sigmKJ^*)$. 


\section{Packing paths}\label{sec:pathpack}
A \emph{path-forest}\index{path-forest} is a non-empty graph each of whose components is a path of length at least~2.
 In this section we show how to pack a collection of path-forests into a quasirandom graph $H$, with several constraints. First, the quasirandom graph will be partitioned into two \emph{sides} $V_{\boxminus}$ and $V_{\boxplus}$ of equal size, and the vertices of each path in each path-forest will be pre-assigned to one side or the other. Furthermore, the path-forests are `anchored', i.e.\ each leaf is pre-embedded to a specific vertex of $H$, and each path-forest is required to avoid a set of `used' vertices given for that path-forest. Helpfully, the path-forests we pack will always contain only paths of bounded length, they will never be spanning (in fact they will always miss a positive proportion of the unused vertices in each part of $H$), and the assignment of sides is such that we finish the packing with a positive density of edges in each side and between the two sides of $H$ unused.
 
 Such a packing is not generally possible; however we will assume that the anchors and used vertices are quasirandomly distributed, and under this additional condition we show the following simple random algorithm works. We go through the path-forests in turn and in each path-forest pack one (randomly chosen) path after another, in each case choosing one of the possible paths consistent with the assigned sides uniformly at random. Our aim is to show that the resulting packing `looks random' in the sense that in each path-forest the vertices used, and in total the edges in and between sides used, look like uniform random sets of the appropriate sizes.
The main result of this section is Lemma~\ref{lem:pathpack}.
 
\subsection{Embedding one path-forest} 
 
 To begin with, we formally state the algorithm we use to embed one path-forest, and give the analysis of this algorithm.

\begin{algorithm}[ht]\index{\RandomPathEmbedding{}} 
	\caption{\RandomPathEmbedding{}}\label{alg:path-embed}
	\SetKwInOut{Input}{Input}
	\Input{
    $\bullet$ a path-forest~$F$ with anchors~$A$, \\
    $\bullet$ a graph~$H$ on $\Vmin\dcup \Vplus$, \\
    $\bullet$ a set $U_0\subset V(H)$ of used vertices, \\
    $\bullet$ an assignment $\chi:V(F)\setminus A\to\{\Vmin,\Vplus\}$, \\
    $\bullet$ an embedding $\phi_0:A\to V(H)$ such that $\im\phi_0\subset U_0$}
  choose a uniform random order $P_1,\dots,P_{t^*}$ of the components of $F$\;
  $\psi_0:=\phi_0$\;
  \For{$t=1$ \KwTo $t^*$}{
    let $x_1,\dots,x_k$ be the path~$P_t$\;
    let $v_1=\phi_0(x_1)$, and $v_k=\phi_0(x_k)$\;
    let $\Paths_t$ be the set of $v_1$-$v_k$-paths
    $v_1,v_2,\dots,v_k$ of length $k-1$ in $H-U_{t-1}$
    such that $v_i\in \chi(x_i)$ for each $1<i<k$\;
    \lIf{$\Paths_t$ is empty}{
			halt with failure}
    choose a path $R$ uniformly at random in $\Paths_t$\;
		$\psi_{t}:=\psi_{t-1}\cup\{P_t\AlgMap R\}$\;
    $U_t:=U_0\cup \im \psi_t$\;
	}
	\Return $\phi=\psi_{t^*}$\;
\end{algorithm}

In order to analyse this algorithm, we only need to specify what precisely we mean by a quasirandom distribution of the anchors. 

\begin{definition}[anchor distribution property]\label{def:anchordistribution}
 Suppose that we are given an $n$-vertex graph $H$ with $V(H)=\Vmin\dcup\Vplus$, a subset $U$ of $V(H)$, a path-forest $F$ whose leaves form the set $A$. For $\gamma\ge 0$, we say an injective map $\phi_0:A\to U$ and assignment $\chi: V(F)\setminus A\to\{\Vmin,\Vplus\}$ have the \emph{$\gamma$-anchor distribution property}\index{anchor distribution property} if the following holds.
 
 Given $a,b,c\in\{\boxminus,\boxplus\}$ let $A_{a,b,c}$ denote the collection of $x\in A$ such that $\phi_0(x)\in V_a$, and such that the neighbour $y$ of $x$ satisfies $\chi(y) = V_b$, and the next vertex $z$ (adjacent to $y$) satisfies $\chi(z) =  V_c$. Define $d_{\boxminus\boxplus}=d_{\boxplus\boxminus}$ to be the bipartite density of $H[V_\boxplus,V_\boxminus]$, $d_{\boxplus\boxplus}$ to be the density of $H[V_\boxplus]$, and $d_{\boxminus\boxminus}$ to be the density of $H[V_\boxminus]$. For each $v\in V_b\setminus U$ we have
 \[\big|\big\{x\in A_{a,b,c}:\phi_0(x)\in \NBH_H(v)\big\}\big|=(1\pm\gamma)d_{ab}|A_{a,b,c}|\pm \tfrac12\gamma n^{0.99}\,.\]
 
 We write $A_{a,b}=A_{a,b,\boxminus}\cup A_{a,b,\boxplus}$, and normally we will only need to refer to $A_{a,b}$; the $\gamma$-anchor distribution property implies
 \[\big|\big\{x\in A_{a,b}:\phi_0(x)\in \NBH_H(v)\big\}\big|=(1\pm\gamma)d_{ab}|A_{a,b}|\pm \gamma n^{0.99}\,.\]
\end{definition}

The following lemma contains all the facts we need about the running of Algorithm~\ref{alg:path-embed}.

\begin{lemma}[Embedding a path-forest]\label{lem:RPE} Given an integer $L\ge 5$ and real $\nu>0$, there is a constant $C>0$ such that the following holds for all $0<\gamma<C^{-1}$ and all sufficiently large $n$ (depending on $\gamma$).

  Let~$F$ be a path-forest with leaves~$A$ (also called anchors), in which each path has at least~$5$ and at most~$L$ vertices. Suppose $\nu n\le |V_\boxplus|,|V_\boxminus|\le n$, and let~$H$ be $(\gamma,L)$-block-quasirandom
  on $\Vmin\dcup \Vplus$ with densities
  $d(H[\Vmin])=d_{\boxminus\boxminus}$,
  $d(H[\Vplus])=d_{\boxplus\boxplus}$,
  $d(H[\Vmin,\Vplus])=d_{\boxminus\boxplus}$, each at least $\nu$.
  Let $\phi_0:A\to V(H)$ be an embedding of the anchors and 
  $\chi:V(F)\to\{\Vmin,\Vplus\}$ be an assignment of the
  vertices of the paths to the parts of~$H$ such that $\phi_0(x)\in\chi(x)$ for each $x\in A$.
  For each $a\in\{\boxminus,\boxplus\}$ let $X_a:=\{x\in V(F)\;:\;\chi(x)=V_a\}$, and suppose that $|X_a\setminus A| \ge \nu n$.
  Let $U\subset V(H)$ be a set of ``used'' vertices such that $\im\phi_0\subset U$, and suppose that for each $a\in\{\boxminus,\boxplus\}$ we have
  \[\big|V_a\setminus U\big|-\big|\big\{x\in V(F)\setminus A:\chi(x)=V_a\big\}\big|\ge \nu n\,.\]
  Suppose that we have the $\gamma$-anchor distribution property, and $(H,U)$ is $(\gamma,L)$-block-diet.
  Suppose that $\big|\big\{(x,y)\,:\,xy\in E(F),x\in
        X_a\setminus A,y\in X_b\setminus A\big\}\big|\ge \nu n$ for each $a,b\in\{\boxminus,\boxplus\}$.
  Let~$V'$ be a set of specified vertices with $|V'|\ge\nu n$ and $V'\subset\Vmin \setminus U$ or
  $V'\subset\Vplus \setminus U$.

  When we execute \RandomPathEmbedding{} (Algorithm~\ref{alg:path-embed}), then
  \begin{enumerate}[label=\itmarab{PE}]
  \item\label{RPE:nofail} with probability at least $1-\exp\big(-n^{0.3}\big)$, we obtain an embedding $\phi:V(F)\to
    \im\phi_0\cup\big(V(H)\setminus U\big)$ of $F$ in $H$ extending $\phi_0$
    and such that $\phi(x)\in\chi(x)$ for each $x\in V(F)\setminus A$.
  \end{enumerate}
  Moreover, if we get such an embedding~$\phi$, the following hold.
  \begin{enumerate}[start=2,label=\itmarab{PE}]
  \item\label{RPE:imageeating} If $a\in\{\boxminus,\boxplus\}$ is such that $V'\subset V_a\setminus U$, then
    \[\Prob\Big[|\im\phi\cap V'|\neq(1\pm C\gamma)|V'|\tfrac{|X_a\setminus A|}{|V_a\setminus U|}\Big]\le \exp\big(-n^{0.3}\big)\,.\]
  \item\label{RPE:probvtx} For each $S_\boxminus\subset\Vmin\setminus U,S_\boxplus\subset \Vplus\setminus U$ with $|S_a|\le L$ for each $a$, we have
    \[\Prob\big[(S_\boxminus\cup S_\boxplus)\cap \im\phi=\emptyset\big]=(1\pm C\gamma)\prod_{a\in\{\boxminus,\boxplus\}}\Big(\tfrac{|V_a\setminus U|-|\{x\in V(F)\setminus A:~\chi(x)=V_a\}|}{|V_a\setminus U|}\Big)^{|S_a|}\,.\]
  \item\label{RPE:probedgefree} For each edge $e\in E(H-U)$, setting $a,b\in\{\boxminus,\boxplus\}$ such
    that $e\in E(H[V_a,V_b])$, we have
    \[\Prob\big[\phi \text{ uses $e$}\big]=(1\pm C\gamma)\frac{\big|\big\{(x,y)\,:\,xy\in E(F),x\in
        X_a\setminus A,y\in X_b\setminus A\big\}\big|}{d_{ab}|V_a\setminus
        U||V_b\setminus U|}\,.\]
  \item\label{RPE:probedgeside}
	For each $a,b\in\{\boxminus,\boxplus\}$,
	each $x\in A$
    such that $x\in X_a$ and the neighbour of $x$ in $F$ is in $X_b$,
    and each 
    $v\in V_b\setminus U$ with $\phi_0(x)v\in E(H)$
    we have 
    \[\Prob\big[\phi \text{ uses $\phi_0(x)v$}\big]=(1\pm C\gamma)\frac{1}{d_{ab}|V_b\setminus U|}\,.\]
  \item\label{RPE:probsec} For each pair of edges $uv,u'v'\in E(H)$, we have
    \[\Prob[\phi \text{ uses $uv$ and $u'v'$}] \le n^{-3/2}\,.\]
  \end{enumerate}
\end{lemma}
\begin{proof} Given an integer $L\ge 5$ and $\nu>0$, we set
\[C:=\exp\big(10^4L^4\nu^{-6}\big)\,,\]
and given $0<\gamma<1/C$, we let $t^*$ be the number of components of $F$, and for each $x\in\RR$ we let
\[\beta_x:=10\gamma\exp\Big(\tfrac{1000L^2x}{\nu^3t^*}\Big)\,.\]

We set $\eps= \gamma^2C^{-1}$, and note that $\beta_{t^*}\le \sqrt{C}\gamma$.

First, we show that after randomly ordering the paths of $F$, each interval of $\eps n$ consecutive paths from this order has the anchor distribution property. This avoids the possibility that neighbours of anchors are distributed in a non-uniform way for most of the process: this imbalance would eventually be corrected, but it would complicate the analysis.

\begin{claim}\label{cl:RPE:orderpaths} In the uniform random order $P_1,\dots,P_{t^*}$ of the components of $F$, with probability at least $1-\exp\big(-n^{0.4}\big)$ for each $1\le i\le t^*+1-\eps n$ the forest $F^*_i$ (with the same $U$, and $\chi$ and $\phi_0$ restricted to $F^*_i$) whose components are $P_i,P_{i+1},\dots,P_{i+\eps n-1}$ has the $\gamma$-anchor distribution property. Furthermore, for each $a\in\{\boxminus,\boxplus\}$, we have $\tfrac{|X_a\cap V(F^*_i)\setminus A|}{|V(F^*_i)\setminus A|}=(1\pm\eps)\tfrac{|X_a\setminus A|}{|V(F)\setminus A|}$.
\end{claim}
\begin{claimproof}[Proof of Claim~\ref{cl:RPE:orderpaths}]
 Since for each $1\le i\le t^*+1-\eps n$, the distribution of paths in $F$ selected to $F^*_i$ is a uniform random set of $\eps n$ paths, by the union bound it suffices to show the above properties of $F^*_i$ hold with probability $1-\exp\big(-n^{0.5}\big)$ for a uniform random set $F^*$ of $\eps n$ paths of $F$. Note that $|X_a|\ge\nu n$ for each $a$, and so in particular $t^*\ge\nu n/L > 2 \eps n$.

 To begin with, we show that the vertices in $X_a\setminus A$ are well distributed. Given $a\in\{\boxminus,\boxplus\}$, we partition the path components of $F$ according to their number of vertices in $X_a\setminus A$. Suppose there are $u_j$ paths with $j$ vertices in $X_a\setminus A$ for each $0\le j\le L$ (and note no path has more than $L$ vertices). Then in expectation we choose $\tfrac{u_j\eps n}{t^*}$ paths with $j$ vertices in $X_a\setminus A$ to $F^*$, and so by Fact~\ref{fact:hypergeometricBasicProperties} the probability we do not pick $\tfrac{u_j\eps n}{t^*}\pm n^{0.98}$ such paths is at most $\exp\big(-n^{0.9}\big)$. Assuming none of these bad events occur, we conclude
 \[|X_a\cap V(F^*)\setminus A|=\sum_{j=0}^L j\cdot \Big(\tfrac{u_j\eps n}{t^*}\pm n^{0.98}\Big)=\tfrac{|X_a\cap V(F)\setminus A|\eps n}{t^*}\pm L^2n^{0.98}\,.\]
So with probability at least $1-\exp\big(-n^{0.8}\big)$ the above estimate holds for each $a\in\{\boxminus,\boxplus\}$, 
which implies
 $\tfrac{|X_a\cap V(F^*)\setminus A|}{|V(F^*)\setminus A|}=(1\pm\eps)\tfrac{|X_a\setminus A|}{|V(F)\setminus A|}$ for each $a\in\{\boxminus,\boxplus\}$.

We now show $F^*$ is likely to have the $\gamma$-anchor distribution property. To that end, fix a vertex $v\in V(H)\setminus U$. For any $a,b,c\in\{\boxminus,\boxplus\}$, consider the set $A_{a,b,c}$. Each path of $F$ may have $0$, $1$, or $2$ endpoints $x\in A_{a,b,c}$ such that $\phi_0(x)\in \NBH_H(v)$, and we partition the paths of $F$ accordingly into sets of sizes $u_0,u_1,u_2$ respectively. 
By the given $\gamma$-anchor distribution property we have
$
\sum_{j=0}^2 ju_j = (1\pm \gamma) d_{ab}|A_{a,b,c}| \pm \tfrac12\gamma n^{0.99}.
$
Moreover, as before, the expected number of paths in the $j$th set we select to $F^*$ is $\tfrac{u_j\eps n}{t^*}$, and by Fact~\ref{fact:hypergeometricBasicProperties} the probability of failing to select $\tfrac{u_j\eps n}{t^*}\pm n^{0.98}$ of these paths is at most $\exp\big(-n^{0.8}\big)$. Assuming none of these bad events occur, we see that the total number of endpoints $x\in A_{a,b,c}$ of paths of $F^*$ such that $\phi_0(x)\in \NBH_H(v)$ is
\[\sum_{j=0}^2j\cdot\Big(\tfrac{u_j\eps n}{t^*}\pm n^{0.98}\Big)=\tfrac{(1\pm\gamma)d_{ab}|A_{a,b,c}|\eps n\pm\tfrac12\gamma\eps n\cdot n^{0.99}}{t^*}\pm 3n^{0.98}.\]
By a similar argument, we have
\[\big|A_{a,b,c}\cap V(F^*)\big|=\tfrac{|A_{a,b,c}|\eps n}{t^*}\pm n^{0.98}\]
with probability at least $1-\exp\big(-n^{0.8}\big)$. It follows that the number of endpoints $x\in A_{a,b,c}$ of paths of $F^*$ such that $\phi_0(x)\in \NBH_H(v)$ is
\[\big(1\pm\gamma\big)d_{ab}\big|A_{a,b,c}\cap V(F^*)\big|\pm \tfrac14\gamma n^{0.99}\pm 5n^{0.98}\]
with probability at least $1-\exp\big(-n^{0.7}\big)$, where we use that $t^\ast>2\eps n$. Taking the union bound over choices of $a$, $b$, $c$ and $v$ we see the probability that this part of the 
$\gamma$-anchor distribution property holds for $F^*$ is at least $1-\exp\big(-n^{0.6}\big)$.
\end{claimproof}

From this point on we will assume that the outcome of the uniform random choice of $P_1,\dots,P_{t^*}$ is fixed and the likely statements of the above claim hold. In particular, this means that (unless Algorithm~\ref{alg:path-embed} fails before time $t$) the number $|V_a\setminus U_{t-1}|$ is fixed for each $a\in\{\boxminus,\boxplus\}$ (although the sets themselves are not). We next argue that vertices are used with a predictable probability in each path.

\begin{claim}\label{cl:RPE:probemb} 
Given any $t\in [t^*]$ and $0<\rho<2^{-4L}$, let $\hist_{t-1}$ denote the event that the embedding of the first $t-1$ paths succeeds and furthermore $(H,U_{t-1})$ is $(\rho,2)$-block-diet. Given any $y\in V(P_t)\setminus A$, let $x$ and $z$ be the neighbours of $y$ on $P_t$; if exactly one of these is in $A$ suppose $x\in A$. Let $a,b,c\in\{\boxminus,\boxplus\}$ be such that $\chi(x) = V_a$, $\chi(y) = V_b$ and $\chi(z) = V_c$. Given any $v\in V_b\setminus U_{t-1}$, we have
 \[\Prob\big[\psi_t(y)=v\big|\hist_{t-1}\big]=\begin{cases}
  (1\pm 5L\rho)\tfrac{1}{|V_b\setminus U_{t-1}|} & \text{if }x,z\notin A \\
  (1\pm 5L\rho)\tfrac{1}{d_{ab}|V_b\setminus U_{t-1}|} & \text{if }x\in A,\,z\notin A,\,v\in \NBH_H\big(\phi_0(x)\big) \\
 \end{cases}\quad.\]
 In all other cases, the probability of embedding $y$ to $v$ is zero.
\end{claim}
\begin{claimproof}[Proof of Claim~\ref{cl:RPE:probemb}]
Observe that since $0<\rho\le 2^{-4L}$ we have $(1\pm2\rho)^L(1\pm2\rho)^{-L}=1\pm 5L\rho$.

The probability zero statement of the claim can be seen as follows. If the first case $x,z\notin A$ does not happen, then by assumption we must have $x\in A$ and $z\notin A$, since $P_t$ has at least 5 vertices. Hence, if we additionally exclude the second case from above, what remains is the case when
$x\in A,\,z\notin A,\,v\notin \NBH_H\big(\phi_0(x)\big)$. 
However, since $y$ must be embedded to a neighbour of the images of $x$ and $z$,
the probability of embedding $y$ to $v$ is zero.
 
 Given any three vertices $p,q,r$ of $V(H)$ and any $3\le k\le \ell-2$, we can count the number of $\ell$-vertex paths from $p$ to $r$ in $H - U_{t-1}$ whose $k$th vertex is $q$ and whose $i$th vertex (for each $1\le i\le \ell$) is in $V_i\in\{\Vmin,\Vplus\}$, for any given choice of the $V_i$. We do this as follows. Mark $p,q,r$ used. We count the number of vertices of $V_2\setminus U_{t-1}$ adjacent to $p$ and not so far used, for each of these the number of vertices of $V_3\setminus U_{t-1}$ adjacent to the first and not so far used, and so on up to the $(k-1)$st vertex which we require to be adjacent to both the $(k-2)$nd and to $q$; we repeat this (avoiding previously used vertices) for the path from $q$ to $r$. 
Denote by $d_{ij}$ the bipartite density $d[V_i,V_j]$
if $V_i\neq V_j$, or the density $d[V_i]$ if $V_i=V_j$. 
Then, by the block-diet condition we have $(1\pm\rho)d_{12}|V_2\setminus U_{t-1}|\pm\ell$ choices for the second vertex, for each of these $(1\pm\rho)d_{23}|V_3\setminus U_{t-1}|\pm\ell$ choices for the third vertex, and so on until we get to $\big(1\pm\rho\big)d_{(k-2)(k-1)}d_{(k-1)k}|V_{k-1}\setminus U_{t-1}|\pm\ell$ choices for the $(k-1)$st vertex, and similarly for the path from $q$ to $r$. Here the $\pm \ell$ comes from the vertices previously used on the path we construct. Since $n$ is sufficiently large, and since $V_a\setminus U_t$ is linear in $n$ for each $a\in\{\boxminus,\boxplus\}$ while $\ell\le L$, we can absorb the $\pm\ell$ into the multiplicative error to obtain $(1\pm 2\rho)d_{12}|V_2\setminus U_{t-1}|$ choices for the second vertex, and so on. Putting this together, the total number of such paths is
 \[(1\pm2\rho)^{\ell-2}\prod_{i=1}^{\ell-1}d_{i(i+1)}\prod_{i=2}^{k-1}|V_i\setminus U_{t-1}|\prod_{i=k+1}^{\ell-1}|V_i\setminus U_{t-1}|\,.\]
 
 If $x,z\notin A$, then using the above with $p$ and $r$ being the embeddings of the ends of $P_t$, letting $q$ range over all vertices of $\chi(y)\setminus U_{t-1}$, and letting $\ell=v(P_t)$ with $y$ being the $k$th vertex of $P_t$ and assignments as specified in the lemma statement, we obtain the total number of paths from $p$ to $r$ from which Algorithm~\ref{alg:path-embed} chooses uniformly at random. We can thus compute the probability of $\psi_t(y)=v$. Most of the terms cancel, leaving
 \[\Prob\big[\psi_t(y)=v\big|\hist_{t-1}\big]=(1\pm2\rho )^{\ell-2}(1\pm2\rho )^{2-\ell}\tfrac{1}{|V_b\setminus U_{t-1}|}\,,\]
 which by choice of $\rho$ is as desired.
 
 Suppose $x\in A$ and $z\notin A$, and that $v\in \NBH_H\big(\phi_0(x)\big)$. By a similar argument, for each vertex $u\in \NBH_H\big(\phi_0(x)\big)\cap V_b$ we can estimate the number of correctly assigned paths in $H-U_{t-1}$ between the embedded anchors of $P_t$. Summing over the $(1\pm\rho )d_{ab}|V_b\setminus U_{t-1}|\pm 3$ choices of $u$ given by the block-diet condition, we find the total number of such paths and hence
 \[\Prob\big[\psi_t(y)=v\big|\hist_{t-1}\big]=(1\pm2\rho )^{\ell-1}(1\pm2\rho )^{1-\ell}\tfrac{1}{d_{ab}|V_b\setminus U_{t-1}|}\,,\]
 which again is as desired.
\end{claimproof}

In light of this claim, given $v\in V(H)\setminus U_0$ and $y\in V(P_t)\setminus A$, it is natural to define a quantity $p_{v,y}$ which approximates the probability that $y$ is embedded to $v$, conditioned on the first $t-1$ paths being embedded in such a way that $v$ is not used and the block-diet condition is maintained. To that end, suppose $\chi(y)=V_b$ and the neighbours of $y$ on $P_t$ are $x$ and $z$; if exactly one of these is in $A$, suppose it is $x$. Suppose $\chi(x)=V_a$ and $\chi(z)=V_c$. We define
\begin{equation}\label{eq:RPE:defpvy}
 p_{v,y}:=\begin{cases}
 \tfrac{1}{|V_b\setminus U_{t-1}|}&\text{if }x,z\notin A\\
 \tfrac{1}{d_{ab}|V_b\setminus U_{t-1}|} & \text{if }x\in A,\,z\notin A,\,v\in \NBH_H\big(\phi_0(x)\big) \\
\end{cases}\quad.
\end{equation}

We also define
\[p^*_{v,t}:=\sum_{t'=t}^{t+\eps n-1}\sum_{\substack{y\in V(P_{t'})\setminus A\\ \chi(y)=V_b}}p_{v,y}\,,\]
which (as we will see) is roughly the probability that $v$ is used in embedding one of the first $\eps n$ paths from $P_t$ onwards, given it was not used earlier. For now, we will use the anchor distribution property established in Claim~\ref{cl:RPE:orderpaths} to evaluate $p^*_{v,t}$.
\begin{claim}\label{cl:RPE:pstarvt} Given $1\le t\le t^*-\eps n+1$, $b\in\big\{\boxminus,\boxplus\big\}$ and $v\in V_b\setminus U_{t-1}$, let $X_{b,t}^{\mathrm{int}}$ denote the set of vertices $y$ which are internal to $P_{t'}$ for some $t\le t'<t+\eps n$ and furthermore satisfy $\chi(y)=V_b$. Then we have 
 \[p^*_{v,t}=(1\pm2\gamma)\tfrac{|X_{b,t}^{\textrm{int}}|}{|V_b\setminus U_{t-1}|}\,.\]
\end{claim}
\begin{claimproof}[Proof of Claim~\ref{cl:RPE:pstarvt}]
We split $X^{\textrm{int}}_{b,t}$ into parts corresponding to the cases of the definition of $p_{v,y}$, including the choice of $a$
in the line where this appears. Let $X^{\textrm{zero}}$ be the subset of $y\in X_{b,t}^{\textrm{int}}$ whose neighbours are not in $A$.
For each $a\in\{\boxminus,\boxplus\}$ let $A_{a,b}$ be the set of $y\in X_{b,t}^{\textrm{int}}$ whose neighbours are $x,z$ with $x\in A$ and $z\notin A$ and where $\chi(x)=V_a$.
Note that these sets partition $X_{b,t}^{\textrm{int}}$.

The contribution to $p^*_{v,t}$ from vertices of $X^{\textrm{zero}}$ is $\tfrac{|X^\textrm{zero}|}{|V_b\setminus U_{t-1}|}$. Now fix $a\in\{\boxminus,\boxplus\}$. The contribution to $p^*_{v,t}$ from vertices of $A_{a,b}$ is, by the $\gamma$-anchor distribution property from Claim~\ref{cl:RPE:orderpaths},
\[\frac{(1\pm \gamma)d_{ab}|A_{a,b}| \pm \gamma n^{0.99}}{d_{ab}|V_b\setminus U_{t-1}|}=(1\pm \gamma)\tfrac{|A_{a,b}|}{|V_b\setminus U_{t-1}|} \pm n^{-0.001}\,,\]
where the equality uses $d_{ab}\ge\nu$ and $|V_b\setminus U_{t-1}|\ge|V_b\setminus U|-|\{x\in V(F)\setminus A:\chi(x)=V_b\}|\ge\nu n$, the last inequality being an assumption of the lemma.
Summing up we obtain
\[p^*_{v,t}=\tfrac{|X^\textrm{zero}|}{|V_b\setminus U_{t-1}|}+(1\pm \gamma)\tfrac{|A_{\boxminus,b}|}{|V_b\setminus U_{t-1}|}+(1\pm \gamma)\tfrac{|A_{\boxplus,b}|}{|V_b\setminus U_{t-1}|} \pm 2n^{-0.001}=(1\pm\gamma)\tfrac{|X_{b,t}^\textrm{int}|}{|V_b\setminus U_{t-1}|}\pm2n^{-0.001}\,.\]
By the `furthermore' part of Claim~\ref{cl:RPE:orderpaths}, and since each path of $F$ has at least $5$ vertices, we have
\[|X^\mathrm{int}_{b,t}|\ge(1-\eps)\tfrac{|X_b\setminus A|}{|V(F)\setminus A|}\cdot 3\eps n\ge\tfrac{\nu n}{2n}\cdot 3\eps n\ge \eps\nu n\,,\]
and this justifies absorbing the $2n^{-0.001}$ additive error to the relative error as claimed.
\end{claimproof}

We are now in a position to estimate the number of vertices embedded to a given subset $X$ of $V_a\setminus U_{t-1}$ when $P_t,\dots,P_{t+\eps n-1}$ are embedded by Algorithm~\ref{alg:path-embed}, where $a\in\{\boxminus,\boxplus\}$.

\begin{claim}\label{cl:RPE:shortint} Given any $1\le t\le t^*-\eps n+1$, let $\hist_{t-1}$ denote the event that embedding of the first $t-1$ paths succeeds and furthermore $(H,U_{t-1})$ is $(\beta_{t-1},2)$-block-diet. Given any $b\in\{\boxminus,\boxplus\}$ and $X\subset V_b\setminus U_{t-1}$ such that $|X|\ge\nu^4 n$, when Algorithm~\ref{alg:path-embed} embeds paths $P_t,\dots,P_{t+\eps n-1}$, we have
\[\Prob\Bigg[\big|X\cap U_{t+\eps n-1}\big|=(1\pm50L\beta_{t-1})|X|\tfrac{|V_b\cap (U_{t+\eps n-1}\setminus U_{t-1})|}{|V_b\setminus U_{t-1}|}\Bigg|\hist_{t-1}\Bigg] \ge 1-\exp\big(-n^{0.5}\big)\,.\]
\end{claim}
\begin{claimproof}[Proof of Claim~\ref{cl:RPE:shortint}]
Suppose that $\hist_{t-1}$ occurs, and consider the running of Algorithm~\ref{alg:path-embed} to embed $P_t$ onwards. Fix $b\in\{\boxminus,\boxplus\}$ and $X\subset V_b\setminus U_{t-1}$.

Observe that in embedding $P_t,\dots,P_{t+\eps n-1}$ we embed in total at most $L\eps n$ vertices. Thus the estimates provided by the $(\beta_{t-1},2)$-block-diet condition for $(H,U_{t-1})$ are changed by at most $L\eps n$ for any $t-1\le t'<t+\eps n$. Since $d_{pq}\ge\nu$ and $|V_p\setminus U_{t-1}|\ge\nu n$ for each $p,q\in\{\boxminus,\boxplus\}$, and by choice of $\eps$, we see that for any $t-1\le t'<t+\eps n$, the pair $(H,U_{t'-1})$ satisfies the $(2\beta_{t-1},2)$-block-diet condition. This in particular implies that Algorithm~\ref{alg:path-embed} cannot fail at any time in this interval, and that the conclusions of Claim~\ref{cl:RPE:probemb} are valid with $\rho=2\beta_{t-1}$ for each of these paths.

Given $v\in X$ and $t\le t'\le t+\eps n-1$, the probability that a vertex of $P_{t'}$ is embedded to $v$, conditioned on the embedding of $P_t,\dots,P_{t'-1}$, is either zero (if a vertex of $P_t,\dots,P_{t'-1}$ has been embedded to $v$) or given by summing, over all internal vertices $y$ of $P_{t'}$ such that $\chi(y)=V_b$, the probability $\big(1\pm 10L\beta_{t-1}\big)p_{v,y}$ given by Claim~\ref{cl:RPE:probemb}. The latter statement holds because the events $\psi_{t'}(y)=v$ for different $y$ in $P_{t'}$ are disjoint. We wish to estimate
\[Z=\sum_{t'=t}^{t+\eps n-1}\Exp\Big[\big|\psi_{t'}\big(V(P_{t'})\big)\cap X\big|\,\Big|\,\text{embeddings of $P_t,\dots,P_{t'-1}$}\Big]\,,\]
and it follows that an upper bound for this sum is $(1\pm 10L\beta_{t-1})\sum_{v\in X}p^*_{v,t}$. Furthermore, the difference between this upper bound and the correct value of the (random variable) $Z$ is that vertices $v$ to which we embed at some time cease to contribute to the sum; there are at most $L\eps n$ such vertices, and so we obtain
\[Z=(1\pm 10L\beta_{t-1})\sum_{v\in X}p^*_{v,t}\, \pm \,L\eps n\cdot 2\tfrac{|X_{b,t}^{\textrm{int}}|}{|V_b\setminus U_{t-1}|}=(1\pm 20L\beta_{t-1})\tfrac{|X|\cdot|X^{\textrm{int}}_{b,t}|}{|V_b\setminus U_{t-1}|}\,,\]
where we use Claim~\ref{cl:RPE:pstarvt} to estimate $p^*_{v,t}$.

Now consider the process of embedding $P_t,\dots,P_{t+\eps n-1}$ following Algorithm~\ref{alg:path-embed}. By Corollary~\ref{cor:freedm}\ref{cor:freedm:tails}, since each path has at most $L$ vertices, we obtain
\[\Prob\bigg[|X\cap U_{t+\eps n-1}\big|\neq(1\pm 40L\beta_{t-1})\tfrac{|X|\cdot|X^{\textrm{int}}_{b,t}|}{|V_b\setminus U_{t-1}|}\bigg|\hist_{t-1}\bigg]\le\exp(-n^{0.5})\,.\]

Since $\big|X^{\textrm{int}}_{b,t}\big|=\big|V_b\cap (U_{t+\eps n-1}\setminus U_{t-1})\big|$, this completes the proof.
\end{claimproof}

Using the above claim, we can show that, given any $1<t'\le t^*$ , any $b\in\{\boxminus,\boxplus\}$ and any large enough $X\subseteq V_b\setminus U_0$, with very high probability the following holds when Algorithm~\ref{alg:path-embed} is run. Either the $(\beta_{t-1},2)$-block-diet condition fails for some $(H,U_{t-1})$ with $1\le t< t'\le t^*$, or $|X\cap U_{t'}|$ is roughly $|X|\tfrac{|V_b\cap (U_{t'}\setminus U_0)|}{|V_b\setminus U_0|}$.

\begin{claim}\label{cl:RPE:longint} Given any $t\ge 1$, any $b\in\{\boxminus,\boxplus\}$ and any $X\subset V_b\setminus U_0$ such that $|X|\ge\nu^2 n$, when Algorithm~\ref{alg:path-embed} embeds paths $P_1,\dots,P_{t}$, we have
\begin{multline*}
 \Prob\Big[\big|X\setminus U_{t}\big|=\big(1\pm\tfrac34\beta_{t}\big)|X|\tfrac{|V_b\setminus U_t|}{|V_b\setminus U_0|}
 \text{ or $\exists t'\in[t-1]\colon (H,U_{t'})$ is not $(\beta_{t'},2)$-block-diet}\Big]\\
 \ge 1-\exp\big(-n^{0.4}\big)\,.
\end{multline*}
\end{claim}
\begin{claimproof}[Proof of Claim~\ref{cl:RPE:longint}]
 The proof of this claim is simply an iteration of Claim~\ref{cl:RPE:shortint}, similar to the proof of~\cite[Lemma~\ref{PAPERpackingdegenerate.lem:dietcon}]{DegPack}. 
For this, observe that $\big|X\setminus U_{t'}\big|=(1\pm\beta_{t'})|X|\tfrac{|V_b\setminus U_{t'}|}{|V_b\setminus U_0|}$
implies $|X\setminus U_{t'}|\ge \nu^4 n$, since $|X|\ge \nu^2n$ and
$|V_b\setminus U_{t'}|\ge \nu n$, and hence we are allowed to apply Claim~\ref{cl:RPE:shortint} iteratively. 
 
Now, given $t\ge1$, $b\in\{\boxminus,\boxplus\}$ and $X\subset V_b\setminus U_0$, we suppose that the likely event of Claim~\ref{cl:RPE:shortint} holds for each time $t'$ with $1\le t'\le t-\eps n+1$, with the input set $X\setminus U_{t'-1}$ (provided this set has size at least $\nu^4 n$). The probability that this likely event does not occur, but $(H,U_{t'})$ is $(\beta_{t'},2)$-block-diet for each such $t'$, is by the union bound at most $n\exp\big(-n^{0.5}\big)<\exp\big(-n^{0.4}\big)$. We claim that outside of the unlikely event the claim statement holds. 

Let $s=\lfloor\tfrac{t}{\eps n}\rfloor$. We require specifically the above likely event for each $t'=1+k\eps n$, where $k$ is an integer between $0$ and $s-1$ inclusive; this yields
\begin{align*}
\big|X\setminus U_{s\eps n}\big|
	& = \big|X\big| \cdot \prod_{k=0}^{s-1}\left(
		1 - (1\pm 50L\beta_{k\eps n})
		 \tfrac{|V_b\cap (U_{(k+1)\eps n}\setminus U_{k\eps n})|}{|V_b\setminus U_{k\eps n}|}\right) \\
	& = \big|X\big| \cdot \prod_{k=0}^{s-1}\left(
		\tfrac{|V_b\setminus U_{(k+1)\eps n}|}{|V_b\setminus U_{k\eps n}|} \pm 50L\beta_{k\eps n}
		 \tfrac{|V_b\cap (U_{(k+1)\eps n}\setminus U_{k\eps n})|}{|V_b\setminus U_{k\eps n}|}\right) \\
	& = \big|X\big| \cdot \prod_{k=0}^{s-1}
		\tfrac{|V_b\setminus U_{(k+1)\eps n}|}{|V_b\setminus U_{k\eps n}|} \left( 1 \pm 50L\beta_{k\eps n}
		 \tfrac{|V_b\cap (U_{(k+1)\eps n}\setminus U_{k\eps n})|}{|V_b\setminus U_{(k+1)\eps n}|}\right)\\
	& = \big|X\big| \cdot \tfrac{|V_b\setminus U_{s\eps n}|}{|V_b\setminus U_{0}|} \cdot \prod_{k=0}^{s-1}
		 \left( 1 \pm 50L\beta_{k\eps n}
		 \tfrac{|U_{(k+1)\eps n}\setminus U_{k\eps n}|}{|V_b\setminus U_{(k+1)\eps n}|}\right).
\end{align*}
Taking logs, we obtain
\begin{align*}
 \log \big|X\setminus U_{s\eps n}\big|&=\log\big|X\big|+
 \log\tfrac{|V_b\setminus U_{s\eps n}|}{|V_b\setminus U_{0}|}+ \sum_{k=0}^{s-1} \log \big(1\pm \tfrac{50L\beta_{k\eps n}|U_{(k+1)\eps n}\setminus U_{k\eps n}|}{|V_b\setminus U_{(k+1)\eps n}|}\big)\\
 &=\log\big|X\big|+\log\tfrac{|V_b\setminus U_{s\eps n}|}{|V_b\setminus U_{0}|}\pm\sum_{k=0}^{s-1}\tfrac{100L\beta_{k\eps n}|U_{(k+1)\eps n}\setminus U_{k\eps n}|}{|V_b\setminus U_{(k+1)\eps n}|}
\end{align*}
where the last line uses the approximation to the logarithm $\log(1\pm x)=\pm 2x$, valid for all sufficiently small $x$. In this case the term is sufficiently small because $\big|U_{(k+1)\eps n}\setminus U_{k\eps n}|\le L \eps n$, because $\eps$ is sufficiently small, and because $|V_b\setminus U_t|\ge \nu n$. Now we have
\begin{align*}
\sum_{k=0}^{s-1}\tfrac{100L\beta_{k\eps n}|U_{(k+1)\eps n}\setminus U_{k\eps n}|}{|V_b\setminus U_{(k+1)\eps n}|}
& \le
\frac{100L}{\nu n}
\sum_{k=0}^{s-1} \beta_{k\eps n} |U_{(k+1)\eps n}\setminus U_{k\eps n}|
\le
\frac{100L}{\nu n}
\sum_{k=0}^{s-1} \beta_{k\eps n} L\eps n  \\
& \le \frac{100L^2}{\nu n}
\int_{-\infty}^{s\eps n}\beta_x\mathrm{d}x
\le \tfrac{\nu^2 t^\ast}{10 n}\beta_{s\eps n}
\le 0.1 \nu^2 \beta_{s\eps n}\,,
\end{align*}
where the third inequality uses the fact that $\beta_x$ is increasing in $x$ and nonnegative. We obtain
\begin{align*}
 \log\big|X\setminus U_{s\eps n}\big|&=\log\big|X\big|+\log\tfrac{|V_b\setminus U_{s\eps n}|}{|V_b\setminus U_{0}|}\pm\tfrac14\beta_{s\eps n}
 =\log\big|X\big|+\log\tfrac{|V_b\setminus U_{s\eps n}|}{|V_b\setminus U_{0}|}+\log\big(1\pm\tfrac12\beta_{s\eps n}\big)\,,
\end{align*}
hence $\big|X\setminus U_{s\eps n}\big|=(1\pm\tfrac12\beta_{s\eps n})|X|\tfrac{|V_b\setminus U_{s\eps n}|}{|V_b\setminus U_{0}|}$.
Since $s\eps n\le t\le (s+1)\eps n$ and since $\beta_x$ is increasing in $x$, we have
\[\big|X\setminus U_{t}\big|=(1\pm\tfrac12\beta_{t})|X|\tfrac{|V_b\setminus U_{t}|\pm L\eps n}{|V_b\setminus U_{0}|}\pm L\eps n\,,\]
which by choice of $\eps$ gives the desired conclusion.
\end{claimproof}

If at some time $t$ Algorithm~\ref{alg:path-embed} gives us a $U_t$ such that the $(\beta_t,2)$-block-diet condition fails for $(H,U_t)$, then there is a first such time $t$. Suppose that $b\in\{\boxminus,\boxplus\}$ and $S$ (which is a set of either one or two vertices) is a witness of the failure. Because $(H,U_0)$ is $(\gamma,2)$-block-diet, we have
\[\big|\NBH_H(S)\setminus U_0\big|=(1\pm\gamma)d_{\boxminus b}^{i_1}d_{\boxplus b}^{i_2}|V_b\setminus U_0|\]
where $i_1$ and $i_2$ are the number of vertices of $S$ in $\boxminus$ and $\boxplus$ respectively. Now by Claim~\ref{cl:RPE:longint}, with probability at least $1-\exp\big(-n^{0.4}\big)$ we have
\[\big|\NBH_H(S)\setminus U_t\big|=\big(1\pm\tfrac34\beta_t\big)(1\pm\gamma)d_{\boxminus b}^{i_1}d_{\boxplus b}^{i_2}|V_b\setminus U_t|\]
which, since $\beta_t\ge10\gamma$, is a contradiction to $b$ and $S$ witnessing failure of $(\beta_t,2)$-block-diet.

Taking the union bound over $1\le t\le t^*$ and witnesses to failure of block-diet (i.e.\ singletons or pairs of vertices and $b\in\{\boxplus,\boxminus\}$), the probability that such a first time $t$ exists is by Claim~\ref{cl:RPE:longint} at most $n^4\exp\big(-n^{0.4}\big)$. In other words, with probability at least $1-\exp\big(-n^{0.3}\big)$, for each $1\le t\le t^*$ the pair $(H,U_t)$ is $(\beta_t,2)$-block-diet (and so Algorithm~\ref{alg:path-embed} does not fail). This establishes~\ref{RPE:nofail} and by Claim~\ref{cl:RPE:longint} also~\ref{RPE:imageeating}.

Next, we find the probability that any one of a given small collection of vertices is embedded to in some interval of time by Algorithm~\ref{alg:path-embed}. To begin with, we deal with intervals of length exactly $\eps n$.

\begin{claim}\label{cl:RPE:shortprob} Given any $1\le t\le t^\ast-\eps n$, let $\hist_{t-1}$ be a history of Algorithm~\ref{alg:path-embed} up to and including the embedding of $P_{t-1}$. Suppose $\hist_{t-1}$ is such that the pair $(H,U_{t-1})$ is $(\beta_{t-1},2)$-block-diet. Let for each $b\in\{\boxminus,\boxplus\}$ the set $S_b\subset V_b\setminus U_{t-1}$ contain at most $L$ elements. Then we have
\[\Prob\big[(S_\boxminus\cup S_\boxplus)\cap  U_{t+\eps n-1}\neq\emptyset\big|\hist_{t-1}\big]=\big(1\pm 30L\beta_{t-1}\big)\sum_{b\in\{\boxminus,\boxplus\}}|S_b|\tfrac{|V_b\cap (U_{t+\eps n-1}\setminus U_{t-1})|}{|V_b\setminus U_{t-1}|}\,.\]
\end{claim}
\begin{claimproof}[Proof of Claim~\ref{cl:RPE:shortprob}] 
Fix a history $\hist_{t-1}$ as in the claim statement; throughout we condition on $\hist_{t-1}$. Let $S=S_\boxminus\cup S_\boxplus$.
We write
\begin{align*}
 \Prob\big[|S \cap U_{t+\eps n-1}|\ge1\big]&=\sum_{t'=t}^{t+\eps n-1}\Prob\big[|S\cap\phi\big(V(P_{t'})\big)|\ge1\text{ and }S\cap\im\psi_{t'-1}=\emptyset\big]\\
 &=\sum_{t'=t}^{t+\eps n - 1}\sum_{\psi_{t'-1}} \Prob\big[|S\cap\phi\big(V(P_{t'})\big)|\ge1\big|\psi_{t'-1}\big]\cdot\Prob[\psi_{t'-1}]\,,
\end{align*}
where the sum over $\psi_{t'-1}$ runs only over embeddings where $S\cap\im\psi_{t'-1}=\emptyset$.  Note that since by choice of $\hist_{t-1}$ the pair $(H,U_{t-1})$ is $(\beta_{t-1},2)$-block-diet, it follows as before that $(H,U_{t'-1})$ is $(2\beta_{t-1},2)$-block-diet for each $t'$ we consider, and in particular the embeddings we assume exist indeed do.
Observe that the expected number of vertices of $P_{t'}$ embedded to $S$ is
\[\sum_{j=1}^{2L}\Prob\big[|S\cap\phi\big(V(P_{t'})\big)|\ge j\,\big|\psi_{t'-1}\big]\,,\]
and rearranging this we get
\[\Prob\big[|S\cap\phi\big(V(P_{t'})\big)|\ge 1\big|\psi_{t'-1}\big]=\Exp\big[\big|S\cap\phi\big(V(P_{t'})\big)\big|\,\big|\psi_{t'-1}\big]\pm |S|\Prob\big[|S\cap\phi\big(V(P_{t'})\big)|\ge 2\big|\psi_{t'-1}\big]\,.\]

For the error term, note that we have $|S|\le 2L$. Moreover, for the expectation above, we have
\[
\Exp\big[\big|S\cap\phi\big(V(P_{t'})\big)\big|\,\big|\psi_{t'-1}\big]
=
\sum_{b\in\{\boxminus,\boxplus\}}\quad \sum_{u\in S_b}
\Prob\big[ u\in\phi\big(V(P_{t'})\big)\,\big|\psi_{t'-1} \big]
\]
where for $u\in S_b$, we have that
$\Prob\big[ u\in\phi\big(V(P_{t'})\big)\,\big|\psi_{t'-1} \big]$ is the (conditional) probability
that a vertex $y\in V(P_{t'})\setminus A$ with $\chi(y)=V_b$ is embedded to $u$.
By Claim~\ref{cl:RPE:probemb} and the definition of $p_{u,y}$ this probability is 
$(1\pm10L\beta_{t-1})\sum p_{u,y}$
where the sum runs over all 
$y\in V(P_{t'})\setminus A$ with $\chi(y)=V_b$. Therefore,
putting everything together, we conclude
\begin{multline*}
 \Prob\big[|S\cap\phi\big(V(P_{t'})\big)|\ge1\big|\psi_{t'-1}\big]\\
 =(1\pm10L\beta_{t-1})\Bigg(\sum_{b\in\{\boxminus,\boxplus\}}\quad \sum_{u\in S_b}\quad \sum_{\substack{y\in V(P_{t'})\setminus A\\ \chi(y)=V_b}}p_{u,y}\Bigg)\pm 2L\Prob\Big[\big|S\cap \phi\big(V(P_{t'})\big)\big|\ge 2\Big|\psi_{t'-1}\Big]\,.
\end{multline*}
We can bound above $2L\Prob\big[\big|S\cap\phi\big(V(P_{t'})\big)\big|\ge2\big|\psi_{t'-1}\big]$ by $20L^5\nu^{-L}n^{-2}<n^{-1.5}$, by similar logic as in Claim~\ref{cl:RPE:probemb}. Since the $p_{u,y}$ do not depend on $\psi_{t'-1}$, we get
\[\Prob\big[|S\cap U_{t+\eps n-1}|\ge1\big]=(1\pm10L\beta_{t-1})\Big(\sum_{u\in S}p^*_{u,t}\pm 2n^{-0.5}\Big)\Prob[S\cap\im\psi_{t'-1}=\emptyset]\,.\]
This expression is trivially bounded above by $(1+10L\beta_{t-1})\big(\sum_{u\in S} p^*_{u,t}+2n^{-0.5}\big)<8L^2\eps\nu^{-1}$, so $\Prob[S\cap \im\psi_{t'-1} = \emptyset]=1\pm8L^2\eps\nu^{-1}$ and we get
\begin{align*}
 \Prob\big[|S\cap U_{t+\eps n-1}|\ge1\big]&=(1\pm10L\beta_{t-1})\Big(\sum_{u\in S}p^*_{u,t}\pm 2n^{-0.5}\Big)\cdot\big(1\pm8L^2\eps\nu^{-1}\big)\\
 &=(1\pm 20L\beta_{t-1})\sum_{u\in S}p^*_{u,t}
 =(1\pm 30L\beta_{t-1})\sum_{b\in\{\boxminus,\boxplus\}}|S_b|\tfrac{|X_{b,t}^{\textrm{int}}|}{|V_b\setminus U_{t-1}|}\,,
\end{align*}
where the last equality uses Claim~\ref{cl:RPE:pstarvt}. Since $|X_{b,t}^{\textrm{int}}|=\big|V_b\cap (U_{t+\eps n-1}\setminus U_{t-1})\big|$, we get the claimed probability.
\end{claimproof}

Using this repeatedly we can deduce similar estimates for general intervals.

\begin{claim}\label{cl:RPE:longprob} Given any $1\le t<t'\le t^*$, let $\hist_{t-1}$ be a history of Algorithm~\ref{alg:path-embed} up to and including embedding $P_{t-1}$. 
Suppose $\hist_{t-1}$ is such that the pair $(H,U_{t-1})$ is $(\beta_{t-1},2)$-block-diet.
Let for each $b\in\{\boxminus,\boxplus\}$ the set $S_b\subset V_b\setminus U_{t-1}$ contain at most $L$ elements. Then we have
\begin{multline*}
 \Prob\big[(S_\boxminus\cup S_\boxplus)\cap U_{t'-1}=\emptyset\big|\hist_{t-1}\big]\\
 =\big(1\pm2\beta_{t'-1}\big)\Bigg(\prod_{a\in\{\boxminus,\boxplus\}}\Big(\tfrac{|V_a\setminus U_{t'-1}|}{|V_a\setminus U_{t-1}|}\Big)^{|S_a|}\Bigg)\pm \eps^{-1}\Prob[\text{block-diet fails}|\hist_{t-1}]\,,
\end{multline*}
where the last probability refers to the event that, at some time $t''$ between $t$ and $t^*$, the pair $(H,U_{t''-1})$ is not $(\beta_{t''-1},2)$-block-diet.
\end{claim}
\begin{claimproof}[Proof of Claim~\ref{cl:RPE:longprob}]
We condition throughout on $\hist_{t-1}$ as given in the claim statement. We allow an error term of size $\eps^{-1}\Prob[\text{block-diet fails}|\hist_{t-1}]$, so we can throughout assume that block-diet does not fail (for the purposes of applying Claim~\ref{cl:RPE:shortprob} at most $\eps^{-1}$ times) since the errors caused by its failure are subsumed in this error term. Because $\eps$ is small, it is enough to prove the desired statements with relative error $1\pm99L\beta_{t'-1}$ under the additional assumption that $t'-t$ is an integer multiple of $\eps n$ (the change in $\beta_{t'-1}$ and the change in the main term caused by increasing $t'$ by a further $\eps n$ is covered by the increase in the relative error). So we may assume $t'=t+k\eps n$ for some integer $k\ge1$. We write $S=S_\boxminus\cup S_\boxplus$.
By Claim~\ref{cl:RPE:shortprob}, we have
\begin{align*}
\Prob\big[S\cap U_{t'-1}=\emptyset\big]
 & =\prod_{\ell=0}^{k-1}
 	\left(
 	1 - \big(1 \pm 30L\beta_{t+\ell\eps n-1} \big)
 	\sum_{a\in\{\boxminus,\boxplus\}}|S_a|\tfrac{|X_{a,t+\ell \eps n}^\mathrm{int}|}{|V_a\setminus U_{t+\ell\eps n-1}|}\right)\\
 & =\prod_{\ell=0}^{k-1}
 	(1 \pm 65L^2\eps\nu^{-1}\beta_{t+\ell\eps n-1} )
 	\left(
 	1 - \sum_{a\in\{\boxminus,\boxplus\}}|S_a|\tfrac{|X_{a,t+\ell \eps n}^\mathrm{int}|}{|V_a\setminus U_{t+\ell\eps n-1}|} \right)\\
 & =\prod_{\ell=0}^{k-1}
 	(1 \pm 70L^2\eps\nu^{-1}\beta_{t+\ell\eps n-1} )
 	\prod_{a\in\{\boxminus,\boxplus\}}\left(
 	1 - \tfrac{|X_{a,t+\ell \eps n}^\mathrm{int}|}{|V_a\setminus U_{t+\ell\eps n-1}|} \right)^{|S_a|}\\
 & =\prod_{\ell=0}^{k-1}\big(1\pm70L^2\eps\nu^{-1}\beta_{t+\ell\eps n-1}\big)\prod_{a\in\{\boxminus,\boxplus\}}\Big(\tfrac{|V_a\setminus U_{t+(\ell+1)\eps n-1}|}{|V_a\setminus U_{t+\ell\eps n-1}|}\Big)^{|S_a|}\\
  & =\Bigg(\prod_{\ell=0}^{k-1}\big(1\pm70L^2\eps\nu^{-1}\beta_{t+\ell\eps n-1}\big)\Bigg)\prod_{a\in\{\boxminus,\boxplus\}}\Big(\tfrac{|V_a\setminus U_{t+k\eps n-1}|}{|V_a\setminus U_{t-1}|}\Big)^{|S_a|}\,.
\end{align*}
To see that the second line is valid, we use $|X_{a,t+\ell \eps n}^\mathrm{int}|\le\eps n$ and $|V_a\setminus U_{t+\ell\eps n-1}|\ge\nu n$. For the third line, observe that expanding out all the brackets in the product on the third line, we get the brackets of the second line plus a collection of at most $2^{2L}$ products in which $\tfrac{|X_{a,t+\ell \eps n}^\mathrm{int}|}{|V_a\setminus U_{t+\ell\eps n-1}|}$ terms appear at least twice. Since $|X_{a,t+\ell \eps n}^\mathrm{int}|\le\eps n$, by choice of $\eps$ all of these terms are tiny compared to $\beta_{0}$ and hence are accounted for by the increase in the error term. The fourth and fifth lines are straightforward equalities.

It remains to show that the product of error terms is sufficiently close to $1$.
For this, note that, using a Taylor expansion, $\log(1\pm x)=\pm x +O(x^2)$ for $x\rightarrow 0$,
and hence $\log(1\pm x)=\pm 1.1x$
for sufficiently small $x$.
Thus, taking logs and using that $\eps$ is sufficiently small, we have
\[\log\prod_{\ell=0}^{k-1}\big(1\pm70L^2\eps\nu^{-1}\beta_{t+\ell\eps n-1}\big)=\pm\sum_{\ell=0}^{k-1}80L^2\eps\nu^{-1}\beta_{t+\ell\eps n-1}=\pm \int_{x=-\infty}^{t'}\tfrac{80L^2\nu^{-1}\beta_{x}}{n}\mathrm{d}x<\beta_{t'-1}\]
where the final equality uses the definition of $\beta_x$ and the fact $t^*\le n$. We obtain
\[\prod_{\ell=0}^{k-1}\big(1\pm70L^2\eps\nu^{-1}\beta_{t+\ell\eps n-1}\big)=\exp\big(\pm\beta_{t'-1}\big)=1\pm2\beta_{t'-1}\,,\]
as desired.
\end{claimproof}

From this point we will no longer require the detailed definition of $\beta_{t-1}$; it will suffice that $\beta_{t-1}\le\beta_{t^*}\le \sqrt{C}\gamma$, and we will simply replace error bounds from the earlier claims with the latter. The remaining parts of Lemma~\ref{lem:RPE} are relatively straightforward at this point. 

As noted earlier, the probability of block-diet failing is
at most $\exp\big(-n^{0.3}\big)$. Now,
taking in particular the case $t=1$ and $t'=t^*$ of Claim~\ref{cl:RPE:longprob}, and observing that $\big|\big\{x\in V(F)\setminus A:\chi(x)=V_a\big\}\big|=|V_a\cap(U_{t^*-1}\setminus U_0)|\pm L$, we obtain~\ref{RPE:probvtx}.

For~\ref{RPE:probedgefree}, fix $a,b\in\{\boxminus,\boxplus\}$ and let $e=uv$ with $u \in V_a$ and $v \in V_b$. It is convenient to give to each directed edge $(y,z)$ of $F$ a weight, which will relate to the probability that $y$ is embedded to $u$ and $z$ to $v$. We do this as follows. If one of $y$ and $z$ is in $A$, we set $q_{y,z}=0$. Otherwise, suppose $x,y,z,w$ are adjacent in $F$ in that order, and suppose $\chi(x)=V_c$, $\chi(y)=V_a$, $\chi(z)=V_b$ and $\chi(w)=V_d$, where $c,d\in\{\boxminus,\boxplus\}$ and $a,b$ are as above.
If $x$ is in $A$, and $\phi_0(x)u$ is an edge of $H$, we set $q_{y,z}=\tfrac{1}{d_{ac}}$. If $w$ is in $A$, and $\phi_0(w)v$ is an edge of $H$, we set $q_{y,z}=\tfrac{1}{d_{bd}}$. Note that we cannot have both $x$ and $w$ in $A$, since all paths have at least~4 edges. If neither $x$ nor $w$ is in $A$, set $q_{y,z}=1$, and in all other cases set $q_{y,z}=0$.

Fix $c\in\{\boxminus,\boxplus\}$ and consider the sum $\sum_{y,z}q_{y,z}$ where the sum runs over pairs $(y,z)$ such that $x\in A$ and $\chi(x)=V_c$. By the $\gamma$-anchor distribution property, there are
\[(1\pm\gamma)d_{ca}|A_{c,a,b}|\pm\tfrac12\gamma n^{0.99}\] summands equal to $\tfrac{1}{d_{ac}}$, and so this sum evaluates to
\[(1\pm\gamma)|A_{c,a,b}|\pm n^{0.999}\,.\]
Similarly, given $d\in\{\boxminus,\boxplus\}$, considering the sum over pairs $(y,z)$ such that $w\in A$ has $\chi(w)=V_d$, we have
\[\sum_{y,z}q_{y,z}=(1\pm\gamma)|A_{d,b,a}|\pm n^{0.999}\,.\]
Finally, the set $\big\{(y,z)\,:\,yz\in E(F),y\in X_a\setminus A, z\in X_b\setminus A\big\}$ splits up into $|A_{c,a,b}|$ pairs $(y,z)$ such that $x\in A$ and $\chi(x)=V_c$, $|A_{d,b,a}|$ pairs $(y,z)$ such that $w\in A$ and $\chi(w)=V_d$, for each $c,d\in\{\boxminus,\boxplus\}$, and the remaining pairs $(y,z)$ each of which has $q_{y,z}=1$. We conclude
\begin{equation}\label{eq:pathemb:pe4}
 \begin{split}
 \sum_{(y,z):yz\in E(F)}q_{y,z}&=(1\pm\gamma)\big|\big\{(y,z)\,:\,yz\in E(F),y\in X_a\setminus A, z\in X_b\setminus A\big\}\big|\pm 4 n^{0.999}\\
 &=(1\pm2\gamma)\big|\big\{(y,z)\,:\,yz\in E(F),y\in X_a\setminus A, z\in X_b\setminus A\big\}\big|\,,
 \end{split}
\end{equation}
where the second equality uses the lemma assumption that there are at least $\nu n$ pairs in the set on the right-hand side of~\eqref{eq:pathemb:pe4}.

We observe that
\[\Prob[\phi\text{ uses }e]=\sum_{(y,z):yz\in E(F)}\Prob[\phi(y)=u,\phi(z)=v]\]
and letting $t$ be such that $y,z\in V(P_t)$ we can write
\[\Prob[\phi(y)=u,\phi(z)=v]=\sum_{\psi_{t-1}}\Prob\big[\phi(y)=u,\phi(z)=v\big|\psi_{t-1}\big]\Prob[\psi_{t-1}]\,.\]
Observe that if $u$ or $v$ is in $\im\psi_{t-1}$, then $\Prob\big[\phi(y)=u,\phi(z)=v\big|\psi_{t-1}\big]=0$. Furthermore, the total probability of $\psi_{t-1}$ such that $(H,U_{t-1})$ is not $(\beta_{t-1},2)$-block-diet is at most $\exp\big(-n^{0.3}\big)$; so we are mainly interested in estimating $\Prob\big[\phi(y)=u,\phi(z)=v\big|\psi_{t-1}\big]$ when $\psi_{t-1}$ is $(\beta_{t-1},2)$-block-diet and $u,v\notin \im\psi_{t-1}$. We can do this by much the same approach as in Claim~\ref{cl:RPE:probemb}; we obtain
\[\Prob\big[\phi(y)=u,\phi(z)=v\big|\psi_{t-1}\big]=\tfrac{(1\pm 5L\beta_{t-1})q_{y,z}}{d_{ab}|V_a\setminus U_{t-1}||V_b\setminus U_{t-1}|}\,.\]

Summing over all $\psi_{t-1}$, we obtain
\begin{multline*}
  \Prob[\phi(y)=u,\phi(z)=v]\\
  \begin{aligned}
 &=\tfrac{(1\pm 5L\beta_{t-1})q_{y,z}}{d_{ab}|V_a\setminus U_{t-1}||V_b\setminus U_{t-1}|}\cdot\big(1\pm 100L\beta_{t-1}\big)\tfrac{|V_a\setminus U_{t-1}||V_b\setminus U_{t-1}|}{|V_a\setminus U_{0}||V_b\setminus U_{0}|}\pm 2\eps^{-1}\exp\big(-n^{0.3}\big)\\
 &=\big(1\pm200L\beta_{t-1}\big)\tfrac{q_{y,z}}{d_{ab}|V_a\setminus U_{0}||V_b\setminus U_{0}|}
  \end{aligned}
\end{multline*}
where the probability $\Prob[u,v\notin\im\psi_{t-1}]$ is estimated using Claim~\ref{cl:RPE:longprob}. Summing over all $(y,z)$ such that $yz\in E(F)$, and using~\eqref{eq:pathemb:pe4}, we obtain~\ref{RPE:probedgefree}. Note that we do not need to be careful estimating a sum of $\beta_t$ here; the error bounds are larger than any of the $\beta_t$ (and this will be the case for the next few properties, too).

For~\ref{RPE:probedgeside}, observe that we use $\phi_0(x)v$, where $x\in V(P_t)$ is in $X_a$, if and only if $v\notin\im\psi_{t-1}$ and then $v$ is selected as the vertex to which we embed the neighbour $y$ of $x$, where $y\in X_b$. The probability $\Prob[v\notin\im\psi_{t-1}]$ is $(1\pm100L\beta_{t-1})\tfrac{|V_b\setminus U_{t-1}|}{|V_b\setminus U_0|}\pm \eps^{-1}\exp\big(-n^{0.3}\big)$ by Claim~\ref{cl:RPE:longprob}. Among $\psi_{t-1}$ in this event, the probability that $(H,U_{t-1})$ is not $(\beta_{t-1},2)$-block-diet is at most $\exp\big(-n^{0.3}\big)$, and for the remaining $\psi_{t-1}$, the probability of $\psi(y)=v$, conditioned on $\psi_{t-1}$, is $(1\pm5L\beta_{t-1})\tfrac{1}{d_{ab}|V_b\setminus U_{t-1}|}$ according to
Claim~\ref{cl:RPE:probemb}. Putting these facts together we get~\ref{RPE:probedgeside}.

Finally, for~\ref{RPE:probsec}, we separate two cases. If $uv$ and $u'v'$ share a vertex, suppose $v=v'$. Observe that we use both $uv$ and $u'v$ if and only if there is some path $P_t$
which uses both edges. There are two subcases to consider. First, one or both of $u$ and $u'$ is in $\phi_0(A)$, in which case there is at most one possible choice of $t$ to consider, and second, neither is in $\phi_0(A)$, in which case we need to consider all paths $P_t$.

We deal with the first subcase first; assume without loss of generality $u\in\phi_0(A)$. If $\psi_{t-1}$ is such that $(H,U_{t-1})$ is $(\beta_{t-1},2)$-block-diet, then the probability of using both $uv$ and $u'v$ is easily estimated to be $O(n^{-2})$ (it is $\Theta(n^{-2})$ if neither $v$ nor $u'$ is in $\im\psi_{t-1}$, but this event need not occur). Moreover, the probability that 
$\psi_{t-1}$ fails to give block-diet is at most
$\exp\big(-n^{0.3}\big)$. Hence we conclude~\ref{RPE:probsec}.

For the second subcase, there are at most $n$ paths $P_t$ to consider and for each one, using both $uv$ and $u'v$ fixes three vertices in $\phi\big(V(P_t)\big)$. When $\psi_{t-1}$ is such that $(H,U_{t-1})$ is $(\beta_{t-1},2)$-block-diet, the probability of this is $O(n^{-3})$, so the probability a given $P_t$ uses $uv$ and $u'v$ is at most $n^{-2.5}$; taking the union bound over at most $n$ choices of $t$ we again get~\ref{RPE:probsec}.

If now $uv$ and $u'v'$ are disjoint edges, we have a few more subcases to consider. If both $uv$ and $u'v'$ intersect $\phi_0(A)$, then this fixes the path or paths that can use $uv$ and $u'v'$. If there is only one path (i.e.\ $uv$ and $u'v'$ between them contain both anchors of that path) then using $uv$ and $u'v'$ fixes two additional vertices on that path, and hence the probability of picking it is $O(n^{-2})$. If there are two such paths (i.e.\ $uv$ and $u'v'$ both contain an anchor but of different paths) then each path has one additional vertex on it fixed and again the probability of choosing such a pair of paths is $O(n^{-2})$. If $uv$ intersects $\phi_0(A)$, but $u'v'$ does not, then there is only one path which can contain $uv$ and doing so fixes one additional vertex on this path. While there are up to $n$ paths which might contain $u'v'$, for any given one of them to do so fixes two additional vertices on that path. By the union bound the probability of choosing any one such path is $O(n^{-1})$, and again the probability of choosing both $uv$ and $u'v'$ is $O(n^{-2})$. The probability of choosing $uv$ and $u'v'$ on the same path is $O(n^{-3})$ since doing so fixes three vertices on that unique path as $uv$ and $u'v'$ are disjoint. Finally, if $uv$ and $u'v'$ are both disjoint from $\phi_0(A)$, then either they can both be chosen on one path, or on two different paths. In the first case, there are $n$ possible paths and for any one to contain both $uv$ and $u'v'$ fixes four vertices, for a probability $O(n^{-3})$ by the union bound. In the second case, there are $n^2$ possible pairs of paths, on each of which two vertices are fixed for a probability $O(n^{-2})$ by the union bound. In all these cases we duly conclude~\ref{RPE:probsec}.
\end{proof}

\subsection{Packing path-forests}

Having analysed Algorithm~\ref{alg:path-embed}, we now state the complete path-forest packing algorithm, which simply runs Algorithm~\ref{alg:path-embed} repeatedly.

\begin{algorithm}[ht]\index{\PathPacking{}}
	\caption{\PathPacking{}}\label{alg:path-pack}
	\SetKwInOut{Input}{Input}
	\Input{
    $\bullet$ path-forests $F_1,\dots,F_{s^*}$, such that~$F_s$ has anchors~$A_s$, \\
    $\bullet$ a graph~$H_0$ on $\Vmin\dcup \Vplus$, \\
    $\bullet$ sets $U_s\subset V(H_0)$ of used vertices for $s\in[s^*]$, \\
    $\bullet$ assignments $\chi_s:V(F_s)\to\{\Vmin,\Vplus\}$, \\
    $\bullet$ embeddings $\phi_s:A_s\to V(H_0)$ such that
    $\im\phi_s\subset U_s$
  }
  Randomly permute the $F_s$ \;
  \For{$s=1$ \KwTo $s^*$}{
		run \RandomPathEmbedding($F_s$,$H_{s-1}$,$U_s$,$\chi_s$,$\phi_s$) to get
		an embedding~$\phi'_s$ of~$F_s$ into~$H_{s-1}$\;
		let $H_{s}$ be the graph obtained from $H_{s-1}$ by removing the
		edges of $\phi'_{s}(F_s)$\;
	}
\end{algorithm}
Abusing notation, we still refer to $F_s$ after the uniform random permutation of the path-forests; what we really mean is the $s$th forest in the random permutation, and we will continue this abuse of notation in our analysis. However the statements we make about the algorithm in Lemma~\ref{lem:pathpack} below are invariant under permutation.

We are now in a position to analyse Algorithm~\ref{alg:path-pack}, and give the important result of this section. Briefly, if we have a collection of anchored path-forests satisfying the conditions of Lemma~\ref{lem:RPE}, and in addition we have the following \emph{pair distribution property}, we can complete the analysis. We should stress that this property depends only weakly on the graph structure of $H$: in particular, if edges are removed from $H$ it cannot cause the pair distribution property to fail.

\begin{definition}[Pair distribution property]\label{def:pairdistprop}
  Given a graph $H$ with $V(H)=\Vmin\dcup\Vplus$ and a collection $F_1,\dots,F_{s^*}$ of path-forests all of whose paths have at least~4 edges, where for each $s\in[s^*]$ the set of endvertices of $F_s$ is $A_s$, and we have a used set $U_s\subset V(H)$ and an injective map $\phi_s:A_s\to U_s$. We say that a collection of maps $(\chi_s: V(F_s)\to\{\Vmin,\Vplus\})_{s\in [s^*]}$ has the \emph{$\gamma$-pair distribution property}\index{pair distribution property} if for each $a,b\in\{\boxminus,\boxplus\}$ and pair $u\in V_a$, $v\in V_b$ of distinct vertices of $H$ the following holds. Define $w_{uv;s}$ by\footnote{\label{foot:double}We emphasise that in the formulae below, a single edge $xy$ can be counted twice, as $(x,y)$ and $(y,x)$.}
\begin{subnumcases}{w_{uv;s}=}
  \label{case:w1}\textstyle
  \frac{|\{(x,y):xy\in E(F_s),\chi_s(x)= V_a,\chi_s(y)= V_b\text{ and }x,y\notin A_s\}|}{|V_a\setminus U_s||V_b\setminus U_s|} \text{ if $u,v\notin U_s$}\\
  \label{case:w2}\textstyle
  \frac{1}{|V_b\setminus U_s|}\text{ if $\phi_s^{-1}(u)=x\in A_s$, $v\notin U_s$, $xy\in E(F_s)$ with $\chi_s(y)=V_b$
  }\\
  \label{case:w3}\textstyle
  \frac{1}{|V_a\setminus U_s|}\text{ if $\phi_s^{-1}(v)=x\in A_s$, $u\notin U_s$, $xy\in E(F_s)$ with $\chi_s(y)=V_a$
  }
\end{subnumcases}
and otherwise $w_{uv;s}=0$. Then if $uv\in E(H)$ we have
 \begin{equation}\label{eq:pairdistrprop}
 \sum_{s=1}^{s^*}w_{uv;s}=(1\pm\gamma) \frac{\sum_{s=1}^{s^*}|\{(x,y):xy\in E(F_s),\chi_s(x)= V_a,\chi_s(y)= V_b\}|}{|V_a||V_b|} \pm\gamma n^{-0.01}\,.
  \end{equation}
\end{definition}

Before stating the lemma, we need a slightly simplified version of the set defined by index-quasirandomness. Given any sets $S_1,S_2\subset V(H)$, and any disjoint $T_1,T_2\subset [s^*]$, and any pairing $\Vmin=\{\boxminus_i\}_{i\in[n/2]}$ and $\Vplus=\{\boxplus_i\}_{i\in[n/2]}$, let \index{$\mathbb{U}_H(S_1,S_2,T_1,T_2)$}$\mathbb{U}_H(S_1,S_2,T_1,T_2)$ denote the set of $i\in[n/2]$ such that $\boxminus_i u\in E(H)$ for each $u\in S_1$, and $\boxplus_i u\in E(H)$ for each $u\in S_2$, and $\boxminus_i\notin U_t$ for each $t\in T_1$, and $\boxplus_i\notin U_t$ for each $t\in T_2$. The simplification as compared to the definition in~\eqref{eq:U5param} is that we have no set $T_3$, that is,
\begin{equation}\label{eq:U4U5}
\mathbb{U}_H(S_1,S_2,T_1,T_2)=\mathbb{U}_H(S_1,S_2,T_1,T_2,\emptyset)\;.
\end{equation}
We now have two definitions of the same symbol $\mathbb{U}_H$; since one is four-parameter and the other five-parameter, no confusion will arise.
We similarly define $\mathbb{U}'_{H'}(S_1,S_2,T_1,T_2)$ by replacing each $U_t$ with $U_t\cup\im\phi'_t$ (where $\phi'_t$ is the embedding of $F_t$ by Algorithm~\ref{alg:path-pack}) and $H$ with $H':=H_{s^*}$ the final graph left by Algorithm~\ref{alg:path-pack}.

\begin{lemma}[Path packing lemma]\label{lem:pathpack}
Given $\nu>0$ and $L\ge 5$, there exists a constant $C'$ such that for all $0<\gamma<1/C'$ the following holds.

Let $n$ be even and sufficiently large. Suppose that $H$ is a graph on $V_\boxminus\dcup V_\boxplus$, where $|V_\boxminus|=|V_\boxplus|=n/2$. 
Suppose $s^*\ge\nu n$. 
\begin{enumerate}[label=\itmarab{PA}]
\item\label{PA:BasicSetting}
For each $s\in [s^*]$, suppose $F_s$ is a path-forest with leaves (anchors) $A_s$ in which each path has between $5$ and $L$ vertices inclusive, and suppose $\phi_s:A_s\to V(H)$ is an embedding of the anchors. 
Suppose that $\chi_s:V(F_s)\to\{V_\boxminus,V_\boxplus\}
$ is an assignment of the vertices of $F_s$ to sides of $H$ such that $\phi_s(x)\in \chi_s(x)$ for each $x\in A_s$. Finally, suppose $U_s$ is a set of vertices containing $\phi_s(A_s)$. 
\item\label{PA:BlockQuasirandom}
Suppose that $H$ (together with the partition $V_\boxminus\dcup V_\boxplus$) is $(\gamma,L)$-block-quasirandom.
\item\label{PA:BlockDiet}
Suppose that for each $s\in[s^*]$ the pair $(H,U_s)$ is $(\gamma,L)$-block-diet.
\item\label{PA:AnchorDistribution} 
Suppose that for each $s\in[s^*]$ we have the $\gamma$-anchor distribution property of $(F_s,A_s,\phi_s,U_s)$ with respect to $H$.
\item\label{PA:Ineq}Suppose that for each $a,b\in\{\boxminus,\boxplus\}$ and each $s\in [s^*]$ we have 
\begin{align*}
|V_a\setminus U_s|-\big|\{x\in V(F_s)\setminus A_s:\chi_s(x)= V_a\}\big| &\ge\nu n\;,\\ \big|\chi_s^{-1}(\{V_a\})\setminus A_s\big|&\ge\nu n\;\mbox{, and}\\ \big|\big\{(x,y)\,:\,xy\in E(F_s),x\in
\chi_s^{-1}(\{V_a\})\setminus A_s,y\in \chi_s^{-1}(\{V_b\})\setminus A_s\big\}\big|&\ge \nu n
\;.
\end{align*}
\item\label{PA:PairDistributionProperty}
         Suppose that we have the $\gamma$-pair distribution property.   
\item\label{PA:LeftOverDensities}
For each $a\in\{\boxminus,\boxplus\}$ define
\begin{align*}
 d'_{aa}&:=\tfrac{1}{\binom{n/2}{2}}\Big(e_H(V_a)-\sum_{s=1}^{s^*}e_{F_s}\big(\chi_s^{-1}(\{V_a\})\big)\Big)\quad\text{and}\\
 d'_{\boxminus\boxplus}&:= \tfrac{4}{n^2}\Big(e_H(V_\boxminus,V_\boxplus)-\sum_{s=1}^{s^*}e_{F_s}\big(\chi_s^{-1}(\{V_\boxminus\}),\chi_s^{-1}(\{V_\boxplus\})\big)\Big)\,.
\end{align*}
	Suppose $d'_{\boxminus\boxminus},d'_{\boxminus\boxplus},d'_{\boxplus\boxplus}\ge\nu$.
\end{enumerate}

Let $\omega:V(H)\to[0,\nu^{-1}]$ be any weight function such that for each $a\in\{\boxminus,\boxplus\}$ we have $\sum_{u\in V_a} \omega(u)$ either equal to~0 or at least $\nu n$. Given any $S_1,S_2\subset V(H)$ and any disjoint $T_1,T_2\subset [s^*]$ such that $|S_1|,|S_2|,|T_1|,|T_2|\le L$, and any pairing $\Vmin=\{\boxminus_i\}_{i\in[n/2]}$ and $\Vplus=\{\boxplus_i\}_{i\in[n/2]}$, let $X\subset \mathbb{U}_H(S_1,S_2,T_1,T_2)$ be any set of size at least $\nu n$. Finally, given sets $S_\boxminus\subset \Vmin$ and $S_\boxplus\subset\Vplus$ with $|S_\boxminus|,|S_\boxplus|\le L$, let $Y$ be any subset of $\{s\in[s^*]:(S_\boxminus\cup S_\boxplus)\cap\im\phi_s=\emptyset\}$ such that $|Y|\ge\nu n$.

Set $d_{\boxminus\boxminus}:=d(H[V_\boxminus])$, $d_{\boxplus\boxplus}:=d(H[V_\boxplus])$, $d_{\boxminus\boxplus}:=d(H[V_\boxminus,V_\boxplus])$.

When we execute \PathPacking{} (Algorithm~\ref{alg:path-pack}), then
\begin{enumerate}[label=\itmarab{PP}]
 \item\label{pathpack:complete} with probability at least $1-\exp\big(-n^{0.2}\big)$, for each $1\le s\le s^*$ we obtain an embedding $\phi'_s$ of $F_s$ in $H\big[\big(V(H)\setminus U_s\big)\cup \phi_s(A_s)\big]$ extending $\phi_s$ and such that $\phi'_s(x)\in\chi_s(x)$ for each $x\in V(F_s)$. Furthermore, if we obtain such embeddings, the following hold for $H':=H_{s^*}$.
 \item\label{pathpack:pack} For each $uv\in E(H)$ there is at most one $s$ such that $\phi'_s$ uses the edge $uv$.
 \item\label{pathpack:quasi} 
 The graph $H'$ is $(C'\gamma,L)$-block-quasirandom on $V_\boxminus\dcup V_\boxplus$ with densities $d_{\boxminus\boxminus}',d_{\boxminus\boxplus}',d_{\boxplus\boxplus}'$.
 \item\label{pathpack:weights} With probability at least $1-\exp\big(-n^{0.2}\big)$, for each $a,b\in\{\boxminus,\boxplus\}$ and $u\in V_a$ we have
  \[\sum_{v\in \NBH_{H'}(u;V_b)}\omega(v)=\big(1\pm C'\gamma\big)\tfrac{d'_{ab}}{d_{ab}}\sum_{v\in \NBH_H(u;V_b)}\omega(v)\,.\]
 \item\label{pathpack:megaqr} With probability at least $1-\exp\big(-n^{0.2}\big)$ we have
 \begin{multline*}
  \big|X\cap \mathbb{U}'_{H'}(S_1,S_2,T_1,T_2)\big|
  =(1\pm C'\gamma)\Big(\frac{d'_{\boxminus\boxminus}}{d_{\boxminus\boxminus}}\Big)^{|S_1\cap \Vmin|}\Big(\frac{d'_{\boxminus\boxplus}}{d_{\boxminus\boxplus}}\Big)^{|S_1\cap \Vplus|+|S_2\cap\Vmin|}\Big(\frac{d'_{\boxplus\boxplus}}{d_{\boxplus\boxplus}}\Big)^{|S_2\cap \Vplus|}\\
  \cdot \Big(\prod_{t\in T_1}\frac{\big|\Vmin\setminus(U_t\cup\im\phi'_t)\big|}{\big|\Vmin\setminus U_t\big|}\Big)\Big(\prod_{t\in T_2}\frac{\big|\Vplus\setminus(U_t\cup\im\phi'_t)\big|}{\big|\Vplus\setminus U_t\big|}\Big) |X|\,.
 \end{multline*}
 \item\label{pathpack:imagecaps} With probability at least $1-\exp\big(-n^{0.2}\big)$ we have
  \[\big|\{s\in Y:(S_\boxminus\cup S_\boxplus)\cap\im\phi'_s=\emptyset\}\big|=(1\pm C'\gamma)\sum_{s\in Y}\prod_{a\in\{\boxminus,\boxplus\}}\Big(\tfrac{|V_a\setminus (U_s\cup\im\phi'_s)|}{|V_a\setminus U_s|}\Big)^{|S_a|}\,.\]
\end{enumerate}
\end{lemma}
The proof of this lemma mainly amounts to collecting the probabilistic estimates of Lemma~\ref{lem:RPE} and our assumptions on the quasirandom distribution of the anchor sets and used sets of the various path-forests to obtain expected values, then applying Corollary~\ref{cor:freedm} to show these expected values are likely all roughly correct. To begin with, we need to show that $H_{s-1}$ is likely to satisfy the required block-diet condition for Lemma~\ref{lem:RPE} with the used set of $F_s$ at each stage $s$ and the required anchor distribution property; for this purpose we will define an exponentially (in $s$) increasing error parameter $\alpha_s$, much as in~\cite{DegPack}. We should note that~\ref{pathpack:megaqr} holds for any pairing, and in particular it can be the case for some choices of $S_i,T_i$ that it holds vacuously (i.e.\ $\mathbb{U}_H(S_1,S_2,T_1,T_2)=\emptyset$). In our applications, we will fix a single pairing and it will be such that there are no vacuous cases.
\begin{proof}
Given $\nu$ and $L$, let $\nu_{\sublem{lem:RPE}}=\nu^{10L}$ and let $C$ be returned by Lemma~\ref{lem:RPE} for input $\nu_{\sublem{lem:RPE}}$ and $L$. We set
\[C'=\exp\big(2000CL\nu^{-20L}\big)\, .\]
Furthermore, for any $0<\gamma<1/C'$, we define 
\[\alpha_x:=10\gamma\exp\Big(\tfrac{1000CL\nu^{-20L}x}{s^*}\Big)\]
for $x\in\RR$, and we set $\eps = \gamma^2 (C')^{-1}$.

We begin with an observation which complements the assumption $s^*\ge\nu n$: the assumptions of the lemma bound $s^*$ from above by a linear function of~$n$ as well. Taking $a=b=\boxminus$ in the third inequality of~\ref{PA:Ineq}, and recalling that each edge of $F_s$ is counted there twice, we obtain $e_{F_s}\big(\chi_s^{-1}(\{\Vmin\})\big)\ge\tfrac{\nu n}{2}$ for each $s\in[s^*]$. On the other hand, $d'_{\boxminus\boxminus}\ge\nu>0$ in~\ref{PA:LeftOverDensities} forces $\sum_{s=1}^{s^*}e_{F_s}\big(\chi_s^{-1}(\{\Vmin\})\big)\le e_H(\Vmin)\le\binom{n/2}{2}$. Combining these two bounds we get
\begin{equation}\label{eq:pathpack:sstarupper}
\nu n\le s^*\le\tfrac{2}{\nu n}\tbinom{n/2}{2}\le\tfrac{n}{4\nu}\;.
\end{equation}

We first argue that the random permutation of the $F_s$ ensures that the pair distribution property holds not just globally but also on short intervals.
\begin{claim}\label{cl:pathpack:pairdist}
 With probability at least $1-\exp\big(-n^{0.5}\big)$, for each $1\le s\le s^*+1-\eps n$, after the uniform permutation, the path-forests $F_s,F_{s+1},\dots,F_{s+\eps n-1}$ have the $2\gamma$-pair distribution property.
\end{claim}
\begin{claimproof}[Proof of Claim~\ref{cl:pathpack:pairdist}]
 By the union bound, it is enough to show that with probability at least $1-\exp\big(-n^{0.55}\big)$ a uniform random choice of $\eps n$ path-forests has the $2\gamma$-pair distribution property, since each interval of $\eps n$ path-forests in the uniform random permutation has this distribution. For this, it is enough to show that a given pair of vertices $u\in V_a,v\in V_b$ with $uv\in E(H)$ is, with probability at least $1-\exp\big(-n^{0.6}\big)$, not a witness of failure of the $2\gamma$-pair distribution property.
 
We partition the set of path-forests according to the weight they assign to $uv$ (as in Definition~\ref{def:pairdistprop}), rounded down to the nearest multiple of $n^{-1.02}$. Since a path-forest has at most $n$ edges, and since $|V_\boxminus\setminus U|,|V_\boxplus\setminus U|\ge\nu n$ for $U$ the used set of any path-forest, the maximum weight assigned to $uv$ by any path-forest is $2\nu^{-2}n^{-1}$ and hence there are at most $3\nu^{-2}n^{0.02}$ parts in our partition. By Fact~\ref{fact:hypergeometricBasicProperties}, with probability at least $1-\exp\big(-n^{0.8}\big)$, in any given part $P$ with $p$ path-forests, the uniform random choice of $\eps n$ path-forests contains $\tfrac{\eps n}{s^*}p\pm n^{0.95}$ of them. Let $Y_P$ denote the random variable which is the sum of weights assigned to $uv$ by the path-forests of $P$ among the uniform random $\eps n$ path-forests; then assuming the above likely event occurs, we have
\[\big|Y_P-\Exp[Y_P]\big|\le n^{-1.02} p+2\nu^{-2} n^{-0.05}\,.\]
With probability at least $1-\exp\big(-n^{0.7}\big)$, by the union bound, the likely event occurs for each $P$, and in this case we therefore have
\[\Big|\sum_P (Y_P-\Exp[Y_P])\Big|\le n^{-1.02} s^*+2\nu^{-2}n^{-0.05}\cdot 3\nu^{-2}n^{0.02}\le 2 n^{-1.02} s^*\,,\]
where the final inequality is since $s^*\ge\nu n$.
Note that $\sum_P\Exp[Y_P]$ is simply $\tfrac{\eps n}{s^*}$ times the total weight assigned to $uv$ by all path-forests, which is as stated in Definition~\ref{def:pairdistprop}. 

By a similar argument we will next conclude that if there are in total $e$ edges between $V_a$ and $V_b$ in all path-forests, 
i.e.~edges $xy\in E(F_s)$ 
such that $\chi_s(x)=V_a$ and $\chi_s(y)=V_b$,
then with probability at least $1-\exp\big(-n^{0.7}\big)$, there are $\tfrac{\eps n}{s^*}e\pm 2n^{0.98}s^\ast$ edges between $V_a$ and $V_b$ in the uniform random $\eps n$ path-forests. Note that by assumption, we have $e\ge \nu^2n^2$,
as each of the $s^\ast\ge \nu n$ path-forests contributes at least $\nu n$ edges.
This time,
we partition the set of path-forests according to the number of edges between $V_a$ and $V_b$, rounded down to the nearest multiple of $n^{0.98}$.
Since a path-forest has at most $n$ edges, there are at most $n^{0.02}$ parts in our partition. Again, by Fact~\ref{fact:hypergeometricBasicProperties}, with probability at least $1-\exp\big(-n^{0.8}\big)$, in any given part $P'$ with $p$ path-forests, the uniform random choice of $\eps n$ path-forests contains $\tfrac{\eps n}{s^*}p\pm n^{0.95}$ of them. Let $Z_{P'}$ denote the random variable counting the number of edges between $V_a$ and $V_b$ in the path-forests of $P'$ among the uniform random $\eps n$ path-forests;
then assuming the above likely event occurs, we have
\[\big|Z_{P'}-\Exp[Z_{P'}]\big|\le n^{0.98} p+n^{1.95}\,.\]
With probability at least $1-\exp\big(-n^{0.7}\big)$, by the union bound, the likely event occurs for each $P'$, and in this case we therefore have
\[\Big|\sum_{P'} (Z_{P'}-\Exp[Z_{P'}])\Big|\le n^{0.98} s^* + n^{1.95}\cdot n^{0.02}\le 2 n^{0.98} s^* .\]
Note that $\sum_{P'}\Exp[Z_{P'}] = \tfrac{\eps n}{s^*}\cdot e$. 

If both this and the previous likely event occur, then we conclude 
\begin{align*}
\sum_{P}Y_P 
& = 
\tfrac{\eps n}{s^*}
 \sum_{s=1}^{s^*}w_{uv;s} \pm 2n^{-1.02}s^\ast \\
\JUSTIFY{by~\ref{PA:PairDistributionProperty}}
& =
\tfrac{\eps n}{s^*}
 \left( (1\pm\gamma) \frac{e}{|V_a||V_b|} \pm \gamma n^{-0.01} \right) \pm 2n^{-1.02}s^\ast \\
& =
(1\pm\gamma)
\frac{\sum_{P'}Z_{P'} \pm 2n^{0.98}s^\ast}{|V_a||V_b|} \pm \tfrac{\eps n}{s^*}\cdot \gamma n^{-0.01}  \pm 2n^{-1.02}s^\ast \\
\JUSTIFY{by~\eqref{eq:pathpack:sstarupper} and $s^\ast\ge \eps n$}& =
(1\pm\gamma)
\frac{\sum_{P'}Z_{P'}}{|V_a||V_b|} \pm 2\gamma n^{-0.01} .
\end{align*}
This means that $uv$ does not witness failure of the $2\gamma$-pair distribution property for the uniform random $\eps n$ path-forests, as desired.
\end{claimproof}

As remarked above, from this point when we write $F_s$ we mean the $s$th path-forest in the chosen permutation (and so on). We assume from this point that the likely event of Claim~\ref{cl:pathpack:pairdist} occurs.

For each $a\in\{\boxminus,\boxplus\}$ we define
\begin{align*}
 d_{aa;s}&:=d_{aa}-\binom{|V_a|}{2}^{-1}\sum_{i=1}^se\big(F_i\big[\chi_i^{-1}(\{V_a\})\big]\big)\quad\text{and}\\
 d_{\boxminus\boxplus;s}&:=d_{\boxminus\boxplus}-|V_\boxminus|^{-1}|V_{\boxplus}|^{-1}\sum_{i=1}^se\big(F_i\big[\chi_i^{-1}(\{V_\boxminus\}),\chi_i^{-1}(\{V_\boxplus\})\big]\big)\,.
\end{align*}
Note that by definition the $d_{ab;s}$ give the densities within and between the sides of $H_s$, assuming this graph is constructed by Algorithm~\ref{alg:path-pack}. Moreover, $d'_{ab}=d_{ab;s^*}$ for each $a,b\in\{\boxminus,\boxplus\}$.

To begin with, much as in the proof of Lemma~\ref{lem:RPE}, we define a quantity $p_{uv;s}$ which approximates the probability that $uv$ is used in the embedding of $F_s$, on the assumption that $uv\in E(H_{s-1})$ and that the conditions of Lemma~\ref{lem:RPE} are met when embedding $F_s$. Suppose $u\in V_a$ and $v\in V_b$, where $a,b\in\{\boxminus,\boxplus\}$.

If neither $u$ nor $v$ is in $U_s$, we set
\[p_{uv;s}:=\frac{|\{(x,y):xy\in E(F_s),\chi_s(x) = V_a,\chi_s(y)= V_b
~ \text{and} ~ x,y\notin A_s\}|}{d_{ab;s-1}|V_a\setminus U_s||V_b\setminus U_s|}\,.\]
If $u=\phi_s(x)$ for some $x\in A_s$, and $v\notin U_s$, and the neighbour $y$ of $x$ has $\chi_s(y)=V_b$, we set
\[p_{uv;s}:=\frac{1}{d_{ab;s-1}|V_b\setminus U_s|}\,,\]
and similarly if 
$v=\phi_s(x)$ for some $x\in A_s$, and $u\notin U_s$, and the neighbour $y$ of $x$ has $\chi_s(y)=V_a$, we set
\[p_{uv;s}:=\frac{1}{d_{ab;s-1}|V_a\setminus U_s|}\,.\]
In all other cases we set $p_{uv;s}:=0$. Note that Lemma~\ref{lem:RPE} parts~\ref{RPE:probedgefree} and~\ref{RPE:probedgeside} state that, conditioning on the running of Algorithm~\ref{alg:path-pack} up to creating $H_{s-1}$ and on the assumption that the conditions of Lemma~\ref{lem:RPE} are satisfied at this point, if $uv$ is any edge of $H_{s-1}$ then the probability that the embedding of $F_s$ uses $uv$ is $(1\pm C\alpha_{s-1}\big)p_{uv;s}$. We also note that, since $d_{ab;s-1}\ge\nu$ (which follows from $d'_{ab} \ge \nu$ in~\ref{PA:LeftOverDensities}) and $|V_a\setminus U_s|,|V_b\setminus U_s|\ge \nu n$ (by~\ref{PA:Ineq}),
and since $e(F_s)\le n$, we always have $p_{uv;s}\le\nu^{-3}n^{-1}$.

Using the pair distribution property~\ref{PA:PairDistributionProperty}, we can estimate the sum $\sum_{i=s+1}^{s+\eps n}p_{uv;i}$ for any $0\le s\le s^*-\eps n$ and any pair $uv\in E(H)$; one should think of this as being the probability that $uv$ is used in the embedding of some graph in this interval, assuming it is in $H_s$ and assuming the process is well-behaved.

\begin{claim}\label{cl:pathpack:sump}
 For each $0\le s\le s^*-\eps n$, each $a,b\in\{\boxminus,\boxplus\}$ and each pair of distinct vertices $u\in V_a$ and $v\in V_b$ with $uv\in E(H)$, we have
 \[\sum_{s'=s+1}^{s+\eps n}p_{uv;s'}=(1\pm3\gamma)\frac{d_{ab;s}-d_{ab;s+\eps n}}{d_{ab;s}} \pm3\gamma n^{-0.01}\,.\]
\end{claim}
\begin{claimproof}[Proof of Claim~\ref{cl:pathpack:sump}]
 Because $d_{ab;s^*}\ge\nu$, and because at most $\eps n^2$ edges are in the path-forests $F_{s+1},\dots,F_{s+\eps n}$, we have $d_{ab;s'}=(1\pm 10\eps\nu^{-1})d_{ab;s}$ for each $s+1\le s'\le s+\eps n$. Because this collection of path-forests has the $2\gamma$-pair distribution property (from~\ref{PA:PairDistributionProperty}), we have
 \[\sum_{s'=s+1}^{s+\eps n}p_{uv;s'}=\tfrac{(1\pm2\gamma)}{(1\pm 10\eps\nu^{-1})d_{ab;s}}\frac{\sum_{s=1}^{s^*}|\{(x,y):xy\in E(F_s),\chi_s(x) = V_a,\chi_s(y)= V_b\}|}{|V_a||V_b|} \pm2\gamma n^{-0.01}\,,\]
 and by choice of $\eps$ and definition of $d_{ab;s}$ we obtain the claim statement. Note that the big fraction in the equation above is exactly $\frac{d_{ab;s}-d_{ab;s+\eps n}}{d_{ab;s}}$ if $a\neq b$; if $a=b$ then it differs in that $|V_a|^2$ is not $2\binom{|V_a|}{2}$, but the ratio between these quantities is $1+O(1/n)$ as $n$ goes to infinity and this is absorbed to the error term.
\end{claimproof}

We now argue that neighbourhoods of sets of not too many vertices in $H_s$ into reasonably large sets decrease roughly as one would expect if edges were removed uniformly at random (rather than as path-forests). For purposes of establishing~\ref{pathpack:megaqr} we need to consider a fixed pairing $\Vmin=\{\boxminus_i:i\in[n/2]\}$ and $\Vplus=\{\boxplus_i:i\in[n/2]\}$. We will not need to assume any special properties of this pairing.
\begin{claim}\label{cl:pathpack:eat}
 Fix each $1\le s<s''\le s^*$, and run Algorithm~\ref{alg:path-pack} up to the point where $H_s$ is defined (i.e.\ $F_s$ has been embedded); condition on this $H_s$. Now choose sets $W_1,W_2\subset V(H)$ with $|W_1|,|W_2|\le L$ and a set $U\subset [n/2]$ with $|U|\ge\nu^{5L}n$ such that for each $i\in U$, $w\in W_1$ and $w'\in W_2$ we have $w\boxminus_i,w'\boxplus_i\in E(H_s)$. We now run Algorithm~\ref{alg:path-pack} further until $H_{s''}$ is defined. With probability at least $1-\exp\big(-n^{0.4}\big)$ one of the following occurs. First, we have
 \begin{equation}\label{eq:pathpack:eat:goal}
  \begin{split}
   \big|\big\{i\in U\,&:\,\forall w\in W_1,w'\in W_2\,\,\text{we have}\,\,w\boxminus_i,w'\boxplus_i\in E(H_{s''})\big\}\big|\\
   &=\big(1\pm\tfrac12\alpha_{s''}\big)|U|\Big(\tfrac{d_{\boxminus\boxminus;s''}}{d_{\boxminus\boxminus;s}}\Big)^{|W_1\cap\Vmin|}\Big(\tfrac{d_{\boxminus\boxplus;s''}}{d_{\boxminus\boxplus;s}}\Big)^{|W_1\cap\Vplus|+|W_2\cap\Vmin|}\Big(\tfrac{d_{\boxplus\boxplus;s''}}{d_{\boxplus\boxplus;s}}\Big)^{|W_2\cap\Vplus|}\,.
  \end{split}
 \end{equation}
 Second, there is some stage $s'$ between $s$ and $s''$ inclusive such that either the conditions of Lemma~\ref{lem:RPE} are not met, or its low-probability event occurs.
\end{claim}
\begin{claimproof}[Proof of Claim~\ref{cl:pathpack:eat}]
 Let $U_s=U$, and for each $s'>s$ let $U_{s'}$ be the set of $i\in U_s$ such that $w\boxminus_i,w'\boxplus_i\in E(H_{s'})$ for all $w\in W_1$ and $w'\in W_2$.
 Observe that it is enough to prove that a given $s''$ is unlikely to be the first $s''$ such that~\eqref{eq:pathpack:eat:goal} fails. So let $\cE$ be the event that the conditions of Lemma~\ref{lem:RPE} are met (which holds initially at $s=1$ due to~\ref{PA:BlockQuasirandom} and~\ref{PA:BlockDiet}), and the low-probability event of that lemma does not occur, and~\eqref{eq:pathpack:eat:goal} does not fail at any stage before $s''$. We aim to show that the probability that $\cE$ occurs and~\eqref{eq:pathpack:eat:goal} fails at stage $s''$ is at most $\exp\big(-n^{0.5}\big)$, since then taking the union bound over choices of $s''$ establishes the claim statement.

 Let $s<s'\le s''$, and condition for a moment on the history of Algorithm~\ref{alg:path-pack} up to and including the embedding of $F_{s'-1}$. Given $i\in U$ such that for each $w\in W_1$ and $w'\in W_2$ we have $w\boxminus_i,w'\boxplus_i\in E(H_{s'-1})$, let $E_i$ denote the edges $w\boxminus_i$ for $w\in W_1$ and $w'\boxplus_i$ for $w'\in W_2$. If $|E_i|=|W_1|+|W_2|$ we say $i$ is \emph{normal}, otherwise it is \emph{special}. Note that if $i$ is special, this means that $\boxminus_i\in W_2$ and $\boxplus_i\in W_1$; there are therefore at most $L$ special $i$.
 
 By Lemma~\ref{lem:RPE}, property \ref{RPE:probedgefree}, the conditional probability that the embedding of $F_{s'}$ uses any given $uv\in E_i$ is $(1\pm C\alpha_{s'-1})p_{uv;s'}$. Furthermore the conditional probability that any given pair of edges in $E_i$ are used is, by~\ref{RPE:probsec}, at most $n^{-3/2}$. It follows that, provided at least one $p_{uv;s'}$ is non-zero with $uv\in E_i$, the conditional probability that at least one edge of $E_i$ is used by the embedding of $F_{s'}$ is
 \[(1\pm C\alpha_{s'-1})\sum_{uv\in E_i} p_{uv;s'}\pm\binom{2L}{2}n^{-3/2}=(1\pm2C\alpha_{s'-1})\sum_{uv\in E_i} p_{uv;s'}\,,\]
 where we justify the equality as follows. If some $p_{uv;s'}$ is non-zero, then by definition of $p_{uv;s'}$ it is at least $\nu n^{-1}$ (here we use the assumption~\ref{PA:Ineq} that at least $\nu n$ edges of $F_{s'}$ are assigned within each of $\Vmin$, $\Vplus$ and between $\Vmin$ and $\Vplus$) and hence the error term, which is $O(n^{-3/2})$, can be absorbed in the relative error.
 If all $p_{uv;s'}$ with $uv\in E_i$ are zero, then the $\binom{2L}{2}n^{-3/2}$ error term can be removed from the left side of the above: the probability that an edge of $E_i$ is used is zero. Thus also in this case the right side is a valid estimate.
 
 By linearity of expectation, we see that the conditional expected number of $i\in U_{s'-1}$ such that at least one edge of $E_i$ is used in embedding $F_{s'}$, is
 \[\big(1\pm 2C\alpha_{s'-1}\big)\sum_{i\in U_{s'-1}}\sum_{uv\in E_i} p_{uv;s'}\,.\]
 
Now since embedding any one graph $F_{s'}$ uses at most two edges at any given vertex (in particular those in $W_1\cup W_2$), there are at most $4L$ values of $i$ such that an edge of $E_i$ is used in embedding $F_{s'}$. So Corollary~\ref{cor:freedm} tells us that with probability at most $\exp\big(-n^{0.5}\big)$ the event $\cE$ occurs and nevertheless we have
\begin{equation*}
\begin{split}\big|U_s\setminus U_{s''}\big|
\neq\sum_{s'=s+1}^{s''}\big(1\pm 3C\alpha_{s'-1}\big)\sum_{i\in U_{s'-1}}\sum_{uv\in E_i} p_{uv;s'}\pm\eps n\,.\end{split}
\end{equation*}
Observe that it is enough to show that within the event $\cE$ we have
\begin{equation}\label{eq:pathpack:weightsconc}
 \begin{split}
  \sum_{s'=s+1}^{s''}&\big(1\pm 3C\alpha_{s'-1}\big)\sum_{i\in U_{s'-1}}\sum_{uv\in E_i} p_{uv;s'}=\\
  &\Bigg(1-\big(1\pm\tfrac14\alpha_{s''}\big)\Big(\tfrac{d_{\boxminus\boxminus;s''}}{d_{\boxminus\boxminus;s}}\Big)^{|W_1\cap\Vmin|}\Big(\tfrac{d_{\boxminus\boxplus;s''}}{d_{\boxminus\boxplus;s}}\Big)^{|W_1\cap\Vplus|+|W_2\cap\Vmin|}\Big(\tfrac{d_{\boxplus\boxplus;s''}}{d_{\boxplus\boxplus;s}}\Big)^{|W_2\cap\Vplus|}\Bigg)|U|\,,
 \end{split}
\end{equation}
since by choice of $\eps$ and $\nu$ the error term $\eps n$ is absorbed into the larger error in~\eqref{eq:pathpack:eat:goal}.

Note that the left side of~\eqref{eq:pathpack:weightsconc} is a random variable; we will now argue that within $\cE$ it is surely in the claimed interval. To that end, we split the sum over $s'$ into disjoint intervals of length $\eps n$, together with a final interval of length less than $\eps n$. We deal first with the final interval: in this interval, we use at most $2\eps n$ edges at any given vertex of $W_1\cup W_2$, hence the sum over this interval is at most $(|W_1|+|W_2|)\cdot 2\eps n\le 4L\eps n$, which we can absorb in the error term. We can thus assume $s''-s$ is a multiple of $\eps n$, and we only need to show
\begin{align*}
 \sum_{s'=s+1}^{s''}&\big(1\pm 3C\alpha_{s'-1}\big)\sum_{i\in U_{s'-1}}\sum_{uv\in E_i} p_{uv;s'}\\
 &=|U|\Bigg(1-\Big(\tfrac{d_{\boxminus\boxminus;s''}}{d_{\boxminus\boxminus;s}}\Big)^{|W_1\cap\Vmin|}\Big(\tfrac{d_{\boxminus\boxplus;s''}}{d_{\boxminus\boxplus;s}}\Big)^{|W_1\cap\Vplus|+|W_2\cap\Vmin|}\Big(\tfrac{d_{\boxplus\boxplus;s''}}{d_{\boxplus\boxplus;s}}\Big)^{|W_2\cap\Vplus|}\Bigg)\pm\tfrac18\alpha_{s''}|U|\nu^{2L}\,.
\end{align*}

Consider the sum
\[\sum_{s'=s+1+\ell\eps n}^{s+(\ell+1)\eps n}\big(1\pm 3C\alpha_{s'-1}\big)\sum_{i\in U_{s'-1}}\sum_{uv\in E_i} p_{uv;s'}\]
where $0\le\ell\le\tfrac{s''-s}{\eps n}-1$. To estimate this sum, fix an $\ell$. We separate the main term
\[\sum_{s'=s+1+\ell\eps n}^{s+(\ell+1)\eps n}\ \sum_{i\in U_{s'-1}}\ \sum_{uv\in E_i} p_{uv;s'}\]
and the error term
\begin{equation}\label{eq:pathpack:eat:error}
 \sum_{s'=s+1+\ell\eps n}^{s+(\ell+1)\eps n}3C\alpha_{s'-1}\sum_{i\in U_{s'-1}}\sum_{uv\in E_i} p_{uv;s'}\le\sum_{s'=s+1+\ell\eps n}^{s+(\ell+1)\eps n}6C\alpha_{s'-1}|U|L\nu^{-3}n^{-1}\,,
\end{equation}
where the inequality uses that each $p_{uv;s'}$ is bounded above by $\nu^{-3}n^{-1}$ and $|E_i|\le 2L$.

To estimate the main term, observe that the set $U_{s'-1}$ decreases, as $s'$ ranges over the interval from $s+1+\ell\eps n$ to $s+(\ell+1)\eps n$, by at most $4L\eps n$ indices. Each normal $i$ which is in $U_{s+(\ell+1)\eps n}$ contributes, by Claim~\ref{cl:pathpack:sump},
\begin{multline*}
 (1\pm3\gamma)\Big(|W_1\cap\Vmin|\frac{d_{\boxminus\boxminus;s+\ell\eps n}-d_{\boxminus\boxminus;s+(\ell+1)\eps n}}{d_{\boxminus\boxminus;s+\ell\eps n}}+\big(|W_1\cap\Vplus|+|W_2\cap\Vmin|\big)\frac{d_{\boxminus\boxplus;s+\ell\eps n}-d_{\boxminus\boxplus;s+(\ell+1)\eps n}}{d_{\boxminus\boxplus;s+\ell\eps n}}\\
 +|W_2\cap\Vplus|\frac{d_{\boxplus\boxplus;s+\ell\eps n}-d_{\boxplus\boxplus;s+(\ell+1)\eps n}}{d_{\boxplus\boxplus;s+\ell\eps n}}\Big) \pm 6L \gamma n^{-0.01}
\end{multline*}
 to the sum. Each special $i$ (of which there are at most $L$) and each normal $i$ which is not in $U_{s+(\ell+1)\eps n}$ contributes at most this amount. Thus our main term is $|U_{s+\ell\eps n}|\pm(4L\eps n+L)$ times the above estimate. Since we are in $\cE$, we know the size of $U_{s+\ell\eps n}$; in particular it is of size at least $\tfrac12\nu^{7L}n$, and so the error $4L\eps n+L$ is tiny compared to $\alpha_0$ by choice of $\eps$. Putting this together, we have
\begin{multline}
  \sum_{s'=s+1+\ell\eps n}^{s+(\ell+1)\eps n}\sum_{i\in U_{s'-1}}\sum_{uv\in E_i} p_{uv;s'}\\
  \begin{aligned}
 \label{eq:pathpack:eat:shortsum}=&\big(1\pm\alpha_{s+\ell\eps n}\big)|U|\Big(\tfrac{d_{\boxminus\boxminus;s+\ell\eps n}}{d_{\boxminus\boxminus;s}}\Big)^{|W_1\cap\Vmin|}\Big(\tfrac{d_{\boxminus\boxplus;s+\ell\eps n}}{d_{\boxminus\boxplus;s}}\Big)^{|W_1\cap\Vplus|+|W_2\cap\Vmin|}\Big(\tfrac{d_{\boxplus\boxplus;s+\ell\eps n}}{d_{\boxplus\boxplus;s}}\Big)^{|W_2\cap\Vplus|}\\
 &\cdot\Bigg(|W_1\cap\Vmin|\tfrac{d_{\boxminus\boxminus;s+\ell\eps n}-d_{\boxminus\boxminus;s+(\ell+1)\eps n}}{d_{\boxminus\boxminus;s+\ell\eps n}}+\big(|W_1\cap\Vplus|+|W_2\cap\Vmin|\big)\tfrac{d_{\boxminus\boxplus;s+\ell\eps n}-d_{\boxminus\boxplus;s+(\ell+1)\eps n}}{d_{\boxminus\boxplus;s+\ell\eps n}}\\
 &+|W_2\cap\Vplus|\tfrac{d_{\boxplus\boxplus;s+\ell\eps n}-d_{\boxplus\boxplus;s+(\ell+1)\eps n}}{d_{\boxplus\boxplus;s+\ell\eps n}}\Bigg) \pm 6L \gamma n^{0.99} 
  \end{aligned}
\end{multline}

To understand the right-hand side of~\eqref{eq:pathpack:eat:shortsum}, consider expanding out
\begin{multline*}
\Big(\tfrac{d_{\boxminus\boxminus;s+\ell\eps n}}{d_{\boxminus\boxminus;s}}\Big)^{|W_1\cap\Vmin|}\Big(\tfrac{d_{\boxminus\boxplus;s+\ell\eps n}}{d_{\boxminus\boxplus;s}}\Big)^{|W_1\cap\Vplus|+|W_2\cap\Vmin|}\Big(\tfrac{d_{\boxplus\boxplus;s+\ell\eps n}}{d_{\boxplus\boxplus;s}}\Big)^{|W_2\cap\Vplus|}\\
\begin{aligned}
=&\Big(\tfrac{d_{\boxminus\boxminus;s+(\ell+1)\eps n}}{d_{\boxminus\boxminus;s}}+\tfrac{d_{\boxminus\boxminus;s+\ell\eps n}-d_{\boxminus\boxminus;s+(\ell+1)\eps n}}{d_{\boxminus\boxminus;s}}\Big)^{|W_1\cap\Vmin|}\\
\cdot &\Big(\tfrac{d_{\boxminus\boxplus;s+(\ell+1)\eps n}}{d_{\boxminus\boxplus;s}}+\tfrac{d_{\boxminus\boxplus;s+\ell\eps n}-d_{\boxminus\boxplus;s+(\ell+1)\eps n}}{d_{\boxminus\boxplus;s}}\Big)^{|W_1\cap\Vplus|+|W_2\cap\Vmin|}\\
\cdot &\Big(\tfrac{d_{\boxplus\boxplus;s+(\ell+1)\eps n}}{d_{\boxplus\boxplus;s}}+\tfrac{d_{\boxplus\boxplus;s+\ell\eps n}-d_{\boxplus\boxplus;s+(\ell+1)\eps n}}{d_{\boxplus\boxplus;s}}\Big)^{|W_2\cap\Vplus|}\,.
\end{aligned}
\end{multline*}
We will split the above expression into three terms: the first term will be obtained by multiplying all the left sides of the brackets; the second term will be given by the sum of the products which are obtained by multiplying
one term from the right sides of the brackets and the remainder from the left;
and the third term will consist of the leftover.

This way, the first term equals
\[\Big(\tfrac{d_{\boxminus\boxminus;s+(\ell+1)\eps n}}{d_{\boxminus\boxminus;s}}\Big)^{|W_1\cap\Vmin|}\cdot \Big(\tfrac{d_{\boxminus\boxplus;s+(\ell+1)\eps n}}{d_{\boxminus\boxplus;s}}\Big)^{|W_1\cap\Vplus|+|W_2\cap\Vmin|}\cdot \Big(\tfrac{d_{\boxplus\boxplus;s+(\ell+1)\eps n}}{d_{\boxplus\boxplus;s}}\Big)^{|W_2\cap\Vplus|}\, .\]
Because $d_{ab;s+(\ell+1)\eps n}\ge\nu$, and because the path-forests $F_{s+1+\ell\eps n},\dots,F_{s+(\ell+1)\eps n}$ have at most $\eps n^2$ edges in total, we have
that $d_{ab;s+(\ell+1)\eps n}$ and $d_{ab;s+\ell\eps n}$ differ at most by a factor $(1\pm 10\eps\nu^{-1})$. Therefore, the second term equals
\begin{multline*}
\Big(\tfrac{d_{\boxminus\boxminus;s+(\ell+1)\eps n}}{d_{\boxminus\boxminus;s}}\Big)^{|W_1\cap\Vmin|}\Big(\tfrac{d_{\boxminus\boxplus;s+(\ell+1)\eps n}}{d_{\boxminus\boxplus;s}}\Big)^{|W_1\cap\Vplus|+|W_2\cap\Vmin|}\Big(\tfrac{d_{\boxplus\boxplus;s+(\ell+1)\eps n}}{d_{\boxplus\boxplus;s}}\Big)^{|W_2\cap\Vplus|}\\ 
 \cdot\Bigg(|W_1\cap\Vmin|\tfrac{d_{\boxminus\boxminus;s+\ell\eps n}-d_{\boxminus\boxminus;s+(\ell+1)\eps n}}{d_{\boxminus\boxminus;s+(\ell+1)\eps n}}+\big(|W_1\cap\Vplus|+|W_2\cap\Vmin|\big)\tfrac{d_{\boxminus\boxplus;s+\ell\eps n}-d_{\boxminus\boxplus;s+(\ell+1)\eps n}}{d_{\boxminus\boxplus;s+(\ell+1)\eps n}}\\
+|W_2\cap\Vplus|\tfrac{d_{\boxplus\boxplus;s+\ell\eps n}-d_{\boxplus\boxplus;s+(\ell+1)\eps n}}{d_{\boxplus\boxplus;s+(\ell+1)\eps n}}\Bigg)\,,
\end{multline*}
which equals
\begin{multline}
\label{eq:pathpack:mixed_term}
(1\pm 30L\eps \nu^{-1}) \Big(\tfrac{d_{\boxminus\boxminus;s+\ell\eps n}}{d_{\boxminus\boxminus;s}}\Big)^{|W_1\cap\Vmin|}\Big(\tfrac{d_{\boxminus\boxplus;s+\ell\eps n}}{d_{\boxminus\boxplus;s}}\Big)^{|W_1\cap\Vplus|+|W_2\cap\Vmin|}\Big(\tfrac{d_{\boxplus\boxplus;s+\ell\eps n}}{d_{\boxplus\boxplus;s}}\Big)^{|W_2\cap\Vplus|}\\ 
\cdot\Bigg(|W_1\cap\Vmin|\tfrac{d_{\boxminus\boxminus;s+\ell\eps n}-d_{\boxminus\boxminus;s+(\ell+1)\eps n}}{d_{\boxminus\boxminus;s+\ell\eps n}}+\big(|W_1\cap\Vplus|+|W_2\cap\Vmin|\big)\tfrac{d_{\boxminus\boxplus;s+\ell\eps n}-d_{\boxminus\boxplus;s+(\ell+1)\eps n}}{d_{\boxminus\boxplus;s+\ell\eps n}}\\ 
+|W_2\cap\Vplus|\tfrac{d_{\boxplus\boxplus;s+\ell\eps n}-d_{\boxplus\boxplus;s+(\ell+1)\eps n}}{d_{\boxplus\boxplus;s+\ell\eps n}}\Bigg)\,,
\end{multline}
which is, except for the factors $(1\pm 30L\eps \nu^{-1})$ and $(1\pm \alpha_{s+\ell\eps n})|U|$, the same as the formula on the right-hand side of~\eqref{eq:pathpack:eat:shortsum}. 

For the third term observe that there are fewer than $2^{2L}$ summands each of which is a product of at least two terms from the right sides of the brackets.
Using again that $d_{ab;s+(\ell+1)\eps n}$ and $d_{ab;s+\ell\eps n}$ differ at most by a factor $(1\pm 10\eps\nu^{-1})$, we see that this third term
amounts to at most $2^{2L}(10\eps\nu^{-1})^2$. 

Putting this together, we see that the term in \eqref{eq:pathpack:mixed_term},
equals
\begin{multline*}
\Big(\tfrac{d_{\boxminus\boxminus;s+\ell\eps n}}{d_{\boxminus\boxminus;s}}\Big)^{|W_1\cap\Vmin|}\Big(\tfrac{d_{\boxminus\boxplus;s+\ell\eps n}}{d_{\boxminus\boxplus;s}}\Big)^{|W_1\cap\Vplus|+|W_2\cap\Vmin|}\Big(\tfrac{d_{\boxplus\boxplus;s+\ell\eps n}}{d_{\boxplus\boxplus;s}}\Big)^{|W_2\cap\Vplus|} \\
 -
\Big(\tfrac{d_{\boxminus\boxminus;s+(\ell+1)\eps n}}{d_{\boxminus\boxminus;s}}\Big)^{|W_1\cap\Vmin|}\cdot \Big(\tfrac{d_{\boxminus\boxplus;s+(\ell+1)\eps n}}{d_{\boxminus\boxplus;s}}\Big)^{|W_1\cap\Vplus|+|W_2\cap\Vmin|}\cdot \Big(\tfrac{d_{\boxplus\boxplus;s+(\ell+1)\eps n}}{d_{\boxplus\boxplus;s}}\Big)^{|W_2\cap\Vplus|}\\
 \pm
2^{2L}(10\eps\nu^{-1})^2\, .
\end{multline*}
Plugging this into \eqref{eq:pathpack:eat:shortsum}, we finally obtain that $\sum_{s'=s+1+\ell\eps n}^{s+(\ell+1)\eps n}\sum_{i\in U_{s'-1}}\sum_{uv\in E_i} p_{uv;s'}$ equals
\begin{multline*}
 \big(1\pm 2\alpha_{s+\ell\eps n}\big)|U|\Bigg(\Big(\tfrac{d_{\boxminus\boxminus;s+\ell\eps n}}{d_{\boxminus\boxminus;s}}\Big)^{|W_1\cap\Vmin|}\cdot \Big(\tfrac{d_{\boxminus\boxplus;s+\ell\eps n}}{d_{\boxminus\boxplus;s}}\Big)^{|W_1\cap\Vplus|+|W_2\cap\Vmin|}\cdot \Big(\tfrac{d_{\boxplus\boxplus;s+\ell\eps n}}{d_{\boxplus\boxplus;s}}\Big)^{|W_2\cap\Vplus|}\\
 -\Big(\tfrac{d_{\boxminus\boxminus;s+(\ell+1)\eps n}}{d_{\boxminus\boxminus;s}}\Big)^{|W_1\cap\Vmin|}\cdot \Big(\tfrac{d_{\boxminus\boxplus;s+(\ell+1)\eps n}}{d_{\boxminus\boxplus;s}}\Big)^{|W_1\cap\Vplus|+|W_2\cap\Vmin|}\cdot \Big(\tfrac{d_{\boxplus\boxplus;s+(\ell+1)\eps n}}{d_{\boxplus\boxplus;s}}\Big)^{|W_2\cap\Vplus|} \\
 \pm 2^{2L}(10\eps\nu^{-1})^2 \Bigg)\,,
\end{multline*}
which in turn equals
\begin{multline*}
 |U|\Bigg(\Big(\tfrac{d_{\boxminus\boxminus;s+\ell\eps n}}{d_{\boxminus\boxminus;s}}\Big)^{|W_1\cap\Vmin|}\cdot \Big(\tfrac{d_{\boxminus\boxplus;s+\ell\eps n}}{d_{\boxminus\boxplus;s}}\Big)^{|W_1\cap\Vplus|+|W_2\cap\Vmin|}\cdot \Big(\tfrac{d_{\boxplus\boxplus;s+\ell\eps n}}{d_{\boxplus\boxplus;s}}\Big)^{|W_2\cap\Vplus|}\\
 -\Big(\tfrac{d_{\boxminus\boxminus;s+(\ell+1)\eps n}}{d_{\boxminus\boxminus;s}}\Big)^{|W_1\cap\Vmin|}\cdot \Big(\tfrac{d_{\boxminus\boxplus;s+(\ell+1)\eps n}}{d_{\boxminus\boxplus;s}}\Big)^{|W_1\cap\Vplus|+|W_2\cap\Vmin|}\cdot \Big(\tfrac{d_{\boxplus\boxplus;s+(\ell+1)\eps n}}{d_{\boxplus\boxplus;s}}\Big)^{|W_2\cap\Vplus|}\Bigg) \\
 \pm 200\nu^{-4L}\eps\alpha_{s+\ell\eps n}|U|\,.
\end{multline*}
At last, we can sum this over $\ell$, together with the error term bounded in~\eqref{eq:pathpack:eat:error}, to obtain
\begin{multline*}
  \sum_{s'=s+1}^{s''}\big(1\pm 3C\alpha_{s'-1}\big)\sum_{i\in U_{s'-1}}\sum_{uv\in E_i} p_{uv;s'}\\
    =|U|\Bigg(\Big(\tfrac{d_{\boxminus\boxminus;s}}{d_{\boxminus\boxminus;s}}\Big)^{|W_1\cap\Vmin|}\cdot \Big(\tfrac{d_{\boxminus\boxplus;s}}{d_{\boxminus\boxplus;s}}\Big)^{|W_1\cap\Vplus|+|W_2\cap\Vmin|}\cdot \Big(\tfrac{d_{\boxplus\boxplus;s}}{d_{\boxplus\boxplus;s}}\Big)^{|W_2\cap\Vplus|}\\
 \hfill -\Big(\tfrac{d_{\boxminus\boxminus;s''}}{d_{\boxminus\boxminus;s}}\Big)^{|W_1\cap\Vmin|}\cdot \Big(\tfrac{d_{\boxminus\boxplus;s''}}{d_{\boxminus\boxplus;s}}\Big)^{|W_1\cap\Vplus|+|W_2\cap\Vmin|}\cdot \Big(\tfrac{d_{\boxplus\boxplus;s''}}{d_{\boxplus\boxplus;s}}\Big)^{|W_2\cap\Vplus|}\Bigg)\\
 \pm\sum_{s'=s}^{s''}(200\nu^{-4L}n^{-1}+6CL\nu^{-3}n^{-1})|U|\alpha_{s'-1}
\end{multline*}
and hence
\begin{multline*}
  \sum_{s'=s+1}^{s''}\big(1\pm 3C\alpha_{s'-1}\big)\sum_{i\in U_{s'-1}}\sum_{uv\in E_i} p_{uv;s'}\\
=|U|\Bigg(1-\Big(\tfrac{d_{\boxminus\boxminus;s''}}{d_{\boxminus\boxminus;s}}\Big)^{|W_1\cap\Vmin|}\cdot \Big(\tfrac{d_{\boxminus\boxplus;s''}}{d_{\boxminus\boxplus;s}}\Big)^{|W_1\cap\Vplus|+|W_2\cap\Vmin|}\cdot \Big(\tfrac{d_{\boxplus\boxplus;s''}}{d_{\boxplus\boxplus;s}}\Big)^{|W_2\cap\Vplus|}\Bigg) \\
 \pm \sum_{s'=s}^{s''} 200CL\nu^{-4L}n^{-1} \alpha_{s'-1} |U|\, .
\end{multline*}
For the last error term, we observe that, using~\eqref{eq:pathpack:sstarupper}, it can be bounded by
\[
\int_{x=-\infty}^{s''} 200CL\nu^{-4L}n^{-1} \alpha_{x} |U| \mathrm{d}x
=0.2\nu^{16L}s^\ast n^{-1} \alpha_{s''}|U|
<\tfrac18\alpha_{s''}|U|\nu^{2L}\, ,
\]
which completes the proof of the claim.
\end{claimproof}

This claim implies (we will see a proof later) that the block-quasirandomness and block-diet requirements of Lemma~\ref{lem:RPE} are satisfied at each stage. What we in addition need is the anchor distribution property, which we will show follows from the next claim.

\begin{claim}\label{cl:pathpack:weights}
 For each $1\le s\le s^*$, any $a,b\in\{\boxminus,\boxplus\}$ and $v\in V_a$, and any weight function $\omega:V(H_0)\to[0,\nu^{-1}]$ such that $\sum_{u\in \NBH_{H_0}(v;V_b)}\omega(u)\ge n^{0.95}$, the following experiment has probability at least $1-\exp\big(-n^{0.4}\big)$ of success. We run Algorithm~\ref{alg:path-pack} up to the point at which $\phi_s$ and $H_s$ are defined. We declare the experiment a success if either
 \begin{equation}\label{eq:pathpack:weightsum}
  \sum_{u\in \NBH_{H_s}(v;V_b)}\omega(u)=(1\pm \tfrac12\alpha_s)\tfrac{d_{ab;s}}{d_{ab}}\sum_{u\in \NBH_{H_0}(v;V_b)}\omega(u)\,,
 \end{equation}
 or at some stage $1\le s'<s$ the conditions of Lemma~\ref{lem:RPE} fail or the low-probability event of that lemma occurs.
\end{claim}
\begin{claimproof}[Proof of Claim~\ref{cl:pathpack:weights}]
Fix a weight function $\omega$, and abbreviate $W:=\sum_{u\in \NBH_{H_0}(v;V_b)}\omega(u)$. As with the previous claim, we show that any given $s$ has probability at most $\exp\big(-n^{0.5}\big)$ of being the first $s$ for which the claim statement fails; the claim then follows by the union bound over $s$. So fix $s$, and let $\cE$ be the event that at each stage from $1$ to $s$ inclusive the conditions of Lemma~\ref{lem:RPE} are met and its low-probability event does not occur, and in addition at each stage $1\le s'\le s-1$ we have
\[\sum_{u\in \NBH_{H_{s'}}(v;V_b)}\omega(u)=(1\pm \tfrac12\alpha_{s'})\tfrac{d_{ab;s'}}{d_{ab}}\sum_{u\in \NBH_{H_0}(v;V_b)}\omega(u)\,.\]

Note that Lemma~\ref{lem:RPE} states that within $\cE$ the probability that the embedding of $F_{s'}$ uses $uv$ is $\big(1\pm C\alpha_{s'-1}\big)p_{uv;s'}$. Let $\hist_{s'-1}$ denote the history of Algorithm~\ref{alg:path-pack} up to and including the embedding of $F_{s'-1}$. Thus by linearity of expectation, within $\cE$ we have
\[\Exp\left[\sum_{u\in \NBH_{H_{s'-1}}(v;V_b)\setminus \NBH_{H_{s'}}(v;V_b)}\omega(u) \Big| \hist_{s'-1}\right]=\big(1\pm C\alpha_{s'-1}\big)\sum_{u\in \NBH_{H_{s'-1}}(v;V_b)}\omega(u)p_{uv;s'}\,.\]

Observe that when we embed any given $F_{s'}$, we use at most two edges leaving $u$, so the maximum change to the weight of the neighbourhood of $u$ is at most $2\nu^{-1}$. Thus by Corollary~\ref{cor:freedm}, the probability that $\cE$ occurs and
\[\sum_{u\in \NBH_{H_0}(v;V_b)\setminus \NBH_{H_s}(v;V_b)}\omega(u)\neq \sum_{s'=1}^s\big(1\pm C\alpha_{s'-1}\big)\sum_{u\in \NBH_{H_{s'-1}}(v;V_b)}\omega(u)p_{uv;s'}\pm n^{0.9}\]
is at most $\exp\big(-n^{0.5}\big)$. As with the previous claim, it thus suffices to show that within $\cE$ we have
\[ \sum_{s'=1}^s\big(1\pm C\alpha_{s'-1}\big)\sum_{u\in \NBH_{H_{s'-1}}(v;V_b)}\omega(u)p_{uv;s'}=\big(1-(1\pm\tfrac14\alpha_s)\tfrac{d_{ab;s}}{d_{ab}}\big)\sum_{u\in \NBH_{H_0}(v;V_b)}\omega(u)\,.\]

As with the previous claim, note that the left side of this is a random variable which in $\cE$ we will see surely lies in the claimed interval; as there, we show this by splitting the sum over $s'$ into intervals of length $\eps n$ together with a final interval of length less than $\eps n$. In the final interval at most $2\eps n$ edges are removed from $v$, and the weight of vertices at these edges is thus at most $2\eps n\nu^{-1}$; this upper bounds the sum of expectations over this interval. Thus it is enough to assume $s$ is a multiple of $\eps n$ and show

\begin{equation}\label{eq:pathpack:weights:finsum}
 \sum_{s'=1}^s\big(1\pm C\alpha_{s'-1}\big)\sum_{u\in \NBH_{H_{s'-1}}(v;V_b)}\omega(u)p_{uv;s'}=\big(1-(1\pm\tfrac18\alpha_s)\tfrac{d_{ab;s}}{d_{ab}}\big)\sum_{u\in \NBH_{H_0}(v;V_b)}\omega(u)\,.
\end{equation}
As before, we split this sum into a main term
\[ \sum_{s'=1}^s\ \sum_{u\in \NBH_{H_{s'-1}}(v;V_b)}\omega(u)p_{uv;s'}\]
and an error term
\begin{align*}
 \sum_{s'=1}^sC\alpha_{s'-1}\sum_{u\in \NBH_{H_{s'-1}}(v;V_b)}\omega(u)p_{uv;s'}&\le C\nu^{-3}n^{-1}W\sum_{s'=1}^s\alpha_{s'-1}\\
 &\le \nu^4 \alpha_s W,
\end{align*}
where the second inequality bounds $\sum_{s'=1}^s\alpha_{s'-1}$ by $\int_{-\infty}^{s}\alpha_x\,\mathrm{d}x=\tfrac{\nu^{20L}\alpha_ss^*}{1000CL}$ and then uses~\eqref{eq:pathpack:sstarupper}.
In the main term, consider, for some integer $\ell\ge0$, the interval $\ell\eps n<s'\le(\ell+1)\eps n$. Each vertex $u$ which is in $\NBH_{H_{s'}}(v;V_b)$ for the entire interval contributes by Claim~\ref{cl:pathpack:sump}
\[\Big((1\pm3\gamma)\tfrac{d_{ab;\ell\eps n}-d_{ab;(\ell+1)\eps n}}{d_{ab;\ell\eps n}}\pm3\gamma n^{-0.01}\Big)\omega(u)\]
to the sum. Thus an upper bound for this interval of the main term is
\[\Big((1\pm3\gamma)\tfrac{d_{ab;\ell\eps n}-d_{ab;(\ell+1)\eps n}}{d_{ab;\ell\eps n}}\pm3\gamma n^{-0.01}\Big)\sum_{u\in \NBH_{H_{\ell\eps n}}(v;V_b)}\omega(u)\,.\]
Let $R_\ell$ be the set of vertices removed from the neighbourhood during this interval. Each of them removes at most its weight times $\nu^{-3}\eps$ (by the general upper bound on $p_{uv;s'}$) from the sum, so this upper bound is off from the correct answer by at most $\nu^{-3}\eps\sum_{u\in R_\ell}\omega(u)$. Furthermore, since we are in $\cE$ we know
\[\sum_{u\in \NBH_{H_{\ell\eps n}}(v;V_b)}\omega(u)=\big(1\pm \tfrac12\alpha_{\ell\eps n}\big)\tfrac{d_{ab;\ell\eps n}}{d_{ab}}\sum_{u\in \NBH_{H_0}(v;V_b)}\omega(u)\,,\]
and we conclude
\begin{multline*}
  \sum_{s'=\ell\eps n+1}^{(\ell+1)\eps n}\ \sum_{u\in \NBH_{H_{s'-1}}(v;V_b)}\omega(u)p_{uv;s'}\\
  \begin{aligned}
 &=\Big((1\pm3\gamma)\tfrac{d_{ab;\ell\eps n}-d_{ab;(\ell+1)\eps n}}{d_{ab;\ell\eps n}}\pm3\gamma n^{-0.01}\Big)\big(1\pm \tfrac12\alpha_{\ell\eps n}\big)\tfrac{d_{ab;\ell\eps n}}{d_{ab}}\sum_{u\in \NBH_{H_0}(v;V_b)}\omega(u)\pm \nu^{-3}\eps\sum_{u\in R_\ell}\omega(u)\\
 &=(1\pm \alpha_{\ell\eps n})\tfrac{d_{ab;\ell\eps n}-d_{ab;(\ell+1)\eps n}}{d_{ab}}\sum_{u\in \NBH_{H_0}(v;V_b)}\omega(u)\pm 2\nu^{-3}\eps\sum_{u\in R_\ell}\omega(u) \\
 &=\tfrac{d_{ab;\ell\eps n}-d_{ab;(\ell+1)\eps n}}{d_{ab}}\sum_{u\in \NBH_{H_0}(v;V_b)}\omega(u)\pm 10\nu^{-1}\eps\alpha_{\ell\eps n}W\pm 2\nu^{-3}\eps\sum_{u\in R_\ell}\omega(u)\,.
  \end{aligned}
\end{multline*}
Summing this over the at most $\eps^{-1}$ choices of $\ell$, we have
\begin{multline*}
  \sum_{s'=1}^s\ \sum_{u\in \NBH_{H_{s'-1}}(v;V_b)}\omega(u)p_{uv;s'} \\
  \begin{aligned}
    &=\sum_{\ell=0}^{\tfrac{s}{\eps n}-1} \tfrac{d_{ab;\ell\eps n}-d_{ab;(\ell+1)\eps n}}{d_{ab}}\sum_{u\in \NBH_{H_0}(v;V_b)}\omega(u)\pm 10\nu^{-1}\eps W\sum_{\ell=0}^{\tfrac{s}{\eps n}-1}\alpha_{\ell\eps n}\pm 2\nu^{-3}\eps W\\
    &=\tfrac{d_{ab}-d_{ab;s}}{d_{ab}}\sum_{u\in \NBH_{H_0}(v;V_b)}\omega(u)\pm\nu^4\alpha_s W\,,
  \end{aligned}
\end{multline*}
where the final line uses our definition of $\alpha_x$, the disjointness of the sets $R_\ell$, and summing the error terms as in previous claims.

Putting back the error term, we find that the left side of~\eqref{eq:pathpack:weights:finsum} is in the interval
\[\tfrac{d_{ab}-d_{ab;s}}{d_{ab}}\sum_{u\in \NBH_{H_0}(v;V_b)}\omega(u)\pm 2\nu^4\alpha_s W\]
which in particular shows~\eqref{eq:pathpack:weights:finsum} holds, proving the claim.
\end{claimproof}

We are now in a position to show it is likely that Algorithm~\ref{alg:path-pack} completes, from which~\ref{pathpack:complete} and~\ref{pathpack:pack} follow. What we need to do is show that it is likely that at each stage $s$ the conditions of Lemma~\ref{lem:RPE} are met, since then it is unlikely (by Lemma~\ref{lem:RPE} and the union bound) that at any stage Lemma~\ref{lem:RPE} fails to produce an embedding. As before, it is enough to show that, for any given $1\le s'\le s^*$, with probability at most $3\exp\big(-n^{0.3}\big)$ stage $s'$ is the first stage at which the conditions of Lemma~\ref{lem:RPE} are not met or its low-probability event occurs. Hence, fix $s'$ and assume that for all $s<s'$ the conditions of Lemma~\ref{lem:RPE} are met and its low-probability event does not occur.

Fix vertices $u_1,\dots,u_p\in V(H_0)$ for some $p\le L$ and $a\in\{\boxminus,\boxplus\}$. The probability that these vertices and side of $H_{s'-1}$ witness a failure of either $(\alpha_{s'-1},L)$-block-quasirandomness of $H_{s'-1}$, or of the $(\alpha_{s'-1},L)$-block-diet condition for $(H_{s'-1},U_{s'})$, is by Claim~\ref{cl:pathpack:eat} at most $2\exp\big(-n^{0.4}\big)$. Here and below, for $a\in\{\boxminus,\boxplus\}$ and $i\in[n/2]$ we write $a_i$ for the $i$th vertex of $V_a$, that is, $a_i=\boxminus_i$ if $a=\boxminus$ and $a_i=\boxplus_i$ if $a=\boxplus$. In both cases we apply the claim with $s=1$ and $s'$, in the first place with $U=\{i\in [n/2]:~ a_i\in\NBH_{H_0}(u_1,\dots,u_p)\cap V_a\}$, and in the second place with $U=\{i\in [n/2]:~ a_i\in\NBH_{H_0}(u_1,\dots,u_p)\cap V_a\setminus U_{s'}\}$; in either case, if $a=\boxminus$ we use $W_1=\{u_1,\dots,u_p\}$ and $W_2=\emptyset$, while if $a=\boxplus$ we use $W_2=\{u_1,\dots,u_p\}$ and $W_1=\emptyset$. Taking the union bound over choices of $p$, $a$ and $u_1,\dots,u_p$ we see that the probability such a witness exists is at most $\exp\big(-n^{0.3}\big)$.

Now fix $a,b,c\in\{\boxminus,\boxplus\}$, a vertex $v\in V_{b}$ and let $\omega$ be the weight function on $V(H_0)$ where
for $u\in V_a$ we set $\omega(u)=1$ if there is $x\in A_{s'}$ such that $\phi_{s'}(x)=u$ and the neighbour $y$ of $x$ in $F_{s'}$ has $\chi_{s'}(y)=V_b$
and the next neighbour $z$ of $y$ in $F_{s'}$ has $\chi_{s'}(z)=V_c$;
and otherwise $\omega(u)=0$. The $\gamma$-anchor distribution property~\ref{PA:AnchorDistribution} states that
\[\sum_{u\in \NBH_{H_0}(v;V_a)}\omega(u)=(1\pm\gamma)d_{ab}\sum_{u\in V_a}\omega(u)\pm\gamma n^{0.99}\,.\]
Suppose first that $\sum_{u\in \NBH_{H_0}(v;V_a)}\omega(u)<n^{0.95}$. The display above then gives
\[d_{ab;s'-1}\sum_{u\in V_a}\omega(u)\le d_{ab}\sum_{u\in V_a}\omega(u)\le 2\gamma n^{0.99}\,,\]
while $\sum_{u\in \NBH_{H_{s'-1}}(v;V_a)}\omega(u)\le n^{0.95}$. Both quantities thus lie in $\big[0,2\gamma n^{0.99}\big]$, so, since $\alpha_{s'-1}\ge10\gamma$, this choice of $a,b,c,v$ does not witness a failure of the $\alpha_{s'-1}$-anchor distribution property. We may therefore assume $\sum_{u\in \NBH_{H_0}(v;V_a)}\omega(u)\ge n^{0.95}$.
Now Claim~\ref{cl:pathpack:weights} tells us that with probability at least $1-\exp\big(-n^{0.4}\big)$ we have
\[\sum_{u\in \NBH_{H_{s'-1}}(v;V_a)}\omega(u)=(1\pm\gamma)(1\pm\tfrac12\alpha_{s'-1})d_{ab;s'-1}\sum_{u\in V_a}\omega(u)\pm\gamma n^{0.99}\,\]
(or in some earlier stage the conditions of Lemma~\ref{lem:RPE} are not met or its low-probability event occurs, which we already  assumed not to happen).
Since $\alpha_{s'-1}\ge10\gamma$, we have $(1\pm\gamma)(1\pm\tfrac12\alpha_{s'-1})=1\pm\alpha_{s'-1}$ and $\gamma n^{0.99}\le\tfrac12\alpha_{s'-1}n^{0.99}$. Hence in particular this choice of $a,b,c,v$ does not witness a failure of the $\alpha_{s'-1}$-anchor distribution property for $F_{s'}$ in $H_{s'-1}$. Taking a union bound over the choices of $a,b,c,v$ we see that with probability at most $\exp\big(-n^{0.3}\big)$ the $\alpha_{s'-1}$-anchor distribution property fails for $H_{s'-1}$ and $F_{s'}$.

Putting these together, the probability that the conditions of Lemma~\ref{lem:RPE} are first not met at stage $s'$ is at most $2\exp\big(-n^{0.3}\big)$, and the probability that the conditions of Lemma~\ref{lem:RPE} are met but Algorithm~\ref{alg:path-pack} fails at stage $s'$ is at most $\exp\big(-n^{0.3}\big)$. Taking the union bound over the choices of $s'$, we conclude that as desired the probability of any of these events occurring is at most $\exp\big(-n^{0.25}\big)$.

By choice of $C'$, we immediately see that if at each stage $s'$, including $s^*$, we have the above claimed quasirandomness then~\ref{pathpack:quasi} holds. Furthermore Claim~\ref{cl:pathpack:weights} states that~\ref{pathpack:weights} holds.

Next we prove~\ref{pathpack:megaqr}. To that end, fix $S_1,S_2$ sets of at most $L$ vertices, and disjoint $T_1,T_2\subset [s^*]$ of size at most $L$, and a pairing $\Vmin=\{\boxminus_i:i\in[n/2]\}$ and $\Vplus=\{\boxplus_i:i\in[n/2]\}$. Let $X$ be any subset of $\mathbb{U}_H(S_1,S_2,T_1,T_2)$ of size at least $\nu n$. Observe that the elements of $T_1\cup T_2$ are natural numbers and hence ordered, and the corresponding path-forests are embedded in that order by Algorithm~\ref{alg:path-pack}. We let $X_0:=X$. We now define some subsets $X'_j$ and $X_j$ of $X$ recursively. For $1\le j\le|T_1|+|T_2|$, let $X'_j$ denote the set of indices $i\in X_{j-1}$ such that $w\boxminus_i,w'\boxplus_i\in E(H_s)$ for each $w\in S_1$ and $w'\in S_2$, where $s$ is the stage immediately before embedding the $j$th element of $T_1\cup T_2$; let $X'_{|T_1|+|T_2|+1}$ be similarly defined, with $s=s^*$. For $1\le j\le |T_1|+|T_2|$, we define $X_j\subset X'_j$ as follows. Let the $j$th element of $T_1\cup T_2$ be $s$. If $s$ is in $T_1$, then $X_j$ is the set of $i\in X'_j$ such that $\boxminus_i\notin\im\phi_s$. If $s$ is in $T_2$, then $X_j$ is the set of $i\in X'_j$ such that $\boxplus_i\notin\im\phi_s$.

We show~\ref{pathpack:megaqr} by finding ratios between sizes of each $X_{j-1}$ and $X'_j$, and each $X'_j$ and $X_j$. Note that $X'_{|T_1|+|T_2|+1}=\mathbb{U}'_{H_{s^*}}(S_1,S_2,T_1,T_2)$ is the set whose size we want to find, by definition. For the former, we use Claim~\ref{cl:pathpack:eat} to estimate the effect of the sequence of embeddings of path-forests strictly between the $(j-1)$st and $j$th elements of $T_1\cup T_2$. We also need to estimate the effect of removing edges when the $(j-1)$st path-forest of $T_1\cup T_2$ is embedded (edges may be removed going to a $\boxplus_i$ while the $(j-1)$st path-forest is in $T_1$, for example, which is not accounted for in $X_{j-1}$) but these can be responsible for removing at most $2L$ vertices from $X_{j-1}$. For the latter, we use~\ref{RPE:imageeating}. In each case, it will be enough to use $\alpha_{s^*}$ to bound error terms, which we will do for simplicity.

We claim inductively that the sets $X_j,X'_j$ are of size at least $\nu^{5L}n$; this justifies the following applications of Claim~\ref{cl:pathpack:eat} and~\ref{RPE:imageeating}. First, for each $1\le j\le |T_1|+|T_2|+1$, we estimate $|X'_j|$. Given $j$, let $s$ be the $(j-1)$st index in $T_1\cup T_2$ (if $j=1$, let $s=0$) and $s''+1$ the $j$th index. Let $Q$ denote the set of $i\in X_{j-1}$ such that $w\boxminus_i,w'\boxplus_i$ are edges of $H_s$ for each $w\in S_1$ and $w'\in S_2$. We have $|Q|=|X_{j-1}|\pm 2L$. We apply Claim~\ref{cl:pathpack:eat}, with input $W_1=S_1,W_2=S_2$, $U=Q$, and $s,s''$. The result is that with probability at least $1-\exp\big(-n^{0.4}\big)$, either we have
\begin{align*}
 |X'_j|&=\big(1\pm\tfrac12\alpha_{s^*}\big)\big(|X_{j-1}|\pm 2L\big)\Big(\tfrac{d_{\boxminus\boxminus;s''}}{d_{\boxminus\boxminus;s}}\Big)^{|S_1\cap\Vmin|}\Big(\tfrac{d_{\boxminus\boxplus;s''}}{d_{\boxminus\boxplus;s}}\Big)^{|S_1\cap\Vplus|+|S_2\cap\Vmin|}\Big(\tfrac{d_{\boxplus\boxplus;s''}}{d_{\boxplus\boxplus;s}}\Big)^{|S_2\cap\Vplus|}\\
 &=\big(1\pm\alpha_{s^*}\big)|X_{j-1}|\Big(\tfrac{d_{\boxminus\boxminus;s''}}{d_{\boxminus\boxminus;s}}\Big)^{|S_1\cap\Vmin|}\Big(\tfrac{d_{\boxminus\boxplus;s''}}{d_{\boxminus\boxplus;s}}\Big)^{|S_1\cap\Vplus|+|S_2\cap\Vmin|}\Big(\tfrac{d_{\boxplus\boxplus;s''}}{d_{\boxplus\boxplus;s}}\Big)^{|S_2\cap\Vplus|}
\end{align*}
or there is some stage $s'$ between $s$ and $s''$ inclusive such that either the conditions of Lemma~\ref{lem:RPE} are not met, or its low-probability event occurs.

Similarly, we can estimate $|X_j|$ for each $1\le j\le |T_1|+|T_2|$. Given $j$, let $s$ be the $j$th index in $T_1\cup T_2$. If $s$ is in $T_1$, let $a=\boxminus$ and $Q=\{\boxminus_i:i\in X'_j\}$; if $s$ is in $T_2$, let $a=\boxplus$ and $Q=\{\boxplus_i:i\in X'_j\}$. We apply Lemma~\ref{lem:RPE} part~\ref{RPE:imageeating} with input $a$ and $V'=Q$; recall that by choice of $X$, we have $Q\cap U_s=\emptyset$. The result is that with probability at least $1-\exp\big(-n^{0.3}\big)$ either we have
\[|X_j|=\big(1\pm C\alpha_{s^*}\big)|X'_j|\frac{|V_a\setminus (U_s\cup \im\phi'_s)|}{|V_a\setminus U_s|}\]
or the conditions of Lemma~\ref{lem:RPE} are not met, or its low-probability event occurs.

If all these equations hold, we obtain precisely the desired~\ref{pathpack:megaqr}, by choice of $C'$. To see this, observe that simply substituting each equation into the next, we are left with the main term of~\ref{pathpack:megaqr}, multiplied by an error term of the form $\big(1\pm C\alpha_{s^*}\big)^{2|T_1|+2|T_2|+1}$, multiplied by a collection of terms of the form $\tfrac{d_{ab;s}}{d_{ab;s-1}}$ where $s\in T_1\cup T_2$ and $a,b\in\{\boxminus,\boxplus\}$. Each of these last terms is of size $1\pm 2\nu^{-1}n^{-1}$, and there are at most $2L(|T_1|+|T_2|)$ of them, so that the product of all these terms is $1\pm\eps$. Since the final estimate is of size at least $\nu^{5L}n$ this justifies the inductive claim on set sizes. So the probability we do not obtain~\ref{pathpack:megaqr} is by the union bound at most $\exp\big(-n^{0.2}\big)$ as required.

Finally, we prove~\ref{pathpack:imagecaps}. Given $S_a\subset V_a$ with $|S_a|\le L$ for each $a\in\{\boxminus,\boxplus\}$, and a set $Y\subset[s^*]$ such that $(S_\boxminus\cup S_\boxplus)\cap U_s=\emptyset$ for each $s\in Y$, with $|Y|\ge\nu n$, we apply Corollary~\ref{cor:freedm} to estimate the number of $s\in Y$ such that $(S_\boxminus\cup S_\boxplus)\cap\im\phi'_s=\emptyset$. We let the good event for Corollary~\ref{cor:freedm} be that at each stage $s$ the conditions of Lemma~\ref{lem:RPE} are met and its low-probability event does not occur. Within this good event, for each $s\in Y$, by~\ref{RPE:probvtx} the probability that $\im\phi'_s\cap (S_\boxminus\cup S_\boxplus)=\emptyset$, conditioned on the embeddings prior to embedding $F_s$, is
\[(1\pm C\alpha_{s^*})\prod_{a\in\{\boxminus,\boxplus\}}\Big(\tfrac{|V_a\setminus (U_s\cup\im\phi'_s)|}{|V_a\setminus U_s|}\Big)^{|S_a|}\,.\]
By Corollary~\ref{cor:freedm}, we conclude that with probability at least $1-\exp(-n^{0.2})$ we have
\[\big|\{s\in Y:(S_\boxminus\cup S_\boxplus)\cap\im\phi'_s=\emptyset\}\big|=(1\pm 2C\alpha_{s^*})\sum_{s\in Y}\prod_{a\in\{\boxminus,\boxplus\}}\Big(\tfrac{|V_a\setminus (U_s\cup\im\phi'_s)|}{|V_a\setminus U_s|}\Big)^{|S_a|}\,,\]
as required for~\ref{pathpack:imagecaps}.
\end{proof}


\section{Stage~A (Proof of Lemma~\ref{lem:StageA})}\label{sec:StageA}

In this stage, the graphs $\left(G_s\right)_{s\in\cG\setminus (\cK\cup\cJ)}\cup \left(G_s^\spadesuit\right)_{s\in\cK}\cup \left(G_s^\parallel\right)_{s\in\cJ}$ are packed into $H$, which is very quasirandom. We do this in two steps. In the first step, we use the algorithm \PackingProcess\ of~\cite{DegPack} to pack the graphs $(G_s)_{s\in\cG\setminus(\cK\cup\cJ)}$ into $H$. We let $H_0$ be the graph of leftover edges from this packing, which we will argue (quoting results of~\cite{DegPack}) is still very quasirandom; we will not need any more facts about this part of the packing.

For the second step, we need to give two algorithms, \RandomEmbedding\ (defined in~\cite{DegPack}) and $\PackingProcess'$ (a modification of \PackingProcess\ from~\cite{DegPack}). To conveniently state these, we need some notation and to define some constants.

\subsection{Notation, constants and algorithms}

Suppose $G$ is a graph whose vertices are ordered by natural numbers,
$V(G)=[\ell]$. Suppose that $i\in V(G)$. We write
$\LNBH(i)=\NBH(i)\cap[i-1]$\index{$\LNBH(i)$} and
$\LEFTDEG(i)=|\LNBH(i)|$\index{$\LEFTDEG(i)$} for the
\emph{left-neighbourhood}\index{left-neighbourhood} and the
\emph{left-degree}\index{left-degree} of $i$.

\begin{definition}
	Let~$G$ be a graph with vertex set $[v(G)]$, and~$H$ be a graph with
	$v(H)\ge v(G)$.	Further, assume $\psi_{j}\colon[j]\rightarrow V(H)$ is a \emph{partial embedding}\index{partial embedding} of $G$ into~$H$ for $j\in[v(G)]$, that is, $\psi_j$ is a graph embedding of $G\big[[j]\big]$ into~$H$. Finally, given $t\in[v(G)]$ we say the \emph{candidate set}\index{candidate set} of $t$ (with respect to~$\psi_j$) is 
	\[\CANDSET_{G\AlgMap H}^{j}(t) 
	=\NBH_{H}\Big(\psi_{j}\big(\LNBH_{G}(t)\cap[j]\big)\Big)\,.
	\index{$\CANDSET_{G\AlgMap H}^{j}(t)$}
	\]
	If $\LNBH_G(t)\cap[j]=\emptyset$, we obtain the common neighbourhood of the empty set, which is $V(H)$. When $j=t-1$, we call $\CANDSET_{G\AlgMap H}^{j}(t)$ the \emph{final candidate set} of $t$\index{final candidate set of $t$}.
\end{definition}

We now define some constants, which we will need in the following analysis and which are identically defined in~\cite{ABCT:PackingManyLeaves} and~\cite{DegPack} (so that we can quote results from both verbatim). For convenience of quoting results, in this section we set $\newD:=2D$ and will work almost exclusively with $\newD$.

We shall work with $\dStageA$ given by~\eqref{eq:dStageA}.
Given $\newD$ and $\deltnonspanning,\gamcore>0$, we choose $0<\gamma<\dStageA$ small enough for two inequalities, which we give after defining the following (copied from~\cite{ABCT:PackingManyLeaves} and~\cite{DegPack}, with the exceptions that $\newD$, $\newdelta$ and $\iniquasi'$ defined below are there $D$, $\delta$ and $\iniquasi$ respectively).

\begin{setting}\label{set:graphs}
	Let $\newD,n\in\NN$ and $\gamma>0$ be given. We define
	\begin{equation}\label{eq:defconsts}\begin{split}
	\eta&=\frac{\gamma^{\newD}}{200\newD}\,,\quad  \newdelta=\frac{\gamma^{10\newD}\eta}{10^{6}\newD^{4}}\,,\quad C=40\newD\exp\big(1000\newD\newdelta^{-2}\gamma^{-2\newD-10}\big)\,,\quad C'=10^4C\newdelta^{-1}\,,\\
	\alpha_x&=\frac{\newdelta}{10^8C\newD}\exp\Big(\frac{10^8C^2\newD^3\newdelta^{-1}\gamma^{-4\newD-6}(x-2n)}{n}\Big) \qquad\text{for each $x\in\RR$},\\
	\eps&= \frac{\alpha_0\newdelta^4\gamma^{10\newD}}{1000C\newD}\,,\quad c= \frac{\newD^{-4}\eps^{4}}{100}\,\quad\text{and}\quad \iniquasi'=\frac{\alpha_0}{100}\,.
	\end{split}\end{equation}
\end{setting}

We require $\gamma$ small enough that $\newdelta\le\tfrac12\deltnonspanning$. In addition, we need
\begin{equation}\label{eq:A:smallerrors}
 10^6CC'\newD\newdelta^{-1}\deltnonspanning^{-1}\gamma^{-1}\alpha_{3n/2}\le\gamcore\,,
\end{equation}
where $\gamcore$ is as defined in~\eqref{eq:CONSTANTS}.
Since several of the terms on the left-hand side of this inequality depend on $\gamma$, it may not be obvious that the left-hand side indeed tends to zero as $\gamma\to0$. We have
\begin{align*}
 10^{6}CC'\newD\newdelta^{-1}\deltnonspanning^{-1}\gamma^{-1}\alpha_{3n/2}&=100C\newdelta^{-1}\deltnonspanning^{-1}\gamma^{-1}\exp\Big(-\tfrac12\cdot10^8C^2\newD^3\newdelta^{-1}\gamma^{-4\newD-6}\Big)\\
 &\le 100C^3\newdelta^{-3}\deltnonspanning^{-1}\exp\Big(-\tfrac12\cdot10^8C\newD^3\newdelta^{-1}\Big)\,.
\end{align*}
Now $\deltnonspanning$ is a fixed constant, while both $C$ and $\newdelta^{-1}$, and so also $C\newdelta^{-1}$, tend to infinity as $\gamma\to 0$. Since any exponential in $C\newdelta^{-1}$ is eventually much larger than any polynomial in $C\newdelta^{-1}$, indeed this quantity tends to $0$ as $\gamma\to 0$, and hence there is a choice of $\gamma>0$ such that~\eqref{eq:A:smallerrors} holds.

We briefly state the purpose of a few of these quantities which are immediately important. We will pack $D$-degenerate graphs (which are also $\newD$-degenerate by definition) on up to $n-\newdelta n$ vertices into an $n$-vertex graph; after completing the packing, at least $\gamma n^2$ edges will remain unused. When we borrow the analysis from~\cite{ABCT:PackingManyLeaves} and~\cite{DegPack}, the quantity $\alpha_s$ (or a multiple of it, usually $C\alpha_s$) will upper bound the errors we have made after packing $s$ graphs; we pack in total at most $\tfrac32n$ graphs, which is why~\eqref{eq:A:smallerrors} is what we need to ensure all the error terms in our final analysis are small compared to $\gamcore$.

We next recall \RandomEmbedding{} from~\cite{DegPack}, which we repeat below. This algorithm takes a graph $G$ with $n-\newdelta n$ vertices, and (if successful) embeds it into an $n$-vertex graph $H$. We will want to use it to embed graphs with at most $n-\deltnonspanning n<n-\newdelta n$ vertices, so we add isolated vertices to each graph we want to embed.

\begin{algorithm}[ht]\index{\RandomEmbedding{}}
	\caption{\RandomEmbedding{}}\label{alg:embed}
	\SetKwInOut{Input}{Input}
	\Input{graphs~$G$ and~$H$, with $V(G)=[v(G)]$ and $v(H)=n$}
	$\psi_0:=\emptyset$\;
	$t^*:=(1-\newdelta)n$\;
	\For{$t=1$ \KwTo $t^*$}{
		\lIf{$\CANDSET_{G\AlgMap H}^{t-1}(t)\setminus\im(\psi_{t-1})=\emptyset$}{
			halt with failure}
		choose $v\in\CANDSET_{G\AlgMap
			H}^{t-1}(t)\setminus\im(\psi_{t-1})$ uniformly at random\;
		$\psi_{t}:=\psi_{t-1}\cup\{t\AlgMap v\}$\;
	}
	\Return $\psi_{t^*}$
\end{algorithm}

In~\cite{DegPack} a packing is obtained by an algorithm \PackingProcess\ which runs \RandomEmbedding\ repeatedly, but which does something in addition in order to allow for packing spanning graphs. Nevertheless, the main task of~\cite{DegPack} is to analyse the repeated running of \RandomEmbedding, and this analysis, and similar analysis in~\cite{ABCT:PackingManyLeaves}, is also valid for the following $\PackingProcess'$ which simply runs \RandomEmbedding\ repeatedly.

\begin{algorithm}[ht]\index{$\PackingProcess'${}}
	\caption{$\PackingProcess'${}}\label{alg:pack}
	\SetKwInOut{Input}{Input}
	\Input{
		$\bullet$ graphs $G_1,\dots,G_{s^*}$, with $G_s$ on vertex set $[v(G_s)]$\\
		$\bullet$ a graph $H_0$ on vertex set of order $n$}
    \For{$s=1$ \KwTo $s^*$}{
		run \RandomEmbedding($G_s$,$H_{s-1}$) to get
		an embedding $\phi_s$ of~$G_s$ into~$H_{s-1}$\;   
		let $H_{s}$ be the graph obtained from $H_{s-1}$ by removing the
		edges of $\phi_{s}(G_s)$\;
	}    
\end{algorithm}

\subsection{Packing and maintaining quasirandomness}\label{sec:A:maintquasi}
 
We first prove that with high probability the process outlined above does return a packing of $\left(G_s\right)_{s\in\cG\setminus (\cK\cup\cJ)}\cup \left(G_s^\spadesuit\right)_{s\in\cK}\cup \left(G_s^\parallel\right)_{s\in\cJ}$, and furthermore after embedding each graph, what is left of $H$ remains quasirandom.

For convenience, we will let $H_0$ be the spanning subgraph of $H$ which we obtain by using~\cite[Theorem~\ref{PAPERpackingdegenerate.thm:maintechGENERAL}]{DegPack} to pack $\left(G_s\right)_{s\in\cG\setminus (\cK\cup\cJ)}$, and we will renumber $\cK\cup\cJ$ to be the integers $[s^*]$. Note that $s^*\le\tfrac32n$, since $|\cK|\le\big|\bigcup_{s\in\cK} \SpecLeaves_s \big|\le(1+2\deltnonspanning)n$ and $|\cJ|\le\sigmKJ n$.

\begin{proof}[Proof of~\ref{A:packing} and $(\alpha_s,2\newD+3)$-quasirandomness of $H_s$]

 Let $\iniquasi'$ be as defined in Setting~\ref{set:graphs}. We choose $\eta^*>0$ such that the graphs of $\left(G_s\right)_{s\in\cG\setminus (\cK\cup\cJ)}$ have at most $e(H)-100\newD(\eta^*)^{1/(8\newD+12)}n^2$ edges, and in addition $\eta^*\le(\iniquasi')^4$. Let $0<\iniquasi\le\iniquasi'$ and $c'>0$ be small enough for \cite[Theorem~\ref{PAPERpackingdegenerate.thm:maintechGENERAL}]{DegPack} with input $\gamma=\eta^*$ and $2\newD$. The constants Lemma~\ref{lem:StageA} returns are then $\iniquasi$ and $\min(c,c')$. For the next paragraph, in which we quote results of~\cite{DegPack}, we need to borrow one further constant. We let $\alpha'$ be the constant defined as $\alpha_x$ with $x=2n$ in~\cite[Setting~\ref{PAPERpackingdegenerate.set:graphs}]{DegPack}, where again we let the input be $\gamma=\eta^*$ and $2\newD$. Note that though this parameter contains a number $n$ in the formula, in the case $x=2n$ the parameter does not depend on $n$. Furthermore, from the formula we have $\alpha'\le\eta^*$.
 
 Given $H$, which is $(\iniquasi,2\newD+3)$-quasirandom (recall $\newD=2D$), we pack $\left(G_s\right)_{s\in\cG\setminus (\cK\cup\cJ)}$ into $H$. To do this, we first need to ensure that $\cG\setminus(\cK\cup\cJ)$ has size at most $2n$. We do this as in~\cite{DegPack}, that is, if we see that $G_s$ and $G_{s'}$ both have $\lfloor n/2 \rfloor$ isolated vertices for some $s\neq s'\in\cG\setminus(\cK\cup\cJ)$, we remove all isolated vertices, take the disjoint union of the two graphs (which has at most $n$ vertices by assumption and continues to respect the maximum degree and degeneracy requirements trivially) and add isolated vertices if necessary to return to $n$ vertices. We repeat this until no two such graphs remain; as in~\cite[Proof of Theorem~\ref{PAPERpackingdegenerate.thm:MAINunbounded}]{DegPack} all but at most one remaining graph have at least $n/4$ edges and hence there are at most $2n$ remaining graphs. Since a packing of this new family trivially gives a packing of the original family, abusing notation we refer to the new family as $\left(G_s\right)_{s\in\cG\setminus(\cK\cup\cJ)}$. We then reorder the vertices of each such $G_s$ such that the final $(2D+1)^{-3}n$ vertices form an independent set all of whose vertices have degree $d$, for some $0\le d\le 2D=\newD$. By~\cite[Lemma~\ref{PAPERpackingdegenerate.lem:degindept}]{DegPack} there exists such an independent set in each $G_s$; we modify the given $D$-degenerate vertex ordering by moving them to the end of the order. For each $G_s$, the new vertex ordering is $\newD$-degenerate. Running \PackingProcess,~\cite[Theorem~\ref{PAPERpackingdegenerate.thm:maintechGENERAL}]{DegPack} with inputs as above states that this packing a.a.s.\ succeeds. Let $H'''$ be the graph of edges of $H$ not used in this packing. Recall that $H'''$ is composed of what remains of the `bulk' and `reservoir'. Now~\cite[Lemma~\ref{PAPERpackingdegenerate.lem:aggregate}]{DegPack} states that the remainder $H'$ of the bulk is $\big(\alpha',2\newD+3\big)$-quasirandom, while the reservoir $H''$ has by~\cite[Lemma~\ref{PAPERpackingdegenerate.lem:firstquasirandomness}]{DegPack} maximum degree $2\eta^* n$ (this follows from our input $\gamma=\eta^*$). The size of $\NBH_{H'''}(v_1,\dots,v_\ell)$ is thus within $\ell\eta^* n$ of $\big|\NBH_{H'}(v_1,\dots,v_\ell)\big|$. Let now $p'$ be such that $e(H')=p'\binom{n}{2}$, and $p'''$ be such that $e(H''')=p'''\binom{n}{2}$. By choice of $\eta^*$, we have $p'''\ge 100\newD\sqrt{\eta^*}$ and $p'''-p'\le 2\eta^*$. Since $H'$ is $(\alpha',2\newD+3)$-quasirandom, and it has at least $50\newD(\eta^*)^{1/(8\newD+12)}n^2$ edges, we have for any $1\le\ell\le 2\newD+3$ and vertices $v_1,\dots,v_\ell$
 \begin{align*}
  \big|\NBH_{H'''}(v_1,\dots,v_\ell)\big|=(1\pm\alpha')(p')^\ell n\pm \ell\eta^* n&=(p''')^\ell n\pm\eta^* n\pm (2\newD+3)\eta^*n\pm (2\newD+3)\frac{4\eta^*}{p'}n\\
  &=(p''')^\ell n\pm (2\newD+4)\eta^* n\pm \sqrt{\eta^*}n\,,
 \end{align*}
 and it follows from our choice of $\iniquasi'$ and since $(p''')^\ell\ge 50\newD(\eta^*)^{1/4}$ that $H'''$ is $(\iniquasi',2\newD+3)$-quasirandom. Having established this fact, we let $H_0:=H'''$ (we avoided referring to this graph as $H_0$ above because the notation $H_0$ is defined and a different graph in~\cite{DegPack}).

 We now run $\PackingProcess'$ to pack the graphs $\big(G_s^\spadesuit\big)_{s\in\cK}\cup \big(G_s^\parallel\big)_{s\in\cJ}$ into $H_0$. For convenience, we renumber these graphs such that $\cK\cup\cJ=[s^*]$. For much of the analysis, we will not want to distinguish between the graphs $G_s^\spadesuit$ and $G_s^\parallel$, so we define graphs $G'_s$ as follows. For $s\in\cK$, we set $G'_s=G_s^\spadesuit$ (and so $v(G'_s)\le (1-\deltnonspanning)n$), and for $s\in\cJ$ we set $G'_s=G_s^\parallel$ and get $v(G'_s)=(1-\deltnonspanning-10\sigmJjedna)n$ for each $s\in \cJ_1$, $v(G'_s)=(1-\deltnonspanning-10\sigmKJ)n$ for each $s\in \cJ_0$, and $v(G'_s)=(1-\deltnonspanning-6\sigmKJ)n$ for each $s\in\cJ_2$. Each $G'_s$ is a subgraph of $G_s$ and hence $D_0$-degenerate; we order its vertices according to a $D_0$-degenerate ordering, so that $\LEFTDEG(y)\le D_0\le D$ for each $y\in V(G'_s)$. In particular the hypotheses of Lemma~\ref{lem:propRandEmb} are met, and in the proofs below we may work with $D$ rather than $\newD$ where the left-degree is what matters.
 
 For each $s$ let $H_s$ be the graph obtained from $H_{s-1}$ by removing the edges used by $\PackingProcess'$ in embedding $G'_s$. Let $\phi_s^\mathbf{A}$ be the embedding of $G'_s$. Observe that~\cite[Lemma~\ref{PAPERpackingdegenerate.lem:dietsimplified}]{DegPack} states that, provided $H_{s-1}$ is $(\alpha_{s-1},2\newD+3)$-quasirandom, \RandomEmbedding\ is unlikely to fail in packing $G'_s$. We claim that the proof of~\cite[Lemma~\ref{PAPERpackingdegenerate.lem:aggregate}]{DegPack} works (trivially, by simply ignoring each mention of the `reservoir') to show that indeed it is unlikely that $H_{s-1}$ is not $(\alpha_{s-1},2\newD+3)$-quasirandom, so that taking a union bound (for each choice of $s$, over the two unlikely events just mentioned) we conclude that $\PackingProcess'$ a.a.s.\ succeeds. However, the reader can also verify that $H_{s-1}$ is unlikely to fail $(\alpha_{s-1},2\newD+3)$-quasirandomness from Lemma~\ref{lem:T} below (taking $T$ to be the entire common neighbourhood, and using the union bound).
\end{proof}

What remains is to prove that the remaining properties of Lemma~\ref{lem:StageA} a.a.s.\ hold, which requires a more careful analysis of $\PackingProcess'$.

\subsection{Behaviour of the random processes}
 We now need to quote some lemmas from~\cite{ABCT:PackingManyLeaves} and~\cite{DegPack} which establish some useful properties of \RandomEmbedding\ and \PackingProcess. We also need one additional concept, the cover condition, which we will see holds for a typical run of \RandomEmbedding.
 
\begin{definition}[cover condition]
Suppose that $G$ and $H$ are two graphs such that $H$ has order $n$, the vertex set of $G$ is $[v(G)]$, and $H$ has density $p$.  Suppose that numbers $\beta,\eps>0$ and $i\in [n-\eps n]$ are given. For each $d\in\NN$ we define
\[X_{i,d}:=\{x\in V(G)\colon i\le x<i+\eps n,|\LNBH(x)|=d\}\,.\]
We say that a partial embedding $\psi$ of $G$ into $H$, which embeds
  $\LNBH(x)$ for each $i\le x<i+\eps n$, satisfies the
  \emph{$(\eps,\beta,i)$-cover condition} if for each $v\in V(H)$ such that $v\notin\im(\psi_{i+\eps n-1})$, and for each $d\in\NN$
 we have
  \[
  \big|\big\{ x\in X_{i,d}:v\in\NBH_{H}\big(\psi(\LNBH(x))\big)\big\}\big|=(1\pm\beta)p^{d}|X_{i,d}|\pm\eps^{2}n\,.
  \]
Note that a corresponding condition for $d=0$ is trivial, even with zero error parameters.
\end{definition}
 
 The following lemma puts together various results from~\cite{ABCT:PackingManyLeaves} and~\cite{DegPack} which describe the typical behaviour of \RandomEmbedding.

\begin{lemma}
  \label{lem:propRandEmb}Assume Setting~\ref{set:graphs}, and let $\alpha\in[\alpha_0,\alpha_{2n}]$ be arbitrary. For each $t\in\RR$ we define
  \[\beta_{t} =2\alpha\exp\left(\tfrac{1000\newD\newdelta^{-2}\gamma^{-2\newD-10}t}{n}\right)\;.\]
  Let~$G$ be a graph on vertex set~$[v(G)]$ with at most $(1-\newdelta)n$ vertices and maximum degree at most
  $cn/\log n$ such that $\LEFTDEG(x)\le \newD$ for each $x\in V(G)$, and let~$H$ be an $(\alpha,2\newD+3)$-quasirandom $n$-vertex
  graph with $p\binom{n}{2}$ edges where $p\ge \gamma$. Fix $k\le 2\newD+3$, a vertex $x\in V(G)$, distinct vertices $v,v',v''\in V(H)$ and distinct neighbours $u_1,\dots,u_k$ of $v$.
  
  Then with probability at
  least $1-2n^{-9}$ all of the following good events hold.
  \begin{enumerate}[label=\abc]
    \item\label{pfA:nofail} When \RandomEmbedding{}
      is run it does not fail and generates a sequence
      $(\psi_{i})_{i\in[v(G)]}$ of partial embeddings of~$G$
      into~$H$. 
    \item\label{pfA:diet} For each $t\in[v(G)]$ the pair $(H,\im\psi_{t})$
      satisfies the $(\beta_t,2\newD+3)$-diet condition.
    \item\label{pfA:cover} The embedding $\psi_{v(G)}$ of~$G$ into~$H$
      satisfies the $(\eps,20\newD\beta_{t-\eps n+2},t-\eps n+2)$-cover
      condition for each $t\in[\eps n-1,n-\newdelta n-1]$.
  \end{enumerate}
  Note that we have $20\newD\beta_{t}\le C\alpha$ for each $0\le t\le n$.
  We have
  \begin{align}
  \label{eq:A:probvtx} \Prob\big[x\AlgMap v\big]&=(1\pm10^4C\alpha \newD\newdelta^{-1})n^{-1}\,,\\
  \label{eq:A:probvtxnoim} \Prob\big[x\AlgMap v\text{ and }v'\notin\im\psi_{v(G)}\big]&=(1\pm10^4C\alpha \newD^2\newdelta^{-1})\big(n-v(G)\big)n^{-2}\qquad\text{and}\\
  \label{eq:A:probvtxnoim2} \Prob\big[x\AlgMap v\text{ and }v',v''\notin\im\psi_{v(G)}\big]&=(1\pm10^4C\alpha \newD^2\newdelta^{-1})\big(n-v(G)\big)^2n^{-3}\,.
  \end{align}
  
  Finally, the probability that there is at least one $u_iv$ to which some edge of $G$ is embedded is
 \begin{equation}\label{eq:A:probstar}
  \hspace{3.6cm}\big(1\pm1000C\alpha\newdelta^{-1}\big)^{4\newD+2}p^{-1}n^{-2}\cdot 2ke(G)\,.\hspace{3.6cm}\qed
 \end{equation}
\end{lemma}

The part of this lemma showing that~\ref{pfA:nofail}--\ref{pfA:cover} are likely is part of~\cite[Lemma~\ref{PAPERmanyleaves.lem:22}]{ABCT:PackingManyLeaves}. The statement~\eqref{eq:A:probvtx} is~\cite[Lemma~\ref{PAPERmanyleaves.lem:vertex}]{ABCT:PackingManyLeaves}. The equations~\eqref{eq:A:probvtxnoim} and~\eqref{eq:A:probvtxnoim2} follow from the method of~\cite[Lemma~\ref{PAPERpackingdegenerate.lem:probedge}]{DegPack} but for completeness we give a proof in the following section. Finally~\eqref{eq:A:probstar} is given by~\cite[Lemma~\ref{PAPERpackingdegenerate.lem:probstar}]{DegPack}. For this last comment, we need to be slightly careful. The statement of~\cite[Lemma~\ref{PAPERpackingdegenerate.lem:probstar}]{DegPack} includes an assumption $e(G)\ge\tfrac{n}{4}$ which we omit here. We now briefly explain why this assumption is not needed. If $e(G)=0$, then~\eqref{eq:A:probstar} is trivial, so we assume $e(G)\ge1$. In the proof of~\cite[Lemma~\ref{PAPERpackingdegenerate.lem:probstar}]{DegPack}, the place where $e(G)\ge \tfrac{n}{4}$ is used is at the end. A probability bound
\[\big(1\pm500C\alpha\newdelta^{-1}\big)^{4\newD+2}p^{-1}n^{-2}\cdot 2ke(G)\pm\binom{2\newD+3}{2}\cdot 2\newD n\Delta(G)\cdot 8\gamma^{-3\newD}\newdelta^{-3}n^{-3}\pm 2n^{-9}\]
is obtained on the probability of embedding an edge of $G$ to at least one $u_iv$. For the term $\binom{2\newD+3}{2}\cdot 2\newD n\Delta(G)\cdot 8\gamma^{-3\newD}\newdelta^{-3}n^{-3}$, we can be a bit more careful in the analysis. The quantity $2\newD n\Delta(G)$ is an upper bound on the number of triples of distinct vertices $x,x',y\in V(G)$ such that $xy,x'y\in E(G)$. However, observe that the number of such triples is
\[\sum_{y\in V(G)}\deg(y)(\deg(y)-1)\le\sum_{y\in V(G)}\deg(y)^2\le\sum_{y\in V(G)}\Delta(G)\deg(y)=2e(G)\Delta(G)\,,\]
and using this tighter bound we obtain the improved probability estimate
\[\big(1\pm500C\alpha\newdelta^{-1}\big)^{4\newD+2}p^{-1}n^{-2}\cdot 2ke(G)\pm\binom{2\newD+3}{2}\cdot 2e(G)\Delta(G)\cdot 8\gamma^{-3\newD}\newdelta^{-3}n^{-3}\pm 2n^{-9}\,,\]
from which~\eqref{eq:A:probstar} follows since (as $\Delta(G)\le cn/\log n$ and $e(G)\ge1$) the two additive error terms are asymptotically smaller than the main term.

We should note that the quantities $\beta_t$ defined in this lemma are needed only for the more precise bounds in~\ref{pfA:diet} and~\ref{pfA:cover}. We will only need these more precise bounds in two places, namely when we state Lemma~\ref{lem:S} below and apply the following lemma from~\cite{ABCT:PackingManyLeaves} in its proof, and when we apply Lemma~\ref{lem:S} to prove~\ref{A:megaquasirandomness1}. Otherwise, when we use~\ref{pfA:diet} and~\ref{pfA:cover} we will need only the weaker bound $C\alpha$.

\begin{lemma}[{\cite[Lemma \ref{PAPERmanyleaves.lem:24}]{ABCT:PackingManyLeaves}}]\label{lem:24}
  We assume Setting~\ref{set:graphs}. Given $\alpha_0\le\alpha\le\alpha_{2n}$, for each $t\in\RR$ we define
  \[\beta_{t} =2\alpha\exp\left(\tfrac{1000\newD\newdelta^{-2}\gamma^{-2\newD-10}t}{n}\right)\;.\]
  Let~$G$ be a graph on vertex set~$[v(G)]$ with at most $(1-\newdelta)n$ vertices and maximum degree at most
  $cn/\log n$ such that $\LEFTDEG(x)\le \newD$ for each $x\in V(G)$, and let~$H$ be an $(\alpha,2\newD+3)$-quasirandom $n$-vertex
  graph with at least~$\gamma\binom{n}{2}$ edges. When we run \RandomEmbedding\ to embed $G$ into~$H$, suppose that it 
  produces a partial embedding~$\psi_j$ such that $(H,\im\psi_j)$ has the
  $(\beta_j,2\newD+3)$-diet condition, and let $T\subset V(H)\setminus\im\psi_j$
  with $|T|\ge\frac12\gamma^{2\newD+3}\newdelta n$. Conditioning on $\psi_j$, with probability at least $1-2n^{-2\newD-19}$, one of the following
  occurs.
  \begin{enumerate}[label=\abc]
  \item\label{e:part1} $\psi_{v(G)}$ does not have the $(\eps,20\newD\beta_j,j)$-cover condition, or
  \item\label{e:part2} $\big|\{x\colon j\le x<j+\eps n, \psi_{j+\eps n}(x)\in T\}\big|=(1\pm40\newD\beta_j)\frac{|T|\eps n}{n-j}$.
  \end{enumerate}
\end{lemma}

\subsection{Further analysis of $\PackingProcess'$}
In this subsection we give some more analysis of $\PackingProcess'$ which is similar to results of~\cite{ABCT:PackingManyLeaves} and~\cite{DegPack}; the proofs are copied and modified appropriately.

We first need a strengthening of~\cite[Lemma~\ref{PAPERpackingdegenerate.lem:probvtx}]{DegPack}, allowing us to estimate accurately the probability that not just one or two but a collection of up to $\newD$ vertices of $H$ are avoided in any given interval of vertices of $G$ that \RandomEmbedding\ embeds. To begin with, we show the statement for intervals of length $\eps n$, a modification of~\cite[Lemma~\ref{PAPERpackingdegenerate.lem:probint}]{DegPack}.

\begin{lemma}\label{lem:A:probint}
Assume Setting~\ref{set:graphs}. The following holds for any $\alpha_0\le\alpha\le\alpha_{2n}$ and all sufficiently large $n$. Suppose
that $G$ is a graph on $[v(G)]$ with $v(G)\le(1-\deltnonspanning)n$ such that $\LEFTDEG(x)\le \newD$ for each $x\in V(G)$, and~$H$ is an $(\alpha,2\newD+3)$-quasirandom graph
with~$n$ vertices and $p\binom{n}{2}$ edges, with $p\ge\gamma$. Suppose that $1\le k\le \newD$ and $u_1,\dots,u_k$ are any distinct vertices of $H$. When \RandomEmbedding{} is run to embed~$G$ into~$H$, for any $1\le t'\le v(G)+1-\eps n$ the following holds.
Suppose the history $\hist_{t'-1}$ up to and including embedding $t'-1$ is such that $u_1,\dots,u_k\notin\im\psi_{t'-1}$, the $(C\alpha,2\newD+3)$-diet condition holds for $(H,\im\psi_{t'-1})$, and 
\[\Prob\big[(H,\psi_{t'+\eps n-1})\text{ does not satisfy the $(\eps,C\alpha,t')$-cover condition}\big|\hist_{t'-1}\big]\le n^{-3}\,.\]
Then we have
\begin{equation*}
\Prob\Big[\big|\{u_1,\dots,u_k\}\cap\im\psi_{t'+\eps n-1}\big|\ge 1 \Big|\hist_{t'-1}\Big]=(1\pm 10C\alpha)\tfrac{k\eps n}{n-t'}\,.
\end{equation*}
\end{lemma}
The proof of this lemma is a straightforward modification of the proof of~\cite[Lemma~\ref{PAPERpackingdegenerate.lem:probint}]{DegPack}, which we include for completeness.
\begin{proof}
We modify \RandomEmbedding\ as follows to obtain \ModRandomEmbedding. At the line where \RandomEmbedding\ chooses $v\in\CANDSET_{G\AlgMap H}^{t-1}(t)\setminus\im(\psi_{t-1})$ uniformly at random, we choose instead
\[w\in\CANDSET_{G\AlgMap H}^{t-1}(t)\setminus\big(\im(\psi_{t-1})\setminus\{u_1,\dots,u_k\}\big)\]
uniformly at random and \emph{report} $w$. If $w\in\im(\psi_{t-1})$, then we choose $v\in\CANDSET_{G\AlgMap H}^{t-1}(t)\setminus\im(\psi_{t-1})$ uniformly at random; if not, we set $v=w$. We then embed $t$ to $v$ and continue as in \RandomEmbedding.

Observe that the distribution over embeddings created by \ModRandomEmbedding\ is identical to that created by \RandomEmbedding. The difference is that \ModRandomEmbedding\ in addition creates a string of reported vertices. For each $t$, let $r(t)$ be the vertex reported by  \ModRandomEmbedding\ at time $t$.  Note that $r(t)$ is the vertex to which $t$ is embedded, except when $r(t)$ is one of the $u_i$ and that $u_i$ is in $\im(\psi_{t-1})$. We shall use the following two auxiliary claims.

Let $t'$ and $\hist_{t'-1}$ be as in the lemma statement. Define $E$ as the random variable counting the times when some $u_i$ is reported by \ModRandomEmbedding\ in the interval $t'\le x<t'+\eps n$,
\[E=\left|\:\big\{x\in [t',t'+\eps n) \::\: r(x-1)\in\{u_1,\dots,u_k\}\big\}\:\right|\;.\]
The probability that \RandomEmbedding\ uses at least one vertex of $\{u_1,\dots,u_k\}$ in the interval $t'\le x<t'+\eps n$, conditioning on $\hist_{t'-1}$, is equal to the probability that  \ModRandomEmbedding\ reports some vertex of $\{u_1,\dots,u_k\}$ at least once in that interval, which probability is at most $\Exp\left(E\:|\:\hist_{t'-1}\right)$ and, by the tail-sum formula for $\Exp\left(E\:|\:\hist_{t'-1}\right)$, at least
\begin{align*}
 \Exp\left(E\:|\:\hist_{t'-1}\right)-\sum_{j=2}^{\eps n}\Prob&\big[\text{vertices of $\{u_1,\dots,u_k\}$ are reported}\\
 &\text{at least $j$ times in the interval $[t',t'+\eps n)$}\big|\hist_{t'-1}\big]\,.
\end{align*}
Our first claim estimates $\Exp\left(E\:|\:\hist_{t'-1}\right)$.
 \begin{claim}\label{cl:probint:E}
 	We have that
 	\[\Exp\left(E\:|\:\hist_{t'-1}\right)=(1\pm 4C\alpha)\frac{\eps k n}{n-t'}\pm 8\newD^2(\newD+1)\eps^2\gamma^{-2\newD}\deltnonspanning^{-2}\;.\]
 \end{claim}

Our second claim is that the sum in the expression above is small.
 
 \begin{claim}\label{cl:probint:sum}
 	We have that
 	\begin{align*}
     \sum_{j=2}^{\eps n}\Prob&\big[\text{vertices of $\{u_1,\dots,u_k\}$ are reported}\\
     &\text{at least $j$ times in the interval $[t',t'+\eps n)$}\big|\hist_{t'-1}\big]\le 8k^2\eps^2\gamma^{-2\newD}\deltnonspanning^{-2}\,.
\end{align*}
 \end{claim}
By choice of $\eps$, we have $16\newD^2(\newD+1)\eps^2\gamma^{-2\newD}\deltnonspanning^{-2}<C\alpha\eps\deltnonspanning$. Thus the two claims  give Lemma~\ref{lem:A:probint}. We now prove the auxiliary Claims~\ref{cl:probint:E} and~\ref{cl:probint:sum}.

\begin{claimproof}[Proof of Claim~\ref{cl:probint:E}]
	Note that since the $(C\alpha,2\newD+3)$-diet condition holds for $(H,\im\psi_{t'-1})$, for each $t'\le x<t'+\eps n$, setting $S=\psi_{x-1}(\LNBH(x))$,  we have
	\begin{equation}\label{eq:probint:setsize}\begin{split}
	\big|\CANDSET_{G\AlgMap H}^{x-1}(x)\setminus\im\psi_{x-1}\big|\pm k&=\big|\NBH_H(S)\setminus \im\psi_{t'-1}\big|\pm \eps n\pm k\\
	&=(1\pm C\alpha)p^{|\LNBH(x)|}(n-t'+1)\pm\eps n\pm k\\
	&=(1\pm 2C\alpha)p^{|\LNBH(x)|}(n-t')\,.
	\end{split}\end{equation}
	
    We first give a simple upper bound on $\Prob\big[\big|\{u_1,\dots,u_k\}\cap\im\psi_{t'+\eps n-1}\big|\ge 1\,\big|\,\hist_{t'-1}\big]$. When we embed any one vertex $x$ with $t'\le x\le t'+\eps n-1$, the probability of embedding $x$ to $\{u_1,\dots,u_k\}$ is at most $k\big|\CANDSET_{G\AlgMap H}^{x-1}(x)\setminus\im\psi_{x-1}\big|^{-1}$. Using~\eqref{eq:probint:setsize} and summing over the $\eps n$ choices of $x$, we have
    \[	\Prob\big[\big|\{u_1,\dots,u_k\}\cap\im\psi_{t'+\eps n-1}\big|\ge 1\,\big|\,\hist_{t'-1}\big]\le\frac{k\eps n}{\tfrac12p^{\newD}\deltnonspanning n}\le 2k\eps\gamma^{-\newD}\deltnonspanning^{-1}\,.\]
	
	By linearity of expectation, we have
	\begin{align}
	\begin{split}\label{eq:obed}
	\Exp\Big[E\:|\:\hist_{t'-1}\Big]&=\sum_{x=t'}^{t'+\eps n-1}\Prob\big[\text{a vertex of $\{u_1,\dots,u_k\}$ is reported at time }x\big|\hist_{t'-1}\big]\\
	&=\sum_{x=t'}^{t'+\eps n-1}\sum_{j=1}^k\Prob\big[u_j\text{ is reported at time }x\big|\hist_{t'-1}\big]\\
	&=\sum_{x=t'}^{t'+\eps n-1}\sum_{j=1}^k\Exp\left(\frac{\mathbbm{1}\{u_j\in\CANDSET_{G\AlgMap H}^{x-1}(x)\}}{|\CANDSET_{G\AlgMap H}^{x-1}(x)\setminus(\im\psi_{x-1}\setminus\{u_1,\dots,u_k\})|}\:\Big|\hist_{t'-1}\right)\\
	&=\sum_{x=t'}^{t'+\eps n-1}\sum_{j=1}^k\Exp\left(\frac{\mathbbm{1}\{u_j\in\CANDSET_{G\AlgMap H}^{x-1}(x)\}}{|\CANDSET_{G\AlgMap H}^{x-1}(x)\setminus \im\psi_{x-1}|\pm k}\:\Big|\hist_{t'-1}\right)\;.
	\end{split}
	\end{align}
	Using~\eqref{eq:probint:setsize}, we get
	\[\Exp\left(E\:|\:\hist_{t'-1}\right)=\sum_{x=t'}^{t'+\eps n-1}\sum_{j=1}^k\frac{\Prob\big[u_j\in\CANDSET_{G\AlgMap H}^{x-1}(x)\big|\hist_{t'-1}\big]}{(1\pm 2C\alpha)p^{|\LNBH(x)|}(n-t')}\,.\]
	Splitting this sum up according to $|\LNBH(x)|$, and again using linearity of expectation, we have
	\[\Exp\left(E\:|\:\hist_{t'-1}\right)=\sum_{j=1}^k\sum_{d=0}^{\newD}\frac{\Exp\left(|\{x\in X_{t',d}:u_j\in\CANDSET_{G\AlgMap H}^{x-1}(x)\}|\big|\hist_{t'-1}\right)}{(1\pm 2C\alpha)p^{d}(n-t')}\,.\]
	
	Now fix $0\le d\le \newD$ and each $j\in[k]$. If the $(\eps,C\alpha,t')$-cover condition holds, and if $u_j\notin\im\psi_{t'+\eps n-1}$, we have $\big|\{x\in X_{t',d}:u_j\in \CANDSET_{G\AlgMap H}^{x-1}(x)\}\big|=\big(1\pm C\alpha\big)p^d|X_{t',d}|\pm\eps^2n$. If the $(\eps,C\alpha,t')$-cover condition fails, or if $u_j\in\im\psi_{t'+\eps n-1}$ (which occur with total probability at most $n^{-3}+2k\eps\gamma^{-\newD}\deltnonspanning^{-1}$), we have $0\le \big|\{x\in X_{t',d}:u_j\in \CANDSET_{G\AlgMap H}^{x-1}(x)\}\big|\le\eps n$. In this case we therefore have
	\[\big|\{x\in X_{t',d}:u_j\in \CANDSET_{G\AlgMap H}^{x-1}(x)\}\big|=\big(1\pm C\alpha\big)p^d|X_{t',d}|\pm\eps^2 n\pm \eps n\,.\]
	Putting these together, we get
	\begin{align*}
		\Exp\left(|\{x\in X_{t',d}:u_j\in\CANDSET_{G\AlgMap H}^{x-1}(x)\}|\big|\hist_{t'-1}\right)&=\big((1\pm C\alpha)p^d|X_{t',d}|\pm\eps^2 n\big)\pm \big(n^{-3}+2k\eps\gamma^{-\newD}\deltnonspanning^{-1}\big)\cdot\eps n\\
		&=(1\pm C\alpha)p^d|X_{t',d}|\pm 4k\gamma^{-\newD}\deltnonspanning^{-1}\eps^2n\,.
	\end{align*}
	Substituting this in, we have
	\begin{align*}
	 \Exp\left(E\:|\:\hist_{t'-1}\right)&=\sum_{j=1}^k\sum_{d=0}^{\newD}\frac{(1\pm C\alpha)p^d|X_{t',d}|\pm 4k\gamma^{-\newD}\deltnonspanning^{-1}\eps^2n}{(1\pm 2C\alpha)p^{d}(n-t')}\\
	 &=(1\pm 4C\alpha)\tfrac{\eps k n}{n-t'}\pm 8\newD^2(\newD+1)\eps^2\gamma^{-2\newD}\deltnonspanning^{-2}\,,
	\end{align*}
	where the last equality uses $k\le \newD$, $p\ge\gamma$ and $n-t'\ge\deltnonspanning n$.
\end{claimproof} 
\begin{claimproof}[Proof of Claim~\ref{cl:probint:sum}]
 Since the $(C\alpha,2\newD+3)$-diet condition holds for $(H,\im\psi_{t'-1})$, since $p\ge\gamma$, and since $n-t'\ge\deltnonspanning n$, for each $x\in[t',t'+\eps n)$, when we embed $x$ we report a uniform random vertex from a set of size at least $\tfrac12\gamma^{\newD}\deltnonspanning n$. The probability of reporting one of $u_1,\dots,u_k$ when we embed $x$ is thus at most $2k\gamma^{-\newD}\deltnonspanning^{-1}n^{-1}$, conditioning on $\hist_{t'-1}$ and any embedding of the vertices $[t',x)$. Since the conditional probabilities multiply, the probability that at each of a given $j$-set of vertices in $[t',t'+\eps n)$ we report a vertex of $\{u_1,\dots,u_k\}$ is at most $2^jk^j\gamma^{-j\newD}\deltnonspanning^{-j}n^{-j}$. Taking the union bound over choices of $j$-sets, we have
 \begin{align*}
 	&\sum_{j=2}^{\eps n}\Prob\big[\text{vertices of $\{u_1,\dots,u_k\}$ are reported at least $j$ times in the interval $[t',t'+\eps n)$}\big|\hist_{t'-1}\big]\\
 	\le&\sum_{j=2}^{\eps n}\binom{\eps n}{j}2^jk^j\gamma^{-j\newD}\deltnonspanning^{-j}n^{-j}\le\sum_{j=2}^{\eps n}\big(2k\eps\gamma^{-\newD}\deltnonspanning^{-1}\big)^j\le\tfrac{4k^2\eps^2\gamma^{-2\newD}\deltnonspanning^{-2}}{1-2k\eps\gamma^{-\newD}\deltnonspanning^{-1}}\le 8k^2\eps^2\gamma^{-2\newD}\deltnonspanning^{-2}\,,
 \end{align*}
 where we use the bound $\binom{\eps n}{j}\le (\eps n)^j$ and sum the resulting geometric series. 
\end{claimproof}
This completes the proof of Lemma~\ref{lem:A:probint}.
\end{proof}

We now deduce a similar result for intervals of any length starting at a given time $t_0\ge 0$, again following~\cite[Lemma~\ref{PAPERpackingdegenerate.lem:probvtx}]{DegPack}. 

 \begin{lemma}\label{lem:A:probvtx} Assume Setting~\ref{set:graphs}. Then the following holds for any $\alpha_0\le\alpha\le\alpha_{2n}$ and all sufficiently large $n$. Suppose
that $G$ is a graph on $[v(G)]$ with $v(G)\le(1-\deltnonspanning) n$ such that $\LEFTDEG(x)\le \newD$ for each $x\in V(G)$, and $H$ is an $(\alpha,2\newD+3)$-quasirandom graph
with $n$ vertices and $p\binom{n}{2}$ edges, with $p\ge\gamma$.
  Let $0\le t_0< t_1\le v(G)$ and let $k\in[\newD]$. Let $\histens$ be a history ensemble of \RandomEmbedding{} up to time $t_0$, and suppose that $\Prob[\histens]\ge n^{-4}$. Then for any distinct vertices $u_1,\dots,u_k\in V(H)$ such that $u_1,\dots,u_k\notin\im\psi_{t_0}$ we have
   \[\Prob[u_1,\dots,u_k\notin\im\psi_{t_1}|\histens]=(1\pm100C\newD\alpha\deltnonspanning^{-1})\big(\tfrac{n-1-t_1}{n-t_0}\big)^k\,.\]
 \end{lemma}
 \begin{proof}
We divide the interval $(t_0,t_1]$ into $j:=\lceil (t_1-t_0)/\eps n\rceil$ intervals, all but the last of length $\eps n$.
  Let $\histens_0:=\histens$. Let, for each $1\le i<j$, the set $\histens_i$ be the embedding histories up to time $t_0+i\eps n$ of \RandomEmbedding{} which extend histories in $\histens_{i-1}$ and are such that $u_1,\dots,u_k\notin \im\psi_{t_0+i\eps n}$. Let $\histens_j$ be the embedding histories up to time $t_1$ extending those in $\histens_{j-1}$ such that $u_1,\dots,u_k\notin\im\psi_{t_1}$. Thus we have
\[\Prob[u_1,\dots,u_k\notin\im\psi_{t_1}]=\Prob[\histens_j]\,.\]
Finally, for each $1\le i\le j$, let the set $\histens'_{i-1}$ consist of all histories in $\histens_{i-1}$ such that the $(C\alpha,2\newD+3)$-diet condition holds for $(H,\im\psi_{t_0+(i-1)\eps n})$ and the probability that the $(\eps,C\alpha,t_0+1+(i-1)\eps n)$-cover condition fails, conditioned on $\psi_{t_0+(i-1)\eps n}$, is at most $n^{-3}$. In other words, $\histens'_i$ is the subset of $\histens_i$ consisting of typical histories, satisfying the conditions of Lemma~\ref{lem:A:probint}. 
  
  We now determine $\Prob[\histens_j]$, and in particular we show inductively that $\Prob[\histens_i]\ge n^{-5}$ for each $i$. Observe that for any time $t$, the probability (not conditioned on any embedding) that either the $(C\alpha,2\newD+3)$-diet condition fails for $(H,\im\psi_\ell)$ for some $\ell\le t$ or that the $(\eps,C\alpha,t+1)$-cover condition has probability greater than $n^{-3}$ of failing, is at most $2n^{-6}$ by Lemma~\ref{lem:propRandEmb}. In other words, for each $i$ we have $\Prob[\histens_i\setminus\histens'_i]\le 2n^{-6}$. Thus by Lemma~\ref{lem:A:probint}
  (with $t'=t_0+(i-1)\eps n + 1$) we have
  \begin{align*}
   \Prob[\histens_i]&=\big(1-(1\pm 10C\alpha)\tfrac{\eps k n}{n-t_0-(i-1)\eps n -1}\big)\Prob[\histens'_{i-1}]\pm 2n^{-6}\\
   &=\big(1-(1\pm 10C\alpha)\tfrac{\eps kn}{n-t_0-(i-1)\eps n -1}\big)\big(\Prob[\histens_{i-1}]\pm 2n^{-6}\big)\pm 2n^{-6}\\
   &=\big(1-(1\pm 20C\alpha)\tfrac{\eps kn}{n-t_0-(i-1)\eps n}\big)\Prob[\histens_{i-1}]\,,
  \end{align*}
  where the final equality uses the lower bound $\Prob[\histens_{i-1}]\ge n^{-5}$. Similarly, we have $\Prob[\histens_j]=\big(1\pm(1+20C\alpha)\tfrac{\eps k n}{n-t_1}\big)\Prob[\histens_{j-1}]$.
  
  Putting these observations together, we can compute $\Prob[\histens_j]$:
  \[\Prob[\histens_j]=\Prob[\histens]\big(1\pm(1+20C\alpha)\tfrac{\eps k n}{n-t_1}\big)\prod_{i=1}^{j-1}\Big(1-(1\pm 20C\alpha)\tfrac{\eps k n}{n-t_0-(i-1)\eps n}\Big)\,.\]
 Observe that the approximation $\log(1+x)=x\pm x^2$ is valid for all sufficiently small $x$. In particular, since $n-t_0-(i-1)\eps n\ge n-t_1\ge \deltnonspanning n$ and by choice of $\eps$, for each $i$ we have
\[\log \Big(1-(1\pm 20C\alpha)\tfrac{\eps k n}{n-t_0-(i-1)\eps n}\Big)=-(1\pm 30C\alpha)\tfrac{\eps k n}{n-t_0-(i-1)\eps n}\,.\]
Thus we obtain
 \begin{align}
\nonumber  \log\Prob[\histens_j]&=\log\Prob[\histens]\pm (1+30C\alpha)\tfrac{\eps k n}{n-t_1}-\sum_{i=1}^{j-1}(1\pm 30C\alpha)\tfrac{\eps k n}{n-t_0-(i-1)\eps n}\\
\nonumber  &=\log\Prob[\histens]\pm 2k\deltnonspanning^{-1}\eps-(1\pm 40C\alpha)\int_{x=t_0}^{t_0+(j-1)\eps n}\tfrac{k}{n-x}\, \mathrm{d}x\\
\nonumber  &=\log\Prob[\histens]\pm 2k\deltnonspanning^{-1}\eps-(1\pm 50C\alpha)k\big(\log (n-t_0)-\log(n-1-t_1)\big)\\
 \label{eq:histensnotsmall} &=\log\Prob[\histens]+\log\big(\tfrac{n-1-t_1}{n-t_0}\big)^k\pm 2k\deltnonspanning^{-1}\eps\pm 50Ck\alpha\log\deltnonspanning^{-1}\,,
 \end{align}
 where we use $t_1\le n-\deltnonspanning n$, and we justify that the integral and sum are close by observing that for each $i$ in the summation, if $t_0+(i-1)\eps n\le x\le t_0+i\eps n$ then we have
 \[\tfrac{1}{n-t_0-(i-1)\eps n}\le\tfrac{1}{n-x}\le\tfrac{1}{n-t_0-i\eps n}\le(1+\alpha)\tfrac{1}{n-t_0-(i-1)\eps n}\,,\]
 where the final inequality uses $n-t_0-i\eps n\ge n-t_1\ge \deltnonspanning n$ and the choice of $\eps$. By choice of $\eps$ and since $k\le \newD$, this gives the lemma statement. Furthermore, \eqref{eq:histensnotsmall}, and the fact $t_1\le n-\deltnonspanning n$, imply that $\Prob[\histens_j]\ge n^{-5}$. Since the $\histens_i$ form a decreasing sequence of events, the same bound holds for each $\histens_i$.
 \end{proof}
 
 Now, we prove~\eqref{eq:A:probvtxnoim} and~\eqref{eq:A:probvtxnoim2}, by copying the proof of~\cite[Lemma~\ref{PAPERpackingdegenerate.lem:probedge}]{DegPack} and modifying it very slightly (simply by removing the vertex $y$ and references to it).

\begin{proof}[Proof of~\eqref{eq:A:probvtxnoim} and~\eqref{eq:A:probvtxnoim2}]
 Let $x\in V(G)$ and three distinct vertices $v,v',v''\in V(H)$ be given. We begin with~\eqref{eq:A:probvtxnoim}.
 Let $z_1,\dots,z_k$ be the vertices of $\LNBH(x)$ in increasing order.
Define time intervals using $z_1,\ldots,z_k,x$ as separators: $I_0=[1,z_1-1]$, $I_1=[z_1+1,z_2-1]$, \ldots, $I_k=[z_k+1,x-1]$, $I_{k+1}=[x+1,v(G)]$.

We now define a nested collection of events, the first being the trivial (always satisfied) event and the last being the event $\{x\AlgMap v,v'\notin\im\psi_{v(G)}\}$, whose probability we wish to estimate. These events are simply that we have not yet (by given increasing times in \RandomEmbedding) made it impossible to have $\{x\AlgMap v,v'\notin\im\psi_{v(G)}\}$. We will see that we can estimate accurately the probability of each successive event, conditioned on its predecessor.

Let $\histens'_{-1}$ be the trivial (always satisfied) event. If $\histens'_{i-1}$ is defined, we let $\histens_i$ be the event that $\histens'_{i-1}$ holds intersected with the event that 
\begin{enumerate}[label=(A\arabic*)]
 \item\label{e:A1} (if $i\le k$:) no vertex of $G$ in the interval $I_i$ is mapped to $v$ or $v'$, or
\item\label{e:A2} (if $i=k+1$:) no vertex of $G$ in the interval $I_{k+1}$ is mapped to $v'$.
\end{enumerate}
In other words, $\histens_i$ is the event that we have not covered $v$ or $v'$ in the interval $I_i$. It turns out that we do not need to know anything else about the embeddings in the interval $I_i$.

If $\histens_i$ is defined, we let $\histens'_i$ be the event that $\histens_i$ holds and that
 \begin{enumerate}[label=(B\arabic*)]
 \item\label{e:Ch1} (if $i<k$:) we have the event $z_{i+1}\AlgMap \NBH_H(v)\setminus\{v'\}$,
 \item\label{e:Ch2} (if $i=k$:) we have the event $x\AlgMap v$.
 \end{enumerate}
Again, in order for $\{x\AlgMap v,v'\notin\im\psi_{v(G)}\}$ to occur we obviously need that a neighbour of $x$ is embedded to a neighbour of $v$ and so on, hence the above conditions.

By definition, we have $\histens_{k+1}=\{x\AlgMap v,v'\notin\im\psi_{v(G)}\}$. Since we have $\histens'_i\subseteq\histens_i\subseteq\histens'_{i-1}$ for each $i$ and $\histens'_{-1}$ is the sure event, we see
\begin{align}
\Prob\left[x\AlgMap v,v'\notin\im\psi_{v(G)}\right]&=\frac{\Prob[\histens_0]}{\Prob[\histens'_{-1}]}\cdot\prod_{i=0}^{k}
\frac{\Prob[\histens'_i]}{\Prob[\histens_i]}\cdot\frac{\Prob[\histens_{i+1}]}{\Prob[\histens'_{i}]}\\
&=\Prob\left[\histens_0\:|\:\histens'_{-1}\right]\prod_{i=0}^{k}
\Prob\left[\histens'_i\:|\:\histens_i\right]\Prob\left[\histens_{i+1}\:|\:\histens'_{i}\right]
\label{eq:nicelygrouped}
\;.
\end{align}
Thus, we need to estimate the factors in~\eqref{eq:nicelygrouped}. This is done in the two claims below. In each claim we assume $\Prob[\histens'_i],\Prob[\histens_i]>n^{-4}$. This assumption is justified, using an implicit induction, since the smallest of all the events we consider is $\histens_{k+1}$, whose probability according to the following~\eqref{eq:probxy} is bigger than $n^{-4}$.

\begin{claim}\label{cl:ng1}
We have
\begin{equation*}
\prod_{i=0}^{k+1} \Prob\left[\histens_i\:|\:\histens'_{i-1}\right]=(1\pm 200C\newD\alpha\newdelta^{-1})^{2k+4}\cdot \frac{(n-x)(n-v(G))}{n^2}\;.
\end{equation*}
\end{claim}
\begin{claimproof}[Proof of Claim~\ref{cl:ng1}]
By definition of~\ref{e:A1}, for each $i=0,\ldots, k$, we have 
\begin{equation}\label{eq:cancelme}
\Prob\left[\histens_i\:|\:\histens'_{i-1}\right]=(1\pm 200C\newD\alpha\newdelta^{-1})\cdot \frac{(n-1-\max(I_{i}))^2}{(n-\min(I_i)+1)^2}
\end{equation}
by the $2$-vertex case of Lemma~\ref{lem:A:probvtx}, with $\histens=\histens'_{i-1}$.   Note that looking at two consecutive indices $i$ and $i+1$ in~\eqref{eq:cancelme} we have cancellation of the former numerator and the latter denominator, $n-1-\max(I_{i})=n-\min(I_{i+1})+1$. Thus,
\begin{equation}\label{eq:cancelmex}
\prod_{i=0}^{k} \Prob\left[\histens_i\:|\:\histens'_{i-1}\right]=(1\pm 200C\newD\alpha\newdelta^{-1})^{2k+2}\cdot \frac{(n-x)^2}{n^2}\;.
\end{equation}
To express $\Prob\left[\histens_{k+1}\:|\:\histens'_{k}\right]$, by definition of~\ref{e:A2} we have to repeat the above replacing the $2$-vertex case of Lemma~\ref{lem:A:probvtx} with the $1$-vertex case. We get that
\begin{equation}\label{eq:cancelmexx}
 \Prob\left[\histens_{k+1}\:|\:\histens'_k\right]=(1\pm 200C\newD\alpha\newdelta^{-1})^2\cdot \frac{n-v(G)}{n-x}\;,
\end{equation}
where we absorb the factor $\tfrac{n-1-v(G)}{n-v(G)}=1\pm\tfrac{1}{\deltnonspanning n}$, coming from the $1$-vertex case of Lemma~\ref{lem:A:probvtx} with $t_0=x$ and $t_1=v(G)$, into the error term; this is possible for large $n$ since $\alpha\ge\alpha_0$. Putting~\eqref{eq:cancelmex} and~\eqref{eq:cancelmexx} together, we get the statement of the claim.
\end{claimproof}

\begin{claim}\label{cl:ng2}
We have
$$
\prod_{i=0}^{k}\Prob\left[\histens'_i\:|\:\histens_i\right]=(1\pm 100C\alpha)^{k+1}\cdot \frac{1}{n+1-x}\;.
$$
\end{claim}
\begin{claimproof}[Proof of Claim~\ref{cl:ng2}]
Suppose that we have embedded up to vertex $\max(I_i)$, and that $\histens_i$ holds. The probability of the event $\histens'_i$ depends on which of the cases~\ref{e:Ch1} and~\ref{e:Ch2} applies. When $\histens'_i$ is defined using~\ref{e:Ch1} then the probability $\Prob[\histens'_i|\histens_i]$ is equal to $\Prob[\{z_{i+1}\AlgMap \NBH_H(v)\setminus\{v'\}\}|\histens_i]$. Let $X:=\NBH_H\big(\psi(\LNBH_G(z_{i+1}))\big)\setminus\im\psi_{z_{i+1}-1}$ be the set of vertices in $H$ to which we could embed $z_{i+1}$, given the embedding of all vertices before $z_{i+1}$. Suppose that the $(C\alpha,2\newD+3)$-diet condition holds for $(H,\im\psi_{z_{i+1}-1})$. Then we have 
\begin{align*}
\Prob\left[z_{i+1}\AlgMap \NBH_H(v)\setminus\{v'\}|\histens_i\right]&=\frac{|(\NBH_H(v)\setminus\{v'\})\cap X|}{|X|}=\frac{|\NBH_H(v)\cap X|\pm 2}{|X|}
\\
&=
\frac{(1\pm C\alpha)p^{1+\LEFTDEG(z_{i+1})}(n-(z_{i+1}-1))\pm 2}
{(1\pm C\alpha)p^{\LEFTDEG(z_{i+1})}(n-(z_{i+1}-1))}
=(1\pm 4C\alpha)p\;,
\end{align*}
where the last line uses the $(C\alpha,2\newD+3)$-diet condition for $(H,\im\psi_{z_{i+1}-1})$ twice, in the denominator with the set $\psi(\LNBH(z_{i+1}))$ and in the numerator with the set $\{v\}\cup\psi(\LNBH(z_{i+1}))$. Recall that we assume the event $\histens_i$, and so we have $v\notin \im \psi_{z_{i+1}-1}$. Therefore, the set $\{v\}\cup \psi(\LNBH_G(z_{i+1}))$ has indeed size $1+\LEFTDEG(z_{i+1})$.

Let us now deal with the term $\Prob\left[\histens'_k\:|\:\histens_k\right]$, which corresponds to~\ref{e:Ch2}. Suppose that $\histens_k$ holds. In particular, $\LNBH(x)$ is embedded to $\NBH_H(v)$.
Suppose first that the $(C\alpha,2\newD+3)$-diet condition for $(H,\im\psi_{x-1})$ holds. With this, conditioning on the embedding up to time $x-1$, the probability of embedding $x$ to $v$ is $(1\pm 2C\alpha)p^{-\LEFTDEG(x)}\tfrac{1}{n+1-x}$.

Thus, letting $\cF$ be the event that the $(C\alpha,2\newD+3)$-diet condition fails for $(H,\im\psi_t)$ for some $t\in[v(G)]$ we have
\[
\prod_{i=0}^{k}\Prob\left[\histens'_i\:|\:\histens_i\right]
=\Big(\big((1\pm 4C\alpha)p\big)^{k} 
\cdot
(1\pm 2C\alpha)p^{-\LEFTDEG(x)}\tfrac{1}{n+1-x}
\Big) \pm \Prob[\cF]\;.
\]
We have $\Prob[\cF]\le 2n^{-9}$ by Lemma~\ref{lem:propRandEmb}. Thus we obtain, using that $\LEFTDEG(x)=|\LNBH(x)|=k$, so that the two powers of $p$ cancel,
\begin{equation*}
\prod_{i=0}^{k}\Prob\left[\histens'_i\:|\:\histens_i\right]
=(1\pm 4C\alpha)^{k+1}\cdot\tfrac{1}{n+1-x}\pm 2n^{-9}\;,
\end{equation*}
and the claim follows, since $\tfrac{1}{n+1-x}\ge\tfrac1n$ and so the additive error $2n^{-9}$ is absorbed by the multiplicative one.
\end{claimproof}

Plugging Claims~\ref{cl:ng1} and~\ref{cl:ng2} into~\eqref{eq:nicelygrouped}, we get 
\begin{equation}\label{eq:probxy}
\Prob\left[x\AlgMap v,v'\notin\im\psi_{v(G)}\right]=(1\pm 500C\newD\alpha\newdelta^{-1})^{3\newD+5}\cdot\tfrac{(n-x)(n-v(G))}{n^2(n-x+1)}\;,
\end{equation}
and~\eqref{eq:A:probvtxnoim} follows.

For~\eqref{eq:A:probvtxnoim2}, we use the same approach. The only difference is that we define events
\begin{enumerate}[label=(A'\arabic*)]
 \item\label{e:A12} (if $i\le k$:) no vertex of $G$ in the interval $I_i$ is mapped to $v$, $v'$ or $v''$, or
\item\label{e:A22} (if $i=k+1$:) no vertex of $G$ in the interval $I_{k+1}$ is mapped to $v'$ or $v''$.
\end{enumerate}
and if $\histens_i$ is defined, we let $\histens'_i$ be the event that $\histens_i$ holds and that
 \begin{enumerate}[label=(B'\arabic*)]
 \item\label{e:Ch12} (if $i<k$:) we have the event $z_{i+1}\AlgMap \NBH_H(v)\setminus\{v',v''\}$,
 \item\label{e:Ch22} (if $i=k$:) we have the event $x\AlgMap v$.
 \end{enumerate}
 The calculations are almost identical: in particular the calculations and results for the B' events are verbatim the same, while for the A' events we use respectively the $3$-vertex and $2$-vertex cases of Lemma~\ref{lem:A:probvtx} rather than the $2$- and $1$-vertex cases; but the calculations go through without further trouble. We omit the details.
\end{proof}

We need a slight strengthening of~\cite[Lemma~\ref{PAPERmanyleaves.lem:T}]{ABCT:PackingManyLeaves}. The idea of this lemma is that we pause $\PackingProcess'$ after it has packed $s$ graphs, and fix a set of vertices $T$ in the common neighbourhood of some $v_1,\dots,v_k$ in the graph $H_s$. We then allow the process to continue up to some point $s'$, and of course the common neighbourhood of $v_1,\dots,v_k$ shrinks; this lemma states that its intersection with $T$ shrinks proportionally. The difference to~\cite[Lemma~\ref{PAPERmanyleaves.lem:T}]{ABCT:PackingManyLeaves} is that we look at the common neighbourhood of a set of vertices rather than one vertex, and we give a different bound on the set size, but the proof is very similar.
For this, let us denote by $p_s$ the density of $H_s$.

\begin{lemma}\label{lem:T}
  Assume Setting~\ref{set:graphs} and let $s,s'\in\cJ\cup\cK$ with $s<s'$. Consider the following experiment. Suppose $v_1,\dots,v_k$ are fixed vertices of $H_0$, with $k\le 2\newD+3$.
  Run $\PackingProcess'$ with input
  $(G'_{s''})_{{s''}\in\cJ\cup\cK}$ and $H_0$ up to and including the
  embedding of~$G'_s$. Then fix $T\subset \NBH_{H_s}(v_1,\dots,v_k)$ with
  $|T|\ge\frac12 \gamma^{\newD}\deltnonspanning^{\newD} n$, and continue $\PackingProcess'$ to perform the embedding
  of $G'_{s+1},\dots,G'_{s'}$.
  
  The probability that $\PackingProcess'$ fails before embedding $G'_{s'}$, or $H_i$ fails to be $(\alpha_i,2\newD+3)$-quasirandom for some $1\le i\le s'$, or we have \[\big|T\cap \NBH_{H_{s'}}(v_1,\ldots,v_k)\big|=(1\pm \alpha_{s'})\big(\tfrac{p_{s'}}{p_s}\big)^k|T|\,,\]
  is at least $1-n^{-C}$.
\end{lemma}
To prove this, we follow the proof of~\cite[{Lemma~\ref{PAPERmanyleaves.lem:T}}]{ABCT:PackingManyLeaves}, using~\eqref{eq:A:probstar} in place of~\cite[{Lemma~\ref{PAPERmanyleaves.lem:probedge}}]{ABCT:PackingManyLeaves}.
\begin{proof}[Proof of Lemma~\ref{lem:T}]
 For $s\le i\le s'$, we define the event~$\cE_i$ that $\PackingProcess'$ does not fail before embedding $G'_i$, and~$H_j$ is $(\alpha_j,2\newD+3)$-quasirandom for each $1\le j\le i$,
  and $|T\cap \NBH_{H_j}(v_1,\dots,v_k)|=(1\pm \alpha_j)\big(\frac{p_j}{p_s}\big)^k|T|$ for each $s\le j\le i$. If the event in the lemma statement fails to occur, then there must exist some $s\le i< s'$ such that $\cE_i$ occurs and
  \[|T\cap \NBH_{H_{i+1}}(v_1,\dots,v_k)|\neq(1\pm \alpha_{i+1})\big(\tfrac{p_{i+1}}{p_s}\big)^k|T|\,.\]
 It suffices to show that each of these bad events occurs with probability at most $n^{-C-2}$, since then the union bound over the at most $\tfrac32 n$ choices of $i$ gives the lemma statement. This is an estimate we can obtain using Corollary~\ref{cor:freedm}. We now fix $s\le i<s'$ and prove the desired estimate.
 
 Suppose $s\le j\le i$, and let $Y_j:=|\NBH_{H_{j}}(v_1,\dots,v_k)\cap T\setminus \NBH_{H_{j+1}}(v_1,\dots,v_k)|$ count the number of stars with leaves $v_1,\dots,v_k$ and centre in $T$, at least one edge of which is used for the embedding of $G'_{j+1}$. Then we have $\big|T\cap \NBH_{H_{i+1}}(v_1,\dots,v_k)\big|=|T|-\sum_{j=s}^iY_j$, and what we want to do is argue that the sum of random variables is concentrated. To that end, suppose $\hist$ is a history of $\PackingProcess'$ up to time $j$ such that $H_j$ is $(\alpha_j,2\newD+3)$-quasirandom and $|T\cap \NBH_{H_j}(v_1,\dots,v_k)|=(1\pm \alpha_j)\big(\frac{p_j}{p_s}\big)^k|T|$. Then we have
 \[\Exp\big[Y_j\,\big|\,\hist\big]=\big(1\pm1000C\alpha_j\newdelta^{-1}\big)^{4\newD+2}p_j^{-1}n^{-2}\cdot 2ke(G'_{j+1})\cdot \big(1\pm\alpha_j\big)\big(\tfrac{p_j}{p_s}\big)^k|T|\,,\]
  where we use linearity of expectation: the first factor is, by~\eqref{eq:A:probstar}, the probability that at least one edge of a given star with leaves $v_1,\dots,v_k$ and centre in $T$ in $H_j$ is used in the embedding of $G'_{j+1}$, and the second factor is the number of such stars. Simplifying, we obtain
  \begin{align*}
   \Exp\big[Y_j\,\big|\,\hist\big]&=\frac{2kp_j^{k-1}e(G'_{j+1})|T|}{p_s^k n(n-1)}\pm\frac{10^5\newdelta^{-1}C\newD^3|T|}{p_s
      n}\alpha_{j}\,,
  \end{align*}
  where for the error term we use the upper bound $e(G'_{j+1})\le \newD n$ and the fact $p_j/p_s\le 1$, and where we have replaced $n^{-2}$ by $\tfrac{1}{n(n-1)}$, which is convenient below and, for large $n$, changes the main term by less than the error term. Let
  \[\tilde{\mu}:=\sum_{j=s}^i\frac{2kp_j^{k-1}e(G'_{j+1})|T|}{p^k_s n(n-1)}\quad\text{and}\quad\tilde\nu:=\sum_{j=s}^i\frac{10^5\newdelta^{-1}C\newD^3|T|}{p_s
      n}\alpha_{j}\]
and observe that $\tilde{\mu}\le k p_s^{-1} |T|\le k \dStageA^{-1} n$ and
$\tilde{\nu} \le \frac{10^5\newdelta^{-1}C\newD^3|T|}{p_s}\alpha_{i} < \frac{|T|}{10^3}$ since $p_s\ge \dStageA$ and
by the definition of $\alpha_j$. 
  
  We trivially have $0\le Y_j\le k\Delta(G'_{j+1})\le  kcn/\log n$. So what Corollary~\ref{cor:freedm}\ref{cor:freedm:tails}, with $\tilde{\rho}=\eps n$, gives us is that
  \[\Prob\left[\cE_i\text{ and }\sum_{j=s}^iY_j\neq \tilde{\mu}\pm(\tilde\nu+\eps n)\right]<2\exp\left(-\tfrac{\eps^2n^2}{4k^2c\dStageA^{-1}n^2/\log n}\right)<n^{-C-2}\,,\]
  where we use the upper bound $\tilde{\mu}+\tilde{\nu}+\tilde{\rho}\le 2k\dStageA^{-1}n$ for the first inequality and the choice of $c$ as well as $\eps < \frac{1}{C}$ for the second. This is the probability bound we wanted. We now simply need to show that if
  \[\sum_{j=s}^iY_j=\tilde{\mu}\pm(\tilde\nu+\eps n)\]
  then we have
  \[|T\cap \NBH_{H_{i+1}}(v_1,\dots,v_k)|=(1\pm \alpha_{i+1})\big(\tfrac{p_{i+1}}{p_s}\big)^k|T|\,.\]
  In order to estimate $\tilde{\mu}$, observe that $e(G'_{j+1})=(p_j-p_{j+1})\binom{n}{2}$, and that for every $x,h\in[0,1]$ we have $(x+h)^k-x^k=kh(x+h)^{k-1}\pm 2^kh^2$. Using the latter equality with $x=p_{j+1}$ and $h=p_j-p_{j+1}$, and using $(p_j-p_{j+1})\binom{n}{2}\le \newD n$, we see \[kp_j^{k-1}e(G'_{j+1})=\big(p_j^k-p_{j+1}^k\big)\binom{n}{2}\pm 2^{k+2}\newD^2\,.\]
  Using this we see
  \[|T|-\tilde{\mu}=|T|\Big(1-\sum_{j=s}^i \big( \big(p_j^k-p_{j+1}^k\big)p_s^{-k}\pm \frac{2^{k+2}\newD^2}{p_s^k\binom{n}{2}} \big)\Big)=\tfrac{p_{i+1}^k}{p_s^k}|T|\pm 2^{k+4}\newD^2p_s^{-k}\,,\]
  since $i+1\le 2n$. So what remains is to argue $\tilde\nu+\eps n+2^{k+4}\newD^2p_s^{-k}<\alpha_{i+1}\left(\tfrac{p_{i+1}}{p_s}\right)^k|T|$. Since $\alpha_j=\frac{\newdelta}{10^8C\newD}\exp\left(\frac{10^8C^2\newD^3\newdelta^{-1}\gamma^{-4\newD-6}(j-2n)}{n}\right)$ is increasing in $j$, we have
  \begin{align}
  \begin{split}
  \label{eq:sum:alpha}
    \sum_{j=s}^i\alpha_j&\le \int_s^{i+1}\alpha_j\,\mathrm{d}j\le\int_{-\infty}^{i+1}\alpha_j\,\mathrm{d}j\\
    &=\Big[\frac{\newdelta}{10^8C\newD}\cdot\frac{n}{10^8C^2\newD^3\newdelta^{-1}\gamma^{-4\newD-6}}\cdot\exp\Big(\frac{10^8C^2\newD^3\newdelta^{-1}\gamma^{-4\newD-6}(j-2n)}{n}\Big)\Big]_{j=-\infty}^{i+1}\\
    &=\frac{\newdelta n}{10^8C^2\newD^3}\gamma^{4\newD+6}\alpha_{i+1}\,.
   \end{split}
  \end{align}
  It follows that
\begin{align*}
	  \tilde\nu+\eps n+2^{k+4}\newD^2p_s^{-k}
	  	& \le \frac{10^5\newdelta^{-1}C\newD^3|T|}{p_s n}\cdot \frac{\newdelta n}{10^8C^2\newD^3}\gamma^{4\newD+6}\alpha_{i+1}+\eps n+2^{k+4}\newD^2p_s^{-k}\\
	  	& \le\tfrac{\alpha_{i+1}}{1000}\gamma^{4\newD+6}\cdot\tfrac{1}{p_s}|T|+\eps n+2^{k+4}\newD^2p_s^{-k}\,.	
\end{align*}	  
    Finally, since $p_{i+1},p_s\ge\gamma$ and $k\le 2\newD+3$, by choice of $\eps$, because $|T|\ge\tfrac12 \deltnonspanning^{\newD}\gamma^{\newD} n$, and because $n$ is sufficiently large, we conclude $\tilde\nu+\eps n+2^{k+4}\newD^2p_s^{-k}\le\alpha_{i+1}\big(\tfrac{p_{i+1}}{p_s}\big)^k|T|$ as desired.
\end{proof}

The next lemma, which is almost the same as~\cite[Lemma~\ref{PAPERmanyleaves.lem:S}]{ABCT:PackingManyLeaves}, shows that in a run of \RandomEmbedding\ (considered as  step $s$ in $\PackingProcess'$) the image of $V(G'_{s})$ covers a predictable amount of any previously given reasonably large vertex set $S$. For completeness, we copy the proof from~\cite{ABCT:PackingManyLeaves} and make the appropriate small modifications.

\begin{lemma}\label{lem:S}
  Assume Setting~\ref{set:graphs} and let $s\in\cJ\cup\cK$.  Run $\PackingProcess'$ with input
  $(G'_{s''})_{{s''}\in\cJ\cup\cK}$ and $H_0$ up to just before the
  embedding of~$G'_s$. Suppose $v(G'_s)\le(1-\deltnonspanning)n$. We define for each $t\in\RR$ 
  \[\beta_{t} =2\alpha_{s-1}\exp\left(\tfrac{1000\newD\newdelta^{-2}\gamma^{-2\newD-10}t}{n}\right)\;.\]
  Then fix any $S\subset V(H_{s-1})$ with $|S|\ge\frac12\gamma^{\newD}\deltnonspanning^{\newD} n$, and
  let $\PackingProcess'$ perform the embedding of~$G'_s$. With
  probability at least $1-n^{-2\newD-18}$, either $H_{s-1}$ is not $(\alpha_{s-1},2\newD+3)$-quasirandom, or \RandomEmbedding\ fails to construct the embedding of $G'_s$, or the embedding of $G'_s$ fails to have the $(\eps,20\newD\beta_t,t)$-cover condition for some $1\le t\le v(G'_s)+1-\eps n$, or for some $1\le t\le v(G'_s)$ the pair $\big(H_{s-1},\phi_s([t])\big)$ does not have the $(\beta_t,2\newD+3)$-diet condition, or we have
  \[|S\setminus\im\phi^\mathbf{A}_{s}|=(1\pm C'\alpha_{s})\tfrac{n-v(G'_s)}{n}|S|\,.\]
\end{lemma}
We will refer to the first four of the events mentioned above as `bad events', and the fifth (the equation) as the good event.
\begin{proof}
  Fix~$s\in\cJ\cup\cK$, and condition on $H_{s-1}$. If $H_{s-1}$ is not $(\alpha_{s-1},2\newD+3)$-quasirandom, or the embedding of $G'_s$ fails, or the embedding of $G'_s$ fails to have the $(\eps,20\newD\beta_t,t)$-cover condition for some $1\le t\le v(G'_s)+1-\eps n$, or for some $1\le t\le v(G'_s)$ the pair $\big(H_{s-1},\phi_s([t])\big)$ does not have the $(\beta_t,2\newD+3)$-diet condition, then one of the bad events of this lemma occurs. So it is enough to show that the probability that none of these bad events occurs but $|S\setminus\im\phi^\mathbf{A}_{s}|\neq(1\pm C'\alpha_{s})\tfrac{n-v(G'_s)}{n}|S|$, conditioned on $H_{s-1}$, is at most $n^{-2\newD-18}$. This is what we will now do, so we suppose that $H_{s-1}$ is $(\alpha_{s-1},2\newD+3)$-quasirandom. Consider the run of \RandomEmbedding{} which finds an embedding $\psi_{v(G'_s)}$ of $G'_s$.
  
Define $S_0=S$, and for $i=1,\dots,\tau$ with
  $\tau=\big\lceil\frac{v(G'_s)}{\eps n}\big\rceil$ set
  $S_i=S_{i-1}\setminus \im \psi_{i\eps n}$. Since $S_\tau\subseteq S\setminus\im\phi^\mathbf{A}_{s}\subseteq S_{\tau-1}$, it is enough to show both $|S_{\tau-1}|$ and $|S_\tau|$ are likely to be in the claimed range. Recall that $v(G'_s)\le(1-\deltnonspanning)n<(1-\newdelta-\eps)n$, so that \RandomEmbedding\ does embed vertices (namely isolated vertices added to $G'_s$) long enough to create $S_\tau$. Since $|S_{\tau-1}|$ and $|S_{\tau}|$ differ by at most $\eps n$, we will focus on estimating $|S_\tau|$.
  
 Given $0\le j\le \tau-1$, either one of the bad events of the lemma occurs, or $\psi_{j\eps n}$ exists and $(H_{s-1},\im\psi_{j\eps n})$ has the $(\beta_{j\eps n},2\newD+3)$-diet condition. Lemma~\ref{lem:24}, with input $T=S_{j}=S\setminus\im\psi_{j\eps n}$, then states that conditioning on $\psi_{j\eps n}$, with probability at least $1-2n^{-2\newD-19}$ either $\psi_{v(G_s')}$ does not have the $(\eps,20\newD\beta_{j\eps n},j\eps n)$-cover condition, or we have
\begin{equation}\label{eq:A:lemS}
    \left| \left\{x:~ j\eps n\le x<(j+1)\eps n:~
    		\psi_{(j+1)\eps n}(x) \in S\setminus \im \psi_{j\eps n} 
    		\right\}\right| =
    		(1\pm 40\newD\beta_{j\eps n})\frac{|S\setminus \im \psi_{j\eps n}|\eps
      n}{n-j\eps n}~ .
\end{equation}
  In particular, since the failure of the cover condition is one of the bad events of the lemma, taking the union bound over the at most $\tau$ choices of $j$, with probability at least $1-2\tau n^{-2\newD-19}\ge 1-n^{-2\newD-18}$, either one of the bad events of the lemma occurs or we have~\eqref{eq:A:lemS} for each $0\le j<\tau$. Let us suppose that the latter is the case. Then we have for each $1\le j\le\tau$
  \[|S_{j}|=|S_{j-1}|\Big(1-\frac{(1\pm 40\newD\beta_{(j-1)\eps n})\eps
      }{1-(j-1)\eps}\Big)\,,\]
  and hence
  \[|S_\tau|=|S|\prod_{i=1}^{\tau}\Big(1-\frac{(1\pm 40\newD\beta_{(i-1)\eps n})\eps
      }{1-(i-1)\eps}\Big)\,.\]
  In order to evaluate this product, observe that
  \[1-\frac{(1\pm 40\newD\beta_{i\eps n})\eps}{1-i\eps}
    =\frac{1-(i+1)\eps}{1-i\eps}\pm\frac{40 \newD\beta_{i\eps n}\eps}{1-i\eps}
    = \frac{1-i\eps-\eps}{1-i\eps}\Big(1\pm\frac{40\newD\beta_{i\eps
        n}\eps }{1-(i+1)\eps}\Big)\,,
  \]
  and therefore
  \[|S_\tau|=|S|\prod_{i=0}^{\tau-1}\frac{1-i\eps-\eps}{1-i\eps}
    \cdot\Big(1\pm\frac{40\newD\beta_{i\eps n}\eps}{1-(i+1)\eps}\Big)
    =|S|(1-\tau\eps) \prod_{i=0}^{\tau-1}\Big(1\pm\frac{40\newD\beta_{i\eps n}\eps}{1-(i+1)\eps}\Big)\,.
  \]
  
By the definition of $\tau$ we have
  $\frac{v(G'_s)}{\eps n}\le\tau\le\frac{v(G'_s)}{\eps n}+1$ and hence $(1-\tau\eps)=\tfrac{n-v(G'_s)}{n}(1\pm\frac\eps\deltnonspanning)$.
Moreover, we obtain that
\begin{align*}
	\sum_{i=0}^{\tau-1}\frac{40\newD\beta_{i\eps n}
		\eps}{1-(i+1)\eps}
	& \le \frac{40\newD\eps }{1-\tau \eps} 
	\sum_{i=0}^{\tau -1} \beta_{i\eps n}
	\le \frac{80\newD\eps n}{n-v(G'_s)} 
	\sum_{i=0}^{\tau -1} \beta_{i\eps n} \\
	&\le\frac{80\newD}{\deltnonspanning n}\int_{0}^{\tau}\eps n\beta_{i\eps n}\,\mathrm{d}i
    \le\frac{80\newD}{\deltnonspanning n}\int_{0}^{\tau\eps n}\beta_{x}\,\mathrm{d}x \\
    &\le
    \frac{80\newD}{\deltnonspanning\cdot 1000
      \newD\newdelta^{-2}\gamma^{-2\newD-10}}\beta_{\tau \eps n}\\
    &\le\frac{80\newD}{\deltnonspanning\cdot 1000
      \newD\newdelta^{-2}\gamma^{-2\newD-10}}\beta_{(1-\deltnonspanning + \eps )n}
    \le\frac12 C'\alpha_{s-1}
\end{align*}  
   since
  $\beta_{(1-\deltnonspanning+\eps)n}=2\alpha_{s-1}\exp(1000\newD\newdelta^{-2}\gamma^{-2\newD-10}(1-\deltnonspanning+\eps))$
  and \[C'=10^4\cdot\frac{40\newD}{\newdelta}\exp(1000\newD\newdelta^{-2}\gamma^{-2\newD-10})\,.\]

So, since $\prod_i(1\pm x_i)=1\pm2\sum_i x_i$
  as long as $\sum_i x_i < \frac{1}{100}$
  and since $\frac{1}{2}C'\alpha_s < \frac{1}{100}$, 
  we get
  \begin{equation*}\begin{split}
      |S_\tau|&=|S|(1-\tau\eps)\Big(1\pm2\sum_{i=0}^{\tau-1}\frac{40\newD\beta_{i\eps
          n}\eps}{1-(i+1)\eps}\Big)
      =|S|(1-\tau\eps)\Big(1\pm\frac{80\newD\eps}{1-\tau\eps}\sum_{i=0}^{\tau-1}\beta_{i\eps
        n}\Big) \\
      &=|S|\Big(1-\tau\eps\pm 80\newD\eps\sum_{i=0}^{\tau-1}\beta_{i\eps
        n}\Big)
      =|S|\tfrac{n-v(G'_s)}{n}\Big(1\pm\frac{\eps}{\deltnonspanning} \pm\frac{80\newD\eps}{\deltnonspanning}\sum_{i=0}^{\tau-1}\beta_{i\eps
        n}\Big) \\
      &=|S|\tfrac{n-v(G'_s)}{n}\Big(1\pm\frac12\alpha_s \pm
      \frac{1}{2}C'\alpha_s\Big)\,,
  \end{split}\end{equation*}
  where for the last equation we use that $\eps\le\alpha_0\newdelta^2\gamma\le\frac12\alpha_s\deltnonspanning$. It follows that  
  \[|S\setminus\im\phi^\mathbf{A}_s|=
  |S_{\tau}|\pm \eps n = 
  |S|\tfrac{n-v(G'_s)}{n}\big(1\pm\tfrac12\alpha_s\pm\tfrac12C'\alpha_s\big)\pm\eps n=\big(1\pm C'\alpha_s\big)\tfrac{n-v(G'_s)}{n}|S|\,,\]
  as desired.
\end{proof}

\subsection{Proofs of the remaining properties of Lemma~\ref{lem:StageA}}
We are now in a position to prove~\ref{A:vertex-covering}. The idea is that first we will show the equation holds with $\HStageA$ replaced by $H_{s-1}$ and $\dStageA$ replaced with the density $d_{s-1}$ of $H_{s-1}$, and then Lemma~\ref{lem:T} tells us it is likely that~\ref{A:vertex-covering} holds. For the first part, we argue (using the diet condition) that each path-end vertex has about a $d_s$ chance of being embedded to the neighbourhood of $v$, conditional on the previous embedding history, and then Corollary~\ref{cor:freedm} gives the conclusion we want.
\begin{proof}[Proof of~\ref{A:vertex-covering}]
 Fix $v\in V(H)$ and $s$. Suppose that $\PackingProcess'$ has not failed before creating $H_{s-1}$, and furthermore $H_{s-1}$ is $(\alpha_{s-1},2\newD+3)$-quasirandom. Let $d_{s-1}$ be the density of $H_{s-1}$.

 Let $X=\bigcup_{P\in \SpecPaths_s}\{\leftpath_0(P), \rightpath_0(P)\}$. Let $\hist_{x-1}$ denote the history of \RandomEmbedding\ embedding $G'_s$ into $H_{s-1}$, up to but not including the embedding of $x$; let $\psi_{x-1}$ denote the partial embedding of $G'_s$ at this time.
 
 Suppose that $\hist_{x-1}$ is such that $(H_{s-1},\im\psi_{x-1})$ satisfies the $(C\alpha_{s-1},2\newD+3)$-diet condition, and furthermore that $v\notin\im\psi_{x-1}$. Then we embed $x$ uniformly to $\NBH_{H_{s-1}}(y_1,\dots,y_\ell)\setminus\im\psi_{x-1}$, where the $y_i$ are the back-neighbours of $x$ in $G'_s$. This is a set of size $(1\pm C\alpha_{s-1})(n+1-x)d_{s-1}^\ell$ by the diet condition. Again by the diet condition, and since $v$ is not one of the $y_i$ (it is not in $\im\psi_{x-1}$), the number of these vertices which are in $\NBH_{H_{s-1}}(v)$ is $(1\pm C\alpha_{s-1})(n+1-x)d_{s-1}^{\ell+1}$. So the probability that $x$ is embedded to $\NBH_{H_{s-1}}(v)$, conditioned on $\hist_{x-1}$, is $(1\pm 3C\alpha_{s-1})d_{s-1}$.
 
 We now apply Corollary~\ref{cor:freedm}\ref{cor:freedm:tails} to estimate the total number of vertices of $X$ embedded to $\NBH_{H_{s-1}}(v)$. We take $\tilde{\eta}=3C\alpha_{s-1}$ and $\tilde{\mu}=|X|d_{s-1}$, $R=1$, and $\cE$ the event that each $\hist_{x-1}$ satisfies the diet condition and $v\notin\im\phi_s^\mathbf{A}$. The conclusion is that with probability at least $1-2\exp\big(-|X|d_{s-1}\cdot9C^2\alpha_{s-1}^2/4\big)\ge 1-n^{-10}$, we have $\bar{\cE}$ or 
 \[\big|\big\{x\in X:\phi_s^\mathbf{A}(x)\in \NBH_{H_{s-1}}(v)\big\}\big|=(1\pm 6C\alpha_{s-1})|X|d_{s-1}\,.\]
 
 By Lemma~\ref{lem:propRandEmb}, the probability that $v\notin\im\phi_s^\mathbf{A}$ and we are not in $\cE$ is at most $2n^{-9}$. So with probability at least $1-3n^{-9}$, either $v\in\im\phi_s^\mathbf{A}$ or
 \[\big|\big\{x\in X:\phi_s^\mathbf{A}(x)\in \NBH_{H_{s-1}}(v)\big\}\big|=(1\pm 6C\alpha_{s-1})|X|d_{s-1}\,.\]
  We condition on this good event. Now if $v\notin\im\phi_s^\mathbf{A}$ then we have $\NBH_{H_{s-1}}(v)=\NBH_{H_s}(v)$, and furthermore $d_{s-1}\binom{n}{2}=d_s\binom{n}{2}\pm \newD n$, so if $v\notin\im\phi_s^\mathbf{A}$ we have
 \[\big|\big\{x\in X:\phi_s^\mathbf{A}(x)\in \NBH_{H_{s}}(v)\big\}\big|=(1\pm 10C\alpha_{s-1})|X|d_{s}\,.\]
 
 We now apply Lemma~\ref{lem:T}, with $T$ the subset of $\NBH_{H_s}(v)$ lying in $X$, and $s'$ the final stage of $\PackingProcess'$. The conclusion is that, with probability at least $1-n^{-C}$, we have
 \[\big|\big\{x\in X:\phi_s^\mathbf{A}(x)\in \NBH_{\HStageA}(v)\big\}\big|=(1\pm\gamma^{-1}\alpha_{s'})(1\pm 10C\alpha_{s-1})|X|\dStageA\,,\]
 which by~\eqref{eq:A:smallerrors} is what we need for~\ref{A:vertex-covering}. Taking the union bound over choices of $s$ and $v$ and the various bad events, we conclude that a.a.s.~\ref{A:vertex-covering} holds.
\end{proof}

Putting Lemmas~\ref{lem:T} and~\ref{lem:S} together, we can show that~\ref{A:megaquasirandomness1} holds with very high probability for any given $S$ and $T$, and so by the union bound holds with high probability.

\begin{proof}[Proof of~\ref{A:megaquasirandomness1}]
 Let $S\subset V(H)$ and $T\subset\cJ\cup \cK$, with $|S|,|T|\le \newD=2D$, be given. Let $T=\{t_1,\dots,t_{|T|}\}$, and for convenience let $t_0=0$ and $t_{|T|+1}=|\cJ|+|\cK|+1$. We aim to show that the following holds for each $0\le i\le|T|$:
 \begin{equation}\label{eq:A:megaquasi}\Big|\NBH_{H_{t_i}}(S)\setminus\bigcup_{j=1}^i\imA(t_j)\Big|=\big(1\pm2(i+1)(1+C')\alpha_{t_i}\big)\Big(\tfrac{e(H_{t_i})}{\binom{n}{2}}\Big)^{|S|} n\Big(\prod_{j=0}^i\big(1-\tfrac{|\imA(t_j)|}{n}\big)\Big)\,.\end{equation}
 Observe that the $i=0$ case of~\eqref{eq:A:megaquasi} is true because $H_0$ is $(\alpha_0,2\newD+3)$-quasirandom. Observe furthermore that for any $i$, the set on the left-hand side of~\eqref{eq:A:megaquasi} is fixed before $\PackingProcess'$ begins (if $i=0$) or immediately after completing the embedding of $G'_{t_i}$ (if $i\ge 1$). Finally, note that the LHS of~\eqref{eq:A:megaquasi} has size at least $\tfrac12\cdot \gamma^{\newD} \deltnonspanning^{\newD}n$, which verifies the size requirement for Lemma~\ref{lem:T}.
 
 It follows that for each $0\le i\le |T|$, we have the setup for Lemma~\ref{lem:T}, which tells us that with probability at least $1-n^{-C}$ either $\PackingProcess'$ fails before embedding $G'_{t_{i+1}-1}$, or $H_j$ fails to be $(\alpha_j,2\newD+3)$-quasirandom for some $1\le j\le t_{i+1}-1$ (these are the bad events of Lemma~\ref{lem:T}), or we have the good event
 \begin{align*}
  \Big|&\NBH_{H_{t_{i+1}-1}}(S)\setminus\bigcup_{j=1}^i\imA(t_j)\Big|\\
  &=(1\pm\alpha_{t_{i+1}-1})\big(1\pm2(i+1)(1+C')\alpha_{t_i}\big)\Big(\tfrac{e(H_{t_{i+1}-1})}{\binom{n}{2}}\Big)^{|S|} n\Big(\prod_{j=1}^i\big(1-\tfrac{|\imA(t_j)|}{n}\big)\Big)\\
  &=\big(1\pm2(i+1)(1+C')\alpha_{t_i}\pm\tfrac32\alpha_{t_{i+1}-1}\big)\Big(\tfrac{e(H_{t_{i+1}-1})}{\binom{n}{2}}\Big)^{|S|} n\Big(\prod_{j=0}^i\big(1-\tfrac{|\imA(t_j)|}{n}\big)\Big)\,.
 \end{align*}
 Observe that if $i=|T|$, by~\eqref{eq:A:smallerrors} this is the desired statement of~\ref{A:megaquasirandomness1}. If $0\le i\le |T|-1$, then we have the setup for Lemma~\ref{lem:S} (as above, the size requirement is given by the fact that the LHS of the above equation is at least $\tfrac12\gamma^{\newD}\deltnonspanning^{\newD} n$), analysing the embedding of $G'_{t_{i+1}}$, which tells us that with probability at least $1-n^{-2\newD-18}$, either one of the following bad events (as listed in Lemma~\ref{lem:S}) occurs, or the final good event occurs. The bad events are: $H_{t_{i+1}-1}$ is not $(\alpha_{t_{i+1}-1},2\newD+3)$-quasirandom, or \RandomEmbedding\ fails to construct the embedding of $G'_{t_{i+1}}$, or $\big(H_{t_{i+1}-1},\phi_{t_{i+1}}^\mathbf{A}([j])\big)$ does not have the $(\beta_j,2\newD+3)$-diet condition for some $1\le j\le v(G'_{t_{i+1}})$, or $\phi^\mathbf{A}_{t_{i+1}}$ fails to have the $(\eps,20\newD\beta_j,j)$-cover condition for some $1\le j\le v(G'_{t_{i+1}})+1-\eps n$. Here $\beta_j=\beta_j\big(\alpha_{t_{i+1}-1}\big)$ is as defined in Lemma~\ref{lem:S}. The good event is that we have
 \begin{align*}
  \Big|&\NBH_{H_{t_{i+1}-1}}(S)\setminus\bigcup_{j=1}^{i+1}\imA(t_j)\Big|\\
  &=(1\pm C'\alpha_{t_{i+1}})\big(1\pm2(i+1)(1+C')\alpha_{t_i}\pm\tfrac32\alpha_{t_{i+1}-1}\big)\Big(\tfrac{e(H_{t_{i+1}-1})}{\binom{n}{2}}\Big)^{|S|} n\Big(\prod_{j=1}^{i+1}\big(1-\tfrac{|\imA(t_j)|}{n}\big)\Big)\\
  &=\big(1\pm2(i+1)(1+C')\alpha_{t_i}\pm\tfrac32\alpha_{t_{i+1}-1}\pm\tfrac32C'\alpha_{t_{i+1}}\big)\Big(\tfrac{e(H_{t_{i+1}-1})}{\binom{n}{2}}\Big)^{|S|} n\Big(\prod_{j=1}^{i+1}\big(1-\tfrac{|\imA(t_j)|}{n}\big)\Big)\,.
 \end{align*}
 Finally, observe that at most $\tfrac{cn}{\log n}$ edges are removed from any vertex of $S$ when $G'_{t_{i+1}}$ is embedded, so we get
 \begin{align*}
  \Big|\NBH_{H_{t_{i+1}}}(S)\setminus\bigcup_{j=1}^{i+1}\imA(t_j)\Big|&=\Big|\NBH_{H_{t_{i+1}-1}}(S)\setminus\bigcup_{j=1}^{i+1}\imA(t_j)\Big|\pm |S|\Delta\\
  &=\big(1\pm2(i+2)(1+C')\alpha_{t_{i+1}}\big)\Big(\tfrac{e(H_{t_{i+1}-1})}{\binom{n}{2}}\Big)^{|S|} n\Big(\prod_{j=1}^{i+1}\big(1-\tfrac{|\imA(t_j)|}{n}\big)\Big)
 \end{align*}
 since $\alpha_j$ is increasing with $j$ and since $n$ is sufficiently large. For $0\le i\le |T|-1$, this shows that with probability at least $1-n^{-C}-n^{-2\newD-18}$, if~\eqref{eq:A:megaquasi} holds for $i$ then it holds for $i+1$ or one of the bad events of Lemma~\ref{lem:T} or~\ref{lem:S} occurs. Putting these together, we conclude that either a bad event of Lemma~\ref{lem:T} or~\ref{lem:S} occurs, or~\ref{A:megaquasirandomness1} holds for the given $S$ and $T$ with probability at least $1-(|T|+1)n^{-C}-|T|n^{-2\newD-18}$, and so taking the union bound over the at most $(2n)^{|S|+|T|}$ choices of $S$ and $T$, we conclude that either one of the bad events of Lemma~\ref{lem:T} or~\ref{lem:S} occurs, or~\ref{A:megaquasirandomness1} holds, with probability at least $1-n^{-15}$. Now we already established in Section~\ref{sec:A:maintquasi} that with high probability $H_{s-1}$ is $(\alpha_{s-1},2\newD+3)$-quasirandom for each $s\in\cK\cup\cJ$, and Lemma~\ref{lem:propRandEmb}~\ref{pfA:nofail}--\ref{pfA:diet} state that the remaining bad events with high probability do not occur. So we conclude that with high probability~\ref{A:megaquasirandomness1} holds, as desired.
\end{proof}

For the two $\SpecLeaves$ conditions, we do not need to obtain accurate estimates and will not try to do so. The first condition,~\ref{A:parleavesSpread}, states that no vertex $v$ of $\HStageA$ can be the embedded neighbour of too many vertices of $\bigcup_{s\in\cK}\SpecLeaves_s$. We separate this estimate into two parts: the contribution due to vertices of some $G_s$ which are adjacent to at most $n^{0.8}$ vertices in $\SpecLeaves_s$, and the rest. The contribution of the first class is not too large by Corollary~\ref{cor:freedm} (the expected contribution is in fact constant, and by definition no one vertex embedding contributes more than $n^{0.8}$), and we will show that it is very unlikely that $20$ or more vertices of the second class (of which there are at most $2\newD n^{0.2}$) are embedded to $v$; since any one vertex contributes at most $\Delta(G_s)\le\tfrac{cn}{\log n}$ we obtain the claimed bound.

\begin{proof}[Proof of~\ref{A:parleavesSpread}]
 Let $\cE$ be the event that $\PackingProcess'$ does not fail and that for each $s$ the graph $H_{s-1}$ is $(\alpha_{s-1},2\newD+3)$-quasirandom and for each time in the embedding of $G'_s$ by \RandomEmbedding\ we have the $(C\alpha_{s-1},2\newD+3)$-diet condition. Note that $\cE$ a.a.s.\ occurs.

 Fix $v\in V(H_0)$. For any given $s\in\cK$ and $y\in V(G'_s)$, let $\hist_y$ denote the history of $\PackingProcess'$ up to immediately before embedding $y$. Suppose that $\hist_y$ is in $\cE$. The probability that we embed $y$ to $v$, conditioning on $\hist_y$, is at most $2\dStageA^{-\newD}\deltnonspanning^{-1}n^{-1}$, since we embed $y$ to a set of size at least $\tfrac12\dStageA^{\newD}\deltnonspanning n$ by the diet condition. Let $\omega(y)$ denote the number of neighbours of $y$ in $\SpecLeaves_s$. Then for~\ref{A:parleavesSpread} we want to show
 \[\sum_y \omega(y)\mathbbm{1}_{y\to v}\le\tfrac{20cn}{\log n}\,.\]
 
 We split this sum into two parts. Let $S$ denote the vertices $y$ such that $\omega(y)\le n^{0.8}$. Note that $\sum_{y\in S}\omega(y)\le 2\newD n$, since $\big|\bigcup_{s\in\cK}\SpecLeaves_s\big|\le 2n$ and
$\DParity\le \newD$.
  So Corollary~\ref{cor:freedm} tells us that the probability that $\cE$ occurs and
 \[\sum_{y\in S}\omega(y)\mathbbm{1}_{y\to v}\ge \frac{4\newD n}{\dStageA^{\newD}\deltnonspanning n}+n^{0.9}\]
 is at most $2\exp\big(-\tfrac{n^{1.8}}{4n^{0.8}n^{0.9}}\big)\le\exp\big(-n^{0.05}\big)$.
 
 Now consider the vertices of $\bar{S}$, i.e.\ those $y$ such that $\omega(y)>n^{0.8}$. There are at most $2\newD n^{0.2}$ such vertices. Since conditional probabilities multiply, the probability that $\cE$ occurs and $20$ or more of the vertices of $\bar{S}$ are mapped to $v$ is at most
 \[\binom{2\newD n^{0.2}}{20}\big(2\dStageA^{-\newD}\deltnonspanning^{-1}n^{-1}\big)^{20}\le n^{-10}\,.\]
 
 We conclude that with probability at least $1-n^{-9}$, either we are not in $\cE$ or
 \[\sum_y \omega(y)\mathbbm{1}_{y\to v}\le\frac{4\newD n}{\dStageA^{\newD}\deltnonspanning n}+n^{0.9}+19\tfrac{cn}{\log n}\le\tfrac{20cn}{\log n}\,.\]
 
 Taking the union bound over choices of $v$, and since as observed $\cE$ a.a.s.\ occurs, we conclude that a.a.s.\ we have~\ref{A:parleavesSpread}.
\end{proof}

For~\ref{A:SpecLeavesCover} we have to work a bit harder. Fix a vertex $v\in V(H_0)$. In order for a vertex $x\in\SpecLeaves_s$ to contribute to the sum in~\ref{A:SpecLeavesCover} (we say $x$ is $\SpecLeaves$ covering $v$; more generally if we have a given embedding of some guest graphs and a given host graph $H$ in the middle of the process, we say $x$ is $\SpecLeaves$ covering $v$ if the neighbours of $x$ are embedded to neighbours of $v$ in the current host graph and $v$ is not in the image of the guest graph containing $x$), two things have to happen. First, we have to embed $\NBH_{G_s}(x)$ to vertices of $\NBH_{H_{s-1}}(v)$ and not embed any vertex of $G'_s$ to $v$ (we say we create $x$ covering $v$ if this happens). Second, we have to not use any of the edges from $v$ to the embedded neighbours of $x$ in the remaining stages of $\PackingProcess'$ (we say we lose $x$ covering $v$ if this happens). It turns out to be rather inconvenient to deal with these two things together.

So what we do is first show that we create a reasonably large number of $\SpecLeaves$ covering $v$ in total. We do this by considering the embedding of each $G'_s$ one after another. It is enough to estimate, for each $s\in\cK$ and $x\in\SpecLeaves_s$, the probability that $\NBH_{G_s}(x)$ is embedded to $\NBH_{H_{s-1}}(v)$ and $v\notin\im\phi_s$, conditioned on the history of $\PackingProcess'$ up to the point where $H_{s-1}$ is defined; this is an easy consequence of the diet condition. Then Corollary~\ref{cor:freedm} gives what we want, because the contribution made by embedding any one $G'_s$ is at most $|\SpecLeaves_s|\le\tfrac{cn}{\log n}$. We should note that it is at this point where we need this bound on $|\SpecLeaves_s|$. But for our approach this bound is necessary: if $|\cK|\le \tfrac1{10}\log n$, it is likely that there would be a vertex $v$ in the image of all the $\phi_s$ for $s\in\cK$, and hence~\ref{A:SpecLeavesCover} would fail.

We then try to argue separately that we do not lose too many of these $\SpecLeaves$ as edges are removed to form $\HStageA$. For this argument, we do not work $G'_s$ by $G'_s$, but rather consider separately each vertex embedding that causes an edge to be removed at $v$. There are two types of such vertex embedding: for each $s\in\cK\cup\cJ$ the vertex embedded to $v$ (if there is one) removes at most $\newD$ edges, and thereafter each neighbour removes one edge. In both cases, by~\ref{A:parleavesSpread} the number of $\SpecLeaves$ we lose is at most $\tfrac{20D cn}{\log n}$, which is small enough for Corollary~\ref{cor:freedm} to give useful probability bounds. For this sketch we ignore the first type, which turns out not to cause much trouble. We will use the diet condition to give an upper bound on the expected number of $\SpecLeaves$ covering $v$ we lose each time we use an edge at $v$, which is proportional to the current number of $\SpecLeaves$ covering $v$.

We separate two cases. First, the number of $\SpecLeaves$ covering $v$ never reaches $\zeta n$ (where $\zeta >0$ is some not-too-small constant). We argue that this case is unlikely to occur. In this case, we can use Corollary~\ref{cor:freedm} and our expectation estimate to show that the total number of $\SpecLeaves$ covering $v$ that we lose is likely small.  But this gives us a contradiction: at the end of the process we have less than $\zeta n$ $\SpecLeaves$ covering $v$, and we lost in total few $\SpecLeaves$, yet we created many $\SpecLeaves$ covering $v$.

Second, there is some first time $s'$ in this process where there are $\zeta n$ $\SpecLeaves$ covering $v$; since the first case is unlikely, this case is likely to occur. We ignore any further creation of $\SpecLeaves$ covering $v$, and focus on how fast these $\zeta n$ $\SpecLeaves$ are lost while the remaining graphs are embedded.  We argue (using Corollary~\ref{cor:freedm}  and our expectation estimate) that the number of edges we have to remove at $v$ in order to lose $\tfrac{\zeta n}{2}$ $\SpecLeaves$ is $\Theta(n)$. Once we removed half of the $\SpecLeaves$, we lose $\SpecLeaves$ at half the rate, so it takes as many edges to remove the next $\tfrac{\zeta n}{4}$ $\SpecLeaves$, and so on. Since we remove less than $n$ edges at $v$ (because $H_0$ is an $n$-vertex graph) this argument necessarily stops after a constant number of iterations, and this gives us~\ref{A:SpecLeavesCover}.

\begin{proof}[Proof of~\ref{A:SpecLeavesCover}]
 Let $\cE$ be the event that $\PackingProcess'$ does not fail, that for each $s$ the graph $H_{s-1}$ is $(\alpha_{s-1},2\newD+3)$-quasirandom, and that at each time in using \RandomEmbedding\ to embed $G'_s$ the $(C\alpha_{s-1},2\newD+3)$-diet condition holds. Note that $\cE$ a.a.s.\ occurs.
 
 Fix $v\in V(H_0)$, $s\in\cK$ and $x\in\SpecLeaves_s$. Suppose that $H_{s-1}$ is $(\alpha_{s-1},2\newD+3)$-quasirandom. We aim to find a lower bound for the probability that $\NBH_{G_s}(x)$ is mapped by \RandomEmbedding\ to $\NBH_{H_{s-1}}(v)$ and furthermore $v\notin\im\phi_s$. Observe that \RandomEmbedding\ embeds in total less than $n$ vertices of $G'_s$. We will now define a coupling with  \RandomEmbedding, then justify that the coupling succeeds if the diet condition holds, and finally explain what this coupling is for.
 
 At each time \RandomEmbedding\ embeds a vertex $y$ of $G'_s$ which is not a neighbour of $x$, we do the following to generate an auxiliary random variable taking value in $\{0,1\}$, with probability $2\dStageA^{-D}\deltnonspanning^{-1}n^{-1}$ of getting $1$. Before embedding $y$, we calculate the probability $p$ that $y\to v$. If $p>2\dStageA^{-D}\deltnonspanning^{-1}n^{-1}$, we say the coupling fails and sample the auxiliary random variable from the desired distribution. Otherwise, we embed $y$, and if $y\to v$ occurs we set the random variable to $1$. If $y$ is not embedded to $v$, we set the random variable to $1$ with probability $q$ such that $p+q(1-p)=2\dStageA^{-D}\deltnonspanning^{-1}n^{-1}$, and otherwise to $0$.
 At each time \RandomEmbedding\ embeds a vertex $y$ of $G'_s$ which is a neighbour of $x$, we generate an auxiliary random variable taking value in $\{0,1\}$ with $\tfrac12\dStageA$ probability of getting a zero, such that when the coupling succeeds the following occurs. If $y\not\to \NBH_{H_{s-1}}(v)$ and $v\notin\im\psi_{y-1}$ then we get $1$. To do this we use the same process as above.
 
 We claim that if at each step we have the $(C\alpha_{s-1},2\newD+3)$-diet condition, then this coupling succeeds. To see that when the diet condition holds, the probability of embedding a non-neighbour $y$ of $x$ to $v$ is at most $2\dStageA^{-D}\deltnonspanning^{-1}n^{-1}$, observe that $y$ has at most $D$ left-neighbours, whose unused common neighbourhood at the time we come to embed $y$ has size at least $\tfrac12\dStageA^{D}\deltnonspanning n$, and at most one of these vertices is $v$. To see that when the diet condition holds the probability of embedding a neighbour $y$ of $x$ to a neighbour of $v$ is at least $\tfrac12\dStageA$, we apply the diet condition twice, once to estimate the number of unused common neighbours of $\psi_{y-1}\big(\LNBH(y)\big)$, and once to estimate the size of the subset of the unused common neighbours which are in addition neighbours of $v$. If $v$ is in the image of $\psi_{y-1}$ there is nothing to prove; if not, then $v$ is not one of the vertices $\psi_{y-1}\big(\LNBH(y)\big)$ and hence we conclude the desired lower bound $\tfrac12\dStageA$ on the probability of getting a zero.
 
  By construction, if in the coupled independent Bernoulli random variables we never get $1$, then necessarily we create $x$ covering $v$. Recalling that at most $n$ vertices of $G'_s$ are embedded, and that $x$ has exactly $\DParity\le D$ neighbours in $G_s$, the probability of none of these Bernoulli random variables taking the value $1$ is at least
 \[\big(1-2\dStageA^{-D}\deltnonspanning^{-1}n^{-1}\big)^n\cdot\big(\tfrac12\dStageA\big)^{D}\ge\exp\big(-4\dStageA^{-D}\deltnonspanning^{-1}\big)2^{-D}\dStageA^{D}\ge\exp\big(-8D\dStageA^{-D}\deltnonspanning^{-1}\big)\,,\]
 and the probability of the coupling failing---that is, of the diet condition at some point not holding---is at most $3n^{-9}$, so the probability that \RandomEmbedding\ creates $x$ covering $v$ is at least $\exp\big(-10D\dStageA^{-D}\deltnonspanning^{-1}\big)$.
 
 By linearity of expectation, the expected number of vertices of $\SpecLeaves_s$ which \RandomEmbedding\ creates covering $v$, conditioning on $H_{s-1}$, is at least $|\SpecLeaves_s|\exp\big(-10D\dStageA^{-D}\deltnonspanning^{-1}\big)$. Summing this over all $s\in\cK$ we obtain at least $n\exp\big(-10D\dStageA^{-D}\deltnonspanning^{-1}\big)$, since there are in total at least $n$ vertices in $\bigcup_{s\in\cK}\SpecLeaves_s$ (cf.~Definition~\ref{def:family}\ref{enu:defSpecLeaves}). Finally Corollary~\ref{cor:freedm} tells us that, since no $\SpecLeaves_s$ has size more than $\tfrac{cn}{\log n}$, with probability at least $1-n^{-C}$ either $\cE$ does not occur, or in total we create at least
 \begin{equation}\label{eq:A:SpecLeavesCover:Lots}
  \tfrac12n\exp\big(-10D\dStageA^{-D}\deltnonspanning^{-1}\big)
 \end{equation}
 $\SpecLeaves$ covering $v$. This completes the first step of the sketch.

 Our next task is to prove the aforementioned expectation estimate. We need to be a little more careful than outlined: we need to estimate the expected number of a given collection of $\SpecLeaves$ created at $v$ which we lose in each step. To that end, suppose that we have completed the embedding of the first $s-1$ guest graphs, and have therefore fixed these embeddings and the graph $H_{s-1}$. We are interested in the vertices $x$ of $\bigcup_{s'\in\cK}\SpecLeaves_{s'}$ with $s'\le s-1$ such that $\NBH_{G_{s'}}(x)$ is embedded to $\NBH_{H_{s-1}}(v)$ and $v\notin\im\phi_{s'}$. Let $X$ be an arbitrary subset of these vertices.
 
 We define two kinds of \emph{loss event}. The first kind of loss event is that we begin embedding $G'_s$, up to and including embedding a vertex to $v$ (if this occurs). The second is that, having already embedded a vertex $x$ of $G'_s$ to $v$, we embed a neighbour $y$ of $x$. We estimate the expectation for each separately, though we will eventually get a single upper bound for both kinds. What will be important is that the total number of loss events is not too large. There are at most $|\cJ\cup \cK|\le 2n$ loss events of the first kind, while every time a loss event of the second kind occurs we use an edge at $v$ and hence there are at most $n-1$ such events.
 
 We need some notation. For $u\in \NBH_{H_{s-1}}(v)$, let
 \[\omega(u;X):=\sum_{x\in X}\mathbbm{1}_{\phi_{s'}^{-1}(u)x\in E(G_{s'})}\]
 where given $x$ in the summation we set $s'$ such that $x\in V(G_{s'})$. Thus removing the edge $vu$ from $H_{s-1}$ means that we lose exactly $\omega(u;X)$ vertices from $X$ (formally we should say that these are $\SpecLeaves$ lost at $v$; we will drop the `at $v$' since $v$ is fixed in what follows); removing several edges $vu_i$ causes us to lose at most $\sum_i\omega(u_i;X)$ vertices from $X$ (it might be that two edges both cause us to lose the same $\SpecLeaves$ vertex; we will not try to account for this). Furthermore we have 
 \[|X|\le\sum_{u\in \NBH_{H_{s-1}}(v)}\omega(u;X)\le D|X|\,,\]
 since each vertex of $X$ has exactly $\DParity$ embedded neighbours, and $1\le\DParity\le D$.
 
 To begin with, we handle the first kind of loss event. Given $H_{s-1}$, suppose that $H_{s-1}$ is $(\alpha_{s-1},2\newD+3)$-quasirandom. Conditioning on $H_{s-1}$, the expected number of vertices of $X$ lost due to the first kind of loss event is at most the expected number of vertices of $X$ lost due to the entire embedding of $G'_s$. By linearity of expectation, this is
\[\sum_{u\in \NBH_{H_{s-1}}(v)}\omega(u;X)\Prob\big[uv\text{ used in embedding $G'_s$}\big|H_{s-1}\big]\,,\]
which we can estimate by~\eqref{eq:A:probstar} with $k=1$; we see that this conditional expectation is at most
\[ 2D|X|\tfrac{2e(G'_s)}{\dStageA n^2}\le 4D^2\dStageA^{-1}|X|n^{-1}\,.\]

We now handle the second kind of loss event. Suppose the vertex $x$ of $G'_s$ has been embedded to $v$; let $y$ be a neighbour of $x$ not yet embedded. We condition on $\psi_{y-1}$, the embedding of all vertices up to $y-1$ of $G'_s$. Suppose that $\im\psi_{y-1}$ satisfies the $(C\alpha_{s-1},2\newD+3)$-diet condition. Then we will embed $y$ uniformly to a set of size at least $\tfrac12\dStageA^{D}\deltnonspanning n$, whose total weight according to $\omega(\cdot;X)$ is at most $D|X|$. In expectation, the weight of our choice for $\psi_y(y)$ is thus at most $\tfrac{2D|X|}{\dStageA^{D}\deltnonspanning n}$.

In either case, provided $\cE$ occurs the expected number of $\SpecLeaves$ in $X$ which we lose in a single loss event (conditioned on the packing up to the previous loss event) is at most
\begin{equation}\label{eq:A:expSpecLeavesLoss}
 4D^2\dStageA^{-D}\deltnonspanning^{-1}n^{-1}|X|\,.
\end{equation}

This is the expectation estimate discussed in the above sketch. We want to use this with Corollary~\ref{cor:freedm}, for which we also need to know the maximum number of $\SpecLeaves$ in $X$ which we lose at $v$ in a single loss event. Observe that when a single edge $uv$ is used, we lose exactly those $\SpecLeaves$ in $X$ which have a neighbour embedded to $u$. This is upper bounded by the total number of $\SpecLeaves$ (whether in $X$ or not) which have a neighbour embedded to $u$, and this is what is upper bounded by~\ref{A:parleavesSpread} (for the vertex $u$) by $\tfrac{20cn}{\log n}$. In a single loss event we use at most $D$ edges at $v$ (in an event of the first kind; one of the second kind always uses one edge) and hence the maximum number of $\SpecLeaves$ in $X$ we lose in a single loss event is at most $\tfrac{20D cn}{\log n}$ as promised.

We have now reached the case separation from the sketch. To begin with, we suppose there is no point in $\PackingProcess'$ when the number of $\SpecLeaves$ covering $v$ reaches
\[\zeta n:=\tfrac1{100}D^{-2}\dStageA^{D}\deltnonspanning\exp\big(-10D\dStageA^{-D}\deltnonspanning^{-1}\big)n\,.\]
In this case, for each loss event we take $X$ to be the set of all $\SpecLeaves$ covering $v$ at the current time. Provided $\cE$ occurs, the sum of conditional expectations from~\eqref{eq:A:expSpecLeavesLoss} is bounded above by
\[3n\cdot 4D^2\dStageA^{-D}\deltnonspanning^{-1}n^{-1}\zeta n=12D^2\dStageA^{-D}\deltnonspanning^{-1}\zeta n\,.\]
By Corollary~\ref{cor:freedm}, with probability at least $1-n^{-C}$ either we are not in $\cE$ or (the good event) the total number of $\SpecLeaves$ lost at $v$ is at most $24D^2\dStageA^{-D}\deltnonspanning^{-1}\zeta n$. However, if the good event occurs, after $\PackingProcess'$ finishes we have at most $\zeta n$ $\SpecLeaves$ covering $v$, and have lost at most $24D^2\dStageA^{-D}\deltnonspanning^{-1}\zeta n$, so at most $25D^2\dStageA^{-D}\deltnonspanning^{-1}\zeta n$ $\SpecLeaves$ covering $v$ can have been created. As we showed in~\eqref{eq:A:SpecLeavesCover:Lots}, the probability of $\cE$ occurring and this few $\SpecLeaves$ being created is at most $n^{-C}$. So the probability of $\cE$ occurring and being in this case is at most $2n^{-C}$. This completes the proof that this case is unlikely.

We now handle the likely case from the sketch: there is some time in $\PackingProcess'$ at which we have at least $\zeta n$ $\SpecLeaves$ covering $v$. Recall that when we embed a guest graph and create $\SpecLeaves$ covering $v$, in particular we do not embed any vertex to $v$ and so do not lose any $\SpecLeaves$. Thus there is a stage $s'$ such that when we have completed embedding $G_{s'}$, we have in total at least $\zeta n$ $\SpecLeaves$ covering $v$. Let $X_0$ be a set of $\zeta n$ $\SpecLeaves$ covering $v$ after stage $s'$. We consider the at most $3n$ loss events from this point on, starting with the embedding of the first vertices of $G_{s'+1}$. Let for each $i$ the set $X_i$ be the set of $\SpecLeaves$ in $X_0$ which are still covering $v$ after the $i$th loss event.

If $\cE$ occurs the expectation of $|X_{i-1}\setminus X_i|$, conditioning on $\PackingProcess'$ up to and including the $(i-1)$st loss event, is at most $4D^2\dStageA^{-D}\deltnonspanning^{-1}n^{-1}|X_{i-1}|$ by~\eqref{eq:A:expSpecLeavesLoss} with $X=X_{i-1}$. To conveniently get an upper bound on this, we split up the loss events into intervals (with random endpoints) as follows: the first interval ends with the first $i$ such that $|X_i|\le \tfrac14\zeta n$, and in general the $j$th interval ends with the first $i$ such that $|X_i|\le 4^{-j}\zeta n$. Thus when $i$ is in the $j$th interval, under $\cE$, the expectation of $|X_{i-1}\setminus X_i|$, conditioned on the packing up to the $(i-1)$st loss event, is at most
\[4D^2\dStageA^{-D}\deltnonspanning^{-1}n^{-1} \cdot 4^{1-j}\zeta n=4D^2\dStageA^{-D}\deltnonspanning^{-1}\cdot 4^{1-j}\zeta\,.\]

We now argue that it is unlikely that any given interval is short and $\cE$ occurs. Deterministically each interval contains at least $\Omega(\log n)>1$ loss events, using the maximum change in each event. However, consider the $\tfrac1{32}D^{-2}\dStageA^{D}\deltnonspanning n$ loss events from the beginning of the $j$th interval. Suppose that $\cE$ occurs; then for each of these loss events $i$ we have $|X_{i-1}|\le 4^{1-j}\zeta n$, and so the sum of the conditional expectations is at most
\[4D^2\dStageA^{-D}\deltnonspanning^{-1}\cdot 4^{1-j}\zeta\cdot \tfrac1{32}D^{-2}\dStageA^{D}\deltnonspanning n=\tfrac{1}{8}\cdot 4^{1-j}\zeta n\,.\]

By Corollary~\ref{cor:freedm}, conditioning on the packing up to and including the end of the $(j-1)$st interval, the probability that $\cE$ occurs and the actual number of $\SpecLeaves$ lost at $v$ in these $\tfrac1{32}D^{-2}\dStageA^{D}\deltnonspanning n$ loss events exceeds $\tfrac16\cdot 4^{1-j}\zeta n$ is at most $n^{-C}$. In particular, if the likely event that it does not exceed this quantity occurs, then all these $\tfrac1{32}D^{-2}\dStageA^{D}\deltnonspanning n$ loss events are in the $j$th interval; we see that the probability within $\cE$ of the $j$th interval containing less than $\tfrac1{32}D^{-2}\dStageA^{D}\deltnonspanning n$ loss events, conditioning on the packing up to the end of the $(j-1)$st interval, is at most $n^{-C}$.

Taking the union bound, the probability that $\cE$ occurs and any of the first $3\cdot 32 D^2\dStageA^{-D}\deltnonspanning^{-1}$ intervals has less than $\tfrac1{32}D^{-2}\dStageA^{D}\deltnonspanning n$ loss events is at most $n^{1-C}$. This is thus a likely event. Assuming this likely event occurs, because there are at most $3n$ loss events $\PackingProcess'$ completes before reaching the end of the last of these intervals. By definition, this means that after the final loss event $t$ we have
\[|X_t|\ge 4^{-3\cdot 32 D^2\dStageA^{-D}\deltnonspanning^{-1}}\zeta n\ge 4^{-200D^2\dStageA^{-D}\deltnonspanning^{-1}} n\,,\]
and since all the vertices of $X_t$ are $\SpecLeaves$ covering $v$ at the end of $\PackingProcess'$, this proves~\ref{A:SpecLeavesCover}.
\end{proof}

To prove~\ref{A:typicality} we use Lemma~\ref{lem:A:probvtx} to estimate the probability that the embedding of a given $G_s$ by \RandomEmbedding\ uses a vertex of $S$. Given this estimate,~\ref{A:typicality} follows easily from Corollary~\ref{cor:freedm} and the union bound.

\begin{proof}[Proof of~\ref{A:typicality}]
 Let $\cE$ be the event that $\PackingProcess'$ succeeds and for each $s$ the graph $H_{s-1}$ is $(\alpha_{s-1},2\newD+3)$-quasirandom.
 
 Given $S_1,S_2\subset V(H_0)$ with $1\le|S_i|\le D=\tfrac12\newD$, $i=1,2$, and $s\in\cJ$, we condition on $H_{s-1}$ and assume we are in $\cE$. By Lemma~\ref{lem:A:probvtx}, with $t_0=0$ and $t_1=v(G_s)$, and $\histens$ the trivial history (which has probability $1$) we can estimate the probability that $\im\phi^\mathbf{A}_s$ is disjoint from $S_1\cup Q$ for every subset $Q$ of $S_2$. By inclusion-exclusion, the probability that $\im\phi^\mathbf{A}_s$ is disjoint from $S_1$ and contains $S_2$ is
 \[\sum_{Q\subset S_2}(-1)^{|Q|}\Prob\big[\im\phi^\mathbf{A}_s\cap (S_1\cup Q)=\emptyset\big]\,,\]
 and substituting in the estimates from Lemma~\ref{lem:A:probvtx} we get the probability
\begin{align*}
 &\Big(\sum_{Q\subset S_2}(-1)^{|Q|}\big(\tfrac{n-1-v(G'_s)}{n}\big)^{|S_1|+|Q|}\Big)\pm 100C\newD\alpha_{s-1}\deltnonspanning^{-1}\Big(\sum_{Q\subset S_2}\big(\tfrac{n-1-v(G'_s)}{n}\big)^{|S_1|+|Q|}\Big)\\
=&\big(\tfrac{n-1-v(G'_s)}{n}\big)^{|S_1|}\Big(1-\tfrac{n-1-v(G'_s)}{n}\Big)^{|S_2|}\pm 100C\newD\alpha_{s-1}\deltnonspanning^{-1}2^{\newD}\big(\tfrac{n-1-v(G'_s)}{n}\big)^{|S_1|}\\
=&\big(1\pm\tfrac18\gamcore\big)\cdot\big(\tfrac{n-v(G'_s)}{n}\big)^{|S_1|}\cdot (\tfrac{v(G'_s)}{n})^{|S_2|}\,,
\end{align*}
where the final line uses our choice $\alpha_{s-1}\ll\gamcore$ and the observation $v(G'_s)>n/2$.
 It follows that
 \[\sum_{s\in\cJ_0}\Prob\big[\im\phi^\mathbf{A}_s\cap S_1=\emptyset,  S_2\subseteq \im\phi^\mathbf{A}_s\big|H_{s-1}\big]=\big(1\pm\tfrac18\gamcore\big)|\cJ_0|\cdot\big(\deltnonspanning+10\sigmKJ\big)^{|S_1|}\cdot \big(1-\deltnonspanning-10\sigmKJ\big)^{|S_2|}\,,\]
 since each $G'_s$ with $s\in\cJ_0$ has $(1-\deltnonspanning-10\sigmKJ)n$ vertices. We obtain similar estimates for $\cJ_1$ and $\cJ_2$; then by Corollary~\ref{cor:freedm} and the union bound over $\{\cJ_0,\cJ_1,\cJ_2\}$ and choices of $S_1$ and $S_2$,~\ref{A:typicality} follows.
\end{proof}

To prove~\ref{A:PathsWithOneVertex}, given $v$, we need to estimate the probability that the embedding of a single $G'_s$ with $s\in\cJ$ puts a $\SpecPaths$ end-vertex on $v$. Given such an estimate, Corollary~\ref{cor:freedm} then tells us that with high probability we have all the three desired equalities for each $v$. To estimate the probability of embedding some path end-vertex of $G_s$ to $v$, it is enough (since the events are disjoint) to estimate the probability of a specific path end-vertex $x\in V(G_s)$ being embedded to $v$, which is provided by~\eqref{eq:A:probvtx}. The proof of~\ref{A:PathsOneVtxIm} is essentially the same, except that~\eqref{eq:A:probvtx} is replaced with~\eqref{eq:A:probvtxnoim}, and we prove both statements together.

\begin{proof}[Proof of~\ref{A:PathsWithOneVertex} and~\ref{A:PathsOneVtxIm}]
 Fix $u,v\in V(H)$, and let $\cE$ be the event that $\PackingProcess'$ succeeds and for each $s$ the graph $H_{s-1}$ is $(\alpha_{s-1},2\newD+3)$-quasirandom.
 
 We first prove~\ref{A:PathsWithOneVertex}. 
 For a given $s\in\cJ$, let $X_s:=\{\leftpath_0(P),\rightpath_0(P):P\in\SpecPaths_s\}$. Note that $|X_s|=2|\SpecPaths_s|$. Conditioning on $H_{s-1}$, and assuming we are in $\cE$, for any $x\in X_s$ we have by~\eqref{eq:A:probvtx}
  \[\Prob\big[x\AlgMap v\big]=\big(1\pm10^4C\alpha_{s-1} \newD\newdelta^{-1}\big)n^{-1}\,.\]
  Since these events are disjoint, the conditional probability that some vertex of $X_s$ is mapped to $v$ is
  \[\big(1\pm 10^4C\alpha_{s-1}\newD\newdelta^{-1}\big)\cdot |X_s|n^{-1}=\big(1\pm\tfrac18\gamcore\big)|X_s|n^{-1}\,.\]
  
  For each $\cJ_i$ the sum of these conditional probabilities is easy to compute. For $\cJ_1$, we have $|\cJ_1|$ summands, in each of which $|X_s|=2\sigmJjedna n$. For $\cJ_0$ and $\cJ_2$, we have $|X_s|=2\sigmKJ n$. Applying Corollary~\ref{cor:freedm}, we see that with high probability either $\cE$ does not occur, or we have~\ref{A:PathsWithOneVertex}. The former is unlikely, so with high probability we have~\ref{A:PathsWithOneVertex}.
  
  The proof of~\ref{A:PathsOneVtxIm} is almost identical, so we only mention the difference. For~\ref{enu:PathsOneVtxIm0}--\ref{enu:PathsOneVtxIm2}, we use~\eqref{eq:A:probvtxnoim} to obtain the probability that the embedding of $G'_s$ maps a given $x\in X_s$ to $v$ and also $u\notin\im\phi_s^\mathbf{A}$, and we note that $\big(n-v(G'_s)\big)n^{-1}=\deltnonspanning+10\sigmJjedna$ for $s\in\cJ_1$, $\big(n-v(G'_s)\big)n^{-1}=\deltnonspanning+10\sigmKJ$ for $s\in\cJ_0$, and $\big(n-v(G'_s)\big)n^{-1}=\deltnonspanning+6\sigmKJ$ for $s\in\cJ_2$. For~\ref{enu:PathsTwoVtxIm0}, we use~\eqref{eq:A:probvtxnoim2} to obtain the probability that the embedding of $G'_s$ maps a given $x\in X_s$ to $v$ and also $u,u'\notin\im\phi_s^\mathbf{A}$.
\end{proof}

Finally we prove~\ref{A:PathsNotSquashed}. Here we do not need an accurate estimate. Given $v_1,v_2$, we find an upper bound for the probability of the embedding of a given $G'_s$ contributing to the count for this pair, and apply Corollary~\ref{cor:freedm} to obtain the desired result.

\begin{proof}[Proof of~\ref{A:PathsNotSquashed}]
 Fix $v_1,v_2\in V(H)$. Let $\cP=\bigcup_{s\in\cJ}\SpecPaths_s$; we aim to count the number of paths in $\cP$ whose ends are embedded to $\{v_1,v_2\}$.
 
 Observe that the embedding of any given $G'_s$ for $s\in\cJ$ contributes at most one to this count, so it is enough to find a good upper bound on the probability that a given $G'_s$ with $s\in\cJ$ does contribute one. Consider the embedding of $G'_s$ into $H_{s-1}$ by \RandomEmbedding. Suppose that at some point a vertex $x$, which is the first end of a path $P\in\cP$, is embedded to one of $\{v_1,v_2\}$ and the other of $\{v_1,v_2\}$ is not yet used in the embedding. Then the embedding of $G'_s$ contributes only if the other end $y$ of $P$ is embedded to the other of $\{v_1,v_2\}$. If the $(C\alpha_{s-1},2\newD+3)$-diet condition holds at the time when $y$ is embedded, then $y$ is embedded to a set of size at least $\tfrac12\dStageA^{\newD}\deltnonspanning n$, so the probability of its being embedded to the other of $\{v_1,v_2\}$ is at most $\tfrac{2}{\dStageA^{\newD}\deltnonspanning n}$ in this case. If $H_{s-1}$ is $(\alpha_{s-1},2\newD+3)$-quasirandom, then by Lemma~\ref{lem:propRandEmb} the probability that the diet condition fails is at most $2n^{-9}$.
 
 So we conclude that if $H_{s-1}$ is $(\alpha_{s-1},2\newD+3)$-quasirandom, the probability that $G_s$ embeds the ends of a path of $\cP$ to $\{v_1,v_2\}$ is at most $\tfrac{3}{\dStageA^{\newD}\deltnonspanning n}$. Let $\cE$ be the event that $H_{s-1}$ is $(\alpha_{s-1},2\newD+3)$-quasirandom for each $s$; then Corollary~\ref{cor:freedm}\ref{cor:freedm:tails}, with $\tilde\mu=\tilde\nu=3\gamma^{-\newD}\deltnonspanning^{-1}$ and $\tilde\rho=n^{0.2}$, says that the probability that $\cE$ occurs and yet more than $2n^{0.2}$ paths of $\cP$ have their ends embedded to $\{v_1,v_2\}$, is at most $\exp\big(-n^{0.1}\big)$.
 
 Taking the union bound over choices of $v_1$ and $v_2$, a.a.s.\ either $\cE$ fails or we have~\ref{A:PathsNotSquashed}. Since we proved that $\cE$ a.a.s.\ does not fail, a.a.s.\ we have~\ref{A:PathsNotSquashed}.
\end{proof}

\section{Stage~B (Proof of Lemma~\ref{lem:StageB})}\label{sec:NEWStageB}
The main purpose of this stage is to embed $(\SpecLeaves_s)_{s\in\cK}$ so that the parities of the degrees of the vertices of the host graph are prepared for later stages of the embedding process. This key parity condition is~\ref{B:parity}.
We begin by finding a pairing of vertices in $\HStageA$ with some vertices of $\bigcup_{s\in\cK}\SpecLeaves_s$ such that embedding any collection of the paired $\bigcup_{s\in\cK}\SpecLeaves_s$ to their paired vertices of $\HStageA$ gives a valid packing. Recall that each graph $G_s^\spadesuit$ with $s\in\cK$ has by~\eqref{eq:familyOrderKEquality} $(1-\deltnonspanning)n-|\SpecLeaves_s|$ vertices. Since $|\SpecLeaves_s|\le\tfrac{cn}{\log n}$ by Definition~\ref{def:family}\ref{enu:defSpecLeaves}, this means $v(G^\spadesuit_s)=\big(1-\deltnonspanning+o(1)\big)n$ and the $o(1)$ will always be absorbed into the error terms in the following arguments.

\begin{lemma}
 There is an injective map $\pi:V(\HStageA)\to\bigcup_{s\in\cK}\SpecLeaves_s$ such that if $\pi(v)=x\in V(G_s)$, then $v\notin\imA(s)$, the pair $v\phi^{\mathbf{A}}_s(y)$ is an edge of $\HStageA$ for each $y\in \NBH_{G_s}(x)$, and no edge of $\HStageA$ is obtained twice in this way (i.e.\ no embedded neighbour of $\pi\big(\phi^{\mathbf{A}}_s(y)\big)$ is embedded to $v$).
\end{lemma}
\begin{proof}
 Let $\overrightarrow{H}$ be a uniform random orientation of $E(\HStageA)$. Given any $\ell\le 2D$ and vertices $v_1,\dots,v_\ell$ of $\HStageA$, and any subset $T$ of $\NBH_{\HStageA}(v_1,\dots,v_\ell)$ such that $|T|\ge\gamcore n$, the expected number of vertices of $T$ which are out-neighbours of each of $v_1,\dots,v_\ell$ in $\overrightarrow{H}$ is $2^{-\ell}|T|$, and by Chernoff's inequality the probability that the actual number is not $(1\pm\gamcore)2^{-\ell}|T|$ is at most $\exp\big(-n^{0.1}\big)$. In particular with high probability, by the union bound, for each of the polynomially many sets defined by~\ref{A:megaquasirandomness1} with $\ell$ vertices, a $(1\pm\gamcore)2^{-\ell}$-fraction are out-neighbours of each of the $\ell$ vertices in $\overrightarrow{H}$. We will need two specific cases of this: if $x\in\SpecLeaves_s$, then
 \begin{equation}\label{eq:B:outdeg}
  \big|\NBHout_{\overrightarrow{H}}\big(\phi^{\mathbf{A}}_s(\NBH_{G_s}(x))\big)\setminus\imA(s)\big|=(1\pm 3\gamcore)2^{-\DParity}\dStageA^{\DParity} \deltnonspanning n\,,
 \end{equation}
 and if in addition $y\in\SpecLeaves_{s'}$, with $s'\neq s$, is such that $\phi^{\mathbf{A}}_s(\NBH_{G_s}(x))$ and $\phi^{\mathbf{A}}_{s'}(\NBH_{G_{s'}}(y))$ are disjoint, we have
 \begin{equation}\label{eq:B:outcodeg}
 \big|\NBHout_{\overrightarrow{H}}\big(\phi^{\mathbf{A}}_s(\NBH_{G_s}(x))\cup \phi^{\mathbf{A}}_{s'}(\NBH_{G_{s'}}(y))\big)\setminus\big(\imA(s)\cup\imA(s') \big)\big|=(1\pm 3\gamcore)2^{-2\DParity}\dStageA^{2\DParity} \deltnonspanning^2 n\,.
 \end{equation}
 
 Furthermore, from~\ref{A:SpecLeavesCover}, with high probability for each $v\in V(\overrightarrow{H})$ there are at least \begin{equation}\label{eq:B:highdegF}
 2^{-\DParity-1}\cdot 4^{-200 D^{2}\dStageA^{-D}\deltnonspanning^{-1}}n
\end{equation} vertices of $\bigcup_{s\in\cK}\SpecLeaves_s$ which we could map to $v$ using only edges going to $v$ and without destroying the packing (i.e.\ $v$ is not in the image of that particular graph).
Let us justify this last statement, which needs a little care because the relevant events are not independent: by~\ref{A:parleavesSpread} a single edge at $v$ can be needed by as many as $\tfrac{20cn}{\log n}$ different vertices of $\bigcup_{s\in\cK}\SpecLeaves_s$.
Fix $v$, let $\cA$ be the set of those $x\in\bigcup_{s\in\cK}\SpecLeaves_s$ which are $\SpecLeaves$ covering $v$, and set $M:=|\cA|$, so that $M\ge 4^{-200 D^{2}\dStageA^{-D}\deltnonspanning^{-1}}n$ by~\ref{A:SpecLeavesCover}.
For $x\in\cA\cap V(G_s)$ let $E_x$ be the set of the $\DParity$ edges joining $v$ to $\phi^{\mathbf{A}}_s\big(\NBH_{G_s}(x)\big)$, so that what we have to bound from below is the number of $x\in\cA$ all of whose edges of $E_x$ are oriented towards $v$.
Fix an arbitrary order $e_1,\dots,e_d$ of the edges of $\HStageA$ at $v$, let $\cF_i$ be generated by the orientations of $e_1,\dots,e_i$, and for $0\le j\le\DParity$ let $X^{(j)}$ be the number of $x\in\cA$ whose first $j$ edges in $E_x$ (with respect to this order) are all oriented towards $v$.
Thus $X^{(0)}=M$, while $X^{(\DParity)}$ is the quantity we want to bound.
Let $Y^{(j)}_i$ be the number of $x\in\cA$ whose $j$-th edge in $E_x$ is $e_i$ and whose first $j$ edges in $E_x$ are all oriented towards $v$.
Then $X^{(j)}=\sum_{i=1}^dY^{(j)}_i$, each $Y^{(j)}_i$ is $\cF_i$-measurable with $0\le Y^{(j)}_i\le\tfrac{20cn}{\log n}$ by~\ref{A:parleavesSpread}, and, since the first $j-1$ edges in $E_x$ precede the $j$-th one,
\[\sum_{i=1}^d\Exp\big[Y^{(j)}_i\big|\cF_{i-1}\big]=\tfrac12X^{(j-1)}\,.\]
Put $\eta_j:=2^{j}4^{-D}$. Applying Corollary~\ref{cor:freedm}\ref{cor:freedm:tails} with $\tilde\mu=2^{-j}M$, $\tilde\nu=\tilde\rho=\eta_{j-1}2^{-j}M$, $R=\tfrac{20cn}{\log n}$ and $\cE$ the event that $X^{(j-1)}=(1\pm\eta_{j-1})2^{-j+1}M$, we obtain inductively that $X^{(j)}=(1\pm\eta_j)2^{-j}M$ for each $j\le\DParity$, with probability at least $1-2\DParity\exp\big(-\Omega(c^{-1}\log n)\big)$, where the implicit constant depends only on $D$, $\dStageA$ and $\deltnonspanning$.
Since $c$ is small compared to these, and since $\eta_{\DParity}\le 2^{-D}\le\tfrac12$, the union bound over the choices of $v$ now gives~\eqref{eq:B:highdegF}.
 
 We draw an auxiliary bipartite graph $F$ whose parts are $V(\overrightarrow{H})$ and $\bigcup_{s\in\cK}\SpecLeaves_s$. We put an edge from $x\in\SpecLeaves_s$ to $v\in V(\overrightarrow{H})$ if and only if $v\notin\imA(s)$ and $v$ is an out-neighbour of $\phi_s^{\mathbf{A}}(y)$ for each $y\in \NBH_{G_s}(x)$. 
 \begin{claim}\label{cl:HallStB}
 There exists a matching in $F$ covering $V(\overrightarrow{H})$.
 \end{claim}
Such a matching gives the desired $\pi$; the orientation of $\overrightarrow{H}$ ensures that no edge is used twice.
\begin{claimproof}[Proof of Claim~\ref{cl:HallStB}]
We verify Hall's condition for $F$. To this end, consider an arbitrary non-empty set $T\subset V(\overrightarrow{H})$. Let $S\subset\bigcup_{s\in\cK}\SpecLeaves_s$ be the vertices of $F$ adjacent to at least one member of $T$. 

Firstly, applying~\eqref{eq:B:highdegF} to an arbitrary vertex of $T$, we see that $|S| \ge 2^{-\DParity-1}\cdot 4^{-200 D^{2}\dStageA^{-D}\deltnonspanning^{-1}}n$. In particular, the only cases for which it remains to prove Hall's condition are when 
\begin{equation}\label{eq:TsThatRemain}
|T|> 2^{-\DParity-1}\cdot 4^{-200 D^{2}\dStageA^{-D}\deltnonspanning^{-1}}n \;.
\end{equation}
We claim that in such a case, we actually have 
 \begin{equation}\label{eq:failfailfail}
|S|\ge \big|\bigcup_{s\in\cK}\SpecLeaves_s\big|-20\gamcore n\;.
 \end{equation} This obviously verifies Hall's condition since
 $$\big|\bigcup_{s\in\cK}\SpecLeaves_s\big|-20\gamcore n>n\ge |T|\;.$$
 
 So, suppose for a contradiction that~\eqref{eq:failfailfail} fails, and for the set $S^*:=\bigcup_{s\in\cK}\SpecLeaves_s\setminus S$ we have 
 \begin{equation}\label{eq:mylittlelowerbound}
 |S^*|\ge 20\gamcore n\;.
  \end{equation}
  Let $T^*:=V(\overrightarrow{H})\setminus T$. Let us count the number $N$ of triples $(u,v,u')$ with $u,u'\in S^*$ and $v\in T^*$ such that $uv,u'v\in E(F)$.
 
By~\eqref{eq:B:outdeg}, we have that the average degree of a vertex in $T^*$ to $S^*$ is
 \[ \frac{e_F(T^*,S^*)}{|T^*|}=\frac{\sum_{x\in S^*} \deg_F(x)}{|T^*|}
 \ge 
 \frac{(1-3\gamcore)\dStageA^{\DParity}2^{-\DParity}\deltnonspanning n |S^*|}{|T^*|}\;.\]
 Hence by Jensen's inequality the average number of ordered pairs of $F$-neighbours a vertex $t\in T^*$ has in $S^*$ is at least
 \[\frac{|S^*|^2}{|T^*|^2}(1-3\gamcore)^2\dStageA^{2\DParity}2^{-2\DParity}\deltnonspanning^2 n^2\,.\]
 It follows that
 \begin{equation}\label{eq:tgb}
 N\ge \frac{|S^*|^2}{|T^*|}(1-3\gamcore)^2\dStageA^{2\DParity}2^{-2\DParity}\deltnonspanning^2 n^2\,.
  \end{equation}
 
On the other hand, we bound $N$ from above by counting starting from pairs $(u,u')$ of (not necessarily distinct) vertices of $S^*$. 
Let $\indexofvertex(u)$ denote the index such that $u\in V(G_{\indexofvertex(u)})$.
We say that $(u,u')$ is of \emph{Type~1} if $\indexofvertex(u)=\indexofvertex(u')$. 
We say that $(u,u')$ is of \emph{Type~2} if $\indexofvertex(u)\neq \indexofvertex(u')$ and the embeddings of the neighbourhoods of $u$ and $u'$ are not disjoint. The remaining pairs $(u,u')$ are \emph{Type~3}. 
To bound the number $M_1$ of pairs of Type~1, we have
\begin{equation}\label{eq:IHM1}
M_1= \sum_{x\in S^*}|S^*\cap \SpecLeaves_{\indexofvertex(x)}|\leBy{D\ref{def:family}\ref{enu:defSpecLeaves}} |S^*| \cdot \tfrac{cn}{\log n}\;.
\end{equation}
 
To bound the number $M_2$ of pairs of Type~2, we consider an arbitrary vertex $x\in S^*$. By~\ref{A:parleavesSpread}, we know that for each $v\in \phi^{\mathbf{A}}_{\indexofvertex(x)}(\NBH_{G_{\indexofvertex(x)}}(x))$ there are at most $\frac{20cn}{\log n}$ vertices $x'\in S^*$ whose embedded neighbourhood touches $v$. That means that there are at most $\tfrac{20\DParity cn}{\log n}$ vertices $x'\in S^*$ whose embedded neighbourhood is not disjoint from the embedded neighbourhood of $x$. Hence,
\begin{equation}\label{eq:IHM2}
M_2\le \frac{20\DParity cn}{\log n} \cdot |S^*|\;.
\end{equation}
 For the number $M_3$ of pairs of Type~3, we use the trivial bound
\begin{equation}\label{eq:IHM3}
M_3\le |S^*|^2\;.
\end{equation} 
We now need to bound the number of ways a given pair $(u,u')$ can be extended to a triple $(u,v,u')$ as above. We shall use the trivial bound that the number of such extensions is at most $|T^*|$ for pairs of Type~1 and~2. For pairs of Type~3, we use~\eqref{eq:B:outcodeg} and see that the number of the extensions is at most
 \[(1+3\gamcore)\dStageA^{2\DParity}2^{-2\DParity}\deltnonspanning^2 n\,.\]
Using~\eqref{eq:IHM1}, \eqref{eq:IHM2} and \eqref{eq:IHM3}, we obtain
\begin{align}
\begin{split}\label{eq:bgt}
N&\le  M_1 \cdot |T^*|+M_2 \cdot |T^*|+ M_3\cdot (1+3\gamcore)\dStageA^{2\DParity}2^{-2\DParity}\deltnonspanning^2 n\\
&\le |S^*|^2(1+4\gamcore)\dStageA^{2\DParity}2^{-2\DParity}\deltnonspanning^2 n\,,
\end{split}
\end{align}
where the argument that the contribution from Types~1 and~2 is negligible compared to Type~3 uses~\eqref{eq:mylittlelowerbound}.

Comparing~\eqref{eq:tgb} and~\eqref{eq:bgt}, we get
 \[n(1-3\gamcore)^2\le(1+4\gamcore)|T^*|\,,\]
which is a contradiction to~\eqref{eq:TsThatRemain}.

 This verifies Hall's condition.
\end{claimproof}
The existence of the map $\pi$ follows from the discussion above.
\end{proof}

Now let
\[S:=\big\{v\in V(\HStageA):\deg_{\HStageA}(v)\not\equiv \PathTerm(v)+\OddOut(v)\mod 2\big\}\,.\]
Observe that $X:=\bigcup_{s\in\cK}\SpecLeaves_s\setminus \pi(S)$ is a set of size at least $\gamcore n$. Let us now show that $|X|$ is even. Indeed, writing $\equiv$ for equality modulo~2, we have
\begin{align*}
|X|
&\equiv\sum_{s\in \cK}|\SpecLeaves_s|+\sum_{v\in V(H)}\left(\deg_{\HStageA}(v)-\OddOut(v)-\PathTerm(v)\right)\\
&\equiv \sum_{v\in V(H)}\left(\deg_{\HStageA}(v)-\PathTerm(v)\right)\\
&\equiv 2e(\HStageA)-2\left|\bigcup_{s\in \cJ}\SpecPaths_s\right|\equiv 0\;.
\end{align*}

We draw a graph $F'$ on vertex set $X$, putting an edge between $u,v\in X$ if $\indexofvertex(u)\neq\indexofvertex(v)$ and their embedded neighbours are disjoint. By~\ref{A:parleavesSpread} and Definition~\ref{def:family}\ref{enu:defSpecLeaves}, every vertex of $X$ is adjacent to all but at most $(\DParity+1)c n$ others. In particular, $F'$ has minimum degree at least $|X|/2$, so by Dirac's theorem $F'$ contains a perfect matching, say $M$.

We now construct the $\phi_s^\mathbf{B}$ for $s\in\cK$. First, for each $v\in S$, suppose $\pi(v)\in V(G_s)$ and extend $\phi_s^\mathbf{A}$ by mapping $\pi(v)$ to $v$. Let $H'$ be the graph of edges remaining from $\HStageA$ (removing the edges used by this embedding). Having embedded exactly one vertex to each $v\in S$, we greedily embed pairs of vertices of $X$ as given by $M$. Suppose $xy\in M$, with $x\in V(G_s)$ and $y\in V(G_{s'})$. Recall that $s\neq s'$. We choose $v\in V(H)$ which is not in $\imA(s)\cup\imA(s')$, and which was not used to embed any vertex of $\SpecLeaves_s\cup\SpecLeaves_{s'}$, and which is a common neighbour in $H'$ of the vertices $\phi_s^\mathbf{A}(\NBH_{G_s}(x))\cup\phi_{s'}^\mathbf{A}(\NBH_{G_{s'}}(y))$, and which was not chosen more than $c n$ times in this process. We embed both $x$ and $y$ to $v$, and remove from $H'$ the edges used by these embeddings. We claim this process succeeds, and defines the required embeddings $\phi_s^\mathbf{B}$ for $s\in\cK$ and remaining graph $\HStageB$. Note that if the process succeeds, then we automatically have~\ref{B:one},~\ref{B:two}.

Observe that in this construction, there are two ways we use edges at a given vertex $v$. We use an edge at $v$ for each vertex of $\bigcup_s\SpecLeaves_s$ for which $v$ is an embedded neighbour, and we use edges at $v$ whenever we embed a vertex of $\bigcup_s\SpecLeaves_s$ to $v$. By~\ref{A:parleavesSpread}, we use at most $\tfrac{20cn}{\log n}$ edges of the first type. Since each choice of $v$ to embed one or two $\SpecLeaves$ uses at most $2\DParity$ edges at $v$, and since we choose $v$ at most $cn$ times, there are at most $6Dcn$ edges of the second type. So in total the degree of $v$ decreases by at most $8Dcn$ from $\HStageA$ to $\HStageB$, hence obtaining~\ref{B:degrees}. Now, the only way in which our construction can fail is that at some point we need to embed some matched pair $xy\in M$ and all the vertices to which we can embed them have already been used $cn$ or more times. By~\ref{A:megaquasirandomness1}, the set of vertices to which we could embed both $x$ and $y$ before the start of this stage has size at least $(1-\gamcore)\dStageA^{2D}\deltnonspanning^2 n$. Of these vertices, at most $\tfrac{2cn}{\log n}$ may have been used for embedding $\SpecLeaves_s$ and $\SpecLeaves_{s'}$ (where $x\in V(G_s)$ and $y\in V(G_{s'})$) and hence have become unavailable. A further at most $16D^2cn$ may have become unavailable since an edge from one of the at most $2D$ embedded neighbours of $x$ or $y$ has been used. However this leaves at least $cn$ vertices; they cannot all have been used $cn$ times since the size of $\bigcup_s\SpecLeaves_s$ is at most $2n$. We conclude that the construction process succeeds.

Let us now turn to the key property~\ref{B:parity}. Given a vertex $v\in V(\HStageA)$, each edge $vw$ used in Stage~B was of one of the following two types:
\begin{enumerate}[label=\abc]
\item for some $x\in\SpecLeaves_s$, and some $y\in \NBH_{G_s}(x)$, $y$ was embedded to $v$ (in Stage~A) and $x$ was embedded to $w$ (in Stage~B),
\item for some $x\in\SpecLeaves_s$, and some $y\in \NBH_{G_s}(x)$, $y$ was embedded to $w$ (in Stage~A) and $x$ was embedded to $v$ (in Stage~B). 
\end{enumerate}
The number of edges as in~(a) is $\OddOut(v)$. As for edges as in~(b) in the case $v\in S$, we recall that when we placed $\pi(v)$ we used exactly $\DParity$ many such edges, which is odd (cf.~Definition~\ref{def:family}). Additional edges as in~(b) come exclusively by placing pairs $\{x,y\}\subset \bigcup_s\SpecLeaves_s$ (here, $xy\in M$) on $v$; placement of each such pair uses exactly $2\DParity$ edges, which is even. In the case $v\notin S$, the edges incident with~$v$ considered in~(b) come only from placing pairs, and hence their total number is even. Put together, we see that if $v\in S$ the total number of edges incident to~$v$ used in Stage~B is equal modulo~2 to $\OddOut(v)+1$, while if $v\notin S$, it is equal modulo~2 to $\OddOut(v)$. This establishes~\ref{B:parity}.

Observe that~\ref{B:megaquasirandomness1} follows from~\ref{A:megaquasirandomness1} and~\ref{B:degrees}: the ranges of $S$ and $T$ in the two statements agree, we have $\imB(s)=\imA(s)$ for each $s\in\cJ$, and $\dStageB=\dStageA\pm O(n^{-1})$ since Stage~B embeds only $O(n)$ edges. Finally, by~\ref{B:degrees} the set whose size we estimate shrinks by at most $|S|\cdot 8Dcn\le16D^2cn$ between $\HStageA$ and $\HStageB$, which is negligible compared to $\gamcore\dStageA^{2D}\deltnonspanning^{D} n$ because $c\ll\gamcore$.

\section{Stage~C (Proof of Lemma~\ref{lem:StageCNew})}\label{sec:StageC}

If $n$ is odd, we pick an arbitrary vertex of $H$ to be $\boxdot$. We then need to embed all the paths in $\bigcup_{s\in\cJ}\SpecPaths_s$ which are anchored at $\boxdot$ (of which there can be at most one per $s\in\cJ$) and in addition some more paths in order that all the edges at $\boxdot$ are used, and exactly one path from each $\SpecPaths_s$ for $s\in\cJ_0$ is embedded.

Let us first provide an overview of the proof. Set\index{$\cJ_{\boxdot}$} \[\cJ_{\boxdot}:=
 \left\{s\in \cJ\::\: \exists P\in \SpecPaths_s, 
 \{\phi_s^{\mathbf{B}}(\leftpath_0(P)), \phi_s^\mathbf{B}(\rightpath_0(P))\}\cap \{\boxdot\}\neq \emptyset\right\}\;.\] Observe that for each $s\in \cJ_{\boxdot}$, there is exactly one such path, which we denote by \index{$P^\boxdot_s$}$P^\boxdot_s$. We need to use $|\cJ_\boxdot|$ edges leaving $\boxdot$ to embed these paths, after which we are left with $\deg_{\HStageB}(\boxdot)-|\cJ_\boxdot|$ edges remaining to cover. By~\ref{B:parity}, $q:=\frac{\deg_{\HStageB}(\boxdot)-|\cJ_\boxdot|}{2}$ is an integer. By~\ref{A:megaquasirandomness1}, with $S=\{\boxdot\}$ and $T=\emptyset$, we have $\deg_{\HStageA}(\boxdot)=(1\pm\gamcore)\dStageA n$ and by~\ref{B:degrees} we therefore have
 \begin{equation*}
 \deg_{\HStageB}(\boxdot)=(1\pm2\gamcore)\dStageA n\eqByRef{eq:dStageAapprox}(14\pm 0.2)\sigmKJ^2 n\;.
 \end{equation*}
By~\ref{A:PathsWithOneVertex} we then have 
\begin{align}
\label{eq:sizeJBoxdot}
|\cJ_{\boxdot}|&=(2\pm \gamcore)\big(\sigmKJ|\cJ_0|+\sigmJjedna|\cJ_1|+\sigmKJ|\cJ_2|\big)\eqByRef{eq:CONSTANTS}(2\pm 0.1)\sigmKJ^2 n\;\mbox{, and}\\
\label{eq:q}
q&=(6\pm 0.3)\sigmKJ^2 n\;.
\end{align}
 
\begin{figure}
	\includegraphics[width=\textwidth]{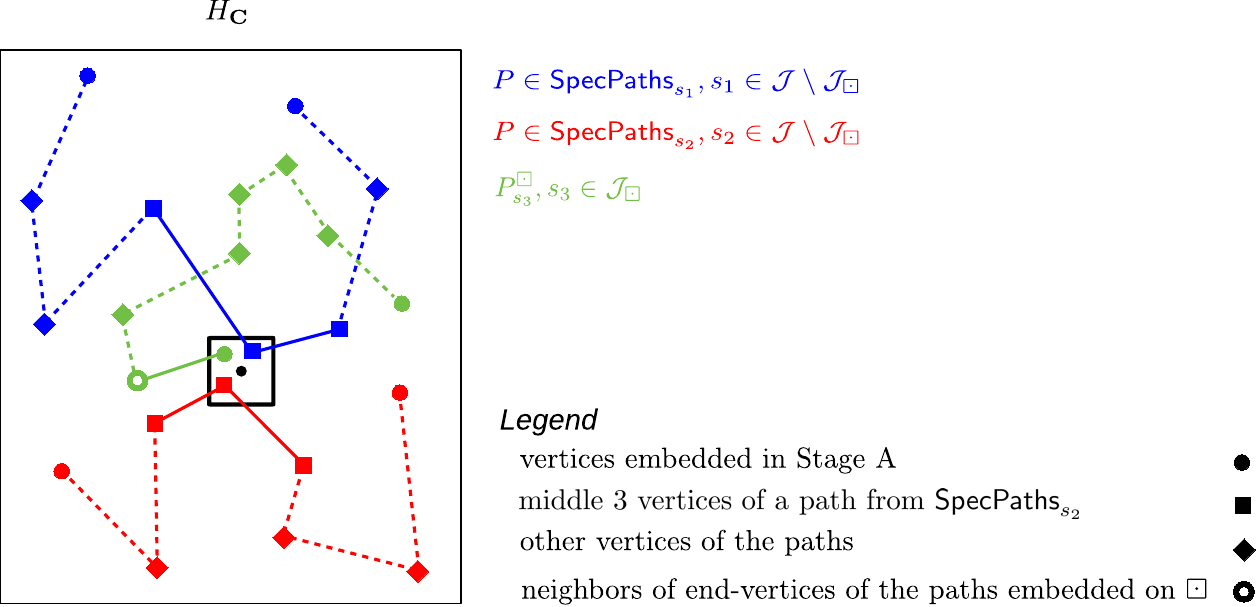}
	\caption{Embedding in Stage~C. Vertex $\boxdot$ enlarged for clarity.}
	\label{fig:StageCPaths}
\end{figure}

The only part of this embedding where we need to be careful is when we embed the edges at $\boxdot$. To do this, we use a matching argument. Recall that since in Stage~B we embed only edges of graphs $G_s$ for $s\in\cK$, we have for each $s\in\cJ$ that $\imA(s)=\imB(s)$. We define an auxiliary balanced bipartite graph $B$. One part of $B$ is $X:=\{P^\boxdot_s:s\in\cJ_\boxdot\}\cup\{x_i,y_i:i\in[q]\}$, where $x_i$ and $y_i$ are auxiliary elements. The other part of $B$ is $\NBH_{\HStageB}(\boxdot)$. We draw an edge in this graph from each $P^\boxdot_s$ to each $u\in \NBH_{\HStageB}(\boxdot)$ such that $u\notin\imB(s)$. To define the remaining edges, we work as follows. We define
 \[Z:=\big\{s\in\cJ\setminus\cJ_\boxdot\,:\,\boxdot\notin\imB(s)\big\}=\big\{s\in\cJ\,:\,\boxdot\notin\imB(s)\big\}\;.\]
From~\ref{A:typicality} with $S_1=\{\boxdot\}$ and $S_2=\emptyset$, we see that 
\begin{equation}\label{eq:sizeZ}
0.5 \deltnonspanning|\cJ|\le |Z|\le 2\deltnonspanning|\cJ|\;.
\end{equation}
In particular, $|Z|\ge q$. Hence, we can choose uniformly at random distinct indices $s_1,\dots,s_q$ from $Z$, and for each $i\in[q]$ we draw an edge from each of $x_i$ and $y_i$ to each $u\in \NBH_{\HStageB}(\boxdot)$ such that $u\notin\imB(s_i)$. The point of this graph is the following: a perfect matching in $B$ tells us how to embed the first edge of each $P^\boxdot_s$ (and we will be able to greedily finish off the embedding) and for each $i\in[q]$, two edges which will be in the middle of a path of $\SpecPaths_{s_i}$ (we will be able to choose which path, and complete the embedding greedily). See Figure~\ref{fig:StageCPaths}.
 
 \begin{lemma}\label{lem:C:match} With positive probability, there is a perfect matching $M$ in $B$.
 \end{lemma}
 \begin{proof}
Since $B$ is balanced, the desired property is equivalent to the existence of a matching covering $X$. To that end, we shall use Hall's condition. Let us note some degree and codegree properties of $B$. For any vertex $a$ in $X$, by~\ref{A:megaquasirandomness1}, with $S=\{\boxdot\}$ and $T$ being either $s$ where $a=P^\boxdot_s$ or $s_i$ where $a\in\{x_i,y_i\}$, and by~\ref{B:degrees}, we see that
 \[\deg_B(a)=(1\pm\gamcore)\dStageA n\big(1-\tfrac{|\imA(s)|}{n}\big)\pm 8Dcn=(1\pm 2\gamcore)\dStageA n\big(1-\tfrac{|\imA(s)|}{n}\big)\,,\]
 where the first term comes from~\ref{A:megaquasirandomness1} and the second term uses~\ref{B:degrees} to estimate the number of these edges removed in Stage~B.
 
 Similarly, given distinct $a,a'\in X$ which are not of the form $x_i,y_i$ for any $i\in[q]$, again using~\ref{A:megaquasirandomness1} with $S=\{\boxdot\}$ and this time with $T=\{s,s'\}$  where $s$ is the index corresponding as above to $a$, and $s'$ (which is by assumption not $s$) is the index corresponding to $a'$, we see that
 \begin{equation}\label{eq:C:pairdeg}\deg_B(\{a,a'\})=(1\pm 2\gamcore)\dStageA n\big(1-\tfrac{|\imA(s)|}{n}\big)\big(1-\tfrac{|\imA(s')|}{n}\big)\,.\end{equation}
 
 Finally, given $v\in \NBH_{\HStageB}(\boxdot)$, using~\ref{A:typicality} with $S_1=\{\boxdot,v\}$ and $S_2=\emptyset$, we see that for the number $f$ of indices $s\in\cJ$ such that neither $\boxdot$ nor $v$ is in $\imB(s)$ we have $f\ge (1-\gamcore)\deltnonspanning^2|\cJ|$. Hence the number of indices $s_i$ with $v\notin \imB(s_i)$ has the hypergeometric distribution with parameters $|Z|$, $f$, $q$. Hence, the expectation of this random variable is 
$$
\frac{fq}{|Z|}\geBy{\eqref{eq:q},\eqref{eq:sizeZ}}\frac{(1-\gamcore)\deltnonspanning^2|\cJ|\cdot 5\sigmKJ^2 n}{2\deltnonspanning|\cJ|}
\geByRef{eq:CONSTANTS} 40\gamcore n\;.
$$
Fact~\ref{fact:hypergeometricBasicProperties} tells us that the probability that we choose fewer than $30\gamcore n$ of these indices among $s_1,\ldots,s_q$ is less than $n^{-2}$. In particular, taking the union bound over the at most $n$ choices of $v$, we see $\deg_B(v)\ge 30\gamcore n$ holds for all $v\in \NBH_{\HStageB}(\boxdot)$ with probability at least $1-n^{-1}$. Suppose that this likely event occurs.
 
Now given any non-empty $X'\subset X$, we verify Hall's condition. We shall distinguish three cases based on the size of $|X'|$.

If $|X'|<\tfrac12\dStageA\deltnonspanning n$, then choosing any $x\in X'$ we have that $\deg_B(x)>|X'|$ by the degree calculation at the beginning of this proof. This verifies Hall's condition for $X'$.

Suppose next that $|X'|>(1-20\gamcore)\dStageA n$. Then using that 
$|X|=\deg_{\HStageB}(\boxdot)=(1\pm2\gamcore)\dStageA n$ and
$\deg_B(v)\ge 30\gamcore n$ holds for all $v\in \NBH_{\HStageB}(\boxdot)$, it follows that the joint neighbourhood of $X'$ in $B$ is the entire part $\NBH_{\HStageB}(\boxdot)$. Hence, Hall's condition for $X'$ holds again.

Suppose now $X'$ is any set with $\tfrac12\dStageA\deltnonspanning n\le|X'|\le (1-20\gamcore)\dStageA n$. Let the joint neighbourhood of $X'$ in $B$ be $Y$. We count the total number of triples $(a,a',y)$ with $a,a'\in X'$ and $y\in Y$ such that $ay$ and $a'y$ are edges of $B$. On the one hand, this is equal to
 \[\sum_{a,a'\in X'}\deg_B(\{a,a'\})\,.\]
We can use the upper bound~\eqref{eq:C:pairdeg} for all but at most $3|X'|$ pairs $a,a'$.
 On the other, it is also equal to
 \[\sum_{y\in Y}\big|\NBH_B(y)\cap X'\big|^2\ge\tfrac{1}{|Y|}\Big(\sum_{y\in Y}\big|\NBH_B(y)\cap X'\big|\Big)^2=\tfrac{1}{|Y|}\Big(\sum_{a\in X'}\deg_B(a)\Big)^2=\tfrac{1}{|Y|}\sum_{a,a'\in X'}\deg_B(a)\deg_B(a')\]
 by Jensen's inequality, where the middle equality holds since $\NBH_B(a)\subseteq Y$ for each $a\in X'$. Putting these two together, and for convenience writing $s$ and $s'$ respectively for the indices corresponding to $a$ and $a'$ in summations, we have
 \begin{multline*}
  \sum_{a,a'\in X'}(1+ 2\gamcore)\dStageA n\big(1-\tfrac{|\imA(s)|}{n}\big)\big(1-\tfrac{|\imA(s')|}{n}\big)+3|X'|n\\
  \ge \tfrac{1}{|Y|}\sum_{a,a'\in X'}(1-2\gamcore)^2\dStageA^2n^2\big(1-\tfrac{|\imA(s)|}{n}\big)\big(1-\tfrac{|\imA(s')|}{n}\big)\,.
 \end{multline*}
 Now the $3|X'|n$ term on the first line is, because $|X'|$ is sufficiently large, small compared to the main term; in particular we have
 \begin{multline*}
  \sum_{a,a'\in X'}(1+ 3\gamcore)\dStageA n\big(1-\tfrac{|\imA(s)|}{n}\big)\big(1-\tfrac{|\imA(s')|}{n}\big)\\
  \ge \tfrac{1}{|Y|}\sum_{a,a'\in X'}(1-2\gamcore)^2\dStageA^2n^2\big(1-\tfrac{|\imA(s)|}{n}\big)\big(1-\tfrac{|\imA(s')|}{n}\big)\,.
 \end{multline*}
 Now most of the terms cancel (in particular the sum over $\big(1-\tfrac{|\imA(s)|}{n}\big)\big(1-\tfrac{|\imA(s')|}{n}\big)$ is identical on both sides) and we get $(1+3\gamcore)\ge\tfrac{1}{|Y|}(1-2\gamcore)^2\dStageA n$. In particular, we have $|Y|>|X'|$ as required. This completes the verification of Hall's condition.
 \end{proof}
 
So, suppose that our choice of indices $s_1,\dots,s_q$ yields a perfect matching in Lemma~\ref{lem:C:match}. We now explain how to use this perfect matching $M$ to embed complete paths. The condition we need to maintain is that we do not use any vertex other than $\boxdot$ more than $n^{0.5}$ times in these embeddings. We assign to each vertex $v$ in $V(H)\setminus\{\boxdot\}$ a number of uses, which initially is equal to the number of paths in $\{P^\boxdot_s:s\in\cJ_\boxdot\}$ which are anchored at $v$ (i.e.~their endvertex not mapped by $\phi^\mathbf{A}_s$ to $\boxdot$ is mapped to $v$), plus one if $v\in \NBH_{\HStageB}(\boxdot)$. By~\ref{A:PathsNotSquashed}, the initial number of uses of any given $v$ is at most $n^{0.3}+1$. Whenever we use a vertex in the following greedy embeddings, we add one to its number of uses. When a vertex has number of uses $n^{0.5}$, we mark it as forbidden; we claim that at no stage do we have more than $20n^{0.5}$ forbidden vertices. Indeed, initially summing the number of uses over all vertices $v$ we have at most $n$ (since at most $n/2$ paths are anchored at $\boxdot$ by~\ref{A:PathsWithOneVertex}, and since the number of `plus ones' is $\deg_{\HStageB}(\boxdot)\le n/2$). The total number of paths we embed greedily is at most $n$ (since $|\cJ|\le n$ and we embed at most one path of any given graph), and all these paths have at most $12$ vertices, hence the total number of uses cannot exceed $20n$.
 
 To begin with, we choose for each $i\in[q]$ a path $Q_{s_i}\in\SpecPaths_{s_i}$\index{$Q_s$}, and for each $s\in\cJ_0\setminus\cJ_\boxdot$ such that $s\neq s_i$ for $i\in[q]$ a path $Q_s\in\SpecPaths_s$. We choose these greedily, in each step choosing a path in $G_s$ whose endvertices are both mapped by $\phi^\mathbf{A}_s$ to vertices which are not forbidden; we then add one to the number of uses of each endvertex and update the set of forbidden vertices. Since the paths of $\SpecPaths_s$ are vertex-disjoint and $\phi^\mathbf{A}_s$ is injective, each forbidden vertex is the anchor image of at most one of them; as $\SpecPaths_s$ contains at least $\sigmJjedna n$ paths and at most $20n^{0.5}$ vertices are forbidden, this is always possible.
 
  We now choose greedily embeddings of the paths $P^\boxdot_s$ for $s\in\cJ_\boxdot$. For $s\in\cJ_{\boxdot}$, we start with the embedding $\phi^{\mathbf{A}}_s$, and extend it to $\phi^\mathbf{C-}_s$ by mapping the vertex $u_2$ of $P^\boxdot_s$ adjacent to $\big(\phi^\mathbf{A}_s\big)^{-1}(\boxdot)$ to the vertex $v_2$ paired in the matching $M$ with $P^\boxdot_s$. Let $v'=\phi^\mathbf{A}_s(u')$ where $u'$ is the other endvertex of $P^\boxdot_s$. We then for each $i\ge 3$ in succession map the $i$th unembedded vertex $u_i$ of $P^\boxdot_s$ to a vertex $v_i$ in $\NBH_{\HStageB}(v_{i-1},v')\setminus\imA(s)$. Note that the condition of $v_i$ being a neighbour of $v'$ is actually only required for the last unembedded vertex $u_i$; we however keep this stronger requirement to avoid case distinction in the following.
Moreover, we choose $v_i$ which is not forbidden, which is not one of the at most $10$ previous vertices used for $G_s$ in Stage~C, and such that the edges $v_{i-1}v_i$ and $v_{i-1}v'$ have not been used for any previous embedding in Stage~C. By construction, the result is an embedding of $P^\boxdot_s$ if each step is possible. We claim that each step is indeed possible: by~\ref{A:megaquasirandomness1}, with $S=\{v_{i-1},v'\}$ and $T=\{s\}$, there are at least $(1-\gamcore)\dStageA^2\deltnonspanning n$ vertices which are adjacent in $\HStageA$ to both $v_{i-1}$ and $v'$ and which are not in $\imA(s)$. At most $16Dcn$ of these vertices are no longer in the common neighbourhood in $\HStageB$ by~\ref{B:degrees}; at most $20n^{0.5}$ are forbidden; at most $10$ were previously used in Stage~C for $G_s$; and at most $2n^{0.5}$ of these vertices cannot be used because the edge to $v_{i-1}$ or $v'$ was previously used in Stage~C. Thus by~\eqref{eq:CONSTANTS} there are at least $(1-2\gamcore)\dStageA^2\deltnonspanning n$ valid choices.
  
  Next, we embed the paths $Q_{s_i}$ for $i\in[q]$. We use the same strategy, except that for each $i$ we embed separately a path from the $M$-neighbour of $x_i$ to one endvertex of $Q_{s_i}$, and a path from the $M$-neighbour of $y_i$ to the other endvertex of $Q_{s_i}$, choosing these two paths with the same number of internal vertices up to an error at most one such that the total number of vertices is $v(Q_{s_i})$. Note that this means each path has at least one internal vertex, and hence by the estimates above this too is always possible.
  
  Finally, we embed the remaining paths $Q_s$, for which we have $s\in\cJ_0\setminus (\cJ_{\boxdot}\cup \{s_1,\ldots,s_q\} )$. We use the same strategy again, and the same estimates show that this is possible.
  
  If $n$ is odd, we have now a collection of embeddings $\phi^\mathbf{C-}_s$ for some $s\in\cJ$. Whether $n$ is odd or even, if we have not yet defined $\phi^\mathbf{C-}_s$ we let it be equal to $\phi^\mathbf{A}_s$, and we let $\HStageCMinus$ be the subgraph of $\HStageB$ obtained by removing all the edges used in these embeddings, and the vertex $\boxdot$. For each $s\in\cJ$, we let $\SpecPaths^{**}_s$ be the set of paths in $\SpecPaths_s$ which are not embedded by $\phi^\mathbf{C-}_s$, and we set $\imCMinus(s)=\im (\phi^\mathbf{C-}_s)$.
   
The critical properties we should note are the following. We have \[|\SpecPaths^{**}_s|\ge|\SpecPaths_s|-1\,.\] The density $\dStageCMinus$ of $\HStageCMinus$ is $\dStageA\pm O(n^{-1})$. By~\ref{B:degrees} and by construction of the embeddings above, for every vertex $v\in V(\HStageCMinus)$ we have $\deg_{\HStageCMinus}(v)\ge\deg_{\HStageA}(v)-10Dcn$. For each $s\in\cJ$, we have $\big|\imCMinus(s)\setminus\imA(s)\big|\le 10$. Finally, for each $v\in V(\HStageCMinus)$, the number of $s\in\cJ$ such that $v\in\imCMinus(s)\setminus\imA(s)$ is at most $10Dcn$. To see that the last of these is true, note that if $v$ is used in Stage~B or~C$^-$ by a given graph $G_s$, then we also use at least one edge at $v$ to embed an edge of $G_s$ in one of these stages, and we use at most $10Dcn$ edges at $v$ in these two stages. What we obtain from these observations and~\eqref{eq:CONSTANTS} is the following, which we will use repeatedly in the rest of this stage.
  \begin{enumerate}[label=($\ast$)]
   \item\label{C:sameasA} For each condition from Lemma~\ref{lem:StageA}, the same statement holds after Stage~C$^-$, replacing $\gamcore$ with $2\gamcore$.
  \end{enumerate}
 
We next label $V(\HStageCMinus)=\{\boxminus_i,\boxplus_i:i\in[\lfloor n/2\rfloor]\}$. Consider such a labelling chosen uniformly at random. We claim that this gives a quasirandom setup.
  
\begin{lemma}\label{lem:C:split}
 With positive probability, the uniform random labelling $V(\HStageCMinus)=\{\boxminus_i,\boxplus_i:i\in[\lfloor n/2\rfloor]\}$, together with $\cJ_i$ for $i=0,1,2$, with the path-forest $F^\mathbf{C-}_s$ consisting of the paths of $\SpecPaths^{**}_s$ for each $s\in\cJ$, with the used sets $U^\mathbf{C-}_s=\im\phi^\mathbf{C-}_s$, with the anchor sets $A^\mathbf{C-}_s$ given by $\phi^\mathbf{C-}_s$ applied to the leaves of $\SpecPaths^{**}_s$ and with $I^{\mathbf{C-}}_s=\{j\in[\lfloor n/2\rfloor]:\boxminus_j,\boxplus_j\notin \imCMinus(s)\}$, is a $(20\LC\gamcore,\LC,\dStageCMinus,\dStageCMinus)$-quasirandom setup. Furthermore, for each $u\in V(\HStageCMinus)$ and $j\in[\lfloor n/2\rfloor]$ such that $u\notin\{\boxminus_j,\boxplus_j\}$ we have
 \[\big|\{s\in\cJ_0:u\in A^\mathbf{C-}_s,j\in I^\mathbf{C-}_s\}\big|=(1\pm 2\gamcore)|\cJ_0|\cdot 2\sigmKJ(\deltnonspanning+10\sigmKJ)^2\,,\]
 and for each $s\in\cJ_0$ we have $|I^\mathbf{C-}_s|=(1\pm 10\gamcore)\tfrac{n}{2}\big(1-\tfrac{|U^\mathbf{C-}_s|}{n}\big)^2$.
 \end{lemma}
\begin{proof}

We split the proof of this lemma into the following claims.

The following claim takes care of~\ref{quasi:PathsWithOneVertex}.
\begin{claim}\label{cl:C:PathWithOneVertexSplit}
For every $v\in V(\HStageCMinus)$, $a\in\{\boxminus,\boxplus\}$ and $i\in\{1,2\}$, let $X_{v,V_a,\cJ_i}$ be a random variable counting the quantity from~\ref{quasi:PathsWithOneVertex}, i.e.~the number of $s\in\cJ_i$ such that $F^\mathbf{C-}_s$ has a path with leaves $x$ and $y$ such that $\phi^{\mathbf{C-}}_s(x)=v$ and $\phi^{\mathbf{C-}}_s(y)\in V_a$. Then $\Exp[X_{v,V_a,\cJ_i}]=(1\pm 2\gamcore)\tfrac{1}{2n} \sum_{s\in\cJ_i}|A^\mathbf{C-}_s|$.

Further, $\Prob\big[\big|X_{v,V_a,\cJ_i}-\Exp[X_{v,V_a,\cJ_i}]\big|>n^{0.9}\big]<\exp(-n^{0.1})$.
\end{claim}
\begin{claimproof}[Proof of Claim~\ref{cl:C:PathWithOneVertexSplit}]
  Given $v$, $a$ and $i$, recall that by~\ref{A:PathsWithOneVertex} and~\ref{C:sameasA}, we have that the number of $P\in\bigcup_{s\in \cJ_i}\SpecPaths^{**}_s$ for which $v\in \{\phi^\mathbf{C-}_s(\leftpath_0(P)),\phi^\mathbf{C-}_s(\rightpath_0(P))\}$ is equal to
  \[(1\pm \gamcore)\tfrac{1}{n} \sum_{s\in\cJ_i}|A^\mathbf{C-}_s|\,.\]
  Each such path is counted in $X_{v,V_a,\cJ_i}$ if and only if the unique vertex of
  \[\Big\{\phi^\mathbf{C-}_s\big(\leftpath_0(P)\big),\phi^\mathbf{C-}_s\big(\rightpath_0(P)\big)\Big\}\setminus\{v\}\]
  ends up in the set $V_a$; the probability of this is $\tfrac12$. Hence, the statement about the expectation follows. 

Let us now turn to the statement about the concentration. To this end, we model the random splitting in a particular way. We consider a space $\Omega=[0,1]^{V(\HStageCMinus)}$ (equipped with the Lebesgue measure), together with a fixed (arbitrary) linear order on $V(\HStageCMinus)$ used only to break ties. Suppose that $\vec{x}\in\Omega$. Order $V(\HStageCMinus)$ by increasing value of $\vec{x}(v)$, breaking any ties using the fixed order; based on $\vec{x}$ we then define $\Vmin$ to be the first $\lfloor \frac n2 \rfloor$ vertices $v\in V(\HStageCMinus)$ in this order and $\Vplus$ to be the remaining $\lfloor \frac n2 \rfloor$ vertices.
Taking $\vec{x}$ at random yields the right distribution of $\Vmin$ and $\Vplus$; moreover, with this convention $\Vmin$ and $\Vplus$ are defined for \emph{every} $\vec{x}\in\Omega$, with no exceptional set required. This way, we can think of $X_{v,V_a,\cJ_i}$ as a random variable on $\Omega$.
Observe that when we change a variable $\vec{x}\in \Omega$ in one coordinate, exactly two vertices exchange sides, and so $X_{v,V_a,\cJ_i}(\vec{x})$ changes by at most $2n^{0.3}$ by~\ref{A:PathsNotSquashed}. Hence, Lemma~\ref{lem:McDiarmidInequ} gives that $\Prob\big[\big|X_{v,V_a,\cJ_i}-\Exp[X_{v,V_a,\cJ_i}]\big|>n^{0.9}\big]<2\exp\left(-\frac{2(n^{0.9})^2}{(2n^{0.3})^2 n}\right)<\exp(-n^{0.1})$.
\end{claimproof}

The following claim takes care of~\ref{quasi:distributionOfAnchors} and~\ref{quasi:anchorsets}.
\begin{claim}\label{cl:C:distributionOfAnchors}
For all $s\in \cJ$, for all $v\notin \imCMinus(s)$ and for all $a,b\in\{\boxminus,\boxplus\}$, let $A^{\mathbf{C-},b}_s$ denote the subset of $u\in A^\mathbf{C-}_s$ such that the path of $\SpecPaths^{**}_s$ with leaves $\big(\phi^\mathbf{C-}_s\big)^{-1}(u)$ and $y$ has $\phi^\mathbf{C-}_s(y)\in V_b$. With probability at least $1-\exp(-n^{0.1})$ we have that
\[|V_a\cap A^{\mathbf{C-},b}_s\cap \NBH_{\HStageCMinus}(v)|=(1\pm 4\gamcore)\dStageCMinus\frac{|A^\mathbf{C-}_s|}{4}\quad\text{ and }\quad|V_a\cap A^{\mathbf{C-},b}_s|=(1\pm 3\gamcore)\frac{|A^\mathbf{C-}_s|}{4}\,.\]
\end{claim}
\begin{claimproof}[Proof of Claim~\ref{cl:C:distributionOfAnchors}]
By~\ref{A:vertex-covering} and~\ref{C:sameasA}, $|A^\mathbf{C-}_s\cap \NBH_{\HStageCMinus}(v)|=(1\pm 2\gamcore)\dStageCMinus |A^\mathbf{C-}_s|$. For each $x\in A^\mathbf{C-}_s\cap \NBH_{\HStageCMinus}(v)$, the probability that $x\in A^{\mathbf{C-},b}_s\cap V_a$ is $\tfrac14\big(1\pm\tfrac2n\big)$: this event asks that two distinct vertices of $V(\HStageCMinus)$, namely $\phi^\mathbf{C-}_s(x)$ and the image of the other leaf of its path, fall into the prescribed sides of a uniform random balanced split of $V(\HStageCMinus)$, which has probability $\tfrac14\big(1\pm\tfrac{1}{|V(\HStageCMinus)|-1}\big)$. Hence the expectation of $|A^{\mathbf{C-},b}_s\cap \NBH_{\HStageCMinus}(v)\cap V_a|$ is $(1\pm 3\gamcore)\dStageCMinus|A^\mathbf{C-}_s|/4$. Similarly $|V_a\cap A^\mathbf{C-,b}_s|$ has expectation $(1\pm 2\gamcore)\tfrac14|A^\mathbf{C-}_s|$. The concentration follows as in Claim~\ref{cl:C:PathWithOneVertexSplit}. In particular, observe that a change of one coordinate in a variable $\vec{x}\in \Omega$ changes $|V_a\cap A^{\mathbf{C-},b}_s\cap \NBH_{\HStageCMinus}(v)|$ or $|V_a\cap A^{\mathbf{C-},b}_s|$ by at most~$2$.
\end{claimproof}

The following claim is a first step towards~\ref{quasi:indexquasi}.
\begin{claim}
	\label{cl:C-quasirandom-firstStage}
	Let $a\in \{\boxminus, \boxplus\}$ and $S\subseteq V(\HStageCMinus)$, and $T\subseteq \cJ$ be such that $|S|\le 2\LC$ and $|T|\le 2\LC$. Then the following holds with probability at least $1-\exp(-n^{0.55})$.
	\begin{equation}\label{eq:C-quasi1}
	\left|V_a\cap \NBH_{\HStageCMinus}(S)\setminus \bigcup_{s\in T}U^\mathbf{C-}_s\right|
	=(1\pm 4\gamcore)\dStageCMinus^{|S|}\prod_{s\in T}\left(1-\frac{|U^\mathbf{C-}_s|}{n}\right)\cdot \frac{n}{2}\;.
	\end{equation}
\end{claim}

\begin{claimproof}[Proof of Claim~\ref{cl:C-quasirandom-firstStage}]
We think of the random labelling described above as follows. Firstly, we split $V(\HStageCMinus)$ into $\Vmin$, $\Vplus$ (where $|\Vmin|=|\Vplus|=\lfloor \frac n2\rfloor$) at random.  Secondly, we consider random labellings $\{\boxminus_i\}_{i\in [\lfloor \frac n2\rfloor]}$ of $\Vmin$ and $\{\boxplus_i\}_{i\in [\lfloor \frac n2\rfloor]}$ of $\Vplus$. Note that after the first stage (without considering the random labelling), we are already able to determine~\eqref{eq:C-quasi1}. By~\ref{A:megaquasirandomness1} and~\ref{C:sameasA} we have 
\[
\left|\NBH_{\HStageCMinus}(S) \setminus \bigcup_{s\in T}U^\mathbf{C-}_s\right|
=(1\pm 2\gamcore)\dStageCMinus^{|S|} n \prod_{s\in T}\left(1-\frac{|U^\mathbf{C-}_s|}{n}\right)\;.\] 
Since each vertex of $\NBH_{\HStageCMinus}(S) \setminus \bigcup_{s\in T}U^\mathbf{C-}_s$ has probability $\tfrac12$ of being in $V_a$, we obtain
\begin{align*}
\Exp\left[\left|V_a\cap \NBH_{\HStageCMinus}(S)\setminus \bigcup_{s\in T}U^\mathbf{C-}_s\right|\right]&=(1\pm 2\gamcore)\dStageCMinus^{|S|} \frac n2 \prod_{s\in T}\left(1-\frac{|U^\mathbf{C-}_s|}{n}\right)\;.	
\end{align*}

Furthermore, the random variable $|V_a\cap \NBH_{\HStageCMinus}(S)\setminus \bigcup_{s\in T}U^\mathbf{C-}_s|$ is exponentially concentrated by Fact~\ref{fact:hypergeometricBasicProperties}. 
\end{claimproof}

We can now finish establishing~\ref{quasi:indexquasi} and our bounds on $|I^\mathbf{C-}_s|$ for $s\in\cJ_0$ by proving the following claim.

\begin{claim}
	\label{cl:C-megaquasirandomness}
Let $S_1,S_2\subseteq \Vmin\cup\Vplus$,  and $T_1,T_2,T_3\subseteq \cJ$ be
pairwise disjoint sets of size at most $\LC$, as in Definition~\ref{def:index-quasirandom}. Then with probability at least $1-\exp(-n^{0.5})$ we have that the set
\[
X:=\left\{j\in [\lfloor\tfrac{n}{2}\rfloor]: \boxminus_j\in \NBH_{\HStageCMinus}(S_1)\setminus \bigcup_{s\in T_1\cup T_3}U^\mathbf{C-}_s, \boxplus_j\in\NBH_{\HStageCMinus}(S_2)\setminus\bigcup_{s\in T_2\cup T_3}U^\mathbf{C-}_s   \right\}
\]
has size \[(1\pm 10\gamcore)\dStageCMinus^{|S_1|+|S_2|}\prod_{s\in T_1\cup T_2}\left(1-\frac{|U^\mathbf{C-}_s|}{n}\right)\cdot \prod_{s\in T_3}\left(1-\frac{|U^\mathbf{C-}_s|}{n}\right)^2\cdot \frac{n}{2}\;.
\]
\end{claim}
\begin{claimproof}[Proof of Claim~\ref{cl:C-megaquasirandomness}]
We consider the same random experiment as in Claim~\ref{cl:C-quasirandom-firstStage}. With probability at least $1-\exp(-n^{0.5})$, $V(\HStageCMinus)$ is split between $\Vmin$ and $\Vplus$  so that~\eqref{eq:C-quasi1} holds simultaneously for all choices of $a\in\{\boxminus,\boxplus\}$, $S$, and $T$. 

Suppose now that $S_1,S_2,T_1,T_2,T_3$ are given.
We expose the labellings $\{\boxminus_j\}_{j\in[\lfloor\tfrac{n}{2}\rfloor]}$ and  $\{\boxplus_j\}_{j\in[\lfloor\tfrac{n}{2}\rfloor]}$. 
 Consider the sets 
\begin{align*}
Y_\boxminus&:=\left\{j\in [\lfloor\tfrac{n}{2}\rfloor]: \boxminus_j\in  \NBH_{\HStageCMinus}(S_1)\setminus\bigcup_{s\in T_1\cup T_3}U^\mathbf{C-}_s\right\} \;
\mbox{and}\\
Y_\boxplus&:=\left\{j\in [\lfloor\tfrac{n}{2}\rfloor]: \boxplus_j\in  \NBH_{\HStageCMinus}(S_2)\setminus\bigcup_{s\in T_2\cup T_3}U^\mathbf{C-}_s\right\} \;.
\end{align*}
We have that $X=Y_\boxminus\cap Y_\boxplus$. Thus, $|X|$ has hypergeometric distribution with parameters $(\lfloor\tfrac{n}{2}\rfloor,|Y_\boxminus|,|Y_\boxplus|)$. Note that the sizes of $Y_\boxminus$ and $Y_\boxplus$ do not depend on the exposed labelling but only on the sets $\Vmin$ and $\Vplus$. We apply~\eqref{eq:C-quasi1} with $a=\boxminus$, $S:=S_1$, $T:=T_1\cup T_3$ and get
\[|Y_\boxminus|=(1\pm 4\gamcore)\dStageCMinus^{|S_1|}\cdot \frac{n}{2}\prod_{s\in T_1\cup T_3}\left(1-\frac{|U^\mathbf{C-}_s|}{n}\right)\;. \]
Similarly, 
\[|Y_\boxplus|=(1\pm 4\gamcore)\dStageCMinus^{|S_2|}\cdot \frac{n}{2}\prod_{s\in T_2\cup T_3}\left(1-\frac{|U^\mathbf{C-}_s|}{n}\right)\;. \]
Now, Fact~\ref{fact:hypergeometricBasicProperties} yields the statement.
\end{claimproof}

Recall that for $s\in\cJ_0$, the set $I^\mathbf{C-}_s$ is the set of $i\in[\lfloor n/2\rfloor]$ such that $\boxminus_i,\boxplus_i\notin U^\mathbf{C-}_s$. Hence, in Claim~\ref{cl:C-megaquasirandomness}, $I^\mathbf{C-}_s$ equals $X$ when we choose $S_1,S_2,T_1,T_2$ to be the empty set and $T_3=\{s\}$.
The good event of this claim then gives $|I^\mathbf{C-}_s|=(1\pm 10\gamcore)\tfrac{n}{2}\big(1-\tfrac{|U^\mathbf{C-}_s|}{n}\big)^2$.
Moreover, having $S_1,S_2\subseteq \Vmin\cup\Vplus$  and $T_1,T_2,T_3\subseteq \cJ$ as in Definition~\ref{def:index-quasirandom}, each of size at most $\LC$, we then get that the set $X$ described in Claim~\ref{cl:C-megaquasirandomness} is likely to be of size
\begin{align*}
& (1\pm 10\gamcore)\dStageCMinus^{|S_1|+|S_2|}\frac{n}{2}\prod_{s\in T_1\cup T_2}\left(1-\frac{|U^\mathbf{C-}_s|}{n}\right)\cdot \prod_{s\in T_3} \frac{|I^\mathbf{C-}_s|}{(1\pm 10\gamcore)\tfrac{n}{2}} \\
& = (1\pm 20\LC\gamcore)\dStageCMinus^{|S_1|+|S_2|}\frac{n}{2}\prod_{s\in T_1\cup T_2}\left(1-\frac{|U^\mathbf{C-}_s|}{n}\right)\cdot \prod_{s\in T_3} \frac{|I^\mathbf{C-}_s|}{n/2}
\end{align*}
as required for~\ref{quasi:indexquasi}.

Next, we turn to verifying~\ref{quasi:termVtxNbs}.
\begin{claim}\label{cl:C:IsInNeighbourhoods}
Suppose that $s\in\cJ_0$ and $v\in V(\HStageCMinus)\setminus U^\mathbf{C-}_s$. Then with probability at least $1-\exp(-n^{0.4})$, we have
\[\big|\{\boxminus_i\in \NBH_{\HStageCMinus}(v):i\in I^\mathbf{C-}_s\}\big|,\big|\{\boxplus_i\in \NBH_{\HStageCMinus}(v):i\in I^\mathbf{C-}_s\}\big|=(1\pm 25\gamcore)\dStageCMinus|I^\mathbf{C-}_s|\,.\]
\end{claim}
\begin{claimproof}[Proof of Claim~\ref{cl:C:IsInNeighbourhoods}]
We prove the statement for the $\boxplus_i$; that for the $\boxminus_i$ is symmetric. By Claim~\ref{cl:C-quasirandom-firstStage}, after revealing the sets $\Vmin$ and $\Vplus$, we have $(1\pm 4\gamcore)\dStageCMinus\big(1-\tfrac{|U^\mathbf{C-}_s|}{n}\big)\tfrac{n}{2}$ vertices in $\NBH_{\HStageCMinus}(v)\cap (\Vplus\setminus U^\mathbf{C-}_s)$, and $(1\pm 4\gamcore)\big(1-\tfrac{|U^\mathbf{C-}_s|}{n}\big)\tfrac{n}{2}$ vertices in $\Vmin\setminus U^\mathbf{C-}_s$. We then consider the uniform random labelling of the $\boxminus_i$ and $\boxplus_i$ within $\Vmin$ and $\Vplus$, as in Claim~\ref{cl:C-megaquasirandomness}. We have $\boxplus_i$ in the desired set if and only if $\boxminus_i\in \Vmin\setminus U^\mathbf{C-}_s$ and $\boxplus_i\in \NBH_{\HStageCMinus}(v)\cap (\Vplus\setminus U^\mathbf{C-}_s)$, so the size $\big|\{\boxplus_i\in \NBH_{\HStageCMinus}(v):i\in I^\mathbf{C-}_s\}\big|$ has hypergeometric distribution with parameters
\[\Big(\lfloor\tfrac{n}{2}\rfloor,(1\pm 4\gamcore)\big(1-\tfrac{|U^\mathbf{C-}_s|}{n}\big)\tfrac{n}{2},(1\pm 4\gamcore)\dStageCMinus\big(1-\tfrac{|U^\mathbf{C-}_s|}{n}\big)\tfrac{n}{2}\Big)\,.\]
Now, using Fact~\ref{fact:hypergeometricBasicProperties} 
and the likely size of $I^\mathbf{C-}_s$, we conclude that the lemma statement holds.
\end{claimproof}
 
For~\ref{quasi:imagecaps}, let $S_1,S_2\subset V(\HStageCMinus)$ be disjoint, and additionally let $T\subset[\lfloor n/2\rfloor]$ be such that for each $j\in T$ we have $\boxminus_j,\boxplus_j\notin S_1\cup S_2$, with each set of size at most $\LC$. 
For $i=0$ we apply~\ref{A:typicality} with input $S_1\cup\{\boxminus_j,\boxplus_j\,:\,j\in T\}$, $S_2$ and $\cJ_0$ together with~\ref{C:sameasA}, and obtain
\begin{align*}   
   & \big|\{s\in\cJ_0\,:\,S_1\cap U_s^{\mathbf{C}-}=\emptyset, S_2\subset U_s^{\mathbf{C}-},T\subset I_s^{\mathbf{C}-}\}\big| \\
   & = \big|\{s\in\cJ_0\,:\, (S_1\cup\{\boxminus_j,\boxplus_j\,:\,j\in T\})\cap U_s^{\mathbf{C}-}=\emptyset, S_2\subset U_s^{\mathbf{C}-}\}\big| \\
   & =
(1\pm 2\gamcore)(\deltnonspanning+10\sigmKJ)^{|S_1|+2|T|}(1-\deltnonspanning-10\sigmKJ)^{|S_2|}|\cJ_0|    \\
  & 
  =(1\pm 20\LC\gamcore)\sum_{s\in\cJ_0}\frac{(n-|U^\mathbf{C-}_s|)^{|S_1|}|U^\mathbf{C-}_s|^{|S_2|}(2|I^\mathbf{C-}_s|)^{|T|}}{n^{|S_1|+|S_2|+|T|}}
\end{align*}
where the last equation is obtained by
using that
$|U^\mathbf{C-}_s|=(1-\deltnonspanning-10\sigmKJ + o(1))n$ and again
$|I^\mathbf{C-}_s|=(1\pm 10\gamcore)\tfrac{n}{2}\big(1-\tfrac{|U^\mathbf{C-}_s|}{n}\big)^2
=(1\pm 11\gamcore)\tfrac{n}{2}\big(\deltnonspanning + 10\sigmKJ\big)^2\,.
$
This is precisely the condition required for $i=0$ in~\ref{quasi:imagecaps}.
For $i\in \{1,2\}$ we use input $S_1$, $S_2$, $\cJ_i$ instead, and analogously we obtain the remaining conditions  required for~\ref{quasi:imagecaps}.

For~\ref{quasi:NumAnchors}, given $i\in\{0,1,2\}$ and $v\in V(\HStageCMinus)$, we apply~\ref{A:PathsWithOneVertex} and~\ref{C:sameasA} which gives~\ref{quasi:NumAnchors}.
For~\ref{quasi:NumAnchorsNoIm}, given $i\in\{0,1,2\}$ and $u\neq v\in V(\HStageCMinus)$, we apply~\ref{A:PathsOneVtxIm} and~\ref{C:sameasA} which gives~\ref{quasi:NumAnchorsNoIm}.

Finally, we prove the required bounds on $\big|\{s\in\cJ_0:u\in A^\mathbf{C-}_s,j\in I^\mathbf{C-}_s\}\big|$ for $u\in V(\HStageCMinus)$ and $j$ an index such that $u\notin\{\boxminus_j,\boxplus_j\}$. By~\ref{A:PathsOneVtxIm}\ref{enu:PathsTwoVtxIm0} and~\ref{C:sameasA}, for any three distinct vertices $u,v,v'$ of $V(\HStageCMinus)$, we have
\[\big|\{s\in\cJ_0:u\in A^\mathbf{C-}_s,v,v'\notin U^\mathbf{C-}_s\}\big|=(1\pm 2\gamcore)|\cJ_0|\cdot 2\sigmKJ(\deltnonspanning+10\sigmKJ)^2\,.\]
We apply this with $v=\boxminus_j$ and $v'=\boxplus_j$, and note that by definition $j\in I^\mathbf{C-}_s$ if and only if $\boxminus_j,\boxplus_j\notin U^\mathbf{C-}_s$.

Since the error probabilities in Claims~\ref{cl:C:PathWithOneVertexSplit}--\ref{cl:C:IsInNeighbourhoods} are superpolynomially small and since there are only polynomially many instances of~\ref{quasi:indexquasi}--\ref{quasi:NumAnchorsNoIm}, we conclude that the random labelling satisfies all of them with positive probability. In particular, there exists one such labelling, completing the proof of Lemma~\ref{lem:C:split}.
\end{proof}

To complete the proof of Lemma~\ref{lem:StageCNew}, we need to embed a few more paths, from the graphs $G_s$ with $s\in\cJ_2$, such that all edges of the form $\boxminus_i\boxplus_i$ are used. We will ensure that we use only few vertices of each $G_s$, and that we embed in total few edges at any given vertex, in this step. This ensures that we preserve our quasirandom setup with only a little degradation of the parameters.

\begin{lemma}\label{lem:C:killtermedges}
	There are sets $\SpecPaths_s^\diamond\subseteq\SpecPaths_s^{**}$, $s\in \cJ_2$ and mappings $(\phi^\diamond_s)_{s\in \cJ_2}$ such that 
	\begin{enumerate}[label=\rom]
		\item \label{C:en:middleedges-embedding}
		For each $s\in \cJ_2$, the mapping $\phi^\diamond_s$ is an embedding of the paths $\SpecPaths^\diamond_s$ into~$\HStageCMinus$ that agrees with $\phi^{\mathbf{C-}}_s$ on the leaves of $\SpecPaths^{\diamond}_s$. For each edge $uv\in E(\HStageCMinus)$, there is at most one $\phi^\diamond_s$ which uses $uv$, and we let $\HStageC$ be the graph obtained from $\HStageCMinus$ by removing all the edges used by $(\phi^\diamond_s)_{s\in\cJ_2}$.
		\item \label{C:en:middleedges-erased} For each $i\in \left[\lfloor\frac{n}{2}\rfloor\right]$ we have $\boxminus_i\boxplus_i\notin E(\HStageC)$.
		\item \label{C:en:middleedges-keepquasirandomness}
		For each $v\in \Vmin\cup\Vplus$ we have $\deg_{\HStageC}(v)>\deg_{\HStageCMinus}(v)-n^{0.6}$, and for each $s\in\cJ_2$ we have $|\SpecPaths^\diamond_s|\le n^{0.5}$.
	\end{enumerate}
\end{lemma}

\begin{proof}
First we specify the sets $\SpecPaths^\diamond_s, s\in \cJ_2$. Initially, set $\SpecPaths^\diamond_s=\emptyset$ for all $s\in \cJ_2$. In order to comply with~\ref{C:en:middleedges-keepquasirandomness}, we need to keep track of how much we map something to a vertex, as any such mapping will decrease its degree by $2$. If we embed too many paths anchored at the same vertex, this also will decrease the degree of this vertex too much, and thus we should choose the sets $\SpecPaths^\diamond_s, s\in \cJ_2$ accordingly. 
For each vertex $v\in V(H)$, let us use a counter $c(v)$ initially set to $2$ if $v\in\{\boxminus_i,\boxplus_i\}$ such that $\boxminus_i\boxplus_i\in E(\HStageCMinus)$ and otherwise to $0$. For each path included in $\SpecPaths^\diamond_s, s\in \cJ_2$ anchored at $v$, we shall increase $c(v)$ by one and for any internal vertex of a path from $\SpecPaths^\diamond_s, s\in \cJ_2$ mapped to $v$, we shall raise $c(v)$ by $2$. Let $U:=\{v\in V(H)\,:\,c(v)>\sqrt{n}\}$ be the set of vertices that get dangerously overloaded. Observe that initially $U=\emptyset$. We shall embed at most $n/2$ paths in this lemma, each of which raises the total of the counters by $1+1+2\cdot 6=14$, while the counters start with total at most $2n$;
observe therefore that the size of $U$ will never exceed $9n/\sqrt{n}=9\sqrt{n}$.

	For each $i\in [\lfloor\frac{n}{2}\rfloor]$,  let $C(i):= \{s\in \cJ_2\,:\,\im\phi^{\mathbf{C-}}_s\cap \{\boxminus_i, \boxplus_i\}=\emptyset\}$. By Lemma~\ref{lem:C:split} and~\ref{quasi:imagecaps}, setting $S_1=\{\boxminus_i, \boxplus_i\}$ and $S_2:=\emptyset$, we have
	\[|C(i)|\ge (1-20\LC\gamcore)(\deltnonspanning+6\sigmKJ)^2|\cJ_2|\geByRef{eq:CONSTANTS} \deltnonspanning^2\sigmKJ n\;.\]
	 For each $i\in  [\lfloor\frac{n}{2}\rfloor]$ such that $\boxminus_i\boxplus_i\in E(\HStageCMinus)$ in succession, pick any $s_i\in C(i)$ such that $|\SpecPaths^\diamond_{s_i}|<\sqrt{n}-1$. Observe that the number of $s\in \cJ_2$ such that $|\SpecPaths^\diamond_s|\ge \sqrt{n}-1$ is fewer than $n/\sqrt{n}=\sqrt{n}<|C(i)|$, so that this is always possible.
	 Take any path $P_i\in \SpecPaths^{**}_{s_i}\setminus\SpecPaths^\diamond_{s_i}$ that is not anchored in $U$. This is possible because the paths of $\SpecPaths^{**}_{s_i}$ are vertex-disjoint and $\phi^\mathbf{C-}_{s_i}$ is injective, so each vertex of $U$ is the anchor image of at most one of them; since $|U|,|\SpecPaths^\diamond_{s_i}|\le 9\sqrt{n}$ and $\SpecPaths^{**}_{s_i}$ contains at least $n^{0.9}$ paths, such a $P_i$ exists.
	 
	 Add $P_i$ to $\SpecPaths^\diamond_{s_i}$ and increase the counter $c(v)$ of its anchors by one. Recall that $P_i$ has 6 internal vertices $x_1,\dots,x_6$ together with leaves $x_0,x_7$. We construct an embedding of $P_i$ as follows. We embed $\phi^\diamond_{s_i}(x_2)=\boxminus_i$ and $\phi^\diamond_{s_i}(x_3)=\boxplus_i$. We choose $\phi^\diamond_{s_i}(x_1)$ arbitrarily from $\NBH_{\HStageCMinus}\big(\phi^\mathbf{C-}_{s_i}(x_0),\boxminus_i\big)\setminus (\im\phi^\mathbf{C-}_{s_i} \cup U)$ such that the edges to $\phi^\mathbf{C-}_{s_i}(x_0)$ and $\boxminus_i$ have not previously been used in this lemma and which has not previously been used in this lemma for $G_{s_i}$. By Lemma~\ref{lem:C:split} and~\ref{quasi:indexquasi}, we have
	 \[\Big|\NBH_{\HStageCMinus}\big(\phi^\mathbf{C-}_{s_i}(x_0),\boxminus_i\big)\setminus \im\phi^\mathbf{C-}_{s_i}\Big|\ge(1-20\LC\gamcore)\dStageCMinus^2(\deltnonspanning+6\sigmKJ)n\,,\] 
	 of which at most $4n^{0.6}$ vertices are disallowed due to being in $U$ or the edges to $\boxminus_i$ or $\phi^\mathbf{C-}_{s_i}(x_0)$, used previously in this lemma, or due to being a vertex used for $G_{s_i}$ in this lemma. Thus the choice of $\phi^\diamond_{s_i}(x_1)$ is possible; and we increase $c\big(\phi^\diamond_{s_i}(x_1)\big)$ by $2$. We now choose similarly for each $j=4,5,6$ in succession a vertex $\phi^\diamond_{s_i}(x_j)$ from
	 \[\NBH_{\HStageCMinus}\big(\phi^\diamond_{s_i}(x_{j-1}),\phi^\mathbf{C-}_{s_i}(x_7)\big)\setminus (\im\phi^\mathbf{C-}_{s_i} \cup U)\]
	  such that the edges to $\phi^\diamond_{s_i}(x_{j-1})$ and $\phi^\mathbf{C-}_{s_i}(x_7)$ have not previously been used in this lemma and the vertex has not been used for $G_{s_i}$ in this lemma. By the same calculation as above, this is always possible. We increase the counter $c\big(\phi^\diamond_{s_i}(x_j)\big)$ by $2$.
\end{proof}

We now complete the proof of Lemma~\ref{lem:StageCNew}. For each $s\in\cJ_2$, we set $\SpecPaths^*_s:=\SpecPaths^{**}_s\setminus\SpecPaths^\diamond_s$, and we set $\phi^\mathbf{C}_s$ to be the union of $\phi^\mathbf{C-}_s$ and $\phi^\diamond_s$; note that these two maps have some common vertices in their domains (the leaves of the $\SpecPaths^\diamond_s$) but by construction they agree on these vertices. For each $s\in\cJ_0\cup\cJ_1$ we set $\SpecPaths^*_s=\SpecPaths^{**}_s$ and $\phi^\mathbf{C}_s=\phi^\mathbf{C-}_s$. Let $\dStageC=e(\HStageC)/\binom{n}{2}$, and observe that $\dStageC=\dStageCMinus\pm o(1)$.

We claim that $\HStageC$ as returned by Lemma~\ref{lem:C:killtermedges}, with the splitting into $\Vmin\cup\Vplus$ from Lemma~\ref{lem:C:split}, with for each $s\in\cJ$ the path-forest $F^\mathbf{C}_s$ consisting of the paths of $\SpecPaths^*_s$, with $U^\mathbf{C}_s=\im\phi^\mathbf{C}_s$, with $A^\mathbf{C}_s$ the image under $\phi^\mathbf{C}_s$ of the leaves of $F^\mathbf{C}_s$ and with the corresponding bijection, and with $I^\mathbf{C}_s$ the set of $i\in [\lfloor n/2\rfloor ]$ such that $\boxminus_i,\boxplus_i\notin \im\phi^\mathbf{C}_s$, is a $(100\LC\gamcore,\LC,\dStageC,\dStageC)$-quasirandom setup.

This follows from Lemma~\ref{lem:C:split}, which guarantees a $(20\LC\gamcore,\LC,\dStageCMinus,\dStageCMinus)$-quasirandom setup after stage $\mathbf{C^-}$, together with Lemma~\ref{lem:C:killtermedges} which tells us that all the quantities in the definition of a quasirandom setup have changed between these two stages by at most $10\LC n^{0.6}$.

Finally, the size bounds relating to $\cJ_0$ in Lemma~\ref{lem:StageCNew} have not changed from after stage $\mathbf{C^-}$ and hence hold by Lemma~\ref{lem:C:split}.

\section{Stage~D (Proof of Lemma~\ref{lem:StageDNew})}\label{sec:StageD}
In this stage, we need to accommodate all vertices of $P$, but the two middle ones for all $s\in\cJ_0$ and $P\in\SpecPaths^*_s$; this corresponds to finding the placement for the solid lines in Figure~\ref{fig:embedJ0} (the outermost vertices of $P$ were embedded already in Stage~A). Recall that the aim here is that after this stage, what remains to pack of each graph in $\cJ_0$ will be a path-forest in which each path has $2$ unembedded vertices, and the two leaves are embedded to a \emph{terminal pair}.

 We do that in two steps. First, in Lemma~\ref{lem:reroot-assign} we decide, for each $P\in\bigcup_{s\in\cJ_0}\SpecPaths^*_s$, what its terminal pair will be. To that end, given such a $P$ we let $y_{P}^{\mathsf{left}}$ denote the fifth vertex of $P$, and $y_{P}^{\mathsf{right}}$ the eighth vertex of $P$.\index{$y_{P}^{\mathsf{left}}$}\index{$y_{P}^{\mathsf{right}}$} That is, $y_P^\mathsf{left}$ is an endvertex of $P^{\mathsf{left}}:=\leftpath_4(P)$, and $y_P^\mathsf{right}$ is an endvertex of $P^{\mathsf{right}}:=\rightpath_4(P)$. We will choose an index $r_P\in[\lfloor n/2\rfloor]$ and embed $y_P^\mathsf{left}$ to $\boxminus_{r_P}$, and $y_P^\mathsf{right}$ to $\boxplus_{r_P}$. For this to be possible, for each $P\in\SpecPaths^*_s$ of course we need
 \[r_P\in I^\mathbf{C}_s:=\{i\in[\lfloor n/2\rfloor]\,:\,\boxminus_i,\boxplus_i\notin \imC(s)\}\,.\]
 Recall from~\ref{D:three} that, for each $s\in\cJ_0$, the set $\SpecShortPaths_s$ consists of the paths $\middlepath_3(P)$ with $P\in\SpecPaths^*_s$; equivalently, $\SpecShortPaths_s$ is the set of subpaths of $\SpecPaths^*_s$ with endvertices $y_P^\mathsf{left}$ and $y_P^\mathsf{right}$. That is, each path in $\SpecShortPaths_s$ is the middle 3~edges (4~vertices) of a path of $\SpecPaths^*_s$; this is what will remain to pack from $(G_s)_{s\in\cJ_0}$ at the conclusion of this stage.
 
 The content of Lemma~\ref{lem:reroot-assign} is then that we can make the assignment of indices $r_P$ and still have a quasirandom setup; in addition, we need one property regarding the set $A^\mathbf{C}_s$, which is the set of anchors after Stage~C. We refer to the situation provided by Lemma~\ref{lem:reroot-assign} as Stage~D$^-$. After Stage~D$^-$, we use Lemma~\ref{lem:pathpack} to argue that we can pack all the $P^\mathsf{left}$ and $P^\mathsf{right}$ and again obtain a quasirandom setup, which completes the proof of Lemma~\ref{lem:StageDNew}.

\begin{lemma}\label{lem:reroot-assign}
There exist indices $(r_P\in[\lfloor\tfrac{n}2\rfloor])_{s\in\cJ_0,P\in\SpecPaths^*_s}$ with the following properties.
\begin{enumerate}[label=\itmit{D\textsuperscript{-}\roman*}]
\item\label{D:reroot:rerootBas1} For each $s\in\cJ_0$, the indices $(r_P)_{P\in\SpecPaths^*_s}$ are distinct elements of $I^\mathbf{C}_s$.

\item\label{D:reroot:DdefimD-} For $s\in\cJ_1\cup\cJ_2$ we set $\phi^\mathbf{D-}_s=\phi^\mathbf{C}_s$. For $s\in\cJ_0$ we extend $\phi^\mathbf{C}_s$ to $\phi^\mathbf{D-}_s$ by setting $\phi^{\mathbf{D}^-}_s(y_{P}^{\mathsf{left}})
 	=\boxminus_{r_P}$ and
 	 $\phi^{\mathbf{D}^-}_s(y_{P}^{\mathsf{right}})
 	 	=\boxplus_{r_P}$.

\item\label{D:reroot:quasisetup} The graph $\HStageC$ with its labelling as $\Vmin\cup\Vplus$, together with the path-forests $F^\mathbf{D-}_s:=\SpecPaths^*_s$ for $s\in\cJ_1\cup\cJ_2$ and $F^\mathbf{D-}_s:=\SpecShortPaths_s$ for $s\in\cJ_0$, with used sets $U^\mathbf{D-}_s:=\imDMinus(s)$ for each $s\in\cJ$, with anchor sets $A^\mathbf{D-}_s$ being the image under $\phi^\mathbf{D-}_s$ of the leaves of $F^\mathbf{D-}_s$ (and with the corresponding bijection) and with $I_s:=I^\mathbf{D-}_s$ where $$I^\mathbf{D-}_s:=\{i\in[\lfloor n/2\rfloor]\,:\,\boxminus_i,\boxplus_i\in A^\mathbf{D-}_s\}$$ for each $s\in\cJ_0$, gives a $(600 \LD 4^{6\LD} \deltnonspanning^{-4\LD}(\sigmKJ^*)^{-3\LD} \LC \gamcore,2\LD,\dStageC,\dStageC)$-quasirandom setup. For convenience, we say this is the \emph{setup after Stage~D$^-$}.

\item\label{D:reroot:quasioldA} For each $u,v\in V(\HStageC)$ such that there is no $j$ with $u,v\in\{\boxminus_j,\boxplus_j\}$, the number of $s\in\cJ_0$ such that $u\in A^\mathbf{C}_s$ and $v\notin U^\mathbf{D-}_s$ is $(1\pm 1000\gamcore)\sum_{s\in\cJ_0}\tfrac{|A^\mathbf{C}_s|(n-|U^\mathbf{D-}_s|)}{n^2}$.
\end{enumerate}
\end{lemma}
\begin{proof}
Suppose that $s\in\cJ_0$ is fixed. After Stage~C we have by Lemma~\ref{lem:StageCNew}\ref{C:four} a $(100\LC\gamcore,\LC,\dStageC,\dStageC)$-quasirandom setup. Further, by Lemma~\ref{lem:StageCNew}\ref{C:five}, for each $s\in\cJ_0$ we have 
\begin{equation}\label{eq:sizeICs}
|I^\mathbf{C}_s|=(1\pm 10\gamcore)\tfrac{n}{2}\big(1-\tfrac{|U^\mathbf{C}_s|}{n}\big)^2
\eqByRef{eq:familyOrderJ0Equality}
(1\pm 11\gamcore)\tfrac{n}{2}(\deltnonspanning+10\sigmKJ)^2\;.
\end{equation}
By~\eqref{eq:CONSTANTS}, \eqref{eq:sizeSpecPathStar} and~\eqref{eq:sigmKJsigmKJ} this means
\begin{equation}\label{eq:followCONST1}
|I^\mathbf{C}_s|>\sigmKJ^* n\;,
\end{equation} and so we can select a set $J_s\subset I^\mathbf{C}_s$ uniformly at random of size $\sigmKJ^* n$. For each $s$, we choose an arbitrary pairing of the paths of 
$F^\mathbf{C}_s = \SpecPaths^*_s$ to the set $J_s$, pairing path $P$ to some index $r_P\in J_s$. Hence, \ref{D:reroot:rerootBas1} is satisfied. In particular we then set 
$\phi^{\mathbf{D}^-}_s(y_{P}^{\mathsf{left}})
 	=\boxminus_{r_P}$ and
 	 $\phi^{\mathbf{D}^-}_s(y_{P}^{\mathsf{right}})
 	 	=\boxplus_{r_P}$, as is described in~\ref{D:reroot:DdefimD-}. We obtain
 $I_s^{\mathbf{D-}}=\{r_P: P\in F^\mathbf{C}_s\}=J_s$.

We now aim to verify that~\ref{D:reroot:quasisetup} holds with high probability. We go through the properties of Definition~\ref{def:quasisetup} one by one.

For~\ref{quasi:indexquasi}, the critical observation is the following. By Fact~\ref{fact:hypergeometricBasicProperties}, for any given set $Y\subset I^\mathbf{C}_s$ of indices with $|Y|\ge cn$, when we reveal the choice of $J_s$, the probability that we do not have
\begin{equation}\label{eq:reroot:quasi}|Y\cap J_s|=|Y|\frac{|J_s|}{|I^\mathbf{C}_s|}\pm n^{0.9}\end{equation}
is at most $\exp(-n^{0.5})$. If $S_1,S_2,T_1',T_2',T_3'$ are any sets as in Definition~\ref{def:index-quasirandom}, each of size at most $4\LD$, by $(\LC,100\LC\gamcore,\dStageC,\dStageC)$-index-quasirandomness 
with respect to the sets $(U^\mathbf{C}_s)_{s\in\cJ}$ and $(I^\mathbf{C}_s)_{s\in\cJ_0}$
after Stage~C we have
\[\big|\mathbb{U}_{\HStageC}(S_1,S_2,T_1',T_2',T_3')\big|=(1\pm100\LC\gamcore)\dStageC^{|S_1|+|S_2|}\tfrac{n}{2}\prod_{\ell\in T_1'\cup T_2'}\tfrac{n-|\imC(\ell)|}{n}\prod_{\ell\in T_3'}\tfrac{|I^\mathbf{C}_\ell|}{n/2}\,.\]
Applying~\eqref{eq:reroot:quasi} iteratively for each $\ell\in T_3'$ successively and taking the union bound over choices of all these sets, we see that with probability at least $1-\exp(-n^{0.4})$, for each choice of sets we have
\begin{equation}\label{eq:reroot:quasi2}
 \Big|\mathbb{U}_{\HStageC}(S_1,S_2,T_1',T_2',T_3')\cap \bigcap_{\ell\in T_3'}J_\ell\Big|
 =(1\pm200\LC\gamcore)\dStageC^{|S_1|+|S_2|}\tfrac{n}{2}\prod_{\ell\in T_1'\cup T_2'}\tfrac{n-|\imC(\ell)|}{n}\prod_{\ell\in T_3'}\tfrac{|I^\mathbf{D-}_\ell|}{n/2}\,.
\end{equation}

Suppose this likely event occurs. In order to verify the index-quasirandomness after Stage~D$^-$ with respect to
$(U^{\mathbf{D}-}_s)_{s\in\cJ}$ and $(J_s)_{s\in\cJ_0}$, 
we have to look at the set
\begin{equation}\label{eq:quasi:stageDminus}
\Big(\mathbb{U}_{\HStageC}(S_1,S_2,T_1,T_2,T_3)\cap\bigcap_{\ell\in T_3}J_\ell\Big)\setminus\big(\cup_{\ell\in (T_1\cup T_2)\cap\cJ_0}J_\ell\big)
\end{equation}
for any sets $S_1,S_2,T_1,T_2,T_3$ as required by
Definition~\ref{def:index-quasirandom},
where the subtraction of the set
$\big(\cup_{\ell\in (T_1\cup T_2)\cap\cJ_0}J_\ell\big)$ is needed due to the fact that 
$U_\ell^{\mathbf{D}-}=U_\ell^{\mathbf{C}}\cup \{\boxminus_i,\boxplus_i:i\in J_\ell\}$
for $\ell\in \cJ_0$.
For short, let $T:=(T_1\cup T_2)\cap \cJ_0$
and $\mathbb{U}_{\HStageC} := \mathbb{U}_{\HStageC}(S_1,S_2,T_1,T_2,T_3)$. By inclusion-exclusion
we see that the size of the set in~\eqref{eq:quasi:stageDminus} is
\begin{align*}
& \Big|\Big(\mathbb{U}_{\HStageC}\cap \bigcap_{\ell\in T_3}J_\ell\Big)\setminus\big(\cup_{\ell\in T} J_\ell \big)\Big| \\
& 
= 
\Big|\mathbb{U}_{\HStageC}\cap \bigcap_{\ell\in T_3}J_\ell \Big| -
\Big|\Big(\mathbb{U}_{\HStageC}\cap \bigcap_{\ell\in T_3}J_\ell\Big)\cap \big(\cup_{\ell\in T} J_\ell \big)\Big| \\
&
=
\Big|\mathbb{U}_{\HStageC}\cap \bigcap_{\ell\in T_3}J_\ell\Big| -
\sum_{\emptyset \neq \mathcal{K}\subset T} (-1)^{|\mathcal{K}|-1}
\Big|\mathbb{U}_{\HStageC}\cap \bigcap_{\ell\in T_3}J_\ell \cap \bigcap_{\ell\in \mathcal{K}}J_\ell \Big| \\
&
=
\sum_{\mathcal{K}\subset T} (-1)^{|\mathcal{K}|} 
\Big|\mathbb{U}_{\HStageC}\cap \bigcap_{\ell\in T_3\cup \mathcal{K}}J_\ell\Big| .
\end{align*}
Moreover, note that 
\[
\mathbb{U}_{\HStageC}(S_1,S_2,T_1,T_2,T_3)\cap \bigcap_{\ell\in T_3\cup \mathcal{K}}J_\ell =
\mathbb{U}_{\HStageC}(S_1,S_2,T_1\setminus \mathcal{K},T_2\setminus \mathcal{K}, T_3\cup \mathcal{K})\cap \bigcap_{\ell\in T_3\cup \mathcal{K}}J_\ell
\]
by the fact that on both sides of the equation we restrict to indices
$i\in \bigcap_{\ell\in \mathcal{K}}J_\ell$, 
and for these we have $\boxminus_i,\boxplus_i \notin \bigcup_{\ell\in\mathcal{K}} U_{\ell}^{\mathbf{C}}$
as well as $i\in \bigcap_{\ell\in\mathcal{K}} I_{\ell}^\mathbf{C}$.  
Putting everything together we see that the size of the set in~\eqref{eq:quasi:stageDminus} is
\begin{align*}
\sum_{r=0}^{|T|} (-1)^r \sum_{\mathcal{K}\in \binom{T}{r}}
  \Big|\mathbb{U}_{\HStageC}(S_1,S_2,T_1\setminus \mathcal{K},T_2\setminus \mathcal{K}, T_3\cup \mathcal{K})\cap \bigcap_{\ell\in T_3\cup \mathcal{K}}J_\ell\Big|\,.
\end{align*}
Hence, applying~\eqref{eq:reroot:quasi2}
with $T_1'=T_1\setminus \mathcal{K}$, $T_2'=T_2\setminus \mathcal{K}$ and $T_3'=T_3\cup \mathcal{K}$ we get that the above sum equals
\begin{align*}
\sum_{r=0}^{|T|} (-1)^r \sum_{\mathcal{K}\in \binom{T}{r}}
(1\pm200\LC\gamcore)\dStageC^{|S_1|+|S_2|}\tfrac{n}{2}
\prod_{\ell\in (T_1\cup T_2)\setminus \mathcal{K}}\tfrac{n-|\imC(\ell)|}{n}\prod_{\ell\in T_3\cup \mathcal{K}}\tfrac{|I^\mathbf{D-}_\ell|}{n/2}\,.
\end{align*}
Now, taking the factor 
$(1\pm200\LC\gamcore)\dStageC^{|S_1|+|S_2|}\tfrac{n}{2}
\prod_{\ell\in (T_1\cup T_2)\setminus \cJ_0}\tfrac{n-|\imC(\ell)|}{n}\prod_{\ell\in T_3}\tfrac{|I^\mathbf{D-}_\ell|}{n/2}$ out, what remains is
\begin{align*}
\sum_{r=0}^{|T|} (-1)^r &\sum_{\mathcal{K}\in \binom{T}{r}}
\prod_{\ell\in T\setminus \mathcal{K}}\tfrac{n-|\imC(\ell)|}{n}\prod_{\ell\in \mathcal{K}}\tfrac{|I^\mathbf{D-}_\ell|}{n/2} 
= 
(1\pm \gamcore) \sum_{r=0}^{|T|} \binom{|T|}{r} (-1)^r 
(\deltnonspanning + 10\sigmKJ)^{|T|-r} 
(2\sigmKJ^*)^r \\
& = 
(1\pm \gamcore) 
\left(
\deltnonspanning + 10\sigmKJ - 2\sigmKJ^\ast
\right)^{|T|} =
(1\pm 2\gamcore) 
\prod_{\ell\in(T_1\cup T_2)\cap\cJ_0}\big(\tfrac{n-|\imC(\ell)|}{n}-\tfrac{|J_\ell|}{n/2}\big)
\end{align*}
where the first and the third equation use that $|\imC(\ell)|=(1-\deltnonspanning -10\sigmKJ)n \pm 11$ and
$|I^\mathbf{D-}_\ell|=|J_\ell|=\sigmKJ^\ast n$ for every $\ell\in \cJ_0$.
Putting everything together, we conclude that the set in~\eqref{eq:quasi:stageDminus} has size
\[
(1\pm 300\LC\gamcore)
\dStageC^{|S_1|+|S_2|}\tfrac{n}{2}\prod_{\ell\in (T_1\cup T_2)\setminus\cJ_0}\tfrac{n-|\imC(\ell)|}{n}\prod_{\ell\in T_3}\tfrac{|I^\mathbf{D-}_\ell|}{n/2}\prod_{\ell\in(T_1\cup T_2)\cap\cJ_0}\big(\tfrac{n-|\imC(\ell)|}{n}-\tfrac{|J_\ell|}{n/2}\big)\,.\]
Since we have $n-|\imC(s)|-2|J_s|=n-|\imDMinus(s)|$ for each $s\in\cJ_0$, and since for each $s\notin\cJ_0$ we have $\imDMinus(s)=\imC(s)$, this formula is the same as appears in the definition of index-quasirandomness after Stage~D$^-$. Hence, we conclude that after Stage~D$^-$ we have the required index-quasirandomness, verifying~\ref{quasi:indexquasi}.

For~\ref{quasi:anchorsets}, observe that if $s\in\cJ_1\cup\cJ_2$ then $A^b_s = A^{\mathbf{D-},b}_s$ has not changed from after Stage~C, and hence the desired statement holds. For $s\in\cJ_0$, by construction exactly half of $A^\mathbf{D-}_s$ is in each of $\Vmin$ and $\Vplus$. Similarly, for~\ref{quasi:PathsWithOneVertex}, nothing has changed for these statements from after Stage~C.

We now prove~\ref{quasi:termVtxNbs} and~\ref{quasi:distributionOfAnchors}. 
For $s\in\cJ_1\cup\cJ_2$, nothing has changed since after Stage~C, so we only need to consider $s\in\cJ_0$.
After Stage~C, by~\ref{quasi:termVtxNbs}, for any given $u\in V(\HStageC)$ and $s\in\cJ_0$ we have
\[\big|\{\boxminus_i\in \NBH_{\HStageC}(u):i\in I^\mathbf{C}_s\}\big|,\big|\{\boxplus_i\in \NBH_{\HStageC}(u):i\in I^\mathbf{C}_s\}\big|=(1\pm 100\LC\gamcore)\dStageC|I^\mathbf{C}_s|\,.\]
When we select $J_s$ uniformly at random from $I^\mathbf{C}_s$, by Fact~\ref{fact:hypergeometricBasicProperties}, with probability at least $1-\exp(-n^{0.5})$ we therefore have
\[\big|\{\boxminus_i\in \NBH_{\HStageC}(u):i\in J_s\}\big|,\big|\{\boxplus_i\in \NBH_{\HStageC}(u):i\in J_s\}\big|=(1\pm 200\LC\gamcore)\dStageC|J_s|\, ,\]
and taking a union bound, we get the above holds with high probability
for every $s\in \cJ_0$ and $u\in V(\HStageC)$. This immediately gives~\ref{quasi:termVtxNbs} after Stage~D$^-$.
Moreover, observe that for $s\in\cJ_0$ we have 
$A^\boxminus_s\cap\Vmin=A^\boxplus_s\cap\Vplus=\emptyset$
and hence it suffices to verify 
\eqref{eq:quasi:distAnchorsa1} and \eqref{eq:quasi:distAnchorsa2} for $a=\boxplus$, and \eqref{eq:quasi:distAnchorsa3} and \eqref{eq:quasi:distAnchorsa4} for $a=\boxminus$. Now, all of these conditions follow directly from the already checked property~\ref{quasi:termVtxNbs},
since
the size of
$A^{\mathbf{D-},\boxplus}_s\cap \Vmin=\{\boxminus_i:i\in I^\mathbf{D-}_s\}$ and $A^{\mathbf{D-},\boxminus}_s\cap\Vplus=\{\boxplus_i:i\in I^\mathbf{D-}_s\}$ equals $|J_s|$, and
since
\[
\{\boxminus_i\in \NBH_{\HStageC}(u):i\in J_s\} =
\{\boxminus_i\in \NBH_{\HStageC}(u) \cap A_s^{\mathbf{D-},\boxplus} \}
\]
and analogously with $\boxplus$ and $\boxminus$ exchanged. In particular, this identification shows that, for $s\in\cJ_0$, properties~\ref{quasi:termVtxNbs} and~\ref{quasi:distributionOfAnchors} are simply two equivalent ways of stating the same fact about $\{\boxminus_i,\boxplus_i:i\in I_s\}\cap\NBH_H(u)$: each can be derived from the other. In later stages, when we say that~\ref{quasi:termVtxNbs} holds ``by construction'' (or vice versa), this equivalence, which depends only on the anchoring structure of paths in $\SpecShortPaths_s$ for $s\in\cJ_0$ and so persists unchanged through the later stages, is what we mean.

We next move to~\ref{quasi:imagecaps}. Since nothing has changed from after Stage~C for $\cJ_1$ and $\cJ_2$, we only need to prove the statement for $\cJ_0$.
Since $|J_s|,|I^\mathbf{C}_s\setminus J_s|=\Omega(n)$, we see by a hypergeometric distribution that for any given disjoint sets $Q,Q'$ of indices of size at most $6\LD$, we have that $\Prob[Q\subset J_s\mid Q'\subset J_s]$ is equal to either 
\begin{align*}
\frac{\binom{|I_s^{\mathbf{C}}|-|Q'|-|Q|}{|J_s|-|Q'|-|Q|}}{\binom{|I_s^{\mathbf{C}}|-|Q'|}{|J_s|-|Q'|}}
& =
\prod_{j=0}^{|Q|-1} \frac{|J_s|-|Q'|-j}{|I_s^{\mathbf{C}}|-|Q'|-j} \\
& =
\prod_{j=0}^{|Q|-1} \frac{\sigmKJ^* n \pm 12 \LD}{(1\pm 12\gamcore)\tfrac{n}{2}(\deltnonspanning+10\sigmKJ^*)^2 \pm 12 \LD }
=
\big((1\pm 100\gamcore) \tfrac{2\sigmKJ^*}{(\deltnonspanning+10\sigmKJ^*)^2}\big)^{|Q|}
\end{align*} 
if $Q\subset I^\mathbf{C}_s$, and $0$ otherwise. We use this approximate independence to prove~\ref{quasi:imagecaps} inductively. Our induction statement is the following $(\ast)_k$, where $0\le k\le 6\LD$ is an integer. As we shall show, the inductive statement holds with high probability.

\begin{InductiveStatementAst}
Let $S$ be a set of vertices, and let $T,T'$ be pairwise disjoint sets of indices. Suppose in addition that there is no $i$ such that $\{\boxminus_i,\boxplus_i\}\subset S$, and that there is no $i\in T\cup T'$ such that $\{\boxminus_i,\boxplus_i\}\cap S\neq\emptyset$. Finally suppose that $|T'|\le k$, and that $|S|+|T|+|T'|\le 6\LD$. Let
\[X(S,T,T'):=\big\{s\in\cJ_0:S\subset\imC(s),T\subset I^\mathbf{C}_s,T'\subset J_s\big\}\,.\]
 Then we have
\begin{equation}\label{D:XST}
 \big|X(S,T,T')\big|=|\cJ_0|(\deltnonspanning+10\sigmKJ^*)^{2|T|}(1-\deltnonspanning-10\sigmKJ^*)^{|S|}(2\sigmKJ^*)^{|T'|}\pm 100(k+1)\LC\gamcore n\,,
\end{equation}
and we reveal only the intersection of $T'$ with the sets $J_s$ to show this.
\end{InductiveStatementAst}
\begin{claimproof}[The base case $(\ast)_0$] We have $T'=\emptyset$. This case is covered by~\ref{quasi:imagecaps} after Stage~C, applied with $S_1=\emptyset$, 
$S_2=S$ and $T$,
\begin{align*}
 \big|X(S,T,\emptyset)\big| 
 & = (1\pm 100\LC\gamcore) 
 \sum_{s\in \cJ_0} 
 \frac{|U_s^{\mathbf{C}}|^{|S|}\cdot (2|I_s^\mathbf{C}|)^{|T|}}{n^{|S|+|T|}} \\
 & = (1\pm 100\LC\gamcore)(1\pm 12\gamcore)^{|T|}\cdot |\cJ_0|(\deltnonspanning+10\sigmKJ^*)^{2|T|}(1-\deltnonspanning-10\sigmKJ^*)^{|S|} \\
  & = |\cJ_0|(\deltnonspanning+10\sigmKJ^*)^{2|T|}(1-\deltnonspanning-10\sigmKJ^*)^{|S|} 
 \pm 100\LC\gamcore n\,.
\end{align*}
\end{claimproof}
\begin{claimproof}[Case $(\ast)_k$ implies with high probability case $(\ast)_{k+1}$]
 For a given $k\le 6\LD$, suppose that we have established $(\ast)_k$; we want to establish a given case $S,T,T'$ of $(\ast)_{k+1}$.
Let $i\in T'$. The set $X(S,T,T')$ consists of those $s\in X(S,T\cup\{i\},T'\setminus\{i\})$ such that $i\in J_s$. As above, $(\ast)_k$ gives us the size of the former set, and then by Chernoff's inequality (Fact~\ref{fact:hypergeometricBasicProperties}), with probability at least $1-\exp(-n^{0.5})$ we have
\begin{align*}
 &|X(S,T,T')|\\
 &=(1\pm 100\gamcore) \tfrac{2\sigmKJ^*}{(\deltnonspanning+10\sigmKJ^*)^2} \left(|\cJ_0|(\deltnonspanning+10\sigmKJ^*)^{2|T|+2}(1-\deltnonspanning-10\sigmKJ^*)^{|S|}(2\sigmKJ^*)^{|T'|-1}\pm 100(k+1)\LC\gamcore n \right) \pm n^{0.9}
\end{align*}
which, rearranging, gives us the desired case of $(\ast)_{k+1}$. This finishes the induction step.
\end{claimproof}

We now suppose that all the above polynomially many good events occur, i.e.\ $(\ast)_{6\LD}$ is true.

\begin{claim}\label{cl:defineY}
 Suppose that $|S|\le4\LD$ and $|T|\le 2\LD$. Moreover, suppose again that there is no $i$ such that $\{\boxminus_i,\boxplus_i\}\subset S$, so that for every $Q\subseteq S$ we may let $I_Q$ be the set of all $i$ for which either $\boxplus_i\in Q$ or $\boxminus_i\in Q$, and then $|I_Q|=|Q|$. Define $Y(S,T):=\big\{s\in\cJ_0:S\subset\imDMinus(s),T\subset J_s\big\}$. Then we have 
\[
\big|Y(S,T)\big|
=
|\cJ_0|(1-\deltnonspanning-8\sigmKJ^*)^{|S|}(2\sigmKJ^*)^{|T|}\pm 600\LD 2^{6\LD}\LC\gamcore n
\;.
\]
\end{claim}
\begin{claimproof}[Proof of Claim~\ref{cl:defineY}]
To begin, we claim that
\begin{equation}\label{eq:UnionY}
	Y(S,T) = \dbigcup_{Q\subset S} X(S\setminus Q,\emptyset, T\cup I_Q)\;.
\end{equation}
First, we argue the disjointness in the union in~\eqref{eq:UnionY}. Assume for contradiction that there exist subsets $Q_1\neq Q_2$ of $S$ such that 
$X(S\setminus Q_1,\emptyset, T\cup I_{Q_1})\cap X(S\setminus Q_2,\emptyset, T\cup I_{Q_2})\neq \emptyset$, and let $s$ be an element of this intersection. Since $Q_1\neq Q_2$, at least one of the sets $Q_1\setminus Q_2$ and $Q_2\setminus Q_1$ is non-empty; as $Q_1$ and $Q_2$ play symmetric roles, we may assume that $Q_1\setminus Q_2\neq\emptyset$. Fix $v\in Q_1\setminus Q_2$ and let $i$ be the index with $v\in\{\boxminus_i,\boxplus_i\}$, so that $i\in I_{Q_1}$. On the one hand, $v\in S\setminus Q_2$, and $s\in X(S\setminus Q_2,\emptyset,T\cup I_{Q_2})$ tells us that $S\setminus Q_2\subseteq\imC(s)$; hence $v\in\imC(s)$. On the other hand, $s\in X(S\setminus Q_1,\emptyset,T\cup I_{Q_1})$ tells us that $I_{Q_1}\subseteq J_s$, so $i\in J_s\subseteq I^\mathbf{C}_s$, and therefore $\boxminus_i,\boxplus_i\notin\imC(s)$ by the definition of $I^\mathbf{C}_s$; in particular $v\notin\imC(s)$. This is a contradiction.

Next, we argue the inclusion $\subseteq$ in~\eqref{eq:UnionY}. Let $s\in Y(S,T)$.
Then observe that if $S\subseteq \imDMinus(s)$, there must exist  a minimal subset $Q'_s$ of $S$ such that $S\setminus Q'_s\subseteq \imC(s)$.
Then $I_{Q'_s}\subseteq J_s$, as all of $Q'_s$ must be used in the embedding of Stage~D$^-$. Hence, $s\in X(S\setminus Q,\emptyset, T\cup I_Q)$ with $Q=Q'_s$.
For the other inclusion $\supseteq$, assume that 
$s\in X(S\setminus Q,\emptyset, T\cup I_Q)$ for some $Q\subset S$, i.e.\ we have
$S\setminus Q\subseteq \imC(s)$, $T\subset J_s$ and $I_Q\subset J_s$. The last condition tells us that all of $Q$ is used in the embedding of Stage~D$^-$. Hence
$Q\subset \imDMinus(s)$ and thus
$S\subset \imDMinus(s)$, leading to 
$s\in Y(S,T)$. This proves~\eqref{eq:UnionY}.

We use~\eqref{eq:UnionY} and see that $\big|Y(S,T)\big|= \sum_{Q\subset S} \left| X(S\setminus Q,\emptyset, T\cup I_Q) \right|$. We use~\eqref{D:XST} to express the terms $\left| X(S\setminus Q,\emptyset, T\cup I_Q) \right|$. There is an error term in~\eqref{D:XST}. Summing this error term up over the $2^{|S|}\le 2^{4\LD}$ subsets $Q\subseteq S$, we see that the cumulative error is $e=\pm 600\LD 2^{6\LD}\LC\gamcore n$. In the rest, we focus on summing up the main terms in~\eqref{D:XST},
\begin{align*}
		\big|Y(S,T)\big|-e&=  \sum_{Q\subseteq S} |\cJ_0|(1-\deltnonspanning-10\sigmKJ^*)^{|S\setminus Q|}(2\sigmKJ^*)^{|T\cup I_Q|} 
		\\
	& =
	|\cJ_0| (1-\deltnonspanning-10\sigmKJ^*)^{|S|}(2\sigmKJ^*)^{|T|} \sum_{q=0}^{|S|} \binom{|S|}{q} \cdot \left( \frac{2\sigmKJ^*}{1-\deltnonspanning-10\sigmKJ^*} \right)^{q}\\
	\JUSTIFY{binomial theorem}& =
	|\cJ_0| (1-\deltnonspanning-10\sigmKJ^*)^{|S|}(2\sigmKJ^*)^{|T|} \left( 1+  \frac{2\sigmKJ^*}{1-\deltnonspanning-10\sigmKJ^*} \right)^{|S|} \\
	& = 
	|\cJ_0|(1-\deltnonspanning-8\sigmKJ^*)^{|S|}(2\sigmKJ^*)^{|T|}\;.
\end{align*}
This proves the claim.
\end{claimproof}
We can now verify~\ref{quasi:imagecaps}. Given sets $S_1,S_2,T$ as in~\ref{quasi:imagecaps}, each of size at most $2\LD$, we have
\begin{align*}
	& |\{s\in\cJ_0:S_1\cap\imDMinus(s)=\emptyset,S_2\subset\imDMinus(s),T\subset J_s\}| \\
	& = \left|Y(S_2,T) \setminus \bigcup_{x\in S_1} Y(S_2\cup \{x\},T)\right| 
	= \sum_{Q\subset S_1} (-1)^{|Q|}|Y(S_2\cup Q,T)|\;,
\end{align*}
by the inclusion-exclusion principle. We can now apply Claim~\ref{cl:defineY},
\begin{align*}
	& |\{s\in\cJ_0:S_1\cap\imDMinus(s)=\emptyset,S_2\subset\imDMinus(s),T\subset J_s\}| \\
	& = 
	|\cJ_0|(1-\deltnonspanning-8\sigmKJ^*)^{|S_2|}(2\sigmKJ^*)^{|T|} 
	\sum_{Q\subset S_1} (-1)^{|Q|}(1-\deltnonspanning-8\sigmKJ^*)^{|Q|}
	\pm 2^{2\LD}\cdot 600\LD 2^{6\LD}\LC\gamcore n\;.
\end{align*}
Using the binomial theorem, we continue as follows:
\begin{align*}
	\sum_{Q\subset S_1} (-1)^{|Q|}(1-\deltnonspanning-8\sigmKJ^*)^{|Q|}
	& =
	(1-\deltnonspanning-8\sigmKJ^*)^{|S_1|}
	\sum_{Q\subset S_1} (-1)^{|Q|}\left(\frac{1}{1-\deltnonspanning-8\sigmKJ^*}\right)^{|S_1|-|Q|} \\
	& =
	(1-\deltnonspanning-8\sigmKJ^*)^{|S_1|}
	\sum_{q=0}^{|S_1|} \binom{|S_1|}{q} (-1)^{q}\left(\frac{1}{1-\deltnonspanning-8\sigmKJ^*}\right)^{|S_1|-q} \\
	& =
	(1-\deltnonspanning-8\sigmKJ^*)^{|S_1|}\cdot\left(-1 + \frac{1}{1-\deltnonspanning-8\sigmKJ^*}\right)^{|S_1|}
	= (\deltnonspanning + 8\sigmKJ^*)^{|S_1|} \, ,
\end{align*}
hence giving that $\{s\in\cJ_0:S_1\cap\imDMinus(s)=\emptyset,S_2\subset\imDMinus(s),T\subset J_s\}$ has size
\[|\cJ_0|(\deltnonspanning+8\sigmKJ^*)^{|S_1|}(1-\deltnonspanning-8\sigmKJ^*)^{|S_2|}(2\sigmKJ^*)^{|T|}\pm 600\LD 4^{6\LD}\LC\gamcore n .\]

This guarantees~\ref{quasi:imagecaps} after Stage~D$^-$ with error term $600\LD 4^{6\LD} \deltnonspanning^{-4\LD}(\sigmKJ^*)^{-3\LD} \LC\gamcore$.

For~\ref{quasi:NumAnchors}, we only need to consider $\cJ_0$ since nothing has changed for $\cJ_1$ or $\cJ_2$. Given $v\in V(\HStageC)$, let $i$ be such that $v\in\{\boxminus_i,\boxplus_i\}$. After Stage~C, by~\ref{quasi:imagecaps} with $S_1=S_2=\emptyset$ and $T=\{i\}$, the number of $s\in\cJ_0$ such that $i\in I^\mathbf{C}_s$ is $(1\pm 200\LC\gamcore)(\deltnonspanning+10\sigmKJ^*)^2|\cJ_0|$. Now for a given $s$ from this collection, the probability that we pick $i\in J_s$ is $|J_s|/|I^\mathbf{C}_s|=(1\pm 100\gamcore)\tfrac{2\sigmKJ^*}{(\deltnonspanning+10\sigmKJ^*)^2}$. By Chernoff's inequality (Fact~\ref{fact:hypergeometricBasicProperties}) we conclude that the number of $s\in\cJ_0$ such that $i\in J_s$ is with probability at least $1-\exp(-n^{0.5})$ equal to $(1\pm 1000\LC\gamcore)(2\sigmKJ^*)|\cJ_0|$, as required for~\ref{quasi:NumAnchors}.

For~\ref{quasi:NumAnchorsNoIm}, we again need only consider $\cJ_0$ since nothing has changed since Stage~C for $\cJ_1$ and $\cJ_2$. Given $u,v\in V(\HStageC)$ with $u\neq v$ and $\{u,v\}$ not a terminal pair, let $i$ be such that $u\in\{\boxminus_i,\boxplus_i\}$ and let $j$ be such that $v\in\{\boxminus_j,\boxplus_j\}$. We want to estimate the number of $s\in\cJ_0$ such that $i\in J_s$ and both $v\notin\imC(s)$ and $j\notin J_s$ occur. Note that we can split these $s$ up according to whether $j\in I^\mathbf{C}_s$. If $j\in I^\mathbf{C}_s$, then automatically $v\notin\imC(s)$. 
By~\ref{quasi:imagecaps} after Stage~C with
$S_1=S_2=\emptyset$ and $T=\{i,j\}$ 
the number of $s\in\cJ_0$ such that $i,j\in I^\mathbf{C}_s$ is
equal to $(1\pm 100\LC\gamcore)\sum_{s\in\cJ_0}\big(2|I^\mathbf{C}_s|/n\big)^2$, and the probability, when we reveal $J_s$, that any given one of these has $i\in J_s$ and $j\notin J_s$ is 
$\frac{|J_s|(|I^\mathbf{C}_s|-|J_s|)}{|I^\mathbf{C}_s|(|I^\mathbf{C}_s|-1)} = \frac{|J_s|(|I^\mathbf{C}_s|-|J_s|)}{|I^\mathbf{C}_s|^2}\pm n^{-0.5}$. So by Chernoff's inequality (Fact~\ref{fact:hypergeometricBasicProperties}), with probability at least $1-\exp(-n^{0.5})$, there are
\[(1\pm 200\LC\gamcore)\sum_{s\in\cJ_0} 4|J_s|(|I^\mathbf{C}_s|-|J_s|)n^{-2}\]
$s\in\cJ_0$ such that $i\in J_s$ and $j\in I^\mathbf{C}_s\setminus J_s$.
Next, we consider those $s\in\cJ_0$ such that $i\in I^\mathbf{C}_s$, $v\notin\imC(s)$ and $j\notin I^\mathbf{C}_s$. Using~\ref{quasi:imagecaps} 
with $S_1=\{v\}$, $S_2=\emptyset$ and $T=\{i\}$
we can count those $s\in\cJ_0$ such that $i\in I^\mathbf{C}_s$ and $v\notin\imC(s)$ holds,
and using~\ref{quasi:imagecaps} 
with $S_1=S_2=\emptyset$ and $T=\{i,j\}$
we can count those $s\in\cJ_0$ such that $i,j\in I^\mathbf{C}_s$, and therefore $v\notin\imC(s)$ (since
$j\in I^\mathbf{C}_s$ implies $v\notin\imC(s)$). Subtracting both
results we see that there are
\[\sum_{s\in\cJ_0}\tfrac{2|I^\mathbf{C}_s|}{n}\cdot \big(\tfrac{n-|U^\mathbf{C}_s|}{n}-\tfrac{2|I^\mathbf{C}_s|}{n}\big)\pm 200\LC\gamcore n\]
indices $s\in\cJ_0$ such that $i\in I^\mathbf{C}_s$, $v\notin\imC(s)$ and $j\notin I^\mathbf{C}_s$.
Of these, when we reveal $J_s$, the probability that any one such $s$ has $i\in J_s$ is $\tfrac{|J_s|}{|I^\mathbf{C}_s|}$, and so by Chernoff's inequality (Fact~\ref{fact:hypergeometricBasicProperties}), with probability at least $1-\exp(-n^{0.5})$, there are
\[\sum_{s\in\cJ_0}\tfrac{2|J_s|}{n}\cdot \big(\tfrac{n-|U^\mathbf{C}_s|}{n}-\tfrac{2|I^\mathbf{C}_s|}{n}\big)\pm 300\LC\gamcore n\]
such $s$. Summing up, we have
\begin{align*}
 (1\pm 200\LC\gamcore)&\sum_{s\in\cJ_0} 4|J_s|(|I^\mathbf{C}_s|-|J_s|)n^{-2}+\sum_{s\in\cJ_0}\tfrac{2|J_s|}{n}\cdot \big(\tfrac{n-|U^\mathbf{C}_s|}{n}-\tfrac{2|I^\mathbf{C}_s|}{n}\big)\pm 300\LC\gamcore n\\
&=\sum_{s\in\cJ_0}\tfrac{2|J_s|}{n}\cdot\tfrac{n-|U^\mathbf{C}_s|-2|J_s|}{n} \pm 500\LC\gamcore n\,.
\end{align*}
Note that $|U^\mathbf{D-}_s|$ is precisely $|U^\mathbf{C}_s|+2|J_s|$, so that this formula guarantees~\ref{quasi:NumAnchorsNoIm} after Stage~D$^-$ with error $600 \LD 4^{6\LD} \deltnonspanning^{-4\LD}(\sigmKJ^*)^{-3\LD} \LC \gamcore$. 

Finally, for~\ref{D:reroot:quasioldA}, we use a similar argument as for~\ref{quasi:NumAnchorsNoIm}. Given $u,v$, if $u\in A^\mathbf{C}_s$ and $v\notin U^\mathbf{D-}_s$, then we need $v\notin\imC(s)$ and $v\notin\{\boxminus_j,\boxplus_j\}$ for all $j\in J_s$. From~\ref{quasi:NumAnchorsNoIm} after Stage~C, we have
\[
  \big|\{s\in\cJ_0:u\in A^\mathbf{C}_s,v\notin\imC(s)\}\big|=(1\pm100\LC\gamcore)|\cJ_0|\cdot 2\sigmKJ^*(\deltnonspanning+10\sigmKJ^*)\,.
\]
By Lemma~\ref{lem:StageCNew}\ref{C:five}, letting $j$ be the index such that $v\in\{\boxminus_j,\boxplus_j\}$, we have
\[\big|\{s\in\cJ_0:u\in A^\mathbf{C}_s,j\in I^\mathbf{C}_s\}\big|=(1\pm 100\LC\gamcore)\sum_{s\in\cJ_0}\tfrac{2|A^\mathbf{C}_s||I^\mathbf{C}_s|}{n^2}=(1\pm 200\LC\gamcore)|\cJ_0|\cdot 2\sigmKJ^*(\deltnonspanning+10\sigmKJ^*)^2\,.\]
Note that the probability that a given $s$ in this last set has $j\in J_s$ is $\tfrac{|J_s|}{|I^\mathbf{C}_s|}$, so by \eqref{eq:sizeICs} and Chernoff's inequality, with probability at least $1-\exp(-n^{0.5})$ we have
\[\big|\{s\in\cJ_0:u\in A^\mathbf{C}_s,j\in J_s\}\big|=(1\pm 300\LC\gamcore)|\cJ_0|\cdot (2\sigmKJ^*)^2\,.\]
Now, if $v\in\{\boxminus_j,\boxplus_j\}$ with $j\in J_s$, then immediately $v\notin \imC(s)$, and hence putting everything together, we have
\begin{align*}
	\big|\{s\in\cJ_0:u\in A^\mathbf{C}_s,v\notin U^\mathbf{D-}_s\}\big|
	& = 
	\big|\{s\in\cJ_0:u\in A^\mathbf{C}_s,v\notin \imC(s)\}\big|
	-
	\big|\{s\in\cJ_0:u\in A^\mathbf{C}_s,j\in J_s\}\big| \\
	& =(1\pm 1000\LC\gamcore)|\cJ_0|\cdot 2\sigmKJ^*(\deltnonspanning+10\sigmKJ^*-2\sigmKJ^*) 
\end{align*}
as required.
\end{proof}

We are now in a position to provide the proof of Lemma~\ref{lem:StageDNew}. As we indicated at the beginning of the section, we first use Lemma~\ref{lem:reroot-assign} to get us to Stage~D$^-$. Next, we would like to use Lemma~\ref{lem:pathpack} to pack $P^\mathsf{left}$ and $P^\mathsf{right}$. Lemma~\ref{lem:pathpack} has a number of hypotheses, and checking these occupies the bulk of the proof.

\begin{proof}[Proof of Lemma~\ref{lem:StageDNew}]\label{pageref:proofOfLem:StageDNewProperties}
	We set $\nu=\tfrac{1}{1000}2^{-4\LD}\deltnonspanning^{4\LD}(\sigmKJ^*)^{5\LD}$. Let $C'$ be the constant returned by Lemma~\ref{lem:pathpack} for parameters $\nu$ and $2\LD$. Set $C'':=2^{10\LD}C'$. Set
	\[\gamma'=1200\LD\cdot 4^{6\LD}\deltnonspanning^{-4\LD}(\sigmKJ^*)^{-3\LD}\LC\gamcore\;.\]
	The smallest constant in the hierarchy~\eqref{eq:CONSTANTS} involved in the definition of the input pair $(\nu,2\LD)$ for Lemma~\ref{lem:pathpack} is $\frac{1}{\LD}$. That means that the entire chain in~\eqref{eq:CONSTANTS} following the `$\frac{1}{\LD}\gg\ldots$' part can be made much smaller than $(C')^{-1}$. In particular, we have $(C'')^{-1}\gg\gamma',\gamAnchor$.
	Lemma~\ref{lem:reroot-assign} guarantees that after Stage~D$^-$ we have a $(\gamma'/2,2\LD,\dStageC,\dStageC)$-quasirandom setup.
	
	Let $H=\HStageC$, which has $2\lfloor\tfrac{n}{2}\rfloor$ vertices. We will therefore be applying Lemma~\ref{lem:pathpack} with a number $2\lfloor\tfrac{n}{2}\rfloor$ of vertices (which is either $n$ or $n-1$). This change (between the $n$ here and $2\lfloor\tfrac{n}{2}\rfloor$ for Lemma~\ref{lem:pathpack}) introduces a relative $1\pm O(n^{-1})$ error in various calculations in the following, which in all cases is absorbed into the error terms.
	
	We let $s^*=|\cJ_0|$. By~\eqref{eq:sizeJ0J1}, we have $s^*=\sigmJnula n\ge \nu n$. Temporarily abusing notation, we relabel the graphs $G_s$ such that $\cJ_0=[s^*]$. For each $s\in\cJ_0$, we do the following.
	Let $F_s$ denote the path-forest whose components are $\leftpath_4(P)$ and $\rightpath_4(P)$ for each $P\in\SpecPaths^*_s$.
	We let $U_s=U^\mathbf{D-}_s=\imDMinus(s)$, and we let $\phi_s$ be the restriction of $\phi^{\mathbf{D-}}_s$ to the leaves of $F_s$.
	Moreover, we let $A_s$ be the images of the leaves of $F_s$. Note that for each $s$, the path-forest $F_s$ is \emph{not} $F^\mathbf{D-}_s$, and the set $A_s$ is \emph{not} the set $A^\mathbf{D-}_s$ from the quasirandom setup after Stage~D$^-$; rather $A_s=A^\mathbf{C}_s\cup A^\mathbf{D-}_s$, where $A^\mathbf{C}_s$ is the set of anchors after Stage~C, and we recall $A^\mathbf{D-}_s=\{\boxminus_i,\boxplus_i:i\in J_s\}$. By construction $A^\mathbf{C}_s$ and $A^\mathbf{D-}_s$ are disjoint.
	
	For each vertex $x\in \phi_s^{-1}(A_s)$, we let $\chi_s(x)$ be whichever of $\Vmin$ and $\Vplus$ contains $\phi_s(x)$. For each vertex $y\in V(F_s)\setminus \phi_s^{-1}(A_s)$, we choose $\chi_s(y)\in\{\Vmin,\Vplus\}$ uniformly at random.
	
	We now want to use Lemma~\ref{lem:pathpack} to finish Stage~D. We shall check that with high probability (with respect to the random choice of the maps $\{\chi_s\}_{s\in \cJ_0}$), all the assumptions of Lemma~\ref{lem:pathpack} with input parameters $\nu$, $2\LD$, and $\gamma=(2\LD+3)\gamma'$ are satisfied. We structure this as follows: The basic requirements of~\ref{PA:BasicSetting} hold by the above constructions. In Claim~\ref{UMDS:BlockQuasiDiet}, we justify assumptions~\ref{PA:BlockQuasirandom} and~\ref{PA:BlockDiet}. In Claim~\ref{UMDS:anchor}, we justify assumption~\ref{PA:AnchorDistribution}. In Claim~\ref{UMDS:Ineq}, we justify the inequalities of~\ref{PA:Ineq}. In Claim~\ref{UMDS:pairdistribution}, we justify assumption~\ref{PA:PairDistributionProperty}. Finally, in Claim~\ref{UMDS:densities}, we verify the leftover densities assumption~\ref{PA:LeftOverDensities}.
	
	\begin{claimInLemmaD}\label{UMDS:BlockQuasiDiet}
		The graph $H$ and the sets $U_s$ satisfy the assumptions~\ref{PA:BlockQuasirandom} and~\ref{PA:BlockDiet}.
	\end{claimInLemmaD}
	\begin{claimproof}[Proof of Claim~\ref{UMDS:BlockQuasiDiet}]
		Recall that by Lemma~\ref{lem:reroot-assign}\ref{D:reroot:quasisetup} the setup after Stage~D$^-$ is $(\gamma'/2,2\LD,\dStageC,\dStageC)$-quasirandom; we use~\ref{quasi:indexquasi} in particular. We stress that it is this setup, rather than the $(100\LC\gamcore,\LC,\dStageC,\dStageC)$-quasirandom setup after Stage~C given by Lemma~\ref{lem:StageCNew}\ref{C:four}, that we need here: the sets $U_s=U^\mathbf{D-}_s=\imDMinus(s)$ for which we want the block-diet condition are the used sets of the setup after Stage~D$^-$, and for $s\in\cJ_0$ these are strictly larger than the used sets $U^\mathbf{C}_s=\imC(s)$ after Stage~C, the difference being the newly assigned anchors $\{\boxminus_i,\boxplus_i\,:\,i\in J_s\}$. Applying Lemma~\ref{lem:indexquasi-implies-blockquasi} with $L=2\LD$ and $\gamma=\gamma'/2$, we conclude that $H$ is $\big((2\LD+3)\tfrac{\gamma'}{2},2\LD\big)$-block-quasirandom, and that for each $s\in\cJ_0$ the pair $(H,U_s)$ is $\big((2\LD+3)\tfrac{\gamma'}{2},2\LD\big)$-block-diet. In particular $H$ is $((2\LD+3)\gamma',2\LD)$-block-quasirandom, and $(H,U_s)$ is $((2\LD+3)\gamma',2\LD)$-block-diet for each $s\in\cJ_0$, as required for the input $\gamma=(2\LD+3)\gamma'$ of Lemma~\ref{lem:pathpack}. Finally, letting $d_{\boxminus\boxminus},d_{\boxplus\boxplus}$ denote the density of $H[\Vmin]$ and $H[\Vplus]$ respectively, and $d_{\boxminus\boxplus}$ denote the bipartite density of $H[\Vmin,\Vplus]$, the density conclusion of Lemma~\ref{lem:indexquasi-implies-blockquasi} gives that each of these three parameters is $(1\pm\gamma')\dStageC$, and in particular $(1\pm(2\LD+3)\gamma')\dStageC$. 
	\end{claimproof}
	
	\begin{claimInLemmaD}\label{UMDS:anchor}
		With high probability, for each $s\in\cJ_0$, the path-forest $F_s$ has the $10\gamma'$-anchor distribution property with respect to $H$, justifying~\ref{PA:AnchorDistribution}. 
	\end{claimInLemmaD}
	\begin{claimproof}[Proof of Claim~\ref{UMDS:anchor}]
		Given $a,b,c\in\{\boxminus,\boxplus\}$ and $v\in V_b\setminus U_s$, we have 
		\begin{align*}
			\big|A_s\cap \NBH_H(v)\cap V_a\big|&=\big|A^\mathbf{C}_s\cap \NBH_H(v)\cap V_a\big|+\big|A^\mathbf{D-}_s\cap \NBH_H(v)\cap V_a\big|
			\\
			\JUSTIFY{\ref{quasi:distributionOfAnchors} after Stage~C and after Stage~D$^-$}&=(1\pm 100\LC\gamcore)\dStageC \big|A^\mathbf{C}_s\cap V_a\big|+(1\pm \gamma'/2)\dStageC\big|A^\mathbf{D-}_s\cap V_a\big|\\
			\JUSTIFY{\ref{quasi:anchorsets} after Stage~C and after Stage~D$^-$}&=(1\pm \gamma')\cdot \frac{\dStageC}2\cdot (|A^\mathbf{C}_s|+|A^\mathbf{D-}_s|)\;.
		\end{align*}
		As each path in $\SpecPaths_s^*$ contributes~$4$ to $|A^\mathbf{C}_s|+|A^\mathbf{D-}_s|$, we conclude by \eqref{eq:sizeSpecPathStar} that
		\[
		\big|A_s\cap \NBH_H(v)\cap V_a\big|=(1\pm\gamma')\dStageC\cdot 2\sigmKJ^* n\,.
		\]
		Given a path $P$ of $F_s$ and a leaf $x$ of that path anchored in $\NBH_H(v)\cap V_a$, the probability that the neighbour $y$ and second neighbour $z$ of $x$ are assigned to sets $\chi_s(y)=V_b$ and $\chi_s(z)=V_c$ is $\tfrac14$. Since the choice of $\chi_s$ on $P$ affects at most the two leaves of $P$, we conclude by Lemma~\ref{lem:McDiarmidInequ} that with probability at least $1-\exp(-n^{0.1})$ we have $(1\pm2\gamma')\dStageC\cdot\tfrac12\sigmKJ^* n$ leaves $x$ of $F_s$ embedded to $\NBH_H(v)\cap V_a$, such that their neighbour and second neighbour are assigned to $V_b$ and $V_c$ respectively by $\chi_s$.
		Similarly, by~\ref{quasi:anchorsets} applied both after Stage~C to $A^\mathbf{C}_s$ and after Stage~D$^-$ to $A^\mathbf{D-}_s$, the number of leaves of $F_s$ anchored in $V_a$ in total is $(1\pm\gamma')\cdot 2\sigmKJ^* n$. By a similar application of Lemma~\ref{lem:McDiarmidInequ}, with probability at least $1-\exp(-n^{0.1})$, the number of these leaves whose neighbour is assigned by $\chi_s$ to $V_b$ and second neighbour is assigned to $V_c$ is $(1\pm 2\gamma')\cdot\tfrac12\sigmKJ^* n$. Taking the union bound over choices of $a$, $b$, $c$, $s$ and $v$, with high probability we have the $10\gamma'$-anchor distribution property for each $s\in\cJ_0$.
	\end{claimproof}
	
	\begin{claimInLemmaD}\label{UMDS:Ineq}
		With high probability, the assignments satisfy the inequalities of~\ref{PA:Ineq}.
	\end{claimInLemmaD}
	\begin{claimproof}[Proof of Claim~\ref{UMDS:Ineq}]
		Observe that for each $a\in\{\boxminus,\boxplus\}$, by~\ref{quasi:indexquasi} applied with $S_1=S_2=T_3=\emptyset$ and either $T_1=\{s\},T_2=\emptyset$ or vice versa (where our setting $T_3=\emptyset$ allows us to use~\eqref{eq:U4U5}), we have
		$|V_a\setminus U_s|=(1\pm\gamma')\cdot\tfrac12(\deltnonspanning+8\sigmKJ^*)n$, while $\big|V(F_s)\setminus \phi_s^{-1}(A_s)\big|=6\sigmKJ^* n$. Even if $\chi_s$ assigned all these vertices to $V_a$, by choice of $\nu$ we have that $\big|V_a\setminus U_s\big|-\big|\{x\in V(F_s)\setminus \phi_s^{-1}(A_s)\,:\,\chi_s(x)=V_a\}\big|\ge\nu n$. Furthermore, each path of $F_s$ has 5 vertices, of which the middle 3 are assigned uniformly at random to $\Vmin$ or $\Vplus$ by $\chi_s$. By the Chernoff bound and the union bound over the at most $n$ choices of $s$, in each $F_s$ we have $\big|\chi_s^{-1}(\{V_a\})\setminus \phi_s^{-1}(A_s)\big|\ge \sigmKJ^* n\ge\nu n$ with probability at least $1-\exp(-n^{0.1})$. 
		Similarly, in each path of $F_s$, the two internal edges each have probability at least $\tfrac14$ of being assigned by $\chi_s$ to any of within $\Vmin$, within $\Vplus$, or between $\Vmin$ and $\Vplus$, and by Lemma~\ref{lem:McDiarmidInequ}, the total number of internal edges of $F_s$ assigned by $\chi_s$ to each of within $\Vmin$, within $\Vplus$, and between $\Vmin$ and $\Vplus$, is at least $\tfrac14\sigmKJ^* n\ge\nu n$ with probability at least $1-\exp(-n^{0.1})$.
	\end{claimproof}
	
\begin{claimInLemmaD}\label{UMDS:pairdistribution}
	With high probability, the random choice of the $\chi_s$ gives the $7\gamma'$-pair distribution property, fulfilling~\ref{PA:PairDistributionProperty}. 
\end{claimInLemmaD}
\begin{claimproof}[Proof of Claim~\ref{UMDS:pairdistribution}]
	To that end, fix distinct vertices $u\in V_a$ and $v\in V_b$, where $a,b\in\{\boxminus,\boxplus\}$. Suppose that we have $uv\in E(\HStageC)$, and in particular by Lemma~\ref{lem:StageCNew}\ref{C:three} we do not have $\{u,v\}=\{\boxminus_i,\boxplus_i\}$ for any $i$. We want to know the various $w_{uv;s}$ as $s$ ranges over $\cJ_0$, as in Definition~\ref{def:pairdistprop}. Recall that there are three different possibilities for $w_{uv;s}>0$. If $u,v\notin U_s$, then we claim it is likely that we have
	\[w_{uv;s}=\frac{(1\pm\gamma')2\sigmKJ^* n}{(1\pm\gamma')^2\tfrac14(\deltnonspanning+8\sigmKJ^*)^2n^2}=(1\pm 4\gamma')\frac{8\sigmKJ^*}{(\deltnonspanning+8\sigmKJ^*)^2 n}\,.\]
	The denominator uses estimates on $|V_a\setminus U_s|$ from~\ref{quasi:indexquasi} as above. For the numerator we observe that $F_s$ has $4\sigmKJ^* n$ internal edges (and so $8\sigmKJ^* n$ pairs span an internal edge), each of which has probability $\tfrac14$ of being assigned by $\chi_s$ to respectively $V_a$ and $V_b$. The random choice of $\chi_s$ on each component of $F_s$ affects the sum by at most $4$. Hence by Lemma~\ref{lem:McDiarmidInequ}, with probability at least $1-\exp(-n^{0.1})$ we obtain the numerator. It remains to estimate the number of $s\in\cJ_0$ such that $u,v\notin U_s$. By~\ref{quasi:imagecaps}, with $S_1=\{u,v\}$ and $S_2=T=\emptyset$, we see that this number is $(1\pm \gamma')(\deltnonspanning+8\sigmKJ^*)^2|\cJ_0|$. Thus the contribution to $\sum_{s\in\cJ_0}w_{uv;s}$ from $s$ falling into this category is
	\[(1\pm 4\gamma')\frac{8\sigmKJ^*}{(\deltnonspanning+8\sigmKJ^*)^2n}\cdot (1\pm\gamma')(\deltnonspanning+8\sigmKJ^*)^2|\cJ_0|=(1\pm 6\gamma')\frac{8\sigmKJ^*|\cJ_0|}{n}\,.\]
	
	Next, if $u\in A_s$ and $v\notin U_s$, then we can have $w_{uv;s}>0$, and if we do, we have
	\[w_{uv;s}=\tfrac{1}{(1\pm\gamma')(\deltnonspanning+8\sigmKJ^*)n/2}\]
	by our estimate on $|V_a\setminus U_s|$ from~\ref{quasi:indexquasi}.  This occurs when the vertex $x$ of $F_s$ anchored to $u$ has neighbour $y$ such that $\chi_s(y)=V_b$. There are two ways this can happen: either $u\in A^\mathbf{C}_s$, $v\notin U_s$, and $\chi_s(v)=V_b$, or $u\in A^\mathbf{D-}_s$, $v\notin U_s$, and $\chi_s(v)=V_b$. For the first of these, by~\ref{D:reroot:quasioldA} we see that there are
	$(1\pm 1000\gamcore)|\cJ_0|2\sigmKJ^*(\deltnonspanning+8\sigmKJ^*)$
	choices of $s$ for which $u\in A^\mathbf{C}_s$ and $v\notin U_s$, and of these by Chernoff's inequality with probability at least $1-\exp(-n^{0.5})$ we see
	\[(1\pm 2000\gamcore)|\cJ_0|\sigmKJ^*(\deltnonspanning+8\sigmKJ^*)\]
	choices of $s$ for which all three conditions hold. 
	For the second, we use~\ref{quasi:NumAnchorsNoIm} with the setup after Stage~D$^-$ (which has error parameter $\gamma'/2$) and Chernoff's inequality to obtain the bound $(1\pm\gamma')|\cJ_0|\sigmKJ^*(\deltnonspanning+8\sigmKJ^*)$. Putting these together, the contribution to $\sum_{s\in\cJ_0}w_{uv;s}$ from $s$ falling into this category is
	\[\tfrac{1}{(1\pm\gamma')(\deltnonspanning+8\sigmKJ^*)n/2}\cdot \big(2\pm 2\gamma'\big)|\cJ_0|\sigmKJ^*(\deltnonspanning+8\sigmKJ^*)=(1\pm 4\gamma')\tfrac{4\sigmKJ^*|\cJ_0|}{n}\,.  \]
	
	Finally, we can have $w_{uv;s}>0$ if $v\in A_s$ and $u\notin U_s$. This is symmetric to the previous case, and we obtain the same bounds. We conclude that
	\[\sum_{s\in\cJ_0}w_{uv;s}=(1\pm 6\gamma')\frac{16\sigmKJ^*|\cJ_0|}{n}\,,\]
	which implies the $7\gamma'$-pair distribution property as required by Definition~\ref{def:pairdistprop}. Indeed, each $F_s$ has $8\sigmKJ^* n$ edges, hence $16\sigmKJ^* n$ pairs $(x,y)$ of vertices such that $xy\in E(F_s)$, of which by the random choice of the $\chi_s$ and random anchor distribution $(1\pm\gamcore)4\sigmKJ^* n$ pairs $(x,y)$ have $\chi_s(x)=V_a$ and $\chi_s(y)=V_b$. Comparing our sum with this target count introduces an additional $(1\pm\gamcore)$ error factor, resulting in a total relative error well covered by the $7\gamma'$ bound.
\end{claimproof}
	
	\begin{claimInLemmaD}\label{UMDS:densities}
		With high probability, the leftover densities satisfy $d'_{ab} \ge \nu$ for each $a,b\in\{\boxminus,\boxplus\}$, fulfilling~\ref{PA:LeftOverDensities}.
	\end{claimInLemmaD}
	\begin{claimproof}[Proof of Claim~\ref{UMDS:densities}]
		We define $d'_{ab}$ for $a,b\in\{\boxminus,\boxplus\}$ as in Lemma~\ref{lem:pathpack} to be the densities of the graph obtained by successfully packing all the path-forests $F_s$ into $\HStageC$ according to $\chi_s$. We now aim to show that each of these quantities is
	\begin{equation}\label{eq:whatthequantitiesreallyare}
		(1\pm(2\LD+4)\gamma')\big(\dStageC-\tfrac{16\sigmKJ^* |\cJ_0|}{n}\big)\,.
	\end{equation}
To that end, we claim that with high probability, for each $s\in\cJ_0$, the number of edges $uv\in E(F_s)$ such that $\chi_s(u),\chi_s(v)=\Vmin$ is $(1\pm 1000\gamcore)\cdot 2\sigmKJ^*n$. The choice of $\chi_s$ on different components of $F_s$ is independent, and the effect of the choice on any given component is to change this count by at most $4$, so by Lemma~\ref{lem:McDiarmidInequ} the probability that the actual number of edges $uv\in E(F_s)$ such that $\chi_s(u),\chi_s(v)=\Vmin$ differs from its expectation by $n^{0.7}$ or more is at most $\exp(-n^{0.1})$. Now the probability that an edge of a path $P$ which is not adjacent to a leaf of $P$ has both vertices assigned to $\Vmin$ is exactly $\tfrac14$. If $uv$ is an edge of $P$ such that $u$ is a leaf vertex (so $\chi_s(u)$ is fixed) then we have $\chi_s(u),\chi_s(v)=\Vmin$ with probability either $\tfrac12$ (if $\chi_s(u)=\Vmin)$ or $0$. By~\ref{quasi:anchorsets}, applied after Stage~C to $A^\mathbf{C}_s$, we see
		\[|A_s\cap\Vmin|=(1\pm100\gamcore)\sigmKJ^* n+\sigmKJ^* n\,,\]
		and putting these together, the expected number of edges of $F_s$ both of whose ends are assigned to $\Vmin$ is $(1\pm100\gamcore)\cdot 2\sigmKJ^* n$, from which the claim follows. 
		By the same argument, the same bound holds replacing $\Vmin$ with $\Vplus$, and taking the union bound, with high probability both bounds hold for each $s\in\cJ_0$. Hence, taking into account that the initial block densities of $\HStageC$ have error $(2\LD+3)\gamma'$ by our earlier claim, we obtain
		\[d'_{ab}=(1\pm(2\LD+4)\gamma')\dStageD\,,\quad\text{where}\quad \dStageD=\dStageC-\tfrac{2\sigmKJ^* n |\cJ_0|}{\binom{\lfloor n/2 \rfloor}{2}}\]
		for each $a,b\in\{\boxminus,\boxplus\}$ as required for~\eqref{eq:whatthequantitiesreallyare}.
		
Now it remains to show that the quantity in~\eqref{eq:whatthequantitiesreallyare} is at least $\nu$. Indeed, substituting $|\cJ_0| = \sigmJnula n$ (recall~\eqref{eq:sizeJ0J1}) and using the hierarchy of constants~\eqref{eq:CONSTANTS} and our discussion in Section~\ref{sssection:EvolutionOfDensities} (which implies $\nu \ll 16\sigmKJ^*\sigmJnula \ll \dStageC$), we have
\[
(1\pm(2\LD+4)\gamma')\big(\dStageC-\tfrac{16\sigmKJ^* |\cJ_0|}{n}\big)
\ge (1-(2\LD+4)\gamma')\big(\dStageC-16\sigmKJ^*\sigmJnula\big)
\ge \frac{\dStageC}{2}
\ge \nu\;,
\]	
as was needed.
\end{claimproof}
	
	Suppose that the random choice of the $\chi_s$ is such that all the above high probability events hold. Then by Lemma~\ref{lem:pathpack}, with input $\nu,2\LD,\gamma=(2\LD+3)\gamma'$, when Algorithm~\ref{alg:path-pack} is executed, with probability at least $1-\exp\big(-n^{0.2}\big)$ we obtain for each $s\in\cJ_0$ an embedding $\phi'_s$ of $F_s$ into $\HStageC\big[(V(H)\setminus U_s)\cup\phi_s(A_s)\big]$ which extends $\phi_s$, such that $\phi'_s(x)\in\chi_s(x)$ for each $x\in V(F_s)$, and such that for each $uv\in E(H)$ there is at most one $s$ such that $\phi'_s$ uses the edge $uv$. Suppose that this likely event occurs; we will in addition list polynomially many further likely events, which we presume all occur, and fix such an outcome.
	
	We now define the setup after Stage~D. We let $\HStageD$ be the graph $H'$ returned by Lemma~\ref{lem:pathpack}. For $s\in \cJ_1\cup\cJ_2$, we set $\phi^\mathbf{D}_s=\phi^\mathbf{C}_s$, with the same anchor and used sets and the same path-forests to embed. For $s\in\cJ_0$, we let $\phi^\mathbf{D}_s$ be the union of $\phi^\mathbf{D-}_s$ and the above obtained $\phi'_s$ (these overlap on the set $A_s$, but they agree on these vertices). We write $\imD(s)=\im(\phi^\mathbf{D}_s)$.\index{$\imD(s)$}
	We let $F^\mathbf{D}_s$ be the path-forest with paths $\SpecShortPaths_s$, with anchor set $A^\mathbf{D}_s$ the image under $\phi^\mathbf{D}_s$ of its leaves and $U^\mathbf{D}_s=\imD(s)$. We let $I^\mathbf{D}_s=\{i\in[\lfloor n/2\rfloor]:\boxminus_i,\boxplus_i\in A^\mathbf{D}_s\}$. 
	
	We now argue that this is a $(\gamAnchor,\LD,\dStageD,\dStageD)$-quasirandom setup, as required for Lemma~\ref{lem:StageDNew}. In each case, we use $\gamma'\ll\gamAnchor$ provided by~\eqref{eq:CONSTANTS}. To ensure clarity, we structure the verification into separate claims: Claim~\ref{UMDS:Quasi1} verifies~\ref{quasi:indexquasi}, Claim~\ref{UMDS:Quasi234} verifies~\ref{quasi:anchorsets}, \ref{quasi:PathsWithOneVertex}, and \ref{quasi:termVtxNbs}, Claim~\ref{UMDS:Quasi5} verifies~\ref{quasi:distributionOfAnchors}, Claim~\ref{UMDS:Quasi6} verifies~\ref{quasi:imagecaps}, and finally Claim~\ref{UMDS:Quasi78} verifies~\ref{quasi:NumAnchors} and \ref{quasi:NumAnchorsNoIm}.
	
	\begin{claimInLemmaD}\label{UMDS:Quasi1}
		The setup satisfies index-quasirandomness as required by~\ref{quasi:indexquasi}.
	\end{claimInLemmaD}
	\begin{claimproof}[Proof of Claim~\ref{UMDS:Quasi1}]
		Fix $S_1,S_2,T_1,T_2,T_3$ as in Definition~\ref{def:index-quasirandom}, each of size at most $\LD$, with families of sets $(U^\mathbf{D-}_s)_{s\in\cJ}$ and $(J_s)_{s\in\cJ_0}$, and let $X=\mathbb{U}_{\HStageC}(S_1,S_2,T_1,T_2,T_3)$. Then we have $|X|\ge\nu n$: since the setup after Stage~D$^-$ is $(\gamma'/2,2\LD,\dStageC,\dStageC)$-quasirandom and $|S_1|+|S_2|\le2\LD$, $|T_1\cup T_2|\le2\LD$, $|T_3|\le\LD$, Definition~\ref{def:index-quasirandom} gives
		\begin{equation}\label{eq:XgeqNuN}|X|=(1\pm\gamma'/2)\dStageC^{|S_1|+|S_2|}\cdot\prod_{\ell\in T_1\cup T_2}\Big(1-\tfrac{|U^\mathbf{D-}_\ell|}{n}\Big)\cdot\prod_{\ell\in T_3}\tfrac{|J_\ell|}{n/2}\cdot\tfrac n2\ge\tfrac12\dStageC^{2\LD}\deltnonspanning^{2\LD}(2\sigmKJ^*)^{\LD}\cdot\tfrac n2\,,\end{equation}
		where we used that every $G_s$, $s\in\cJ$, has exactly $(1-\deltnonspanning)n$ vertices (so $|U^\mathbf{D-}_\ell|\le(1-\deltnonspanning)n$, giving $1-|U^\mathbf{D-}_\ell|/n\ge\deltnonspanning$, for each $\ell$) and $|J_\ell|=\sigmKJ^* n$ by construction, and that a product of at most $2\LD$ (resp.\ $\LD$) factors each in $(0,1)$ is minimised by taking all $2\LD$ (resp.\ $\LD$) factors at their smallest value. Since $\dStageC\ge13(\sigmKJ^*)^2$ for $n$ sufficiently large (Section~\ref{sssection:EvolutionOfDensities}), this is at least $\tfrac{n}{4}\cdot13^{2\LD}2^{\LD}\deltnonspanning^{2\LD}(\sigmKJ^*)^{5\LD}\ge\tfrac{n}{1000}2^{-4\LD}\deltnonspanning^{4\LD}(\sigmKJ^*)^{5\LD}=\nu n$, as claimed. Now~\ref{pathpack:megaqr} is stated in terms of the true densities $d_{ab}:=d(\HStageC[V_a])$ (resp.\ $d(\HStageC[\Vmin,\Vplus])$) and $d'_{ab}$ from~\ref{PA:LeftOverDensities}, not directly in terms of $\dStageC,\dStageD$. By Claims~\ref{UMDS:BlockQuasiDiet} and~\ref{UMDS:densities} we have $d_{ab}=(1\pm(2\LD+3)\gamma')\dStageC$ and $d'_{ab}=(1\pm(2\LD+4)\gamma')\dStageD$ for each $a,b$, so $\tfrac{d'_{ab}}{d_{ab}}=(1\pm(8\LD+13)\gamma')\tfrac{\dStageD}{\dStageC}$; raising this to a power of at most $|S_1|+|S_2|\le2\LD$ multiplies the relative error by at most $2\LD$, i.e.\ contributes at most $(32\LD^2+52\LD)\gamma'$, which is dominated by (and hence absorbed into) $C'(2\LD+3)\gamma'$ since $C'$ depends only on $\nu,2\LD$ and is larger than any fixed polynomial in $\LD$. The number of tuples $(S_1,S_2,T_1,T_2,T_3)$ of size at most $\LD$ each is at most $n^{5\LD}$, and the failure probability $\exp(-n^{0.2})$ in~\ref{pathpack:megaqr} is uniform over the choice of $S_1,S_2,T_1,T_2,X$ (it does not depend on which are chosen), so a union bound over all $n^{5\LD}=o(\exp(n^{0.2}))$ such tuples shows that with probability $1-o(1)$ the conclusion of~\ref{pathpack:megaqr} holds simultaneously for every such tuple. So by~\ref{pathpack:megaqr} it is likely that we have (recall that $\mathbb{U}'$ is the variant of $\mathbb{U}$ introduced before Lemma~\ref{lem:pathpack}, formed with the enlarged used sets $U^\mathbf{D}_t=U^\mathbf{D-}_t\cup\im\phi'_t$, so that here $\mathbb{U}'_{\HStageD}$ is exactly $\mathbb{U}$ taken with respect to the setup after Stage~D)
		\begin{equation}\label{eq:XcapUHD}\big|X\cap \mathbb{U}'_{\HStageD}(S_1,S_2,T_1\cap\cJ_0,T_2\cap\cJ_0)\big|=(1\pm 2C'(2\LD+3)\gamma')\big(\tfrac{\dStageD}{\dStageC}\big)^{|S_1|+|S_2|}\prod_{t\in T_1}\frac{|\Vmin\setminus U^\mathbf{D}_t|}{|\Vmin\setminus U^\mathbf{D-}_t|}\prod_{t\in T_2}\frac{|\Vplus\setminus U^\mathbf{D}_t|}{|\Vplus\setminus U^\mathbf{D-}_t|}|X|\,.\end{equation}
		Using~\ref{quasi:indexquasi} as before (with the setup after Stage~D$^-$ being $(\gamma'/2,2\LD,\dStageC,\dStageC)$-quasirandom), we see that for $a\in \{\boxplus,\boxminus\}$, we have
		$|V_a\setminus U^\mathbf{D-}_t|=(1\pm \tfrac{\gamma'}2)\frac{|V\setminus U^\mathbf{D-}_t|}{2}$. Writing $N_a(t)$ for the number of vertices of $V(F_t)\setminus\phi_t^{-1}(A_t)$ assigned to $V_a$ by $\chi_t$, we have $|V_a\setminus U^\mathbf{D}_t|=|V_a\setminus U^\mathbf{D-}_t|-N_a(t)$ and $|V\setminus U^\mathbf{D}_t|=|V\setminus U^\mathbf{D-}_t|-|V(F_t)\setminus\phi_t^{-1}(A_t)|$. Since $\chi_t$ assigns the vertices of $V(F_t)\setminus\phi_t^{-1}(A_t)$ essentially independently, each changing $N_a(t)$ by at most $1$, Lemma~\ref{lem:McDiarmidInequ} gives $N_a(t)=\tfrac12|V(F_t)\setminus\phi_t^{-1}(A_t)|\pm n^{0.6}$ with probability at least $1-\exp(-n^{0.1})$, and this a.a.s.\ holds simultaneously for all $t\in \cJ_0$ by the union bound over the at most $n$ choices of $t$ (recall that the random assignment $\chi_s(y)\in\{\Vmin,\Vplus\}$ for $y\in V(F_s)\setminus \phi_s^{-1}(A_s)$ determines how many vertices of $V_a$ get used in the second step of Stage~D). Combining these,
		\[|V_a\setminus U^\mathbf{D}_t|=(1\pm \tfrac{\gamma'}2)\frac{|V\setminus U^\mathbf{D-}_t|}{2}-\tfrac12|V(F_t)\setminus\phi_t^{-1}(A_t)|\pm n^{0.6}=(1\pm\tfrac{\gamma'}2)\frac{|V\setminus U^\mathbf{D}_t|}{2}\pm n^{0.6}\,,\]
		and since $|V\setminus U^\mathbf{D}_t|=\Theta(n)$, the additive $n^{0.6}$ term is at most $\tfrac{\gamma'}2\cdot\tfrac{|V\setminus U^\mathbf{D}_t|}{2}$ for $n$ sufficiently large, giving
		$|V_a\setminus U^\mathbf{D}_t|=(1\pm \gamma')\frac{|V\setminus U^\mathbf{D}_t|}{2}$ for all $t\in \cJ_0$, as required.
		We now substitute the formula for $|X|$ given by~\eqref{eq:XgeqNuN} into~\eqref{eq:XcapUHD}. The power $\dStageC^{|S_1|+|S_2|}$ from that formula combines with $\big(\tfrac{\dStageD}{\dStageC}\big)^{|S_1|+|S_2|}$ to give $\dStageD^{|S_1|+|S_2|}$, and the $T_3$ factors $\prod_{\ell\in T_3}\tfrac{|J_\ell|}{n/2}\cdot\tfrac n2$ are unaffected, exactly as required by Definition~\ref{def:index-quasirandom} for $\mathbb{U}_{\HStageD}(S_1,S_2,T_1,T_2,T_3)$. It remains to check the $T_1,T_2$ factors: for each $t\in T_1$ we must show
		\[\Big(1-\tfrac{|U^\mathbf{D-}_t|}n\Big)\cdot\frac{|\Vmin\setminus U^\mathbf{D}_t|}{|\Vmin\setminus U^\mathbf{D-}_t|}=(1\pm3\gamma')\Big(1-\tfrac{|U^\mathbf{D}_t|}n\Big)\,,\]
		(and likewise with $\Vplus$ in place of $\Vmin$ for $t\in T_2$). Since $|V|=n$ we have $|V\setminus U_t|=n\big(1-|U_t|/n\big)$ for either used set, so the two identities established above give $|\Vmin\setminus U^\mathbf{D-}_t|=(1\pm\tfrac{\gamma'}2)\cdot\tfrac n2\big(1-|U^\mathbf{D-}_t|/n\big)$ and $|\Vmin\setminus U^\mathbf{D}_t|=(1\pm\gamma')\cdot\tfrac n2\big(1-|U^\mathbf{D}_t|/n\big)$; substituting these into the displayed ratio, the factor $\tfrac n2$ occurs in both the numerator and the denominator and so cancels exactly, while the factor $\big(1-|U^\mathbf{D-}_t|/n\big)$ in the denominator cancels against the leading factor of the left-hand side, leaving $\tfrac{1\pm\gamma'}{1\pm\gamma'/2}=(1\pm3\gamma')$ as claimed. Since $|T_1\cup T_2|\le2\LD$, taking the product over $t\in T_1\cup T_2$ multiplies this error by at most $2\LD$, again dominated by $C'(2\LD+3)\gamma'$. Combining this with the~\ref{pathpack:megaqr} error and the $\gamma'/2$ error from~\ref{quasi:indexquasi} after Stage~D$^-$, and observing that the resulting set is precisely $\mathbb{U}_{\HStageD}(S_1,S_2,T_1,T_2,T_3)$ with families of sets $(U^\mathbf{D}_s)_{s\in\cJ}$ and $(J_s)_{s\in\cJ_0}$, this is what is required for $(\LD,3C'(2\LD+3)\gamma',\dStageD,\dStageD)$-index-quasirandomness.
	\end{claimproof}
	
	\begin{claimInLemmaD}\label{UMDS:Quasi234}
		The setup satisfies~\ref{quasi:anchorsets}, \ref{quasi:PathsWithOneVertex}, and \ref{quasi:termVtxNbs}.
	\end{claimInLemmaD}
	\begin{claimproof}[Proof of Claim~\ref{UMDS:Quasi234}]
		For~\ref{quasi:anchorsets}, observe that for $\cJ_1$ and $\cJ_2$ nothing has changed from after Stage~D$^-$, while for $\cJ_0$ the required equality holds even with no error by construction. Similarly, for~\ref{quasi:PathsWithOneVertex} nothing has changed from after Stage~D$^-$. For~\ref{quasi:termVtxNbs}, it suffices by construction to establish~\ref{quasi:distributionOfAnchors}.
	\end{claimproof}
	
	\begin{claimInLemmaD}\label{UMDS:Quasi5}
		The setup satisfies~\ref{quasi:distributionOfAnchors}.
	\end{claimInLemmaD}
	\begin{claimproof}[Proof of Claim~\ref{UMDS:Quasi5}]
		Fix $s\in\cJ$, $a\in\{\boxminus,\boxplus\}$ and $u\in V(\HStageD)\setminus U^\mathbf{D}_s$. By~\ref{quasi:distributionOfAnchors} after Stage~D$^-$, writing $A^a_s$ for the set of anchors $v\in A^\mathbf{D}_s$ such that the path end not anchored to $v$ is anchored in $V_a$, we have
		\begin{align*}
			\big|\NBH_{\HStageC}(u)\cap A^a_s\cap\Vmin\big|&=(1\pm\gamma')\dStageC |A^a_s\cap\Vmin|\quad\text{and}\\
			\big|\NBH_{\HStageC}(u)\cap A^a_s\cap\Vplus\big|&=(1\pm\gamma')\dStageC |A^a_s\cap\Vplus|\,.
		\end{align*}
		For each $s\in\cJ$ and $a\in\{\boxminus,\boxplus\}$, we define a weight function $\omega=\omega_{s,a}$ on $V(\HStageC)$ by setting $\omega(v)=1$ if $v\in A^a_s$, and otherwise $\omega(v)=0$; note that $\sum_{v\in V_b}\omega(v)=|A^a_s\cap V_b|$ is either equal to $0$ or at least $\nu n$, for every such $\omega_{s,a}$: for $s\in\cJ_1\cup\cJ_2$, we have $|A^a_s\cap V_b|=(1\pm\gamma')\tfrac14|A_s|$ by~\ref{quasi:anchorsets}, which is positive (and comfortably at least $\nu n$, since $|A_s|=\Theta(n)$); for $s\in\cJ_0$, by construction each path of $F_s$ is anchored at a matched pair $\{\boxminus_i,\boxplus_i\}$ with one endpoint on each side, so an anchor's partner always lies on the side opposite the anchor itself, giving $|A^a_s\cap V_a|=0$ and $|A^a_s\cap V_{a'}|=\tfrac12|A_s|\ge\nu n$ for $a'$ the side other than $a$. The failure probability $\exp(-n^{0.2})$ in~\ref{pathpack:weights} is uniform over all admissible weight functions (it does not depend on which $\omega$ is chosen), so applying~\ref{pathpack:weights} separately with $\omega=\omega_{s,a}$ for each of the $2|\cJ|=O(n)$ pairs $(s,a)$ and taking the union bound, with probability at least $1-n\exp(-n^{0.2})=1-o(1)$ we have, simultaneously for all $s\in\cJ$, $a\in\{\boxminus,\boxplus\}$, and $u\in V(\HStageD)\setminus U^\mathbf{D}_s$,
		\begin{align*}
			\big|\{\boxminus_i\in \NBH_{\HStageD}(u)\cap A^a_s\}\big|&=(1\pm\gamma')\dStageC |A^a_s\cap\Vmin|(1\pm 2C'(2\LD+3)\gamma')\tfrac{\dStageD}{\dStageC}\quad\text{and}\\
			\big|\{\boxplus_i\in \NBH_{\HStageD}(u)\cap A^a_s\}\big|&=(1\pm\gamma')\dStageC |A^a_s\cap\Vplus|(1\pm 2C'(2\LD+3)\gamma')\tfrac{\dStageD}{\dStageC}\,,
		\end{align*}
		where, as in the proof of Claim~\ref{UMDS:Quasi1}, we used that~\ref{pathpack:weights} is stated in terms of $d'_{ab}/d_{ab}$ rather than $\dStageD/\dStageC$, and that the resulting extra error of $(8\LD+13)\gamma'$ (by Claims~\ref{UMDS:BlockQuasiDiet} and~\ref{UMDS:densities}) is absorbed by doubling the coefficient. Combining the two factors $(1\pm\gamma')$ and $(1\pm2C'(2\LD+3)\gamma')$ above (and discarding the resulting negligible $\gamma'^2$ term), we get
		\[\big|\{\boxminus_i\in \NBH_{\HStageD}(u)\cap A^a_s\}\big|=(1\pm3C'(2\LD+3)\gamma')\dStageD |A^a_s\cap\Vmin|\,,\]
		and likewise with $\boxplus$ in place of $\boxminus$, which is as required for~\ref{quasi:distributionOfAnchors} with $\gamma=3C'(2\LD+3)\gamma'$ (matching the final error obtained in Claim~\ref{UMDS:Quasi1}).
	\end{claimproof}
	
	\begin{claimInLemmaD}\label{UMDS:Quasi6}
		The setup satisfies~\ref{quasi:imagecaps}.
	\end{claimInLemmaD}
	\begin{claimproof}[Proof of Claim~\ref{UMDS:Quasi6}]
		Fix $S_1,S_2,T$ as in~\ref{quasi:imagecaps}, each of size at most $\LD$. By~\ref{quasi:imagecaps} after Stage~D$^-$, for each $Q\subset S_2$ we have (recalling that the setup after Stage~D$^-$ has $I_s=I^\mathbf{D-}_s$ by~\ref{D:reroot:quasisetup}, and that $I^\mathbf{D-}_s=J_s$ by construction, as shown in the proof of Lemma~\ref{lem:reroot-assign})
		\[\big|\{s\in\cJ_0:(S_1\cup Q)\cap U^\mathbf{D-}_s=\emptyset,T\subset I^\mathbf{D-}_s\}\big|=(1\pm\gamma')|\cJ_0|(\deltnonspanning+8\sigmKJ^*)^{|S_1|+|Q|}(2\sigmKJ^*)^{|T|}\,,\]
		and by~\ref{pathpack:imagecaps}, applied with $Y:=\{s\in\cJ_0:(S_1\cup Q)\cap U^\mathbf{D-}_s=\emptyset,T\subset I^\mathbf{D-}_s\}$ (of size at least $\nu n$ by the display above) and $S_\boxminus:=(S_1\cup Q)\cap\Vmin$, $S_\boxplus:=(S_1\cup Q)\cap\Vplus$, we get, with probability at least $1-\exp(-n^{0.2})$ (uniform over the choice of underlying sets, so absorbed into the union bound as above, using $U_s^\mathbf{D}=U_s^\mathbf{D-}\cup\im\phi'_s$),
		\begin{equation}\label{eq:Quasi6imagecapsraw}\big|\{s\in Y:(S_1\cup Q)\cap U^\mathbf{D}_s=\emptyset\}\big|=(1\pm C'\gamma)\sum_{s\in Y}\prod_{a\in\{\boxminus,\boxplus\}}\Big(\frac{|V_a\setminus U^\mathbf{D}_s|}{|V_a\setminus U^\mathbf{D-}_s|}\Big)^{|S_a|}\,,\end{equation}
		where $\gamma=(2\LD+3)\gamma'$ as before. We now show every factor $|V_a\setminus U^\mathbf{D}_s|/|V_a\setminus U^\mathbf{D-}_s|$ above is within $(1\pm5\gamma')$ of the same value $(\deltnonspanning+2\sigmKJ^*)/(\deltnonspanning+8\sigmKJ^*)$, regardless of $s\in Y$ or $a$. By~\ref{quasi:indexquasi}, exactly as in the proof of Claim~\ref{UMDS:Ineq}, we have $|V_a\setminus U^\mathbf{D-}_s|=(1\pm\gamma')\tfrac n2(\deltnonspanning+8\sigmKJ^*)$ for every $s\in\cJ$ and $a\in\{\boxminus,\boxplus\}$; summing over $a$ gives $|V\setminus U^\mathbf{D-}_s|=(1\pm\gamma')n(\deltnonspanning+8\sigmKJ^*)$, and subtracting the (deterministic) count $|V(F_s)\setminus\phi_s^{-1}(A_s)|=6\sigmKJ^* n$ of newly-used vertices gives $|V\setminus U^\mathbf{D}_s|=(1\pm3\gamma')n(\deltnonspanning+2\sigmKJ^*)$ for $s\in\cJ_0$, using $\sigmKJ^*\ll\deltnonspanning$ to absorb this subtraction into a slightly larger relative error, exactly as in the derivation of~\eqref{eq:XcapUHD}. Combined with $|V_a\setminus U^\mathbf{D}_s|=(1\pm\gamma')\tfrac{|V\setminus U^\mathbf{D}_s|}2$ (established above, for all $s\in\cJ_0$), this gives $|V_a\setminus U^\mathbf{D}_s|=(1\pm4\gamma')\tfrac n2(\deltnonspanning+2\sigmKJ^*)$, and so
		\[\frac{|V_a\setminus U^\mathbf{D}_s|}{|V_a\setminus U^\mathbf{D-}_s|}=(1\pm5\gamma')\frac{\deltnonspanning+2\sigmKJ^*}{\deltnonspanning+8\sigmKJ^*}\]
		as claimed, for every $s\in Y\subset\cJ_0$ and every $a$. Substituting this into~\eqref{eq:Quasi6imagecapsraw}, each term of the sum over $s\in Y$ becomes $(1\pm5\gamma')^{|S_\boxminus|+|S_\boxplus|}\big(\tfrac{\deltnonspanning+2\sigmKJ^*}{\deltnonspanning+8\sigmKJ^*}\big)^{|S_1|+|Q|}$ (using $|S_\boxminus|+|S_\boxplus|=|S_1\cup Q|=|S_1|+|Q|$), independent of $s$, so the sum is just $|Y|$ copies of this common value; bounding $(1\pm5\gamma')^{|S_1|+|Q|}$, with $|S_1|+|Q|\le2\LD$, by $(1\pm10\LD\gamma')$, and combining with the~\ref{pathpack:imagecaps} error $C'\gamma=C'(2\LD+3)\gamma'$ (which dominates, since $C'$ is larger than any fixed polynomial in $\LD$), and finally substituting $|Y|=(1\pm\gamma')|\cJ_0|(\deltnonspanning+8\sigmKJ^*)^{|S_1|+|Q|}(2\sigmKJ^*)^{|T|}$ from the display above (and, in the set on the left below, writing $I^\mathbf{D}_s$ for the analogous index set after Stage~D, using $I^\mathbf{D}_s=I^\mathbf{D-}_s=J_s$: the second-step embedding in Stage~D only adds vertices to $U_s$, not to the anchor set $A_s$, so $A^\mathbf{D}_s=A^\mathbf{D-}_s$, and $I_s$ depends only on $A_s$), it is likely that we have
		\begin{multline}\label{eq:Quasi6imagecapsSQ}
			\big|\{s\in\cJ_0:(S_1\cup Q)\cap U^\mathbf{D}_s=\emptyset,T\subset I^\mathbf{D}_s\}\big|\\
			\begin{aligned}
				&=(1\pm\gamma')|\cJ_0|(\deltnonspanning+8\sigmKJ^*)^{|S_1|+|Q|}(2\sigmKJ^*)^{|T|}(1\pm 2C'(2\LD+3)\gamma')\big(\tfrac{\deltnonspanning+2\sigmKJ^*}{\deltnonspanning+8\sigmKJ^*}\big)^{|S_1|+|Q|}\\
				&=|\cJ_0|(\deltnonspanning+2\sigmKJ^*)^{|S_1|+|Q|}(2\sigmKJ^*)^{|T|}\pm 3C'(2\LD+3)\gamma'n\,.
			\end{aligned}
		\end{multline}
		Since $\ONE[S_2\subset U^\mathbf{D}_s]=\prod_{x\in S_2}\big(1-\ONE[x\notin U^\mathbf{D}_s]\big)=\sum_{Q\subset S_2}(-1)^{|Q|}\ONE[Q\cap U^\mathbf{D}_s=\emptyset]$, and $S_1\cap U^\mathbf{D}_s=\emptyset$ together with $Q\cap U^\mathbf{D}_s=\emptyset$ is the same as $(S_1\cup Q)\cap U^\mathbf{D}_s=\emptyset$, summing over $s\in\cJ_0$ with $S_1\cap U^\mathbf{D}_s=\emptyset$ and $T\subset I^\mathbf{D}_s$ gives
		\[\big|\{s\in\cJ_0:S_1\cap U^\mathbf{D}_s=\emptyset,S_2\subset U^\mathbf{D}_s,T\subset I^\mathbf{D}_s\}\big|=\sum_{Q\subset S_2}(-1)^{|Q|}\big|\{s\in\cJ_0:(S_1\cup Q)\cap U^\mathbf{D}_s=\emptyset,T\subset I^\mathbf{D}_s\}\big|\,.\]
		Substituting~\eqref{eq:Quasi6imagecapsSQ} for each of the (at most $2^{\LD}$, since $|S_2|\le\LD$) terms in this sum, and using the binomial theorem $\sum_{Q\subset S_2}(-1)^{|Q|}(\deltnonspanning+2\sigmKJ^*)^{|Q|}=\big(1-(\deltnonspanning+2\sigmKJ^*)\big)^{|S_2|}$ to simplify the main terms, while crudely bounding the sum of the $2^{|S_2|}\le2^{\LD}$ error terms $\pm3C'(2\LD+3)\gamma'n$ by their number times their size, we obtain
		\begin{align*}
		&\big|\{s\in\cJ_0:S_1\cap U^\mathbf{D}_s=\emptyset,S_2\subset U^\mathbf{D}_s,T\subset I^\mathbf{D}_s\}\big|\\
		&=|\cJ_0|(\deltnonspanning+2\sigmKJ^*)^{|S_1|}(1-\deltnonspanning-2\sigmKJ^*)^{|S_2|}(2\sigmKJ^*)^{|T|}\pm 2^{\LD}\cdot3C'(2\LD+3)\gamma'n\,.
		\end{align*}
		To convert this to the multiplicative form required by~\ref{quasi:imagecaps}, note that since $|S_1|,|S_2|,|T|\le\LD$ and each of $\deltnonspanning+2\sigmKJ^*$, $1-\deltnonspanning-2\sigmKJ^*$, $2\sigmKJ^*$ lies in $(0,1)$, the main term above is at least
		\[|\cJ_0|\big[\deltnonspanning(1-\deltnonspanning-2\sigmKJ^*)(2\sigmKJ^*)\big]^{\LD}=\sigmJnula n\big[\deltnonspanning(1-\deltnonspanning-2\sigmKJ^*)(2\sigmKJ^*)\big]^{\LD}\,,\]
		so the additive error above is at most a $\gamma$-fraction of the main term, for
		\[\gamma:=\frac{2^{\LD}\cdot3C'(2\LD+3)}{\sigmJnula\big[\deltnonspanning(1-\deltnonspanning-2\sigmKJ^*)(2\sigmKJ^*)\big]^{\LD}}\gamma'\,.\]
		Since $C'=\exp(2000CL\nu^{-20L})$ for $L=2\LD$ (where $\nu$ is itself exponentially small in $\LD$), $C'$ dominates the merely exponential-in-$\LD$ prefactor above, so $\gamma\le C'^2(2\LD+3)\gamma'$; hence
		\[\big|\{s\in\cJ_0:S_1\cap U^\mathbf{D}_s=\emptyset,S_2\subset U^\mathbf{D}_s,T\subset I^\mathbf{D}_s\}\big|=\big(1\pm C'^2(2\LD+3)\gamma'\big)|\cJ_0|(\deltnonspanning+2\sigmKJ^*)^{|S_1|}(1-\deltnonspanning-2\sigmKJ^*)^{|S_2|}(2\sigmKJ^*)^{|T|}\,,\]
		which is as required for~\ref{quasi:imagecaps} with $\gamma=C'^2(2\LD+3)\gamma'$.
	\end{claimproof}
	
	\begin{claimInLemmaD}\label{UMDS:Quasi78}
		The setup satisfies~\ref{quasi:NumAnchors} and~\ref{quasi:NumAnchorsNoIm}.
	\end{claimInLemmaD}
	\begin{claimproof}[Proof of Claim~\ref{UMDS:Quasi78}]
		For~\ref{quasi:NumAnchors}, observe that nothing has changed since after Stage~D$^-$. The same holds for~\ref{quasi:NumAnchorsNoIm} for $i=1,2$. For $\cJ_0$, fix distinct $u,v\in V(\HStageD)$ such that for no $j$ we have $\{u,v\}=\{\boxminus_j,\boxplus_j\}$. Let $Z=\{s\in\cJ_0:u\in A^\mathbf{D}_s,v\notin U^{\mathbf{D-}}_s\}$. Note that $A^\mathbf{D}_s=A^\mathbf{D-}_s$. By~\ref{quasi:NumAnchorsNoIm} after Stage~D$^-$, we have $|Z|=(1\pm\gamma')|\cJ_0|\cdot 2\sigmKJ^*(\deltnonspanning+8\sigmKJ^*)$.
		 By~\ref{pathpack:imagecaps}, applied with $Y:=Z$ and (writing $v\in V_b$) $S_b:=\{v\}$, $S_{\bar b}:=\emptyset$, we get, exactly as in the proof of Claim~\ref{UMDS:Quasi6}, with probability at least $1-\exp(-n^{0.2})$ (again absorbed into a union bound over the at most $n^2$ choices of $(u,v)$),
		\[\big|\{s\in Z:v\notin U^\mathbf{D}_s\}\big|=(1\pm C'\gamma)\sum_{s\in Z}\frac{|V_b\setminus U^\mathbf{D}_s|}{|V_b\setminus U^\mathbf{D-}_s|}\,,\]
		where $\gamma=(2\LD+3)\gamma'$. Since $|V_b\setminus U^\mathbf{D}_s|/|V_b\setminus U^\mathbf{D-}_s|=(1\pm5\gamma')(\deltnonspanning+2\sigmKJ^*)/(\deltnonspanning+8\sigmKJ^*)$ for every $s\in\cJ_0$ (established in the proof of Claim~\ref{UMDS:Quasi6}), the summand does not depend on $s$, so the sum is just $|Z|$ copies of this common value, giving
		\[\big|\{s\in Z:v\notin U^\mathbf{D}_s\}\big|=(1\pm 2C'(2\LD+3)\gamma')\frac{\deltnonspanning+2\sigmKJ^*}{\deltnonspanning+8\sigmKJ^*}\cdot (1\pm\gamma')|\cJ_0|\cdot 2\sigmKJ^*(\deltnonspanning+8\sigmKJ^*)\;.\]
		The factors $\deltnonspanning+8\sigmKJ^*$ cancel exactly, and combining the two remaining error terms (as in Claim~\ref{UMDS:Quasi1}) gives the simplified form
		\[\big|\{s\in Z:v\notin U^\mathbf{D}_s\}\big|=(1\pm 3C'(2\LD+3)\gamma')\cdot 2\sigmKJ^*(\deltnonspanning+2\sigmKJ^*)|\cJ_0|\;,\]
		as required for~\ref{quasi:NumAnchorsNoIm} with $\gamma=3C'(2\LD+3)\gamma'$.
	\end{claimproof}
	Each error term obtained in Claims~\ref{UMDS:Quasi1}--\ref{UMDS:Quasi78} above is $\gamma'$ multiplied by a quantity depending only on $\LD$ and $C'$ (itself a function of $\nu$ and $2\LD$ alone), so, exactly as in the bound $(C'')^{-1}\gg\gamma',\gamAnchor$ noted at the start of this section, the hierarchy~\eqref{eq:CONSTANTS} lets us choose $\gamcore$ (and hence $\gamma'$) small enough that all of these error terms are $\ll\gamAnchor$. This confirms the $(\gamAnchor,\LD,\dStageD,\dStageD)$-quasirandom setup claimed above.
\end{proof}

\section{Stage~E (Proof of Lemma~\ref{lem:StageENew})}\label{sec:StageE2}

In this stage, we need to pack the paths $\SpecPaths^*_s$ for $s\in\cJ_2$. Recall that each of these paths has 7~edges. For each $s\in\cJ_2$, the set $\SpecPaths_s$ contains exactly $\sigmKJ n$ paths; in Stage~C we packed at most $n^{0.6}$ of these paths to leave $\SpecPaths^*_s$. The goal of this section is to pack a specified number of edges of these paths into each of $\Vmin$ and $\Vplus$, and between the sides $\Vmin$ and $\Vplus$. Due to our choice $\sigmKJ\gg\sigmJnula,\sigmJjedna$ in~\eqref{eq:CONSTANTS}, the density of the graph $\HStageE$ after this stage will be much smaller than the density of $\HStageD$. Thus we need to embed roughly half the edges of the paths $\bigcup_{s\in\cJ_2}\SpecPaths^*_s$ crossing between $\Vmin$ and $\Vplus$, with the other roughly half about equally split between sides. To begin with, we assign to each path a number of edges it should use within $\Vmin$, crossing between $\Vmin$ and $\Vplus$, and within $\Vplus$ (Lemma~\ref{lem:E2decideParts}). Observe that a path with one anchor on each side necessarily has an \emph{odd} number of crossing edges, while a path with both anchors on the same side has an \emph{even} number of crossing edges. Having done this, we apply the path-packing lemma, Lemma~\ref{lem:pathpack}, to prove Lemma~\ref{lem:StageENew}.

We will use the following definition with $X$ a set of size $\big|\bigcup_{s\in\cJ_2}\SpecPaths^*_s\big|$, where our goal is to use $N_1$ edges within $\Vmin$, $N_2$ edges crossing between the sides, and $N_3$ edges within $\Vplus$. The map $m$ identifies which paths are odd and which even, while the occupancy assignment $\aleph(x)$ says how many edges the path $x$ needs to use in $\Vmin$, crossing between sides, and in $\Vplus$ in that order.

\begin{definition}\label{def:occupancies}
Suppose that $X$ is a finite set, and $m:X\rightarrow\{\mathrm{odd},\mathrm{even}\}$ is a map. Then we say that a triple $(N_1,N_2,N_3)$ of non-negative integers is a \emph{moderately balanced occupancy requirement}\index{moderately balanced occupancy requirement} if
\begin{enumerate}[label=\abc]
	\item\label{en:occ:seven} $N_1+N_2+N_3=7 |X|$,
	\item\label{en:occ:balanced} $3|X| \le N_2\le 4|X|$ 
		and $N_1,N_3\ge |X|$, and
	\item\label{en:occ:reqPar} the parity of $N_2$ is the same as the parity of $|m^{-1}(\mathrm{odd})|$.
\end{enumerate}
Given a moderately balanced occupancy requirement $(N_1,N_2,N_3)$, we say that a map $\aleph:X\rightarrow \ZZ^3$ is an \emph{occupancy assignment}\index{occupancy assignment} for $(N_1,N_2,N_3)$ if the following hold.
\begin{enumerate}[resume,label=\abc]
	\item\label{en:occ:nonneg} for each $x\in X$, we have that 
	$\aleph(x)_1\ge 1$, $\aleph(x)_3 \ge 1$ and $\aleph(x)_2\ge 2$,
	\item\label{en:occ:sum7} for each $x\in X$, we have that $\sum_{i=1}^3\aleph(x)_i=7$,
	\item\label{en:occ:assPar} for each $x\in X$, we have that $\aleph(x)_2$ is odd if and only if $m(x)=\mathrm{odd}$, 
	\item\label{en:occ:sumN} for each $i\in\{1,2,3\}$, we have that $N_i=\sum_{x\in X}\aleph(x)_i$.
\end{enumerate}
We call vectors in the image of $\aleph$ \index{occupancy vector}\emph{occupancy vectors}.
\end{definition}
In the definition of a moderately balanced occupancy requirement,~\ref{en:occ:seven} says that in total we assign all $7|X|$ edges of the paths and~\ref{en:occ:balanced} that the desired assignment is not too far from the quarter-half-quarter distribution mentioned above.

The first condition~\ref{en:occ:nonneg} of an occupancy assignment is something we will need to help embed paths later in this stage; the remaining conditions simply say that we assign all 7~edges of each path, that precisely the odd paths have an odd number of crossing edges assigned, and that in total we assign the desired number of edges to each side and crossing. What an occupancy assignment does \emph{not} specify is which vertices of each path will be assigned to which side; we will do that later in the stage.

The next easy lemma says that given a moderately balanced occupancy requirement, we can find an occupancy assignment.
\begin{lemma}\label{lem:OccAssign}
Suppose that $X$ is a finite set, that $m:X\rightarrow\{\mathrm{odd},\mathrm{even}\}$ is a map, and that $(N_1,N_2,N_3)$ is a moderately balanced occupancy requirement for $X$ and $m$. Then there exists an occupancy assignment for $(N_1,N_2,N_3)$.
\end{lemma}
\begin{proof}
Let a set $X=\{x_1,x_2,\ldots,x_{|X|}\}$, a map $m:X\rightarrow \{\mathrm{odd},\mathrm{even}\}$
and a moderately balanced occupancy requirement be given. 
We will show that an occupancy assignment $\aleph$ can be found by Algorithm~\ref{alg:occ} below.

  \begin{algorithm}[ht]\label{alg:occ}
    \caption{Generating an occupancy assignment}    
\SetKwInOut{Input}{Input}
\SetKwInOut{Output}{Output}
\Input{~a set $X=\{x_1,x_2,\ldots,x_{|X|}\}$, map $m:X\rightarrow \{\mathrm{odd},\mathrm{even}\}$ and moderately balanced occupancy requirement~$(N_1,N_2,N_3)$}
\Output{~an occupancy assignment $\aleph:X \rightarrow \ZZ^3$}
    \For{$k=1,2,\ldots |X|$}{
    \quad\mbox{} $\aleph(x_k)_1=\aleph(x_k)_3=1$ and $\aleph(x_k)_2=2$\;
         \lIf{$m(x_k)=\mathrm{odd}$\\ 
         \quad\mbox{}}{
        $\aleph(x_k)_2=3$}}
    \While{$\sum_k \aleph(x_k)_2 < N_2$}{
    \quad\mbox{} $t = \min \{ \ell\in [|X|]:~ \aleph(x_{\ell})_2 \le \aleph(x_k)_2~ \text{for all $k$}\}$\;
	\quad\mbox{} $\aleph(x_t)_2 \leftarrow \aleph(x_t)_2+2$\;}
    \While{$\sum_{i,k} \aleph(x_k)_i < \sum_i N_i$}{
    \quad\mbox{} 
    		$t = \min \{ \ell\in [|X|]:~ \sum_i \aleph(x_{\ell})_i < 7\}$\;
    	\quad\mbox{} 
		$j = \min \{ i\in [3]:~ \sum_{\ell} \aleph(x_{\ell})_i < N_i\}$\;
	\quad\mbox{} 
 		$\aleph(x_t)_j \leftarrow \aleph(x_t)_j+1$;       
     }
  \end{algorithm}

The for loop sets initial values for $\aleph(x_k)_i$
for every $k\in [|X|]$ and $i\in [3]$.
Right at the end of that loop we have
$\sum_k \aleph(x_k)_1 = |X| \le N_1$,
$\sum_k \aleph(x_k)_3 = |X| \le N_3$
and
$\sum_k \aleph(x_k)_2\le 3|X|$. 
By the if line, the parity of $\sum_k\aleph(x_k)_2$ is the same as the parity
of $|m^{-1}(\mathrm{odd})|$ and hence, by~\ref{en:occ:balanced} and~\ref{en:occ:reqPar}, $N_2-\sum_k \aleph(x_k)_2\ge 0$
is even. Further, 
$|\aleph(x_{k_1})_2-\aleph(x_{k_2})_2|\le 2$
for every $k_1,k_2\in [|X|]$.

During the subsequent while loop, the value of
$\sum_k \aleph(x_k)_2$ gets increased by 2
in each iteration, hence forcing the loop to stop
eventually with $N_2=\sum_k \aleph(x_k)_2$.
By the choice of $x_t$ in each such iteration,
we also ensure that
$|\aleph(x_{k_1})_2-\aleph(x_{k_2})_2|\le 2$
is maintained for every $k_1,k_2\in [|X|]$.
In particular, using $N_2\le 4|X|$ from~\ref{en:occ:balanced},
we know that $\aleph(x_{k})_2\le 5$ needs to hold
for every $k\in [|X|]$ by the end of that while loop.
Hence, at the end of this while loop, we have 
$\sum_i \aleph(x_{\ell})_i\le 7$ for every $\ell\in [|X|]$,
and $\sum_{\ell} \aleph(x_\ell)_i\le N_i$ for every $i\in [3]$.

Afterwards, during the last loop, $\aleph(x_k)_2$
is not updated for any $k\in [|X|]$.
However, as long as $\sum_{i,k} \aleph(x_k)_i < \sum_i N_i$
holds, the following is done. We fix
$t = \min \{ \ell\in [|X|]:~ \sum_i \aleph(x_{\ell})_i < 7\}$
which must exist, since otherwise we would
have $\sum_i \aleph(x_{\ell})_i \ge 7$ for every $\ell \in [|X|]$,
which leads to $\sum_{i,k} \aleph(x_k)_i \ge 7|X| = \sum_i N_i$,
in contradiction to the condition for processing an iteration of the while loop.
We also fix
$j = \min \{ i\in [3]:~ \sum_{\ell} \aleph(x_{\ell})_i < N_i\}$
whose existence is ensured analogously. Then we update $\aleph$ by increasing
$\aleph(x_t)_j$ by 1.
Therefore, since $\sum_{i,k} \aleph(x_k)_i$ is increased with every iteration, the while loop and hence the whole algorithm must end eventually.
Inductively, using the choice of $t$ and $j$, one verifies easily that throughout the while loop we always maintain
$\sum_i \aleph(x_{\ell})_i\le 7$ for every $\ell\in [|X|]$,
and $\sum_{\ell} \aleph(x_\ell)_i\le N_i$ for every
$i\in [3]$. Using all this, we will show now that the output $\aleph$ satisfies Properties~\ref{en:occ:nonneg}--\ref{en:occ:sumN}.

Property~\ref{en:occ:nonneg} already holds at the end of the
for loop. Afterwards it is kept to be true, since
the values $\aleph(x_k)_i$ never decrease throughout the algorithm.

Suppose that Property~\ref{en:occ:sum7} fails. Then, by the observations for the second while loop, there would have to exist some $\ell\in [|X|]$ with
$\sum_{i} \aleph(x_\ell)_i < 7$ while
$\sum_i \aleph(x_{k})_i\le 7$ for every $k\neq \ell$.
However, this would lead to
$\sum_{i,k} \aleph(x_k)_i < 7|X| = \sum_i N_i$ by~\ref{en:occ:seven},
in contradiction to the while loop having stopped.
Property~\ref{en:occ:sumN} is proved analogously.

Finally, Property~\ref{en:occ:assPar} holds,
since for every $k\in [|X|]$, 
the initial value for $\aleph(x_k)_2$ in the for loop
is chosen such that it satisfies~\ref{en:occ:assPar}, 
and since in the subsequent while loops 
we never change the parity of $\aleph(x_k)_2$.
\end{proof}

The following lemma is the main work of this section. It returns a map $\chi$ which assigns to each vertex $x$ of $\bigcup_{s\in\cJ_2}\SpecPaths^*_s$ a side (i.e.\ $\Vmin$ or $\Vplus$) to which we will embed $x$. Of course, $\chi$ has to assign the embedded endvertices to the sides on which they are embedded (property~\ref{en:kaficko} of the lemma). Observe that once this map is fixed, we know exactly how many edges will be used in Stage~E within $\Vmin$, within $\Vplus$ and between $\Vmin$ and $\Vplus$. Property~\ref{en:Etotaledges} of the lemma will give us a $\chi$ which guarantees that the number of edges used in each set is correct (i.e.\ that the requirement on the number of edges~\eqref{eq:E:edges} of Lemma~\ref{lem:StageENew} in each place will be met).
What complicates the lemma is that in order to actually pack the paths and maintain quasi\-random\-ness we need some additional properties of $\chi$ (parts~\ref{en:medunkovycaj},~\ref{en:nocheins},~\ref{en:abstractpscatt} and~\ref{en:pairdistrpropE} of the lemma).

For compactness, in this lemma and its proof we use the following abbreviations, where $a\in\{\boxminus,\boxplus\}$ and $s\in\cJ_2$:
\begin{align*}
 E^*&:=\bigcup_{s\in \cJ_2} E(\SpecPaths^*_s)\index{$E^*$}\quad\text{and}\\
 X_s^a&:=\big\{x\in V(\SpecPaths^*_s):\chi(x)=V_a\text{ and }x\text{ is not an endvertex of }\SpecPaths^*_s\big\}\,.
\end{align*}
Moreover, each path $P\in\bigcup_{s\in\cJ_2}\SpecPaths^*_s$ lies in $\SpecPaths^*_s$ for exactly one $s\in\cJ_2$; in this lemma and its proof we therefore suppress this index and write $\phi^\mathbf{A}(x)$ and $\phi^\mathbf{D}(x)$ for $\phi^\mathbf{A}_s(x)$ and $\phi^\mathbf{D}_s(x)$ whenever $x$ is a vertex of such a path~$P$. Note that for $s\in\cJ_2$ the maps $\phi^\mathbf{A}_s$ and $\phi^\mathbf{D}_s$ agree on the endvertices of $\SpecPaths^*_s$, since those vertices were already embedded in Stage~A and $\phi^\mathbf{D}_s=\phi^\mathbf{C}_s$ by~\ref{D:two}; we use whichever of the two is more natural in context.

\begin{lemma}\label{lem:E2decideParts}
There exists a map $\chi:\bigcup_{s\in \cJ_2}V\left(\SpecPaths^*_s\right)\rightarrow\{\Vmin,\Vplus\}$ with the following properties.
\begin{enumerate}[label=\rom]
\item\label{en:kaficko} For each $P\in \bigcup_{s\in \cJ_2}\SpecPaths^*_s$ and each $x\in \{\leftpath_0(P),\rightpath_0(P)\}$ we have 
$\phi^\mathbf{A}(x) \in \chi(x)$.
\item\label{en:medunkovycaj} For each $a\in\{\boxminus,\boxplus\}$ and $s\in\cJ_2$, we have $\big|X^a_s\big|\ge \tfrac15\sigmKJ n$.
\item\label{en:nocheins} For each $a,b\in\{\boxminus,\boxplus\}$ and $s\in\cJ_2$, we have
\[\big|\big\{(x,y)\,:\,xy\in E(\SpecPaths^*_s),x\in X^a_s,y\in X^b_s\big\}\big|\ge \tfrac1{16}\sigmKJ  n\,.\]
\item\label{en:Etotaledges} We have
\begin{align}
\label{dminmin}
\sum_{xy\in E^*}\ONE_{\chi(x),\chi(y)=\Vmin} &=e_{\HStageD}(\Vmin)-\left(\sigmJnula\sigmKJ^* n^2+2 j_{\boxminus\boxminus}+3(j_{\boxminus\boxplus}+ j_{\boxplus\boxplus})\right)\;,\\
\label{dplusplus}
\sum_{xy\in E^*}\ONE_{\chi(x),\chi(y)=\Vplus}
&=e_{\HStageD}(\Vplus)-\left(\sigmJnula\sigmKJ^* n^2+2 j_{\boxplus\boxplus}+ 3(j_{\boxminus\boxplus}+ j_{\boxminus\boxminus})\right) \;,\\	
\label{dminplus}
\sum_{xy\in E^*}\ONE_{\{\chi(x),\chi(y)\}=\{\Vmin,\Vplus\}}
&=e_{\HStageD}(\Vmin,\Vplus)-\left(\sigmJnula\sigmKJ^* n^2+6( j_{\boxminus\boxminus}+ j_{\boxplus\boxplus})+5 j_{\boxminus\boxplus}\right)
\;.
\end{align}
\item\label{en:abstractpscatt} For every $s\in \cJ_2$, the pair $\left(\phi^\mathbf{D}_{\restriction\textrm{endvertices of $\SpecPaths^*_s$}},
\chi_{\restriction\textrm{inner vertices of $\SpecPaths^*_s$}}\right)$ has the $6\gamAnchor$-anchor distribution property with respect to the path-forest $\SpecPaths^*_s$, partition $(\Vmin,\Vplus)$ and the graph $\HStageD$.
\item\label{en:pairdistrpropE}
The collection of path-forests $(\SpecPaths^*_s)_{s\in \cJ_2}$ together with used sets $(\imD(s))_{s\in \cJ_2}$, anchorings $(\phi^\mathbf{D}_{\restriction\textrm{endvertices of  $\SpecPaths^*_s$}})_{s\in \cJ_2}$ and assignments $(\chi_{\restriction\SpecPaths^*_s})_{s\in \cJ_2}$ has the $20\gamAnchor$-pair distribution property.
\end{enumerate}
\end{lemma}

The proof of this lemma goes as follows. First, we check that the numbers of edges we need to use in $\Vmin$, crossing, and in $\Vplus$ according to property~\ref{en:Etotaledges} give us a moderately balanced occupancy requirement (Claim~\ref{cl:IsModeratelyBalanced}).

We then split the paths $\bigcup_{s\in\cJ_2}\SpecPaths^*_s$ into three sets $\mathcal{Q}_{\boxminus\boxminus}$, $\mathcal{Q}_{\boxminus\boxplus}$, $\mathcal{Q}_{\boxplus\boxplus}$ according to the sides to which the endvertices are embedded by the maps $\phi^{\mathbf{D}}_s$. We then create three corresponding auxiliary sets $X_{\boxminus\boxminus}$, $X_{\boxminus\boxplus}$ and $X_{\boxplus\boxplus}$ with $|X_{ab}|=|\mathcal{Q}_{ab}|$ for each choice of $a$ and $b$, and use Lemma~\ref{lem:OccAssign} to find an occupancy assignment $\aleph$ on the auxiliary sets.
The purpose of this is the following: the number of $x\in X_{\boxminus\boxminus}$ which are assigned (for example) $(4,2,1)$ by $\aleph$ is the number of paths in $\mathcal{Q}_{\boxminus\boxminus}$ for which we want to use $4$ edges in $\Vmin$, $2$ crossing, and $1$ in $\Vplus$. We apply Lemma~\ref{lem:OccAssign} to the auxiliary sets rather than to the sets of paths themselves because we want to choose the specific paths assigned $(4,2,1)$ randomly in the next step.

We make the following random choices. First, we pick a bijective map $\pi_{ab}:\mathcal{Q}_{ab}\to X_{ab}$ uniformly at random for each choice of $a$ and $b$. Second, we partition $X=X^\triangledown\dcup X^\blacktriangle$ uniformly at random. These two random choices determine $\chi$ as follows. For a given $P\in\mathcal{Q}_{ab}$, we check whether $\pi_{ab}(P)$ is in $X^\triangledown$ or $X^\blacktriangle$, and look at one of six tables according to the three possible values of $ab$ and choice of $X^\triangledown$ or $X^\blacktriangle$. We look in the given table at the row for the occupancy vector $\aleph\big(\pi_{ab}(P)\big)$ and this row defines $\chi$ on $P$.

One can verify that properties~\ref{en:kaficko},~\ref{en:medunkovycaj} and~\ref{en:Etotaledges} of the lemma hold deterministically, directly from these tables. We will show that~\ref{en:nocheins} holds with high probability over the random choice of $X^\triangledown$ and $X^\blacktriangle$. Finally, we will show that~\ref{en:abstractpscatt} and~\ref{en:pairdistrpropE} each hold with high probability over all the random choices.
\begin{proof}[Proof of Lemma~\ref{lem:E2decideParts}]
Let us write $D_{\boxminus\boxminus}$, $D_{\boxplus\boxplus}$, and $D_{\boxminus\boxplus}$ for the quantities on the right-hand sides of~\eqref{dminmin}, \eqref{dplusplus}, and~\eqref{dminplus}, respectively. Define a map $m:\bigcup_{s\in \cJ_2}\SpecPaths^*_s\rightarrow\{\mathrm{odd},\mathrm{even}\}$ as follows. If $P\in \bigcup_{s\in \cJ_2}\SpecPaths^*_s$ is such that $\phi^\mathbf{A}(\leftpath_0(P))$ and $\phi^\mathbf{A}(\rightpath_0(P))$ both belong to $\Vmin$, or both to $\Vplus$, then we set $m(P):=\mathrm{even}$. Otherwise, set $m(P):=\mathrm{odd}$.

\begin{claim}\label{cl:IsModeratelyBalanced}
The triple $(D_{\boxminus\boxminus},D_{\boxminus\boxplus},D_{\boxplus\boxplus})$ is a moderately balanced occupancy requirement for the set $\bigcup_{s\in \cJ_2}\SpecPaths^*_s$ and the map $m$. 
\end{claim}
\begin{claimproof}[Proof of Claim~\ref{cl:IsModeratelyBalanced}]
Let us first verify Definition~\ref{def:occupancies}\ref{en:occ:seven}. By summing up the right-hand sides of \eqref{dminmin}--\eqref{dminplus}, we get 
\begin{align}
\begin{split}\label{horriblesoundsbybrucespringsteen}
D_{\boxminus\boxminus}+D_{\boxplus\boxplus}+D_{\boxminus\boxplus}&=e(\HStageD)-3\sigmJnula\sigmKJ^* n^2-11(j_{\boxminus\boxminus}+j_{\boxplus\boxplus}+j_{\boxminus\boxplus})\\
&=e(\HStageD)-3\sum_{s\in\cJ_0}|\SpecShortPaths_s|-11\sum_{s\in\cJ_1}|\SpecPaths^*_s|\;.
\end{split}
\end{align}
Recall that we plan to eventually use every edge of $H$ for the packing (cf.~Section~\ref{ssec:IneqEqual} and Definition~\ref{def:family}\ref{enu:familyEachG}). Some edges of $H$ were used prior to Stage~E, $7\sum_{s\in \cJ_2}|\SpecPaths^*_s|$ edges will be used in Stage~E, $11\sum_{s\in\cJ_1}|\SpecPaths^*_s|$ edges will be used in Stage~F and $3\sum_{s\in\cJ_0}|\SpecShortPaths_s|$ edges will be used in Stage~G. This gives
\begin{align*}
e(H)&=(e(H)-e(\HStageD))\\
&+7\sum_{s\in \cJ_2}|\SpecPaths^*_s|+3\sum_{s\in\cJ_0}|\SpecShortPaths_s|+11\sum_{s\in\cJ_1}|\SpecPaths^*_s|\;.
\end{align*}
Plugging this back into~\eqref{horriblesoundsbybrucespringsteen}, we get
$$D_{\boxminus\boxminus}+D_{\boxplus\boxplus}+D_{\boxminus\boxplus}=7\sum_{s\in \cJ_2}|\SpecPaths^*_s|\;,$$
as was needed.

Next, we verify Definition~\ref{def:occupancies}\ref{en:occ:reqPar}. Here we use the parity correction established in Stage~B. We first consider the edges at $\Vmin$ and the paths anchored in $\Vmin$ after Stage~B. By~\ref{B:parity}, we have $\deg_{\HStageB}(v)\equiv\PathTerm(v)\mod 2$ for every $v\in V(\HStageB)$. Summing over $v\in\Vmin$, we have
\[\sum_{v\in\Vmin}\deg_{\HStageB}(v)\equiv\sum_{v\in\Vmin}\PathTerm(v)\mod 2\,.\]
Observe that on the left-hand side, every edge in $\Vmin$ is counted twice, and on the right-hand side every path with both terminals in $\Vmin$ is counted twice. Thus we see that the parity of the number of edges in $\HStageB$ leaving $\Vmin$ is equal to the parity of the number of paths in $\bigcup_{s\in\cJ}\SpecPaths_s$ with exactly one anchor in $\Vmin$; we call these \emph{crossing paths}. We claim that this property is maintained in Stages~C and~D: that is, the parity of the number of edges in $\HStageD$ leaving $\Vmin$ is equal to the parity of the number of unembedded crossing paths. More formally, we claim
\begin{equation}\label{eq:E:parityedges}
 e_{\HStageD}\big(\Vmin,\Vplus\big)\equiv \big|\bigcup_{s\in\cJ_0}\SpecShortPaths_s\big|+j_{\boxminus\boxplus}+\big|m^{-1}(\mathrm{odd})\big|\mod 2\,.
\end{equation}
To see that this is true, consider the embeddings in Stages~C and~D. In Stage~C, we embed some complete paths. When we embed a path with both anchors in $\Vmin$, or both anchors not in $\Vmin$, we do not change the number of crossing paths, and we necessarily use an even number of edges leaving $\Vmin$, so that the parities of the two quantities remain the same. When we embed a crossing path, the number of crossing paths drops by one and so the parity changes, but we use an odd number of edges leaving $\Vmin$ and so that parity also changes to match. Thus after Stage~C, the parity of the number of edges going from $\Vmin$ to $V(\HStageB)\setminus\Vmin$ is equal to the parity of the number of unembedded crossing paths. Note that after Stage~C, there are no edges of $\HStageC$ at $\boxdot$ (if it exists), nor are any unembedded paths anchored there. Thus we can ignore $\boxdot$, and we get that the parity of $e_{\HStageC}\big(\Vmin,\Vplus\big)$ is equal to the parity of the number of unembedded crossing paths, i.e.\ the number of paths in $\bigcup_{s\in\cJ}\SpecPaths_s^*$ with exactly one anchor in $\Vmin$. This last number is equal to the number of paths in $\bigcup_{s\in\cJ_0}\SpecPaths_s^*$ with exactly one anchor in $\Vmin$, plus $j_{\boxminus\boxplus}+|m^{-1}(\mathrm{odd})|$ (which account for the paths belonging to $\cJ_1\cup\cJ_2$).

In Stage~D, we embed parts of the paths $\bigcup_{s\in\cJ_0}\SpecPaths_s^*$. Consider some $P\in\SpecPaths^*_s$. Each of the two pieces of $P$ that we embed starts at an anchor of $P$ and terminates at an anchor of the corresponding path $P'$ in $\SpecShortPaths_s$. Note that $P'$ is a crossing path by construction. If both anchors of $P$ are in $\Vmin$, or both are in $\Vplus$, then when we embed the two pieces of $P$ we use in total an odd number of edges leaving $\Vmin$ (we use an even number in one piece and an odd number in the other) and increase the number of crossing paths by one; if $P$ is a crossing path, then we embed in total an even number of edges leaving $\Vmin$ (either both pieces use an even number of edges, or both pieces use an odd number of edges, leaving $\Vmin$) and the number of crossing paths does not change. This, with the observation that all paths in $\bigcup_{s\in\cJ_0}\SpecShortPaths_s$ are crossing, gives~\eqref{eq:E:parityedges}.

Recall that $D_{\boxminus\boxplus}$ is the right-hand side of~\eqref{dminplus} and hence
\[D_{\boxminus\boxplus}\equiv e_{\HStageD}(\Vmin,\Vplus)-\sigmJnula\sigmKJ^*n^2-j_{\boxminus\boxplus}\mod 2\,.\]
Since $\big|\bigcup_{s\in\cJ_0}\SpecShortPaths_s\big|=\sigmJnula\sigmKJ^*n^2$, we see from~\eqref{eq:E:parityedges} that $D_{\boxminus\boxplus}\equiv |m^{-1}(\mathrm{odd})|\mod 2$, as desired for~\ref{en:occ:reqPar}.

Finally, we verify Definition~\ref{def:occupancies}\ref{en:occ:balanced}. By Lemma~\ref{lem:indexquasi-implies-blockquasi}, we have $d\big(\HStageD[\Vmin]\big)=(1\pm2\gamAnchor)\dStageD$. Thus we have
\begin{align*}
D_{\boxminus\boxminus}&=e_{\HStageD}(\Vmin)-\left(\sigmJnula\sigmKJ^* n^2+2 j_{\boxminus\boxminus}+3(j_{\boxminus\boxplus}+ j_{\boxplus\boxplus})\right)\\
\JUSTIFY{by choice of $\sigmJjedna$}&=e_{\HStageD}(\Vmin)\pm 2\sigmJnula\sigmKJ^* n^2\\
&=\big(1\pm2\gamAnchor\big)\dStageD\binom{|\Vmin|}{2}\pm 2\sigmJnula\sigmKJ^* n^2\\
&=\dStageD
\cdot \frac{n^2}8 \pm 3\sigmJnula\sigmKJ^* n^2\eqByRef{eq:dStageD}
14\sigmKJ^2\cdot \frac{n^2}8 \pm 5\sigmJnula\sigmKJ n^2\;.
\end{align*}
(Expanding~\eqref{eq:dStageD} contributes a further $\sigmKJ\sigmJnula n^2$ to the error here: the $-8\sigmKJ\sigmJnula$ term inside its parentheses becomes exactly $-\sigmKJ\sigmJnula n^2$ once multiplied by $\frac{n^2}8$, while the remaining $\pm\iniquasi$ factor and the $\sigmJjedna$-terms of~\eqref{eq:dStageD} are negligible by comparison and are absorbed into the same bound.)
Similarly, we obtain $D_{\boxplus\boxplus}=14\sigmKJ^2\cdot \frac{n^2}8 \pm 5\sigmJnula\sigmKJ n^2$, and $D_{\boxminus\boxplus}=28\sigmKJ^2\cdot \frac{n^2}8 \pm 5\sigmJnula\sigmKJ n^2$. Definition~\ref{def:occupancies}\ref{en:occ:balanced} then follows from~\eqref{eq:CONSTANTS}.
\end{claimproof}

Having verified that we have a moderately balanced occupancy requirement, we can use Lemma~\ref{lem:OccAssign} to obtain an occupancy assignment. However, we cannot directly use the output of Lemma~\ref{lem:OccAssign} to decide the number of edges that each path of $\bigcup_{s\in\cJ_2}\SpecPaths_s^*$ will embed to each side and crossing: this output could be imbalanced between the different $s\in\cJ_2$ and lead to failure of~\ref{en:abstractpscatt} or~\ref{en:pairdistrpropE}. We put in an extra step and randomisation to deal with this as follows.

We split $\bigcup_{s\in \cJ_2}\SpecPaths^*_s$ into three sets $\mathcal{Q}_{\boxminus\boxminus}$, $\mathcal{Q}_{\boxminus\boxplus}$, and $\mathcal{Q}_{\boxplus\boxplus}$,
\begin{align*}
\mathcal{Q}_{\boxminus\boxminus}&:=\left\{P\in \bigcup_{s\in \cJ_2}\SpecPaths^*_s: \phi^\mathbf{A}(\leftpath_0(P)),\phi^\mathbf{A}(\rightpath_0(P))\in \Vmin\right\} \;,\\
\mathcal{Q}_{\boxplus\boxplus}&:=\left\{P\in \bigcup_{s\in \cJ_2}\SpecPaths^*_s: \phi^\mathbf{A}(\leftpath_0(P)),\phi^\mathbf{A}(\rightpath_0(P))\in \Vplus\right\} \;,\\
\mathcal{Q}_{\boxminus\boxplus}&:=\bigcup_{s\in \cJ_2}\SpecPaths^*_s\setminus (\mathcal{Q}_{\boxminus\boxminus}\cup \mathcal{Q}_{\boxplus\boxplus}) \;.
\end{align*}

We define auxiliary sets $X=X_{\boxminus\boxminus}\dcup X_{\boxminus\boxplus}\dcup X_{\boxplus\boxplus}$, where the sizes of $X_{\boxminus\boxminus}$, $X_{\boxminus\boxplus}$, and $X_{\boxplus\boxplus}$ are $|\mathcal{Q}_{\boxminus\boxminus}|$, $|\mathcal{Q}_{\boxminus\boxplus}|$, and $|\mathcal{Q}_{\boxplus\boxplus}|$, respectively. The reason we introduce this abstract set, rather than applying Lemma~\ref{lem:OccAssign} directly to $\bigcup_{s\in\cJ_2}\SpecPaths^*_s$, is that Lemma~\ref{lem:OccAssign} gives no control over how the resulting occupancy vectors correlate with the identity of the path they get assigned to; as noted above, we need this correlation to be as if uniformly random, so that it does not interact badly with $s$ when we later verify~\ref{en:abstractpscatt} and~\ref{en:pairdistrpropE}. So instead we first apply Lemma~\ref{lem:OccAssign} to the `unlabelled' set $X$, obtaining an occupancy assignment on $X$ with no reference yet to the actual paths, and only then, via an independent uniformly random bijection $\pi_{ab}$ from each $\mathcal{Q}_{ab}$ to the correspondingly sized $X_{ab}$ (chosen a few lines below), do we transfer this to the paths $\bigcup_{s\in\cJ_2}\SpecPaths^*_s$. It is this independent randomness in the bijection, not anything intrinsic to $X$, that gives the required lack of correlation.

We apply Lemma~\ref{lem:OccAssign} with the set $X$, the map sending each element of $X_{\boxminus\boxminus}\cup X_{\boxplus\boxplus}$ to $\mathrm{even}$ and each element of $X_{\boxminus\boxplus}$ to $\mathrm{odd}$ (this is simply the map $m$ defined above, transferred from $\bigcup_{s\in\cJ_2}\SpecPaths^*_s$ to $X$ via the bijections $\pi_{ab}$ chosen below), and the tuple $(D_{\boxminus\boxminus},D_{\boxminus\boxplus},D_{\boxplus\boxplus})$. Let $\aleph$ be the occupancy assignment we obtain.

We pick a uniform random (and independent of $\aleph$) partition of $X$ into sets $X^\triangledown$ and $X^\blacktriangle$.
Let $\pi_{\boxminus\boxminus}:\mathcal{Q}_{\boxminus\boxminus}\rightarrow X_{\boxminus\boxminus}$, $\pi_{\boxminus\boxplus}:\mathcal{Q}_{\boxminus\boxplus}\rightarrow X_{\boxminus\boxplus}$, and $\pi_{\boxplus\boxplus}:\mathcal{Q}_{\boxplus\boxplus}\rightarrow X_{\boxplus\boxplus}$ be independent uniformly random bijections. We write $\pi:\bigcup_{s\in \cJ_2}\SpecPaths^*_s\rightarrow X$ for the union of these three bijections.

We can now define $\chi$ on each set $V(P)$, $P\in \SpecPaths^*_s$, $s\in \cJ_2$ as follows. Let $P=x_1x_2\cdots x_8$, where in case $P\in \mathcal{Q}_{\boxminus\boxminus}\cup\mathcal{Q}_{\boxplus\boxplus}$ we take $x_1=\leftpath_0(P)$ and in case $P\in \mathcal{Q}_{\boxminus\boxplus}$ we take the orientation of $P$ such that $\phi^\mathbf{A}(x_1) \in \Vmin$. If $\pi(P)\in X_{\boxminus\boxminus}\cap X^\triangledown$ then we define $\chi_{\restriction V(P)}$ according to Table~\ref{tab:VminVminL}, whereas if $\pi(P)\in X_{\boxminus\boxminus}\cap X^\blacktriangle$, we define $\chi_{\restriction V(P)}$ according to Table~\ref{tab:VminVminR}.
\begin{table}[!ht]
\begin{tabular}{|c|c|c|c|}
	\hline
	value of $\aleph(\pi_{\boxminus\boxminus}(P))$ & the map $\chi$ on $V(P)$ & left pattern & right pattern\\ \hline
	$(1,2,4)$ &  $x_1,x_2,x_8\mapsto \Vmin$, $x_3,x_4,\ldots,x_7\mapsto \Vplus$ & $(\boxminus,\boxminus,\boxplus)$ & $(\boxminus,\boxplus,\boxplus)$\\
	$(1,4,2)$ &  $x_1,x_2,x_4,x_8\mapsto \Vmin$, $x_3,x_5,x_6,x_7\mapsto \Vplus$ & $(\boxminus,\boxminus,\boxplus)$ & $(\boxminus,\boxplus,\boxplus)$\\
	$(2,2,3)$ &  $x_1,x_2,x_3,x_8\mapsto \Vmin$, $x_4,x_5,x_6,x_7\mapsto \Vplus$ &$(\boxminus,\boxminus,\boxminus)$ & $(\boxminus,\boxplus,\boxplus)$\\
	$(2,4,1)$ &  $x_1,x_2,x_3,x_5,x_8\mapsto \Vmin$, $x_4,x_6,x_7\mapsto \Vplus$ & $(\boxminus,\boxminus,\boxminus)$ & $(\boxminus,\boxplus,\boxplus)$\\
	$(3,2,2)$ &  $x_1,x_2,x_3,x_4,x_8\mapsto \Vmin$, $x_5,x_6,x_7\mapsto \Vplus$& $(\boxminus,\boxminus,\boxminus)$ & $(\boxminus,\boxplus,\boxplus)$\\
	$(4,2,1)$ &  $x_1,x_2,\ldots,x_5,x_8\mapsto \Vmin$, $x_6,x_7\mapsto \Vplus$& $(\boxminus,\boxminus,\boxminus)$ & $(\boxminus,\boxplus,\boxplus)$
	\\ \hline
\end{tabular}
\caption{Defining $\chi_{\restriction V(P)}$ when $\pi(P)\in X_{\boxminus\boxminus}\cap X^\triangledown$.}
\label{tab:VminVminL}
\end{table}

\begin{table}[!ht]
\begin{tabular}{|c|c|c|c|}
	\hline
	value of $\aleph(\pi_{\boxminus\boxminus}(P))$ & the map $\chi$ on $V(P)$ & left pattern & right pattern\\ \hline
	$(1,2,4)$ &  $x_1,x_7,x_8\mapsto \Vmin$, $x_2, x_3,\ldots,x_6\mapsto \Vplus$ & $(\boxminus,\boxplus,\boxplus)$ & $(\boxminus,\boxminus,\boxplus)$\\
	$(1,4,2)$ &  $x_1,x_5, x_7,x_8\mapsto \Vmin$, $x_2,x_3,x_4,x_6\mapsto \Vplus$ & $(\boxminus,\boxplus,\boxplus)$ & $(\boxminus,\boxminus,\boxplus)$\\
	$(2,2,3)$ &  $x_1,x_6,x_7,x_8\mapsto \Vmin$, $x_2,x_3,x_4,x_5\mapsto \Vplus$ &$(\boxminus,\boxplus,\boxplus)$ & $(\boxminus,\boxminus,\boxminus)$\\
	$(2,4,1)$ &  $x_1,x_4, x_6, x_7,x_8\mapsto \Vmin$, $x_2,x_3,x_5\mapsto \Vplus$ & $(\boxminus,\boxplus,\boxplus)$ & $(\boxminus,\boxminus,\boxminus)$\\
	$(3,2,2)$ &  $x_1,x_5,x_6,x_7,x_8\mapsto \Vmin$, $x_2,x_3,x_4\mapsto \Vplus$& $(\boxminus,\boxplus,\boxplus)$ & $(\boxminus,\boxminus,\boxminus)$\\
	$(4,2,1)$ &  $x_1,x_4,\ldots,x_8\mapsto \Vmin$, $x_2,x_3\mapsto \Vplus$& $(\boxminus,\boxplus,\boxplus)$ & $(\boxminus,\boxminus,\boxminus)$
	\\ \hline
\end{tabular}
\caption{Defining $\chi_{\restriction V(P)}$ when $\pi(P)\in X_{\boxminus\boxminus}\cap X^\blacktriangle$. Note that the table is the same as Table~\ref{tab:VminVminL} if the paths are followed in the reversed direction.}
\label{tab:VminVminR}
\end{table}

The point of Tables~\ref{tab:VminVminL} and~\ref{tab:VminVminR} is that they satisfy requirement~\ref{en:kaficko} of the lemma, and the numbers of edges of $P$ that are mapped into $\Vmin$, between $\Vmin$ and $\Vplus$, and into $\Vplus$, are $\aleph(\pi_{\boxminus\boxminus}(P))_1$, $\aleph(\pi_{\boxminus\boxminus}(P))_2$, and $\aleph(\pi_{\boxminus\boxminus}(P))_3$, respectively.
Let us now explain the `left pattern' and the `right pattern' columns in Tables~\ref{tab:VminVminL} and~\ref{tab:VminVminR}. The left pattern is the sequence which says where the three leftmost vertices of $P$ are mapped. Note that this pattern can be read off from the map $\chi$, but it is convenient to have it in this form for verification of~\ref{en:abstractpscatt} and~\ref{en:pairdistrpropE}. The right pattern is the sequence which says where the three rightmost vertices of $P$ are mapped (starting with $x_8$).

By symmetry, there is a similar assignment $\chi_{\restriction V(P)}$ if $\pi(P)\in X_{\boxplus\boxplus}\cap X^\triangledown$ and $\pi(P)\in X_{\boxplus\boxplus}\cap X^\blacktriangle$, which for completeness we give in Table~\ref{tab:VplusVplusL} and Table~\ref{tab:VplusVplusR}, respectively.
\begin{table}[!ht]
	\begin{tabular}{|c|c|c|c|}
		\hline
		value of $\aleph(\pi_{\boxplus\boxplus}(P))$ & the map $\chi$ on $V(P)$ & left pattern & right pattern\\ \hline
		$(4,2,1)$ &  $x_1,x_2,x_8\mapsto \Vplus$, $x_3,x_4,\ldots,x_7\mapsto \Vmin$ & $(\boxplus,\boxplus,\boxminus)$ & $(\boxplus,\boxminus,\boxminus)$\\
		$(2,4,1)$ &  $x_1,x_2,x_4,x_8\mapsto \Vplus$, $x_3,x_5,x_6,x_7\mapsto \Vmin$ & $(\boxplus,\boxplus,\boxminus)$ & $(\boxplus,\boxminus,\boxminus)$\\
		$(3,2,2)$ &  $x_1,x_2,x_3,x_8\mapsto \Vplus$, $x_4,x_5,x_6,x_7\mapsto \Vmin$ &$(\boxplus,\boxplus,\boxplus)$ & $(\boxplus,\boxminus,\boxminus)$\\
		$(1,4,2)$ &  $x_1,x_2,x_3,x_5,x_8\mapsto \Vplus$, $x_4,x_6,x_7\mapsto \Vmin$ & $(\boxplus,\boxplus,\boxplus)$ & $(\boxplus,\boxminus,\boxminus)$\\
		$(2,2,3)$ &  $x_1,x_2,x_3,x_4,x_8\mapsto \Vplus$, $x_5,x_6,x_7\mapsto \Vmin$& $(\boxplus,\boxplus,\boxplus)$ & $(\boxplus,\boxminus,\boxminus)$\\
		$(1,2,4)$ &  $x_1,x_2,\ldots,x_5,x_8\mapsto \Vplus$, $x_6,x_7\mapsto \Vmin$& $(\boxplus,\boxplus,\boxplus)$ & $(\boxplus,\boxminus,\boxminus)$
		\\ \hline
	\end{tabular}
	\caption{Defining $\chi_{\restriction V(P)}$ when $\pi(P)\in X_{\boxplus\boxplus}\cap X^\triangledown$.}
	\label{tab:VplusVplusL}
\end{table}

\begin{table}[!ht]
	\begin{tabular}{|c|c|c|c|}
		\hline
		value of $\aleph(\pi_{\boxplus\boxplus}(P))$ & the map $\chi$ on $V(P)$ & left pattern & right pattern\\ \hline
		$(4,2,1)$ &  $x_1,x_7,x_8\mapsto \Vplus$, $x_2, x_3,\ldots,x_6\mapsto \Vmin$ & $(\boxplus,\boxminus,\boxminus)$ & $(\boxplus,\boxplus,\boxminus)$\\
		$(2,4,1)$ &  $x_1,x_5, x_7,x_8\mapsto \Vplus$, $x_2,x_3,x_4,x_6\mapsto \Vmin$ & $(\boxplus,\boxminus,\boxminus)$ & $(\boxplus,\boxplus,\boxminus)$\\
		$(3,2,2)$ &  $x_1,x_6,x_7,x_8\mapsto \Vplus$, $x_2,x_3,x_4,x_5\mapsto \Vmin$ &$(\boxplus,\boxminus,\boxminus)$ & $(\boxplus,\boxplus,\boxplus)$\\
		$(1,4,2)$ &  $x_1,x_4, x_6, x_7,x_8\mapsto \Vplus$, $x_2,x_3, x_5\mapsto \Vmin$ & $(\boxplus,\boxminus,\boxminus)$ & $(\boxplus,\boxplus,\boxplus)$\\
		$(2,2,3)$ &  $x_1,x_5,x_6,x_7,x_8\mapsto \Vplus$, $x_2,x_3,x_4\mapsto \Vmin$& $(\boxplus,\boxminus,\boxminus)$ & $(\boxplus,\boxplus,\boxplus)$\\
		$(1,2,4)$ &  $x_1,x_4,\ldots,x_8\mapsto \Vplus$, $x_2,x_3\mapsto \Vmin$& $(\boxplus,\boxminus,\boxminus)$ & $(\boxplus,\boxplus,\boxplus)$
		\\ \hline
	\end{tabular}
	\caption{Defining $\chi_{\restriction V(P)}$ when $\pi(P)\in X_{\boxplus\boxplus}\cap X^\blacktriangle$. As with Tables~\ref{tab:VminVminL} and~\ref{tab:VminVminR}, the paths in this table are reversed copies from Table~\ref{tab:VplusVplusL}.}
	\label{tab:VplusVplusR}
\end{table}

Next, we define $\chi_{\restriction V(P)}$ in the case when $P\in \mathcal{Q}_{\boxminus\boxplus}$. To this end we use Tables~\ref{tab:VminVplus} and~\ref{tab:VminVplus2}.
\begin{table}[!ht]
\begin{tabular}{|c|c|c|c|}
	\hline
	value of $\aleph(\pi_{\boxminus\boxplus}(P))$ & the map $\chi$ on $V(P)$&left pattern & right pattern\\ \hline
	$(1,3,3)$ & $x_1,x_2,x_4\mapsto \Vmin$, $x_3,x_5,x_6,x_7,x_8\mapsto \Vplus$ & $(\boxminus,\boxminus,\boxplus)$
	&$(\boxplus,\boxplus,\boxplus)$\\
	$(1,5,1)$ & $x_1,x_2,x_4,x_6\mapsto \Vmin$, $x_3,x_5,x_7,x_8\mapsto \Vplus$&$(\boxminus,\boxminus,\boxplus)$&
	$(\boxplus,\boxplus,\boxminus)$\\
	$(2,3,2)$ & $x_1,x_2,x_3,x_5\mapsto \Vmin$, $x_4,x_6,x_7,x_8\mapsto \Vplus$&$(\boxminus,\boxminus,\boxminus)$
	&$(\boxplus,\boxplus,\boxplus)$\\
	$(3,3,1)$ & $x_1,\ldots,x_4,x_6\mapsto \Vmin$, $x_5,x_7,x_8\mapsto \Vplus$&$(\boxminus,\boxminus,\boxminus)$
	&$(\boxplus,\boxplus,\boxminus)$\\
	\hline
\end{tabular}
\caption{Defining $\chi_{\restriction V(P)}$ when $\pi(P)\in X_{\boxminus\boxplus}\cap X^\triangledown$.}
\label{tab:VminVplus}
\end{table}

\begin{table}[!ht]
	\begin{tabular}{|c|c|c|c|}
		\hline
		value of $\aleph(\pi_{\boxminus\boxplus}(P))$ & the map $\chi$ on $V(P)$&left pattern & right pattern\\ \hline
		$(1,3,3)$ & $x_1,x_6,x_7\mapsto \Vmin$, $x_2,\ldots,x_5,x_8\mapsto \Vplus$ & $(\boxminus,\boxplus,\boxplus)$
		&$(\boxplus,\boxminus,\boxminus)$\\
		$(1,5,1)$ & $x_1,x_4,x_6, x_7\mapsto \Vmin$, $x_2,x_3,x_5,x_8\mapsto \Vplus$&$(\boxminus,\boxplus,\boxplus)$&
		$(\boxplus,\boxminus,\boxminus)$\\
		$(2,3,2)$ & $x_1,x_5,x_6,x_7\mapsto \Vmin$, $x_2,x_3,x_4,x_8\mapsto \Vplus$&$(\boxminus,\boxplus,\boxplus)$
		&$(\boxplus,\boxminus,\boxminus)$\\
		$(3,3,1)$ & $x_1,x_4, \ldots,x_7\mapsto \Vmin$, $x_2,x_3,x_8\mapsto \Vplus$&$(\boxminus,\boxplus,\boxplus)$
		&$(\boxplus,\boxminus,\boxminus)$\\
		\hline
	\end{tabular}
	\caption{Defining $\chi_{\restriction V(P)}$ when $\pi(P)\in X_{\boxminus\boxplus}\cap X^\blacktriangle$.}
	\label{tab:VminVplus2}
\end{table}

Again, it is obvious that Tables~\ref{tab:VminVplus} and~\ref{tab:VminVplus2} obey requirement~\ref{en:kaficko} of the lemma, and the numbers of edges of $P$ that are mapped into $\Vmin$, between $\Vmin$ and $\Vplus$, and into $\Vplus$, are $\aleph(\pi_{\boxminus\boxplus}(P))_1$, $\aleph(\pi_{\boxminus\boxplus}(P))_2$, and $\aleph(\pi_{\boxminus\boxplus}(P))_3$, respectively. Thus requirement~\ref{en:Etotaledges} is inherited from Definition~\ref{def:occupancies}\ref{en:occ:sumN}.

By inspecting Tables~\ref{tab:VminVminL} and~\ref{tab:VminVminR}, we observe that in all entries one of $x_2$ and $x_7$ is allocated to $\Vmin$, and the other is allocated to $\Vplus$. Thus each path of $\cQ_{\boxminus\boxminus}\cap\SpecPaths^*_s$ contributes at least one vertex to $X^a_s$ for each $a\in\{\boxminus,\boxplus\}$. Since each path of $\cQ_{\boxminus\boxminus}$ accounts for two of the anchors counted by $A^\boxminus_s\cap V_\boxminus$, we have by~\ref{quasi:anchorsets} that
\[\big|\cQ_{\boxminus\boxminus}\cap\SpecPaths^*_s\big|=\tfrac12\big|A^\boxminus_s\cap V_\boxminus\big|=(1\pm\gamAnchor)\tfrac14\big(\sigmKJ n\pm n^{0.6}\big)\ge\tfrac15\sigmKJ n\,,\]
and we conclude~\ref{en:medunkovycaj}.

Finally, inspecting Table~\ref{tab:VminVplus2}, we see that in all entries in this table the edge $x_6x_7$ is assigned to $\Vmin$, and $x_2x_3$ is assigned to $\Vplus$, and (since at least three edges are assigned between $\Vmin$ and $\Vplus$) there is an edge other than $x_1x_2$ and $x_7x_8$ which is assigned between $\Vmin$ and $\Vplus$. Thus each path $P$ of $\SpecPaths^*_s$ which is in $\cQ_{\boxminus\boxplus}$ and has $\pi(P)\in X^\blacktriangle$ contributes at least one pair to each of the three sets considered in~\ref{en:nocheins} (observe that the case $a=\boxminus,b=\boxplus$ is identical to $a=\boxplus,b=\boxminus$ by symmetry). As above, by~\ref{quasi:anchorsets} we have $\big|\cQ_{\boxminus\boxplus}\cap\SpecPaths^*_s\big|\ge\tfrac14\sigmKJ n$, and by Chernoff's inequality (Fact~\ref{fact:hypergeometricBasicProperties}), with probability at least $1-\exp(-n^{0.5})$ at least $\tfrac1{16}\sigmKJ n$ of these paths $P$ have $\pi(P)\in X^\blacktriangle$. Taking the union bound over the at most $n$ choices of $s$, we conclude~\ref{en:nocheins}.

\smallskip

What remains is to check~\ref{en:abstractpscatt} and~\ref{en:pairdistrpropE}. 
To this end, we use that $\pi_{\boxminus\boxminus}$, $\pi_{\boxminus\boxplus}$ and $\pi_{\boxplus\boxplus}$ are random, and that we split $X=X^\triangledown\dcup X^\blacktriangle$ randomly.

Let us first focus on proving~\ref{en:abstractpscatt}, which is the following claim.
\begin{claim}\label{cl:E:anchordist}
With high probability for every $s\in \cJ_2$, the pair
\[\left(\phi^\mathbf{D}_{\restriction\textrm{endvertices of $\SpecPaths^*_s$}},
\chi_{\restriction\textrm{inner vertices of $\SpecPaths^*_s$}}\right)\]
has the $6\gamAnchor$-anchor distribution property with respect to the path-forest $\SpecPaths^*_s$, partition $(\Vmin,\Vplus)$ and the graph $\HStageD$.
\end{claim}
\begin{claimproof}[Proof of Claim~\ref{cl:E:anchordist}]
We fix $s\in\cJ_2$, and aim to prove that the $6\gamAnchor$-anchor distribution property holds for $\SpecPaths^*_s$ with sufficiently high probability to take a union bound over the at most $n$ choices of $s$.

As in Definition~\ref{def:anchordistribution}, given $a,b,c\in\{\boxminus,\boxplus\}$ let $A_{a,b,c}$ denote the collection of those endvertices $x$ of paths of $\SpecPaths^*_s$ such that $\phi^\mathbf{D}_s(x)\in V_a$, and such that the neighbour $y$ of $x$ satisfies $\chi(y) = V_b$, and the next vertex $z$ satisfies $\chi(z) =  V_c$.
For each $v\in V_b\setminus \imD(s)$ we claim that with probability at least $1-\exp(-n^{0.4})$ we have
\begin{equation}\label{eq:provemeanchor}
\big|\big\{x\in A_{a,b,c}:\phi^\mathbf{D}_s(x)\in \NBH_{\HStageD}(v)\big\}\big|=(1\pm 3\gamAnchor)\dStageD |A_{a,b,c}|\pm \gamAnchor n^{0.99}\,.
\end{equation}
Note that in this setting all the densities in Definition~\ref{def:anchordistribution} are $(1\pm2\gamAnchor)\dStageD$ by the fact that from Lemma~\ref{lem:StageDNew} we have a $(\gamAnchor,\LD,\dStageD,\dStageD)$-quasirandom setup, and using Lemma~\ref{lem:indexquasi-implies-blockquasi} with~\ref{quasi:indexquasi}, so that establishing~\eqref{eq:provemeanchor} for each pattern $a,b,c$ and each $v\in V_b\setminus\imD(s)$ gives the required $6\gamAnchor$-anchor distribution property. To simplify notation, we suppose $a=\boxminus$; the other case is symmetric.

Observe that $|A_{\boxminus,b,c}|$ is the total number of left and right patterns (as in Tables~\ref{tab:VminVminL}--\ref{tab:VminVplus2}) whose consecutive values are $\boxminus$, $b$, and $c$ which are assigned to paths of $\SpecPaths^*_s$ by $\chi$. This can be either a left or a right pattern of a path from $\cQ_{\boxminus\boxminus}$, or a left pattern of a path from $\cQ_{\boxminus\boxplus}$. We split up $A_{\boxminus,b,c}$ into sets $A^{p}_{\boxminus,b,c}$, where $p\in\{\boxminus,\boxplus\}$, and a vertex $x$ of $A^{p}_{\boxminus,b,c}$ is in a path of $\SpecPaths^*_s$ whose other endvertex is anchored in $V_p$.

Recall that the partition $X=X^\triangledown\dcup X^\blacktriangle$ and the bijection $\pi_{\boxminus\boxminus}$ were chosen randomly. Whatever the partition $X=X^\triangledown\dcup X^\blacktriangle$ turns out to be, it, together with $\aleph$, determines the total number $N_\boxminus$ of left and right patterns of paths in $\cQ_{\boxminus\boxminus}$ which are $\boxminus,b,c$, and also the total number $N_\boxplus$ of left patterns of paths in $\cQ_{\boxminus\boxplus}$ which are $\boxminus,b,c$. Note that no path of $\cQ_{\boxminus\boxminus}$ has its left and right pattern equal: in Table~\ref{tab:VminVminL} the right pattern is $(\boxminus,\boxplus,\boxplus)$ in every row while the left pattern always assigns $x_2$ to $\Vmin$, and symmetrically in Table~\ref{tab:VminVminR}. Hence $N_\boxminus$ is also the number of elements of $X_{\boxminus\boxminus}$ carrying a pattern $\boxminus,b,c$, and each path whose $\pi_{\boxminus\boxminus}$-image is such an element contributes exactly one vertex to $A^\boxminus_{\boxminus,b,c}$. So, when we reveal $\pi_{\boxminus\boxminus}$, the size $\big|A_{\boxminus,b,c}^\boxminus\big|$ follows the hypergeometric distribution with parameters $\big(\big|\cQ_{\boxminus \boxminus}\big|,N_\boxminus,\big|\cQ_{\boxminus \boxminus}\cap\SpecPaths^*_s\big|\big)$. (Recall that this is the random experiment of choosing, from a set of size $\big|\cQ_{\boxminus \boxminus}\big|$ of which $N_\boxminus$ are marked, a subset uniformly at random of size $\big|\cQ_{\boxminus \boxminus}\cap\SpecPaths^*_s\big|$, and seeing how many are marked.) For $\big|A_{\boxminus,b,c}^\boxminus\cap \NBH_{\HStageD}(v)\big|$ we argue in the same way, splitting according to whether the pattern $\boxminus,b,c$ occurs as a left or as a right pattern: writing $N_\boxminus=N^{\mathrm{L}}_\boxminus+N^{\mathrm{R}}_\boxminus$ accordingly, and $A_{v,\boxminus}=A^{\mathrm{L}}_{v,\boxminus}+A^{\mathrm{R}}_{v,\boxminus}$ for the number of vertices in $V_\boxminus\cap \NBH_{\HStageD}(v)$ which are anchors of paths in $\cQ_{\boxminus\boxminus}\cap\SpecPaths_s^*$, split according to whether the anchor is the image of $\leftpath_0$ or of $\rightpath_0$ of its path, each of the two corresponding counts is hypergeometric, with means $N^{\mathrm{L}}_\boxminus A^{\mathrm{L}}_{v,\boxminus}\big/\big|\cQ_{\boxminus\boxminus}\big|$ and $N^{\mathrm{R}}_\boxminus A^{\mathrm{R}}_{v,\boxminus}\big/\big|\cQ_{\boxminus\boxminus}\big|$. Now Tables~\ref{tab:VminVminL} and~\ref{tab:VminVminR} are mirror images of each other, and the choice between them is made by the uniform random partition $X=X^\triangledown\dcup X^\blacktriangle$; hence $N^{\mathrm{L}}_\boxminus$ and $N^{\mathrm{R}}_\boxminus$ are the two sides of a uniform random split of a set determined by $\aleph$ alone, so that, by Fact~\ref{fact:hypergeometricBasicProperties}, $N^{\mathrm{L}}_\boxminus=\tfrac12N_\boxminus\pm n^{1.3}$ (and hence also $N^{\mathrm{R}}_\boxminus=\tfrac12N_\boxminus\mp n^{1.3}$) with probability at least $1-\exp(-\sqrt n)$. Since $A^{\mathrm{L}}_{v,\boxminus},A^{\mathrm{R}}_{v,\boxminus}\le n$ while $\big|\cQ_{\boxminus\boxminus}\big|=\Theta(n^2)$, the sum of the two means is therefore $N_\boxminus A_{v,\boxminus}\big/2\big|\cQ_{\boxminus\boxminus}\big|\pm n^{0.4}$. By Fact~\ref{fact:hypergeometricBasicProperties}, with probability at least $1-\exp(-\sqrt{n})$ we get
\[\big|A_{\boxminus,b,c}^\boxminus\big|=\frac{N_\boxminus\cdot \big|\cQ_{\boxminus \boxminus}\cap\SpecPaths^*_s\big|}{\big|\cQ_{\boxminus \boxminus}\big|}\pm n^{0.9}\quad\text{and}\quad \big|A_{\boxminus,b,c}^\boxminus\cap \NBH_{\HStageD}(v)\big|=\frac{N_\boxminus\cdot A_{v,\boxminus}}{2\big|\cQ_{\boxminus \boxminus}\big|}\pm 3n^{0.9}\,.\]
By similar logic, recalling that paths of $\cQ_{\boxminus\boxplus}$ have only one end anchored in $\Vmin$, we get
\[\big|A_{\boxminus,b,c}^\boxplus\big|=\frac{N_\boxplus\cdot \big|\cQ_{\boxminus \boxplus}\cap\SpecPaths^*_s\big|}{\big|\cQ_{\boxminus \boxplus}\big|}\pm n^{0.9}\quad\text{and}\quad \big|A_{\boxminus,b,c}^\boxplus\cap \NBH_{\HStageD}(v)\big|=\frac{N_\boxplus\cdot A_{v,\boxplus}}{\big|\cQ_{\boxminus \boxplus}\big|}\pm n^{0.9}\,.\]

Suppose that all four of these equalities hold. By~\ref{quasi:distributionOfAnchors} (specifically~\eqref{eq:quasi:distAnchorsa1} or~\eqref{eq:quasi:distAnchorsa2}, depending on whether $v\in\Vmin$ or $\Vplus$), we have $A_{v,\boxminus}=(1\pm\gamAnchor)\dStageD\cdot 2\big|\cQ_{\boxminus \boxminus}\cap\SpecPaths^*_s\big|$, and similarly $A_{v,\boxplus}=(1\pm\gamAnchor)\dStageD\cdot \big|\cQ_{\boxminus \boxplus}\cap\SpecPaths^*_s\big|$. Substituting these into the above equations, we get
\begin{align*}
 \big|A_{\boxminus,b,c}^\boxminus\cap \NBH_{\HStageD}(v)\big|&=(1\pm\gamAnchor)\dStageD\big|A_{\boxminus,b,c}^\boxminus\big|\pm n^{0.95}\\
 \text{and}\quad\big|A_{\boxminus,b,c}^\boxplus\cap \NBH_{\HStageD}(v)\big|&=(1\pm\gamAnchor)\dStageD\big|A_{\boxminus,b,c}^\boxplus\big|\pm n^{0.95}\,.
\end{align*}

Summing these equations, since $A_{\boxminus,b,c}=A_{\boxminus,b,c}^\boxminus\dcup A_{\boxminus,b,c}^\boxplus$, we obtain the desired~\eqref{eq:provemeanchor}. Taking the union bound over the four equalities and the choices of $a,b,c,v$ we see that it holds for all these choices with probability at least $1-\exp(-n^{0.4})$ as required. Finally the union bound over the choices of $s$ completes the proof of the claim.
\end{claimproof}

Finally, we shall prove that with high probability the map $\chi$ has the $20\gamAnchor$-pair distribution property, as required by~\ref{en:pairdistrpropE} in Lemma~\ref{lem:E2decideParts}. 

Recall that the graphs $(G_s)_{s\in\cJ_2}$ were (partially) packed in Stage~A and Stage~C. For such an $s\in \cJ_2$, the embedding $\phi^\mathbf{A}_s$ used exactly $(1-\deltnonspanning-6\sigmKJ)n$ vertices of $H$ by~\eqref{eq:familyOrderJ2Equality}, so that after Stage~A there were $(\deltnonspanning+6\sigmKJ)n$ vertices of $H$ left unused. In Stage~C, between $0$ and $6n^{0.6}$ additional vertices of $G_s$ were embedded. Hence, 
\begin{equation}\label{eq:D-sizeOfFreeSpace}
|V(H)\setminus \imD(s)|=(\deltnonspanning+6\sigmKJ)n\pm 6n^{0.6}\;.
\end{equation}
Similarly, the number of endvertices of paths in $\SpecPaths_s^*$ is
\begin{equation}\label{eq:sizeOfAs}
2|\SpecPaths^*_s|=2\sigmKJ n\pm 2n^{0.6} \;.
\end{equation}

Recall that $E^*:=\bigcup_{s\in \cJ_2} E(\SpecPaths^*_s)$, and let $E^{**}\subset E^*$ be the edges that do not contain endvertices (i.e.\ the $5$ edges in the middle of each path). We refer to the edges of $E^{**}$ as \emph{internal} and those of $E^*\setminus E^{**}$ as \emph{peripheral}. For $a,b\in \{\boxminus,\boxplus\}$ let
\begin{equation}\label{eq:f}
f_{ab}:=\sum_{(x,y):xy\in E^{**}}\ONE_{\chi(x)=V_a,\chi(y)=V_b}\;,
\end{equation}
where, as in Definition~\ref{def:pairdistprop} of the pair distribution property, an edge $xy$ with $\chi(x)=\chi(y)$ counts twice in this sum, once as $(x,y)$ and once as $(y,x)$. The same applies to subsequent similar sums in this section.
For $s\in\cJ_2$, let
\begin{equation}\label{eq:fs}
f_{ab,s}:=\sum_{(x,y):xy\in E^{**}\cap E(\SpecPaths^*_s)}\ONE_{\chi(x)=V_a,\chi(y)=V_b}\;.
\end{equation}
We define $B_{\boxminus\boxminus}:=2D_{\boxminus\boxminus}$, and $B_{\boxplus\boxplus}:=2D_{\boxplus\boxplus}$, and $B_{\boxminus\boxplus},B_{\boxplus\boxminus}:=D_{\boxminus\boxplus}$. Thus $B_{ab}$ is the number of edges in total, counted with multiplicity as in the above sums, that we want to embed in this stage with one end in $V_a$ and the other in $V_b$. Recall that the $D_{ab}$ form a moderately balanced occupancy requirement for $\bigcup_{s\in\cJ_2}\SpecPaths^*_s$, a set of size $|X|$, so by Definition~\ref{def:occupancies}\ref{en:occ:balanced} we have $D_{\boxminus\boxminus},D_{\boxplus\boxplus}\ge|X|$ and $D_{\boxminus\boxplus}\ge3|X|$; hence $B_{\boxminus\boxminus}=2D_{\boxminus\boxminus}\ge2|X|$, $B_{\boxplus\boxplus}=2D_{\boxplus\boxplus}\ge2|X|$, and $B_{\boxminus\boxplus}=B_{\boxplus\boxminus}=D_{\boxminus\boxplus}\ge3|X|>2|X|$, so in every case $B_{ab}\ge2|X|$. Since $|X|=\sum_{s\in\cJ_2}\big|\SpecPaths^*_s\big|\ge|\cJ_2|\big(\sigmKJ n-n^{0.6}\big)\ge\tfrac78\sigmKJ n|\cJ_2|$, we in particular conclude
\begin{equation}\label{eq:E:pairdist:lowerB}
 B_{ab}\ge \tfrac74\sigmKJ n|\cJ_2|=\tfrac74\sigmKJ (\sigmKJ-\sigmJnula-\sigmJjedna)n^2\,,
\end{equation}
using~\eqref{eq:sizeJ2} for the final equality.
We now prove the following claim.
\begin{claim}\label{claim:chacha}
For $s\in \cJ_2$, with probability at least $1-2\exp(-\sqrt{n})$ we have for $a,b\in \{\boxminus, \boxplus\}$
\begin{align*}
f_{ab,s}&=(1\pm9\gamAnchor)\left(\frac{B_{ab}}{(\sigmKJ-\sigmJnula-\sigmJjedna)n}-\sigmKJ n\right)\;.
\end{align*}
\end{claim}
\begin{claimproof}[Proof of Claim~\ref{claim:chacha}]
For $a,b\in \{\boxminus, \boxplus\}$, we define $B_{ab,s}$ as the total number of pairs $(x,y)$, $xy\in E(\SpecPaths_s^*)$ with $\chi(x)=V_a$ and $\chi(y)=V_b$. We claim that with probability at least $1-\exp(-\sqrt{n})$, for every $s\in \cJ_2$ we have
\begin{equation}\label{eq:Dabs}
B_{ab,s}=(1\pm 4\gamAnchor)\frac{B_{ab}}{(\sigmKJ-\sigmJnula-\sigmJjedna)n}\ge \tfrac32\sigmKJ n\;.
\end{equation}
The inequality in~\eqref{eq:Dabs} is an immediate consequence of~\eqref{eq:E:pairdist:lowerB} and the choice of $\gamAnchor$, so the difficulty is to prove the equality. To prove the equality in~\eqref{eq:Dabs}, fix $a,b\in\{\boxminus,\boxplus\}$ and $s\in\cJ_2$ and let $\mathfrak{v}$ be in the range of $\aleph$. By definition for each $c,d\in\{\boxminus,\boxplus\}$ there are $\big|\aleph^{-1}(\mathfrak{v})\cap X_{cd}\big|$ elements of $X_{cd}$ which are assigned $\mathfrak{v}$. The number of these elements which are in $\pi_{cd}\big(\SpecPaths^*_s\cap\cQ_{cd}\big)$ follows the hypergeometric distribution with parameters
\[\big(|X_{cd}|\;,\;|X_{cd}\cap\aleph^{-1}(\mathfrak{v})|\;,\;|\SpecPaths^*_s\cap \mathcal{Q}_{cd}|\big)\]
and by Fact~\ref{fact:hypergeometricBasicProperties}, with probability at least $1-\exp(-n^{0.6})$ we have
\begin{equation}\label{eq:E:pairdist:chacha}
 \big|\aleph^{-1}(\mathfrak{v})\cap X_{cd}\cap\pi_{cd}\big(\SpecPaths_s^*\cap\cQ_{cd}\big)\big|=\frac{|\SpecPaths^*_s\cap \mathcal{Q}_{cd}||X_{cd}\cap\aleph^{-1}(\mathfrak{v})|}{|X_{cd}|}\pm n^{0.9}\,.
\end{equation}
Suppose that~\eqref{eq:E:pairdist:chacha} holds for each choice of $\mathfrak{v}$ and of $c,d,s$; there are at most $100n$ such choices in total, so that this occurs with probability at least $1-\exp(-\sqrt{n})$.

Suppose for a moment that $a=b=\boxminus$. Consider the following weighted sum of~\eqref{eq:E:pairdist:chacha} (and note that the intersection with $X_{cd}$ is redundant since $\pi_{cd}\big(\SpecPaths_s^*\cap\cQ_{cd}\big)\subset X_{cd}$). We have
\begin{equation}\label{eq:E:pairdist:cha2}
\sum_{\mathfrak{v}}2\mathfrak{v}_1\cdot \big|\aleph^{-1}(\mathfrak{v})\cap\pi_{cd}\big(\SpecPaths_s^*\cap\cQ_{cd}\big)\big|=\sum_{\mathfrak{v}}2\mathfrak{v}_1\cdot\tfrac{|\SpecPaths^*_s\cap \mathcal{Q}_{cd}||X_{cd}\cap\aleph^{-1}(\mathfrak{v})|}{|X_{cd}|}\pm n^{0.91}\,.
\end{equation}
The left-hand side of this equality is precisely the contribution to $B_{ab,s}$ made by paths in $\cQ_{cd}$, while on the right-hand side, $\sum_{\mathfrak{v}}2\mathfrak{v}_1|X_{cd}\cap\aleph^{-1}(\mathfrak{v})|$ is the contribution to $B_{ab}$ made by paths in $\cQ_{cd}$. By~\ref{quasi:anchorsets}, if $c=d$ then we have
\[|\SpecPaths^*_s\cap \mathcal{Q}_{cd}|=(1\pm\gamAnchor)\tfrac14|\SpecPaths^*_s|\eqByRef{eq:sizeOfAs}(1\pm\gamAnchor)\tfrac14\big(\sigmKJ n\pm n^{0.6}\big)\]
and, since $\cQ_{cd}=\bigsqcup_{s\in\cJ_2}\big(\SpecPaths^*_s\cap\cQ_{cd}\big)$, summing the above bound over the $|\cJ_2|$ choices of $s$ gives $|\cQ_{cd}|=|\cJ_2|\cdot(1\pm\gamAnchor)\tfrac14\big(\sigmKJ n\pm n^{0.6}\big)$. If $c\neq d$, we replace the fraction $\tfrac14$ by $\tfrac12$. In either case, we have
\[\frac{|\SpecPaths^*_s\cap \mathcal{Q}_{cd}|}{|\cQ_{cd}|}=\frac{1\pm3\gamAnchor}{|\cJ_2|}\,,\]
where we used that $\tfrac{1+\gamAnchor}{1-\gamAnchor}\le1+3\gamAnchor$ for the two independent $(1\pm\gamAnchor)$ factors, and absorbed the $n^{0.6}$ terms into the same error.
Plugging these observations in, and summing~\eqref{eq:E:pairdist:cha2} over choices of $c$ and $d$, we get
\begin{equation}\label{eq:DabsSharp}
B_{ab,s}=(1\pm3\gamAnchor)\frac{B_{ab}}{|\cJ_2|}\pm n^{0.95}\,,
\end{equation}
from which~\eqref{eq:Dabs} follows, since $|\cJ_2|=(\sigmKJ-\sigmJnula-\sigmJjedna)n$ and, by~\eqref{eq:E:pairdist:lowerB}, the additive error $n^{0.95}$ is absorbed by increasing the relative error from $3\gamAnchor$ to $4\gamAnchor$. We stress that the estimates below use the form~\eqref{eq:DabsSharp}, in which the additive error is kept separate, rather than the rounded form~\eqref{eq:Dabs}. To deal with the case $a=b=\boxplus$, we use the same argument but replace $2\mathfrak{v}_1$ with $2\mathfrak{v}_3$ throughout; to deal with $\{a,b\}=\{\boxminus,\boxplus\}$ we replace $2\mathfrak{v}_1$ with $\mathfrak{v}_2$.

The difference between $B_{ab,s}$ and the required $f_{ab,s}$ is precisely that peripheral edges count towards the former, but not the latter.

Before treating the individual cases, we record a general observation that we will use for each of them: for any $a,b\in\{\boxminus,\boxplus\}$, writing $M_{ab}:=\tfrac{B_{ab}}{(\sigmKJ-\sigmJnula-\sigmJjedna)n}$, by~\eqref{eq:E:pairdist:lowerB} we have $M_{ab}\ge\tfrac74\sigmKJ n$, so $M_{ab}-\sigmKJ n\ge\tfrac34\sigmKJ n$, i.e.\ $\sigmKJ n\le\tfrac43(M_{ab}-\sigmKJ n)$. Consequently
\begin{align*}
3\gamAnchor M_{ab}+\gamAnchor\sigmKJ n
&=3\gamAnchor(M_{ab}-\sigmKJ n)+4\gamAnchor\cdot\sigmKJ n\le\Big(3+4\cdot\tfrac43\Big)\gamAnchor(M_{ab}-\sigmKJ n)\\
&=\tfrac{25}3\gamAnchor(M_{ab}-\sigmKJ n)\le 9\gamAnchor(M_{ab}-\sigmKJ n)\,,
\end{align*}
and since $M_{ab}-\sigmKJ n=\Theta(n)$ while the additive polynomial error terms appearing below are $O(n^{0.95})=o(n)$, for $n$ sufficiently large this bound of $9\gamAnchor(M_{ab}-\sigmKJ n)$ absorbs those too; this is exactly what lets us combine the two error terms $3\gamAnchor M_{ab}$ (from~\eqref{eq:DabsSharp}) and $\gamAnchor \sigmKJ n$ (from the peripheral-edge count) below into the single error $9\gamAnchor(M_{ab}-\sigmKJ n)$.

Let us look first at $f_{\boxminus\boxminus, s}$. Here, peripheral edges count twice towards $B_{\boxminus\boxminus,s}$ (once for each orientation) and they are necessarily peripheral edges of paths either in $\cQ_{\boxminus\boxminus}\cap\SpecPaths^*_s$, or in $\cQ_{\boxminus\boxplus}\cap\SpecPaths^*_s$. The former is easy to deal with: by~\ref{quasi:anchorsets} and~\eqref{eq:sizeOfAs} we have $\big|\cQ_{\boxminus\boxminus}\cap\SpecPaths^*_s\big|=\tfrac14(1\pm\gamAnchor)\big(\sigmKJ n\pm n^{0.6}\big)$, and consulting Tables~\ref{tab:VminVminL} and~\ref{tab:VminVminR} we see that in all cases, in each such path, exactly one of the two peripheral edges is assigned by $\chi$ within $\Vmin$. Thus the peripheral edge contribution from these paths to $B_{\boxminus\boxminus,s}$ is $\tfrac12(1\pm\gamAnchor)\big(\sigmKJ n\pm n^{0.6}\big)$. For a path $P\in\cQ_{\boxminus\boxplus}$, on the other hand, only one of the two peripheral edges can possibly lie in $\Vmin$. From Tables~\ref{tab:VminVplus} and~\ref{tab:VminVplus2}, we see that whether it does so depends on whether $\pi_{\boxminus\boxplus}(P)$ is in $X^\blacktriangle$ or $X^\triangledown$. These two cases are equally likely and chosen independently for all paths in $\cQ_{\boxminus\boxplus}$. We have by~\ref{quasi:anchorsets} and~\eqref{eq:sizeOfAs} that $\big|\cQ_{\boxminus\boxplus}\cap\SpecPaths^*_s\big|=\tfrac12(1\pm\gamAnchor)(\sigmKJ n\pm n^{0.6})$, so by the Chernoff bound (Fact~\ref{fact:hypergeometricBasicProperties}), with probability at least $1-\exp(-\sqrt{n})$ the contribution from these paths to $B_{\boxminus\boxminus,s}$ is $\tfrac12(1\pm\gamAnchor)(\sigmKJ n\pm n^{0.6})\pm n^{0.9}$. Putting this together, we have
\begin{align*}
 f_{\boxminus\boxminus,s}&=B_{\boxminus\boxminus,s}-(1\pm\gamAnchor)(\sigmKJ n\pm n^{0.6})\pm n^{0.9}\\
 &\eqByRef{eq:DabsSharp} (1\pm3\gamAnchor)M_{\boxminus\boxminus}-\sigmKJ n\mp \gamAnchor\sigmKJ n\pm(n^{0.95}+n^{0.6}+n^{0.9})=(1\pm 9\gamAnchor)\Big(M_{\boxminus\boxminus}-\sigmKJ n\Big)\,,
\end{align*}
using the general observation above (with $a=b=\boxminus$) for the last step. This proves the Claim for $f_{\boxminus\boxminus,s}$. The same argument, replacing $\boxminus\boxminus$ with $\boxplus\boxplus$ throughout (including in the general observation above, which applies equally to $a=b=\boxplus$) and looking at Tables~\ref{tab:VplusVplusL} and~\ref{tab:VplusVplusR}, proves the Claim for $f_{\boxplus\boxplus,s}$.

What remains is to prove the Claim for $f_{\boxminus\boxplus,s}=f_{\boxplus\boxminus,s}$. Note that peripheral edges count only once towards $B_{\boxminus\boxplus,s}$. By exactly the same argument as above, the number of peripheral edges from paths of $\cQ_{\boxminus\boxminus}\cap\SpecPaths^*_s$ which contribute to $B_{\boxminus\boxplus,s}$ is $\tfrac14(1\pm\gamAnchor)\big(\sigmKJ n\pm n^{0.6}\big)$, and we obtain the same bound for paths of $\cQ_{\boxplus\boxplus}\cap\SpecPaths^*_s$. For paths $P$ of $\cQ_{\boxminus\boxplus}\cap\SpecPaths^*_s$, of which there are as above $\tfrac12(1\pm\gamAnchor)(\sigmKJ n\pm n^{0.6})$, from Tables~\ref{tab:VminVplus} and~\ref{tab:VminVplus2} we see that either both or neither peripheral edges contribute to $B_{\boxminus\boxplus,s}$, according to whether $\pi_{\boxminus\boxplus}(P)$ is in $X^\blacktriangle$ or $X^\triangledown$. Much as above, we see that with probability at least $1-\exp(-\sqrt{n})$ the contribution from these paths to $B_{\boxminus\boxplus,s}$ is $\tfrac12(1\pm\gamAnchor)(\sigmKJ n\pm n^{0.6})\pm 2n^{0.9}$, and summing we get
\begin{align*}
 f_{\boxminus\boxplus,s}&=B_{\boxminus\boxplus,s}-(1\pm\gamAnchor)(\sigmKJ n\pm n^{0.6})\pm 2n^{0.9}\\
 &\eqByRef{eq:DabsSharp} (1\pm3\gamAnchor)M_{\boxminus\boxplus}-\sigmKJ n\mp \gamAnchor\sigmKJ n\pm(n^{0.95}+n^{0.6}+2n^{0.9})=(1\pm 9\gamAnchor)\Big(M_{\boxminus\boxplus}-\sigmKJ n\Big)\,,
\end{align*}
again using the general observation above (with $a=\boxminus,b=\boxplus$) for the last step, and where the (very minor) change from $n^{0.9}$ to $2n^{0.9}$ makes no difference since both are $o(n)$ and hence absorbed exactly as before, as required.
\end{claimproof}

As in the definition of the pair distribution property (Definition~\ref{def:pairdistprop}) suppose that we are given $a,b\in \{\boxminus,\boxplus\}$ and an edge $uv\in E(\HStageD)$, with $u\in V_a$ and $v\in V_b$. We claim that with probability at least $1-4\exp(-\sqrt{n})$ we have that
\begin{equation}\label{eq:prob-pair-distr-prop}
\sum_{s\in \cJ_2}w_{uv;s}=(1\pm 20\gamAnchor)
\frac{\sum_{(x,y):xy\in E^*}\ONE_{\chi(x)=V_a, \chi(y)=V_b}}{|V_a||V_b|}\pm 20\gamAnchor n^{-0.01}\;,
\end{equation}
where the quantity $w_{uv;s}$ is as defined in Definition~\ref{def:pairdistprop}.
Assuming we have this, then by the union bound over all possible edges $uv\in E(\HStageD)$ we obtain that with high probability, $\chi$ has the $20\gamAnchor$-pair distribution property. So we need to prove~\eqref{eq:prob-pair-distr-prop} holds with the claimed probability. We split the sum into the contribution coming from~\eqref{case:w2} and~\eqref{case:w3} (i.e.\ respectively $u$ or $v$ is an anchor of some path in $\SpecPaths^*_s$, and the other is not in $\im\phi^\mathbf{D}_s$) and that from~\eqref{case:w1} ($u,v\notin\im\phi^{\mathbf{D}}_s$). These cases are not exhaustive, but for $s\in\cJ_2$ where none of the three cases occurs, we have $w_{uv;s}=0$.

We first deal with the case~\eqref{case:w1}. The numerator of $w_{uv;s}$ in this case~\eqref{case:w1} counts (with multiplicity) the number of pairs which are internal edges of $\SpecPaths^*_s$ with one end assigned in $V_a$ and the other in $V_b$, in other words it is $f_{ab,s}$. The denominator is by~\ref{quasi:indexquasi} (with $S_1=S_2=T_3=\emptyset$ and one of $T_1$ and $T_2$ being $\{s\}$ and the other $\emptyset$) and~\eqref{eq:D-sizeOfFreeSpace} equal to
\[(1\pm\gamAnchor)^2\tfrac14\big((\deltnonspanning+6\sigmKJ )n\pm 6n^{0.6}\big)^2\,,\]
and therefore
\begin{align*}
 w_{uv;s}&=\frac{f_{ab,s}}{(1\pm\gamAnchor)^2\tfrac14\big((\deltnonspanning+6\sigmKJ )n\pm 6n^{0.6}\big)^2}\eqBy{Claim~\ref{claim:chacha}}\frac{(1\pm9\gamAnchor)\big(\frac{B_{ab}}{(\sigmKJ-\sigmJnula-\sigmJjedna)n}-\sigmKJ n\big)}{(1\pm\gamAnchor)^2\tfrac14\big((\deltnonspanning+6\sigmKJ )n\pm 6n^{0.6}\big)^2}\\
 &=4(1\pm 12\gamAnchor)\frac{B_{ab}-\sigmKJ (\sigmKJ-\sigmJnula-\sigmJjedna)n^2}{(\deltnonspanning+6\sigmKJ)^2(\sigmKJ-\sigmJnula-\sigmJjedna)n^3}\,.
\end{align*}
Thus we simply need to know the number of $s\in\cJ_2$ such that $u,v\notin\im\phi^\mathbf{D}_s$. By~\ref{quasi:imagecaps} (with $S_1=\{u,v\}$ and $S_2=\emptyset$) this number is
\[(1\pm\gamAnchor)|\cJ_2|\big(\deltnonspanning+6\sigmKJ\pm n^{-0.4}\big)^2=(1\pm2\gamAnchor)(\sigmKJ-\sigmJnula-\sigmJjedna)n(\deltnonspanning+6\sigmKJ)^2\,,\]
and so we get
\begin{equation}\label{eq:contr5a}
\sum_{\substack{s\in\cJ_2\\ u,v\notin\im\phi^\mathbf{D}_s}}w_{uv;s}=4(1\pm15\gamAnchor)\big(B_{ab}-\sigmKJ (\sigmKJ-\sigmJnula-\sigmJjedna)n^2\big)n^{-2}\,.
\end{equation}

We now deal with the case~\eqref{case:w2}. The denominator of~\eqref{case:w2} is $(1\pm 2\gamAnchor)\tfrac12(\deltnonspanning+6\sigmKJ)n$ by~\ref{quasi:indexquasi} as above, so we have
\[w_{uv;s}=2(1\pm3\gamAnchor)(\deltnonspanning+6\sigmKJ)^{-1}n^{-1}\,.\]
As before, what we need to do is to estimate the number of $s\in\cJ_2$ such that $u=\phi^\mathbf{D}_s(x)$ is an anchor of some path in $\SpecPaths_s^*$ whose first two vertices are $xy$, and $v\notin\im\phi^\mathbf{D}_s$, and $\chi(y)=V_b$. To do this, first we define 
$$S_{uv}:=\{s\in\cJ_2\::\:  \big(\phi^\mathbf{D}_s\big)^{-1}(u)\in \mbox{endvertices of }\SpecPaths^*_s, v\notin \im\phi^\mathbf{D}_s\}\;.$$
What we want, then, is to estimate the number of $s\in S_{uv}$ such that $\chi(y)=V_b$, where $y$ is the neighbour of $\big(\phi^\mathbf{D}_s\big)^{-1}(u)$ in $\SpecPaths_s^*$. To begin with, substituting~\eqref{eq:D-sizeOfFreeSpace}, ~\eqref{eq:sizeOfAs} and \eqref{eq:sizeJ2} into~\ref{quasi:NumAnchorsNoIm}, we have
\begin{align*}
|S_{uv}|&=(1\pm\gamAnchor)|\cJ_2|\cdot \frac{(2\sigmKJ n\pm 2n^{0.6})((\deltnonspanning+6\sigmKJ)n\pm 6n^{0.6})}{n^2}\\
\nonumber
 &=(1\pm2\gamAnchor)
 \cdot 2\sigmKJ(\deltnonspanning+6\sigmKJ)(\sigmKJ-\sigmJnula-\sigmJjedna)n
 \;.
\end{align*}
Now, let $s\in S_{uv}$ be arbitrary. Let $x\in V(G_s)$ be the anchor that is mapped on $u$, and let $y$ be its neighbour in the path $P\in\SpecPaths^*_s$ containing $x$. Whatever $\aleph\big(\pi(P)\big)$ is, observe from Tables~\ref{tab:VminVminL}--\ref{tab:VminVplus2} that $\chi(y)=V_b$ is true for exactly one of $\pi(P)\in X^\blacktriangle$ and $\pi(P)\in X^\triangledown$. In other words, when the random map $\pi$ is fixed, the event $\chi(y)=V_b$ holds with probability $\tfrac12$, and these events are independent for the various $s\in S_{uv}$. Thus with probability at least $1-\exp(-\sqrt{n})$ the number of indices $s\in \cJ_2$ counted in~\eqref{case:w2} is $\frac12|S_{uv}|\pm n^{0.9}$. Under this assumption, we get
\begin{align}
\label{eq:contr5b}
 \sum_{s\in\cJ_2: \text{\eqref{case:w2} applies}}w_{uv;s}&=(\tfrac12|S_{uv}|\pm n^{0.9})\cdot 2(1\pm3\gamAnchor)(\deltnonspanning+6\sigmKJ)^{-1}n^{-1}\\
\nonumber &=(1\pm 6\gamAnchor)
 \cdot 2\sigmKJ(\sigmKJ-\sigmJnula-\sigmJjedna)
 \;,
\end{align}
and similarly we obtain that with probability at least $1-\exp(-\sqrt{n})$,
\begin{align}
\label{eq:contr5c}
 \sum_{s\in\cJ_2: \text{\eqref{case:w3} applies}}w_{uv;s}=(1\pm 6\gamAnchor)
 \cdot 2\sigmKJ(\sigmKJ-\sigmJnula-\sigmJjedna)
 \;.
\end{align}
Summing up~\eqref{eq:contr5a}, \eqref{eq:contr5b} and~\eqref{eq:contr5c}, and writing $K:=\sigmKJ(\sigmKJ-\sigmJnula-\sigmJjedna)n^2$, the total error is at most $4\big(15\gamAnchor(B_{ab}-K)+6\gamAnchor K\big)n^{-2}=4\big(15\gamAnchor B_{ab}-9\gamAnchor K\big)n^{-2}\le 4\cdot 20\gamAnchor B_{ab}n^{-2}$, since $K\ge0$; that is,
\[\sum_{s\in\cJ_2}w_{uv;s}=4(1\pm 20\gamAnchor)\frac{B_{ab}}{n^2}\,.\]
Since $\sum_{(x,y):xy\in E^*}\ONE_{\chi(x)=V_a,\chi(y)=V_b}=B_{ab}$ by the definition of $B_{ab}$ as the occupancy target for the pair $(a,b)$, and $|V_a||V_b|=\big(1\pm O(n^{-1})\big)\tfrac14 n^2$, this is precisely the statement of~\eqref{eq:prob-pair-distr-prop}. Taking the union bound over the various choices, we see that this holds for all $u$ and $v$, giving the desired $20\gamAnchor$-pair distribution property.
\end{proof}

By combining Lemma~\ref{lem:E2decideParts} with Lemma~\ref{lem:pathpack} we can prove Lemma~\ref{lem:StageENew}. This proof is broadly similar to the proof for Stage~D, Lemma~\ref{lem:StageDNew}.

\begin{proof}[Proof of Lemma~\ref{lem:StageENew}]
	We set $\nu=\tfrac{1}{1000}2^{-4\LE}\deltnonspanning^{4\LE}(\sigmKJ^*)^{5\LE}$. Let $C'$ be such that $\tfrac{1}{20(\LD+2)}\cdot 2^{-10\LE}C'$ is returned by Lemma~\ref{lem:pathpack} for input $\nu$ and $\LD$. Note that $C'\gamAnchor\ll \gamNew$ by~\eqref{eq:CONSTANTS}. 
	Let $H=\HStageD$, which has $2\lfloor\tfrac{n}{2}\rfloor$ vertices. As before, the $1\pm O(n^{-1})$ difference between $n$ and $|V(\HStageD)|$ will be absorbed in the error bounds in what follows.
	From Lemma~\ref{lem:E2decideParts} we obtain a map $\chi$. 
	
	We now translate our setting into the setup for Lemma~\ref{lem:pathpack} with $\nu$ and $L=\LD$, and with $\gamma=20(\LD+2)\gamAnchor$. We use graph $H$ with the partition $\Vmin\dcup\Vplus$. Temporarily abusing notation, let $s^*:=|\cJ_2|\ge\nu n$ and suppose $\cJ_2=[s^*]$. Given $s\in\cJ_2$ we obtain a path-forest $F_s=\SpecPaths^*_s$, whose leaves $A_s$ are anchored by the embedding $\phi_s:=\phi_s^{\mathbf{D}}\big|_{A_s}$, with used set $U_s:=\im\phi_s^{\mathbf{D}}$. Recall that each path in $F_s$ has 8 vertices, which is less than $\LD$. The required assignment of sides is provided by $\chi_s$, defined as the restriction of $\chi$ (as defined in Lemma~\ref{lem:E2decideParts}) to $\SpecPaths_s^*$.
	
	We now check that all the assumptions of Lemma~\ref{lem:pathpack} are satisfied. The structural requirements of~\ref{PA:BasicSetting} hold by the above constructions. In Claim~\ref{UMES:BlockQuasiDiet}, we justify assumptions~\ref{PA:BlockQuasirandom} and~\ref{PA:BlockDiet}. In Claim~\ref{UMES:anchor}, we justify assumption~\ref{PA:AnchorDistribution}. In Claim~\ref{UMES:pairdistribution}, we justify assumption~\ref{PA:PairDistributionProperty}. In Claim~\ref{UMES:Ineq}, we justify the inequalities of~\ref{PA:Ineq}. Finally, in Claim~\ref{UMES:densities}, we verify the leftover densities requirement~\ref{PA:LeftOverDensities}.
	
	\begin{claimInLemmaE}\label{UMES:BlockQuasiDiet}
		The graph $H$ and the sets $U_s$ satisfy the requirements of~\ref{PA:BlockQuasirandom} and~\ref{PA:BlockDiet}.
	\end{claimInLemmaE}
	\begin{claimproof}[Proof of Claim~\ref{UMES:BlockQuasiDiet}]
		Recall that by~\ref{quasi:indexquasi} the graph $\HStageD$ with partition $(\Vmin,\Vplus)$ is $(\LD,\gamAnchor,\dStageD,\dStageD)$-index-quasirandom with respect to $\big(\imD(s)\big)_{s\in\cJ}$ and $\big(I_s\big)_{s\in\cJ_0}$, with $I_s:= \{i\in [\lfloor\frac{n}{2}\rfloor]: \boxminus_i,\boxplus_i\in A_s\}$. By Lemma~\ref{lem:indexquasi-implies-blockquasi}, we conclude that $\HStageD$ is $((\LD+3)\gamAnchor,\LD)$-block-quasirandom and $((\LD+3)\gamAnchor,\LD)$-block-diet with respect to each set $U_s$ for $s\in\cJ_2$. By choice of $\gamma$ we have $\gamma\ge(\LD+3)\gamAnchor$.
	\end{claimproof}
	
	\begin{claimInLemmaE}\label{UMES:anchor}
		For each $s\in\cJ_2$, the path-forest $F_s$ has the $\gamma$-anchor distribution property with respect to $H$, fulfilling~\ref{PA:AnchorDistribution}.
	\end{claimInLemmaE}
	\begin{claimproof}[Proof of Claim~\ref{UMES:anchor}]
		By Lemma~\ref{lem:E2decideParts}\ref{en:abstractpscatt}, for each $s\in \cJ_2$ the forest $F_s$ with anchors $\phi_s^{\mathbf{D}}|_{A_s}$ has the $6\gamAnchor$-anchor distribution property with respect to $\HStageD$. Since $\gamma \ge 20(\LD+2)\gamAnchor > 6\gamAnchor$, the $\gamma$-anchor distribution property follows immediately.
	\end{claimproof}
	
	\begin{claimInLemmaE}\label{UMES:pairdistribution}
		The collection of path-forests with anchors and assignment $\chi$ has the $\gamma$-pair distribution property, fulfilling~\ref{PA:PairDistributionProperty}.
	\end{claimInLemmaE}
	\begin{claimproof}[Proof of Claim~\ref{UMES:pairdistribution}]
		By Lemma~\ref{lem:E2decideParts}\ref{en:pairdistrpropE} the collection of forests with anchors and assignment $\chi$ has the $20\gamAnchor$-pair distribution property. Again, by our choice of $\gamma$, this implies the $\gamma$-pair distribution property.
	\end{claimproof}
	
	\begin{claimInLemmaE}\label{UMES:Ineq}
		The assignments satisfy the inequalities of~\ref{PA:Ineq}.
	\end{claimInLemmaE}
	\begin{claimproof}[Proof of Claim~\ref{UMES:Ineq}]
		Lemma~\ref{lem:E2decideParts}\ref{en:medunkovycaj} gives us $|\chi_s^{-1}(\{V_a\})\setminus A_s|\ge\tfrac15\sigmKJ n\ge \nu n$ for each $a\in\{\boxminus,\boxplus\}$ and $s\in\cJ_2$, and Lemma~\ref{lem:E2decideParts}\ref{en:nocheins} gives us
		\[\big|\big\{(x,y)\in E(F_s):x\in \chi_s^{-1}(\{V_a\})\setminus A_s, y\in \chi_s^{-1}(\{V_b\})\setminus A_s\big\}\big|\ge\tfrac{1}{16}\sigmKJ n\ge\nu n\]
		for each $a,b\in\{\boxminus,\boxplus\}$ and $s\in\cJ_2$.
		
		For each $s\in\cJ_2$, recall from~\eqref{eq:D-sizeOfFreeSpace} (using $U_s=\imD(s)$) that $n-|U_s|=(\deltnonspanning+6\sigmKJ)n\pm6n^{0.6}$. Hence by~\ref{quasi:indexquasi} we have, for each $a\in\{\boxminus,\boxplus\}$, $|V_a\setminus U_s|=(1\pm\gamAnchor)\cdot\tfrac12\big(n-|U_s|\big)=(1\pm\gamAnchor)\cdot\tfrac12\big((\deltnonspanning+6\sigmKJ)n\pm6n^{0.6}\big)\ge\tfrac14\deltnonspanning n+6\sigmKJ n$, where we use the fact that $\sigmKJ\ll\deltnonspanning$. Since $|\SpecPaths_s^*|\le\sigmKJ n$ and each path has $6$ interior vertices, even if all the vertices of the paths $\SpecPaths_s^*$ were assigned to $V_a$ we obtain $|V_a\setminus U_s|-\big|\{x\in V(F_s)\setminus A_s:\chi_s(x)=V_a\}\big|\ge\tfrac14\deltnonspanning n\ge \nu n$. These bounds confirm all inequalities of~\ref{PA:Ineq}.
	\end{claimproof}
	
	\begin{claimInLemmaE}\label{UMES:densities}
		The leftover densities satisfy $d'_{\boxminus\boxminus}, d'_{\boxminus\boxplus}, d'_{\boxplus\boxplus} \ge \nu$, fulfilling~\ref{PA:LeftOverDensities}. Moreover, $d'_{\boxminus\boxminus}=d'_{\boxplus\boxplus}=(1\pm3\gamAnchor)\dStageE^*$ and $d'_{\boxminus\boxplus}=(1\pm3\gamAnchor)\dStageEbar$.
	\end{claimInLemmaE}
	\begin{claimproof}[Proof of Claim~\ref{UMES:densities}]
		We define $d'_{\boxminus\boxminus},d'_{\boxminus\boxplus},d'_{\boxplus\boxplus}$ to be the densities obtained by a successful packing of the $F_s$ according to $\chi_s$; explicitly, $d'_{\boxminus\boxminus}:=e(\HStageE[\Vmin])/\binom{\lfloor n/2\rfloor}{2}$, $d'_{\boxplus\boxplus}:=e(\HStageE[\Vplus])/\binom{\lfloor n/2\rfloor}{2}$, and $d'_{\boxminus\boxplus}:=e(\HStageE[\Vmin,\Vplus])/\lfloor n/2\rfloor^2$. By Lemma~\ref{lem:E2decideParts}\ref{en:Etotaledges}, the number of edges of $E^*$ assigned by $\chi$ entirely within $\Vmin$ is exactly $\sum_{xy\in E^*}\ONE_{\chi(x),\chi(y)=\Vmin}=e_{\HStageD}(\Vmin)-\big(\sigmJnula\sigmKJ^* n^2+2j_{\boxminus\boxminus}+3(j_{\boxminus\boxplus}+j_{\boxplus\boxplus})\big)$; since these are precisely the edges of $\HStageD[\Vmin]$ consumed in packing the $F_s$, this gives
		\[e(\HStageE[\Vmin])=e_{\HStageD}(\Vmin)-\sum_{xy\in E^*}\ONE_{\chi(x),\chi(y)=\Vmin}=\sigmJnula\sigmKJ^* n^2+2j_{\boxminus\boxminus}+3(j_{\boxminus\boxplus}+j_{\boxplus\boxplus})\,,\]
		and analogously $e(\HStageE[\Vplus])=\sigmJnula\sigmKJ^* n^2+2j_{\boxplus\boxplus}+3(j_{\boxminus\boxplus}+j_{\boxminus\boxminus})$ and $e(\HStageE[\Vmin,\Vplus])=\sigmJnula\sigmKJ^* n^2+6(j_{\boxminus\boxminus}+j_{\boxplus\boxplus})+5j_{\boxminus\boxplus}$, matching~\eqref{eq:E:edges} exactly. In particular, unlike in the proof of Claim~\ref{UMDS:densities} for Stage~D, no probabilistic estimate is needed here, since~\ref{en:Etotaledges} gives these edge counts exactly rather than only with high probability. Since $j_{\boxminus\boxminus},j_{\boxminus\boxplus},j_{\boxplus\boxplus}\ge0$ (being counts of paths), each of the three displayed quantities is at least $\sigmJnula\sigmKJ^* n^2$, and so
		\[d'_{\boxminus\boxminus},d'_{\boxminus\boxplus},d'_{\boxplus\boxplus}\ge \frac{\sigmJnula\sigmKJ^* n^2}{\binom{\lfloor n/2\rfloor}{2}}=(1+o(1))\cdot8\sigmJnula\sigmKJ^*\;.\]
		By~\eqref{eq:CONSTANTS} we have $\sigmJnula\sigmKJ^*\gg\nu$, so for $n$ sufficiently large the right-hand side is at least $\nu$, giving $d'_{\boxminus\boxminus},d'_{\boxminus\boxplus},d'_{\boxplus\boxplus}\ge\nu$ as required.

		For the second part, we compute $j_{\boxminus\boxminus},j_{\boxminus\boxplus},j_{\boxplus\boxplus}$ precisely. Recall from~\ref{quasi:anchorsets} (applied to $\HStageD$, which holds for $s\in\cJ_1$ in particular) that for each $a,b\in\{\boxminus,\boxplus\}$ and $s\in\cJ_1$ we have $|A_s^b\cap V_a|=(1\pm\gamAnchor)\tfrac14|A_s|$, where $A_s^b$ denotes the anchors of $F_s=\SpecPaths^*_s$ whose \emph{other} endvertex lies in $V_b$. A path of $F_s$ with both anchors in $V_\boxminus$ contributes two elements to $A_s^\boxminus\cap V_\boxminus$ (one for each of its anchors), while a path with anchors split between $\Vmin$ and $\Vplus$ contributes exactly one element to $A_s^\boxplus\cap V_\boxminus$ (and one to $A_s^\boxminus\cap V_\boxplus$). Hence, summing over $s\in\cJ_1$,
		\[j_{\boxminus\boxminus}=\tfrac12\sum_{s\in\cJ_1}|A_s^\boxminus\cap V_\boxminus|=\tfrac18(1\pm\gamAnchor)\sum_{s\in\cJ_1}|A_s|\,,\quad j_{\boxplus\boxplus}=\tfrac18(1\pm\gamAnchor)\sum_{s\in\cJ_1}|A_s|\,,\quad j_{\boxminus\boxplus}=\tfrac14(1\pm\gamAnchor)\sum_{s\in\cJ_1}|A_s|\,,\]
		using that the sum of quantities each within a $(1\pm\gamAnchor)$-factor of a common multiple of a nonnegative base value is itself within the same $(1\pm\gamAnchor)$-factor of the sum of those base values. Since each path of $F_s$ contributes two anchors, $|A_s|=2|\SpecPaths^*_s|=2(\sigmJjedna n\pm n^{0.6})$ for $s\in\cJ_1$, and $|\cJ_1|=\sigmJjedna n$ exactly by~\eqref{eq:sizeJ0J1}, so
		\[\sum_{s\in\cJ_1}|A_s|=2\sigmJjedna n(\sigmJjedna n\pm n^{0.6})=2\sigmJjedna^2 n^2\pm2\sigmJjedna n^{1.6}\;.\]
		Since $\sigmJjedna$ is a fixed constant, the additive term $2\sigmJjedna n^{1.6}$ is $o(\gamAnchor\sigmJjedna^2n^2)$ for $n$ sufficiently large, so for such $n$ we obtain (absorbing this lower-order term by increasing the multiplicative error from $\gamAnchor$ to $2\gamAnchor$)
		\[j_{\boxminus\boxminus}=j_{\boxplus\boxplus}=\tfrac14\sigmJjedna^2n^2(1\pm2\gamAnchor)\,,\qquad j_{\boxminus\boxplus}=\tfrac12\sigmJjedna^2n^2(1\pm2\gamAnchor)\;.\]
		(As a sanity check, $j_{\boxminus\boxminus}+j_{\boxplus\boxplus}+j_{\boxminus\boxplus}\approx\sigmJjedna^2n^2\approx\sum_{s\in\cJ_1}|\SpecPaths^*_s|$, as it must, since every path of $F_s$, $s\in\cJ_1$, falls into exactly one of these three categories.)
		Substituting into~\eqref{eq:E:edges}, and using that $2j_{\boxminus\boxminus}+3(j_{\boxminus\boxplus}+j_{\boxplus\boxplus})$ and $6(j_{\boxminus\boxminus}+j_{\boxplus\boxplus})+5j_{\boxminus\boxplus}$ are each sums of nonnegative terms with the common multiplicative error $(1\pm2\gamAnchor)$, we get
		\[e_{\HStageE}(\Vmin)=\sigmJnula\sigmKJ^*n^2+\tfrac{11}4\sigmJjedna^2n^2(1\pm2\gamAnchor)\,,\qquad e_{\HStageE}(\Vmin,\Vplus)=\sigmJnula\sigmKJ^*n^2+\tfrac{11}2\sigmJjedna^2n^2(1\pm2\gamAnchor)\;,\]
		and likewise for $e_{\HStageE}(\Vplus)$ in place of $e_{\HStageE}(\Vmin)$. Dividing by $\binom{\lfloor n/2\rfloor}2=(1+O(n^{-1}))\tfrac{n^2}8$ (resp.\ $\lfloor n/2\rfloor^2=(1+O(n^{-1}))\tfrac{n^2}4$) gives
		\[d'_{\boxminus\boxminus},d'_{\boxplus\boxplus}=\big(1+O(n^{-1})\big)\dStageE^*\pm44\gamAnchor\sigmJjedna^2\big(1+O(n^{-1})\big)\,,\qquad d'_{\boxminus\boxplus}=\big(1+O(n^{-1})\big)\dStageEbar\pm44\gamAnchor\sigmJjedna^2\big(1+O(n^{-1})\big)\;,\]
		using the definitions $\dStageE^*=8(\sigmKJ^*\sigmJnula+\tfrac{11}4\sigmJjedna^2)$ and $\dStageEbar=4(\sigmKJ^*\sigmJnula+\tfrac{11}2\sigmJjedna^2)$. Since both $\dStageE^*,\dStageEbar\ge22\sigmJjedna^2$ (dropping the nonnegative $\sigmKJ^*\sigmJnula$ term from each), the additive error $44\gamAnchor\sigmJjedna^2$ is at most $2\gamAnchor$ times either target, so for $n$ sufficiently large (absorbing the $O(n^{-1})$ term as before) we conclude
		\[d'_{\boxminus\boxminus}=d'_{\boxplus\boxplus}=(1\pm3\gamAnchor)\dStageE^*\,,\qquad d'_{\boxminus\boxplus}=(1\pm3\gamAnchor)\dStageEbar\;,\]
		as required.
	\end{claimproof}
	
	This is the setup for Lemma~\ref{lem:pathpack}, and the conclusion of that lemma is that \PathPacking\ with high probability succeeds in packing each $F_s$ as specified. That is, we obtain for each $s\in\cJ_2$ an embedding $\phi'_s$ of $F_s$ into $\HStageD\big[(V(H)\setminus U_s)\cup\phi_s(A_s)\big]$ which extends $\phi_s$, such that $\phi'_s(x)\in\chi_s(x)$ for each $x\in V(F_s)$, and such that for each $uv\in E(H)$ there is at most one $s$ such that $\phi'_s$ uses the edge $uv$. Suppose that this likely event occurs, and as before we will list polynomially many further events which we presume all occur, and fix such an outcome. Note that the existence of a packing according to $\chi$, by Lemma~\ref{lem:E2decideParts}~\ref{en:Etotaledges}, in particular gives us the required edge counts~\eqref{eq:E:edges}.
	
	We now define the setup after Stage~E and check it satisfies the conditions of Lemma~\ref{lem:StageENew}. We let $\HStageE$ be the graph $H'$ returned by Lemma~\ref{lem:pathpack}. For $s\in\cJ_0\cup\cJ_1$, we set $\phi^\mathbf{E}_s=\phi^\mathbf{D}_s$, with the same anchor and used sets and the same path-forests to embed. Recall that the conclusion of Lemma~\ref{lem:StageENew} has path-forests indexed by $\cJ_0\cup\cJ_1$, with $\emptyset$ in place of $\cJ_2$. We now argue that this is a $(\gamNew,\LE,\dStageE^*,\dStageEbar)$-quasirandom setup, as required for Lemma~\ref{lem:StageENew}. To ensure clarity, we structure the verification into separate claims.
	
	\begin{claimInLemmaE}\label{UMES:Quasi1}
		The setup satisfies index-quasirandomness as required by~\ref{quasi:indexquasi}.
	\end{claimInLemmaE}
	\begin{claimproof}[Proof of Claim~\ref{UMES:Quasi1}]
		Fix $S_1,S_2,T_1,T_2,T_3$ as in Definition~\ref{def:index-quasirandom}, each of size at most $\LE$, with families of sets $(U^\mathbf{E}_s)_{s\in\cJ_0\cup\cJ_1}$ and $(J_s)_{s\in\cJ_0}$, and let $X=\mathbb{U}_{\HStageD}(S_1,S_2,T_1,T_2,T_3)$. Then we have $|X|\ge\nu n$: since the setup after Stage~D is $(\LD,\gamAnchor,\dStageD,\dStageD)$-index-quasirandom (and $\LE\le\LD$ by~\eqref{eq:CONSTANTS}) with $|S_1|+|S_2|\le2\LE$, $|T_1\cup T_2|\le2\LE$, $|T_3|\le\LE$, Definition~\ref{def:index-quasirandom} gives
		\[|X|=(1\pm\gamAnchor)\dStageD^{|S_1|+|S_2|}\cdot\prod_{\ell\in T_1\cup T_2}\Big(1-\tfrac{|U^\mathbf{E}_\ell|}{n}\Big)\cdot\prod_{\ell\in T_3}\tfrac{|J_\ell|}{n/2}\cdot\tfrac n2\ge\tfrac12\dStageD^{2\LE}\deltnonspanning^{2\LE}(2\sigmKJ^*)^{\LE}\cdot\tfrac n2\,,\]
		where we used that every $G_s$, $s\in\cJ$ (so in particular $s\in\cJ_0\cup\cJ_1$), has exactly $(1-\deltnonspanning)n$ vertices (so $|U^\mathbf{E}_\ell|\le(1-\deltnonspanning)n$, giving $1-|U^\mathbf{E}_\ell|/n\ge\deltnonspanning$, for each $\ell$) and $|J_\ell|=\sigmKJ^* n$ by construction, and that a product of at most $2\LE$ (resp.\ $\LE$) factors each in $(0,1)$ is minimised by taking all $2\LE$ (resp.\ $\LE$) factors at their smallest value. Since $\dStageD\ge13(\sigmKJ^*)^2$ for $n$ sufficiently large (Section~\ref{sssection:EvolutionOfDensities}, using $\sigmKJ^*\le\sigmKJ$ and $\sigmJnula,\sigmJjedna\ll\sigmKJ$ in~\eqref{eq:dStageD}), this is at least $\tfrac{n}{4}\cdot13^{2\LE}2^{\LE}\deltnonspanning^{2\LE}(\sigmKJ^*)^{5\LE}\ge\tfrac{n}{1000}2^{-4\LE}\deltnonspanning^{4\LE}(\sigmKJ^*)^{5\LE}=\nu n$, as claimed. Now~\ref{pathpack:megaqr} is stated in terms of the true densities $d_{ab}:=d(\HStageD[V_a])$ (resp.\ $d(\HStageD[\Vmin,\Vplus])$) and $d'_{ab}$ from~\ref{PA:LeftOverDensities}, not directly in terms of $\dStageD,\dStageE^*,\dStageEbar$. By Claim~\ref{UMES:BlockQuasiDiet} (via Lemma~\ref{lem:indexquasi-implies-blockquasi} applied there) we have $d_{ab}=(1\pm2\gamAnchor)\dStageD$ for each $a,b$, and by Claim~\ref{UMES:densities} we have $d'_{\boxminus\boxminus}=d'_{\boxplus\boxplus}=(1\pm3\gamAnchor)\dStageE^*$ and $d'_{\boxminus\boxplus}=(1\pm3\gamAnchor)\dStageEbar$; substituting these (and dividing out the common, now-justified, factor of $\dStageD$ up to a $(1\pm5\gamAnchor)$-type error, which is dominated by and absorbed into the coefficient $C'$ of the main error term below, since $C'$ is larger than any fixed constant) gives $\tfrac{d'_{\boxminus\boxminus}}{d_{\boxminus\boxminus}}=\tfrac{d'_{\boxplus\boxplus}}{d_{\boxplus\boxplus}}=(1\pm5\gamAnchor)\tfrac{\dStageE^*}{\dStageD}$ and $\tfrac{d'_{\boxminus\boxplus}}{d_{\boxminus\boxplus}}=\tfrac{d'_{\boxplus\boxminus}}{d_{\boxplus\boxminus}}=(1\pm5\gamAnchor)\tfrac{\dStageEbar}{\dStageD}$. So by~\ref{pathpack:megaqr} of Lemma~\ref{lem:pathpack} it is likely that we have
		\[
		\big|X\cap \mathbb{U}'_{\HStageE}(S_1,S_2,\emptyset,\emptyset)\big|=
		(1\pm 20C'(\LD+2)\gamAnchor)\big(\tfrac{\dStageE^*}{\dStageD}\big)^{|S_1\cap \Vmin|+|S_2\cap\Vplus|}\big(\tfrac{\dStageEbar}{\dStageD}\big)^{|S_1\cap \Vplus|+|S_2\cap\Vmin|}|X|\,.
		\]
		Plugging in the size $|X|$ given by~\ref{quasi:indexquasi} after Stage~D, and observing that the above set is precisely $\mathbb{U}_{\HStageE}(S_1,S_2,T_1,T_2,T_3)$ with families of sets $(U^\mathbf{E}_s)_{s\in\cJ_0\cup\cJ_1}$ and $(J_s)_{s\in\cJ_0}$, this is what is required for $(\LE,\gamNew,\dStageE^*,\dStageEbar)$-index-quasirandomness.
	\end{claimproof}
	
	\begin{claimInLemmaE}\label{UMES:Quasi234}
		The setup satisfies~\ref{quasi:anchorsets}, \ref{quasi:PathsWithOneVertex}, and \ref{quasi:termVtxNbs}.
	\end{claimInLemmaE}
	\begin{claimproof}[Proof of Claim~\ref{UMES:Quasi234}]
		For~\ref{quasi:anchorsets}, observe that for $\cJ_1$ nothing has changed from after Stage~D, and we no longer consider $\cJ_2$. Similarly, for~\ref{quasi:PathsWithOneVertex} nothing has changed from after Stage~D for $\cJ_1$. For~\ref{quasi:termVtxNbs}, it suffices by construction to establish~\ref{quasi:distributionOfAnchors} as in Stage~D.
	\end{claimproof}
	
	\begin{claimInLemmaE}\label{UMES:Quasi5}
		The setup satisfies~\ref{quasi:distributionOfAnchors}.
	\end{claimInLemmaE}
	\begin{claimproof}[Proof of Claim~\ref{UMES:Quasi5}]
		Fix $s\in\cJ_0\cup\cJ_1$, $a\in\{\boxminus,\boxplus\}$ and $u\in V(\HStageE)\setminus U^\mathbf{E}_s$. By~\ref{quasi:distributionOfAnchors} after Stage~D, writing $A^a_s$ for the set of anchors $v\in A^\mathbf{E}_s$ such that the path end not anchored to $v$ is anchored in $V_a$, we have
		\begin{align*}
			\big|\Vmin\cap \NBH_{\HStageD}(u)\cap A^a_s\big|&=(1\pm\gamAnchor)\dStageD |A^a_s\cap\Vmin|\quad\text{and}\\
			\big|\Vplus\cap \NBH_{\HStageD}(u)\cap A^a_s \big|&=(1\pm\gamAnchor)\dStageD |A^a_s\cap\Vplus|\,.
		\end{align*}
		For each $s\in\cJ_0\cup\cJ_1$ and $a\in\{\boxminus,\boxplus\}$, we define a weight function $w=w_{s,a}$ on $V(\HStageD)$ by setting $w(v)=1$ if $v\in A^a_s$, and otherwise $w(v)=0$; note that $\sum_{v\in V_b}w(v)=|A^a_s\cap V_b|$ is either equal to $0$ or at least $\nu n$, for every such $w_{s,a}$: for $s\in\cJ_1$, we have $|A^a_s\cap V_b|=(1\pm\gamAnchor)\tfrac14|A_s|$ by~\ref{quasi:anchorsets}, which is positive (and comfortably at least $\nu n$, since $|A_s|=\Theta(n)$); for $s\in\cJ_0$, by construction each path of $F_s$ is anchored at a matched pair $\{\boxminus_i,\boxplus_i\}$ with one endpoint on each side, so an anchor's partner always lies on the side opposite the anchor itself, giving $|A^a_s\cap V_a|=0$ and $|A^a_s\cap V_{a'}|=\tfrac12|A_s|\ge\nu n$ for $a'$ the side other than $a$. The failure probability $\exp(-n^{0.2})$ in~\ref{pathpack:weights} is uniform over all admissible weight functions (it does not depend on which $w$ is chosen), so applying~\ref{pathpack:weights} separately with $w=w_{s,a}$ for each of the $2|\cJ_0\cup\cJ_1|=O(n)$ pairs $(s,a)$ and taking the union bound, with probability at least $1-n\exp(-n^{0.2})=1-o(1)$ the conclusion of~\ref{pathpack:weights} holds simultaneously for all $s\in\cJ_0\cup\cJ_1$, $a\in\{\boxminus,\boxplus\}$, and $u\in V(\HStageE)\setminus U^\mathbf{E}_s$. Now~\ref{pathpack:weights} is stated in terms of the true densities $d_{ab}$ and $d'_{ab}$, not directly in terms of $\dStageD,\dStageE^*,\dStageEbar$; as in the proof of Claim~\ref{UMES:Quasi1}, Claim~\ref{UMES:BlockQuasiDiet} gives $d_{ab}=(1\pm2\gamAnchor)\dStageD$ for each $a,b$, and Claim~\ref{UMES:densities} gives $d'_{\boxminus\boxminus}=d'_{\boxplus\boxplus}=(1\pm3\gamAnchor)\dStageE^*$ and $d'_{\boxminus\boxplus}=(1\pm3\gamAnchor)\dStageEbar$; substituting these (with the resulting extra $(1\pm5\gamAnchor)$-type error again absorbed into $C'$, exactly as before) gives $\tfrac{d'_{\boxminus\boxminus}}{d_{\boxminus\boxminus}}=\tfrac{d'_{\boxplus\boxplus}}{d_{\boxplus\boxplus}}=(1\pm5\gamAnchor)\tfrac{\dStageE^*}{\dStageD}$ and $\tfrac{d'_{\boxminus\boxplus}}{d_{\boxminus\boxplus}}=\tfrac{d'_{\boxplus\boxminus}}{d_{\boxplus\boxminus}}=(1\pm5\gamAnchor)\tfrac{\dStageEbar}{\dStageD}$. Thus, if $u\in\Vmin$ it is likely that we have
		\begin{align*}
			\big|\Vmin\cap \NBH_{\HStageE}(u)\cap A^a_s \big|&=(1\pm\gamAnchor)\dStageD |A^a_s\cap\Vmin|(1\pm 20C'(\LD+2)\gamAnchor)\tfrac{\dStageE^*}{\dStageD}\quad\text{and}\\
			\big|\Vplus\cap \NBH_{\HStageE}(u)\cap A^a_s \big|&=(1\pm\gamAnchor)\dStageD |A^a_s\cap\Vplus|(1\pm 20C'(\LD+2)\gamAnchor)\tfrac{\dStageEbar}{\dStageD}\,.
		\end{align*}
		Similarly if $u\in\Vplus$ it is likely that we have
		\begin{align*}
			\big|\Vmin \cap \NBH_{\HStageE}(u)\cap A^a_s\big|&=(1\pm\gamAnchor)\dStageD |A^a_s\cap\Vmin|(1\pm 20C'(\LD+2)\gamAnchor)\tfrac{\dStageEbar}{\dStageD}\quad\text{and}\\
			\big|\Vplus \cap \NBH_{\HStageE}(u)\cap A^a_s\big|&=(1\pm\gamAnchor)\dStageD |A^a_s\cap\Vplus|(1\pm 20C'(\LD+2)\gamAnchor)\tfrac{\dStageE^*}{\dStageD}\,.
		\end{align*}
		This is as required for~\ref{quasi:distributionOfAnchors}. 
	\end{claimproof}
	
	\begin{claimInLemmaE}\label{UMES:Quasi678}
		The setup satisfies~\ref{quasi:imagecaps}, \ref{quasi:NumAnchors}, and~\ref{quasi:NumAnchorsNoIm}.
	\end{claimInLemmaE}
	\begin{claimproof}[Proof of Claim~\ref{UMES:Quasi678}]
		For~\ref{quasi:imagecaps}, observe that nothing has changed for $\cJ_0\cup\cJ_1$ and that we no longer consider $\cJ_2$.
		For~\ref{quasi:NumAnchors}, observe that nothing has changed since after Stage~D for $\cJ_0$ and $\cJ_1$, while we no longer consider $\cJ_2$. The same holds for~\ref{quasi:NumAnchorsNoIm}.
	\end{claimproof}
	
\end{proof}

\section{Stage~F (Proof of Lemma~\ref{lem:StageFNew})}\label{sec:StageF}

In this stage we embed the paths of $\bigcup_{s\in\cJ_1}\SpecPaths^*_s$ to precisely adjust the degree condition of the vertices. Recall from~\eqref{eq:t(v)} that $t(v)$ denotes the number of paths from $\{\SpecShortPaths_s,s\in\cJ_0\}$ anchored at $v$. We need to obtain~\eqref{eq:stageF:degtarget}, that is
\[ \deg_{\HStageF}(\boxminus_i;\Vplus)=\deg_{\HStageF}(\boxminus_i;\Vmin)-t(\boxminus_i)\,\text{ and }\,\deg_{\HStageF}(\boxplus_i;\Vmin)=\deg_{\HStageF}(\boxplus_i;\Vplus)-t(\boxplus_i)\]
for each $i\in[\lfloor n/2\rfloor]$. We write $\imE(s)=\im(\phi^\mathbf{E}_s)$\index{$\imE(s)$} and $\imF(s)=\im(\phi^\mathbf{F}_s)$.\index{$\imF(s)$}

For each $s\notin\cJ_1$ we set $\phi^\mathbf{F}_s:=\phi^\mathbf{E}_s$. We will construct embeddings $\phi^\mathbf{F}_s$ extending $\phi^\mathbf{E}_s$ for each $s\in\cJ_1$ that do not use any edge of $H$ multiple times and give us~\eqref{eq:stageF:degtarget} for each $i$. First, we argue that any such maps give us the required quasirandom setup. In this, we use the fact that $\gamNew,\sigmJjedna\ll\gamJ\ll \sigmJnula$ from~\eqref{eq:CONSTANTS}. Recall that $\dStageE^*$ and $\dStageEbar$ denote the (approximate) densities within the sets $\Vmin$ and $\Vplus$, and between them, respectively (see~\eqref{eq:dStageE-def}), and that $\dStageEbar,\dStageE^*\ge 4\sigmJnula\sigmKJ^*\ge\sigmJnula^2$ from Fact~\ref{fact:min:d_E}.

We begin with~\ref{quasi:indexquasi}.
For an arbitrary set $T\subseteq \cJ_0$ we have $\bigcup_{s\in T}\imF(s)=\bigcup_{s\in T}\imE(s)$. 
As we are extending the mappings only on paths from $\bigcup_{s\in \cJ_1}\SpecPaths^*_s$, and since for each $s\in \cJ_1$ the mapping $\phi^{\mathbf{F}}_s$ is an embedding, observe that for each $v\in V(H)$, we have 
\begin{equation}\label{eq:lk98}
|\NBH_{\HStageE}(v)\setminus \NBH_{\HStageF}(v)|\le 2|\cJ_1|=2\sigmJjedna n\;.
\end{equation}
We can also already fix the densities $\dStageF^*$ and $\dStageFbar$ promised by Lemma~\ref{lem:StageFNew}. Each of the paths that remain to be embedded in Stage~F will use prescribed numbers of edges within $\Vmin$, within $\Vplus$, and between the two sides; this is guaranteed by~\ref{it:F_Discr1}--\ref{it:F_Discr3} of Lemma~\ref{lem:repairDegreeDiscrepancy} and by the \emph{Moreover} part of Lemma~\ref{lem:neutralembedding} below. Summing those numbers over the $j_{\boxminus\boxminus}$, $j_{\boxminus\boxplus}$ and $j_{\boxplus\boxplus}$ paths of the three types, and subtracting from the edge counts~\eqref{eq:E:edges} of Lemma~\ref{lem:StageENew}, every term involving the $j_{ab}$ cancels and we are left with
\[e_{\HStageF}(\Vmin)=e_{\HStageF}(\Vplus)=e_{\HStageF}(\Vmin,\Vplus)=\sigmJnula\sigmKJ^* n^2\;.\]
Accordingly we set\index{$\dStageF^*$}\index{$\dStageFbar$}
\[\dStageF^*:= 8\sigmKJ^*\sigmJnula \quad\text{and}\quad \dStageFbar:= 4\sigmKJ^*\sigmJnula\;,\]
which are exactly the densities $d(\HStageF[\Vmin])=d(\HStageF[\Vplus])$ and $d(\HStageF[\Vmin,\Vplus])$ up to a factor $1+O(n^{-1})$, absorbed into the error terms as before. In particular the densities in $\HStageE$ and in $\HStageF$ are nearly identical: comparing with~\eqref{eq:dStageE-def} we get $\dStageE^*-\dStageF^*=\dStageEbar-\dStageFbar=22\sigmJjedna^2\le 100\sigmJjedna^2$.

Hence, for sets $S\subseteq V(\HStageE)$ and $T\subseteq\cJ_0$, we have 
\begin{align}
\left|\NBH_{\HStageF}(S)\setminus \bigcup_{s\in T}\imF(s)\right|
\eqByRef{eq:lk98}\label{eq:F-jointneighbourhood}
\left|\NBH_{\HStageE}(S)\setminus \bigcup_{s\in T}\imE(s)\right|\pm |S|\cdot 2\sigmJjedna n\;.
\end{align}

Assume that we are given sets $S_1,S_2\subset \Vmin\cup\Vplus$ and pairwise disjoint sets $T_1,T_2,T_3\subseteq\cJ_0$ with $|S_i|\le \LF$ and $|T_j|\le \LF$ for $i\in [2]$ and $j\in [3]$; recall that $\LF\le\LE$ by~\eqref{eq:CONSTANTS}. 
Combining~\eqref{eq:F-jointneighbourhood} with the fact that 
$(\HStageE, \Vmin,\Vplus)$  is $(\LE,\gamNew, \dStageE^*, \dStageEbar)$-index-quasirandom with respect to the collections $(\imE(s))_{s\in \cJ_0\cup \cJ_1}$ and $(I_s)_{s\in \cJ_0}$, with $I_s:= \{i\in [\lfloor\frac{n}{2}\rfloor]: \boxminus_i,\boxplus_i\in A_s\}$ (see~\ref{quasi:indexquasi}), we obtain
\begin{multline*}
  |\mathbb{U}_{\HStageF}(S_1,S_2,T_1,T_2,T_3)|= |\mathbb{U}_{\HStageE}(S_1,S_2,T_1,T_2,T_3)|\pm (|S_1|+|S_2|)\cdot 2\sigmJjedna n\\
  \begin{aligned}
    &\eqBy{Definition~\ref{def:index-quasirandom}} (1\pm\gamNew)(\dStageE^*)^{|S_1\cap \Vmin|+|S_2\cap \Vplus|}\cdot (\dStageEbar)^{|S_1\cap \Vplus|+|S_2\cap \Vmin|} \\
    & \hspace{8em} \cdot \prod_{s\in T_1\cup T_2}\left(1-\frac{|\imE(s)|}{n}\right)\cdot (2\sigmKJ^*)^{|T_3|}\cdot \frac{n}{2}\pm (|S_1|+|S_2|)2\sigmJjedna n\\
    &=(1\pm \gamJ)(\dStageF^*)^{|S_1\cap \Vmin|+|S_2\cap \Vplus|}(\dStageFbar)^{|S_1\cap \Vplus|+|S_2\cap \Vmin|}\cdot (\deltnonspanning+2\sigmKJ^*)^{|T_1|+|T_2|}\cdot(2\sigmKJ^*)^{|T_3|}\cdot \frac{n}{2}\;,
  \end{aligned}
  \end{multline*}
proving that $(\HStageF, \Vmin,\Vplus)$ is $(\LF, \gamJ, \dStageF^*, \dStageFbar)$-index-quasirandom with respect to the collections $(\imF(s))_{s\in \cJ_0}$ and $(I_s)_{s\in \cJ_0}$, with $I_s:= \{i\in [\lfloor\frac{n}{2}\rfloor]: \boxminus_i,\boxplus_i\in A_s\}$, as required for~\ref{quasi:indexquasi}.

For~\ref{quasi:anchorsets},~\ref{quasi:PathsWithOneVertex},~\ref{quasi:imagecaps},~\ref{quasi:NumAnchors} and~\ref{quasi:NumAnchorsNoIm}, either nothing has changed from Stage~E, or the condition is vacuous since we ask for a quasirandom setup with index sets $\cJ_0,\emptyset,\emptyset$. It remains to verify~\ref{quasi:distributionOfAnchors}, which follows from the above calculation and the observation that for a given $u$ and $v$, the change in any of~\eqref{eq:quasi:distAnchorsa1}--\eqref{eq:quasi:distAnchorsa4} in going from Stage~E to Stage~F is at most $\deg_{\HStageE}(u)-\deg_{\HStageF}(u)$ or $\deg_{\HStageE}(v)-\deg_{\HStageF}(v)$ respectively. By the above calculation, this change is absorbed in the larger error term $\gamJ$. As in Stage~D and Stage~E, \ref{quasi:termVtxNbs} follows by construction from~\ref{quasi:distributionOfAnchors}, so it is established by the same argument, with the same (small) change absorbed into $\gamJ$. This completes the proof that we obtain the required quasirandom setup for Lemma~\ref{lem:StageFNew}.

\smallskip

We now construct the $\phi^\mathbf{F}_s$ for $s\in\cJ_1$. Let $\cP_0:=\bigcup_{s\in\cJ_1}\SpecPaths^*_s$.
Let $v\in \Vmin\cup \Vplus$. Let $\{V_v, \bar{V}_v\}= \{\Vmin,\Vplus\}$\index{$V_v, \bar{V}_v$} be such that $v\in V_v$ and $v\notin \bar{V}_v$. Embedding each path of $\cP_0$ anchored at $v$ will decrease the degree of $v$ by $1$; all other paths will decrease the degree by either $2$ or $0$, depending on whether we use $v$ in the embedding or not. It turns out to be convenient to use the following rule: if a path of $\cP_0$ is anchored at $v$ and its other endpoint is also anchored in $V_v$, then we embed the edge from $v$ in $E(\Vmin,\Vplus)$, while if the other endpoint is anchored in $\bar{V}_v$, then we embed the edge from $v$ within $V_v$. To keep track of how many edges this uses, for a given set of paths $\cP^*\subseteq \cP_0$ we define $j(v, \cP^*)$ to be the number of paths $P\in \cP^*$ anchored at $v$ whose second anchor is in $V_v$, and $\bar{j}(v, \cP^*)$ to be the number of paths $P\in \cP^*$ anchored at $v$ whose second anchor is in $\bar{V}_v$.

To satisfy~\eqref{eq:stageF:degtarget}, by packing $\cP_0$ we need to achieve that $\deg_{\HStageF}(v;V_v)-t(v)=\deg_{\HStageF}(v;\bar{V}_v)$. We denote the \emph{discrepancy degree}\index{discrepancy degree} (defined with respect to a host graph $H^*\subseteq H$ and a set $\cP^*\subseteq \cP_0$) of a vertex $v\in \Vmin\cup\Vplus$  by\index{$d_{H^*,\cP^*}(v)$} \[d_{{H^*},\cP^*}(v)=\deg_{H^*}(v;V_v)-t(v)-\deg_{H^*}(v;\bar{V}_v)+a(v, \cP^*)\;,\] 
where $a(v, \cP^*)=j(v, \cP^*)-\bar{j}(v, \cP^*)$. Observe that by~\ref{quasi:PathsWithOneVertex} we have
\begin{equation}\label{eq:j(v,P_0)}
\begin{split}
j(v,\cP_0)=(1\pm \gamNew)\sigmJjedna^2n\qquad\text{and}\qquad \bar{j}(v, \cP_0)=(1\pm \gamNew)\sigmJjedna^2n\;.
\end{split}
\end{equation} Hence, our ultimate goal is to obtain that $d_{\HStageF, \emptyset}(v)=0$ for all $v\in V(\HStageF)$. 

The following three facts state that the discrepancy degree is $0$ on average, that it does not deviate much from this value, and that it is always even. We will then embed paths one by one, preserving the condition that the average discrepancy degree is zero.

\begin{fact}\label{fact:average-degree-discrepancy}
 We have 
	\[\sum_{v\in \Vmin}d_{\HStageE, \cP_0}(v)=\sum_{v\in \Vplus}d_{\HStageE, \cP_0}(v)=0\;.\]
\end{fact}

\begin{proof}We prove the statement for $\Vmin$; the case of $\Vplus$ is analogous. In the third line below we use the edge counts of~\eqref{eq:E:edges}; that $\sum_{v\in \Vmin}t(v)=\sum_{s\in\cJ_0}|\SpecShortPaths_s|=\sigmJnula\sigmKJ^* n^2$, since by~\ref{D:two} each path of $\SpecShortPaths_s$ is anchored at a terminal pair and so has exactly one anchor in $\Vmin$, while $|\SpecShortPaths_s|=|\SpecPaths^*_s|=\sigmKJ^* n$ by~\eqref{eq:sizeSpecPathStar} and $|\cJ_0|=\sigmJnula n$ by~\eqref{eq:sizeJ0J1}; and that $\sum_{v\in \Vmin}\big(j(v, \cP_0)-\bar{j}(v, \cP_0)\big)=2j_{\boxminus\boxminus}-j_{\boxminus\boxplus}$, since a path of $\cP_0$ with both anchors in $\Vmin$ contributes $1$ to $j(v,\cP_0)$ for each of its two anchors $v$, whereas a path with one anchor in each of $\Vmin$ and $\Vplus$ contributes $1$ to $\bar{j}(v,\cP_0)$ for its anchor $v\in\Vmin$. 
	\begin{align*}
		\sum_{v\in \Vmin}&d_{\HStageE,\cP_0}(v)=\sum_{v\in \Vmin}\left(\deg_{\HStageE}(v;\Vmin)-t(v)-\deg_{\HStageE}(v;\Vplus)+a(v,  \cP_0)\right)\\
		&=2e_{\HStageE}(\Vmin) -e_{\HStageE}(\Vmin,\Vplus)-\sum_{v\in \Vmin}t(v)+\sum_{v\in \Vmin}\left(j(v, \cP_0)-\bar{j}(v, \cP_0)\right)\\
		&\eqBy{Lemma~\ref{lem:StageENew}}2(\sigmJnula\sigmKJ^* n^2+2 j_{\boxminus\boxminus}+3(j_{\boxminus\boxplus}+ j_{\boxplus\boxplus}))-(\sigmJnula\sigmKJ^* n^2+6( j_{\boxminus\boxminus}+ j_{\boxplus\boxplus})+5 j_{\boxminus\boxplus})-\sigmJnula\sigmKJ^* n^2\\
		&+(2j_{\boxminus\boxminus}-j_{\boxminus\boxplus})\\
		&=0\;.\qedhere
	\end{align*}	
\end{proof}

\begin{fact}\label{fact:discrepancy-concentration}
	For each $v\in \Vmin\cup \Vplus$ we have $|d_{\HStageE, \cP_0}(v)|\le \gamNew n$.
\end{fact}

\begin{proof}
	Fix an arbitrary $v\in \Vmin\cup \Vplus$. By definition of $d_{\HStageE, \cP_0}(\cdot)$, we have
	\begin{equation}\label{eq:ladygrey}
	|d_{\HStageE, \cP_0}(v)|=\big|\deg_{\HStageE}(v;V_v)-t(v)-\deg_{\HStageE}(v;\bar{V}_v)+j(v, \cP_0)-\bar{j}(v,\cP_0)\big|\,.
	\end{equation}
	Now, after Stage~E we have a $(\gamNew,\LE,\dStageE^*,\dStageEbar)$-quasirandom setup. By~\eqref{eq:j(v,P_0)}, we have $j(v,\cP_0)-\bar{j}(v,\cP_0)=\pm 2\gamNew \sigmJjedna^2 n$. Also by~\ref{quasi:NumAnchors} we have $t(v)=(1\pm\gamNew)\cdot 2\sigmKJ^*\sigmJnula n$. Finally by~\ref{quasi:indexquasi}, we have $\deg_{\HStageE}(v;V_v)=(1\pm\gamNew)\dStageE^* n/2$ and $\deg_{\HStageE}(v;\bar{V}_v)=(1\pm\gamNew)\dStageEbar n/2$. Plugging this into~\eqref{eq:ladygrey}, and using the definitions of $\dStageE^*$ and $\dStageEbar$, we get
	\begin{align*}
	 &|d_{\HStageE, \cP_0}(v)|=\big|(1\pm\gamNew)\dStageE^*\tfrac{n}{2}-(1\pm\gamNew)2\sigmKJ^*\sigmJnula n-(1\pm\gamNew)\dStageEbar\tfrac{n}{2}\pm 2\gamNew \sigmJjedna^2 n\big|\\
	 =&\big|(1\pm\gamNew)\big(8(\sigmKJ^*\sigmJnula+\tfrac{11}{4}\sigmJjedna^2)\big)\tfrac{n}{2}-(1\pm\gamNew)2\sigmKJ^*\sigmJnula n-(1\pm\gamNew)\big(4(\sigmKJ^*\sigmJnula+\tfrac{11}{2}\sigmJjedna^2)\big)\tfrac{n}{2}\pm 2\gamNew \sigmJjedna^2 n\big|\\
	 \le& \gamNew n\,.\qedhere
	\end{align*}
\end{proof}

\begin{fact}\label{fact:discrepancy-parity}
	For each $v\in \Vmin\cup \Vplus$ the number $d_{\HStageE, \cP_0}(v)$ is even.
\end{fact}

\begin{proof}
	For a leftover graph $H^*$ arising after one of the Stages~B--E, write $p_{H^*}(v)$ for the number of paths still to be embedded at that point which are anchored at $v$; thus $p_{\HStageB}(v)=\PathTerm(v)$, and~\ref{B:parity} says that
	\[\deg_{\HStageB}(v)\equiv p_{\HStageB}(v)\mod 2\,.\]
	We claim that this congruence is maintained in Stages~C, D and~E. Indeed, consider the embedding of a single path, or (in Stage~D) of a single one of the two outer pieces of a path $P\in\SpecPaths^*_s$ with $s\in\cJ_0$, and fix $v\in V(H)$. If $v$ is an internal vertex of the embedded piece, then exactly two edges at $v$ are used and the collection of unembedded paths anchored at $v$ does not change. If $v$ is an endvertex of the embedded piece, then exactly one edge at $v$ is used and $p(v)$ also changes by one: it drops by one if $v$ was an anchor of the path being embedded, and it rises by one if $v$ instead becomes an anchor of the leftover piece of $P$, which is what happens at the terminal pair $\{\boxminus_i,\boxplus_i\}$ in Stage~D. In all remaining cases neither quantity changes.

	Since after Stage~E the paths still to be embedded are precisely those of $\bigcup_{s\in\cJ_0}\SpecShortPaths_s$ and those of $\cP_0$, we conclude that $\deg_{\HStageE}(v)\equiv p_{\HStageE}(v)=t(v)+j(v,\cP_0)+\bar{j}(v,\cP_0)\mod 2$. Using that $\deg_{\HStageE}(v;V_v)+\deg_{\HStageE}(v;\bar{V}_v)=\deg_{\HStageE}(v)$, this gives
	\[d_{\HStageE,\cP_0}(v)\equiv \deg_{\HStageE}(v)+t(v)+j(v,\cP_0)+\bar{j}(v,\cP_0)\equiv 0\mod 2\,.\qedhere\]
\end{proof}

The mappings $(\phi^{\mathbf{F}}_s)_{s\in \cJ_1}$ will be defined in roughly $\sigmJjedna^2n^2$ steps, where in each step $i$ we define a mapping $\phi_i$ of one path $P_i\in \cP_0$. The reason for `roughly' is that in Stage~C we packed some (between $0$ and $n^{0.6}$) paths from each $\SpecPaths_s$ with $s\in\cJ_1$. Let $m:=\sum_{s\in\cJ_1}|\SpecPaths^*_s|$ denote the precise number of paths we need to pack in this stage; then we have $\tfrac12\sigmJjedna^2n^2\le m\le\sigmJjedna^2n^2$.

Denote by $\cP_i$ the set of paths of $\cP_0$ that have not yet been packed after step~$i$, in particular  $\cP_{m}=\emptyset$.
Let $\im_0(s):=\imE(s)$ be the image of the partial embedding of $G_s$ at the beginning of Stage~F. Let us assume that we are in step $i\in [m]$ and are about to embed a path $P_i\in \cP_{i-1}$. For $s\in \cJ_1$, denote by $\im_i(s)$ the image of $G_s$  just after step~$i$, i.e., $\im_i(s)=\im_{i-1}(s)$ if $P_i\notin \SpecPaths^*_s$ and $\im_i(s)=\im_{i-1}(s)\cup \phi_i(V(P_i))$, otherwise. Let $H_i$ be the host graph just after step~$i$ (i.e., the union of edges not yet used by the packing after step~$i$, the graph $H_0$ being $\HStageE$). 
 For the first (at most) $\gamNew n^2$ steps, we shall be correcting the degree discrepancy of the vertices. During these steps, we shall be careful not to use any one vertex more than $\gamNew\big(1+20(\sigmJnula\sigmKJ^*)^{-2}\deltnonspanning^{-1}\big)n$ times.
  After that, we shall just pack the rest of $\cP_0$ in such a way as not to change the degree discrepancy of any vertex.

First, we state an auxiliary fact that will ensure the existence of the mappings $\phi_i$, $i\in [m]$, that we shall define in the proofs of Lemmas~\ref{lem:repairDegreeDiscrepancy} and~\ref{lem:neutralembedding}. 

	\begin{fact}\label{fact:repair-largedegree}
	For any $i\le \sigmJjedna^2n^2$, for any $s\in \cJ_1$, for any set $S\subseteq V(H)$ with $|S|\le 2$, and for $V_*\in \{\Vmin,\Vplus\}$, we have 
	$|\NBH_{H_{i-1}}(S)\cap V_*\setminus \im_{i-1}(s)|>(\sigmJnula\sigmKJ^*)^2\deltnonspanning n$.
\end{fact}

\begin{proof}
	From~\ref{quasi:indexquasi}, by setting $S_1:=S$, $S_2:=\emptyset$, $T_1:=\{s\}$, and  $T_2=T_3:=\emptyset$, if $V_*=\Vmin$, and setting $S_1:=\emptyset$, $S_2:=S$, $T_1:=\emptyset$, $T_2:=\{s\}$, and $T_3:=\emptyset$, if $V_*=\Vplus$,  we obtain that 
	\[|\NBH_{\HStageE}(S)\cap V_*\setminus \imE(s)|\ge \frac{1}{2}\cdot\min\{\dStageE^*,\dStageEbar\}^2\cdot \frac{n}{2}\cdot \deltnonspanning \ge 4 (\sigmJnula\sigmKJ^*)^2\deltnonspanning n\,,\]
	where the last inequality comes from Fact~\ref{fact:min:d_E}.
	
	Observe that for an arbitrary $v\in \Vmin\cup \Vplus$, the difference between $|\NBH_{\HStageE}(v)|$ and $|\NBH_{H_{i-1}}(v)|$ is at most $2\sigmJjedna n$, as $|\cJ_1|=\sigmJjedna n$ and, for each $s\in\cJ_1$, at most $2$ edges incident to $v$ are used in Stage~F for $G_s$, since $v$ is the image of at most one vertex of $G_s$. Also observe that the difference between $|\imE(s)|$ and $|\im_{i-1}(s)|$ is at most $10 \sigmJjedna n$: these are exactly the vertices of paths of $\SpecPaths^*_s$ embedded so far in Stage~F, each of these paths has $10$ vertices to embed, and $|\SpecPaths^*_s|\le|\SpecPaths_s|=\sigmJjedna n$. 
Therefore, we have by~\eqref{eq:CONSTANTS}
\begin{align*}
	|\NBH_{H_{i-1}}(S)\cap V_*\setminus \im_{i-1}(s)|&\ge |\NBH_{\HStageE}(S)\cap V_*\setminus \imE(s)|-4\sigmJjedna n-10\sigmJjedna n\\
	&\ge 4 (\sigmJnula\sigmKJ^*)^2\deltnonspanning n - 14\sigmJjedna n>(\sigmJnula\sigmKJ^*)^2\deltnonspanning n\;.\qedhere
\end{align*}	
\end{proof}

In Lemma~\ref{lem:repairDegreeDiscrepancy} below we claim that we can embed a path from $\cP_{i-1}$ in such a way that we can improve by $2$ the degree discrepancy of two vertices, without changing the degree discrepancy of any other vertex (and thus keeping the average degree discrepancy equal to zero as it was at the start, see Fact~\ref{fact:average-degree-discrepancy}). By Fact~\ref{fact:discrepancy-parity} the degree discrepancy of each vertex is even, so that repeating this procedure can indeed terminate with all vertices having degree discrepancy zero.

\begin{lemma}\label{lem:repairDegreeDiscrepancy}
Let $i\le \gamNew n^2$.		
Suppose we have two vertices $u,w\in V_a$ where $a\in \{\boxminus, \boxplus\}$ such that $d_{H_{i-1},\cP_{i-1}}(u)>0$ and $d_{H_{i-1}, \cP_{i-1}}(w)<0$. 

	Then there exists a path $P_i\in \cP_{i-1}$ and there exists a mapping $\phi_i$ of $P_i$ in $H_{i-1}$ such that $d_{H_{i}, \cP_i}(u)=d_{H_{i-1}, \cP_{i-1}}(u)-2$, $d_{H_{i},\cP_i}(w)=d_{H_{i-1}, \cP_{i-1}}(w)+2$, and $d_{H_{i},\cP_i}(v)=d_{H_{i-1}, \cP_{i-1}}(v)$ for all $v\in V(H)\setminus \{u,w\}$ and $\phi_i$ satisfies the following. 
	\begin{enumerate}[label=\abc]
	\item \label{it:F_Discr1}If one anchor of $P_i$ is in $\Vmin$ and the other in $\Vplus$, then the image of $\phi_i$ uses $3$ edges in $\Vmin$, $3$ edges in $\Vplus$, and $5$ edges between $\Vmin$ and $\Vplus$.
	\item \label{it:F_Discr2} If both anchors are in $\Vmin$, then $\phi_i$ uses $2$ edges in $\Vmin$, $3$ edges in $\Vplus$, and $6$ edges between $\Vmin$ and $\Vplus$.
	\item \label{it:F_Discr3}If both anchors are in $\Vplus$, then $\phi_i$ uses $2$ edges in $\Vplus$, $3$ edges in $\Vmin$, and $6$ edges between $\Vmin$ and $\Vplus$.
	\end{enumerate}
Moreover, throughout the at most $\gamNew n^2$ applications of this lemma, no vertex is used more than $\gamNew(1+20(\sigmJnula\sigmKJ^*)^{-2}\deltnonspanning^{-1})n$ times for the embedding of these paths. 
\end{lemma}

\begin{proof} The next auxiliary claim helps us to pick a suitable path $P_i$ from the set $\cP_{i-1}$.
	\begin{claim}\label{cl:existence-Pi}
	Suppose that $i\le \gamNew n^2$. Then there exists a path $P_i\in \cP_{i-1}\cap \SpecPaths^*_s$ with $s$ such that $\{u,w\}\cap \im_{i-1}(s)=\emptyset$.
	\end{claim}
\begin{claimproof}[Proof of Claim~\ref{cl:existence-Pi}]
	From~\ref{quasi:imagecaps}, by setting $S_1:=\{u,w\}$, $S_2=\emptyset$ and looking in $\cJ_1$ (where for each $s\in\cJ_1$ we have $n-|\imE(s)|\ge\deltnonspanning n$), we obtain that 
	\[|\{s\in \cJ_1\::\: \imE(s)\cap \{u,w\}= \emptyset\}|\ge (1-2\gamNew)\deltnonspanning^2\cdot \sigmJjedna n\ge\tfrac12\deltnonspanning^2\cdot\sigmJjedna n\;.\]
	At step $i\le \gamNew n^2$, each of $u$ and $w$ has been used at most $\gamNew\big(1+20(\sigmJnula\sigmKJ^*)^{-2}\deltnonspanning^{-1}\big)n$ times in Stage~F, and in particular each of them lies in $\im_{i-1}(s)\setminus\imE(s)$ for at most that many indices $s$. Since $|\SpecPaths^*_s|\ge\tfrac12\sigmJjedna n$ for each $s\in\cJ_1$, and since we have embedded at most $\gamNew n^2$ paths so far, we have at least $\big(\tfrac12\deltnonspanning^2\sigmJjedna n-2\gamNew\big(1+20(\sigmJnula\sigmKJ^*)^{-2}\deltnonspanning^{-1}\big)n\big)\cdot\tfrac12\sigmJjedna n-\gamNew n^2>0$ paths to choose $P_i$ from.
\end{claimproof}
	
	Assume that $u,w\in \Vmin$. The other case is dealt with symmetrically and we omit the details. 
Let $P_i=x_1x_2\cdots x_{12}$ be a path given by Claim~\ref{cl:existence-Pi}, let $s$ be such that $P_i\in \SpecPaths_s^*$, and let $v_1:=\phi^{\mathbf{E}}_s(x_1)$ and $v_{12}:=\phi^{\mathbf{E}}_s(x_{12})$. 
	
	First consider the case when $v_1\in \Vmin$ and $v_{12}\in \Vplus$ (or the other way around, which is done analogously). We set $v_5=\phi_i(x_5):=w$ and $v_9=\phi_i(x_9):=u$. Each vertex $x_j\in V(P_i)\setminus\{x_1,x_5,x_9,x_{12}\}$ (which we call the \emph{connection vertices})
	is successively mapped to $v_j$ chosen in  
	\begin{equation}\label{eq:repair-simpleneighbourhood}
\NBH_{H_{i-1}}(v_{j-1})\cap V_*\setminus \left(\im_{i-1}(s)\cup\{v_2, \ldots, v_{j-1}\}\cup \{u,w\}\right)
	\end{equation}
	or in 
	\begin{equation}\label{eq:repair-doubleneighbourhood}
	\NBH_{H_{i-1}}(v_{j-1})\cap \NBH_{H_{i-1}}(v_{j+1})\cap V_*\setminus \left(\im_{i-1}(s)\cup\{v_2, \ldots, v_{j-1}\}\cup \{u,w\}\right), 
	\end{equation} if $v_{j+1}$ is already defined, where $V_*=\Vmin$ for $j\in \{2,8,10\}$ and $V_*=\Vplus$ for $j\in \{3,4,6,7,11\}$. We insist additionally (and likewise in the two cases below) that at each of these steps we choose a vertex to which fewer than $20\gamNew(\sigmJnula\sigmKJ^*)^{-2}\deltnonspanning^{-1} n$ connection vertices have so far been embedded in Stage~F. By Fact~\ref{fact:repair-largedegree}, before this additional restriction the set from which $v_j$ is chosen has more than $(\sigmJnula\sigmKJ^*)^2\deltnonspanning n$ elements. The total number of paths we embed while applying this lemma is at most $\gamNew n^2$ and each has $10$ vertices to embed, so at most $10\gamNew n^2$ connection vertices are embedded in total, and hence the restriction excludes at most $\tfrac12(\sigmJnula\sigmKJ^*)^2\deltnonspanning n$ vertices. Excluding in addition the at most $10$ vertices of $\{v_2,\ldots,v_{j-1}\}\cup\{u,w\}$, the set from which $v_j$ is chosen is non-empty. Note that when embedding paths in this lemma, we use vertices either as connection vertices or because their degree discrepancy is non-zero; by Fact~\ref{fact:discrepancy-concentration} we use any vertex at most $\gamNew n$ times for the latter, and this implies that we use any vertex in total at most $\gamNew\big(1+ 20(\sigmJnula\sigmKJ^*)^{-2}\deltnonspanning^{-1} \big)n$ times while applying this lemma. This satisfies the \emph{Moreover} part of the lemma.

Observe that
\begin{align*}
  \deg_{\phi_i(P_i)}(u;\Vmin)&=2\,, \quad \deg_{\phi_i(P_i)}(u;\Vplus)=0\,, \\
  \deg_{\phi_i(P_i)}(w;\Vmin)&=0\,, \quad \deg_{\phi_i(P_i)}(w;\Vplus)=2\,, \\
  \deg_{\phi_i(P_i)}(v_j;\Vmin)&=\deg_{\phi_i(P_i)}(v_j;\Vplus)=1\,,
\end{align*}
for $j\in \{2,3,4,6,7,8,10,11\}$, and $\deg_{\phi_i(P_i)}(v_1;\Vmin)=1$, $\deg_{\phi_i(P_i)}(v_{12};\Vplus)=1$, while there is one fewer path anchored at $v_1\in \Vmin$ and $v_{12}\in \Vplus$ in $\cP_i$ than in $\cP_{i-1}$. This ensures that the requirements for the degree discrepancy $d_{H_i, \cP_i}(v)$ are fulfilled for all $v\in \Vmin\cup \Vplus$. 
	 Observe that $\phi_i$ satisfies~\ref{it:F_Discr1}, as it uses the following $5$ edges between $\Vmin$ and $\Vplus$: $\{v_2,v_3\}, \{v_4,v_5\}, \{v_5, v_6\}, \{v_7, v_8\}$, and $\{v_{10},v_{11}\}$, 
	 the following $3$ edges inside $\Vmin$: $\{v_1,v_2\},\{v_8,v_9\}$, and $\{v_9,v_{10}\}$, 
	 and the following $3$ edges inside $\Vplus$: $\{v_3,v_4\}, \{v_6,v_7\}$, and $\{v_{11},v_{12}\}$.

	Second, consider the case when $v_1,v_{12}\in \Vmin$. We set $v_4=\phi_i(x_4):=w$ and $v_8=\phi_i(x_8):=u$.  Each vertex $x_j\in V(P_i)\setminus\{x_1,x_4,x_8,x_{12}\}$ 
	is successively mapped to $v_j$ chosen  as in~\eqref{eq:repair-simpleneighbourhood} or~\eqref{eq:repair-doubleneighbourhood} where $V_*=\Vmin$ for $j\in \{7,9\}$ and $V_*=\Vplus$ for $j\in\{2,3,5,6,10,11\}$. As above, Fact~\ref{fact:repair-largedegree} ensures that we can define $v_j$ without any vertex having too many connection vertices embedded to it, and the assignment of $V_*$ ensures that the requirements for $d_{H_i, \cP_i}(\cdot)$ are satisfied. Observe that $\left\{\{v_1,v_2\},\{v_3,v_4\},\{v_4,v_5\},\{v_6,v_7\},\{v_9,v_{10}\},\{v_{11},v_{12}\}\right\}\subset \Vmin\times\Vplus$, $\left\{\{v_7,v_8\},\{v_8,v_9\}\right\}\subseteq \binom{\Vmin}{2}$, and $\left\{\{v_2,v_3\},\{v_5,v_6\},\{v_{10},v_{11}\}\right\}\subseteq \binom{\Vplus}{2}$, satisfying~\ref{it:F_Discr2}.

	Finally, consider the case when $v_1,v_{12}\in \Vplus$. We set $v_3=\phi_i(x_3):=u$ and $v_7=\phi_i(x_7):=w$. Each vertex $x_j\in V(P_i)\setminus \{x_1,x_3,x_7,x_{12}\}$ is successively mapped to $v_j$ chosen  in~\eqref{eq:repair-simpleneighbourhood} or in~\eqref{eq:repair-doubleneighbourhood}, where $V_*=\Vmin$ for $j\in \{2,4,10,11\}$ and $V_*=\Vplus$ for $j\in\{5,6,8,9\}$. Again we can define $v_j$ due to  Fact~\ref{fact:repair-largedegree} and the requirements of $d_{H_i, \cP_i}(\cdot)$ are satisfied. This establishes~\ref{it:F_Discr3}: $\phi_i$ uses $6$ edges between $\Vmin$ and $\Vplus$, $3$ edges inside $\Vmin$, and $2$ edges inside $\Vplus$.
	\end{proof}

\begin{lemma}\label{lem:neutralembedding}
	Let $i\in [m]$.		
	For any path $P_i\in \cP_{i-1}$, there exists a mapping $\phi_i$ of $P_i$ in $H_{i-1}$ such that $d_{H_{i},\cP_i}(v)=d_{H_{i-1}, \cP_{i-1}}(v)$ for all $v\in \Vmin\cup\Vplus$.
	Moreover, $\phi_i$ satisfies the following. 
	\begin{itemize}
		\item If one anchor of $P_i$ is in $\Vmin$ and the other in $\Vplus$, then the image of $\phi_i$ uses $3$ edges in $\Vmin$, $3$ edges in $\Vplus$, and $5$ edges between $\Vmin$ and $\Vplus$.
		\item If both anchors are in $\Vmin$, then $\phi_i$ uses $2$ edges in $\Vmin$, $3$ edges in $\Vplus$, and $6$ edges between $\Vmin$ and $\Vplus$.
		\item If both anchors are in $\Vplus$, then $\phi_i$ uses $2$ edges in $\Vplus$, $3$ edges in $\Vmin$, and $6$ edges between $\Vmin$ and $\Vplus$.
	\end{itemize}
\end{lemma}

\begin{proof}
 Let $P_i=x_1x_2\cdots x_{12}$ and let $s$ be such that $P_i\in \SpecPaths^*_s$ and let $v_1:=\phi^{\mathbf{E}}_s(x_1)$ and $v_{12}:=\phi^{\mathbf{E}}_s(x_{12})$. 
 
 	First, consider the case when $v_1\in \Vmin$ and $v_{12}\in \Vplus$ (or the other way around, which is done analogously). Each vertex $x_j\in V(P_i)\setminus\{x_1,x_{12}\}$ 
 	is successively mapped to $v_j$ chosen in  
 	\eqref{eq:repair-simpleneighbourhood}, or in \eqref{eq:repair-doubleneighbourhood} if $j=11$ (in both cases with the term $\{u,w\}$ omitted, since no vertices of $P_i$ are prescribed here), where $V_*=\Vmin$ for $j\in \{2,5,6,9,10\}$ and $V_*=\Vplus$ for $j\in \{3,4,7,8,11\}$.
 	By Fact~\ref{fact:repair-largedegree} the set from which $v_j$ is chosen is non-empty, as $\phi_i$ is eventually defined on $10$ vertices. As in the proof of Lemma~\ref{lem:repairDegreeDiscrepancy}, it is easy to observe that the requirements on $d_{H_i, \cP_i}(\cdot)$ are satisfied. For the \emph{Moreover} part, observe that $\phi_i$ uses the edges $\{v_2,v_3\}, \{v_4,v_5\}, \{v_6,v_7\}, \{v_8,v_9\}$, and $\{v_{10},v_{11}\}$ from $\Vmin\times\Vplus$. The edges $\{v_1,v_2\}, \{v_5,v_6\}$, and $\{v_9,v_{10}\}$ lie inside $\Vmin$ and the edges $\{v_3,v_4\}, \{v_7,v_8\}$, and $\{v_{11},v_{12}\}$ lie inside $\Vplus$.

 	Second, consider the case when $v_1, v_{12}\in \Vmin$ (the case when $v_1, v_{12}\in \Vplus$ is done analogously). Each vertex $x_j\in V(P_i)\setminus\{x_1,x_{12}\}$ 
 	is successively mapped to $v_j$ chosen in  
 	\eqref{eq:repair-simpleneighbourhood}, or in \eqref{eq:repair-doubleneighbourhood} if $j=11$, where $V_*=\Vmin$ for $j\in \{4,5,8,9\}$ and $V_*=\Vplus$ for $j\in \{2,3,6,7,10,11\}$.
 	By Fact~\ref{fact:repair-largedegree} the set from which $v_j$ is chosen is non-empty, as $\phi_i$ is eventually defined on $10$ vertices. As in the proof of Lemma~\ref{lem:repairDegreeDiscrepancy}, it is easy to observe that the requirements on $d_{H_i, \cP_i}(\cdot)$ are satisfied. As just above, one can easily check that $\phi_i$ uses $2$ edges in $\Vmin$, $3$ edges in $\Vplus$, and $6$ edges between $\Vmin$ and $\Vplus$.
\end{proof}

\begin{fact}\label{fact:discrepancy-invariant}
Throughout the process of repeatedly applying Lemma~\ref{lem:repairDegreeDiscrepancy}, every discrepancy $d_{H_i,\cP_i}(v)$ remains even, and the discrepancies sum to $0$ on each of $\Vmin$ and $\Vplus$ separately. Consequently, as long as $d_{H_i,\cP_i}(v)\ne0$ for some $v$, Lemma~\ref{lem:repairDegreeDiscrepancy} is applicable.
\end{fact}
\begin{proof}
At $i=0$, every discrepancy is even by Fact~\ref{fact:discrepancy-parity}, and by Fact~\ref{fact:average-degree-discrepancy} the discrepancies sum to $0$ on each of $\Vmin$ and $\Vplus$ separately. Each application of Lemma~\ref{lem:repairDegreeDiscrepancy} changes two same-side discrepancies by $\mp2$, which preserves both properties. Given these, if $d_{H_i,\cP_i}(v)\ne0$ for some $v\in V_a$, then, as the discrepancies on $V_a$ sum to $0$, some vertex of $V_a$ must have discrepancy of the opposite sign to $v$, giving the hypothesis of Lemma~\ref{lem:repairDegreeDiscrepancy}.
\end{proof}

Once we have $d_{H_i, \cP_i}(v)=0$ for all $v\in \Vmin\cup \Vplus$, we use Lemma~\ref{lem:neutralembedding} to pack the leftover paths from $\cP_i$. So it remains to prove that we manage to correct all the degree discrepancies within the first $\gamNew n^2$ steps by using Lemma~\ref{lem:repairDegreeDiscrepancy}. By Fact~\ref{fact:discrepancy-invariant}, Lemma~\ref{lem:repairDegreeDiscrepancy} remains applicable throughout this process, and since discrepancies stay even, the process indeed terminates with $d_{H_i,\cP_i}(v)=0$, rather than merely small, for every $v$. Observe that by Fact~\ref{fact:discrepancy-concentration} we have
\[\sum_{v\in \Vmin\cup \Vplus}|d_{H_0,\cP_0}(v)|\le \gamNew n^2\;.\]
Each application of Lemma~\ref{lem:repairDegreeDiscrepancy} at a step~$i$ gives $|d_{H_i, \cP_i}(v)|=|d_{H_{i-1}, \cP_{i-1}}(v)|-2$ for the two vertices $u,w$ to which the lemma is applied at that step, and $d_{H_i, \cP_i}(v)=d_{H_{i-1}, \cP_{i-1}}(v)$ for every other $v$; here we use that $d_{H_{i-1}, \cP_{i-1}}(u)>0>d_{H_{i-1}, \cP_{i-1}}(w)$ and that both are even by Fact~\ref{fact:discrepancy-invariant}, so that both absolute values indeed decrease. Hence
\[\sum_{v\in \Vmin\cup \Vplus}|d_{H_i, \cP_i}(v)|=\sum_{v\in \Vmin\cup \Vplus}|d_{H_{i-1}, \cP_{i-1}}(v)|-4\;.\]
Together with the previous display, it follows that after at most $\tfrac14\gamNew n^2$ steps this sum vanishes, that is, $d_{H_i, \cP_i}(v)=0$ for every $v\in \Vmin\cup \Vplus$. This completes the proof of Lemma~\ref{lem:StageFNew}.

\section{Stage~G (Proof of Lemma~\ref{lem:StageG})}\label{sec:StageG}
In this section we shall prove Lemma~\ref{lem:StageG} by verifying the hypotheses of Proposition~\ref{prop:Designs-oursetting} for an auxiliary multigraph $\mathcal M$.

Proposition~\ref{prop:Designs-oursetting} with input parameters $d_{P\ref{prop:Designs-oursetting}}:= \sigmKJ\sigmJnula$ and $\sigma_{P\ref{prop:Designs-oursetting}}:= \sigmJnula$ outputs parameters $L_{P\ref{prop:Designs-oursetting}}$, $n_{0,P\ref{prop:Designs-oursetting}}$ and $\gamma_{P\ref{prop:Designs-oursetting}}$.
By~\eqref{eq:CONSTANTS}, we have $n\ge 2n_{0,P\ref{prop:Designs-oursetting}}$, $\LF\ge  L_{P\ref{prop:Designs-oursetting}}$, $5\LF\gamJ\le\gamma_{P\ref{prop:Designs-oursetting}}$ and $\gamJ\le \sigmJnula\sigmKJ\deltnonspanning/8$.

From the leftover graph~$\HStageF$ and the paths of $\bigcup_{s\in \cJ_0}\SpecShortPaths_s$, we define
\[\mathcal M=([\lfloor n/2 \rfloor]\dcup \cJ_0; \vec{E}_1, E_2, E_3, E_4, E_5, E_6)\]
with $\vec{E}_1\subseteq [\lfloor n/2 \rfloor]^2$, 
$E_2, E_3\subseteq \binom{[\lfloor n/2 \rfloor]}{2}$, 
and $E_4, E_5, E_6\subseteq [\lfloor n/2 \rfloor]\times \cJ_0$, as follows. Each pair $\boxminus_i, \boxplus_i$ 
corresponds to the vertex $i\in [\lfloor n/2 \rfloor]$. For each edge $\{\boxminus_i, \boxplus_j\}\in E(\HStageF)$ insert an oriented edge $(i,j)\in \vec{E}_1$. For each edge $\{\boxminus_i, \boxminus_j\}\in E(\HStageF)$ insert an edge $\{i,j\}\in E_2$. For each edge $\{\boxplus_i, \boxplus_j\}\in E(\HStageF)$ insert an edge $\{i,j\}\in E_3$. Insert an edge $\{i,s\}\in E_6$ if there is a path $P\in \SpecShortPaths_s$ with $(\phi^\mathbf{F}_s)^{-1}(\{\boxminus_i,\boxplus_i\})=\{\leftpath_0(P),\rightpath_0(P)\}$; since by~\ref{D:two} every path of $\SpecShortPaths_s$ is anchored at a terminal pair, this is the case precisely when $i\in I_s$.
Insert an edge $\{i,s\}\in E_4$ if $\boxminus_i\notin \imF(s)$, and an edge $\{i,s\}\in E_5$ if $\boxplus_i\notin \imF(s)$.

As $\HStageF$ is a simple graph, none of the edge sets $E_i$, $i\in\{2,\ldots,6\}$, nor $\vec{E}_1$ has multiple edges that would be oriented in the same direction. As none of the edges $\{\boxminus_i, \boxplus_i\}$ with $i\in [\lfloor n/2 \rfloor]$ is present in $E(\HStageF)$ (by~\ref{C:three} every such edge of $H$ was already used in Stage~C), there is no loop in $\vec{E}_1$. Hence, the multigraph $\mathcal M$ is a chest in the sense of Definition~\ref{def:chest}.

Moreover, the edge counts $e_{\HStageF}(\Vmin)=e_{\HStageF}(\Vplus)=e_{\HStageF}(\Vmin,\Vplus)=\sigmJnula\sigmKJ^* n^2$ established in Section~\ref{sec:StageF} give
\[|\vec{E}_1|=|E_2|=|E_3|=|E_6|=\Big|\bigcup_{s\in \cJ_{0}}\SpecShortPaths_s\Big|=\sigmJnula\sigmKJ^* n^2\,,\]
as Proposition~\ref{prop:Designs-oursetting} requires. Recall also that $|\imF(s)|=(1-\deltnonspanning-2\sigmKJ^*)n$ for every $s\in\cJ_0$, since $\phi^\mathbf{F}_s=\phi^\mathbf{D}_s$ leaves exactly the $2|\SpecPaths^*_s|=2\sigmKJ^*n$ middle vertices of the paths $\SpecPaths^*_s$ unembedded. Hence the densities of Definition~\ref{def:quasichest} satisfy
\[d_1=\big(1+O(n^{-1})\big)\dStageFbar\,,\qquad d_2=d_3=\big(1+O(n^{-1})\big)\dStageF^*\,,\qquad d_6=\big(1+O(n^{-1})\big)2\sigmKJ^*\,,\]
while $d_4,d_5=(1\pm\gamJ)\big(1+O(n^{-1})\big)(\deltnonspanning+2\sigmKJ^*)$, the extra factor $1\pm\gamJ$ arising because the free space of $G_s$ is split between $\Vmin$ and $\Vplus$ only up to that factor. In particular $d_1=\big(1+O(n^{-1})\big)4\sigmKJ^*\sigmJnula>\sigmKJ\sigmJnula$ and $d_4,d_5\ge d_6>\sigmKJ\sigmJnula$. Finally, $|V(\mathcal M)|=\lfloor n/2\rfloor+|\cJ_{0}|>n_{0,P\ref{prop:Designs-oursetting}}$ and $\min\{\lfloor n/2\rfloor,|\cJ_{0}|\}=\sigmJnula n>\sigmJnula|V(\mathcal M)|$.

\begin{lemma}\label{lem:chest-degree}
The chest $\mathcal M$ fulfils the degree conditions from Proposition~\ref{prop:Designs-oursetting}.
\end{lemma}
\begin{proof}
Fix any vertex $i\in \left[ \lfloor\frac{n}{2}\rfloor \right]$. Condition~\ref{itm:designs:i} of Proposition~\ref{prop:Designs-oursetting} that $\deg_{E_2}(i)=\degout_{\vec{E}_1}(i)+\deg_{E_6}(i)$ translates to $\deg_{\HStageF}(\boxminus_i; \Vmin)=\deg_{\HStageF}(\boxminus_i; \Vplus)+|\{s\in \cJ_0 : i\in I_s\}|$, with $I_s$ as in Lemma~\ref{lem:StageFNew}\ref{F:four}. This is given by~\eqref{eq:stageF:degtarget} holding as stated in Lemma~\ref{lem:StageFNew}. Condition~\ref{itm:designs:ii} of Proposition~\ref{prop:Designs-oursetting} that $\deg_{E_3}(i)=\degin_{\vec{E}_1}(i)+\deg_{E_6}(i)$ translates to $\deg_{\HStageF}(\boxplus_i;\Vplus)=\deg_{\HStageF}(\boxplus_i; \Vmin)+|\{s\in \cJ_0 : i\in I_s\}|$, again given by~\eqref{eq:stageF:degtarget}. 
As for Conditions~\ref{itm:designs:iii} and~\ref{itm:designs:iv} of Proposition~\ref{prop:Designs-oursetting}, observe that, using the notation of Definition~\ref{def:index-quasirandom}, $\deg_{\HStageF}(\boxminus_i;\Vplus)=|\mathbb{U}_{\HStageF}(\emptyset,\{\boxminus_i\},\emptyset,\emptyset,\emptyset)|$ and $\deg_{\HStageF}(\boxplus_i;\Vmin)=|\mathbb{U}_{\HStageF}(\{\boxplus_i\},\emptyset,\emptyset,\emptyset,\emptyset)|$. Hence, since by~\ref{quasi:indexquasi} the graph $\HStageF$ with the partition $(\Vmin,\Vplus)$ is $(\LF,\gamJ, \dStageF^*, \dStageFbar)$-index-quasirandom with respect to the collections $(\imF(s))_{s\in \cJ_0}$ and $(I_s)_{s\in \cJ_0}$, we get
\[\deg_{\HStageF}(\boxminus_i;\Vplus)\,,\ \deg_{\HStageF}(\boxplus_i;\Vmin)\ \le\ (1+\gamJ)\dStageFbar\tfrac{n}{2}=(1+\gamJ)\cdot2\sigmKJ^*\sigmJnula n<\tfrac14\deltnonspanning\sigmJnula n\;,\]
using $\dStageFbar=4\sigmKJ^*\sigmJnula$ and $\sigmKJ^*\le\sigmKJ\ll\deltnonspanning$ from~\eqref{eq:CONSTANTS}.

We need to show
\[|\{s\in \cJ_0 : \boxminus_i\notin \imF(s)\}|\ge \deg_{\HStageF}(\boxminus_i;\Vplus) +\sigmKJ^5 \sigmJnula^7|\cJ_0| \eqByRef{eq:sizeJ0J1} \deg_{\HStageF}(\boxminus_i;\Vplus) +\sigmKJ^5 \sigmJnula^8n\;,\]
and the analogous statement with $\boxminus_i,\Vplus$ replaced by $\boxplus_i,\Vmin$. From~\ref{quasi:imagecaps}, applied with $S_1=\{\boxminus_i\}$ (resp.\ $S_1=\{\boxplus_i\}$) and $T=\emptyset$, each left-hand side is at least $\frac 12 (\deltnonspanning+2\sigmKJ^*)\sigmJnula n\ge\tfrac12\deltnonspanning\sigmJnula n$. By~\eqref{eq:CONSTANTS}, this exceeds the bound $\tfrac14\deltnonspanning\sigmJnula n$ on $\deg_{\HStageF}(\boxminus_i;\Vplus)$ (resp.\ $\deg_{\HStageF}(\boxplus_i;\Vmin)$) established above by more than $\sigmKJ^5\sigmJnula^8n$, giving Conditions~\ref{itm:designs:iii} and~\ref{itm:designs:iv} of Proposition~\ref{prop:Designs-oursetting}.

Now fix any vertex $s\in \cJ_0$. Observe that, using the notation of Definition~\ref{def:index-quasirandom},
 \begin{align*}
 \deg_{E_4}(s)&=|\Vmin\setminus \imF(s)|=|\mathbb{U}_{\HStageF}(\emptyset, \emptyset, \{s\}, \emptyset, \emptyset)|\,,\\
 \deg_{E_5}(s)&=|\Vplus\setminus \imF(s)|=|\mathbb{U}_{\HStageF}(\emptyset, \emptyset, \emptyset, \{s\}, \emptyset)|\,,\text{ and}\\
 \deg_{E_6}(s)&=|I_s|=|\mathbb{U}_{\HStageF}(\emptyset, \emptyset,\emptyset,  \emptyset, \{s\})|\,.
 \end{align*}
 Hence, by~\ref{quasi:indexquasi} again, we have
 $|\Vmin\setminus \imF (s)|=(1\pm \gamJ)(\deltnonspanning+2\sigmKJ^*)\frac{n}{2}$, $|\Vplus\setminus \imF (s)|=(1\pm \gamJ)(\deltnonspanning+2\sigmKJ^*)\frac{n}{2}$, and we have $|I_s|=(2\sigmKJ^*)\frac{n}{2}$. Together with~\eqref{eq:CONSTANTS}, these imply~\ref{itm:designs:v} and~\ref{itm:designs:vi} of Proposition~\ref{prop:Designs-oursetting}.
\end{proof}

To apply Proposition~\ref{prop:Designs-oursetting}, it remains to prove the quasirandomness of the chest $\mathcal M$ (see Definition~\ref{def:quasichest}).
\begin{lemma}\label{lem:chest-quasirandomness}
 The chest $\mathcal{M}$ is $(5\LF\gamJ,\LF)$-quasirandom.
\end{lemma}
\begin{proof}
The $(5\LF\gamJ,\LF)$-quasirandomness of the chest corresponds to the following conditions in the leftover graph $\HStageF$.
 Suppose we are given mutually disjoint sets $X_1, X'_1, X_2, X_3, X'_4, X'_5, X'_6\subseteq [\lfloor n/2 \rfloor]$ and $X_4,X_5,X_6\subseteq \cJ_0$ of total size at most $\LF$. Let $d_i$ be the density of $E_i$ as in Definition~\ref{def:quasichest}.  We need to show that 
\begin{align}
\label{eq:conjquasi-V}| \NBHout_{\vec{E}_1}(X_1)\cap \NBHin_{\vec{E}_1}(X'_1)\cap \bigcap_{i=2}^6\NBH_{E_i}(X_i)|&=(1\pm 5\LF\gamJ)d_1^{|X_1|+|X'_1|}\cdot \prod_{i=2}^{6} d_i^{|X_i|}\frac{n}{2}\;\mbox{, and}\\
\label{eq:conjquasi-U}|\cJ_0 \cap \bigcap_{i=4}^6\NBH_{E_i}(X'_i)|&=(1\pm 5\LF\gamJ)\prod_{i=4}^{6}d_i^{|X'_i|}|\cJ_0|\;.
\end{align}
 Define sets $Q_1, Q_2, Q'_4\subseteq \Vmin$ and $Q'_1, Q_3, Q'_5\subseteq \Vplus$ by $Q_1=\{\boxminus_i:i\in X_1\}$, $Q'_1=\{\boxplus_i:i\in X'_1\}$, $Q_2=\{\boxminus_i:i\in X_2\}$, $Q_3=\{\boxplus_i:i\in X_3\}$, $Q'_4=\{\boxminus_i:i\in X'_4\}$, and $Q'_5=\{\boxplus_i:i\in X'_5\}$. 
 
  Then 
 \begin{itemize}
 	\item $\NBHout_{\vec{E}_1}(X_1)$ corresponds to $\{i\in [\lfloor n/2 \rfloor]: \boxplus_i\in \NBH_{\HStageF}(Q_1)\}$,
 	\item $\NBHin_{\vec{E}_1}(X'_1)$ corresponds to $\{i\in [\lfloor n/2 \rfloor]: \boxminus_i\in \NBH_{\HStageF}(Q'_1)\}$,
 	\item $\NBH_{E_2}(X_2)$ corresponds to $\{i\in [\lfloor n/2 \rfloor]: \boxminus_i\in\NBH_{\HStageF}(Q_2)\}$,
 	\item $\NBH_{E_3}(X_3)$ corresponds to $\{i\in [\lfloor n/2 \rfloor]: \boxplus_i\in\NBH_{\HStageF}(Q_3)\}$,
 	\item $\NBH_{E_4}(X_4)$ corresponds to $\{i\in [\lfloor n/2 \rfloor]: \boxminus_i\notin \bigcup_{s\in X_4}\imF(s)\}$,
 	\item $\NBH_{E_5}(X_5)$ corresponds to $\{i\in [\lfloor n/2 \rfloor]: \boxplus_i\notin \bigcup_{s\in X_5}\imF(s)\}$,
 	\item $\NBH_{E_6}(X_6)$ corresponds to $\bigcap_{s\in X_6}I_s$,
 	\item $\NBH_{E_4}(X'_4)$ corresponds to $\{s\in \cJ_0: \imF (s)\cap Q'_4=\emptyset\}$,
 	\item $\NBH_{E_5}(X'_5)$ corresponds to $\{s\in \cJ_0: \imF (s)\cap Q'_5=\emptyset\}$,
 	\item $\NBH_{E_6}(X'_6)$ corresponds to $\{s\in \cJ_0: X'_6\subseteq I_s\}$.
\end{itemize}

Hence,~\eqref{eq:conjquasi-V} and~\eqref{eq:conjquasi-U} translate as
\begin{multline}
\label{eq:conjquasi-V(H)}
\bigg|\bigg\{i\in [\lfloor\tfrac{n}{2}\rfloor] \cap \bigcap_{s\in X_6}I_s: \Big(\boxminus_i\in \NBH_{\HStageF}(Q'_1\cup Q_2)\setminus \bigcup_{s\in X_4}\imF(s)\Big) \\
\text{and} \quad \Big(\boxplus_i\in \NBH_{\HStageF}(Q_1\cup Q_3)\setminus \bigcup_{s\in X_5}\imF(s)\Big)\bigg\}\bigg|\\
=(1\pm 5\LF\gamJ)d_1^{|Q_1|+|Q'_1|}\cdot\prod_{i=2}^{3} d_i^{|Q_i|}\cdot \prod_{i=4}^{6} d_i^{|X_i|}\cdot \frac{n}{2}
\end{multline}
and
\begin{equation}
\label{eq:conjquasi-J_0}
|\{s\in \cJ_0: X'_6\subseteq I_s\text{ and }\imF (s)\cap ( Q'_4\cup Q'_5)=\emptyset\}|=(1\pm 5\LF\gamJ)d_4^{|Q'_4|}d_5^{|Q'_5|}d_6^{|X'_6|}\sigmJnula n\;,
\end{equation}
with the densities $d_1,\ldots,d_6$ as computed above.

Condition~\eqref{eq:conjquasi-V(H)} follows from~\ref{quasi:indexquasi}, by which $\HStageF$ is index-quasirandom with parameters $(\LF,\gamJ,\dStageF^*,\dStageFbar)$ (see Definition~\ref{def:index-quasirandom}), on setting $S_1:=Q_1'\cup Q_2$, $S_2:= Q_1\cup Q_3$, $T_1:= X_4$, $T_2:= X_5$ and $T_3:=X_6$: this gives the required count with relative error $\gamJ$, but with $\dStageF^*$, $\dStageFbar$, $\deltnonspanning+2\sigmKJ^*$ and $2\sigmKJ^*$ in place of $d_2=d_3$, $d_1$, $d_4=d_5$ and $d_6$, respectively. Replacing each of these by the corresponding density costs a further factor $1\pm2\gamJ$ per occurrence by the computation above, and there are at most $\LF$ occurrences in total, so the relative error stays within $5\LF\gamJ$. Condition~\eqref{eq:conjquasi-J_0} follows in the same way from~\ref{quasi:imagecaps}, by setting $S_1:=Q_4'\cup Q_5'$, $S_2=\emptyset$ and $T:=X_6'$.
\end{proof}

By Lemmas~\ref{lem:chest-degree} and~\ref{lem:chest-quasirandomness}, together with the parameter bounds observed above and the fact that $5\LF\gamJ\le\gamma_{P\ref{prop:Designs-oursetting}}$ and $\LF\ge L_{P\ref{prop:Designs-oursetting}}$, Proposition~\ref{prop:Designs-oursetting} applies and gives us a diamond core-decomposition of the chest $\mathcal M$.
Next, we shall explain how we use the diamond core-decomposition given by Proposition~\ref{prop:Designs-oursetting} to pack the leftover path-forest.

For each diamond of the core-decomposition, its edge in $E_6$ tells us which path $P$ we shall pack using this diamond: if that edge lies between $s\in \cJ_0$ and $i\in [\lfloor\frac{n}{2}\rfloor]$, then $P$ is the only path in $\SpecShortPaths_s$ anchored at $\{\boxminus_i, \boxplus_i\}$. The edges of the diamond in $E_2$, $\vec{E}_1$ and $E_3$, say between $i$, $j$ and $k$, tell us where the three edges of the path $P$ are to be mapped. So if the $E_2$-edge between $i$ and $j$ ensures there is an edge in $\HStageF$ between $\boxminus_i$ and $\boxminus_j$, the $\vec{E}_1$-edge between $j$ and $k$ ensures there is an edge in $\HStageF$ between $\boxminus_j$ and $\boxplus_k$, and the $E_3$-edge between $k$ and $i$ ensures there is an edge in $\HStageF$ between $\boxplus_k$ and $\boxplus_i$. The edges from $E_4$ and $E_5$ ensure that $\boxminus_j$ and $\boxplus_k$ are not already used by $G_s$. 
	
	As we have a diamond core-decomposition, each edge of $\vec{E}_1\dcup E_2\dcup E_3$ (that is, each edge of $\HStageF$) lies in exactly one diamond and is therefore used exactly once, and each edge of $E_6$ (that is, each path of $\bigcup_{s\in \cJ_{0}}\SpecShortPaths_s$) is packed exactly once. Furthermore, each edge of $E_4\dcup E_5$ lies in at most one diamond, which is what guarantees that the resulting map is injective on each $G_s$: for a fixed $s\in\cJ_0$, no two of the packed paths of $\SpecShortPaths_s$ use the same vertex of $\Vmin$, and no two of them use the same vertex of $\Vplus$. These maps $(\phi^\mathbf{G}_s)_{s\in\cG}$ therefore extend $(\phi^\mathbf{F}_s)_{s\in\cG}$ and form a packing of $(G_s)_{s\in\cG}$ into~$H$, which completes the proof of Lemma~\ref{lem:StageG}. 

\section{Concluding remarks}
\label{sec:concl}

Theorem~\ref{thm:maintechnical} allows us to perfectly pack large degenerate graphs
of maximum degree $O(n/\log n)$, most of which may be spanning, and some of
which are non-spanning and give us quadratically many bare paths and linearly
many odd-degree vertices. As discussed, this implies the tree packing conjecture
for large~$n$ and trees with maximum degree $O(n/\log n)$. Hence, to resolve
this conjecture in full (for large~$n$), it only remains to consider the case of
trees with large degree vertices. As we explained in Section~\ref{ssec:opt}, our methods crucially need
the maximum degree assumption, witnessed by the fact that they also work for
pseudorandom graphs, in which the maximum degree assumption is sharp (up to
constants). It would be very interesting to solve this last remaining case of
the tree packing conjecture.

One could ask more generally which families of trees pack into~$K_n$. Obviously, the number of edges in such a family of trees may not exceed $\binom{n}{2}$. We further argued in Section~\ref{ssec:opt} that $\Omega(n/\log n)$ of the trees must be at least $\Omega(n/\log n)$-far from spanning. In Section~9.2 of~\cite{BHPT} it was explained that if only an aggregate bound on the maximum degrees is imposed, then this bound cannot be less than $n/2$. This motivates the following conjecture.

\begin{conjecture}
  There exist $\deltnonspanning>0$ and $n_{0}\in\NN$ such that for each $n\ge
  n_{0}$ any family of trees $\left(T_s\right)_{s\in[N]}$ with $\Delta(T_s)\le
  n/2$ and $v(T_s)\le n$ for all $s\in[N]$, with $\deltnonspanning n\le
  v(T_s)\le(1-\deltnonspanning)n$ for all $s\in[\deltnonspanning n]$, and with
  $\sum_{s\in[N]}e(T_s)\le\binom{n}{2}$, packs into $K_{n}$.
\end{conjecture}

We further remark that we use the required odd-degree vertices only for
parity corrections. In scenarios where this is not needed, our methods would
also work without them. Let us give an example. The
Oberwolfach problem on packings of cycle-factors was solved by Glock, Joos,
Kim, K\"uhn, and Osthus~\cite{GloJooKimKueOst}, and generalised by Keevash and
Staden~\cite{KS:Oberwolfach}, who proved that any large and sufficiently
quasirandom $2r$-regular graph can be decomposed into any family of~$r$
two-factors. Our methods would give a variant of this result where we do not
require the host graph to be regular but need some non-spanning $2$-regular
graphs in the family. Namely, we could pack any family of $2$-regular graphs
into any sufficiently quasirandom graph all of whose vertex degrees are even,
provided that $\deltnonspanning n$ of the graphs in the family have at least
$\deltnonspanning n$ and at most $(1-\deltnonspanning)n$ vertices and contain
$\deltnonspanning n$ cycles of length at least~12. This last restriction on
cycle lengths comes from the way we work with bare paths and could probably be removed with some extra work.

The motivation for the name `Oberwolfach problem' is that conference attendees sit around circular tables for their meals. The Oberwolfach problem is to determine when it is possible that every mathematician sits next to every other one exactly once at meals. This is exactly the problem of packing $K_n$ with $n$-vertex $2$-regular graphs. Due to the Schwarzwald scenery, however, some mathematicians go on hikes, so that at some meals the number of vertices in the corresponding $2$-factor to pack is smaller than $n$. We therefore propose the following `Oberwolfach problem with hikes' as a more practical version of the original. Given $n$ and a collection of $2$-regular graphs $G_1,\dots,G_t$ each on at most $n$ vertices with $\sum_{i=1}^te(G_i)=\binom{n}{2}$, when is it true that $G_1,\dots,G_t$ pack into $K_n$?

\smallskip

Concerning the packing of $r$-regular graphs with $r\ge 2$, Glock, Joos, Kim,
K\"uhn, and Osthus~\cite{GloJooKimKueOst} formulated the following conjecture,
which already for $r=3$ seems challenging.

\begin{conjecture}
  Given $r$, there exists $n_{0}\in\NN$ such that for each
  $n\ge n_{0}$ any family of $n$-vertex graphs $\left(G_s\right)_{s\in[N]}$ such that $G_s$ is $r_s$-regular
  with $r_s\le r$ for all $s\in[N]$ and such that
  $\sum_{s\in[N]} r_s=n-1$, packs into $K_{n}$.
\end{conjecture}

Considering, more generally, $D$-degenerate graphs, the situation is certainly
more complex. It is clear that an appropriate analogue of the tree packing
conjecture with trees replaced by $D$-degenerate graphs cannot hold: if one of
the graphs in the family is $K_{2,n-2}$, its embedding would isolate an edge in
the host graph. However, an analogue of Ringel's conjecture could still hold.

\begin{conjecture}
  Given $D$, there exists $n_{0}\in\NN$ such that for each
  $n\ge n_{0}$ the following holds. Let~$G$ be any $D$-degenerate graph on~$n$ vertices and with $m\ge n-1$ edges.
  Then $2m+1$ copies of~$G$ pack into $K_{2m+1}$.
\end{conjecture}

\section{Acknowledgements}
This work was commenced during the Southwestern German Workshop on Graph Theory in August 2018. We thank the organisers, Maria Axenovich,
Yury Person, and Dieter Rautenbach, we thank Felix Joos for his explanations concerning~\cite{GloJooKimKueOst} at this workshop,
and we (PA, JB, JH, DP) thank Pavel Hladk\'y for babysitting our children during that time. Four anonymous referees provided very detailed and helpful comments on the manuscript. Google Gemini~3.6, Claude Sonnet~5 and Claude Opus~5 were used during the revision of this paper in 2026 to improve language, and to incorporate the referees' comments; all mathematical content was verified by the authors.

\appendix
\section{Deducing Theorem~\ref{thm:DesignsIntermediate} \texorpdfstring{from~\cite{KeeColouredDirected}}{}}\label{appendix:designs}

In this appendix we deduce Theorem~\ref{thm:DesignsIntermediate} from \textcolor{blue}{Theorem~19}
of~\cite{KeeColouredDirected}. This deduction is, on the one hand, fairly pedestrian;
on the other hand, it requires understanding some of the rather complicated theory
of~\cite{KeeColouredDirected}. Our paper and~\cite{KeeColouredDirected} use
similar terms (such as ``divisibility''), but our definitions are tailored to a
more specialised setting.  In this appendix, when referring to terminology or
theorem, definition, and page numbers from~\cite{KeeColouredDirected}, we use
\textcolor{blue}{blue text}. For the convenience of those readers who are using
the arXiv version of~\cite{KeeColouredDirected} instead, we also occasionally add (in
brackets and in \textcolor{red}{red}) the corresponding theorem, definition,
and page numbers from this arXiv version~\cite{KeeColouredDirected-arXiv}.  We
remark that \cite{KeeColouredDirected} contains numerous sample applications of
\textcolor{blue}{Theorem~19} that are similar to our application here.  Two
further applications of \textcolor{blue}{Theorem~19} are used
in~\cite{KS:Oberwolfach,KS:Ringel} to deduce the main results there.

We now explain how to deduce Theorem~\ref{thm:DesignsIntermediate}.  Let $q,
D\in \mathbb{N}$ and $\sigma>0$ be given and let $q'$ be as defined in
Theorem~\ref{thm:DesignsIntermediate}. Let $n_0$ and $\omega_0$
be the output of \textcolor{blue}{Theorem~19} (\textcolor{red}{Theorem~7.4}) with input
$\textcolor{blue}{q}:=q',D$ and for uniformity $\textcolor{blue}{r}:=2$.
It is easy to check that the constants $h$ and $\delta$
given on \textcolor{blue}{page~283} (\textcolor{red}{page~4}) are consistent with those
in the statement of Theorem~\ref{thm:DesignsIntermediate}.  Observe that $\sigma
n\ge \frac{n}{h}$ by choice of~$q'$ and~$h$.

Let $n>n_0$, set $\textcolor{blue}{n_1}:=n+q'-q$, and let $\omega\in (n^{-\delta},
\omega_0)$, as is required in our Theorem~\ref{thm:DesignsIntermediate} but also
in \textcolor{blue}{Theorem~19}. Let $\mathcal P=(P_1, \ldots, P_t)$ be a
partition of $[q]$ and $\mathcal P'=(P'_1, \ldots, P'_t)$ be a partition of
$[n]$ with $|P'_i|\ge \sigma n\ge \frac{n}{h}$ for each $i\in [t]$ and trivially $|P'_i|\le n\le \textcolor{blue}{n_1}$.  We extend
$\mathcal P$ by adding $q'-q$ isolated vertices to $P_t$. Abusing notation,
we call this new partition also $\mathcal P$.

Let $\mathcal H^*$ be a simple $\mathcal P$-canonical family of
$[D]$-edge-coloured digraphs on $[q]$ supplied by Theorem~\ref{thm:DesignsIntermediate}. We extend each $H^*\in \mathcal H^*$ by adding
the $q'-q$ new vertices of $P_t$ as isolated vertices and call this new family
$\mathcal H$. This step is only necessary to ensure that each digraph $H\in\mathcal{H}$ is on
vertex set $[q']$. The reason for moving from vertex set $[q]$ to vertex set
$[q']$ in turn is that~\cite{KeeColouredDirected} requires that $h=2^{50(q')^3}$
where $[q']$ is the vertex set of the digraphs we want to decompose into, but
here we need that~$h$ is big in order to guarantee $\sigma n\ge
\frac{n}{h}$.

Let $G$ be a general $[D]$-edge-coloured multi-digraph on $[n]$ such that
$(G,\mathcal P')$
is $(\mathcal H^*,\mathcal
  P)$-divisible, $(\mathcal H^*,\omega^{h^{20}},2h^{q'-q}\omega)$-regular and $(\mathcal
H^*,\mathcal P,\sqrt[q'h]{\omega}, h)$-vertex-extendable. Note that the notion of
divisibility and vertex-extendability is not changed by adding the
$q'-q$ isolated vertices to digraphs of $\mathcal H^*$, so $(G,\mathcal P')$
is $(\mathcal H,\mathcal
  P)$-divisible and $(\mathcal
H,\mathcal P,\sqrt[q'h]{\omega}, h)$-vertex-extendable. We claim it is also $(\mathcal H,\omega^{h^{20}},\omega)$-regular.

To see this, suppose that for each coloured copy $H'$ of a graph in $\mathcal{H}^*$ in $G$ we have given a weight $w^*_{H'}$ witnessing $(\mathcal H^*,\omega^{h^{20}},2h^{q'-q}\omega)$-regularity. In other words, each $w^*_{H'}$ is in the range $[2h^{q'-q}\omega n^{2-q},(2h^{q'-q}\omega)^{-1} n^{2-q}]$, and, for any edge $e\in G$, when we sum the weights of all $H'$ (over all choices of graph in $\mathcal{H}^*$) using $e$ we obtain $1\pm\omega^{h^{20}}$. We define the falling factorial $K=\big(|P'_t|-|P_t|+q'-q\big)_{q'-q}$, which is the number of ways to extend a given coloured copy $H'$ in $G$ of some graph $H^*\in \mathcal{H}^*$ to a copy of the corresponding $H\in \mathcal{H}$ (observe that by construction this number is independent of the choice of $H^*$ and of the specific coloured copy $H'$), with the $q'-q$ added vertices being mapped to $P'_t$. Now let $H^\dagger$ be a coloured copy in $G$ of $H\in\mathcal{H}$ (we use this notation for coloured copies of graphs of $\mathcal H$ throughout this appendix), and let $H'$ be the induced coloured copy in $G$ of the corresponding $H^*\in\mathcal{H}^*$. We define $w_{H^\dagger}:=K^{-1}w^*_{H'}$. We claim that this collection of weights witnesses that $(G,\mathcal{P}')$ is $(\mathcal H,\omega^{h^{20}},\omega)$-regular. Indeed, for every edge $e\in E(G)$ the sum of weights $w_{H^\dagger}$ over copies $H^\dagger$ of graphs of $\mathcal{H}$ using $e$ is exactly equal to the sum of weights $w^*_{H'}$ over copies $H'$ of graphs of $\mathcal{H}^*$ using $e$, and in particular it is $1\pm \omega^{h^{20}}$ as required. Furthermore, for sufficiently large $n$ we have $\tfrac12\sigma^{q'-q}n^{q'-q}<K<n^{q'-q}$ because $\sigma n\le |P'_t|\le n$. It follows that for each $H^\dagger$ we have
\[\omega n^{2-q'}\le (n^{q'-q})^{-1}\cdot 2h^{q'-q}\omega n^{2-q}\le w_{H^\dagger}\le (\tfrac12\sigma^{q'-q}n^{q'-q})^{-1}\cdot (2h^{q'-q}\omega)^{-1} n^{2-q}\le\omega^{-1}n^{2-q'}\]
where the outermost inequalities use $1\ge \sigma\ge h^{-1}$. This verifies the range for $w_{H^\dagger}$ required for Definition~\ref{def:regularity(designs)}, completing the check that $(G,\mathcal{P}')$ is $(\mathcal H,\omega^{h^{20}},\omega)$-regular.

From this point on we aim to find an $\mathcal{H}$-decomposition of $G$. Since the graphs of $\mathcal H$ differ from the corresponding graphs of $\mathcal H^*$ only by isolated vertices, such a decomposition is at the same time an $\mathcal{H}^*$-decomposition of $G$, which is what Theorem~\ref{thm:DesignsIntermediate} asserts.

We shall check that $G$ satisfies the conditions of
\textcolor{blue}{Theorem~19}. First of all, this theorem requires~$\mathcal H$
to be canonical in the sense of \textcolor{blue}{Definition~35}
(\textcolor{red}{Definition~7.1}). Observe that without loss of generality, we may
assume that every edge in any graph from $\mathcal H$ between two parts~$P_i$
and~$P_j$ with $i<j$ is oriented towards~$P_j$. (If this is not the case, that
is, there is some colour~$d$ with colour location $(j,i)$, then simply reverse
all the edges of colour~$d$ in~$\mathcal H$ and~$G$, and accordingly change the
colour location of~$d$ to $(i,j)$.) This assumption is necessary so that we are
consistent with \textcolor{blue}{Definition~35}.  This definition now requires
us to specify an index vector $\textcolor{blue}{\mathbf {i}^d}\in \NN_0^t$ for
each colour~$d$ such that $\textcolor{blue}{\mathbf {i}^d}$ encodes the colour
location of $d$ as follows. Let~$d$ be a colour with colour location $(i,j)$. If
$i=j$, that is, edges of colour~$d$ lie within $P_i$, then we let
$\textcolor{blue}{\mathbf {i}^d}$ be the vector with a~$2$ at coordinate $i=j$
and~$0$ everywhere else. If $i\neq j$, on the other hand, that is, all edges of
colour $d$ run between $P_i$ and $P_j$, then we let $\textcolor{blue}{\mathbf
  {i}^d}$ be the vector with a~$1$ at position~$i$ and a~$1$ at position~$j$,
and~$0$ everywhere else. \textcolor{blue}{Definition~35} then defines a
partition
\textcolor{blue}{$R(\mathbf{i}^d)=(R(\mathbf{i}^d)_1,\dots,R(\mathbf{i}^d)_t)$}
of $[r]=[2]$ for each of these vectors
$\textcolor{blue}{\mathbf {i}^d}$. In our case, if the colour location $(i,j)$
of~$d$ satisfies $i=j$, then we get
$\textcolor{blue}{R(\mathbf{i}^d)_i}=\{1,2\}$ and
$\textcolor{blue}{R(\mathbf{i}^d)_{i'}}=\emptyset$ for all $i'\neq i$; if $i<j$,
then we get $\textcolor{blue}{R(\mathbf{i}^d)_i}=\{1\}$,
$\textcolor{blue}{R(\mathbf{i}^d)_j}=\{2\}$ and
$\textcolor{blue}{R(\mathbf{i}^d)_{i'}}=\emptyset$ for all $i'\neq i,j$.
Further, \textcolor{blue}{Definition~35} requires us to specify for each
colour~$d$ and each $j'\in[t]$ a permutation group
$\textcolor{blue}{\Lambda_{j'}^d}$ on $\textcolor{blue}{R(\mathbf{i}^d)_{j'}}$.
The purpose of these permutation groups is to allow for a mixture of
(generalisations of) partly directed and undirected edges within parts; for us they will all be
either trivial or the symmetric group $\mathbb{S}_2$ on elements $\{1,2\}$.  Indeed, we
let $\textcolor{blue}{\Lambda_{j'}^d}$ be the trivial group $\operatorname{Id}$
on $\textcolor{blue}{R(\mathbf{i}^d)}_{j'}$ only containing the identity
whenever $\textcolor{blue}{R(\mathbf{i}^d)}_{j'}$ is empty, contains only one
element, or if~$d$ is an oriented colour. If~$d$ is an unoriented colour with
colour location $(j',j')$, on the other hand, then we let
$\textcolor{blue}{\Lambda_{j'}^d}$ be $\mathbb{S}_2$.  As in 
\textcolor{blue}{Definition~35}, we set
$\textcolor{blue}{\Lambda^d}=\prod_{j'}\textcolor{blue}{\Lambda_{j'}^d}$ and
$\textcolor{blue}{\Lambda}:=(\textcolor{blue}{\Lambda^d}:d\in [D])$.

We now claim that $\mathcal H$ is \textcolor{blue}{$(\mathcal
  P,\Lambda)$-canonical} as defined in \textcolor{blue}{Definition~35}.  Indeed,
\textcolor{blue}{Definition~35i} just says that for each colour~$d$ edges
of colour~$d$ respect the colour location
of~$d$. \textcolor{blue}{Definition~35ii} translates to the following: For each
$H\in\mathcal{H}$ individually,
\begin{enumerate}[label=\abc]
\item (``\textcolor{blue}{for $\theta'\in Bij([r],Im(\theta))$ we have
  $\theta'\notin H\setminus H^d$}'') $H$ does not have parallel or antiparallel
  edges of different colours, and
\item(``\textcolor{blue}{$\theta'\in H^d$ iff $\theta^{-1}\theta'\in
  \Lambda^d$}'') that edges of colour~$d$ with colour location $(j,j)$ come in
  pairs forming directed $2$-cycles if and only if the corresponding permutation
  group $\textcolor{blue}{\Lambda^d_j}$ equals $\mathbb S_2$.
\end{enumerate}
Hence, \textcolor{blue}{Definition~35} is consistent with our
Definition~\ref{def:canonical}: loops are not allowed (see the definition of an $r$-digraph in \textcolor{blue}{Definition~28} (\textcolor{red}{Definition 6.1})) and in Definition~\ref{def:canonical} we consider the special case of simple digraphs, i.e., we forbid parallel edges of the same colour, and also antiparallel edges of any given oriented colour (which would form $2$-cycles).

Let us now turn to the interpretation of \textcolor{blue}{Definition~36} (\textcolor{red}{Definition~7.2}). This definition involves the set $\textcolor{blue}{I^r_{q'}}$ of injective maps from $[r]$ to $[q']$. In~\cite{KeeColouredDirected}, such injections are referred to as \textcolor{blue}{edges}. To avoid confusion with other uses of the term ``edge'', we refer to such injections as \emph{$r$-tuples}, and more generally we use the term \emph{$r$-tuple} for any injective map from $[r]$ to a given ground set. 
So, $\textcolor{blue}{I^r_{q'}}$ is the set of all
possible $r$-tuples on $[q']$.
Here, we work with $r=2$. Further,
\textcolor{blue}{Definition~36} works with a vector formalism in which each
$H\in \mathcal{H}$ is encoded as $H\in (\NN_0^D)^{\smash{I^2_{q'}}}$ in the
following way. The element $H_{d,f}$ with $d\in [D]$ and $f\in I^2_{q'}$ equals
$1$ if there is an edge of colour $d$ starting in vertex $f(1)$ and ending in
vertex $f(2)$, and equals~$0$ otherwise. The group~$\textcolor{blue}{\Sigma}$
in \textcolor{blue}{Definition~36} is the permutation group on $[q']$ containing
all permutations~$\sigma$ that preserve the partition~$\mathcal P$,
that is, $\sigma(P_i)=P_i$ for all $i\in[t]$.
The \textcolor{blue}{$\Sigma$-adapted $[q']$-complex $\Phi$} on vertex set $[n]$
is used to encode where copies of $H\in\mathcal H$ are potentially allowed to be
mapped in a graph on~$[n]$. For us this just means that the partitions~$\mathcal
P$ and~$\mathcal P'$ have to be respected. We let \textcolor{blue}{$\Phi$} be
the complete $\mathcal P$-partite $[q']$-complex on~$[n]$ with
partition~$\mathcal P'$, that is, \textcolor{blue}{$\Phi$} firstly contains
every tuple $(v_1,\dots,v_{q'})$ with pairwise distinct vertices
such that $v_i\in P'_j$ whenever $i\in P_j$; and secondly contains all tuples in the down-closure of any such tuple $(v_1,\dots,v_{q'})$. Note that $\Phi$ is in fact \textcolor{blue}{exactly $\Sigma$-adapted}: given $B,B'\subseteq[q']$, a map $\varphi\in\Phi_B$ and a bijection $\tau$ from $B'$ to $B$, we have $\varphi\circ\tau\in\Phi_{B'}$ if and only if $\tau$ maps $B'\cap P_j$ into $P_j$ for each $j\in[t]$, which in turn holds if and only if $\tau$ is the restriction of a permutation of $[q']$ preserving~$\mathcal P$. The
(directed) $q'$-graph \textcolor{blue}{$\Phi_{q'}$} then consists of all the
$q'$-tuples of \textcolor{blue}{$\Phi$}. Similarly, \textcolor{blue}{$\Phi_{2}$}
contains all underlying $2$-tuples, that is, all $2$-tuples
within each part~$P'_i$ and all $2$-tuples between any two parts~$P'_i$
and~$P'_j$ with $i<j$ that are directed towards~$P'_j$. Further,
\textcolor{blue}{$\mathcal H(\Phi)$} denotes the set of coloured labelled
copies~$H^\dagger$ on~$[n]$ of digraphs $H\in\mathcal H$ corresponding to
$q'$-tuples $\phi\in\Phi_{q'}$ (such a copy is written $\textcolor{blue}{\phi H}$ in~\cite{KeeColouredDirected}), that is, vertices of~$H$ from~$P_i$ are mapped to vertices
from~$P'_i$ in~$H^\dagger$ (and edges of~$H$ are mapped to $2$-tuples in $\Phi_2$).
The host graph~$G$ is then defined as a coloured multi-digraph, in which every
$2$-tuple from $\Phi_{2}$ may appear as an edge multiple times and in different colours,
which is consistent with our Definition~\ref{def:general-digraph}.
 
So now it remains to verify the three conditions of \textcolor{blue}{Theorem 19},
namely that~$G$ is \textcolor{blue}{$\mathcal H$-divisible} in~$\Phi$ and
\textcolor{blue}{$(\mathcal H, c, \omega)$-regular} in~$\Phi$, and that all
$(\Phi,G^H)$ are \textcolor{blue}{$(\omega,h)$-extendable}. These definitions
are related to Definitions~\ref{def:divisibility},
\ref{def:regularity(designs)}, and~\ref{def:vertex-extensibility}, but we need
to make the correspondence precise.

For \textcolor{blue}{$\mathcal H$-divisibility} we first explain the notion of
the \textcolor{blue}{degree vector $G(\psi)^*\in \NN_0^{[D]\times I^i_2}$}. In
our setting, \textcolor{blue}{$\psi$} is an injective map from $[i]$ to $[n]$,
where $i=0,1,2$. For $i=0$, the only choice is the trivial one,
$\textcolor{blue}{\psi}=\emptyset_n$, and also \textcolor{blue}{$I^0_2$} only
contains the trivial map $\emptyset_2$.  Hence,
$\textcolor{blue}{G(\emptyset_n)^*}$ has only one coordinate
$\textcolor{blue}{G(\emptyset_n)^*}_{d,\emptyset_2}$ for each colour $d\in[D]$,
which equals the total number of edges of colour $d$ in $G$. For $i=1$, a given
\textcolor{blue}{$\psi$} corresponds to a choice of a vertex $v\in [n]$ and
$I^1_2$ contains the two elements $1\mapsto 1$ and $1\mapsto 2$. So in this
case, \textcolor{blue}{$G(v)^*$} is a list of $2D$ numbers encoding the in- and
out-degree of each colour at vertex $v$ in the~$D$ colours. For $i=2$, a given
\textcolor{blue}{$\psi$} corresponds to a choice of a $2$-tuple $uv$ of
vertices of~$[n]$ and $I^2_2$ contains the two possible permutations on~$2$ elements.
In this case, \textcolor{blue}{$G(uv)^*$} is a list of $2D$
numbers encoding the number of edges going between $u$ and $v$ in each colour
and each of the two directions. The degree vector \textcolor{blue}{$H(\theta)^*$} is defined
analogously for $H\in \mathcal H$ and $\theta$ an injective map from $[i]$
to $[q']$ with $i=0,1,2$.

Let us clarify one minor point in the definition of \textcolor{blue}{$\mathcal
  H$-divisibility} in \textcolor{blue}{Definition~36}:
\begin{itemize}[leftmargin=*]
\item
\textcolor{blue}{
For $\mathbf{i'} \in \NN_0^t$ we let 
$H\langle\mathbf{i'}\rangle = \langle H(\theta)^* : i_{\mathcal{P}}(\theta)=\mathbf{i'}\rangle$.
We say $G$ is $\mathcal{H}$-divisible (in $\Phi$) if
$G(\psi)^* \in H\langle \mathbf{i'} \rangle$ whenever $i_{\mathcal{P}'}(\psi)=\mathbf{i'}$.
}
\end{itemize}

\noindent
should read

\begin{itemize}[leftmargin=*]
\item
\textcolor{blue}{
For $\mathbf{i'} \in \NN_0^t$ we let 
$\mathcal{H}\langle\mathbf{i'}\rangle = \langle H(\theta)^* : H\in\mathcal{H}, i_{\mathcal{P}}(\theta)=\mathbf{i'}\rangle$.
We say $G$ is $\mathcal{H}$-divisible (in $\Phi$) if
$G(\psi)^* \in \mathcal{H}\langle \mathbf{i'} \rangle$ whenever $i_{\mathcal{P}'}(\psi)=\mathbf{i'}$.
}
\end{itemize}

What this says is the following.  For a fixed index vector
\textcolor{blue}{$\mathbf{i'} \in \NN_0^t$}, specifying how many vertices we want
to fix in each part of our partition, we let~$i$ be the sum $\sum\textcolor{blue}{\mathbf{i'}}$ of the entries of
\textcolor{blue}{$\mathbf{i'}$}.  We then consider all injective maps~$\theta$
from $[i]$ to $[q']$ with this index vector, that is, $\theta$ sends exactly
\textcolor{blue}{$\mathbf{i'}_j$} vertices to $P_j$. Each of these~$\theta$
(specifying which vertices we are fixing) and each $H\in\mathcal H$ then gives
us a degree vector $H(\theta)^*$. The lattice generated by all these vectors is
denoted by \textcolor{blue}{$\mathcal{H}\langle\mathbf{i'}\rangle$}. The
graph~$G$ is then \textcolor{blue}{$\mathcal{H}$-divisible} if for any mapping
$\psi$ from some $[i]$ to $[n]$ (fixing $i$ vertices of~$G$) the following
holds. We let $\textcolor{blue}{\mathbf{i'}}\in\NN_0^t$ be the index vector of
$\psi$. We then require that the degree vector \textcolor{blue}{$G(\psi)^*$} is
in \textcolor{blue}{$\mathcal{H}\langle\mathbf{i'}\rangle$}.

Accordingly, this definition has to be checked for each
``$\textcolor{blue}{\mathbf{i'}\in\NN_0^t}$'' and each \textcolor{blue}{$\psi$}
with ``$\textcolor{blue}{i_{\mathcal P'}(\psi)=\mathbf{i'}}$''. For us, there are three
cases to distinguish, depending on whether $\sum\textcolor{blue}{\mathbf{i'}}$
is 0, 1, or 2; for higher values of $\sum\textcolor{blue}{\mathbf{i'}}$ the
condition in \textcolor{blue}{Definition~36} is void. We now show that these
three cases correspond to 0-divisibility, 1-divisibility, and 2-divisibility in
Definition~\ref{def:divisibility}, respectively.
\begin{itemize}[leftmargin=*]
\item First suppose that $\sum\textcolor{blue}{\mathbf{i'}}=0$. Then the
  condition requires that \textcolor{blue}{$G(\emptyset_n)^*$} is a linear
  combination of the degree vectors $\textcolor{blue}{H(\emptyset_{q'})^*}$,
  $H\in \mathcal H$, with integer coefficients $m_H$.  This means exactly that
  the number of edges of colour~$d$ in~$G$ is $\sum_{H\in\mathcal H}m_H
  c_{d,H}$, where $c_{d,H}$ is the number of edges of colour~$d$ in~$H$, which is 0-divisibility.
\item If $\sum\textcolor{blue}{\mathbf{i'}}=1$, the condition mandates the
  following for each choice of vertex $v\in [n]$. Let $j\in[t]$ be such that
  $v\in P'_j$.  The degree vector $\textcolor{blue}{G(v)^*}$ has to be a linear
  combination with integer coefficients of the degree vectors $H(\theta)^*$,
  with $H\in \mathcal H$ and $\theta:[1]\rightarrow [q'] $ such that
  $\theta(1)\in P_j$. This is exactly 1-divisibility.
\item If $\sum\textcolor{blue}{\mathbf{i'}}=2$, then the condition asks that for
  each choice of a $2$-tuple $uv$ from $[n]$ the following holds. Let $j,j'\in [t]$
  be such that $u\in P'_j$, $v\in P'_{j'}$. The degree vector
  $\textcolor{blue}{G(uv)^*}$ has to be a linear combination with integer
  coefficients of the degree vectors $H(\theta)^*$ with $H\in \mathcal H$ and
  $\theta:[2]\rightarrow [q']$ such that $\theta(1)\in P_j$ and $\theta(2)\in
  P_{j'}$. Because our graphs $H\in\mathcal{H}$ are simple, the entries of the degree vector $H(\theta)^*$ are always $0$ or $1$. This means that an edge of colour~$d$ can appear in~$G$ between
  $P'_j$ and $P'_{j'}$ only if there is an edge of colour~$d$ between $P_j$ and
  $P_{j'}$ in some $H\in\mathcal{H}$. Further, since for unoriented colours~$d$
  with colour location $(j,j)$ the edges of colour~$d$ in each $H\in\mathcal H$
  come in pairs that form directed $2$-cycles, the same has to be true for
  edges of colour~$d$ in~$G$. This is exactly 2-divisibility.
\end{itemize}

Let us now turn to \textcolor{blue}{$(\mathcal H, c, \omega)$-regularity}, which
requires that for each $H\in\mathcal H$ and each coloured copy
\textcolor{blue}{$\phi H$} of~$H$ in~$G$ with
\textcolor{blue}{$\phi\in\Phi_{q'}$} (that is, $\phi$ maps vertices from $P_j$
to $P'_j$; in the notation introduced above, $\textcolor{blue}{\phi H}$ is the copy $H^\dagger$ determined by~$\phi$) we can define a weight \textcolor{blue}{$y_\phi^H\in[\omega
    n^{2-q'},\omega^{-1} n^{2-q'}]$} so that for each edge $e\in E(G)$, when we
sum these weights over all copies \textcolor{blue}{$\phi H$} which use~$e$ we
get $(1\pm c)$.
 Since in our setting, any
coloured copy of any $H\in \mathcal H$ maps vertices from~$P_j$ to~$P'_j$, this
is identical to Definition~\ref{def:regularity(designs)}. 

The last notion is that of \textcolor{blue}{extendability} for all $H\in\mathcal
H$. By a remark before \textcolor{blue}{Definition~37} (\textcolor{red}{Definition 7.3}), the required
\textcolor{blue}{$(\omega,h)$-extendability} of
\textcolor{blue}{$(\Phi,G^H)$} follows from
\textcolor{blue}{$(\sqrt[q'h]{\omega},h)$-vertex-extendability} of
\textcolor{blue}{$(\Phi,G)$}. According to \textcolor{blue}{Definition~37}
\textcolor{blue}{$(\sqrt[q'h]{\omega},h)$-vertex-extendability} of
\textcolor{blue}{$(\Phi,G)$} mandates the
following.
For every vertex $x\in [q']$ and for any choice of vertex-disjoint sets
$\textcolor{blue}{A_i}$ with $i\in[q']\setminus\{x\}$ of size at most~$h$ in $G$
such that for every collection of $v_i\in A_i$ we have that
\textcolor{blue}{$(i\mapsto v_i:i\in[q']\setminus\{x\})\in\Phi$} (which for us
just means that $v_i\in P'_j$ where~$j$ is such that $i\in P_j$, so overall this condition simply means that $A_i\subset P'_j$), there have to
be at least $\sqrt[q'h]{\omega}n$ vertices $v\in\textcolor{blue}{\Phi^\circ_x}$
(for us $\textcolor{blue}{\Phi^\circ_x}$ is the part $P'_{j'}$ such that $x\in
P_{j'}$) satisfying the following properties. Firstly,
\textcolor{blue}{$(i\mapsto v_i:i\in[q'])\in\Phi$ whenever $v_x=v$ and $v_i\in
  A_i$ for $i\neq x$}, which, since $v\in\textcolor{blue}{\Phi^\circ_x}$, is always
satisfied. Secondly, for every colour~$d$ and every
edge~\textcolor{blue}{$\theta$} of~$H$ which is of colour~$d$ and contains~$x$,
we have that all $2$-tuples \textcolor{blue}{$(i\mapsto v_i:i\in[2])$} are edges
of~$G$ in colour~$d$ whenever \textcolor{blue}{$v_j=v$ for $j=\theta^{-1}(x)$
  and $v_i\in A_{\theta(i)}$ for all $i\neq j$}.  This means that for any edge
$xy$ in~$H$ of colour~$d$ we have that all edges $vw$ with $w\in A_y$ of
colour~$d$ are present in~$G$; and that for any edge $yx$ in~$H$ of colour~$d$
we have that all edges $wv$ with $w\in A_y$ of colour~$d$ are present
in~$G$.  This corresponds exactly to Definition~\ref{def:vertex-extensibility}.

Hence the host graph $G$ from Theorem~\ref{thm:DesignsIntermediate} satisfies the
conditions of \textcolor{blue}{Theorem~19} and thus $G$ has an $\mathcal
H$-decomposition as claimed.

\printindex

\bibliographystyle{amsplain_yk}
\bibliography{packing}	 
	
\end{document}